\documentclass[reqno,table]{amsart}
 \usepackage[utf8]{inputenc}
\usepackage[table]{xcolor}
\usepackage{fancyhdr}
\usepackage[backref=page,bookmarksopen,bookmarksdepth=2,colorlinks=true]{hyperref}
\definecolor{darkergreen}{rgb}{0.0, 0.8, 0.0}  
\hypersetup{
    linkcolor=blue,      
    citecolor=magenta,      
    urlcolor=darkergreen,       
    }
\usepackage{verbatim}
\usepackage{graphicx}
\usepackage{accents}
\usepackage{subcaption}
\usepackage{amssymb,amsmath,amsfonts,amsthm,epsfig,amscd,stmaryrd,mathrsfs}
\usepackage{stmaryrd}
\usepackage[all,cmtip,poly]{xy}
\usepackage{color} 
\usepackage{enumitem} 
\def\QQQQ{\mathrm{Den}}
\def\MMM{\mathcal{M}}

\newcommand{\iso}{\xrightarrow{\simeq}}

\usepackage{ifthen}

\newboolean{commenton}
\setboolean{commenton}{false} 
\newcommand{\silentcomment}[1]{
  \ifthenelse{\boolean{commenton}}{\textcolor{red}{#1}}{\ignorespaces}
}

\newenvironment{quietcomment}{
  \ifthenelse{\boolean{commenton}}{\color{red}}
     {\expandafter\comment\expandafter}}{
  \ifthenelse{\boolean{commenton}}
    {}
    {\expandafter\endcomment}}

\def\tareesidedbox#1{\setbox0=\hbox{$#1$}\dimen0=\wd0 \advance\dimen0 by3pt\rlap{\hbox{\vrule height9pt width.4pt depth2pt \kern-.4pt\vrule height9.4pt width\dimen0 depth-9pt\kern-.4pt \vrule height9pt width.4pt depth2pt}} \relax \hbox to\dimen0{\hss$#1$\hss}}
\def\ho#1{\tareesidedbox{#1}}

\SetSymbolFont{stmry}{bold}{U}{stmry}{b}{n}

\def\NwithzeroB{\mathbf{N}}
\def\NwithzeroA{\NwithzeroB}
\def\NwithoutzeroA{\mathbf{N}_{> 0}}

\usepackage{silence}
\def\Slit{\mathrm{Slit}}
\def\Res{\mathrm{Res}}

\def\AA{\widehat{A}}
\def\BB{\widehat{B}}

\def\LL{\mathcal{L}}

\def\GGG{\psi}
\def\HHH{\mathrm{Gob}}

\def\hT{\widehat{T}}
\def\hA{\widehat{A}}
\def\wV{\widetilde{\mathcal{V}}}
\def\A{\mathbb{A}}
\def\oA{\accentset{\circ}{A}}
\def\oT{\accentset{\circ}{T}}

\def\hatG{\widehat{G}}

 \def\cR{\mathcal{R}}
 \def\cQ{\mathcal{Q}}

\renewcommand\setminus{\smallsetminus}
\def\spec{\mathrm{Spec} \,}
\def\spf{\mathrm{Spf} \,}

\def\li{\mathrm{Li}}
\def\HH{\mathcal{H}}
\def\im{\mathrm{im}}

\def\chip{\chi_{-3}}
\def\Li{\mathrm{Li}}
\def\rank{\mathrm{rank} \,}
\def\Gr{\mathrm{Gr}}

\def\den{\mathrm{den}}

\usepackage{graphicx}

\newmuskip\pFqmuskip

\newcommand*\pFq[6][8]{
  \begingroup 
  \pFqmuskip=#1mu\relax
  \mathcode`\,=\string"8000
  \begingroup\lccode`\~=`\,
  \lowercase{\endgroup\let~}\pFqcomma
  {}_{#2}F_{#3}{\left[\genfrac..{0pt}{}{#4}{#5};#6\right]}
  \endgroup
}
\newcommand{\pFqcomma}{\mskip\pFqmuskip}

\def\mv{\mu_{\mathrm{Haar}}}
\def\PP{\mathbb{P}}

\def\T{\mathbf{T}}

\def\F{\mathbf{F}}
\def\Z{\mathbf{Z}}

\def\R{\mathbf{R}}

\def\OL{\mathcal{O}}

\def\Q{\mathbf{Q}}

\def\Qbar{\overline{\Q}}

\def\D{\mathbf{D}}
\def\Db{\overline{\mathbf{D}}}
\def\Cspan{\mathrm{Span} }
\def\ardeg{\widehat{\deg} \,}

\def\bx{\mathbf{x}}
\def\bk{\mathbf{k}}

\def\bz{\mathbf{z}}

\def\bs{\mathbf{s}}

\newif\iffinalrun
\iffinalrun
\else
\fi

\iffinalrun
  \newcommand{\need}[1]{}
  \newcommand{\mar}[1]{}
\else
  \newcounter{margcount}
  \newcommand{\need}[1]{{\tiny *** #1}}
  \newcommand{\mar}[1]{\stepcounter{margcount}\marginpar{\raggedright\tiny \themargcount.FIXME #1}}
\fi

\renewcommand\mathbb{\mathbf}

\numberwithin{equation}{subsection}
\def\numequation{\addtocounter{subsubsection}{1}\begin{equation}}

\counterwithin{figure}{equation}

\makeatletter\let\c@equation\c@subsubsection\makeatother

\newtheorem{theorem}[equation]{Theorem}
\newtheorem{cor}[equation]{Corollary}
\newtheorem{proposition}[equation]{Proposition}
\newtheorem{fact}[equation]{Fact}

\newtheorem{problem}[equation]{Problem}
\newtheorem{conjecture}[equation]{Conjecture}
\newtheorem{question}[equation]{Question}

\newtheorem{thm}[equation]{Theorem}
\newtheorem{lemma}[equation]{Lemma}

\newtheorem{iprob}{Problem}

\newtheorem{ithm}[iprob]{Theorem}

\newtheorem{icor}[iprob]{Corollary}
\newtheorem*{BombieriInequality}{Bombieri's Inequality}

\theoremstyle{definition}
\newtheorem{remark}[equation]{Remark}
\newtheorem{basicremark}[equation]{Basic Remark}
\newtheorem{example}[equation]{Example}

\newcommand\xqed[1]{
  \leavevmode\unskip\penalty9999 \hbox{}\nobreak\hfill
  \quad\hbox{#1}}
\def\endofremark{\xqed{$\triangle$}}

\newtheorem{df}[equation]{Definition}

\usepackage{tikz}
\usetikzlibrary{decorations.markings}
\usetikzlibrary{positioning}

\def\C{\mathbf{C}}

\def\Gal{\mathrm{Gal}}

\def\SL{\mathrm{SL}}

\def\Sym{\mathrm{Sym}}

\def\vol{\mathrm{vol}}

\newmuskip\pFqmuskip

\usepackage{tocbasic}
\makeatletter
\def\l@figure{\@dottedtocline{1}{1.5em}{3em}}
\makeatother

\def\H{\mathbf{H}}
\def\ardeg{\widehat{\deg} \,}

\def\cC{\mathcal{C}}

\def\cX{\mathcal{X}}
\def\cL{\mathcal{L}}
\def\ovcL{\overline{\mathcal{L}}}

\def\bn{\mathbf{n}}
\def\bi{\mathbf{i}}
\def\bt{\mathbf{t}}
\def\be{\mathbf{e}}
\def\bb{\mathbf{b}}
\def\P{\mathbf{P}}
\def\rk{\mathrm{rk}}
\def\slope{\widehat{\mu}}
\def\Spf{\mathrm{Spf} \,}
\def\bzo{\mathbf{0}}
\def\bm{\mathbf{m}}
\def\ord{\mathrm{ord}}
\def\ovE{\overline{E}}
\def\br{\mathbf{r}}
\def\ovtau{\tau^{\flat\kern-1.2pt\flat}}

\ifdefined\NewCommandCopy
  \NewCommandCopy\latexbibitem\bibitem
\else
  \let\latexbibitem\bibitem
  \usepackage{xparse}
\fi

\renewcommand{\bibitem}[2][]{
  \def\mykey{#1}
  \def\cmpkeyA{Gau21}
  \def\cmpkeyB{Eul35}
  \def\cmpkeyC{Jac59}
  \def\cmpkeyD{Lin82}
  \def\cmpkeyE{Leg94}
  \def\cmpkeyF{Dir37}
  \def\cmpkeyG{Wei85}
  \def\cmpkeyH{Nie09}
  \def\cmpkeyI{Cat82}
  \def\cmpkeyJ{Her74}
  \def\cmpkeyK{Her93}
  \def\cmpkeyL{Her17}
  \def\cmpkeyM{Bor14}
  \def\cmpkeyN{Ang19}
  \def\cmpkeyO{Sie21}
  \def\cmpkeyP{P{\'o}l23}
    \def\cmpkeyQ{Had99}
    \def\cmpkeyR{Mai06}
  \ifx\mykey\cmpkeyA
    \latexbibitem[Gau1821]{#2}
  \else\ifx\mykey\cmpkeyB
    \latexbibitem[Eul1735]{#2}
  \else\ifx\mykey\cmpkeyC
    \latexbibitem[Jac1859]{#2}
  \else\ifx\mykey\cmpkeyD
    \latexbibitem[Lin1882]{#2}
  \else\ifx\mykey\cmpkeyE
    \latexbibitem[Leg1794]{#2}
  \else\ifx\mykey\cmpkeyF
    \latexbibitem[Dir1837]{#2}
  \else\ifx\mykey\cmpkeyG
    \latexbibitem[Wei1885]{#2}
  \else\ifx\mykey\cmpkeyH
    \latexbibitem[Nie1909]{#2}
  \else\ifx\mykey\cmpkeyI
    \latexbibitem[Cat1882]{#2}
  \else\ifx\mykey\cmpkeyJ
    \latexbibitem[Her1874]{#2}
  \else\ifx\mykey\cmpkeyK
    \latexbibitem[Her1893]{#2}
  \else\ifx\mykey\cmpkeyL
    \latexbibitem[Her1917]{#2}
  \else\ifx\mykey\cmpkeyM
    \latexbibitem[Bor1914]{#2}
  \else\ifx\mykey\cmpkeyN
    \latexbibitem[Ang1919]{#2}
  \else\ifx\mykey\cmpkeyO
    \latexbibitem[Sie1921]{#2}
      \else\ifx\mykey\cmpkeyP
    \latexbibitem[P{\'o}l1923]{#2}
         \else\ifx\mykey\cmpkeyQ
    \latexbibitem[Had1899]{#2}
             \else\ifx\mykey\cmpkeyR
    \latexbibitem[Mai1906]{#2}
  \else
    \latexbibitem[#1]{#2}
  \fi\fi\fi\fi\fi\fi\fi\fi\fi\fi\fi\fi\fi\fi\fi\fi\fi
}

\title{The linear independence of $1$, $\zeta(2)$, and $L(2,\chi_{-3})$}

\author[F. Calegari]{Frank Calegari}
 \email{fcale@math.uchicago.edu} \address{The University of Chicago,
5734 S University Ave,
Chicago, IL 60637, USA}

\author[V. Dimitrov]{Vesselin Dimitrov}
 \email{dimitrov@caltech.edu} \address{Department of Mathematics, 
 California Institute of Technology,
Pasadena, CA 91125, USA}

\author[Y. Tang]{Yunqing Tang}
\email{yungqing.tang@berkeley.edu} \address{Department of Mathematics, University of California, Berkeley, Evans Hall, Berkeley, CA 94720, USA}

\begin{document}

\begin{abstract}
We prove the irrationality of the classical Dirichlet $L$-value 
$$
L(2,\chi_{-3}) = \frac{1}{1^2} - \frac{1}{2^2} + \frac{1}{4^2} - \frac{1}{5^2} + \frac{1}{7^2} - \frac{1}{8^2}  +  \ldots.
$$ 
The argument applies a new kind of arithmetic holonomy bound to a well-known construction of Zagier~\cite{Zagier}. In fact our work  also
establishes the $\Q$-linear independence
of $1, \zeta(2)$ and $L(2,\chi_{-3})$.
We also give a number of other applications of our method
to other problems in irrationality.
\end{abstract}

\maketitle

{\footnotesize
\setcounter{tocdepth}{2}
\tableofcontents
}

\newpage

\section{Introduction}
\label{intro}

\subsection{Dirichlet~\texorpdfstring{$L$}{L}-values}
The values of the Riemann zeta function~$\zeta(k)$  for positive integers~$k$, and
more generally the Dirichlet~$L$-values
$$L(k,\chi) = \sum_{n=1}^{\infty} \frac{\chi(n)}{n^k}$$
  for 
quadratic  characters~$\chi$, 
 have long been  a source of interest to mathematicians. 
Suppose that~$\chi$ is a primitive quadratic character of conductor~$D$, where we use the convention
that the sign of~$D$ is the sign of~$\chi(-1)$.
Starting with work of
Euler~\cite{Euler} and Dirichlet~\cite{Dirichlet} (or even far before that in the
special case of~$D=-4$ and~$k=1$~\cite{MR1081274}), we know that, for positive integers~$k$:
\begin{equation}
\label{dirichlet}
L(k,\chi) \in \begin{cases}
\pi^k \cdot \sqrt{D} \cdot \Q^{\times}, & \text{$k$ even and~$\chi(-1) = 1$}, \\
\pi^k \cdot \sqrt{-D} \cdot \Q^{\times}, & \text{$k$ odd and~$\chi(-1) = -1$}, \\
\sqrt{D} \log \left(| \Qbar^{\times} | \setminus \{1\} \right), & \text{$k = 1$, $\chi(-1)  = 1$, and~$D \ne 1$}.
 \end{cases}
\end{equation}
Combined with Lindemann's theorem~\cite{Lindemann} that~$\pi$ is transcendental and
Weierstrass's extension~\cite{Weierstrass} to the transcendence of the natural logarithms of
algebraic numbers other than~$0$ or~$1$,
one knows all of these values to be transcendental.
The remaining~$L$-values are far less well understood.
Indeed, in the (approximately) last 140 years since~\cite{Weierstrass}
only a \emph{single} further explicit number~$L(k,\chi)$ 
has been shown to be irrational,
namely Ap\'{e}ry's unexpected
1978 proof that~$\zeta(3)$ is irrational~\cite{AperyHistoric,Apery,CohenApery}.
In this paper, we establish the irrationality of a new~$L$-value~$L(k,\chi)$;
in some sense the ``simplest'' open case corresponding to~$k = 2$ and the character~$\chi = \chi_{-3}$
of smallest possible conductor:

\begin{ithm} \label{mainA}
The period
$$L(2,\chi_{-3}) = \sum_{n=0}^{\infty} \left( \frac{1}{(3n+1)^2} - \frac{1}{(3n+2)^2} \right) = 0.7813024128964862968\ldots $$
$$ =   \iint_{1 \ge y \ge x \ge 0}  \frac{dx dy}{ y (1 + x + x^2)}
= -  \int_{0}^{1} \frac{ \log(x) dx}{1+ x+x^2}$$
is irrational. More generally, the three periods~$1, \pi^2, L(2,\chi_{-3})$
are linearly independent over~$\Q$.
\end{ithm}

The formula above exhibits~$L(2,\chi_{-3})$  as a period in the sense of Kontsevich--Zagier~\cite{KontsevichZagier}. 
There are a panoply of other more complicated expressions for $L(2,\chi_{-3})$ as an integral, or an infinite sum,
 for example, the following sum  of hypergeometric type~(\cite[\S3]{MR2900448}):
$$L(2,\chi_{-3}) =  \frac{1}{27} \sum_{n=1}^{\infty} \frac{(4 - 15 n) (-27)^n}{n^3 \displaystyle{
\binom{2n}{n}^2 \binom{3n}{n}}
},
$$
or, more serendipitously,
 in terms of the sum of the inverse squares of the entries greater than one in Pascal's triangle~\cite{Pascal}:
$$L(2,\chi_{-3}) = - \frac{1}{3} + \sum_{n >m > 0} \frac{1}{
\displaystyle{
\binom{n}{m}^2
}
}.
$$
The constant
$3\sqrt{3} \, L(2,\chi_{-3})/4  
=  \mathrm{Im}(\Li_2(e^{ \pi i/3}))  = \frac{3}{2} \cdot  \mathrm{Im}(\Li_2(e^{2 \pi i/3}))
 = 1.014941\ldots$ is the volume of the regular ideal hyperbolic
tetrahedron (the one of the maximal volume), and is also the  volume
of the non-compact hyperbolic manifold with minimal volume~\cite{Adams} (the Gieseking manifold, whose
orientable double cover is the complement of the figure~$8$ knot~\cite{ThurstonBook}). 
It is an open problem to show that the volumes of hyperbolic~$3$-manifolds
are not all rationally related
(see~\cite[Problem~23]{Thurston}, and~\cite{Milnor2,Milnor1}). While our result does not have any direct implications for
this question, it is the first unconditional 
 result to make  contact with the arithmetic nature of these volumes.
Another appearance of~$L(2,\chi_{-3})$ is in Smyth's formula~\cite[Prop.~3.4]{BrunaultZudilin}
\begin{equation}
\label{smyth}
\frac{3\sqrt{3}}{4\pi} L(2,\chi_{-3}) = m(1+x+y ) :=  \int_0^1 \int_0^1 \log{| 1  + e^{2\pi i s} + e^{2 \pi i t} |} \, ds \, dt
\end{equation}
linking $L(2,\chi_{-3})$ to the \emph{Mahler measure} of the simplest essentially bivariate polynomial $1+x+y$, or equivalently in the 
language of Diophantine and Arakelov geometry~\cite{Philippon,BostGilletSoule},
to the canonical height of the subvariety $1+x+y = 0$ of the linear algebraic torus~$\mathbb{G}_m^2$. 
Unfortunately, while the nonvanishing Mahler measures of the integer univariate polynomials are all known to be transcendental by the Hermite--Lindemann--Weierstrass theorem~\cite{Hermite,Lindemann,Weierstrass}, our result has no direct bearing on the conjectured irrationality
of any such 
 canonical heights! (We do, incidentally, also prove the irrationality of the Mahler measure of the \emph{rational} coefficients bivariate polynomial $(1+x+y)^4/3$, which is not a canonical height. This is in Theorem~\ref{mixed}.)

An immediate consequence of Theorem~\ref{mainA} is 
the irrationality of the following values of the ``trigamma''
function~$\psi_1(z) =\displaystyle{ \sum_{n=0}^{\infty} \frac{1}{(z+n)^2} = \zeta(2,z)}$:

\begin{icor}
The following numbers are irrational:
$$\begin{aligned}
\frac{1}{1^2} + \frac{1}{4^2} + \frac{1}{7^2} + \frac{1}{10^2} +  \ldots = & \
\quad \frac{L(2,\chi_{-3})}{2} + \frac{2 \pi^2}{27}, \\
\frac{1}{2^2} + \frac{1}{5^2} + \frac{1}{8^2} + \frac{1}{11^2} + \ldots = & \
-\frac{L(2,\chi_{-3})}{2} + \frac{2 \pi^2}{27}, \\
\frac{1}{1^2} + \frac{1}{7^2} + \frac{1}{13^2} + \frac{1}{19^2} + \ldots = & \
\quad \frac{5 L(2,\chi_{-3})}{8} + \frac{\pi^2}{18}, \\
\frac{1}{5^2} + \frac{1}{11^2} + \frac{1}{17^2} + \frac{1}{23^2} + \ldots = & \
- \frac{5 L(2,\chi_{-3})}{8} + \frac{\pi^2}{18}. \end{aligned}
$$
\end{icor}

Note that since~$\psi_1(x+1) - \psi_1(x) + 1/x^2 = 0$,   it 
also follows from Theorem~\ref{mainA}
(together with the fact that~$\psi_1(1) = \pi^2/6$ and~$\psi_1(1/2) = \pi^2/2$)
that~$\psi_1(n/6)$ is irrational for any~$n \in \NwithoutzeroA$.
These are the first new irrationality results for~$\psi_1$ since Legendre's
proof in 1794~\cite{Legendre} that~$\pi^2$ is irrational!

As another application of what turns out to be exactly the same argument, we also prove the following irrationality result for
certain products of two logarithms (see Theorem~\ref{logsmainagain} in~\S~\ref{sec:logs}).

\begin{ithm} \label{logsmain}
Let~$m, n \in \Z \setminus \{-1,0\}$ be  integers such that
$\displaystyle{\left| \frac{m}{n} - 1 \right| < \frac{1}{10^6}}$.
Then
\begin{equation*} 
\log \left(1 + \frac{1}{m} \right) \log \left( 1 + \frac{1}{n} \right)
\end{equation*}
is irrational.
Moreover,
for~$m \ne n$, the following are linearly independent over~$\Q$:
\begin{equation*}
1, \quad \log \left(1 + \frac{1}{m} \right), \quad  \log \left( 1 + \frac{1}{n} \right), \quad 
\log \left(1 + \frac{1}{m} \right) \log \left( 1 + \frac{1}{n} \right). 
\end{equation*}
\end{ithm}

\subsection{Comparisons to the work of Ap\'{e}ry}

Ap\'{e}ry's proof~\cite{Apery} consisted of finding an explicit sequence of rational approximations which
converged  ``sufficiently quickly'' to~$\zeta(3)$ to prove that~$\zeta(3)$ is irrational.
Ever since Ap\'{e}ry's result, considerable effort has been expended in searching
for analogous
sequences which demonstrate the irrationality of other~$L$-values~$L(k,\chi)$ beyond those
of the form~(\ref{dirichlet}). 
Unfortunately, despite enormous efforts, no such sequences have ever been 
found.\footnote{
One significant step which is not directly related to the irrationality of specific~$L$-values is
the theorem of Ball and Rivoal~\cite{Rivoal1,BallRivoal} that infinitely many odd zeta values are irrational.
One refinement by Zudilin~\cite{Zudilin57911} 
proves that at least one of the values~$\zeta(5)$,
$\zeta(7)$, $\zeta(9)$ and~$\zeta(11)$ is irrational.}
In particular, in this paper, we do \emph{not} find (directly)
any new convergent sequences to~$L(2,\chi_{-3})$.
Instead, we show how one can exploit
the arithmetic nature of \emph{known} approximations 
(found by Ap\'{e}ry and others) in a more subtle way
using both methods from transcendental number theory
and complex analysis. 

In order to introduce our main idea, we begin with an exposition
and then a reformulation of some of the key features of Ap\'{e}ry's proof. 
The first remark to make is that Ap\'{e}ry's proof uses very little
number theory;
indeed the only number theoretic input
is a (weak form) of the prime number theorem and
 the following elementary lemma:

\begin{lemma} \label{simple} If there is a~$\delta > 0$ and a sequence
of rational numbers~$p_n/q_n \ne \beta$ with~$q_n \rightarrow \infty$
such that
$$\left| \beta - \frac{p_n}{q_n} \right| < \frac{1}{q^{1+\delta}_n} \quad
n = 1,2, \ldots,$$
then~$\beta$ is irrational.
\end{lemma}

This lemma is true even with the weaker hypothesis that~$\left| \beta - p_n/q_n  \right| = o(1/q_n)$.
Ap\'{e}ry writes down a pair of power series~$A(x), B(x) \in \Q\llbracket x \rrbracket$ and the linear 
combination
$$P(x) =  B(x) - \zeta(3) A(x) = \sum_{n=0}^{\infty} x^n \left(b_n - \zeta(3) a_n \right).$$
 The coefficients~$a_n$ and~$b_n$ are rational numbers, and more precisely: 
 $$
 a_n \in \Z, \qquad [1,2,3,\ldots,n]^3 b_n \in \Z.
 $$
  Here and throughout 
 our paper, we follow the conventional notation~$[1,2,\ldots,n]$ for the lowest common multiple of the first~$n$ integers. The prime number theorem 
 determines the growth rate of these denominators:
$$
\log [1,2,3,\ldots,n] = n + o(n). 
$$
At the same time, Ap\'{e}ry proves that~$A(x)$ (and~$B(x)$) have radius
of convergence~$(\sqrt{2}-1)^4$ whereas~$P(x)$ has radius of convergence exactly~$(\sqrt{2}+1)^4$.
Now one exploits the inequality
\begin{equation}
\label{lucky}
4 \log (\sqrt{2}+1) > 3
\end{equation} to deduce, by Lemma~\ref{simple} with~$p_n/q_n = b_n/a_n$
and
$$\delta = \frac{4 \log (\sqrt{2}+1) - 3}{4 \log (\sqrt{2}+1) + 3},$$
 that~$\zeta(3) \notin \Q$.
 
There are a number of other situations where one can construct functions~$A(x)$, $B(x)$
 of a similar flavour so that a particular linear combination~$P(x) = B(x) - \eta A(x)$
 has extra convergence properties, and
 where~$\eta = L(k,\chi)$ turns out to be the unique complex number characterized in this way.
 But the analogous inequality~(\ref{lucky}) always seems to fail,\footnote{except in 
one notable example found by Ap\'{e}ry himself with~$\eta = \zeta(2)$.}
and one can draw no consequences about the arithmetic of the corresponding~$L$-value. 
(For one
particularly interesting study of sequences of the form considered
by Ap\'{e}ry, see~\cite{Zagier}. In our proof of Theorem~\ref{mainA}, we will make a central use of some of the sequences (re-)discovered by Zagier
in his search.)

The starting point of our investigation is that, even when the analogue of~(\ref{lucky}) fails as it usually does, the functions~$P(x)$ arising in these constructions have
more structure which has not previously been exploited. Ap\'{e}ry's functions~$A(x)$
and~$B(x)$ turn out to
satisfy a linear ordinary differential equation (ODE) with coefficients in~$\Z[x]$
which only has (regular) singular points at~$x=0$, $\infty$,
and~$(\sqrt{2} \pm 1)^4$. The function~$P(x)$ arises as the unique (up to scalar)
linear combination of the two dimensional space of solutions to this ODE
which are holomorphic at~$0$ with the additional
property that it is also holomorphic at~$(\sqrt{2}-1)^4$.
This implies that~$P(x)$, for example, is not merely holomorphic on the disc of radius~$(\sqrt{2}+1)^4$, but extends to a holomorphic function on all 
of~$\C \setminus [(\sqrt{2}+1)^4,\infty)$, or (more relevantly for our ultimate
purposes, but less important for the introduction) to a
function on the universal cover of
$$\mathbf{P}^1 \setminus \{0,(\sqrt{2}-1)^4,(\sqrt{2}+1)^4,\infty\}$$
which is holomorphic at~$x=0$
and overconverges beyond the first singularity~$(\sqrt{2}-1)^4$.
All Ap\'{e}ry uses is that~$P(x)$ is holomorphic on the disc of radius~$(\sqrt{2}+1)^4$.

Now imagine an analogous situation where~$A(x)$ and~$B(x)$
are holomorphic (at~$x=0$) solutions to an ODE\footnote{By this we mean: a \emph{linear} ODE over~$\Q(x)$, as always in this paper.} with regular singular
points at~$0$, $\infty$ and a pair of real numbers~$0<\alpha<\beta$,
and~$P(x) = A(x) - \eta B(x)$ is a linear combination which is also
holomorphic at~$\alpha$, but whose analytic continuation has a singularity
at~$\beta$ and does not analytically continue to a meromorphic
function at~$x = \beta$.
 But now suppose   ---   taking into account the denominators of~$a_n$ and~$b_n$
   ---   that the constant~$\beta$
is not large enough to imply that the corresponding convergents~$p_n/q_n$
satisfy Lemma~\ref{simple}. Is there a way to exploit
the fact that not only is~$P(x)$ a  holomorphic function on the disc
of radius~$\beta$, but also extends to a holomorphic
function on~$\C \setminus [\beta,\infty)$ and on the universal
cover of~$\mathbf{P}^1 \setminus \{0,\alpha,\beta,\infty\}$?

To make things simpler (too simple, in fact --- we will return to the issue of the necessity of denominators), let us momentarily suppose that the~$a_n$
and~$b_n$ are actually integers.  To run Ap\'{e}ry's proof scheme via Lemma~\ref{simple},
it then would have sufficed that~$\beta > 1$.
Lemma~\ref{simple} in this case has the following alternate formulation: 
\begin{lemma}[Simple Lemma] \label{simpletwo}
 An integer power series in~$\Z \llbracket x \rrbracket$ that defines a holomorphic function on the disc $|x| < R$
of a radius~$R>1$ is a polynomial. 
\end{lemma}
In our running example, if~$\beta < 1$, we can deduce
nothing from Lemma~\ref{simpletwo}, but we can still derive that~$P(x)$
is holomorphic on~$\C \setminus [\beta,\infty)$. 
There is an entire subject devoted to more subtle extensions of Lemma~\ref{simpletwo},
beginning with the Theorem of Borel--P\'{o}lya~\cite[Chapter 5]{Amice}, which allows one to make conclusions about~$P(x) \in \Z \llbracket x \rrbracket$ from weaker analytic hypotheses
than simply converging on a disc of sufficiently large radius. We recall (a special case of) this theorem now.
If~$\Omega \subset \C$ is a simply connected open region containing~$0$, then, from the Riemann mapping theorem, there exists  a biholomorphic map~$\varphi: \D \rightarrow \Omega$
with~$\varphi(0) = 0$. The map~$\varphi$ is unique up to biholomorphisms of the unit disc fixing~$0$, which are all given by rotations. In particular, 
the invariant~$|\varphi'(0)|$ does not depend on the choice of~$\varphi$,  and (by definition) is equal to the conformal radius~$\rho(\Omega,0)$ 
of~$\Omega$ at~$0$. 
The conformal radius of the disc~$\D_R  = D(0,R)$ is equal to~$R$ (via the map~$\varphi(z) = Rz$), but the conformal radius of any other~$\Omega$ is strictly
larger than the radius of the largest disc contained in~$\Omega$ and centered at~$z=0$. 
We have:

\begin{theorem}[Borel--P\'{o}lya,~\cite{PolyaOriginal}] \label{borelpolya}
A power series~$P(x) \in \Z \llbracket x \rrbracket$ that continues analytically to a simply connected open region~$0 \in \Omega \subset \C$ of conformal radius~$\rho(\Omega,0) = |\varphi'(0)| > 1$ is necessarily a rational function:~$P(x) \in \Q(x)$. 
\end{theorem}
For example, the biholomorphic map
$$\varphi: \D \rightarrow \C \setminus [\beta,\infty), \quad
z \mapsto \frac{4 \beta z}{(1 + z)^2}$$
shows that~$\C \setminus [\beta,\infty)$
has conformal radius~$|\varphi'(0)| = 4 \beta$. It follows from Theorem~\ref{borelpolya}
 in our imagined example
above that~$P(x)$ is a rational
function as soon as~$4\beta > 1$, 
contradicting the assumption that~$P(x)$ was not meromorphic at~$\beta$,
and implying that~$\eta$ is irrational.
This  is already clearly an improvement
on the condition that~$\beta > 1$. (The basic idea for this special case of P\'olya's theorem is sketched in Remark~\ref{remark on the slit plane} at the end of this introduction.)

Even beyond Theorem~\ref{borelpolya} (as we shall discuss in Section~\ref{seriousintro} below), there are algebraicity theorems of Andr\'{e} and others with even weaker hypotheses that allow one to deduce that~$P(x)$ is algebraic over~$\Q(x)$ (see~\cite{Andre,AndreG} and~\cite{UDC}), which can often be ruled out directly
in practice for any particular~$P(x)$.

The main thrust of our paper is in adapting and honing up the methods of Borel, P\'olya, and Andr\'e to fit into the Ap\'{e}ry irrationality
proofs context. The algebraicity criteria as such do not apply, because
the power series~$A(x)$ and~$B(x)$ of relevance to Ap\'{e}ry style proofs
never (both) have integral coefficients. 
And indeed, when one introduces denominators (even of some controlled flavour), it turns out that algebraicity is no longer the right property to consider.
To begin with, a theorem of Eisenstein~\cite[\S11.4]{BombieriGubler} states that the power series expansion of any algebraic function in~$\overline{\Q(x)} \cap \Q\llbracket x \rrbracket$ has~$\Z[1/S]$ coefficients for some~$S \in \NwithoutzeroA$.
But from the point of view of the various proofs of Borel's theorem (and its variations),
 if~$P(x) \in \Z \llbracket x \rrbracket$,
then so too are all of its powers; but if~$P(x)$ has denominators, then the powers of~$P(x)$ typically have \emph{worse} denominators. An example to keep in mind is
$$-\log(1-x) = \sum_{n=1}^{\infty} \frac{x^n}{n}.$$
This function has the property that multiplication by~$[1,2,3,\ldots,n]$ (of order~$e^n$) simultaneously clears the denominators of the first~$n$ coefficients. But in order to clear
the denominators of the first~$n$ coefficients of~$\log^m(1-x)$, one has to multiply by a denominator of order
$$[1,2,\ldots,n] \times [1,\ldots, n/2 ] \times
\cdots \times [1,\ldots, n/m ] = \exp \left(n \left(1 + \frac{1}{2} + \ldots + \frac{1}{m} \right) + o(n) \right).$$
Here and throughout the paper by~$[1,2,3, \ldots, bn]$ for~$b \in \R^{>0}$ we mean
by abuse of notation~$[1,2,3, \ldots, \lfloor bn \rfloor]$.
On the other hand, \emph{differentiation} does  preserve the property of controlled denominator growth. Hence, instead of an algebraicity theorem, one should expect an \emph{arithmetic holonomy} bound,
where one bounds the dimension of a~$\Q(x)$-vector space generated by functions with certain denominator growth and analytic properties,
and which is closed under differentiation. This in particular
implies that, in the appropriate generalization of the Borel--P\'olya conditions, the \emph{solutions}  ---  a precise formulation is given by Corollary~\ref{holonomic criterion}, to be 
discussed in detail in the next section  ---  are~$G$-functions in the sense of Siegel
(see~\cite{Siegel1929SNS} and~\cite[\S VIII.1]{Dwork}, see also~Definition~\ref{gfunctiondef}).
Moreover, one can hope to give a --- good enough --- explicit bound on the order of the~$G$-function
(that is, the rank of the~$\Q(x)$-module generated by~$G$ and its derivatives),
which  ---  in a given situation  ---  contradicts the structure of some explicit approximation
function~$P(x) = B(x) - \eta A(x)$. When this is achieved, the ultimate contradiction is in the supposition that 
$P(x) \in \Q\llbracket x \rrbracket$, that is that $\eta \in \Q$.

These arithmetic holonomy bounds are ultimately the main concern of this paper,
and we take up a detailed introduction to them in our next section~\S~\ref{seriousintro}. 
 
\begin{remark}  \label{remark on the slit plane}
 In a very special case, a hint in this direction has 
 been previously proposed
 (although without any application to a new irrationality proof) by Zudilin~\cite{ZudilinDet}, who isolated a condition on the linear forms $c_n = b_n - \eta a_n$ which implies an analytic continuation of the generating function
 $P(x) = \sum_{n=0}^{\infty} c_n x^n$ to a slit plane $\C \setminus [\beta,\infty)$. Whereas Ap\'ery's use of the convergence radius focused on the decay rate $\beta^{-n+o(n)}$
 of the coefficients $c_n$, Zudilin highlights the improved 
decay rate $(4\beta)^{-n^2+o(n^2)}$ of the sequence of \emph{Hankel determinants} $\det(c_{i+j})_{i,j=0}^{n}$; this is indeed a well-known consequence~\cite{PolyaDet,PommerenkeDet} 
of the analyticity of $P(x)$ on  $\C \setminus [\beta,\infty)$. The latter is in fact closely linked to the proof of Theorem~\ref{borelpolya} in this particular case of $\Omega = \C \setminus [\beta,\infty)$ with 
$\beta > 1/4$: if all $c_n \in \Z$, the Hankel determinants are also rational integers, therefore they vanish from some point onward if they decay at an exponential rate smaller than one, 
and finally this means $P(x) \in \Q(x)$ by the rationality criterion of Kronecker. Quantifying the denominator of the Hankel determinant in the case $c_n \in \Q$  leads as well to Zudilin's determinantal criterion for $\eta \notin \Q$. See~\S~\ref{univalent det} for a precise formulation and a generalization.~\endofremark
\end{remark}

\begin{remark}[A remark on exposition]
We take the point of view that the readership of this paper might include mathematicians not familiar
with either the details of our previous paper~\cite{UDC} or the methods of Diophantine analysis more broadly.
At the risk of interrupting the flow of the exposition, we have included a number of
expositional asides denoted by ``basic remarks'' 
 throughout
the paper which are intended to help orient the reader less familiar with this material; the expert should
feel free to skip over these.
\end{remark}

\begin{remark}[A remark on notation]
We shall use~$\NwithzeroB = \{0,1,2, \ldots\}$ to denote the natural numbers with zero,
and~$\NwithoutzeroA$  to denote the positive integers. Depending on the context, $\P^1$ will signify either a scheme isomorphic
to~$\mathrm{Proj}{ \, \Z[T_0,T_1]}$ (the projective line over~$\Z$), or 
the complex manifold~$\P^1(\C)$ of its $\C$-valued points (the Riemann sphere with 
coordinate~$z := T_0/T_1$, 
elsewhere commonly denoted~$\widehat{\C}$ or~$\C\P^1$). A similar ambiguity is adopted for the modular stacks~$Y_0(2)$ and~$Y(2)$.  The complex disc~$D(0,R)$ of radius~$R \in (0,\infty]$ in the relevant 
coordinate (always clear by the context, but most frequently denoted~$z$) will be denoted by~$\D_R$, and we shall write $\D := \D_1$ for the unit radius disc and~$\Db$ for its closure in~$\C$. The unit circle~$\partial \Db = \{e^{2\pi i \theta} : \theta \in [0,1]\}$ is denoted~$\T$ and its uniform measure~$d\theta$ is denoted~$\mv$. For a connected complex manifold~$M$, we shall denote by $\mathcal{O}(M)$ and $\mathcal{M}(M)$, respectively, the ring of holomorphic functions the field of meromorphic 
functions on~$M$. The notation~$\mathcal{O}(\Db)$ and~$\mathcal{M}(\Db)$ is used for the corresponding functions on some unspecified open neighborhood of the closed unit disc~$\Db \subset \C$. 
Throughout our paper, we will usually write $q := e^{\pi i \tau}$ for $\tau$ belonging to the upper half plane $\H$, although
we will occasionally write~$q = e^{2 \pi i \tau}$.  (As noted in~\cite{UDC}, this is forced upon us
by historical convention, but we always use the first choice unless explicitly stated otherwise.) By a mild and harmless notational abuse, 
the \emph{modular lambda function} 
\begin{equation}
\label{lambdadef}
\lambda(q) := \frac{ \left(
\displaystyle{ \sum_{n \in 1+2\Z} q^{n^2/4} }\right)^4  }{\left(
\displaystyle{ \sum_{n \in 2\Z} q^{n^2/4}} \right)^4} = 16 q \prod_{n=1}^{\infty} \left( \frac{1+q^{2n}}{1+q^{2n-1}} \right)^8: 
\, \H \to \C \setminus \{0, 1\}, \  \{ q = 0\} \mapsto 0
\end{equation}
will be written in the cusp-filling coordinate $q \in \D  := \{ |q| < 1\}$ rather than $\tau = \log(q) / (\pi i)$. 
The letter~$e$ is generally reserved for the Euler constant $e \approx 2.718281$.
Finally, we admit a minor notational abuse by adopting the convention of writing~$X \setminus \{A,B\} := X \setminus (A \cup B)$ for any subsets~$A$ and~$B$ of a set~$X$. 
\end{remark}

\subsection{The paths to Theorems~\ref{mainA} and~\ref{logsmain}, and an outline of the paper}
\label{sec:leitfaden}
The following leitfaden (Figure~\ref{leit}) gives in  summary the logical structure of our paper.
Here the pair of dotted lines indicates that
there are two alternate paths to Theorems~\ref{mainA} and~\ref{logsmain},
either through~\S~\ref{fine section} (by multivariable methods, based on measure concentration) or~\S~\ref{new slopes} (by single variable methods, 
based on some Arakelov theory and Bost's inequality on evaluation heights).
We also omit~\S~\ref{app:PerelliZannier}, which is most closely
related (though there is no dependency in either direction)
to~\S~\ref{new slopes}.
There could be some (modest) economy if we restricted ourselves
to the shortest possible proof of Theorem~\ref{mainA}. However,
with a view to both  future developments and applications, we 
felt it was better
to include all these new ideas. In many ways this reflects our
experience with our
 previous
paper~\cite{UDC} which included three proofs
of the main holonomicity theorem~\cite[Theorem~2.0.1]{UDC}.
One of the referees  of that paper
recommended removing one particular proof of this theorem
whose
ideas subsequently proved essential  for the advances in this paper.

{\small
\begin{center}
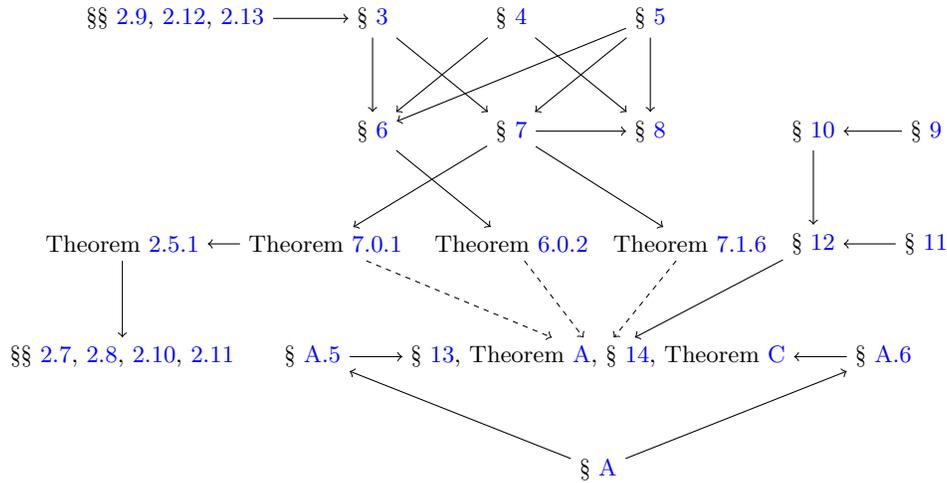
\begin{figure} 
\begin{tikzpicture}[node distance = 1.5cm, auto]
\node (A){\S\S~\ref{local univalent leaves}, \ref{sec:ideas outline}, \ref{sec_refined}};
\node (B)[right=1cm of A]{\S~\ref{functional transcendence}};
\node (C)[right=1.2cm of B]{\S~\ref{concentration of measure}	}; 
\node (CC)[right=1.2cm of C]{\S~\ref{integration cost} };
\node (H)[below of=B]{ \S~\ref{fine section} };
\node (I)[below of=C]{ \S~\ref{new slopes}	};
\node (J)[below of=CC]{	\S~\ref{slopes}};
\node (X)[below of=I]{Theorem~\ref{main:elementary form}};
\node (Y)[left=0.2cm of X]{Theorem~\ref{main:BC form} };
\node (Z)[right=0.1cm of X]{Theorem~\ref{main:BC conv discrete}};
\node (K)[right=0.1cm of Z]{\S~\ref{sec:lindep} };
\node (E)[above of=K]{\S~\ref{sec:purefunctions}};
\node (F)[right of=E]{	\S~\ref{sec:YtoY0(2)}	 };
\node (G)[right of=K]{ \S~\ref{Zagier local system}		 };
\node (MM)[below of=Z]{};
\node (M)[left=-1.5cm of MM]{ \S~\ref{sec:proofA}, Theorem~\ref{mainA}, \S~\ref{sec:logs}, Theorem~\ref{logsmain}};
\node (N)[left=0.7cm of M]{\S~\ref{numerical integration}};
\node (O)[right=0.7cm of M]{ \S~\ref{contour choiceC}};
\node (P)[below of=M]{\S~\ref{contour choiceA}};
\node (Q)[left=2.9cm of X]{Theorem~\ref{basic main}};
\node (R)[below of= Q]{ \S\S~\ref{sec:logarithmchar}, \ref{sec:beyondlog},
\ref{sec:first irrationality proofs}, \ref{bivalent app}};
\draw[->] (C) -- (H);
\draw[->] (C) -- (J);
\draw[->] (A) -- (B);
\draw[->] (B) -- (I);
\draw[<-] (E) -- (F);
\draw[->] (G) -- (K);
\draw[->] (B) -- (H);
\draw[->] (CC) -- (I);
\draw[->] (CC) -- (J);
\draw[->] (I) -- (J);
\draw[->] (E) -- (K);
\draw[->] (K) --(M);
\draw[->] (Y)--(Q);
\draw[->] (CC)--(H);
\draw[->] (I)--(Z);
\draw[->] (I)--(Y);
\draw[->] (H)--(X);
\draw[dash pattern=on 2pt off 2pt, ->] (X) -- (M);
\draw[dash pattern=on 2pt off 2pt, ->] (Z) -- (M);
\draw[dash pattern=on 2pt off 2pt, ->] (Y) -- (M);
\draw[->] (N)-- (M);
\draw[->] (O)--(M);
\draw[->] (P) -- (N);
\draw[->] (P) -- (O);
\draw[->] (Q) -- (R);
\end{tikzpicture}
\caption{Leitfaden: paths to Theorems~\ref{mainA} and~\ref{logsmain}}
\label{leit}
\end{figure}
\end{center}
}
We now very briefly outline the paper.
Section~\S~\ref{seriousintro} is mainly introductory,
although~\S\S~\ref{sec:holbasic}--\ref{bivalent app} present a basic form
of our main results (which will not be proved until~\S~\ref{fine section}) together
with some applications, and~\S\S~\ref{sec:ideas outline}--\ref{sec_refined}
outlines our approach
to proving holonomy bounds, which is followed up in precise detail
in~\S~\ref{functional transcendence}, with elements of functional transcendence theory. In~\S~\ref{functional transcendence}, we
also include some further exposition of related material in its proper historical context,  intended to help place our ideas into a broader context.
In~\S~\ref{concentration of measure}, we collect some basic  facts
concerning large deviations and the concentration of measure phenomenon in high dimensions.
In~\S~\ref{integration cost}, we introduce the idea (possibly
counterintuitive in light of the discussions in~\S~\ref{seriousintro})
that it can sometimes be useful to  integrate our putative functions despite introducing new denominators.
Also included are some technical computations related
to extra denominators arising from integrations, which follow from
 the prime number theorem.
In~\S~\ref{fine section}, we prove
our first main holonomy bound Theorem~\ref{main:elementary form}.
In~\S~\ref{new slopes}, based on the work of Bost and Charles~\cite{BostCharles}, we prove our second main holonomy bound (or more precisely, several closely related bounds) using Bost's slopes method framework. In particular, Theorem~\ref{main:BC form} is essentially the bound of Bost--Charles in~\cite{BostCharles} incorporated with our treatment of denominators in~\S~\ref{fine section}; Theorem~\ref{main:BC conv discrete} is a further improvement of Theorem~\ref{main:BC form} using the convexity property of a growth characteristic function which is closely related to the Bost--Charles bound and behaves similarly to a Nevanlinna characteristic function. In~\S~\ref{slopes}, we unify our methods from~\S\S~\ref{fine section}--\ref{new slopes} and obtain, with an eye to future applications, the sharpest holonomy bound in our paper. 
In sections~\S\S~\ref{sec:YtoY0(2)}--\ref{sec:lindep} we return to a discussion
of specific \emph{templates} (situations
in which the denominator types and singularities are fixed)
in order to prepare for the application
of our holonomy bounds to our main irrationality results.
In~\S~\ref{sec:YtoY0(2)}, we use the map of modular curves~$Y(2) \rightarrow Y_0(2)$
to relate two templates over~$\PP^1 \setminus \{0,1,\infty\}$ and~$\PP^1 \setminus \{0,4,\infty\}$,
respectively.
In~\S~\ref{sec:purefunctions}, we discuss some~$G$-functions
on~$\PP^1 \setminus \{0,1,\infty\}$ with simple denominator types (most of them
well-known, but also one which was surprising to us),
and in~\S~\ref{Zagier local system}
we introduce certain local systems arising in~\cite{Zagier} which, contingent
on a hypothetical linear dependence of~$1$, $\pi^2$, and~$L(2,\chi_{-3})$, give rise
to more~$G$-functions.
 \S~\ref{sec:lindep} is concerned with proving
the linear independence of all these functions over~$\C(x)$.
In~\S\S~\ref{sec:proofA} and~\ref{sec:logs} we
give the proofs of Theorems~\ref{mainA} and~\ref{logsmain} respectively,
using some explicit computations which are explained in detail
in~\S~\ref{contour choiceA}. (For a proof of Theorem~\ref{logsmain} only, a number of
 subsections, including all of~\S~\ref{Zagier local system}, can also be omitted.)
Finally, the short \S~\ref{app:PerelliZannier} is intended as a showcase of the basic proof scheme, and can serve most particularly as an introduction to~\S~\ref{new slopes}.

\section{The main arithmetic holonomy bound}
\label{seriousintro}

We begin with a discussion of the dimension bounds in their simplest case. 

\subsection{The algebraic case} \label{sec:algebraic}
Our solution~\cite{UDC} of the ``unbounded denominators'' conjecture was based on the following  dimension upper estimate on a certain $\Q(x)$-linear space
of algebraic functions. We called this type of result an \emph{arithmetic holonomy bound}, and while our reason for this name remained obscure in~\cite{UDC}, 
we hope it should be vindicated by our present paper where we treat more general holonomic functions whose analytic continuations generate an infinite and nonsolvable monodromy group. 
Given a holomorphic mapping $\varphi : \Db \to \C$ on some neighborhood of 
the closed unit disc~$\Db \subset \C$ and taking $\varphi(0)=  0$ with $|\varphi'(0)| > 1$, we established~\cite[Theorem~2.1]{UDC} the dimension upper bound
\begin{equation} \label{original UDC bound}
\dim_{\Q(x)} \mathcal{H}(\varphi) \le e \cdot \frac{\displaystyle{\int_{\T} \log^+{|\varphi|} \, \mv}}{\log{|\varphi'(0)|}}
\end{equation}
on the $\Q(x)$-linear span~$\mathcal{H}(\varphi)$  of the $\Z\llbracket x \rrbracket$ formal power series  $f(x)$ whose pullback $f(\varphi(z))$ also converges on 
a neighborhood of~$\Db$. 
Here~$\T$ is the unit circle, $\log^{+}|x|$ is defined to be~$\max(0,\log |x|)$, and~$\mv$ is just the usual Haar measure on~$\T$, so that the integral
in the numerator can equally be written as~$\int_{0}^{1} \max \left(0,\log |\varphi(e^{2 \pi i t})|\right) dt$. 

\begin{basicremark}  \label{Polya and Andre}
Suppose that~$f(x)$ is a power series
which extends to a holomorphic function on a domain~$\Omega \subset \C$ containing the origin,
of conformal mapping radius~$\rho(\Omega,0)  > 1$. (See the beginning of \S~\ref{sec:logarithmchar} for a precise definition of conformal mapping radius.) By definition, there  consequently exists a biholomorphic
map~$\varphi: \D \rightarrow \Omega$ with~$\varphi(0) = 0$ and~$|\varphi'(0)| = \rho(\Omega, 0) > 1$.
In turn, the holomorphy of~$f(x)$ on~$\Omega$ means exactly that the pulled back power series~$f(\varphi(z)) \in \C \llbracket z \rrbracket$  converges on~$\D$.
The bound~\eqref{original UDC bound} then implies that the~$\Q(x)$-vector space
 generated by~$f(x) \in \Z\llbracket x \rrbracket$ \emph{and its powers}
is finite dimensional (since the powers of~$f(x)$ also lie in~$\HH(\varphi)$), and thus~$f(x)$ is algebraic (of some explicitly bounded
degree). However, under these assumptions, one can already deduce the
\emph{rationality} of~$f(x)$ from the Borel--P\'olya Theorem~\ref{borelpolya}. 
So the bound~(\ref{original UDC bound})  is more interesting when~$\varphi$ is
\emph{not} univalent. (We will eventually find that the Borel--P\'{o}lya theorem 
too will be completely subsumed into holonomy bounds finer than~\eqref{original UDC bound}, 
such as the bound~\eqref{BC integral} below, which is due to Bost and Charles~\cite{BostCharles}, and ultimately
our main new holonomy bound~\eqref{new bound} in this paper.)

A non-univalent example  is as follows.
Suppose that~$f(x) \in \Z\llbracket x \rrbracket$ can be analytically continued on any path from~$0$ in~$\C$
avoiding both~$0$ and some fixed real number~$\alpha > 0$.
For example, take~$f(x) = (1-4x)^{-1/2} = \sum \binom{2n}{n} x^n$, and~$\alpha = 1/4$.
 Then one can take~$\varphi$
to be any holomorphic function with~$\varphi(0) = 0$ but which has no
other preimages of either~$0$ or~$\alpha$. One such function
is
$$\varphi(z) = \alpha \lambda(z)$$
 where~$\lambda$
is the modular~$\lambda$ function as given in~(\ref{lambdadef}).
In this case, we have~$|\varphi'(0)| = 16 \alpha$.
Hence, if~$\alpha >  1/16$, we deduce that~$f(x)$ is algebraic (with some degree explicitly
bounded by~(\ref{original UDC bound}) over~$\Q(x)$). This example is already due to Andr\'{e}.
The paper~\cite{UDC} is concerned with the case when~$\alpha = 1/16$, where there are infinitely many~$\Q(x)$-linearly
independent algebraic examples including
\[f(x) = \sum_{n=0}^{\infty} \binom{4n}{2n} x^n = \sqrt{ 
\frac{1 + \sqrt{1-16x}}{2-32x}  };\]
 but also  the algebraicity fails without any additional hypothesis,
as can be seen from the hypergeometric example~$f(x) = \sum_{n=0}^{\infty} \binom{2n}{n}^2 x^n$, a case used in~\cite{AndreGtranscendence}
and further discussed in~\cite[Appendix~A]{Andre}. 
We refer any further discussion of the algebraic case to~\cite{UDC}.   
\endofremark
\end{basicremark}

\subsection{Denominators} \label{sec:introdenominators}
 For linear independence proofs,  as suggested
by the examples in~\S~\ref{intro}, we need holonomy bounds on functions in $\Q\llbracket x \rrbracket$ rather than $\Z\llbracket x \rrbracket$. Indeed, the holonomic
coefficients of interest  ---  such as $\eta = L(2,\chi_{-3})$ as our primary focus here  ---  are conjecturally transcendental, and so any realization as numbers in a period matrix must necessarily involve 
a local system with an infinite global monodromy group. 
On the other hand, if~$P(x) \in \Z\llbracket x \rrbracket$ lies in a holonomic module~$\HH(\varphi)$ attached to some $\varphi : \Db \to \C$ with $\varphi(0) = 0$ and $|\varphi'(0)| > 1$, then (as noted previously)
it would follow that~$P(x)$ is algebraic.
The
 \emph{Grothendieck--Katz $p$-curvature conjecture}~\cite{katzpcurvature,Andre} (proved by Katz in many of the cases that are of geometric origin) 
 informally equates  the infinitude of the global monodromy 
 group of an integrable connection with the nonvanishing of the
 \emph{$p$-curvature operator}   ---   the local obstruction to integrability modulo~$p$   ---   for a positive density  of the primes~$p$. 
 But we remind the reader that, even for an irreducible linear homogeneous ODE, a single $\Z[1/S]\llbracket x \rrbracket$ solution does not imply vanishing of the~$p$-curvatures;
rather, a \emph{basis} of $\Z[1/S]\llbracket x^{1/h} \rrbracket$ solutions does. 
 As an example, the function~$A(x) \in \Z\llbracket x \rrbracket$ in Ap\'{e}ry's
 argument (discussed in~~\S~\ref{intro}) has integral coefficients but is not algebraic. To square this example with the remarks about~$P(x)$ above, remember that
 the holonomy bounds are never being applied to~$A(x) \in \Z\llbracket x \rrbracket$ itself, but rather to a (supposed for the contradiction!) $\Q$-linear combination $P(x) = B(x) - \eta A(x)$ of~$A(x)$ and some other solution~$B(x) \in \Q\llbracket x \rrbracket$ of the same ODE. This second solution~$B(x)$ does indeed have denominators involving
 infinitely many primes. 
 
From~\cite{FischlerRivoal}, we have a conjectural\footnote{This is an unconditional theorem for the case of $G$-functions that ``arise from geometry,'' based on
the existence of an $F$-crystal structure at all  but the finitely many primes of bad reduction, cf.~\cite[\S~V app.]{AndreG}.
One can also be more precise: if~$\LL(f) = 0$ for some nonzero~$r^{\mathrm{th}}$-order Fuchsian operator having for~$x=0$ local exponents rational numbers
with denominators dividing~$b$, then the denominators form of~$f$ may be taken as $ A^{n+1} [1, \ldots, b n + b_0]^{r-1}$ for
some~$A \in \NwithoutzeroA$ and~$b_0 \in \Z$. } understanding of the denominator types of Taylor series $P(x) = \sum a_n x^n \in \Q \llbracket x \rrbracket$ arising from~$G$-functions:
 there should exist~$A \in \NwithoutzeroA$, $b \in \Q_{>0}$, and~$\sigma \in \NwithzeroA$ such that
\begin{equation}
\label{integrality}
 a_n A^{n+1} [1, \ldots, b n]^{\sigma} \in \Z \qquad \forall n \in \NwithzeroA;
\end{equation}
here and throughout our paper (as noted previously in the introduction), $[1,\ldots,n]$ is used to denote the lowest common multiple of the first~$\lfloor n \rfloor$ positive integers.

The most basic example is the~$G$-function~$\log(1-x)$.  It has the type~\eqref{integrality} with~$A=1$, $b=1$, and~$\sigma = 1$, but that form can in this case  clearly
 be improved: only an~$n$ is needed out of the $[1,\ldots,n]$ clearance, in reflection of the fact that 
 \[\log(1-x) = \int \frac{dx}{x-1}\]
  is an integral of an algebraic function.   
  It  turns out, cf.~\S~\ref{integration cost}, that it will ultimately be important to exploit such refinements from integrals.
 In any case, the necessity of at least the~$[1,\ldots,bn]$ denominators 
forces us to venture outside of the proper\footnote{A natural framework would be the construction and comparison of integrable connections over formal-analytic arithmetic varieties and their algebraizations. We do not attempt to get into such a concept in the present paper, apart 
from raising one specific finiteness problem
in~\S~\ref{integral bounds}, but we  do hope to turn to it on another occasion.} scope 
 of the theory of formal-analytic arithmetic surfaces~\cite{BostBook,BostCharles}.

 With the presence of denominators, for given holomorphic $\varphi: \Db \rightarrow \C$ and parameters $b\in \Q_{\geq 0}$ and $\sigma 
 \in \NwithzeroA$, we define the \emph{holonomic module} $\mathcal{H}(\varphi; b; \sigma)$ to be the $\Q(x)$-linear span
of all the formal functions of the form
\begin{equation}
f(x) = \sum_{n=0}^{\infty} a_n \, \frac{x^n}{[1,\ldots,bn]^{\sigma}} \in \Q\llbracket x \rrbracket, \qquad a_n \in \Z \quad \forall n \in \NwithzeroA
\end{equation}
whose $\varphi$-pullbacks $f(\varphi(z))$ converge on $\D$.
 The proofs in~\cite[\S~2]{UDC} extend routinely to establish a first result in this direction: 
\begin{equation} \label{founding hol}
\dim_{\Q(x)} \mathcal{H}(\varphi; b; \sigma) \leq e \cdot \frac{\displaystyle{\int_{\T} \log^+{|\varphi|} \, \mv}}{\log{|\varphi'(0)| - \tau}},
\end{equation}
where $\tau := b \sigma$ and we now assume that $\varphi$ has the \emph{conformal size} $|\varphi'(0)| > e^{\tau}$.

Unfortunately, the holonomy bound~\eqref{founding hol}, which worked nicely in the asymptotic framework of~\cite{UDC} where the absolute numerical coefficient was immaterial,  is now far too crude to prove the irrationality of $L(2,\chi_{-3})$. It is then of interest to know the least possible
value that may take the place of the constant~$e$ in the bound~\eqref{founding hol}. Progress was made by Bost and Charles~\cite[Corollary~8.3.5]{BostCharles} who, in the original $\sigma = 0$ case of~\cite{UDC}, established the finer bound by
\begin{equation} \label{BC integral}
\frac{\displaystyle{ \iint_{\T^2}  \log{|\varphi(z) - \varphi(w)  |} \, \mv(z) \mv(w)}  }{\log{|\varphi'(0)|}}.
\end{equation}
In 2023, in response to our question about a similar dimension bound for the general holonomic modules $\mathcal{H}(\varphi; b; \sigma)$, Charles explained to us how their proof can be directly generalized to obtain
\begin{equation} \label{BC integral 2}
\dim_{\Q(x)} \mathcal{H}(\varphi; b; \sigma) \leq \frac{ \displaystyle{ \iint_{\T^2}  \log{|\varphi(z) - \varphi(w)  |} \, \mv(z) \mv(w)  }}{\log{|\varphi'(0)| - b\sigma}}.
\end{equation}
This in particular implies (see, for example, Corollary~\ref{crude indeed}) that the coefficient $e$ in~\eqref{original UDC bound} and~\eqref{founding hol} can be taken down to the better constant~$2$.
Bost and Charles's work has been a major stimulus for our exploration of the applications to irrationality. 
Inspired by~\cite{BostCharles}, but going
outside of their framework of formal-analytic arithmetic surfaces and incorporating an idea of Perelli and Zannier~\cite{PerelliZannier}, we prove in \S~\ref{app:PerelliZannier} the reduction $e \rightsquigarrow 2$ in~\eqref{founding hol}. In~\S~\ref{new slopes}, we carry this further based on
some of Bost and Charles's results from~\cite{BostCharles} re-interpreted for analytic purposes into  Bost's prior method of evaluation heights,  in order to generalize~\eqref{BC integral} and~\eqref{BC integral 2} to incorporate a refined denominator term; see Theorem~\ref{basic main} for a special case of our bounds from \S~\ref{new slopes}.
Our companion treatise~\cite{zeta5} of the irrationality of 
the $2$-adic zeta value $\zeta_2(5)$ explores these bounds in a wider context.

\subsection{A preview of the various holonomy bounds}\label{sec_preview}
The core of our present paper consists of refined holonomy bounds that improve~\eqref{BC integral 2} and unify the proof methods behind~\eqref{founding hol} and~\eqref{BC integral}. One aspect of these bounds is to improve the $\tau = b\sigma$ term in~\eqref{founding hol}; here the high-dimensional methods ultimately yield a more precise information, although the difference is invisible to all our applications in this paper. The other aspect is to carry out a more refined complex analytic estimate (see \S~\ref{high dimensional analysis} and~\S~\ref{varia} for a summary of ideas) to further improve the double integral in~\eqref{BC integral}; here the improvements are the same in the single variable as in the high-dimensional treatments.
One technical novelty is a probabilistic input from large deviations theory which accommodates the $e \rightsquigarrow 2$ lowering in~\eqref{original UDC bound} even in the 
elementary multivariable framework of our original analysis in~\cite[\S~2]{UDC}.
    This is established through a Diophantine approximation argument in~$d$ auxiliary variables, and the point of achieving the $e \rightsquigarrow 2$ coefficient improvement
in precisely this way is that the high-dimensional geometric features of the $d \to \infty$ asymptotic make an additional room for further independent improvements. The sharpest holonomy bound (Theorem~\ref{high dim BC convexity}) that we have in this paper is a product of the measure concentration feature in the high-dimensional evaluation module. 

For the applications in this paper, including Theorems~\ref{mainA} and~\ref{logsmain}, the finest improvement concerning the general denominators does not make a difference.  We have 
two general simplified lines to these theorems. 
One is via Theorem~\ref{main:elementary form} using the high-dimensional techniques in a basic Siegel lemma framework, but another is via
Theorem~\ref{main:BC form} and alternatively Theorem~\ref{main: easy convexity}
using single variable methods.
  (See~\S~\ref{sec:leitfaden} for more details on the dependencies between
different sections of this paper, and the various paths to Theorems~\ref{mainA} and~\ref{logsmain}.)
For the application to Theorem~\ref{mainA}, Theorem~\ref{main:BC form} gives the weakest passable bound 
(sufficient by only a narrow margin)
compared to these two other theorems; while its ``convexity refinement,'' Theorem~\ref{main:BC convexity}, gives a stronger bound than either of them. 
The proof of Theorem~\ref{main:BC form} is a direct combination of the work of Bost and Charles, together with our improvement of the $\tau=b\sigma$ term in a relatively simple setting (see Theorem~\ref{basic main}; this simple setting allows us to get the optimal improvement of $\tau$ even without a high dimensional method), and a computation in~\S~\ref{integration cost} to accommodate added powers of~$n$ in the denominator types~\eqref{den type integrated}. 

To get the stronger bounds that handle Theorem~\ref{mainA} by  a more comfortable margin,
we use more refined complex analytic estimates to prove Theorems~\ref{main:elementary form} and~\ref{main:BC conv discrete}. In the case of the former, the 
large deviations input 
is used  not only to reach the improvement of the denominators rate term~$\tau$, but also to obtain a replacement of the Bost--Charles double integral by a more elementary \emph{rearrangement integral}\footnote{A \emph{rearrangement integral} here 
refers to a more general set of functions than the~$\log{|\varphi|}$ in the integrand of $\int_0^1 2t \cdot (\log{|\varphi(e^{2\pi i t})|})^* \, dt$. The latter, as we will see in~\S~\ref{rmk_Nazarov}, is larger than the Bost--Charles double integral. In general, we replace the $\varphi$ inside this integrand by a piecewise weighted combination of the functions $z \mapsto \varphi(rz)$, using a suitable set of radii $r$ that facilitate our refined complex analytic estimate.} 
which we introduce in \S~\ref{sec:variation};
   the proofs here 
are fully independent of
\cite{BostCharles}. On the other hand, based on \cite{BostCharles} and Theorem~\ref{main:BC form}, we undertake a closer study of the optimal archimedean estimates for the heights of the evaluation maps in Bost's slopes framework, and employ these improvements to prove Theorem~\ref{main:BC conv discrete}. This is what we dub the \emph{improvement from convexity}, a choice of terminology that refers to a classical theorem in the value distribution theory of meromorphic functions: the Nevanlinna characteristic~$T(r,\varphi)$ of a meromorphic function~$\varphi$ is a convex increasing function of~$\log{r}$. Further, if we choose a certain heuristically optimal Hermitian structure on the evaluation module of auxiliary polynomials, 
the argument of Theorem~\ref{main:BC conv discrete} leads to a heuristically optimal bound which we formulate as Theorem~\ref{main:BC fullconv}, still using single variable methods. In the basic denominators capping such as we introduce already  in Theorem~\ref{basic main} further down in this introduction, 
Theorem~\ref{main:BC fullconv} is the same as Theorem~\ref{high dim BC convexity} (cf. Remark~\ref{equalityoftaus}),
and we expect (cf. Remark~\ref{stationary choice}) both to give a stronger bound than Theorem~\ref{main:elementary form}.

\medskip

We proceed now to describe some of these basic improvements, and then state a first form of our new holonomy bounds. 

\subsection{Variants of the Nevanlinna growth characteristic}  \label{sec:variation}
From the starting bound~\eqref{founding hol}, on further pursuing~\cite[Remark 2.3.3]{UDC}, the multiple variables naturally  improve Nevanlinna's growth characteristic term $ \int_0^1 \log^+{|\varphi(e^{2\pi i t})|} \, dt$ 
 to the manifestly smaller \emph{rearrangement integral}
$\int_0^1 t \cdot (\log{|\varphi(e^{2\pi i t})|})^* \, dt$; here and throughout our paper, we follow the classical analysis custom to designate by 
\begin{equation} \label{incr}
g^*(t) := \inf_{s \in \R} \left\{ \mathbb{P}\left( x \in (0,1) \, : \, g(x) > s \right) \leq t \right\} =  \inf_{s \in \R} \left\{  \int_0^1 \chi_{g^{-1}([s,\infty))} \, dt \leq t   \right\}
\end{equation}
 the \emph{increasing rearrangement}
of a measurable function $g : (0,1) \to \R$. (See Basic Remark~\ref{increasingrearrangement}.)
This is the unique\footnote{Up to functions vanishing outside of a set of measure zero. Some authors prefer to use the term \emph{nondecreasing rearrangement function}. } nondecreasing measurable function that has the same distribution function as~$g$. We thus have
\begin{equation}
\label{maxintegral}
\int_0^1 2t \cdot g^*(t) \, dt = \int_0^1 \int_0^1 \max(g(s),g(t)) \, ds \, dt,
\end{equation}
inviting a comparison to the Bost--Charles double integral term from~\eqref{BC integral}. We will see in~\S~\ref{rmk_Nazarov} that the latter
is always, and in practice only slightly, smaller than the former. 

 It is, however, the left-hand
side of~\eqref{maxintegral} that arises naturally in the probabilistic character of our new argument. 
For our discussion here it suffices to note the trivial inequality 
\begin{equation}
\label{trivialnequality}
\int_0^1 2t \cdot g^*(t) \, dt \leq \int_0^1 2 \max(g^*(t),0) \, dt = \int_0^1 2 \max(g(t),0)^* \, dt 
= 2 \int_0^1 \max(g(t),0) \, dt
\end{equation}
 for any measurable function~$g$, and so in particular this recovers the $e \rightsquigarrow 2$ coefficient reduction from a genuinely high-dimensional perspective following~\cite[\S~2]{UDC} which is in some sense an approach ``orthogonal'' to the single variable analyses of either~\cite{BostCharles} or~\S~\ref{new slopes}.
 (These latter approaches have an Arakelovian character, and carry their own and different refinement of the Nevanlinna growth characteristic, which in~\S~\ref{sec:BC convexity} we dub the \emph{Bost--Charles characteristic}.) 
 Now the point is that the $d \to \infty$ argument further allows for an analogous ``denominator increasing rearrangements'' improvement of the 
 term~$\tau = b\sigma$ in the extension~\eqref{founding hol} to $\Q\llbracket x \rrbracket$ functions. 
 Some such improvement 
  is essential for all our proofs of Theorems~\ref{mainA} and~\ref{logsmain}. We also do give a single variable treatment in~\S~\ref{new slopes} of the main results of~\S~\ref{fine section}. The high dimensional method, on the other hand, leads to an even more precise bound in the denominators aspect, a refinement that could be useful in further developments or applications of our method.

 \begin{basicremark} \label{increasingrearrangement}
 A basic way to understand the definition of~$g^*(t)$ in equation~(\ref{incr}) is as follows. Assume that~$g(t)$ is a continuous (and hence bounded) function on~$[0,1]$.
 If~$g(t)$ is monotonically increasing, then~$g^*(t) = g(t)$. For~$g(t)$ arbitrary, let~$g_n(t)$ for~$n \ge 1$ denote the piecewise constant step function which
takes the value~$g(k/n)$ on the interval~$I_k = [(k-1)/n,k/n)$ for~$k=1,\ldots,n$ (extending the final interval~$I_n$ to include~$1$). The functions~$g_n(t)$ converge uniformly to~$g(t)$ as~$n \rightarrow \infty$.
 Now let~$g^*_n(t)$ denote the  step function which
 is also constant in the~$n$ intervals~$I_1$, $\ldots$, $I_n$, except now taking the~$n$ respective values
 \[\{g(1/n), g(2/n), g(3/n), \ldots g(n/n)\}\]
 \emph{rearranged} in \emph{increasing order} (hence the name). Then~$g^*_n(t)$ is the increasing rearrangement of~$g_n(t)$, and the functions~$g^*_n(t)$ converge uniformly
 to~$g^*(t)$.  
 \endofremark
 \end{basicremark}

 \subsection{Arithmetic holonomy bounds, basic form} \label{sec:holbasic}
Our first main result is the following simultaneous strengthening of all the \emph{holonomy bounds} or arithmetic rationality or algebraicity criteria that we have explicitly stated so far.

\begin{thm} \label{basic main}
Consider two positive integers~$m,r \in \NwithoutzeroA$ and an $m \times r$ rectangular array of nonnegative real numbers $\mathbf{b} := \big( b_{i,j} \big)_{\substack{ 1 \leq i \leq m, \,
1 \leq j \leq r }}$, all of whose columns are of the form: 
$$
0 = b_{1,j} =\cdots = b_{u_j,j} < b_{u_j+1,j}= \cdots = b_{m,j}=: b_j,  \qquad \forall j = 1, \ldots, r,
$$
for some $u_j \in \{0, 1, \ldots, m\}$.
Let
$$
\sigma_i := b_{i,1} +\ldots + b_{i,r}, \qquad i = 1, \ldots, m
$$
be the $i$-th row sum, and define
\begin{equation} \label{equitau}
\tau(\mathbf{b}) := \frac{1}{m^2} \sum_{i=1}^{m} (2i-1) \sigma_i   =\sigma_m -\frac{1}{m^2}\sum_{j=1}^r u_j^2 b_j \in [0, \sigma_m]. 
\end{equation}
Further, consider a holomorphic mapping $\varphi : (\Db, 0) \to (\C,0)$ with derivative ({\it conformal size}) satisfying $|\varphi'(0)| > e^{\sigma_m}$.
 
Suppose there exists an $m$-tuple $f_1, \ldots, f_m \in \Q \llbracket x \rrbracket$ of $\Q(x)$-linearly independent formal functions with denominator types of the form
\begin{equation}   \label{den type}
f_i(x) = \sum_{n=0}^{\infty} a_{i,n} \frac{x^n}{[ 1, \ldots, b_{i,1} \cdot n] \cdots [1,\ldots, b_{i,r} \cdot n]}, \qquad a_{i,n} \in \Z, 
\end{equation}
such that $f_i(\varphi(z)) \in \C \llbracket z \rrbracket$ is the germ of a meromorphic function on $|z| < 1$, for all $i = 1, \ldots, m$. 
 Then we have the bound
\begin{equation} \label{new bound}
m  \leq  \frac{\displaystyle{
 \iint_{\T^2} \log{|\varphi(z) - \varphi(w)|} \, \mv(z) \mv(w)} }{  \log{|\varphi'(0)|} - \tau(\mathbf{b}) }. 
\end{equation}
If, moreover, all functions~$f_i$ are \emph{a priori} assumed to be holonomic, the 
condition $|\varphi'(0)| > e^{\sigma_m}$  
can be relaxed to  $|\varphi'(0)| > e^{\tau(\mathbf{b})}$.
\end{thm} 

With elementary methods based on the phenomenon of measure concentration in high dimensions, we prove directly in~\S~\ref{fine section} the following variant using the increasing rearrangement function: 
\begin{equation} \label{new bound 2}
m  \leq  \frac{\displaystyle{
 \iint_{\T^2} \log\left( \max (|\varphi(z)|, \varphi(w)| \right) \, \mv(z) \mv(w)} }{  \log{|\varphi'(0)|} - \tau(\mathbf{b}) }  
 =  \frac{\displaystyle{
 \int_0^1 2t \cdot ( \log{|\varphi(e^{2\pi i t})} )^* \,  dt}}{  \log{|\varphi'(0)|} - \tau(\mathbf{b}) }.   
\end{equation}

To highlight the similarity of $\tau(\bb)$ with the increasing rearrangement function that emerges from a simple probabilistic consideration, let us note (writing $\sigma_0 := 0$) that the weighted average in~\eqref{equitau} can be expressed in a form rather similar to~\eqref{maxintegral}:
\begin{equation} \label{finitary rearrangement}
\tau(\mathbf{b}) 
 =
\sum_{i=1}^{m} \sigma_i \int_{(i-1)/m}^{i/m} 2t \, dt = \int_0^1 2t \cdot \left( \sum_{i=1}^{m} (\sigma_i-\sigma_{i-1}) \chi_{[0,i/m]}(t)  \right) \, dt.
\end{equation}
Noting the monotonicity of the step function $\sum_{i=1}^{m} (\sigma_i-\sigma_{i-1}) \chi_{[0,i/m]}(t) $,  our requirement that $\mathbf{b}$ is column-wise nondecreasing serves as the counterpart for denominators of the increasing rearrangement function
$(\log{|\varphi|})^*$. As remarked above, we will see in~\S~\ref{rmk_Nazarov} that the Bost--Charles integral in~\eqref{new bound} can be tightly majorized by $\int_0^1 2t \cdot (\log{|\varphi(e^{2\pi i t})|})^* \, dt$, with the effect that the bound~\eqref{new bound} implies the bound~\eqref{new bound 2}. 
 For either of the rearrangement integrals~\eqref{new bound 2} and~\eqref{finitary rearrangement}, the integration weight~$2t$ arises as the cumulative distribution function of $\left( [0,1], \mu_{\mathrm{Lebesgue}} \right)$. 
One mechanism for both these improvements over~\eqref{founding hol} is held by the concentration of measure phenomenon~\S~\ref{concentration of measure}; it is explained in~\S~\ref{horizontal integration}. 

For the rather rudimentary shape of the denominator type form~\eqref{den type} in our statement of Theorem~\ref{basic main}, a single variable proof is nevertheless also possible, as we discover with the slopes method in \S~\ref{one variable slopes}. 
In that context, both the denominators rate term~$\tau(\bb)$  and the Bost--Charles double integral term in~\eqref{new bound} emerge from the computation of the  covolume of the Euclidean lattice of auxiliary polynomial functions  chosen in the usual Diophantine analysis proof scheme: the former as the minimizer of a multivariable quadratic form arising from a basic template sought for the integral structure, and the latter as the ``infinite part'' based on a combination (due to Bost and Charles~\cite[\S~5]{BostCharles}) of the Poincar\'e--Lelong formula in complex analysis and the arithmetic Hilbert--Samuel formula in Arakelov theory.

\subsection{Siegel's  \texorpdfstring{$G$}{G}-functions}
As discussed above, Theorem~\ref{basic main} has a crude qualitative corollary which we may read as an arithmetic holonomicity criterion. 
It is due to Andr\'e~\cite[\S~VIII 1.6]{AndreG} (where the set of places $V$ in \emph{loc. cit.} must be assumed to be finite); 
in a slightly different context, the first holonomicity result of such a kind is probably the one
discovered by Perelli and Zannier~\cite[Thm. 1~B]{PerelliZannier}. 

\begin{cor}  \label{holonomic criterion}
If a formal function $f \in \Q \llbracket x \rrbracket$ has rational coefficients of the form
\begin{equation} \label{gen den form}
f(x) = \sum_{n=0}^{\infty} a_{n} \frac{x^n}{[ 1, \ldots, b_{1} n] \cdots [1,\ldots, b_{r}  n]}, \qquad a_{n} \in \Z
\end{equation}
and admits an analytic mapping $\varphi : (\D,0) \to (\C,0)$ with conformal size
$|\varphi'(0)| > e^{b_1 +\ldots + b_r}$ and such that the composite function 
germ $f(\varphi(z)) \in \C \llbracket z \rrbracket$ is the germ of a meromorphic function on~$\D$, then 
$f(x)$ is a \emph{holonomic function}: there exists a nonzero linear differential operator $\LL$ with $\Q[x]$-coefficients
that satisfies $\LL(f) = 0$.
\end{cor}
In this paper, we will exhibit and exploit such $f \in \Q \llbracket x \rrbracket$ whose holonomicity can be recognized
by this criterion. The special form of the denominators~\eqref{gen den form} then situates us
more specifically into the context of Siegel's theory of $G$-functions; in particular, see Remark~\ref{remark on global nilpotence} for 
a discussion,
 the linear differential operator~$\LL$ can {\it a posteriori} be taken to be of the Fuchsian class
 with only regular singular points and with rational exponents~\cite[III 6.1, VII 2.1, and VIII 1.5]{Dwork}.
 A major open question, which is closely related to the discussion of~\S~\ref{integral bounds} with implications to irrationality proofs and effective Siegel integral points problems, is to control the possible 
 \emph{apparent} singularities of the linear differential operator~$\LL$ in a minimal-order \emph{inhomogeneous} ODE
 $\LL(f) \in \Q[x]$.

 \begin{basicremark} \label{basiclog}
 A simplest example is $f(x) = \log(1-x)$, with type given by $(b_1,\ldots,b_r) = (1)$ and minimal differential operator
 $\LL := (1-x) (d/dx)^2 - (d/dx)$, varying holonomically on the domain $\Omega = \C \setminus \{1\}$ to define 
 a rank-$2$ local system
 $$
 \Cspan_{\C} \big\{1, \log(1-x)\big\}
 $$
 on~$\Omega$ with monodromy the infinite cyclic group generated by the unipotent matrix
 $$
T := \left( \begin{array}{cc}  1 & 0 \\ -2 \pi i & 1 \end{array} \right). 
 $$
 This expresses the fact that the analytic continuation process~$T$  ---  the \emph{local monodromy operator}  ---  for $\log(1-x)$ under the counterclockwise direction along a simple closed loop 
 encircling the singularity $\{1\}$ leaves $f_1 := 1$ invariant but adds to $f_2 := \log(1-x)$ the \emph{period} $-2\pi i$ times $f_1$:
 $$
 T^k(\log(1-x)) = \log(1-x) - 2\pi i k, \qquad T^k(1) = 1. 
 $$
 This holonomic example is furthermore recognized as a case of the holonomicity criterion Corollary~\ref{holonomic criterion}, for instance with
 the multivalent choice $\varphi(z) := 1 - e^{-Rz}$ for any $R > e$, or  the multivalent choice $\varphi(z) := \lambda(z)$ with $|\varphi'(0)| = 16 > e$, 
 or the univalent choice $\varphi(z) := 4z/(1+z)^2$ with $|\varphi'(0)|  = 4 > e$.  
 
 In Theorem~\ref{basic main}, the denominators type is captured by the~$2 \times 1$ matrix~$\bb = (0,1)^{\mathrm{t}}$,
 with~$\tau(\bb) = (1 \cdot 0 + 3 \cdot 1)/2^2 = 3/4$. For the choice~$\varphi(z) := 4z/(1+z)^2$, the holonomy quotient 
 is~$\log{4} /( \log{4} - 3/4) \approx 2.1787$, an upper bound 
 on the dimension~$m=2$ of this local system. 
 \endofremark
\end{basicremark}

 In~\S~\ref{Zagier local system}, we will make a thorough study of Zagier's holonomic functions~\cite{Zagier} that endow the numbers
 $\zeta(2)$ and $L(2,\chi_{-3})$ similarly as periods in a much more complicated local system spread over the domain $\Omega = \C \setminus \{0, 1/9, 1\} \cong \H / \Gamma_0(6)$. 
 For this local system, which emerged from analyzing the form of the recursion from Ap\'ery's $\zeta(2)$ irrationality proof~\cite{AperyHistoric,CohenApery,Apery} and is based on the theory of Eichler integrals, we will now have
 the main integrality type $x^n/[1,\ldots,n]^2$. We will then reduce the $\Q$-linear independence problem of $1, \zeta(2), L(2,\chi_{-3})$ to a  
 Diophantine analysis problem on the nonexistence of a $G$-function of the type $x^n/[1,\ldots,n]^2$ and with certain analytic properties: specifically, 
 our task becomes to prove that Zagier's local
 system cannot contain a nonzero $\Q\llbracket x \rrbracket$ element which is regular  ---  \emph{overconvergent}  ---  at the singularities $\{0,1/9\}$. 

A direct application of~\eqref{new bound}, see~\S~\ref{bivalent app} further down in this section, suffices for proving 
the irrationality of the mixed period
$$L(2,\chi_{-3})  - \pi \frac{\log{3}}{3\sqrt{3}} = L(2,\chi_{-3}) -    L(1,\chi_{-3}) \log 3.$$
(Another  irrationality result for a mixed period
is Beukers's proof~\cite[Thm~4]{Beukers} using modular forms that~$\zeta(3)  - 5 \sqrt{5} \, L(3,\chi_{5}) \notin \Q(\sqrt{5})$.) 
For the irrationality proof of the pure $L(2,\chi_{-3})$, as discussed in \S~\ref{sec_preview}, we need an even finer result than this to also take into account the integrals of the functions.
These more elaborate versions of Theorem~\ref{basic main}
(including Theorems~\ref{main:elementary form}, ~\ref{main:BC form}, \ref{main:BC conv discrete}, and~\ref{main: easy convexity}) are
deferred to \S~\ref{fine section} and \S~\ref{new slopes} below where they are proved.
The particular application to Theorem~\ref{mainA} is fairly delicate,
and among the many local systems generating $\zeta(2)$ and $L(2,\chi_{-3})$ among their holonomic coefficients, the choice that ends up working for us is highly
reducible (although with nonsolvable monodromy) 
and involves integrations that lead to denominators essentially\footnote{More precisely,
of the form~$n [1,\ldots,2n+3]^2$, but this can more or less be treated as having
the shape~$n [1,\ldots,2n]^2$, by Remark~\ref{rem:overflow}.}
of the form $n[1,\ldots,2n]^2$.

\subsection{Univalent holonomy bounds and an arithmetic characterization
of the logarithm}    \label{sec:logarithmchar}

We now consider the specialization of Theorem~\ref{basic main} to the setting
where the map~$\varphi$ is univalent. We remark that although, for general $\varphi$, we have various improvements of Theorem~\ref{basic main}, such as Theorems~\ref{main:BC conv discrete} and~\ref{main:BC fullconv} (assuming $\be$ in \emph{loc.~cit.} is $\mathbf{0}$), in the case of univalent~$\varphi$, all these  reduce to the same Theorem~\ref{univalent case} below.

For  $\Omega \subset \C$  a contractible domain containing~$0$, the Riemann mapping theorem supplies a biholomorphic map $\varphi : \D 
  \iso  \Omega$ with $\varphi(0) = 0$, which by Schwarz's lemma is uniquely defined up to pre-composing by a circle rotation. That makes the absolute value $|\varphi'(0)| \in (0,\infty]$ well-defined; we denote it by $\rho(\Omega,0)$ and call it \emph{the conformal mapping radius of the pointed contractible domain~$(\Omega,0)$}. The holomorphic mapping $\varphi : \D \to \C$ is said to be \emph{univalent} if it is biholomorphic onto its image, or equivalently, if $\varphi : \D \hookrightarrow \C$ is injective.

\begin{theorem}[Univalent holonomy bound] \label{univalent case}
Under the notations and assumptions of Theorem~\ref{basic main}, consider  $\Omega \subset \C$ a contractible domain with $0 \in \Omega$ and having a conformal
mapping radius $\rho(\Omega,0) > e^{\tau(\bb)}$. For any $m$-tuple of $\Q(x)$-linearly independent formal functions of the type~\eqref{den type} and meromorphic in~$\Omega$, the following holonomy bound holds: 
$$
m \leq \frac{\log{\rho(\Omega,0)}}{ \log{\rho(\Omega,0)} - \tau(\bb)}. 
$$
\end{theorem}

\begin{proof}  This follows directly from Theorem~\ref{basic main}.
The point to observe is that the Bost--Charles double integral term satisfies the inequality
 $$
 \iint_{\T^2} \log{|\varphi(z) - \varphi(w)|} \, \mv(z) \mv(w) \geq \log{|\varphi'(0)|}, 
 $$
\emph{with equality if and only if $\varphi : \D \hookrightarrow \C$ is univalent on the open disc}. To see this, simply observe that the 
univalence is equivalent to having the bivariate holomorphic function
$$
\frac{\varphi(z) - \varphi(w)}{z-w} = \varphi'(0)   + O\left(|z|+|w|\right)  \in \mathcal{O}(\D^2)
$$
to be nonvanishing throughout the unit polydisc. Hence the function
$$
G(z,w) := \log{\left|\frac{\varphi(z) - \varphi(w)}{z-w} \right| }  \, : \,  \D^2 \to \R \cup \{-\infty\}
$$
is plurisubharmonic, and harmonic if and only if $\varphi$ is univalent.
Both claims now follow upon remarking that $G(0,0) = \log{|\varphi'(0)|}$ while 
$$
\iint_{\T^2} G(z,w) \, \mv(z)\mv(w) =  \iint_{\T^2} \log{|\varphi(z) - \varphi(w)|} \, \mv(z) \mv(w),
$$
by the basic integral $\iint_{\T^2} \log\frac{1}{|z-w|} \, \mv(z)\mv(w) = 0$. 
\end{proof}

We note the following application of Theorem~\ref{univalent case} to the logarithm
function, which  is the example of Basic Remark~\ref{basiclog}.

\begin{theorem}  \label{logcharacterization}
Suppose $f(x) = \sum_{n=0}^{\infty} a_n x^n \in \Q\llbracket x \rrbracket$ is
a power series such that:  
\begin{enumerate}
\item \label{logdens} $[1,\ldots,n] a_n \in \Z$ for all $n \in \NwithzeroA$. 
\item $f(x)$ is holomorphic on $\C \setminus [1,\infty)$.
\end{enumerate}
Then 
$$
f(x) = Q_0(x) + Q_1(x) \log(1-x)
$$
for some rational functions $Q_0, Q_1 \in \displaystyle{\Q \left[ x,  \frac{1}{1-x}  \right]} \subset \Q(x)$. 
\end{theorem}

We view Theorem~\ref{logcharacterization} as an \emph{arithmetic
characterization} of the logarithm function.

\begin{proof} We consider the contractible domain $\Omega := \C \setminus [1,\infty)$, of 
conformal mapping radius $\rho(\Omega,0) = 4$ with the Riemann map $\varphi(z) = 4z/(1+z)^2$. Applying Theorem~\ref{univalent case} with $m=3, r = 1$, and $\mathbf{b} = (0,1,1)^{\mathrm{t}}$ with $\tau(\mathbf{b}) = 8/9$, the numerology
\begin{equation}  \label{first numerology}
\frac{\log{4}}{ \log{4} - 8/9 }  = 2.787050\ldots < 3
\end{equation}
proves that there is no third such function $\Q(x)$-linearly independent from the two known examples $f_1 = 1$ and $f_2 = \log(1-x)$ for the type~\eqref{gen den form}
with $(b_1, \ldots, b_r) = (1)$ and analytic on $\Omega = \C \setminus [1,\infty)$. This means that all 
such examples are of the form $Q_0(x) + Q_1(x) \log(1-x)$ with $Q_0(x), Q_1(x) \in \Q(x)$. 

At this point, we know that $f(x)$ is regular (holomorphic) on $\C \setminus  \left( [1,\infty) \cap \overline{\Q} \right)$, and that every point $x \neq 1$ in~$\C$ is at worst a meromorphic pole of~$f(x)$.  It remains to prove two things: 
\begin{enumerate}[label=(\roman*)]
\item \label{x0 part} $Q_0(x)$ and~$Q_1(x)$ are from the subring $\Q\left[ x, \frac{1}{x}, \frac{1}{1-x} \right]$ of~$\Q(x)$. 
\item \label{no x0} It is impossible to have~$Q_0(x),Q_1(x) \in \Q[x,1/x]$ without having both $Q_0(x),Q_1(x) \in \Q[x]$. 
\end{enumerate}
Indeed,~\ref{no x0} gives what we want assuming~\ref{x0 part} and upon changing~$f(x)$ to~$(1-x)^kf(x)$ with a sufficiently 
high power~$k \in \NwithzeroA$ to clear the~$(1-x)$ denominators from $Q_0(x)$ and~$Q_1(x)$.

We first prove~\ref{no x0}.
Suppose~$Q_0(x)$ and~$Q_1(x)$  are not both in~$\Q[x]$.
 If
 \[Q_1(x) \in \Q[x] \subset (1/x)\Q[x],\]
  then~$Q_1(x) \log(1-x)$ is holomorphic at~$x=0$, but then $Q_0(x) = f(x)-Q_1(x) \log(1-x)$ 
  is also holomorphic at~$x=0$   and then~$Q_0(x) \in \Q[x]$. Hence we may assume that~$Q_1(x) \in \Q[x,1/x] \setminus \Q[x]$.
 After multiplying~$f(x)$ by the correct power of~$x$ and a suitably divisible positive 
integer, we may assume that~$Q_1(x) \in \left( \Z[x] + \Z \cdot x^{-1} \right) \setminus \Z[x]$ and~$Q_0(x)$
(which is now holomorphic by the argument above) lies in~$\Z[x,1/x]$ and hence also in~$\Z[x]$,
and still with the denominator property~\eqref{logdens} in place. 
 In turn, upon subtracting from~$f(x)$ a suitable
element of~$\Z[x] + \log(1-x) \Z[x]$, we are left with analyzing the case
\[f(x) =  \frac{q_1 \log(1-x)}{x}\]
 with~$q_1 \in \Z \setminus 0$.
But then the~$x^{p-1}$ coefficient of~$f(x)$ is equal to~$q_1/p$, which, when~$p > |q_1|$ is a prime, is not of the required form~\eqref{logdens}. 
This completes the reduction step~\ref{no x0}.

We now consider~\ref{x0 part}.  
Suppose for contradiction that the rational functions $Q_0(x)$ and~$Q_1(x)$ are not from the subring $\Q\left[ x, \frac{1}{x}, \frac{1}{1-x} \right]$; then at least one of them will have a pole $\alpha \in \overline{\Q} \setminus \{0,1\}$. 

Fix a complex embedding~$\Qbar \hookrightarrow \C$, and consider firstly the case that $\alpha \notin (1,\infty)$ for at least one of the poles of $Q_0(x)$ or~$Q_1(x)$. 
In that case, our assumption that $f(x)$ is holomorphic at~$\alpha$ implies
that
\[Q_0(x) = \frac{V(x)}{(x-\alpha)^k}, \quad Q_1(x) = \frac{-U(x)}{(x-\alpha)^k}\] 
with some \emph{positive} integer $k \in \NwithoutzeroA$ and some rational functions
$U,V \in \overline{\Q}(x)$ regular and nonzero at~$x=\alpha$. Setting $x = \alpha$ in the equation
 \[(x-\alpha)^kf(x) = V(x) - U(x) \log(1-x)\]
  yields
a nontrivial vanishing combination $V(\alpha) - U(\alpha) \log(1-\alpha) = 0$ with nonzero algebraic number coefficients $U(\alpha), V(\alpha) \in \overline{\Q}^{\times}$. 
But this contradicts the Hermite--Lindemann--Weierstrass theorem on transcendental values of the function $\log(1-x)$ on~$\overline{\Q} \setminus \{0, 1\}$. 

It remains to handle the case that \emph{all} poles $\alpha \neq 0,1, \infty$ of $Q_0$ and $Q_1$ 
belong to $\alpha \in (1,\infty) \cap \overline{\Q}$, and that this set of poles is nonempty. Here, our $f \in \mathcal{O}(\C \setminus [1,\infty))$ holomorphy condition does not rule 
out a meromorphic pole at $x=\alpha$, and we need a different argument. As the set of poles in consideration is stable under $\Gal(\overline{\Q}/\Q)$, our assumption implies that all Galois conjugates of $\alpha$ lie in $(1,\infty) \cap \overline{\Q}$.
We then deduce from the product formula that there is a prime~$p$
and a choice of a pole $\alpha \in \overline{\Q} \cap \C_p$ lying within the open disc $|x|_p < 1$. Then we get the same contradiction $p$-adically, upon citing\footnote{The $p$-adic counterpart of the full Hermite--Lindemann--Weierstrass theorem on the algebraic independence of special values of the exponential
function is a well-known and still-unresolved conjecture. We refer to Nesterenko's work~\cite{NesterenkoPadic}, for partial results, and~\cite[\S~2.4]{NesterenkoMahler}, for an overview of the subject and an introduction to Mahler's argument.} Mahler's theorem~\cite{MahlerLindemann} 
 on the transcendence of all convergent values of the $p$-adic exponential function
at nonzero algebraic arguments; which is equivalent to the transcendence of all values of the $p$-adic logarithm function $\log(1-x)$ at the algebraic points of the punctured open unit disc 
$0 < |x|_p < 1$. 
\end{proof} 

\begin{remark}  \label{x=0 special role} Theorem~\ref{logcharacterization} and its proof
also holds with, for example,~\eqref{logdens} relaxed to 
the form
\[\left[1,\ldots,\left(1 + \frac{1}{100}\right)n\right];\]
but then with the weaker conclusion~$Q_0, Q_1 \in \Q[x,1/x,1/(1-x)]$ from step~\ref{x0 part} alone, where indeed~$1/x$ can no longer be removed, 
as instanced by the function~$f(x) = 100! \cdot\log(1-x)/x$.
Here the constant~$1+1/100$ could equally be replaced by any element of~$(1,(3 \log 2)/2) = (1,1.03972\ldots)$.~\endofremark
\end{remark}

\begin{remark} \label{integral module log} In the conclusion of Theorem~\ref{logcharacterization}, we can completely characterize
the possible~$Q_i(x)$. Namely, $f(x) = Q_0(x) + Q_1(x) \log(1-x)$
has the required form if and only if the following two conditions hold:
\begin{enumerate}
\item $Q_1(x) \in \Z[x,1/(1-x)]$.
\item $Q_0(x)$ lies in the~$\Z[x,1/(1-x)]$-module generated by~$x^n/[1,\ldots,n]$ for each~$n$ ---
 equivalently, generated by~$x^{q}/q$ for each prime power~$q$.
 \end{enumerate}
 This gives the full description by generators and relations of the (infinite) $\Z[x,1/(1-x)]$-module
 of solutions in Theorem~\ref{logcharacterization}. 
 
 It is plain that these conditions yield the requirements of Theorem~\ref{logcharacterization}. 
To prove the converse, consider an~$f(x) = Q_0(x) + Q_1(x) \log(1-x)$ in the theorem. 
 The conclusion for~$Q_0(x)$ is clear once we establish the conclusion for~$Q_1(x)$. Without loss of generality
(after multiplying by a power of~$(1-x)$), it suffices to show that if~$Q_i(x) \in \Q[x]$, then~$Q_1(x) \in \Z[x]$.
If~$Q_1(x) \notin \Z[x]$, then there exists a prime~$p$ and a monomial~$q_{1,m} x^m$ of~$Q_1(x)$ such that the $p$-adic
valuation $\mathrm{val}_p(q_{1,m}) < 0$ is negative and minimal amongst the $p$-adic valuations of all coefficients of~$Q_1(x)$.
But now, if~$p^r > \deg(Q_0(x)), \deg(Q_1(x))$, it is easy to check that the $p$-adic valuation
of~$[1,2,\ldots,n] a_n$ is negative for~$n = p^r+m$ and~$a_n$ the coefficient of~$x^n$ in~$f(x)$.~\endofremark 
\end{remark}

\begin{remark} \label{reverse}
The resort to the Hermite--Lindemann--Weierstrass theorem and Mahler's (partial) $p$-adic analog is not accidental in the proof of Theorem~\ref{logcharacterization}. In fact, reversing the logic at least in part, the statement of the theorem implies, for example, the irrationality of $\log(1-1/n)$ for all integers $n \in \Z \setminus \{1\}$; for if this (archimedean) logarithm took a rational value~$p/q$, then
$$
f(x) :=   \frac{q\log(1-x) - p}{1 - nx}  = q \frac{\log(1-x) - \log(1-1/n)}{1 - nx} \in \Q \llbracket x \rrbracket \cap \mathcal{O}\left(\C \setminus [1,\infty)\right)
$$
would meet the integrality and holomorphy constraints in Theorem~\ref{logcharacterization}, but the rational functions~$Q_0(x), Q_1(x) \in \Q(x)$ in the expression
$f(x) = Q_0(x) + Q_1(x) \log(1-x)$ would be singular at $x = 1/n$, and thus definitely not from the ring~$\Q\left[ x, \frac{1}{1-x} \right]$. We shall return to this type of issue in~\S~\ref{integral bounds}.  \endofremark
\end{remark}

\subsubsection{Theorem~\ref{univalent case} as a refinement of the Borel--P\'olya--Zudilin rationality criterion} \label{univalent det}
We make three remarks about Theorem~\ref{univalent case}. First, in the discussion in Basic Remark~\ref{Polya and Andre} we do indeed recover the more precise rationality statement
in the original Borel--P\'olya theorem, for we can have $ \tau(\mathbf{b})  = \tau(\mathbf{0}) = 0$ in that setting. 
Second, 
on a given simply connected domain $\Omega \ni 0$ of the complex plane with conformal mapping radius  $\rho(\Omega,0) > 1$, all transcendental $\Q\llbracket x \rrbracket$ formal function germs with a denominator type of the form
\begin{equation} \label{s den type}
f(x) = \sum_{n=0}^{\infty} a_{n} \frac{x^n}{[ 1, \ldots, b_{1}  n] \cdots [1,\ldots, b_{r} n]}, \qquad a_{n} \in \Z, 
\end{equation}
must meet  a denominator type gap
\begin{equation}  \label{den gap}
b_1 +\ldots + b_r \geq (2/3) \log{\rho(\Omega,0)}. 
\end{equation}
If there are at least~$m \geq 2$ such $\Q(x)$-linearly independent functions, the coefficient~$2/3$ in~\eqref{den gap}
improves to $m/(m+1)$. 

Finally, in the most basic situation of all taking~$\Omega$ to be a round disc $|x| < R$ 
of  a radius~$R > e^{\sigma}$ centered at~$0$ (as considered by Ap\'{e}ry, except that now --- like Borel --- we assume
meromorphy rather than holomorphy), we explain how Theorem~\ref{univalent case} implies that~$f(x) \in \Q(x)$.
Applying Theorem~\ref{univalent case} directly, we deduce to start with that the corresponding~$\Q(x)$-vector
space~$\HH$ generated by such functions is finite dimensional, and in particular
consists of holonomic functions.
However,
 if~$f(x) = \sum a_n x^n \in \Q \llbracket x \rrbracket$
is meromorphic on~$\Omega$, then, with~$\zeta = e^{2 \pi i/m}$,  so are the twists
$$\frac{1}{m} \sum_{i=0}^{m-1} f( \zeta^i x) \zeta^{- i k} = \sum a_{mn+k} x^{nm + k},$$
and those have the same denominator type as~$f(x)$. It follows that~$\HH$ 
is preserved by~$x \mapsto \zeta x$ for any~$m$. 
The (non-apparent) singularities of the corresponding differential equation  cannot be invariant under all 
these rational  rotations unless 
they are a subset of~$\{0,\infty\}$. But this implies that any such~$f(x)$ must be meromorphic on~$\C$,
and (after clearing denominators) we may apply Theorem~\ref{univalent case} again, taking now~$R$ to be arbitrarily large,
to deduce that~$\dim_{\Q(x)}\HH = 1$.
(Note that there do exist  finite-dimensional~$\Q(x)$-vector spaces of  dimension greater than~$1$
which are generated
by $\Z \llbracket x \rrbracket$ holomorphic functions on $\D$  and are invariant under~$x \mapsto \zeta x$ for all rational rotations; for example, the~$\Q(x)$-vector
space generated by~$1$ and~$f(x) = \sum x^{n!}$. The latter, of course, is non-holonomic.)

We can summarize the three remarks by the following refinement of the Borel--P\'olya rationality criterion, 
and also of Zudilin's determinantal criterion~\cite{ZudilinDet}. 

\begin{thm}  \label{cap rationality}
Consider a contractible open domain $\Omega \ni 0$ in the complex plane and a formal power series of the arithmetic type
\begin{equation} \label{den type 3}
f(x) = \sum_{n=0}^{\infty} a_n \,  \frac{x^n}{[1,\ldots,b_1n] \cdots [1,\ldots, b_rn]}, \qquad a_n \in \Z ,\quad \forall n \in \NwithzeroA,
\end{equation} 
which is the $x=0$ germ of a meromorphic function on~$\Omega$. Suppose that either 
\begin{enumerate}[label=(\roman*)]
\item  \label{rounddisc} $\Omega$ is a round disc $|x| < R$ of a radius $R > \exp(b_1 +\ldots + b_r)$; or else that
\item \label{bigconformal} the conformal mapping radius $\rho(\Omega,0)$ of $\Omega$ at the origin exceeds
\[\exp\left( \frac{3}{2} (b_1 +\ldots + b_r) \right).\] 
\end{enumerate}
Then  $f(x) \in \Q(x)$ is the Taylor expansion 
of a rational function.   \hfill{$\square$}
\end{thm}

\subsection{Arithmetic characterizations beyond the logarithm}
\label{sec:beyondlog}

In light of Theorem~\ref{logcharacterization}, it is natural to inquire of arithmetic characterizations of other basic  transcendental functions in terms of
their domains of analyticity and the arithmetic behavior of their power series.
In view of Bely\u{\i}'s theorem~\cite[\S~12.3]{BombieriGubler}, a natural place to start is (as in Theorem~\ref{logcharacterization}) with power series
that can be analytically continued as \emph{multivalued}  holomorphic functions along all paths in~$\mathbf{P}^1 \setminus \{0,1,\infty\}$. 
Going further than the denominator type~$[1,\ldots,n]$ of Theorem~\ref{logcharacterization} requires to use a multivalent map~$\varphi$, but there is still
a local univalence input, discussed in~\S~\ref{local univalent leaves} and formalized in~\S~\ref{def:univalentleaf} and Corollary~\ref{stacky overconvergence}, which is essential for our approach to irrationality proofs. 

In any case, if~$\tau =\tau(\bb)$,  a necessary condition for our methods to have any hope of applying 
is that~$|\varphi'(0)| > e^{\tau}$. A theorem of Carath\'eodory~\cite[(412.8) on page~198]{Caratheodory} shows that $|\varphi'(0)| \leq 16$
for all holomorphic maps $\varphi : \D \to \C \setminus \{1\}$ subject to $\varphi^{-1}(0) = \{0\}$, with equality
holding if and only if $\varphi(z) = \lambda( cz )$ with $|c| = 1$. Hence, in this setting, it is necessary that~$\lambda'(0) = 16 > e^{\tau}$ (see \S~\ref{local univalent leaves} for details on why the specific assumptions in Carath\'eodory's theorem is relevant).
This necessary condition is certainly met by~$\tau = 2$. In particular, Corollary~\ref{holonomic criterion} implies that the~$\Q(x)$-vector
space of type $[1,\ldots,n]^2$ functions holonomic on $\P^1 \setminus \{0, 1, \infty\}$ is finite dimensional.
As we shall explain below,  our method of proof for both Theorems~\ref{mainA} and~\ref{logsmain} can be summarized as making a sufficient way towards the determination 
of that finite-dimensional space. 

\begin{conjecture} \label{conjecture five elements}
The following conditions on a formal power series $f(x) = \sum_{n=0}^{\infty} a_n x^n \in \Q\llbracket x \rrbracket$ convergent in $|x| < 1$ 
are equivalent: 
\begin{enumerate}
\item \label{five char} $f(x)$ is analytically continuable as a holomorphic function to~$\C \setminus [1,\infty)$, and furthermore \emph{as a meromorphic function} along all paths in
$\P^1 \setminus \{0,1,\infty\}$; and there is an $M \in \NwithoutzeroA$ such that $[1,\ldots,n]^2 a_n \in M^{-1} \Z$ for all $n \in \NwithzeroA$.
\item \label{solns five char} There are  rational functions 
$Q_0, \ldots, Q_4 \in \displaystyle{\Q \left[ x, \frac{1}{1-x}  \right]} \subset \Q(x)$ with
\begin{equation*}
\begin{aligned}
f(x) &  = Q_0(x) + Q_1(x) \log(1-x) + Q_2(x) \log^2(1-x) + Q_3(x) \Li_2(x)  \\
& +
\frac{ Q_4(x) }{\sqrt{1-x}} \int_0^x \frac{\log{(1-t)}}{t\sqrt{1-t}} \, dt.  \end{aligned}
\end{equation*}
\end{enumerate}
\end{conjecture}

Here, $\Li_2(x) := -\int_0^x \log(1-t) \, d\log{t} = \sum_{n=1}^{\infty} x^n/n^2$ is the standard dilogarithm function branch, and the hypothetical
solution space is discussed in more depth in~\S~\ref{sec:pureX2}.  
One should compare this conjecture to Theorem~\ref{logcharacterization}. 
In either case, one may consider a bipartite approach. The first part is to devise a setup in Theorem~\ref{basic main} that proves the finite-dimensionality of the~$\Q(x)$-vector
space of such functions; the second part is to give a bound for this space which coincides with the number
of known functions.
The fact that~$16 > e$ (respectively $16 > e^2$) establishes the first claim in either case.
In the second case, however, the best bound on the dimension we can currently establish is~$9$ rather than~$5$ (see Remark~\ref{improvement} and
 Equation~\ref{biggerthannine}).
Ruling any possible further functions out remains a difficult problem currently beyond the reach of our methods in this
paper.

\begin{remark} \label{contains in itself}
Similarly to Remark~\ref{reverse}, the~$\Q[x,1/(1-x)]$ refinement contains, like a hidden particular clause in this form of Conjecture~\ref{conjecture five elements},  the $\Q$-linear independence of the~$x=1/n$ special values of the five functions~$1$, $\log(1-x)$, $\log^2(1-x)$,
$\Li_2(x)$, and~$\frac{1}{\sqrt{1-x}} \int_0^x \frac{\log(1-t)}{ t \sqrt{1-t}} \, dt$ in the statement of the conjecture, for 
every~$n \in \Z \setminus \{0,1\}$. 
If for example $\Li_2(1/n) \in \Q$, then the point is that the function
$$
f(x) := \frac{\Li_2(x) - \Li_2(1/n)}{1-nx}  \in \mathcal{O}(\C \setminus [1,\infty)) \cap \Q\llbracket x \rrbracket
$$
meets all the conditions in~\eqref{five char} of the conjecture, but it is manifestly
not contained in the solution~$\Q[x,1/(1-x)]$-module prescribed by~\eqref{solns five char}. 

For all~$|n| \geq N_0$, where~$N_0$ is some (large, explicitly computable) number, the $\Q$-linear 
independence of the~$x=1/n$ special values of those five functions follows as a very particular case of
the general Theorem~\ref{BombieriChudnovsky} from the theory of special values of $G$-functions. 
For the dilogarithm function, the first such result was proved already by Maier~\cite[\S~8]{WMaier}, in a work that foreshadowed (and directly inspired) Siegel's~1929 paper~\cite{Siegel1929SNS}. 
The irrationality~$\Li_2(1/n) \notin \Q$ has at present only been proved~\cite{HataDilog,RhinViolaDilog1,RhinViolaDilog2} for $n \notin \{ -4, -3, -2, 2, 3, 4, 5 \}$. 
The issue is discussed further in~\S~\ref{integral bounds}.~\endofremark
\end{remark}

In the main spirit of our paper, one could even ask for variations of Conjecture~\ref{conjecture five elements} that allow for further (possible) singularities in the convergence disc $|x| < 1$ of the 
original branch $f(x)$, for example:

\begin{question} \label{noextras}
Do the conclusions of Conjecture~\ref{conjecture five elements} still hold if the meromorphic continuability on~$\P^1 \setminus \{0,1,\infty\}$ in condition~\eqref{five char} is relaxed to a meromorphic continuability on~$\P^1 \setminus \{0, \delta, 1, \infty\}$, for some $\delta \in [-1/2,1/2]$? 

In particular, for a power series $f(x) \in \Q\llbracket x \rrbracket$ convergent on $|x| < 1$ and defining holonomic functions on $\P^1 \setminus \{0, \delta, 1, \infty\}$ of the denominators type condition $a_n [1,\ldots,n]^2 \in \Z$, where $\delta \in [-1/2,1/2]$ is an arbitrary fourth puncture, does~$f(x)$ automatically extend through that fourth puncture~$x =\delta$  to
define a holonomic function on~$\mathbf{P}^1 \setminus \{0,1,\infty\}$? 
\end{question}

While these questions seem rather awkward for our method of rational holonomy bounds as developed in this paper, we are able to fully resolves a sub-problem intermediate in difficulty between Theorem~\ref{logcharacterization}
and Conjecture~\ref{conjecture five elements}, namely when the denominator type has the form~$[1,2,\ldots,n][1,2,\ldots,n/2]$, 
that is ``a case of $\tau = 3/2$'' where the first new function after $\log(1-x)$ pops out, namely, the function $\log^2(1-x)$. Here Conjecture~\ref{conjecture five elements} becomes the $\delta = 0$ case of the following theorem, responding affirmatively to Question~\ref{noextras} for the subcase of $[1,\ldots,n][1,\ldots,n/2]$ types: 

\begin{thm}
\label{thm three elements} \label{noextrassmall}
Suppose $f(x) = \sum_{n=0}^{\infty} a_n x^n \in \Q\llbracket x \rrbracket$ has $[1,\ldots,n][1,\ldots,n/2] a_n \in \Z$ for all $n \in \NwithzeroA$, 
is holomorphic in $\C \setminus [1,\infty)$, 
and is analytically continuable as a meromorphic function along all paths in  $\PP^1 \setminus \{0, \delta, 1, \infty\}$, 
for some $\delta
\in (-\infty,1)$.

 Then
\begin{equation}  \label{three form eq}
f(x) = Q_0(x) + Q_1(x) \log(1-x) + Q_2(x) \log^2(1-x)
\end{equation}
for some rational functions of the form~$Q_0,Q_1, Q_2 \in \Q\left[ x, \frac{1}{1-x} \right] \subset \Q(x)$. 

In particular, 
\begin{equation*}
(\ast) \qquad
\begin{aligned}
&f(x) \text{ continues analytically as a meromorphic function} \\
&\text{along all paths in } \mathbb{P}^1 \setminus \{0,1,\infty\}.
\end{aligned}
\end{equation*}
\end{thm}

Some immediate applications of part~$(\ast)$ of Theorem~\ref{thm three elements} to $\Q$-linear independence proofs are treated in~\S~\ref{mixed examples} further down in this introduction, as a proof-of-concept for our method. It is there that~$(\ast)$ is proved, 
as an application of Theorem~\ref{basic main}. To conclude the full Theorem~\ref{thm three elements} requires a subtler holonomy bound
and it is carried out in~\S~\ref{sec:completion three elements}.

\subsection{Overconvergence and univalent leaves} \label{local univalent leaves}
We now turn to the basic mechanism for irrationality proofs by extending the method of Ap\'ery limits. We will follow this in~\S~\ref{sec:first irrationality proofs} with
some explicit examples, and in~\S~\ref{bivalent app} with a proof-of-concept application to some new $\Q$-linear independence proofs. 

Consider $\Sigma \subset \D_R := \{ x \in \C \, : \, |x| < R \} $ a discrete subset of the open complex disc 
of radius~$R \in (0,\infty]$ (possibly including the disc center $x=0$), and   $f(x) \in \C\llbracket x \rrbracket$ a holomorphic function germ
at the center point that continues analytically as a holomorphic function along all paths in $\D_R \setminus \Sigma$. 
Let us define the subset $\Sigma_f^+ \subset \Sigma$, to necessarily include~$0$ if $0 \in \Sigma$, to 
consist of
those $\beta \in \Sigma$ for which the radial analytic continuation of $f(x) \in \C \llbracket x \rrbracket$ from $x=0$ towards $x = \beta$ remains bounded. We say that the 
power series $f(x)$ is \emph{overconvergent} at~$\Sigma_f^+$ and \emph{extends to~$\D_R \setminus \Sigma$ as a multivalued holomorphic function}. 
 Then the radius of convergence of the initial power series germ~$f(x) \in \C\llbracket x \rrbracket$ is equal to $\min_{\beta \in
 \left(\{R \} \cup \Sigma  \right)\setminus \Sigma_f^+} |\beta|$. The following trivial lemma is crucial for our approach to Theorems~\ref{mainA} and~\ref{logsmain}; we note that this type of statement on compatibility with integrations becomes completely false if we replace \emph{holomorphic} by \emph{meromorphic} everywhere in the previous paragraph. 
 
 \begin{lemma}  \label{integral holomorphy}
 There is an equality $\Sigma_{\int_0^x f(t) \, dt}^+ = \Sigma_f^+$. 
 \end{lemma}
  
    Given now a holomorphic mapping
$\varphi : \D \to \D_R$ with $\varphi(0) = 0$, we can apply the same notion to the pulled-back power series $f(\varphi(z)) \in \C \llbracket z \rrbracket$,
which is a $z=0$ holomorphic function germ that extends to $\D \setminus \varphi^{-1}(\Sigma)$ as a multivalued holomorphic function. 
In general, there is no relationship between the overconvergence sets $\Sigma_f^+ \subset \Sigma \subset \D_R$ and $\left( \varphi^{-1}(\Sigma) \right)_{\varphi^*f}^+ 
\subset \varphi^{-1}(\Sigma) \subset \D$ for $f$ and $\varphi^*f$.

But suppose there is a contractible open neighborhood $0\in \Omega \in \D$ on which $\varphi|_\Omega: \Omega \iso \varphi(\Omega)$ is univalent and, therefore, a conformal isomorphism onto the image open neighborhood~$\varphi(U) \ni 0$. Assume furthermore that $\varphi^{-1}(0) = \{0\}$ and that each 
point in $\varphi(\Omega) \cap \Sigma_f^+$ has exactly one pre-image under the analytic map~$\varphi$, that is: 
$$
\varphi^{-1}\left(\varphi(\Omega) \cap \Sigma_f^+ \right) \subset \Omega.
$$
 Then, in particular, $f(\varphi(z))$ is holomorphic on at least~$\Omega$: 

\begin{equation} \label{uni leaf}
\left( \varphi^{-1} \left( \Sigma \right) \cap \Omega  \right)_{\varphi^* f}^+ = \varphi^{-1}\left( \Sigma_f^+ \right) \cap \Omega = \varphi^{-1}\left( \Sigma_f^+ \cap \varphi(\Omega) \right). 
\end{equation}
This is the univalence input we alluded to. If now, in addition, $\varphi(\D) \cap \Sigma = \Sigma_f^+ \cap \varphi(\Omega)$, 
it follows
at once that the multivalued holomorphic function $f(\varphi(z))$  on $\D \setminus \varphi^{-1}(\Sigma) = \D \setminus  \varphi^{-1}\left( \Sigma_f^+ \cap \varphi(\Omega) \right) 
= \D \setminus \left( \varphi^{-1} \left( \Sigma \right) \cap \Omega  \right)_{\varphi^* f}^+$ is in fact a (single-valued) holomorphic function on the whole disc~$\D$, that is a convergent power series on that disc.

We summarize the basic property that we just proved: 

\begin{proposition}\label{overconvergence}
Let $f \in \C \llbracket x \rrbracket$ be a holomorphic function germ which extends as a multivalued holomorphic function on the Riemann surface $\P^1 \setminus \Sigma$, for some finite set of punctures $\Sigma$ on the Riemann sphere. Consider a disjoint partition $\Sigma = \Sigma^0 \sqcup \Sigma^1$, a holomorphic map $\varphi : \D \to \P^1 \setminus \Sigma^1$ that takes
$\varphi(0) = 0$, and a contractible open neighborhood $0\in \Omega \subset \D$ on which~$\varphi$ restricts as a univalent map (equivalently: $\varphi_{|\Omega} : \Omega \iso \varphi(\Omega)$ is a conformal isomorphism). We assume that $\varphi^{-1}\left( \Sigma^0 \right) \subset \Omega$ and that $f \in \mathcal{O}(\varphi(\Omega))$ is holomorphic on~$\varphi(\Omega)$. 

Then, the pulled-back germ $f(\varphi(z)) \in \C \llbracket z \rrbracket$ converges on the full disc $\D$.
\end{proposition}

\begin{remark}
The assumptions on the triple~$(\varphi,\Omega,\Sigma)$ in Proposition~\ref{overconvergence} can alternatively, and slightly more succinctly, be summarized
by having a holomorphic mapping $\varphi : (\D,0) \to (\C,0)$ that restricts univalently on the contractible open neighborhood~$\Omega \ni 0$,
and such that~$\varphi^{-1}(\Sigma) \subset \Omega$ for the finite puncture set~$\Sigma$. We chose the formulation with~$\Sigma = \Sigma^0 \sqcup \Sigma^1$ to highlight the practical presence of universal maps~$\varphi$ when the singularity type~$(\Sigma^0, \Sigma^1)$ is given but the open neighborhood~$\Omega \ni 0$ is kept unspecified. 
\end{remark}

\subsubsection{The modular lambda map}  \label{the lambda template}
In this general context, the significance of the modular lambda map~\eqref{lambdadef} is in the observation that $\varphi(z) := \lambda(z)$ is the universal map in  Proposition~\ref{overconvergence} 
for the case $\Sigma^0 = \{0\}$ and
$\Sigma^1 = \{1, \infty\}$ (upon keeping fluid the choice of an unspecified open neighborhood~$\Omega \ni 0)$. Its derivative $\lambda'(0) = 16$ therefore maximizes the conformal size of any such map. This can be considered (see Remark~\ref{Landen scheme} for a direct connection) as the multivalent analog of the role of the
domain $\Omega = \C \setminus [1,\infty)$ and the Koebe map $\varphi(z) = 4z/(1+z)^2$ in the proof of Theorem~\ref{logcharacterization}. Concretely, 
 if $f(x) \in \C \llbracket x \rrbracket$ continues analytically along all paths as a holomorphic function on~$\P^1 \setminus \{0, 1, \infty\}$ (an example
 is any balanced hypergeometric series),
then  $f(\lambda(z)) \in \C \llbracket z \rrbracket$ converges on the open unit disc~$z \in \D$. A basic illustration is the classic Jacobi formula
\[\sum_{n=0}^{\infty} \binom{2n}{n}^2 \left( \frac{\lambda(q)}{16}\right)^n = \left( \sum_{n=0}^{\infty} q^{n^2} \right)^2,\]
 where the holonomicity in~$x = \lambda(q)$ 
is an expression of the Picard--Fuchs ODE for the de Rham cohomology of the Legendre elliptic curve over
\[Y(2)_{\C} = \spec{\C}\left[x,1/x,1/(1-x)\right].\]
Our proof of Theorem~\ref{mainA} will involve a similar expression~\S~\ref{X06} of the Picard--Fuchs ODE over the modular curve
 \[Y_0(6)_{\C} \cong \spec{\C} \left[ x, 1/x, 1/(1-x), 1/(1-9x) \right],\]
  in which~$\zeta(2)$ and~$L(2,\chip)$ emerge as the Eichler periods.

\subsection{First irrationality proofs} \label{sec:first irrationality proofs}
In Remark~\ref{contains in itself} on Conjecture~\ref{conjecture five elements}, we observed that the prescribed~$\Q[x,1/(1-x)]$-module
 has direct irrationality implications on special values at points of the form~$x=1/n$. 
However, the method of our present paper only addresses~$\Q(x)$-vector spaces in the framework of
Theorem~\ref{basic main}, but not their integral structures over 
finitely generated $\Q$-algebras intermediate between~$\Q[x]$  and~$\Q(x)$. 
We now explain how even the cruder $\Q(x)$-form of Conjecture~\ref{conjecture five elements} (as enhanced by Question~\ref{noextras})
casts  a method for establishing irrationality proofs. These are now in the form of Ap\'ery limits, as opposed to the straight special values 
of the functions in the relevant holonomic module.

The following expands upon what
we have already discussed in the introduction. 
The ideal situation is as follows. Given an interesting period~$\eta$, one writes down a holonomic function~$f(x)$ with
coefficients in~$\Q(\eta)$. Assuming for the contradiction that~$\eta \in \Q$ and hence~$f(x) \in \Q \llbracket x \rrbracket$, this function (together with its derivatives)
provides a space of holonomic functions of some explicit denominator type
and dimension over~$\Q(x)$. In addition, depending on the circumstances, there will also exist other known
functions in this space. 
Considerations of monodromy (or otherwise) typically
allow one to
show that this space of known functions is~$\Q(x)$ (and even~$\C(x)$)-linearly independent from the 
functions coming from~$f(x)$. If the lower bound coming from  the span of such functions exceeds the
upper bounds from our theorem, we obtain the desired irrationality of~$\eta$.

Consider, for instance, our task to establish the~$\Q$-linear dependence~$1,\zeta(2)$, and $L(2,\chip)$. The ideal scenario would be to use a  putative $\Q$-linear dependence to write down such an~$f(x)$ with
denominators of type~$\tau = [1,2,\ldots,n]^2$ which extends holomorphically along all paths in~$\mathbf{P}^1 \setminus \{0,1,\infty\}$, but such that~$f(x)$ is \emph{not} in the~$\Q(x)$-vector space generated by the 
five functions in Conjecture~\ref{conjecture five elements}. Then Conjecture~\ref{conjecture five elements} would immediately
give a contradiction. This is not possible, but clearly we can get away with something weaker. As mentioned, we \emph{can}
prove a bound of~$9$ on the dimension of such functions. Now such a bound would still be sufficient as long as the span of~$f(x)$
and its derivatives were linearly independent from these five functions and gave a complementary~$\Q(x)$-vector space
of dimension at least~$5$.  In practice, even this fails in two respects. First, the function~$f(x)$
we construct only generates a holonomic module of dimension~$4$. Second, the function~$f(x)$ has additional
singularities at paths in~$\mathbf{P}^1 \setminus \{0,1,\infty\}$ to both~$\delta = 1/9$ and~$\delta = -1/8$.
It turns out that we can still bound the space of functions by~$9$ with these additional singularities, but the
numerology still falls \emph{just} short of our desired application. Instead, we have to additionally also include
\emph{integrations} of these (and other) functions into our story, and this is how we ultimately
achieve the proof of Theorem~\ref{mainA}, which is perhaps the most subtle of our applications.
It seems useful, however, to give examples where the approach as described above works directly,
first by reproving the (known by Lambert in~1761!) irrationality of~$\log{3}$, and then 
(in Theorem~\ref{mixed}) to devising a new irrationality result.

\begin{basicremark} 
\label{log 3 example}
Turning now to the main style of applications of holonomy bounds to irrationality proofs, the following
is a simple example due to Zudilin~\cite[\S~3]{ZudilinDet}, in which case~\ref{bigconformal}, but not case~\ref{rounddisc} of Theorem~\ref{cap rationality} 
provides an irrationality proof of the period $\log{3}$ out of the consideration of the integrals
\begin{equation}
\begin{aligned}
& \qquad  f(x) :=   \frac{1}{\sqrt{1-4x + x^2}} \int_{2-\sqrt{3}}^x \frac{dt}{\sqrt{1-4t + t^2}} \label{prod3}  \\ 
= & \  \frac{1}{2} \sum_{n=0}^{\infty} \left(   b_n   - a_n \log{3}    \right) \cdot   \frac{x^n}{2^n}, \qquad a_n \in \Z, \quad [1,\ldots,n] \, b_n \in \Z, 
\quad \forall n \in \NwithzeroA.
\end{aligned}
\end{equation}
For the second line, we use the binomial expansion of  $(1-4x+x^2)^{-1/2} \in \Z\llbracket x/2 \rrbracket$ 
and the fact that $-\log\left(2-x+\sqrt{1-4x+x^2}\right)$ is an explicit primitive of $1/\sqrt{1-4x+x^2}$. 

Of course, $f(x)$ is not a $G$-function as it has transcendental coefficients from involving~$\log{3}$; rather, it is a  $\C$-linear 
combination of two $G$-functions on $\P^1 \setminus \{ 2 \pm \sqrt{3}, \infty  \}$, and $\log{3} \in \C$ gets characterized as
the unique (holonomic) coefficient in such a combination to give a branch regular (holomorphic) at the smaller singularity
$x = 2-\sqrt{3}$. (This is rather transparently revealed by the fact that both factors in~\eqref{prod3} switch sign after a simple loop going around that singularity $x = 2-\sqrt{3}$, 
and thus their product has no monodromy at $x = 2-\sqrt{3}$.)

But we can turn this around and get an irrationality proof of $\log{3} \notin \Q$ as an application of Theorem~\ref{cap rationality} (and, hence, ultimately 
of the univalent holonomy bound). Proving $\log{3} \notin \Q$ means precisely proving that $f (2x) \notin \Q \llbracket x \rrbracket$. Suppose 
not. Then $f(2x) \in \Q \llbracket x \rrbracket$ has, upon clearing a fixed positive integer denominator, visibly 
the type~\eqref{den type 3} with $r = 1$ and $(b_1, \ldots, b_r) = (1)$. At the same time, by construction we have $f(2x)$
 holomorphic on the domain 
$$
x \in \Omega := \C \setminus \left[   \frac{2+\sqrt{3}}{2}, \infty  \right)
$$
of Riemann mapping radius 
$$
\rho(\Omega,0) = 2(2+\sqrt{3}) = 7.4641\ldots > 4.481689\ldots = e^{3/2}.
$$
Theorem~\ref{cap rationality}~(ii) proves that every such function~$f(2x)$ has to be a rational function. Obviously,
the putative function from~\eqref{prod3} (which would only have existed had $\log{3}$ been a rational number) is not rational, and so contrapositively
this argument gives a proof of the irrationality of~$\log{3}$. 

And yet,  as $2+\sqrt{3} = 3.73205\ldots  <  5.43656\ldots = 2e$, the $2^n [1,\ldots,n]$ denominators growth rate
in these approximations $\log{3} \approx b_n/a_n$ exceeds the reciprocal of the decay rate~$2-\sqrt{3}$ of the 
error $\left| \log{3} - b_n/a_n \right|$ of the approximations. In other words, case~\ref{bigconformal} applies in the theorem, whereas case~\ref{rounddisc} does not.    
Thus we find an irrationality proof, by $G$-function methods, without actually constructing any rapidly convergent explicit (\emph{holonomic}) rational
approximants. (A more complicated construction~\cite{Salikhovlog3,Sorokin} to pass the latter requirement is known in the case of~$\log{3}$, but not for say $\log{p}$
where $p$ is any sufficiently big prime.) 

As we will see, the usefulness of Theorem~\ref{basic main} lies in the possibility of using  ---  instead of domains $\Omega \subset \C$ as on this example  ---  \emph{multivalent} mappings
such as $\varphi(z) := \lambda(z)$,  a holomorphic function on~$\D$ whose derivative $|\varphi'(0)| = 16 > e^2$ fortuitously exceeds the growth
rate of the $[1,\ldots,n]^2$ layer of denominators common to several linear independence problems of interest here (including the case of Theorems~\ref{mainA} and~\ref{logsmain}), and 
which applies to the holonomic functions on $\P^1 \setminus \{0, 1, \infty\}$.   \endofremark
\end{basicremark}

\subsection{The multivalent case: first new linear independence results} \label{bivalent app}  In this section, we prove half of Theorem~\ref{thm three elements}   ---  namely, part~$(\ast)$  ---   using a multivalent map~$\varphi$ in Theorem~\ref{basic main}, and derive as a consequence a first $\Q$-linear independence proof which, unlike with Basic Remark~\ref{log 3 example} to which it is otherwise
entirely similar, is actually a new result. The second half of Theorem~\ref{thm three elements}, which is irrelevant to this application, will be proved in~\S~\ref{sec:completion three elements}. 

A key point to observe is that we shall definitely need a multivalent choice for~$\varphi$.

\begin{basicremark}  \label{bivalent clean}
Koebe's quarter theorem states that $|\varphi'(0)| \leq 4$ for all univalent holomorphic maps of pointed domains $\varphi : (\D, 0) \to (\C \setminus \{1\}, 0)$, 
and that equality holds if and only if $\varphi(z) = G(cz)$ with $|c| = 1$, where 
\[G(z) := \frac{4z}{(1+z)^2} = 1 - \left( \frac{1-z}{1+z}\right)^2\]
 is \emph{Koebe's extremal function}, the Riemann uniformization
map at the origin of the slit complex plane $\C \setminus [1,\infty)$. 
In particular, if we restrict~$\varphi$ to univalent 
maps in
 Theorem~\ref{basic main} then we cannot hope to prove  Theorem~\ref{thm three elements}
 since then 
 \[|\varphi'(0)| \leq 4 < 4.481689\ldots = e^{3/2}.\]
 The Koebe map is $1:1$ on the open unit disc but it extends to a $2:1$ rational map $\C \setminus \{\pm 1\}
\to \C \setminus \{1\}$. Pre-composing this quadratic rational map with the Riemann uniformization map $\D \to \C \setminus \left( (-\infty,-1] \cup [1,\infty) \right)$, which is simply the map $\sqrt{G(z^2)} = 2z/(1+z^2)$, we end up with the bivalent map 
\begin{equation}  \label{Landen second layer}
\varphi : \D \to \C \setminus \{1\}, \qquad  \varphi(z) := G\left(\sqrt{G(z^2)}\right) = \frac{8(z+z^3)}{(1+z)^4} = 1 - \left( \frac{1-z}{1+z} \right)^4.
\end{equation}
In the present section, analogously to the role of Koebe's univalent map for the proof of Theorem~\ref{logcharacterization}, we make a use of the basic properties of 
the bivalent map~\eqref{Landen second layer}. This map bijects $(-1,1)  \iso (-\infty,1)$
and is bivalent on $\D \setminus (-1,1)$, taking either of the two connected halves conformally isomorphically onto $\C \setminus (-\infty,1]$. This 
shows in particular that the case $\varphi : (\D,0) \to (\C \setminus \{1\}, 0)$ in Proposition~\ref{overconvergence} with $\Sigma^1 = \{1,\infty\}$ and an \emph{arbitrary} $\Sigma^0 \subset (-\infty,1)$ can have a derivative as big as $|\varphi'(0)| = 8$. 
 \endofremark
\end{basicremark}

The continuation of this construction explains the central role of the modular lambda map in our paper: 

\begin{remark}  \label{Landen scheme}
We can repeat the process of getting from $G(z) = 4z/(1+z)^2$ to~$G\left(\sqrt{G(z^2)}\right) = 8(z+z^3)/(1+z)^4$ by post-composing next
with the Riemann map of the complement in~$\C$ of the union of the four normal external rays out from the fourth roots of unity $z = \pm i$ (the points that give additional zeros of the map~\eqref{Landen second layer}, that we want to avoid having for $\Sigma^0 = \{0\}$) and $z = \pm 1$ (which give values~$1$ and~$\infty$ for~\eqref{Landen second layer}, which we want to avoid having for $\Sigma^1 = \{1, \infty\}$). But the Riemann map of this $\Z/4$-rotationally symmetrically slit region is just~$\sqrt[4]{G(z^4)}$. The result is the quadrivalent map
\begin{equation}  \label{Landen third layer}
\varphi : \D \to \C \setminus \{1\}, \qquad  \varphi(z) := G\left(\sqrt{G\left(\sqrt{G(z^4)}\right)}\right) = 
\frac{8\sqrt{2} \left( 1+z^2 \right)^2 \sqrt{1+z^4}}{\Big( z\sqrt{2} + \sqrt{1+z^4} \Big)^4},
\end{equation}
which has the bigger derivative $\varphi'(0) = 8\sqrt{2}$ while still serving in Proposition~\ref{overconvergence} for the case
$\Sigma^0 = \{0\}$ and $\Sigma^1 = \{1,\infty\}$. 

Continuing these iterations, we find that the nesting with~$n$ square roots
\begin{equation}  \label{Landen nth layer}
\varphi_n : \D \to \C \setminus \{1\}, \qquad  \varphi_n(z) := G\left(\sqrt{G\left(\sqrt{G\left(  \sqrt{\cdots G(z^{2^{n}}) } \right)}\right)}\right) 
\end{equation}
continues to serve in the $\Sigma^0 = \{0\}, \Sigma^1 = \{1,\infty\}$ case of~\eqref{overconvergence}, while having the derivative
\begin{equation}
\varphi_n'(0) = 4^{ 1 + \frac{1}{2} + \frac{1}{4} +\ldots + \frac{1}{2^n} }    = 16^{1 - 2^{-n-1}}. 
\end{equation}
This constructs a sequence of $2$-solvable algebraic power series $\varphi_0(q), \varphi_1(q), \varphi_2(q),\ldots$ in $\C \llbracket q \rrbracket$ 
starting with the Koebe map $\varphi_0(q) = 4q/(1+q)^2$ and converging coefficients-wise, as well as locally uniformly on $q \in \D$, 
to the modular lambda map~\eqref{lambdadef}. The latter fact was known in essence to Landen, Legendre, and Gauss~\cite[\S~1]{BorweinPiAGM} in
the form of the \emph{arithmetic-geometric mean iteration} $(a,b) \rightsquigarrow \left( \frac{a+b}{2}, \sqrt{ab} \right)$.~\endofremark
\end{remark}

We base our proof of Theorem~\ref{thm three elements} on the bivalent example
$$
\varphi(z) := G\left(\sqrt{G(z^2)}\right) = 8(z+z^3)/(1+z)^4
$$
 from Basic Remark~\ref{bivalent clean}. Crucially, the restriction to~$\D$ of this rational map  has the fairly big derivative $\varphi'(0) = 8$ all the while inducing a bijection $(-1,1) \iso  (-\infty,1)$ and conformal isomorphisms 
 $\D \cap \{ \im(z) > 0 \} \iso  \C \setminus (-\infty,1]$ and $\D \cap \{ \im(z) < 0 \} \iso \C \setminus (-\infty,1]$.

The rationality of this  basic function also allows for an explicit formula of the double integral occurring in the holonomy bound~\eqref{new bound}. 

\begin{lemma}  \label{bivalent BC integral}
The Bost--Charles double integral of the map 
\[\varphi(z) := 8(z+z^3)/(1+z)^4\]
 has the following explicit evaluation: 
\begin{equation}  \label{BC Smyth}
\iint_{\T^2} \log{| \varphi(z)-\varphi(w) |}  \, \mv(z) \mv(w)  = \log{8} + \frac{4G}{\pi},
\end{equation}
where $G := L(2,\chi_{-4})$ is the Catalan constant. 
\end{lemma}

\begin{proof}
This time, to compare with the proof of the univalent case of Theorem~\ref{univalent case}, we have the factorization
\begin{equation}
\frac{\varphi(z)-\varphi(w)}{z-w} = 8 \frac{(1-zw)(1 + i x -iy - xy) (1 - i x + iy - xy) }{(1+z)^4(1+w)^4}, 
\end{equation}
which does have zeros on the unit polydisc $\D^2$: the Bost--Charles \emph{overflow}~\cite[\S~5]{BostCharles} is positive, and it equals
the Mahler measure 
\begin{equation*}
\begin{aligned}
m\left(1+x^2 + y^2 - 4xy + x^2y^2 \right) &  = m(1 + i x -iy - xy) + m(1 - i x + iy - xy) \\
& = 2 m (1 + x +y - xy) = 4G / \pi,
\end{aligned}
\end{equation*} 
where the two integrals make the respective unimodular change of variables $(x,y) \rightsquigarrow (\pm ix, \mp iy)$, and 
the last evaluation is due to Smyth~\cite{BoydMahlerSurvey}. 
\end{proof}

 We will divide the proof of Theorem~\ref{thm three elements} into two parts: property~$(\ast)$ in our statement of the theorem, 
 regarding the meromorphic extendability through~$\delta$ in all analytic continuations; and the derivation of the full form~\eqref{three form eq}
 granting~$(\ast)$. We now prove the first part --- property~$(\ast)$ --- and derive from it a showcase application in~\S~\ref{mixed examples} to $\Q$-linear independence. The second
 part is subtler and will be proved in~\S~\ref{sec:completion three elements} based on the refined holonomy bound Theorem~\ref{main:elementary form}. 

\begin{proof}[Proof of part~$(\ast)$ in Theorem~\ref{thm three elements}] 
 Suppose to the contrary that there exists a $\delta \in (-\infty,1) \setminus \{0\}$ 
such that $f(x)$ converges on $|x| < 1$ while having analytic continuations along all paths in $\PP^1 \setminus \{0,\delta, 1, \infty\}$
and with eventually a  nontrivial local monodromy around $x= \delta$. 
Consider the M\"obius involution $x \mapsto x/(x-1)$ that fixes the origin of the expansions, preserves the $[1,\ldots,n][1,\ldots,n/2]$ denominators type, exchanges the punctures $1 \leftrightarrow \infty$, and maps $\delta \leftrightarrow \delta/(\delta-1)$ to a different puncture (since $\delta \neq 2$), which is also in $(-\infty,1)$.
Then the formal power series $f\big(x/(x-1)\big) \in \Q\llbracket x \rrbracket$ has similar properties to $f(x)$, except now for 
having meromorphic  continuations along all paths in $\PP^1  \setminus \big\{0, \delta/(\delta-1), 1, \infty \big\}$ and with eventually a nontrivial
local monodromy around $x = \delta/(\delta-1)$. As $\delta \neq \delta/(\delta-1)$, it follows at once that the following five functions, all of the $[1,\ldots,n][1,\ldots,n/2]$
denominator type, are $\C(x)$-linearly independent: 
\begin{equation}  \label{li11 case 1}
1, \quad \log(1-x), \quad \log^2(1-x), \quad f(x), \quad f \left( \frac{x}{x-1} \right). 
\end{equation}
We use Theorem~\ref{basic main} with the $5 \times 2$ array 
$$
\mathbf{b} := \left(  \begin{array}{lllll}   0 & 1 & 1 & 1 & 1  \\
0 & 0 & \frac{1}{2} & \frac{1}{2} & \frac{1}{2}   \end{array} \right)^{\mathrm{t}},
$$
corresponding to the denominator types in the ordered list~\eqref{li11 case 1}. We calculate
\begin{equation}
\tau(\mathbf{b}) = \frac{ 1 \cdot 0 + 3 \cdot 1 + (5+7+9) \cdot (3/2) }{5^2} = \frac{69}{50} = 1.38. 
\end{equation}
For the map $\varphi$, we select
$$
\varphi(z) := \frac{8(z+z^3)}{(1+z)^4} = 8 z - 32 z^2 + 88 z^3 - 192 z^4 +\ldots,
$$
whose basic implications we discussed in Basic Remark~\ref{bivalent clean}. This map meets the criteria
in Proposition~\ref{overconvergence} for $\Sigma^1 := \{1,\infty\}$ and $\Sigma^0 := \left\{0,\delta, \delta/(\delta-1) \right\} 
\subset (-\infty,1)$, and with~$f(x)$ replaced by~$Q(x)f(x)$ for a suitable non-zero polynomial~$Q \in \C[x] \setminus \{0\}$
such that~$Q(x) f(x)$ and $Q(x) f(x/(x-1))$ are \emph{holomorphic} (rather than merely meromorphic) 
under analytic continuation along the~$\varphi_*$-images of all paths in~$\D \setminus \{0\}$. 
For the~$\Q(x)$-linear span~$\HH$ of the five functions~\eqref{li11 case 1}, Proposition~\ref{overconvergence}
thus gives~$\varphi^* \HH \subset \mathcal{M}(\D)$, supplying the analyticity hypotheses 
for Theorem~\ref{basic main}. 

 By Lemma~\ref{bivalent BC integral}, the holonomy bound~\eqref{new bound} becomes
$$
5 = m \leq \frac{\log{8} +(4G/\pi)}{\log{8} - 69/50} = 4.640395\ldots, 
$$
a contradiction. 
\end{proof}

\subsubsection{Some mixed periods}  \label{mixed examples} We give an application to irrationality of the theorem we just proved.

\begin{lemma}  \label{around 1/4}
Define
$$
\begin{aligned}
H_A(x) := & \ \frac{1}{\sqrt{1-4x}} \in \Z\llbracket x \rrbracket, \\
H_B(x) := & \ \frac{1}{\sqrt{1-4x}} \int_{0}^{x} \frac{1}{1-t} \frac{1}{\sqrt{1-4t}} \, dt \in \Q \llbracket x \rrbracket, \\
H_C(x) := & \ \frac{1}{\sqrt{1-4x}} \int_{0}^{x} \frac{\log(1-t)}{t\sqrt{1-4t}} \, dt \in   \Q \llbracket x \rrbracket, \\
H_D(x) := & \ \frac{1}{\sqrt{1-4x}} \int_{0}^{x} \frac{\log(1-t)}{1-t} \frac{1}{\sqrt{1-4t}} \, dt \in \Q \llbracket x \rrbracket.
\end{aligned} 
$$
Then~$H_A(x)$, $H_B(x)$, $H_C(x)$, and~$H_D(x)$ have $|x| < 1/4$ for the convergence disc of their Taylor series, and continue as holomorphic
functions along all paths in $\PP^1 \setminus \{0, 1/4, 1, \infty\}$. They have the respective denominator types~$1$, for $H_A$; $[1,2,\ldots,n]$, for $H_B$ and $[1,2,\ldots,n][1,2,\ldots,n/2]$, for $H_C$ and $H_D$. A simple counterclockwise loop encircling the singularity $x=1/4$ induces the following unipotent local monodromy operator: 
\begin{equation}
\begin{aligned}
T(H_A) =  & \  -H_A = H_A - 2H_A, \\
T(H_B) = & \ H_B - 2L(1,\chip) H_A,  \\
T(H_C) = & \ H_C + \frac{\pi^2}{9} H_A, \\
T(H_D) = & \  H_D - 2\left( L(1,\chip) \log{3} - L(2,\chip)  \right) H_A. 
\end{aligned}
\end{equation}
We also have
\begin{equation} \label{L 1 chi}
L(1,\chip) = \frac{\pi}{ 3\sqrt{3}}.
\end{equation}
\end{lemma}
 
\begin{proof}
All are straightforward; we indicate the computation of the $x=1/4$ local monodromy operator~$T$. The first equation, $T(H_A) = -H_A$, is evident as the analytic continuation must be the unique algebraic conjugate. The equation for $T(H_B)$ follows from~\eqref{L 1 chi} and the closed form integration evaluation (and integration by parts 
using $\arctan(1/\sqrt{3}) = \pi/6$) 
\begin{equation} \label{HB t}
H_B(x) = \frac{\pi}{3\sqrt{3}}H_A(x) - \frac{2}{\sqrt{3}} \frac{\arctan{\sqrt{\frac{1-4x}{3}}} }{\sqrt{1-4x}},
\end{equation}
where the second term is regular at $x=1/4$. In general, just like in Basic Remark~\ref{log 3 example}, for any meromorphic function $f(x)$ around~$x=1/4$ we have the meromorphy of $\frac{1}{\sqrt{1-4x}} \int_{1/4}^x \frac{f(t)}{\sqrt{1-4t}} \, dt$ near $x=1/4$ (with both factors switching signs under the monodromy operator~$T$). Consequently, the $x=1/4$ monodromy operator~$T$ acts by
\begin{equation*}
\begin{aligned}
T(H_A)  & = -H_A = H_A - 2H_A, \\
  T(H_B) & = H_B - 2H_A\int_0^{1/4} \frac{1}{1-t} \frac{1}{\sqrt{1-4t}} \, dt, \\
T(H_C) &  = H_C - 2H_A\int_0^{1/4} \frac{\log(1-t)}{t\sqrt{1-4t}} \, dt , \\
 T(H_D) & = H_D - 2H_A\int_0^{1/4}  \frac{\log(1-t)}{1-t} \frac{1}{\sqrt{1-4t}} \, dt. 
 \end{aligned}
\end{equation*}
Fairly straightforward integrations reveal the holonomic coefficients 
$$
\int_0^{1/4} \frac{1}{1-t} \frac{1}{\sqrt{1-4t}} \, dt = \frac{\pi}{3\sqrt{3}} = L(1,\chip),
$$
reaffirming~\eqref{HB t}, and 
$$
\int_0^{1/4} \frac{\log(1-t)}{t\sqrt{1-4t}} \, dt = - \frac{\pi^2}{18}. 
$$
Lastly, an only slightly more involved integration --- or a computing package --- leads to the evaluation of
$$
V(x):=\int \frac{\log(1-x)}{1-x} \frac{1}{\sqrt{1-4x}} \, dx$$
as
\begin{equation*}
\begin{aligned}
& - \frac{2i}{\sqrt{3}} \left(
 \arctan{\sqrt{ \frac{1-4x}{3} }} \left(  \arctan{\sqrt{ \frac{1-4x}{3} }} 
 - i \left( \log{\frac{4(1-x)}{\left(1+\sqrt{ \frac{1-4x}{3} }\right)^2}} \right) \right)  \right) \\
&  +  \frac{2i}{\sqrt{3}}  \left( \Li_2 \left( \frac{\sqrt{ \frac{4x-1}{3} } - 1}{\sqrt{ \frac{4x-1}{3} }+1} \right)  \right),
 \end{aligned}
\end{equation*}
whereupon the familiar formulas 
$$
\begin{aligned}
\arctan{\frac{1}{\sqrt{3}}} &  = \frac{\pi}{6}, \\
L(1,\chip)&  =   \frac{\pi}{3\sqrt{3}}, \\
\Li_2(-1) & =   - \frac{\pi^2}{12}, \\
\Li_2\left(e^{2\pi i /3} \right) & =  -  \frac{\pi^2}{18} + i \frac{\sqrt{3}}{2} L(2,\chip)
\end{aligned}
$$
straightforwardly evaluate the requisite holonomic coefficient 
\begin{equation*}
\int_0^{1/4}  \frac{\log(1-t)}{1-t} \frac{1}{\sqrt{1-4t}} \, dt   = V(1/4) - V(0) =   L(1,\chip) \log{3} - L(2,\chip). 
\end{equation*}
The lemma follows from these integral evaluations. The special value~\eqref{L 1 chi}, which we used in the preceding 
derivation,  is none other than the Dirichlet class number formula for the complex quadratic field~$\Q(\sqrt{-3})$. 
\end{proof}

From part~$(\ast)$ that we already proved in Theorem~\ref{thm three elements} (assuming Theorem~\ref{basic main}, whose treatment is in~\S~\ref{new slopes}), we can thus readily derive a $\Q$-linear independence result out of the circumstance
that the $x=1/4$ local monodromy operator~$T$ simultaneously transforms $H_A, H_B, H_C$, and~$H_D$ by a scalar multiple of
the common function~$H_A$. The linear independence thus sifting through is for the holonomic coefficients in these monodromies: 

\begin{theorem} \label{mixed} The four periods
$$
1, \quad  \frac{ \pi}{\sqrt{3}}, \quad  \pi^2, \quad 3  L(2,\chi_{-3}) -   \frac{\pi}{\sqrt{3}} \log 3
$$ 
are $\Q$-linearly independent. 

In particular, the Mahler measure
\begin{equation} \label{Mahler measure case}
m\left( \frac{(1+x+y)^4}{3} \right) = \int_0^1 \int_0^1 \log{\left|  \frac{(1+e^{2\pi i s}  + e^{2\pi i t})^4}{3} \right|} \, ds \, dt \notin \Q
\end{equation}
is irrational.
\end{theorem}

\begin{proof}
By Lemma~\ref{around 1/4}, 
a $\C$-linear combination 
\begin{equation}
\label{mixedlinear}
f(x) = a H_A(x) + b H_B(x) + c H_C(x) + d H_D(x)
\end{equation}
overconverges at the singularity~$x=1/4$ if and only if
$$
a +  b \frac{\pi}{3\sqrt{3}} - c \frac{\pi^2}{18} +   d \left( \frac{\pi}{3\sqrt{3}}\log{3} -  L(2,\chi_{-3})  \right)  = 0. 
$$
 If this relation 
held with some nonzero integer vector $(a,b,c,d) \in \Z^4 \setminus \{(0,0,0,0)\}$, the combination~\eqref{mixedlinear} would have had
all the requirements  of Theorem~\ref{thm three elements}
  with $\delta := 1/4$. Yet, clearly, $f(x)$ does not vary holonomically on
$\PP^1 \setminus \{0, 1, \infty\}$, only on  $\PP^1 \setminus \{0, \delta, 1, \infty\} = \PP^1 \setminus \{0, 1/4, 1, \infty\}$. 

The irrationality of the Mahler measure~\eqref{Mahler measure case} follows immediately by Smyth's formula~\eqref{smyth}, which
we can rewrite as
\begin{equation*}
\begin{aligned}
m\left( \frac{(1+x+y)^4}{3} \right)  & = 4m(1+x+y) - \log{3} \\ 
& =
  \frac{3\sqrt{3}}{\pi} L(2,\chi_{-3}) - \log 3  =  
 \frac{3 L(2,\chi_{-2}) - \frac{\pi}{\sqrt{3}} \log 3}{ \pi \big/ \sqrt{3}}.
 \end{aligned}
\end{equation*}
This concludes the proof assuming Theorem~\ref{basic main}, which we already proved to imply the requisite part~$(\ast)$ of Theorem~\ref{thm three elements}.

Theorem~\ref{basic main} will be proved in~\S~\ref{new slopes}, and  the full Theorem~\ref{thm three elements} (which we did not need in the preceding argument) will be completed in~\S~\ref{sec:completion three elements}. 
\end{proof}

\subsection{How we prove holonomy bounds}  \label{sec:ideas outline}   We distinguish three principal steps: 
\begin{enumerate}[label=(\roman*)]
\item \label{stepone} Setting up an auxiliary polynomials module $(Q_1,\ldots,Q_m)$, by which we consider auxiliary functions such as $F := \sum_{i=1}^{m} Q_i f_i$ or its multivariable generalizations.
\item \label{dirichletbox} Arranging a Dirichlet box principle or Thue--Siegel lemma for the unknown coefficients of the auxiliary polynomials $Q_i$ 
to have the associated function~$F$ vanish to a high order at~$x=0$. 
\item \label{stepthree} Performing a Diophantine analysis of the lowest order coefficient
of the auxiliary function~$F$. 
\end{enumerate}
In especially favorable circumstances such as Hermite's approximants to the exponential function~\cite{Hermite,Hermite2} or the ensuing approximants
to the logarithm and binomial functions~\cite{ChudnovskyHermite,ChudnovskyThueSiegel}, 
step~\ref{dirichletbox} is replaced by an explicit construction of the requisite polynomials~$Q_i$. Some simplest examples are discussed in~\S~\ref{PerfectPadeExamples}. Such constructions, in the rare occasions that they are possible, usually lead to stronger quantitative results than~\ref{dirichletbox}. For our intricate applications, however, as well
as for the abstract theorems, some form of the Dirichlet box principle is essential. 

The simplest arrangement, which already obtains \emph{some} (rather poor) holonomy bound on the maximal number~$m$
of $\Q(x)$-linearly independent functions, is the following. A commonly used corollary of Siegel's lemma~\cite[Lemma~2.9.1]{BombieriGubler}
states that for a linear homogeneous system of~$M$ linear equations in which the coefficients are rational integers of absolute
values bounded exponentially in a parameter~$\alpha$, while the 
number~$N$ of free variables is no less than twice the number of equations to be solved ($N \geq 2M$), there exist solutions whose components
are rational integers, not all zero, and with absolute values bounded exponentially in~$\max\left(\alpha,\log{N}\right)$. (Cramer's formula constructs
explicitly a nonzero solution of the linear system as soon as the number of free variables strictly exceeds the number 
of equations; but the determinantal expression of this solution gives in general a bound which is exponential in~$M\alpha$ rather
than~$\alpha$; in our setting with $M \asymp \alpha$,  this means  that  the Cramer solution is bounded exponentially in~$\alpha^2$ rather than~$\alpha$.
 For Hermite--Pad\'e approximants to holonomic functions this is not a methodological limitation but actually the correct size in
a majority of naturally occurring cases; cf.~\cite{BombieriCohenPade} for a complete study of the algebraic case.) 

Following Thue~\cite[\S~11]{ThueSelected}, we can improve the upper bound on the solution of the linear system, from exponential in~$\alpha$ to asymptotically subexponential in~$\alpha$,
by using  $N = (1+C)M$ free variables for a large constant~$C$.
 Siegel's lemma 
then supplies nontrivial solutions in rational integers bounded in magnitude by $\exp\big(  O \left( \alpha / C \right) \big)$; this becomes subexponential
in the asymptotic where $C \to \infty$ after $\alpha \to \infty$. 
 Hence, if we have a $\Q(x)$-linearly independent set $f_1, \ldots, f_m$ with denominators of the type $A^{n+1} [1,\ldots, bn]^{\sigma}$
 and with $m$ sufficiently big with regard to $A, b, \sigma$, and the smallest convergence radius of an $f_i(x)$ (this is ultimately handled
 in Lemma~\ref{Siegel}, in a high-dimensional setting that we will need for proving our refined bounds),  Siegel's lemma guarantees
 the existence of a \emph{nonzero} auxiliary function 
 $$
 F(x) := \sum_{i=1}^{m} Q_i(x) f_i (x)  = \beta x^n + O(x^{n+1})   \in \Q \llbracket x \rrbracket, \quad \beta \in \Q^{\times}
 $$
 that vanishes  to some high order~$n$ at $x = 0$, 
  all the while involving integer polynomials $Q_i \in \Z[x]$ whose degrees and coefficients, taken on the logarithmic scale,\footnote{This
  means that all these polynomials have degrees smaller than $cn$ and  rational integer coefficients with absolute values $\ll e^{cn}$.} are
  smaller than an arbitrary desired linear rate~$c n$ in the vanishing order~$n$. 
  
  But the meaning of ``an arbitrary desired linear rate~$cn$''
  is that an arbitrarily small~$c > 0$ is attainable when the number~$m$ of independent functions~$f_i(x)$ is supposed correspondingly large: 
giving a combined number of as many as~$N = mD$ undetermined coefficients for the auxiliary polynomial $m$-tuple $(Q_1, \ldots, Q_m) \in \Z[x]^{\oplus m}_{\deg < D}$. 
  Making this quantitative will ultimately read into holonomy bounds such as~\eqref{new bound}. To explain
  where those derive from, observe that if the functions $f_i(x)$ are of the denominator type $[1,\ldots,b_1n] \cdots [1,\ldots, b_rn]$, 
  then since the auxiliary $Q_i(x)$ have integer coefficients, the lowest order  coefficient $\beta \in \Q^{\times}$
  is some \emph{nonzero} rational number of this denominator, hence
  \begin{equation}  \label{Liouville bit}
    |\beta| \geq \frac{1}{[1, \ldots, b_1 n] \cdots [1, \ldots, b_r n]}. 
  \end{equation}
  By the prime number theorem, this gives a Diophantine lower bound by 
  \[e^{-(b_1 +\ldots + b_r) n + o(n)}\]
   on that
  leading coefficient. Now suppose we have a holomorphic mapping $\varphi : (\Db, 0) \to (\C, 0)$
  of derivative $|\varphi'(0)| > e^{b_1 +\ldots + b_r}$ and turning all $f_i(\varphi(z)) \in \C \llbracket z \rrbracket$
  holomorphic (\emph{convergent}) in a neighborhood of the closed unit disc $z \in \Db$. Then $G(z) := z^{-n}F(\varphi(z))
  = \varphi'(0)^n \beta + O(z)$ is a holomorphic 
  function in a neighborhood of the closed unit disc, but taking an exponentially large value  
  \begin{equation} \label{lower bit}
  |G(0)| = |\varphi'(0)|^n |\beta|   \geq  \exp \left( \left( \log{|\varphi'(0)|} - \sum_{h=1}^{r} b_h  \right)n  + o(n) \right) 
  \end{equation}
  at the center $z = 0$ of that disc. Yet, since by construction the degrees of the polynomials $Q_i(x)$ are smaller than~$cn$
  while their coefficients are smaller than $e^{cn}$, we know in this construction that on the unit circle~$\T$ the 
  holomorphic function~$G(z)$ has the pointwise upper bound
  \begin{equation}  \label{upper bit}
  \sup_{\T} |G| \leq \exp \left(  O \left(  c  \left( \max_{\T} \log{|\varphi|}  \right) \cdot n \right) \right).  
 \end{equation}
 Since we can make the coefficient $c > 0$ arbitrarily small upon assuming~$m$ to be correspondingly big, 
 but the maximum principle for holomorphic functions restrains the left-hand side of~\eqref{lower bit} to be
 not greater than the left-hand side of~\eqref{upper bit}, our assumption of the positive rate in the lower bound~\eqref{lower bit} 
 sets an upper limitation on the maximal number~$m$ of our $\Q(x)$-linearly independent functions $f_i(x)$. 
 This dimension bound only depends on the holomorphic mapping~$\varphi$ and on the positive difference
 $\log{|\varphi'(0)|} - \sum_{h=1}^{r} b_h$ that occurred through~\eqref{lower bit}. We call it an \emph{arithmetic holonomy
 bound} due to the Diophantine way it was proved. 
 
 For simplicity of this sketch, we assumed the $f_i(\varphi(z))$ to be holomorphic rather than meromorphic functions on 
 a neighborhood of the closed unit disc. The general meromorphic case is handled in exactly the same way
 just by changing the definition of the holomorphic function~$G(z)$ to $G(z) := h(z) z^{-n} F(\varphi(z))$, where 
 $h \in \mathcal{O}(\Db)$ is a holomorphic function on a neighborhood of the closed unit disc that 
 has $h(0) = 1$ and all $h(z) f_i(\varphi(z))$ simultaneously holomorphic on that disc. 
 
 In particular, this sketch proves Andr\'e's holonomicity criterion (Corollary~\ref{holonomic criterion}), for 
 by the chain rule, the $\Q(x)$-linear span of all $f(x)$ in Corollary~\ref{holonomic criterion}  
 is closed under the derivation~$d/dx$. This is how holonomy arises out of finiteness theorems. 

\subsection{Refined methods} \label{sec_refined}
This subsection is a deeper and more technical introduction than the rest of \S~\ref{seriousintro}, and it serves as a more detailed summary of the ideas in the proofs of our holonomy bounds. It is not strictly required for the logic of these proofs. The reader might therefore opt to skip any part in the following,   and refer back as needed later. 

The rudimentary proof method we just described in~\S~\ref{sec:ideas outline} is completely 
standard in the subject of Diophantine analysis. It is referred to as Gelfond's method in the works of D\`ebes~\cite{DebesG}
and Andr\'e~\cite[\S~VIII.3]{AndreG}, and found its first applications to arithmetic algebraization in the trailblazing work~\cite{ChudnovskyAlg,ChudnovskyIso} of David 
and Gregory Chudnovsky. Our Appendix~\ref{app:PerelliZannier} refines these ideas through the prism
of Perelli and Zannier's work~\cite{PerelliZannier} to re-derive the bound~\eqref{founding hol} in our context, 
including the $e \rightsquigarrow 2$ coefficient reduction by a single-variable analysis. 
As mentioned in \S~\ref{sec_preview}, for our applications to irrationality, we have two alternative lines of holonomy bounds: one via high-dimensional techniques (Theorem~\ref{main:elementary form}, which implies~\eqref{new bound 2}), and the other via the single variable slopes method (Theorem~\ref{main:BC form}, which implies Theorem~\ref{basic main}, and its strengthening Theorem~\ref{main:BC conv discrete}). We now discuss what can be improved in the preceding scheme to obtain these two lines of refined results separately.
 We begin with the ideas of the proof of Theorem~\ref{main:elementary form}, based on Diophantine approximation
 in several variables. The basic idea can be summarized by saying that our multivariable
 evaluation module will lose none of the simplicity of the essentially one-dimensional features similar to~\S~\ref{app:PerelliZannier}, yet it also has all the added flexibility of the Law of Large Numbers
   inherent in any Diophantine approximation scheme with $d \to \infty$ variables.

\subsubsection{The possible vanishing orders}  \label{filtrations} We can formulate step~\ref{stepone}
 of the preceding scheme differently. 
We do this just as easily in a multivariable framework with $\bx := (x_1,\ldots,x_d)$, which as we will see is ultimately advantageous
for the proofs upon working with the $d$-th Cartesian power of the single-variable evaluation module. 
Given 
\begin{itemize}
\item
a $\Q(\bx)$-linearly independent set  $\left\{ f_{\bi}(\bx) \right\}_{\bi \in I}$ of $\Q\llbracket \bx \rrbracket$ formal power series,
to be indexed by a finite set~$I$ which for our purposes will be taken a subset $I \subset \{1,\ldots,m\}^d$, 
\item and a  bounded Lebesgue-measurable subset $\Omega \subset [0,\infty)^d$,
\end{itemize}
we can express the preceding argument by introducing a parameter~$D$ and taking
$(Q_1, \ldots, Q_m)$, or $(Q_{\bi})_{\bi \in I}$ in this generality, to range from the \emph{auxiliary polynomials module}
$$
E_{D,\Omega}^I := \Cspan_{\Z} \left\{ \bx^{\bk} \, : \, \bk \in (D \cdot \Omega) \cap \Z^d \right\}^{\oplus I},
$$
a free $\Z$-module of rank $R_{D,\Omega}^I = (1+o(1)) (\#I) \, \vol(\Omega) D^d$. For the original case of $d=1$, $I = \{1,\ldots,m\}$,
and $\Omega = [0,1)$, we simplify the notation to $E_D = \Z[x]_{< D}^{\oplus m}$, of rank~$R_{D,[0,1)}^{\{1,\ldots,m\}} = mD$.

The $\Q(\bx)$-linear independence condition on the $f_{\bi}(\bx)$  means 
exactly that, for all~$D$ and $\Omega$, the \emph{evaluation homomorphism}
$$
\psi_D : E_{D,\Omega}^I \hookrightarrow \Q \llbracket \bx \rrbracket, \qquad (Q_{\bi})_{\bi \in I} \mapsto \sum_{\bi \in I} Q_{\bi} f_{\bi} \in \Q \llbracket \bx \rrbracket,
$$
is injective. (We drop $\Omega$ and $I$ from the notation of~$\psi_D$, as they will ultimately be considered fixed throughout the procedure, whereas $D$ will be
the first asymptotic parameter to be let $\to \infty$.) 

In the outline~\S~\ref{sec:ideas outline}, we  considered some power series $F(x) = \beta x^n + O(x^{n+1})$ from the range of this evaluation map (for $d=1$) 
that vanished at $x =0$ to the exact order~$n$. But the possible leading order exponents $\bn \in \NwithzeroB^d$ in any $F = \sum_{\bi} Q_{\bi} f_{\bi} \in E_{D,\Omega}^I$
take up exactly~$R_{D,\Omega}^I = \dim_{\Q}(E_{D,\Omega}^I  \otimes \Q)$ possibilities that depend only on the evaluation module $(E_D, \psi_D)$, and not on 
the specific element $(Q_{\bi}) \in E_{D,\Omega}^I$. These \emph{(vanishing) filtration jumps}\footnote{They may as well be termed the \emph{successive minima} 
of the evaluation module, as in~\cite{BertrandExpAndre} taking an inspiration from the Weierstrass gaps on algebraic curves.} form a size-$R_{D,\Omega}^I$ subset
of $\NwithzeroB^d$, which we formally define in~\S~\ref{vanfiltration}.

For our final results in this paper, we ultimately only consider single variable ODEs. The high-dimensional modules $E_{D,\Omega}^I$ arise from involving
the $d$-fold Cartesian power
\begin{equation} \label{Descartes}
E_{D, [0,1)^d}^{\{1,\ldots,m\}^d} = E_D \times \cdots \times E_D, \quad f_{\bi}(\bx) := \prod_{s=1}^d f_{i_s}(x_s),  \quad \textrm{ of rank } (mD)^d,
\end{equation} 
of the univariate module $E_D = E_{D,[0,1)}^{\{1,\ldots,m\}}$ generated by the functions $f_1,\ldots,f_m \in \Q\llbracket x \rrbracket$
of Theorem~\ref{basic main}, and their suitable submodules --- this is the idea of measure concentration in the $d \to \infty$ limit --- given by restriction to statistically preponderant subsets $\Omega \subset [0,1)^d$ and
$I \subset \{1,\ldots,m\}^d$. 

A basic idea for our new developments here over the results in~\cite[\S~2]{UDC} is a simple lemma (Corollary~\ref{jumps Cartesian}) about 
the commutation in the formations of Cartesian products of evaluation modules and the sets of vanishing filtration jumps. Concretely, the $(mD)^d$ filtration jumps of the evaluation
module~\eqref{Descartes} are at a Cartesian power set of the form $S^d$ for some $S \subset \NwithzeroB$ with $\#S = mD$.

\subsubsection{Methods from differential algebra and functional bad approximability}   \label{holonomic subspace theorem}

In our proofs of the general holonomy bounds, we use a more precise information on the vanishing filtration jumps.  This takes
on the role of the ``zero estimates'' in the traditional transcendence theory proofs. In our context, the latter can be seen as functional analogs
of the Schmidt Subspace theorem on bad approximability. (See, indeed,~\cite{JulieWang} for the case of algebraic functions.)  Easier but cruder versions 
---  analogous rather to Liouville's Diophantine inequality for differential algebra --- include the prototypical Shidlovsky lemma~\cite[\S~3.5, Lemma~8]{Shidlovsky} 
from the historical proof~\cite{Shidlovsky1959} of the Siegel--Shidlovsky theorem on special values of $E$-functions, with its multitude of effectively computable variations~\cite[\S~11]{ChudnovskyShidlovsky}, \cite{BertrandBeukers,BertrandChirskiiYebbou}, \cite[\S~2]{Bertrand} available in the literature. The general bad approximability theorem was known as Kolchin's problem (\cite{Kolchin}, see Problem~\ref{kolchinsproblem}), 
before it was proven independently by David and Gregory Chudnovsky~\cite{Chudnovsky4} and Osgood~\cite{Osgood4},  for the essential case of holonomic $f_1, \ldots, f_m$.
Its statement amounts to saying that the vanishing filtration jumps set~$S$ is close to the generic jumps $\{0,1,\ldots,mD-1\}$, in the sense that 
$$
S \subset \big\{ 0, 1, \ldots,  mD + o_{f_1,\ldots,f_m}(D) \big\},  \qquad \#S = mD.
$$
In an asymptotic sense, this almost determines the vanishing filtration jumps for all the holonomic evaluation modules of relevance to our paper: 
those being the modules $E_{D,\Omega}^I$ with $\Omega \subset [0,1)^d$ of $\vol(\Omega) = 1 - o_{d \to \infty}(1)$; $I \subset \{1,\ldots,m\}^d$
with $\#I = m^d - o_{d \to \infty}(m^d)$, and holonomic $f_1,\ldots,f_m$.  
These improvements are discussed in~\S~\ref{functional bad approximability}.

Technically, for the qualitative linear independence proofs of Theorem~\ref{mainA} and Theorem~\ref{logsmain} (up to replacing the numerical threshold~$10^{-6}$ by a smaller absolute constant), it is actually possible to avoid all recourse to this differential algebra material~\S~\ref{functional bad approximability}. 
 It is however an unnecessarily convoluted route to insist on; moreover, some version of the theorems collected in~\S~\ref{functional bad approximability} is indispensable in pursuing any quantitative refinements to Diophantine measures of linear independence. 
 We choose to use functional bad approximability in our main proofs as well, for the holonomy bounds in \S\S~\ref{fine section}--\ref{slopes}, as that allows
 for cleaner arguments, and is actually (as far as we are aware) necessary for most of the general --- qualitative! --- holonomy bounds in the clean structural form in which we have stated them. 
 We do observe, however, that such structural necessities do not concern the main~$|\varphi'(0)| > e^{\sigma_m}$ 
 case of Theorem~\ref{basic main} itself, which does admit clean proofs not relying on any functional bad approximability theorems. (Remark~\ref{converse remark} shows that the~$|\varphi'(0)| > e^{\tau(\bb)}$ case \emph{must} make some special use of the ODE.) All this is discussed in~\S~\ref{sec_withoutShid}. 
  The reader may compare the situation with the simpler \S~\ref{app:PerelliZannier}, where no special information on the vanishing filtration jumps is relevant to the proof of the qualitative holonomy bound~\eqref{basic basic}.

\subsubsection{Multiple variables unlock the Law of Large Numbers} \label{high dims}  We next discuss how, in the $d \to \infty$ asymptotic modeled
by independent
and identically distributed random variables, we can 
exploit the full-measure subsets $\Omega \subset [0,1)^d$ and $I \subset \{1,\ldots,m\}^d$. Historically, Diophantine approximation by multiple variables
was the key to refining Liouville's bad approximability theorem $|\alpha - p/q| \gg q^{-[\Q(\alpha):\Q]}$ 
to Roth's ``best-possible''  bad approximability measure $|\alpha - p/q| \gg_{\varepsilon} q^{-2 - \varepsilon}$ 
(when the target~$\alpha$ is algebraic and irrational). The purpose of the scheme\footnote{Found by Siegel, and attempted with partial success
by Schneider~\cite{Schneider} prior to Roth's work~\cite{Roth}.} is to make the maximum use of the free parameters count in the application of
Siegel's lemma. Having a multivariable auxiliary function $F(x_1, \ldots, x_d)$ vanish to a high $(D_1, \ldots, D_d)$-weighted order $\geq \xi d$
at a point $(0,\ldots,0)$ means to vanish all monomials $\mathbf{x^n} := x_1^{n_1} \cdots x_d^{n_d}$ with $n_1/D_1 +\ldots + n_d/D_d \geq \xi d$.
 But as $d \to \infty$ and $D_i \to \infty$ with $t_i := n_i/D_i \in [0,1]$, the Law of Large Numbers (in Chernoff's form) for  the sum $\sum_{i=1}^{d} t_i \approx d/2$ of $d \to \infty$ uniform and identically distributed
 random variables $t_i \in [0,1]$ shows that, with an $1- \exp(O(-d \varepsilon^2 ))$ probability, $\xi =  \frac{1}{2} - \varepsilon$ is the correct reasonable weighted vanishing order 
 to attain by the parameter count in the Thue--Siegel lemma. (To contrast, the single variable construction only reaches the Liouville-strength 
vanishing order coefficient $\xi = 1/[\Q(\alpha):\Q]$, and the two-variables construction only reaches a vanishing order coefficient of about~$\xi = 1/(2\sqrt{[\Q(\alpha):\Q]})$, giving 
the exponent in Siegel's sub-Liouville theorem~\cite{SiegelZeitschrift}.) 

Further work of Wirsing~\cite[see \S~4.2]{Wirsing}, aimed at correcting Roth's Corrigendum in~\cite{Roth} regarding approximation 
of an algebraic number target by algebraic number approximants of a fixed degree over~$\Q$, pivoted around a refinement of the
above Law of Large numbers, the {\it measure concentration property of the high-dimensional hypercube $[0,1]^d$}, which states
that not only $\frac{1}{d} \sum_{i=1}^{d} t_i$ converges in probability to the expectation $\mathbf{E}[t] = \int_0^1 t \, dt = \frac{1}{2}$ as $d \to \infty$, but further and more precisely, that
with high asymptotic probability as $d \to \infty$, the random vector $(t_1, \ldots, t_d) \in [0,1]^d$ has uniformly distributed components.
This has a precise meaning in our Theorem~\ref{thm_MeasureConcentration} below refining~\cite[Lemma~13]{Wirsing}: the $\varepsilon$-high discrepancy set (see Definition~\ref{box discrepancy})
$$
B_{\varepsilon}^d :=  \left\{  \mathbf{t} \in [0,1)^d \, : \, \exists [a,b) \subset [0,1), \,  
\Big| (b-a)- \frac{1}{d}  \# \{ i \, : \, t_i \in I \} \Big| \geq \varepsilon  \right\}
$$
has $d$-dimensional Lebesgue measure $\vol(B_{\varepsilon}^d) \leq 100\exp\left( -\varepsilon^4 d/300 \right)$. 
 This is ultimately the statistical property behind 
the rearrangement integral in our bound
\begin{equation}\label{intro_elementary bound}
m \leq \frac{\int_0^1 2 t \cdot (\log{|\varphi(e^{2\pi i t})|})^* \, dt}{  \log{|\varphi'(0)|} - \tau(\mathbf{b})}
\end{equation}
discussed in \S~\ref{sec:variation}; this bound is a special case of Theorem~\ref{main:elementary form} with $\be=\mathbf{0}, l=0, \varphi_0=\varphi$.

 We review in~\S~\ref{concentration of measure} the topic of measure concentration and large statistical deviations. These ideas are used not only to sift through the subsets $\Omega \subset [0,1)^d$ and $I \subset \{1,\ldots,m\}^d$ in the make-up of the evaluation module $E_{D,\Omega}^I$, but also to usefully limit the shape of the monomials~$\bx^{\bn}$ from the leading order jet of the $d$-variate auxiliary function $F \in \sum_{\bi} Q_{\bi} f_{\bi} \in \Q \llbracket \bx \rrbracket \setminus \{0\}$. The former 
 leads to the rearrangement integral; the latter two lead, in particular, to the refined denominators rate~$\tau(\bb)$. 
 We discuss in \S~\ref{oo parameter count} the mechanism for both these improvements. As we have mentioned in \S~\ref{sec_preview}, for the case of  basic denominator types as in Theorem~\ref{basic main} (as well as in all the other holonomy bounds in \S\S~\ref{fine section}--\ref{new slopes}), the exact same denominator saving comes through also by the single variable method (hence no measure concentration) of \S~\ref{new slopes}; while the best general denominator term comes through in Theorem~\ref{high dim BC convexity}, again by measure concentration. Although the proofs of Theorems~\ref{main:elementary form} and~\ref{high dim BC convexity} are described in different languages (one via the Thue--Siegel Lemma, the other via Bost's slopes method), the ideas on treating denominators behind both proofs are the same, as is the scope for the further improvements in the denominators aspect. 

In all three proofs of our holonomy bounds in~\cite[\S~2]{UDC}, we used $d \to \infty$ for its automatic improvement of the Dirichlet exponent --- namely, if~$M$ is  the number of equations and~$N$ is the number of parameters,
then~$M/N = o_{d \to \infty}(1)$, and hence the Dirichlet exponent~$M/(N-M)$ is also~$o_{d \to \infty}(1)$.
 In this paper, this aspect is shown again in~\eqref{small Dirichlet exponent}, but for this particular point, $d \to \infty$ is used only
as a methodological feature of working with the most traditional form of the Thue--Siegel lemma. (Appendix~\ref{app:PerelliZannier} explains 
how we could bypass the auxiliary coefficients size while sticking to the single variable module~$E_D$, in a form similar to the slopes method
treatment in~\S~\ref{new slopes}.) The input from measure concentration is by far the more
essential use of the high dimensions. 

The fine improvements in the numerator and denominator of the fraction~\eqref{founding hol} are however only relevant insofar
that they also come with an $e \rightsquigarrow 2$ overall coefficient reduction. We discuss next how this is achieved by exploiting, in the Dirichlet box principle, Lemma~\ref{jumps Cartesian} on the Cartesian power structure of the vanishing filtration jumps, in the sense described in~\S~\ref{filtrations}. This point will also clarify the employment of the functional bad approximability results that we mentioned in \S~\ref{holonomic subspace theorem}.

\subsubsection{The high-dimensional parameter count} \label{oo parameter count} 
At the outset, to have an auxiliary function $F := \sum Q_{\mathbf{i}} f_{\mathbf{i}}$ in 
the range of the general evaluation module $(E_{\Omega,D}^I,\psi_D)$ to  vanish to an order at least~$\alpha$ at $\mathbf{x=0}$, involves solving
$\binom{ \alpha  + d }{d} \sim \alpha^d / d!$ linear equations in the $\sim \mathrm{vol}(\Omega) (mD)^d$ unknown coefficients of the 
polynomials $Q_{\mathbf{i}}$. By Stirling's asymptotic $d! = d^d/e^{d-o(d)}$, the maximal attainable vanishing order in the high-dimensional asymptotic
$d \to \infty$ appears to be $\alpha \sim mdD/e$. This 
was why in~\cite[\S~2]{UDC} we have the coefficient~$e$ in the holonomy bound that we 
established there with the hypercube choice $\Omega := [0,1)^d$. Had we used instead the simplex choice 
$$
\Omega := \left\{ (t_1, \ldots, t_d) \in [0,1] \, : \, t_1 +\ldots + t_d < 1   \right\},  \qquad \mathrm{vol}(\Omega) = 1/d!, 
$$
we would have entertained an asymptotic vanishing order as high as~$\alpha \sim mD$ (without the number~$e$ entering in as a coefficient); but in this case 
the functions $Q_{\mathbf{i}}(\varphi(z_1),\ldots, \varphi(z_d))$ would be far too big on the unit polycircle $\mathbf{z} = (z_1, \ldots, z_d) \in \T^d$. In such an approach,
 we would have only obtained an inadequately big holonomy bound with  (in the context) an exponentially larger numerator such as $\sup_{\T} \log{|\varphi|}$, instead of the 
  Nevanlinna growth characteristic $T(\varphi) = \int_{\T} \log^+{|\varphi|} \, \mv$. In the present paper, for a similar reason,
we still use the hypercube shape $\Omega = [0,1)^d$, or more precisely, its measure-concentrated subsets $\Omega = P_{\epsilon}^d :=  [0,1)^d \setminus B_{\epsilon}^d$  for the auxiliary monomials exponents range. 
As we discussed above, we do rely on these statistically preponderant parts of the high-dimensional hypercube in order to get the
refined growth integral~\eqref{maxintegral}; moreover, we will explain how we use these statistics to control the shape of lowest order terms in Siegel Lemma construction in order to obtain the denominators counterpart~\eqref{finitary rearrangement} of the refined growth integral. 

The Cartesian power situation is special for enforcing, as discussed in~\S~\ref{filtrations}, a Cartesian power structure (Corollary~\ref{jumps Cartesian}) on the vanishing filtration jumps vectors $\subset \NwithzeroB^d$ of $E_{[0,1)^d,D}$. 
These are in turn brought to exploit a certain automatic vanishing of many of the coefficients of the sought-for auxiliary function~$F$.
The simplest instance of this automatic vanishing is showcased in~\S~\ref{dynamic box}. Instead of directly solving for the vanishing of all the low-degree monomials of~$F$
(which we definitely need for the maximum modulus principle step when we carry out the higher-dimensional extension of step~\ref{stepthree} of~\S~\ref{sec:ideas outline}), we set up the Thue--Siegel lemma differently by focusing on the $mD$ filtration jumps
$$
0 \leq u(1) < u(2) < \cdots < u(mD)
$$
 of the single-variable evaluation module $E_D$. In the single-variable situation the procedure simply reduces to setting to zero the $x^{u(p)}$ coefficient $\beta_{u(p)} = 0$  of~$F(x)$ for $p = 1, \ldots, mD$. 
 In general, for any subset $T \subset [0,mD]^d$, we write
$$
u(T) := \big\{ (u(s_1), \ldots, u(s_d)) \, : \, (s_1, \ldots, s_d) \in T \cap \NwithoutzeroA^d \big\}. 
$$
 Then, in the
measure-concentrated submodule $E_{D,\Omega}^I \subset E_{D,[0,1)^d}^{\{1,\ldots,mD\}^d}$ with $\vol(\Omega) > 1 -  100\exp\left(-\epsilon^4 d/300\right)$, we use our $(mD)^{d-o(d)}$ degrees of freedom in the auxiliary polynomials coefficients to construct a nonzero $F(\mathbf{x}) = \sum_{\bi} Q_{\bi}(\bx) f_{\bi}(\bx) = \sum \beta_{\mathbf{n}} \mathbf{x^n}$ with all auxiliary polynomials~$Q_{\bi}$ having integer coefficients bounded by~$e^{\epsilon dD}$ in absolute value, and in which (with a sufficiently small~$\delta \in (0,\epsilon)$)
\begin{equation} \label{from 1 gen}
\beta_{\bn} = 0 \quad \textrm{for all} \quad
\bn \in   
u \big( [0,(m-\delta)D]^d  \big)
\cup u\big( (m+\delta)D\cdot B_{\epsilon}^d \big),
\end{equation} 
provided~$\delta \in (0,\epsilon)$ is small enough to have~$(m-\delta)/(m+\delta) >   \exp\left(-\epsilon^4/400\right)$.
For Theorem~\ref{basic main} when~$|\varphi'(0)| > e^{\sigma_m}$, and  for some further forms of our bounds that are discussed in~\S~\ref{sec_withoutShid} (which do cover, in particular, the ultimate application to Theorem~\ref{mainA}), it is technically possible to devise a proof directly out of this construction, and without appeal to the ideas of~\S~\ref{holonomic subspace theorem}. 

In any event, for our practical purposes in this paper, if the reader would like to further simplify the essential mental picture,
 it would be very reasonable to imagine at this point that the filtration jumps
 are as simple as possible, namely given by~$u(i) = i - 1$. Such is for example the case with the classical Hermite--Pad\'e systems that we discuss in~\S~\ref{PerfectPadeExamples}.
 The tenor of the functional bad approximability theorems of \S~\ref{functional bad approximability} is that, for the purposes of many applications including ours, such an assumption 
 \silentcomment{(while being, in a certain sense~\cite{BombieriCohenPade}, actually typically false)}
  is not far from being satisfied: the Chudnovsky--Osgood theorem, as we formulated in~\S~\ref{holonomic subspace theorem}, 
  can be stated as the upper bound $u(mD) \leq (m+\varepsilon)D + C(\varepsilon)$, which is at most $(m+\epsilon)D$ for $D\gg 1$ if we assume $\varepsilon < \delta < \epsilon$. We observe as a statistical effect that the difference between $u(i)$ and its lower bound~$i-1$ becomes negligible in the asymptotic analysis
  of $D \to \infty$ followed by~$\epsilon \to 0$. 
 Our gateway to the functional bad approximability theorems is through Andr\'e's holonomicity criterion (Corollary~\ref{holonomic criterion}; unless, as in the applications, the $f_i$ are \emph{a priori} given holonomic), whose proof was outlined in~\S~\ref{sec:ideas outline} and laid out in full in~\S~\ref{app:PerelliZannier}. 
 With the Chudnovsky--Osgood theorem, the previous construction reduces simply to attaining
\begin{equation} \label{from 2 Chudnovsky}
\beta_{\bn} = 0 \quad \textrm{for all} \quad
\bn \in [0,(m-\delta)D]^d  \cup (m+\delta)D \cdot B_{\epsilon}^d.
\end{equation}

Whichever the approach, the routine for step~\ref{stepthree} of~\ref{sec:ideas outline} is to
 examine the possibilities for a lowest-order nonzero coefficient 
 $\beta := \beta_{\bn} \neq 0$. 
On the one hand, as we discussed above, the Cartesian structure 
restrains $\bn$ to be of the form
$$
\bn = \left( u(p_1), u(p_2), \ldots, u(p_d) \right), \qquad \textrm{ for some } p_1, \ldots, p_d \in \{1,\ldots, mD\}. 
$$
Since our 
Thue--Siegel lemma construction disposed of all the multi-indices $(p_1,\ldots,p_d)$ in $(mD+\delta) \cdot B_{\epsilon}^d$ lying in the $\epsilon$-high
discrepancy part of the hypercube, the above tuple $(p_1, \ldots, p_d)$ must belong to the complementary part $(mD+\delta) \cdot P_{\epsilon}^d$ of the statistically
typical points. 
Heuristically speaking, the components $n_j = u(p_j)$ of each lowest-order exponent vector $\mathbf{n}$ in the Taylor series of $F \in \Q \llbracket \bx \rrbracket
\setminus \{0\}$ are close --- as $d \to \infty$ followed by $\epsilon \to 0$ --- to some ordering of the set 
$$
\left\{ u(\lfloor mD/d \rfloor), u(\lfloor 2mD/d \rfloor), 
\ldots, u(\lfloor dmD/d \rfloor)    \right\}.
$$
 In particular, the vanishing order in this auxiliary construction satisfies
\begin{equation*}
\begin{aligned}
\mathrm{ord}_{\mathbf{x=0}} F = |\bn|  &  = \left( 1+ o(1) \right) \sum_{j=1}^{d} u(\lfloor jmD/d \rfloor) \\
& \geq  \left( 1+ o(1) \right) \sum_{j=1}^{d}  jmD/d  = (1+o(1)) mdD/2, 
\end{aligned}
\end{equation*}
a notable improvement of the asymptotic vanishing order parameter $\alpha \sim mdD/e$ in~\cite[\S~2]{UDC}. 

This heuristic lower estimate does indeed match the accurate asymptotic formula from using the Chudnovsky--Osgood theorem and~\eqref{from 2 Chudnovsky}. 
The one (fundamentally minor) technical point in arguing directly from~\eqref{from 1 gen}, for the reader who may desire additionally here to forsake the theorems in \S~\ref{functional bad approximability},  is that --- for the discrepancy theory purposes of our proofs --- the uniform distribution $\{p_1,\ldots,p_d\} \approx \left\{\lfloor mD/d \rfloor, \lfloor 2mD/d \rfloor, \ldots, \lfloor dmD/d \rfloor \right\} $ does not preserve the~$\approx$ relation upon applying~$u$ to both sides. 
This is however irrelevant to the above outline; all that matters is that the facts that $u(i)\geq i-1$ and that $(p_1, \ldots, p_d)$ has  asymptotically uniformly distributed components by themselves
imply  $|\bn| \geq (1+o(1)) mdD/2$.

\subsubsection{Effects on denominators}
We now discuss how to obtain the refined denominators saving in Theorem~\ref{main:elementary form}, and by extension,
in Theorem~\ref{high dim BC convexity}. To illustrate the idea, we use the construction~\eqref{from 2 Chudnovsky} contingent upon the Chudnovsky--Osgood theorem, and we consider the lexicographically minimal term $\beta \, \bx^\bn$ in $F(\bx)$ among all the terms of the minimal total degree $n = |\bn|$. 
Thus $\bn \in (m+\delta)D \cdot P_{\epsilon}^d$, and $F$ is a $\Q[\bx]$-linear combination of $f_{\bi}$, where every $\bi \in I$ is \emph{balanced}, namely each $ i_0 \in \{1,\ldots,m\}$ occurs about $ d/m $ times among all $i_j, 1\leq j \leq d$. The denominator of the nonzero rational number $\beta \in \Q^{\times}$ divides the lowest common multiple of the denominators of the $\bx^\bn$ terms in all formal functions from the modules $\Q[\bx] f_\bi$, as $\bi \in I$ ranges over the balanced multi-index sets. The particular form of the denominators assumed in Theorem~\ref{basic main} --- with the types of $f_1, \ldots, f_m$ being ``from best to worst'' in this order --- implies that said lowest common multiple formally agrees asymptotically with the $\bx^\bn$ coefficient of $f_{\bi_0}$, where $\bi_0$ is a balanced multi-index arranged in nondecreasing order. This observation yields our denominator saving term $\tau(\bb)$ as a ``finite	 rearrangement integral''~\eqref{finitary rearrangement}. In general, there is not a single particular~$\bi$ to make the asymptotic denominator of $\beta$; this ``collective~$\bi_0$'' is rather the formal effect of working only with the balanced~$\bi$, which --- as another effect of the measure concentration\footnote{Here basically amounting to the maximality of the central multinomial coefficients, cf. Lemma~\ref{balanced indices}.} advantage of $d\to \infty$ --- are statistically preponderant in $\{1,\ldots,m\}^d$.

\subsubsection{Complex-analytic tools} \label{high dimensional analysis}  The maximum principle can be replaced 
by the Poisson--Jensen formula (\S~\ref{sec:slopes archimedean} or \cite[\S~2.4]{UDC}) or enhanced by seeking the
optimal quotient representation $\varphi = v/u$ by holomorphic functions $v,u \in \mathcal{O}(\Db)$ with $u(0) =1$
(\S~\ref{iii coeff} and~\cite[\S~2.3]{UDC}). But in the $d \to \infty$ asymptotic we discussed in~\S~\ref{oo parameter count}, a better holomorphic
dampener than $u(z_1)^D \cdots u(z_d)^D$ to use in the multivariable analytic function $F(\varphi(z_1), 
\ldots, \varphi(z_d))$ would be to take a suitable power of the discriminant polynomial $\prod_{1 \leq i < j \leq d} (z_i - z_j)$, which is very small
on the part of the torus $\mathbf{z} = (z_1,\ldots,z_d) \in \T^d$ where the set $\{z_1, \ldots, z_d\}$ has a non-small discrepancy from 
the uniform measure~$\mv$ of the circle~$\T$. This was the \emph{ad hoc} method in~\cite[\S~2.5]{UDC}, which fits here the most naturally into the
cross-variables integration scheme stemming from~\S~\ref{oo parameter count}, ultimately leading into the bound~\eqref{intro_elementary bound}.  This is our treatment in~\S~\ref{equidistribution}.
 In~\S~\ref{sec:slopes archimedean}, we give a second treatment based on the Poisson--Jensen formula. 
 
 The further refinements that we mentioned in~\S~\ref{sec_preview} are based on the following idea. If we consider another holomorphic mapping $\psi: (\Db,0) \rightarrow (\C,0)$ also having all $\psi^* f_i \in \mathcal{M}(\D)$, we may replace a subset of the $\varphi(z_j)$ in $F(\varphi(z_1), \ldots, \varphi(z_d))$ by $\psi(z_j)$, and carry out a similar analysis thus using the combined analytic maps~$\varphi$ and~$\psi$.  To use $\varphi$ for $j \in S_1$ and $\psi$ for $j \in S_2$, for some partitioning $\{1, \ldots, d\} = S_1 \sqcup S_2$ of the indexing set into proportionally large subsets $S_1$ and $S_2$, observe that upon taking our holomorphic dampener to be a suitable power of 
 $$
 \prod_{1 \leq i < j \leq d, i,j \in S_1} (z_i - z_j)\prod_{1 \leq i < j \leq d, i,j \in S_2} (z_i - z_j),
 $$
  the main contribution to the growth of the auxiliary function pullback on~$\T^d$ comes from the part of the torus $\mathbf{z} \in \T^d$ where \emph{both} sequences $(z_j)_{j\in S_1}$ and $(z_j)_{j\in S_2}$ have a small discrepancy from 
the uniform distribution on the circle~$\T$. The point is that, when we estimate the leading-order  $\bx^\bn$ coefficient~$\beta$ by the analytic method, only the variables indexed by $j \in S_1$ use~$\varphi$ while the variables indexed by $j \in S_2$ use~$\psi$. 
We select the partition so as to minimize the upper bound from maximum principle over $\bz \in \T^d$. For a given $\varphi$, we may certainly take our second (or, repeating the procedure, our ``next'') map to be $\psi (z) := \varphi(rz)$ for an arbitrary $0<r<1$. As far as $\varphi$ is not univalent, we prove that depending on the size of $n_j$, one may choose for the variable~$z_j$ a certain optimal radius $r = r(n_j) \in (0,1)$ to obtain a strictly better estimate. This is the idea behind the improvement from the bound~\eqref{intro_elementary bound} to the full Theorem~\ref{main:elementary form}.

\subsubsection{A dynamic box principle or a finer Geometry of Numbers} The basic sketch given in~\S~\ref{sec:ideas outline} was grounded
in a ``static''  Thue--Siegel lemma construction: finding a nonzero auxiliary function $F \in E_D$, then arguing ``by extrapolation'' from putting
together the arithmetic and the analytic properties of the lowest-order nonzero coefficient $\beta \in \Q^{\times}$. This simple-minded procedure is insufficient
for obtaining the original holonomy bound~\eqref{original UDC bound} by a single-variable analysis, because in
the Thue--Siegel lemma of~\S~\ref{sec:ideas outline} it is impossible to attain a small Dirichlet exponent $M/(N-M) < c$ all the while 
having a near-maximal
vanishing order~$M \approx N$. In~\cite[\S~2]{UDC} we exploited the decaying Dirichlet exponent under $d \to \infty$ (in the present
paper, this is the step~\eqref{small Dirichlet exponent}), making the issue go away in the high-dimensional analysis. 
As the Bost--Charles work~\cite{BostCharles} made it abundantly clear, it is possible to prove~\eqref{founding hol}, even with the 
coefficient reduction $e \rightsquigarrow 2$, by one-dimensional methods once the rudimentary Thue--Siegel lemma is replaced by a sufficiently precise
arrangement of the pigeonhole or Minkowski arguments. Our Appendix~\S~\ref{app:PerelliZannier}  gives an essentially elementary such 
proof based on the dynamic box principle technique of Perelli and Zannier~\cite[Lemma~1]{PerelliZannier}. This may be also read as 
an introduction to Bost's slopes method framework, whose idea is very similar but cast into the language of Hermitian vector bundles
over $\spec{\Z}$, and which is the content of~\S~\ref{new slopes}. \endofremark \\

We now discuss the more specific ideas for the proofs of Theorems~\ref{main:BC form},~\ref{main:BC conv discrete}, and~\ref{main: easy convexity} via Bost's method of slopes.
Common ingredients (with the simplification applied to a single variable situation) are \S\S~\ref{holonomic subspace theorem} and~\ref{high dimensional analysis}.

\subsubsection{Bost's slopes method with ingredients from Bost--Charles \cite{BostCharles}}  \label{BC summary}
We adapt the notation from \S~\ref{filtrations} to consider a filtered $\Z$-module $E_D$ and an evaluation homomorphism~$\psi_D$. (In the bulk of~\S~\ref{new slopes}, we opt to rather work with $x^{1-D} E_D$ as that allows for a more natural identification with the global sections of a certain ample line bundle; for simplicity here, we stick to the positive degree monomials, like we do in one of our more elementary slopes method variations in~\S~\ref{sec: proof easy convexity}.) We let $E_D^{(n)} \subset E_D$  to denote the $n^{\mathrm{th}}$ vanishing order filtration, namely the submodule consisting of  those elements whose image under $\psi_D$ vanishes to order at least $n$ at $x=0$. The evaluation homomorphism $\psi_D : E_D \hookrightarrow \Q \llbracket x \rrbracket$ then induces a set of \emph{monomorphisms} $\psi_D^{(n)}$ on the graded quotients. Once one endows $E_D$ with a Euclidean lattice structure, one can define an \emph{arithmetic degree} of the underlying Hermitian vector bundle~$\overline{E}_D$, and the \emph{heights} of the evaluation maps $\psi_D^{(n)}$. (Doing this involves fixing a lattice of $\Q\llbracket x \rrbracket$ and endowing it with a pro-Euclidean structure. This then defines the local and global heights of $\psi_D^{(n)}$ following~\cite[\S~1.4.3]{BostBook}.  We stick to the natural lattice choice, namely $\Z \llbracket x \rrbracket$ with pro-Euclidean 
structure induced from using  $\{x^n\}_{n=0}^{N-1}$ for an orthonormal basis of each finite-dimensional quotient
$\R\llbracket x \rrbracket \big/ x^N \R\llbracket x \rrbracket$ of~$\R\llbracket x \rrbracket$.)  

Bost's slopes inequality~\eqref{slope-inequality} provides an upper bound on the arithmetic degree of $\overline{E}_D$ in terms of the heights of the evaluation maps $\psi_D^{(n)}$. In \cite{BostFoliations, BostGerms}, Bost proved various algebraicity criteria in arithmetic-geometric settings similar to~\S~\ref{sec:algebraic}. His methods combined a crude version of the global arithmetic Hilbert--Samuel formula, used as a lower bound on the arithmetic degree of~$\overline{E}_D$, 
with local complex and $p$-adic analysis tools, employed to devise place-by-place upper bounds on all the local heights of all the evaluation maps $\psi_D^{(n)}$. The algebraization results then sift out from the slopes inequality under the $D\rightarrow \infty$ asymptotic. 
The recent work of Bost and Charles~\cite{BostCharles} is written (in part) under the framework of Bost's theory~\cite{BostBook} of theta invariants of infinite-dimensional Hermitian vector bundles over arithmetic curves, but one can certainly interpret the arguments in the language of the more rudimentary slopes method. We stick to the latter choice because the convexity enhancements in~\S~\ref{sec:convex} seem to be more of an analytic than a geometric nature, and we do not attempt here to include these into the theory of the theta invariant. 

The main ingredients of the proofs of the bounds~\eqref{BC integral} and~\eqref{BC integral 2} are the arithmetic Hilbert--Samuel formula for the exact
asymptotic arithmetic degree of $\overline{E}_D$, and a choice of the Euclidean structure giving rise to $\overline{E}_D$ based upon optimizing the complex analysis of the archimedean local heights of $\psi_D^{(n)}$. The latter relies on the standard tools of the subject: the Poincar\'e--Lelong and Poisson--Jensen formulas. One technical point in Bost and Charles's theory~\cite[\S~4]{BostCharles}, needed for carrying out the arithmetic intersection number computations, is to extend the scope of the classical Arakelov theory to allow for Green functions and Hermitian metrics that are not necessarily smooth but have, in Bost and Charles's terminology, a $\mathcal{C}^{\mathrm{b} \Delta}$ regularity: a condition~\cite[Def.~4.1.1]{BostCharles} related to using continuous Green functions locally of bounded variation. We use this framework in \S\S~\ref{new slopes},~\ref{slopes}.

\subsubsection{Varia} \label{varia}
To prove Theorem~\ref{main:BC form}, we use the same Euclidean norm on $E_D\otimes_{\Z} \R$ as alluded to in the final paragraph of~\S~\ref{BC summary}; and we adapt the same complex analytic estimates on the archimedean local heights of $\psi_D^{(n)}$. On the other hand, based on the denominator type~\eqref{den type integrated}, we choose a new $\Z$-sublattice of $E_D\otimes_{\Z} \Q$ to optimize the comparison between the arithmetic degree of~$\overline{E}_D$ and the combined non-archimedean heights of all the evaluation maps~$\psi_D^{(n)}$ involved in the slopes inequality~\eqref{slope-inequality}. The latter brings out the vanishing filtration jumps sets~\S~\ref{filtrations} of the evaluation module, and this is where the theorems from \S~\ref{functional bad approximability} (as summarized by~\S~\ref{holonomic subspace theorem})
are used in this method also. 

The extra input for Theorem~\ref{main:BC conv discrete} rests on the idea of~\S~\ref{high dimensional analysis} with the multiple holomorphic maps~$\varphi, \psi, \ldots$ for devising sharper estimates on the various archimedean local evaluation heights~$h_{\infty}(\psi_D^{(n)})$ 
in accordance with the range of~$n/D$. We stick to the choice $\psi(z) = \varphi(rz)$ derived from a single holomorphic ambience~$\varphi : \Db \to \C$, where the convexity property in~$\log{r}$ for various Nevanlinna-style growth characteristic functions implies that, in the multivalent case, there is always some improvement from every single intermediate radius~$r \in (0,1]$. Ultimately this leads to the limiting form in Theorem~\ref{main:BC convexity}, where the total convexity saving is presented as a~$dr/r$ integral over~$r \in [0,1]$ of the square of an analog of the Ahlfors--Shimizu covering area function. 
Such a principle applies to a number of variations of the Nevanlinna characteristic of a meromorphic map that could be used for the principal term of the holonomy bounds, including the traditional Ahlfors--Shimizu characteristic figuring in~\cite[Prop.~5.4.5]{BostCharles}; more significantly for us (see Ex.~\ref{examplecompare} for an illustrative comparison), it holds for the Bost--Charles characteristic that we define in~\S~\ref{BC characteristic} (sticking for simplicity to the most basic case that we use of a holomorphic mapping~$\Db \to \C$).

We remark on the other hand that once we involve this new improvement on the archimedean local evaluation height bounds from using a finite set of intermediate radii~$r$, Bost and Charles's choice of Euclidean norm on $E_D\otimes_{\Z} \R$ is no longer the optimal in general (unless $\varphi$ is univalent).
 We propose in Theorem~\ref{main:BC fullconv} a heuristically optimal choice
(see Remark~\ref{stationary choice}) for the Euclidean norm.   \endofremark
 
\section{Filtered evaluation modules and Functional Transcendence} \label{functional transcendence}  

In this section, we develop the structure of the vanishing filtration jumps of the 
multivariable evaluation modules of auxiliary polynomial functions
that we described in~\S~\ref{filtrations}. 
Their formalism and the basic facts are collected in~\S~\ref{evaluation}, where we prove the commutativity in the formation
of Cartesian products and vanishing filtration jump exponent vectors of the evaluation modules. 
In~\S~\ref{functional bad approximability}, we
survey some of the literature on the classical Shidlovsky lemma from the historical proof~\cite{Shidlovsky1959} of the 
Siegel--Shidlovsky theorem on special values of $E$-functions, 
 and the deeper work of Chudnovsky and Osgood on the functional Schmidt subspace theorem --- \emph{Kolchin's problem} --- in differential algebra. 
 This finer information simplifies the statements and proofs of our arithmetic holonomy bounds, and it is furthermore indispensable for any refinements 
to quantitative linear independence measures and Diophantine inequalities. We do nevertheless remark that, following indications  in~\S~\ref{sec_withoutShid},  one could in principle dispense with the 
differential algebra theorems for the particular \emph{qualitative} linear independence proofs in this paper. Finally, just to give a sense of completeness and a proper historical context, we collect in~\S~\ref{PerfectPadeExamples} some of the most basic examples of perfect Pad\'e approximants to holonomic functions, which can be considered as a prototype and an introduction to the functional bad approximability theorems collected in~\S~\ref{functional bad approximability}.

\subsection{The evaluation module for Cartesian products} \label{evaluation}

We formalize the discussion of~\S~\ref{filtrations}. 

\subsubsection{Evaluation module} Consider a bounded Lebesgue-measurable subset $\Omega \subset [0,\infty)^d$
and a finite indexing set~$I$. In practice we will think of~$I$ as a subset of $\{1,\ldots,m\}^d$, and so we will use the boldface notation for the index elements
$\bi \in I$.  We fix for the time being an $I$-tuple of $\Q(\bx)$-linearly independent formal power series
$$
f_{\bi}(\bx) \in \Q \llbracket \bx \rrbracket, \qquad \bx := (x_1,\ldots,x_d), \quad \bi \in I.
$$
{\it The finite-rank free
 $\Z$-modules  in the following will all depend on the given function $f_{\bi}$, which will be considered as fixed and dropped from the notation. }

\begin{df} \label{eval map}
The {\emph{evaluation module}} $(E_{D,\Omega}^I, \psi_D)$ defined by the data 
\[\left( (f_{\bi})_{i \in I}; D; \Omega \right)\]
 is a pair of a finite-rank free $\Z$-module $E_{D,\Omega}^I$ together with
a $\Z$-module homomorphism $\psi_D : E_{D,\Omega}^I \to \Q \llbracket \bx\rrbracket$, constructed as follows. For $E_{D,\Omega}^I$, 
take the $\Z$-linear span of all $I$-tuples of monomials $\bx^{\bk(\bi)}$ with $\bk(\bi) \in (D \cdot \Omega) \cap \Z^d$  for all $\bi \in I$.
For $\psi_D$, we take the Taylor series development map $\Z$-linearly generated on the basis by the Taylor expansion of $\mathbf{x^{k(i)}}  f_{\bi}(\bx)$: 
$$
\psi_D \, : \,  E_{D,\Omega}^I \otimes_{\Z} \Q  \hookrightarrow  \Q \llbracket \bx\rrbracket, \qquad   \mathbf{x^{k(i)}} \mapsto \mathbf{x^{k(i)}}  f_{\bi}(\bx) \in \Q \llbracket \bx \rrbracket.
$$
\end{df}
The evaluation map $\psi_D$ is an \emph{injective} homomorphism, precisely by the $\Q(\bx)$-linear independence we assumed on the
 $I$-tuple $f_{\bi}(\bx) \in \Q\llbracket \bx \rrbracket$.

\subsubsection{Evaluation filtration} \label{vanfiltration}  
Consider the filtration of the infinite-dimensional $\Q$-vector space $\Q \llbracket \bx \rrbracket$ of formal power series in $d$ commuting variables $\bx := (x_1,\ldots,
x_d)$, obtained by firstly grading the monomial basis $\bx^{\bn}$ by the total degree $|\bn| := n_1 +\ldots +n_d$, and then filtering the $\binom{n+d-1}{d-1}$-dimensional $\Q$-vector space homogeneous piece of degree-$n$ elements by the lexicographical ordering of the exponents $\bn = (n_1,\ldots,n_d)$ with $|\bn| = n$:
\[
\mathbf{m} \prec \mathbf{n}  \Longleftrightarrow  \textrm{either $|\mathbf{m}| < |\mathbf{n}|$,
 or $|\mathbf{m}| = |\mathbf{n}|$
and $\mathbf{m}$ precedes $\mathbf{n}$ lexicographically}.  
\]
We denote the successor function of this total ordering by $\bn \mapsto \bn^+$. 
 The resulting filtration 
on 
$$
\Q\llbracket \bx \rrbracket  =: F = \bigcup_{\substack{ \bn \in \NwithzeroB^d, \, \prec }} F^{(\bn)}
$$
 is split and \emph{maximal}: 
 the successor quotient $\Q$-vector spaces $F^{(\bn)}/F^{(\bn^+)} \cong \Q$ are one-dimensional with basis the class of the unique monomial $\bx^{\bn}$ in $F^{(\bn)} \setminus F^{(\bn^+)}$. Using the monomorphism $\psi_D : E_{D,\Omega}^I \hookrightarrow F$, the $(\NwithzeroB^d,\prec)$-filtration on~$F$ induces an
$(\NwithzeroB^d,\prec)$-filtration on the $\Q$-vector space $E_{D,\Omega}^I \otimes_{\Z} \Q$: 
$$
E_{D,\Omega}^{I, (\bn)}  := \psi_D^{-1} \left( F^{(\bn)} \right) \subset E_{D,\Omega}^{I}  \otimes_{\Z} \Q.
$$
The monomorphism $\psi_D : E_{D,\Omega}^{I} \hookrightarrow F$ induces a \emph{monomorphism} on the graded successive quotients: 
$$
\psi_D^{(\bn)} : E_{D,\Omega}^{I, (\bn)}  / E_{D,\Omega}^{I, (\bn^+)} \hookrightarrow F^{(\bn)}/F^{(\bn^+)}.
$$
and since the codomain of this monomorphism is the one-dimensional $\Q$-vector space $F^{(\bn)}/F^{(\bn^+)}$, it follows that
\begin{equation}  \label{van fj}
\forall \bn \in \NwithzeroB^d, \qquad
\dim_{\Q}  \left( E_{D,\Omega}^{I, (\bn)}  / E_{D,\Omega}^{I, (\bn^+)} \right)   \in \{0, 1\}. 
\end{equation}
The sum of all those $\{0,1\}$ dimensions equals $\dim_{\Q} (E_{D,\Omega}^{I} \otimes \Q) = \rank(E_{D,\Omega}^{I})$. It follows that 
the \emph{vanishing filtration jumps} set
$$
\mathcal{V}_{D,\Omega}^I  :=  \left\{ \bn \in \NwithzeroB^d \, : \,  \dim_{\Q}  \left( E_{D,\Omega}^{I, (\bn)}  / E_{D,\Omega}^{I, (\bn^+)} \right)  = 1 \right\}
\subset \NwithzeroB^d
$$
satisfies 
\begin{equation}  \label{cv}
\# \mathcal{V}_{D,\Omega}^I = \mathrm{rank}(E_{D,\Omega}^I), \quad \textrm{ and } \quad E_{D,\Omega}^{I, (\bn^+)} = E_{D,\Omega}^{I, (\bn)}  
\textrm{ for } \bn \notin \mathcal{V}_{D,\Omega}^I. 
\end{equation}
The $\prec$ filtration also shows that $\bn \in \mathcal{V}_{D,\Omega}^I$ are exactly the exponents that occur in the monomials $\mathbf{x^n}$ in the
leading order jet $|\bn| = \mathrm{ord}_{\mathbf{x=0}}(F)$ of some nonzero element $F \in  \psi_D\left(E_{D,\Omega}^{I}) \right) \setminus \{0\}$. 

We have proved: 
\begin{lemma} \label{total jumps}
Under the total ordering $(\NwithzeroB^d, \prec)$ and the above premise of the $\Q(\bx)$-linear independence of the power series $(f_{\bi})_{\bi \in I}$, 
there are precisely $\rank(E_{D,\Omega}^{I})$ exponent vectors $\bn \in \NwithzeroB^d$ for which there exists a nonzero element 
$G \in \psi_D\left(E_{D,\Omega}^{I}) \right) \subset \Q \llbracket \bx \rrbracket$  whose $\prec$-minimal exponent monomial is $c \, \mathbf{x^n}$
for some nonzero $c \in \Q^{\times}$. 

Furthermore, for any $G \in \psi_D\left(E_{D,\Omega}^{I}) \right) \setminus \{0\}$ of $\mathbf{x=0}$ 
vanishing order~$n$, and for any nonzero monomial $c \, \mathbf{x^n}$ in~$G$ of the minimal degree~$|\bn| = n$, there exists 
an $F  \in \psi_D\left(E_{D,\Omega}^{I}) \right)$  with $\prec$-minimal monomial~$\mathbf{x^n}$. 
\end{lemma}

\subsubsection{Cartesian products}    Consider two
bounded Lebesgue-measurable subsets $\Omega_1 \subset [0,\infty)^{d_1}$ and $\Omega_2 \subset [0,\infty)^{d_2}$, a positive integer $D \in \NwithoutzeroA$, 
two respective index sets $I_1$ and $I_2$ as above, and for each $h \in \{1,2\}$, two respective $I_h$-tuples of 
$\Q\left(x^{(h)}_1, \ldots, x^{(h)}_{d_h}\right)$-linearly independent formal power series $\{ f_{\bi} \}_{\bi \in I_1}, \{g_{\mathbf{j}} \}_{\mathbf{j} \in I_2}$. 
These respective data sets define two evaluation modules $\psi_D^{(h)} : E_{D,\Omega_h}^{I_h} \hookrightarrow \Q \llbracket  \mathbf{x}^{(h)} \rrbracket$, 
as well as a $(d_1+d_2)$-variable evaluation module $ E_{D, \Omega_1 \times \Omega_2}^{I_1 \times I_2}$, the \emph{Cartesian product}, defined by 
the $\Q(\bx^{(1)}, \bx^{(2)})$-linearly independent $(I_1 \times I_2)$-tuple of formal power series $f_{\bi}(\bx)g_{\mathbf{j}}(\mathbf{y})$.
There is hence a tautological $\Z$-module isomorphism
\begin{equation}  \label{prod mod}
\begin{aligned} 
E_{D, \Omega_1}^{I_1} \times E_{D, \Omega_2}^{I_2}  &  \iso E_{D, \Omega_1 \times \Omega_2}^{I_1 \times I_2}, \\
  \left( f_{\bi}(\bx), g_{\mathbf{j}}(\mathbf{y}) \right)_{\mathbf{i} \in I_1, 
\mathbf{j} \in I_2} & \mapsto  \left(f_{\bi}(\bx)g_{\mathbf{j}}(\mathbf{y})\right)_{(\mathbf{i,j}) \in I_1 \times I_2},
\end{aligned}
\end{equation}
under which the evaluation map $\psi_D :  E_{D, \Omega_1 \times \Omega_2}^{I_1 \times I_2} \hookrightarrow \Q \llbracket
\mathbf{x}^{(1)}, \mathbf{x}^{(2)} \rrbracket$ of the product commutes with the product of the evaluation maps 
\[(\psi_D^{(1)}, \psi_D^{(2)}) : E_{D,\Omega_1}^{I_1} \times E_{D,\Omega_2}^{I_2}  \hookrightarrow \Q \llbracket
\mathbf{x}^{(1)}, \mathbf{x}^{(2)} \rrbracket.\]
  In combination with Lemma~\ref{total jumps}, this remarks leads to the following key result, which we can most succinctly 
  express by saying that the formation of Cartesian products of evaluation modules commutes with the formation of their vanishing filtration jumps. 
  
  \begin{lemma}  \label{commutation}
  Fix the $\Q(\bx^{(1)})$-linearly independent $I_1$-tuple $\left( f_{\bi}(\bx^{(1)}) \right)_{\bi \in I_1}$ and  the $\Q(\bx^{(2)})$-linearly independent $I_2$-tuple $\left( g_{\mathbf{j}}(\bx^{(2)}) \right)_{\mathbf{j} \in I_2}$. 
  Under the notation and premises of the current~\S~\ref{evaluation}, the vanishing filtration jumps of the associated evaluation modules satisfy
  \begin{equation}  \label{commuting}
  \mathcal{V}_{D,\Omega_1 \times \Omega_2}^{I_1 \times I_2} = \mathcal{V}_{D,\Omega_1}^{I_1} \times \mathcal{V}_{D,\Omega_2}^{I_2},
  \end{equation} 
  as subsets of $\NwithzeroB^{d_1 + d_2} = \NwithzeroB^{d_1} \times \NwithzeroB^{d_2}$. 
  \end{lemma}
  
  In view of the importance of this basic lemma for the sequel, we give two proofs. 
  
  \begin{proof}[First proof of Lemma~\ref{commutation}]
  By~\eqref{cv} and~\eqref{prod mod}, the two sets in the asserted equality~\eqref{commuting} are finite and of the same cardinality 
  \[\rank(E_{D,\Omega_1 \times \Omega_2}^{I_2 \times I_2}) = \rank(E_{D,\Omega_1}^{I_1}) \, \rank(E_{D,\Omega_2}^{I_2}).\] It 
  therefore suffices to prove that one of these sets is contained by the other. But
   \[\mathcal{V}_{D,\Omega_1}^{I_1} \times \mathcal{V}_{D,\Omega_2}^{I_2}
  \subseteq   \mathcal{V}_{D,\Omega_1 \times \Omega_2}^{I_1 \times I_2}\]
   is made clear by the following product construction. 
  For any pair $\bn_1 \in \mathcal{V}_{D,\Omega_1}^{I_1} \subset \NwithzeroB^{d_1}$ and $\bn_2 \in \mathcal{V}_{D,\Omega_2}^{I_2} \subset \NwithzeroB^{d_2}$, there
  exist by definition two auxiliary function evaluations $G_1(\bx^{(1)}) \in \psi_D\left(E_{D,\Omega_1}^{I_1}\right) \subset \Q \llbracket \bx^{(1)} \rrbracket$
  and $G_2(\bx^{(2)}) \in \psi_D\left(E_{D,\Omega_2}^{I_2}\right) \subset \Q \llbracket \bx^{(2)} \rrbracket$ such that $\bn_1$ is the $\prec$-minimal exponent in
  $G_1(\bx^{(1)}) = c_1 \cdot (\bx^{(1)})^{\bn_1} +\ldots$ among the totally ordered exponent set $(\NwithzeroB^{d_1}, \prec)$, and $\bn_2$ is the $\prec$-minimal exponent in
  $G_2(\bx^{(2)}) = c_2 \cdot (\bx^{(2)})^{\bn_2} +\ldots$ among the totally ordered exponent set $(\NwithzeroB^{d_2}, \prec)$. Consider then the product
  $$
  G(\bx^{(1)}, \bx^{(2)}) := G_1(\bx^{(1)})G_2(\bx^{(2)})   \in \Q \llbracket \bx^{(1)}, \bx^{(2)} \rrbracket. 
  $$
  By construction of the Cartesian product,  we see that $G$ belongs to the evaluation range $\psi_D\left(E_{D,\Omega_1 \times \Omega_2}^{I_2 \times I_2} \right)$. 
  It has the nonzero coefficient $c_1c_2 \in \Q^{\times}$ in the multidegree $(\bn_1,\bn_2) \in \NwithzeroB^{d_1} \times \NwithzeroB^{d_2} = 
  \NwithzeroB^{d_1+d_2}$. We claim that for the $(\NwithzeroB^{d_1+d_2}, \prec)$ total ordering of the exponents, this is the minimal multidegree in~$G$.  
  It is certainly of the minimal possible vanishing order $|\bn_1| + |\bn_2|$, for by definition of $\prec$ the factor power series $G_1$
  and $G_2$ have respective vanishing orders $|\bn_1|$ and $|\bn_2|$. Now the monomial degrees $(\mathbf{u,v})$ in $G$ having the minimal possible order $|\mathbf{u}|
  + |\mathbf{v}| = |\bn_1| + |\bn_2|$ 
  have, by $|\mathbf{u}| \geq |\bn_1|$ and $|\mathbf{v}| \geq |\bn_2|$, partial degrees $|\mathbf{u}| = |\mathbf{n_1}|$ and $|\mathbf{v}| = |\mathbf{n_2}|$. 
  We have $\bn_1 \preceq \mathbf{u}$ and $\bn_2 \preceq \mathbf{v}$. 
  By the definition of the lexicographical ordering, if $(\mathbf{u,v}) \neq (\bn_1, \bn_2)$, it follows that  $(\bn_1, \bn_2) \prec (\mathbf{u,v})$. 
  Hence, through the example of~$G(\bx^{(1)},\bx^{(2)}) = c_1c_2 \cdot (\bx^{(1)})^{\bn_1}(\bx^{(2)})^{\bn_2} +\ldots$, we have found that
  $(\bn_1, \bn_2) \in   \mathcal{V}_{D,\Omega_1 \times \Omega_2}^{I_1 \times I_2} $, and 
  in this way we have proved the requisite inclusion $ \mathcal{V}_{D,\Omega_1}^{I_1} \times \mathcal{V}_{D,\Omega_2}^{I_2} 
  \subseteq \mathcal{V}_{D,\Omega_1 \times \Omega_2}^{I_1 \times I_2}$. \end{proof} 
  
    \begin{proof}[Second proof of Lemma~\ref{commutation}]
    We can also directly see the reverse inclusion, 
    \[\mathcal{V}_{D,\Omega_1 \times \Omega_2}^{I_1 \times I_2} 
    \subseteq \mathcal{V}_{D,\Omega_1}^{I_1} \times \mathcal{V}_{D,\Omega_2}^{I_2}.\]
    Let 
    $$
    F(\bx^{(1)}, \bx^{(2)}) = \sum_{\bi \in I_1, \, \mathbf{j} \in I_2} Q_{\mathbf{i,j}}(\bx^{(1)},\bx^{(2)}) 
    f_{\bi}(\bx^{(1)}) g_{\mathbf{j}}(\bx^{(2)}) \in \psi_D \left( E_{D,\Omega_1 \times \Omega_2}^{I_1 \times I_2} \right)
    $$
    be an arbitrary auxiliary function evaluation of the Cartesian product module, with $(\NwithzeroB^{d_1+d_2}, \prec)$ minimal
    monomial $\beta \cdot (\bx^{(1)})^{\bn_1} (\bx^{(2)})^{\bn_2}$. 
    Then
$$
\frac{1}{\bn_2!}\frac{\partial^{|\bn_2|}}{(\partial \bx^{(2)})^{\bn_2}} \Big|_{\bx^{(2)} = \mathbf{0}} \Big\{ F(\bx^{(1)},\bx^{(2)}) \Big\} = \beta \cdot (\bx^{(1)})^{\bn_1} +\ldots 
\textrm{\big[$\prec$-higher terms\big]}, 
$$
for the monomials of this specialization are exactly the $\gamma \cdot (\bx^{(1)})^{\bk}$
such that $\gamma \cdot  (\bx^{(1)})^{\bk} (\bx^{(2)})^{\bn_2}$ are monomials from $F(\bx^{(1)},\bx^{(2)})$. 
  In this way, $\bn_1 \in  \mathcal{V}_{D,\Omega_1}^{I_1}$.  Similarly, $\bn_2 \in  \mathcal{V}_{D,\Omega_2}^{I_2}$, and the 
  requisite inclusion $\mathcal{V}_{D,\Omega_1 \times \Omega_2}^{I_1 \times I_2} 
    \subseteq \mathcal{V}_{D,\Omega_1}^{I_1} \times \mathcal{V}_{D,\Omega_2}^{I_2}$ is proved. 
  \end{proof}

 We record the main corollary we will use. 

\begin{cor} \label{jumps Cartesian}
Consider a positive integer $D \in \NwithoutzeroA$ and an $m$-tuple $f_1, \ldots, f_m$ in $\C \llbracket x \rrbracket$ of $\C(x)$-linearly independent 
formal power series. For these data, there exists a sequence 
$$
0 \leq u(1) < \cdots < u(mD)
$$
 of $mD$ non-negative integers such that the following holds for every $d = 1, 2, 3, \ldots$: 
 
 In every nonzero formal power series of the shape
 $$
F(\mathbf{x}) := \sum_{\bi \in \{1,\ldots,m\}^d} Q_{\bi}(\mathbf{x}) f_{i_1}(x_1) \cdots f_{i_d}(x_d) \in \C \llbracket x_1, \ldots, x_d \rrbracket \setminus \{0\}, 
$$
where $Q_{\bi}(x_1, \ldots, x_d) \in \C [x_1,\ldots,x_d]$ are polynomials having all their partial degrees $\deg_{x_j}Q_{\bi} < D$, all monomials
$\beta \mathbf{x^n}$ with minimal total degree $|\bn| = n_1 +\ldots + n_d$ have
$$
\bn \in \big\{ 0 \leq u(1) < \cdots < u(mD) \big\}^d \subset \NwithzeroB^d.
$$
\end{cor}

\begin{proof}
Take the evaluation module $\psi_D: E_{D,[0,1)}^{\{1,\ldots,m\}} \hookrightarrow \Q \llbracket x \rrbracket$ defined by the $\Q(x)$-linearly 
independent power series $f_1, \ldots, f_ \in \Q \llbracket x \rrbracket$. Clearly, $\rank\left(E_{D,[0,1)}^{\{1,\ldots,m\}} \right) = mD$. 
We define 
$$
\big\{ 0 \leq u(1) < \cdots < u(mD) \big\} = \mathcal{V}_{D,[0,1)}^{\{1,\ldots,m\}}  \subset \NwithzeroB
$$
to be the vanishing filtration jumps for this single-variable evaluation module. By Lemma~\ref{commutation},
it then follows for each $d = 1, 2, 3, \ldots$ that the Cartesian $d$-th power evaluation module $\left( E_{D,[0,1)^d}^{\{1,\ldots,m\}^d} , \psi_D \right)$
defined by the $\Q(x_1,\ldots,x_d)$-linearly independent formal power series $f_{\bi}(\xi) := f_{i_1}(x_1) \cdots f_{i_d}(x_d)$ has vanishing filtration jumps
at exactly the $d$-th Cartesian power set
$$
\mathcal{V}_{D,[0,1)^d}^{\{1,\ldots,m\}^d} =  \big\{ 0 \leq u(1) < \cdots < u(mD) \big\}^d \subset \NwithzeroB^d.
$$
The result now follows by Lemma~\ref{total jumps}. 
\end{proof}

\subsection{Functional bad approximability}  \label{functional bad approximability}
When the functions
 $f_1,\ldots,f_m \in \Q\llbracket x \rrbracket$ are holonomic,  it turns out possible to almost completely determine the~$mD$ vanishing filtration
jumps of the ensuing univariate evaluation module $E_{D} := E_{D,[0,1)}^{\{1,\ldots,m\}}$. This is the content of the Chudnovsky--Osgood theorem~\ref{KolchinSolved}, which can be seen as a functional analog for holonomic functions of the Roth--Schmidt bad approximability theorem. The roots of all of this are in Hermite's memoir (discussion in~\S~\ref{PerfectPadeExamples} below) on the exponential function and the transcendence of the number~$e$.

\begin{basicremark}
For the system $\{f_1, \ldots, f_m\} = \{ e^{\alpha_1x},\ldots,e^{\alpha_m x}\}$ of pairwise distinct exponential functions,
Hermite~\cite{Hermite2}, 
\silentcomment{this is also \cite[vol.~4, pages 357--377]{HermiteOeuvres},}
in a letter published in~1893 (after having published similar formulas already in~\cite{Hermite}), found the explicit $\C[x]$-linear combination of the maximal 
$x=0$ vanishing order for an arbitrary degree vector $(D_1,\ldots,D_m) \in \NwithzeroB^m$: 
\begin{equation} \label{Hermite2}
\begin{aligned}
\int_{|z| = R} \frac{e^{xz} \, \mv(z)}{(z-\alpha_1)^{D_1+1} \cdots (z-\alpha_m)^{D_m+1}}
=: \sum_{i=1}^{m} P_i(x) e^{\alpha_i x}, \qquad R > \max_i{|\alpha_i|}, \\ 
= \frac{1}{(D_1+\ldots+D_m+m)!} \, x^{-1 + \sum_{i=1}^{m} (D_i+1)} + O\left(x^{\sum_{i=1}^{m} (D_i+1)} \right),
\end{aligned}
\end{equation}
This follows upon unfolding the residue calculus of the complex contour integral 
and finding the thus-explicitable polynomials $P_1,\ldots,P_m$ to have the \emph{exact} degrees $\deg{P_i} = D_i$.
The right-hand side of~\eqref{Hermite2} follows by~$D_i$ partial integrations
upon computing the integrand residues at the poles $z = \alpha_1, \ldots, \alpha_m$ in the bounded component of $\C \setminus \{ |z| = R\}$; on the other hand, the $x=0$ exact vanishing order development~\eqref{Hermite2} follows by computing the residue at the unique pole $z = \infty$ in the complementary
component.  Having for these particular polynomials --- the so-called \emph{type I Hermite--Pad\'e approximants} --- the exact degrees $\deg{P_i} = D_i$ (which can furthermore be taken
completely arbitrary), and this
exact vanishing order~\eqref{Hermite2},  proves  by an argument similar to the proof of Lemma~\ref{total jumps} that for \emph{arbitrary} polynomials $Q_1, \ldots, Q_m \in \C[x]$, 
the strongest possible form of (rational) functional bad approximability is in place: 
  \begin{equation} \label{bad rational approach}
\mathrm{ord}_{x=0} \left(  Q_1 f_1 +\ldots + Q_mf_m  \right)   \leq \sum_{i=1}^{m} \left( \deg{Q_i} + 1 \right) - 1. 
  \end{equation}
   Mahler~\cite{MahlerPerfect} termed such systems~\emph{perfect}, and found a few other examples (incidentally obtainable from Hermite's formula 
   by a substitution and a limit~\cite[page~331]{ChudnovskyThueSiegel}), including the binomial system 
   \begin{equation} \label{Binomials system}
   \{(1-x)^{\alpha_1}, \ldots, (1-x)^{\alpha_m}\},  \qquad \textrm{when all } \alpha_i - \alpha_j \notin \Z,
   \end{equation}
    and, under the additional constraint\footnote{Sometimes termed \emph{weak perfection}.} $D_1 \leq D_2 \leq \cdots \leq D_m$, the logarithm system~\cite{MahlerLog,Jager}
    \begin{equation} \label{log Mahler}
    \{ 1, \log(1-x), \log^2(1-x), \ldots, \log^{m-1}(1-x) \}. 
    \end{equation} 
    \silentcomment{At least in the diagonal (restricted) case, the  perfectness of the system $\left\{ 1,   \pFq{2}{1}{a,,1}{c}{x} \right\}$ had been already known to Jacobi~\cite[\S~8]{Jacobi}.}
   In~\cite{MahlerPi}, Mahler used the explicit linear forms for the system~\eqref{log Mahler} to prove
   the explicit inequality
   $\left| \pi - p/q \right| > q^{-42}$ for all positive integers $p,q \geq 2$. Gregory Chudnovsky~\cite{ChudnovskyHermite,ChudnovskyThueSiegel,ChudnovskyExt} has
   used (like Thue, Siegel, and Baker before him, cf.~\S~\ref{binomials discussion}) the systems~\eqref{Binomials system} and~\eqref{log Mahler} to derive excellent effective irrationality exponents for suitable
   roots $\sqrt[n]{b/a}$ from rational numbers, as well as for logarithms of rational numbers. A fairly general class of perfect systems are the \emph{Angelesco--Nikishin systems}~\cite{Angelesco,Nikishin,SorokinN,NikishinSorokin} in the
   theory of the Cauchy transform and orthogonal polynomials; their perfection
   was proven in full generality by Fidalgo Prieto and L\'opez-Lagomasino \cite{PrietoLagomasino,PrietoLagomasino2}. The fact that the Pad\'e approximants to the 
   polylogarithm system $\{f_1, \ldots, f_m\} = \{1, \li_1, \li_2, \ldots, \li_{m-1}\}$ turn out to be Angelesco--Nikishin systems was at the root\footnote{As a starting or inspiration point, even though they ultimately devised a different (but related) function system, see~\cite[Th\'eor\`eme~1]{FischlerRivoalPade} or~\cite[\S~2.4]{FischlerApery}.} of
   Ball and Rivoal's work~\cite{Rivoal1,BallRivoal} on the arithmetic of zeta values. 
   
 In the particular case $D_ 1 = \cdots = D_m = D-1$, we can equivalently express the functional bad approximability property~\eqref{bad rational approach} into 
 the framework of~\S~\ref{evaluation}: it precisely means that 
 the evaluation module~$E_D$ has the vanishing filtration jumps set
\begin{equation}   \label{the generic jumps}
\mathcal{V}_{D,[0,1)}^{\{1,\ldots,m\}} = \{0, 1, \ldots, mD-1\}. 
\end{equation} 
 \endofremark
\end{basicremark}
In differential algebra, as we briefly indicated in~\S~\ref{holonomic subspace theorem}, Kolchin~\cite{Kolchin} proved an analog of
Liouville's Diophantine inequality, and
asked\footnote{From~\cite{Kolchin}:
``It remains to make the obvious remark,
in view of the deep Thue--Siegel--Roth improvement 
of Liouville's theorem (see K.~F.~Roth, Mathematika vol. 2 (1955) pp. 1--20), that it would
be desirable to obtain a similar improvement in the present theorem.''}
 for an analog of the fundamental theorem on algebraic numbers that Roth had proved four years prior:

\begin{problem}[Kolchin's Problem] \label{kolchinsproblem}
 Given a non-rational formal power series solution $f \in \C \llbracket x \rrbracket \setminus \C(x)$ of some linear ODE $L(f) = 0$ over $\C(x)$, to prove that 
$(2+\varepsilon) \max(\deg{P},\deg{Q}) + O_{\varepsilon,f}(1)$ is the highest $x=0$ vanishing order possible for the error $f(x) - P(x)/Q(x)$ in any rational function
approximation.
\end{problem}

In fact, Kolchin's setup was more general and not limited to linear ODEs; he worked in an arbitrary nontrivial valued differential field, and his Liouville inequality~\cite[\S~5]{Kolchin} thus also applied to arbitrary (nonlinear) ODEs over~$\C(x)$. The first such result, weaker than Kolchin's, appears to be Maillet's~\cite[page~266]{Maillet}.

In terms of our evaluation modules in~\S~\ref{evaluation}, Kolchin's Liouville-type result for the case of a formal power series solution to an $r^{\textrm{th}}$ order linear ODE~$\LL(f) = 0$ can be expressed by saying that the~$2D$ vanishing filtration jumps in the module 
defined by $\{f_1, f_2\} = \{1, f\}$ are contained by $\left\{0,1, \ldots, rD+O_{\LL}(1)\right\}$. His (implicit) Roth-type conjecture is that they should in fact be contained  by $\left\{ 0, 1, \ldots, (2+\varepsilon) D + O_{\varepsilon,\LL}(1) \right\}$ for every $\varepsilon > 0$. 

  Independently in the same year,
Shidlovsky~\cite{Shidlovsky1959} (see also~\cite[\S~3.5, Lemma~8]{Shidlovsky}, \cite[\S~VII.3]{Lang}, or~\cite[\S~4]{MahlerBook}) discovered a more accurate  form of the functional
Liouville inequality in the case of linear ODEs over~$\C(x)$, and used it to complete the main results of Siegel's algebraic independence theory~\cite{Siegel1929SNS,SiegelBook} for special values of $E$-functions.  
We reformulate Shidlovsky's lemma into our language of the vanishing filtration jumps. 

\begin{thm}[Shidlovsky]\label{Shidlovsky}
For the case of a system $\{f_1, \ldots, f_m\}$ whose $\Q(x)$-linear span is $m$-dimensional and closed under the derivation $d/dx$, there exists a constant $C = C(f_1, \ldots, f_m)$
such that, for every $D \in \NwithoutzeroA$, the vanishing filtration jumps of the evaluation module $E_{D,[0,1)}^{\{1,\ldots,m\}}$ satisfy
\begin{equation} \label{Shidlovsky jumps}
\mathcal{V}_{D,[0,1)}^{\{1,\ldots,m\}} \subseteq \{0, 1, \ldots, mD + C\}, \qquad \#\mathcal{V}_{D,[0,1)}^{\{1,\ldots,m\}} = mD. 
\end{equation}
\end{thm}
Although Shidlovsky's original work did not supply an effective procedure to compute the constant $C$ out of the rank-$m$ first-order linear differential
system $\mathbf{y}' = A \mathbf{y}$ that has $\mathbf{y} = (f_1,\ldots, f_m)^{\mathrm{t}}$ as a solution, such theorems were eventually obtained, firstly by Chudnovsky~\cite[Corollary 11.3.10]{ChudnovskyShidlovsky} in the Fuchsian case (which is certainly the case we are concerned with, see Remark~\ref{remark on global nilpotence} below), and then eventually in the general case by Bertrand, Beukers, Chirskii, and Yebbou, see~\cite{BertrandBeukers} as complemented by~\cite[\S~3]{BertrandChirskiiYebbou}.  A far-reaching generalization of Shidlovsky's lemma is in Bertrand~\cite[Th\'eor\`eme~2]{Bertrand}. 

 We state only a crude version of Chudnovsky's result on the Fuchsian case.
A brief treatment of this explicit zero estimate is also sketched in Andr\'e's book~\cite[\S~III, Appendix]{Andre}. 

\begin{thm}[Chudnovsky] \label{Chudnovsky Fuchsian form}
Suppose the system $\mathbf{f} := \{f_1, \ldots, f_m\} \in \Q \llbracket x \rrbracket^m \setminus \left( x \Q\llbracket x \rrbracket \right)^m$ of $\Q(x)$-linearly independent formal power series arises as the full component vector of some solution $\mathbf{y} =\mathbf{f}^{\mathrm{t}}$ to a Fuchsian first-order linear differential system $ \mathbf{y}' = A \mathbf{y}$, where 
$A \in M_{m \times m}\left(\Q(x)\right)$. Let $S \subset \PP^1$ be the set of poles in the matrix of rational functions~$A$, and define $h := \sum_{s \in S} \varepsilon_s$, 
where $\varepsilon_s$ is the negative of the smallest real part of any exponent that occurs in the asymptotic development of any one of the functions $f_i(x)$ at the regular singular point $x=s$ of the Fuchsian ODE.

Then, for all $D \in \NwithoutzeroA$, the~$mD$ vanishing filtration jumps of the evaluation module $E_{D,[0,1)}^{\{1,\ldots,m\}}$
are contained by the set
\begin{equation} \label{Chudnovsky jumps}
\mathcal{V}_{D,[0,1)}^{\{1,\ldots,m\}} \subseteq \left\{0, 1, 2, \ldots, mD + (\#S-2) \binom{m}{2}  + mh   \right\};
\  \#\mathcal{V}_{D,[0,1)}^{\{1,\ldots,m\}} = mD. 
\end{equation}
\end{thm}

\begin{remark} \label{remark on global nilpotence}
Thanks to the work on the global nilpotence property by David and Gregory Chudnovsky~\cite{ChudnovskyG}, \cite[\S~VIII]{Dwork}, \cite[\S~VI]{AndreG}, \cite{diVizio} and the theorem of Honda and Katz~\cite[\S~III.6]{Dwork}, \cite[\S~IV.5.3]{AndreG}, the holonomic power series in all our (abstract) theorems in this paper are automatically of the Fuchsian class: they have only regular singular points (with rational exponents). Hence Theorem~\ref{Chudnovsky Fuchsian form} applies to them as an explicit Shidlovsky bound. For the proofs of Theorem~\ref{basic main} and all its generalizations, this theorem is already sufficient under the supplemental assumption --- which is satisfied in all the applications we could conceive of --- that the~$\Q(x)$-linear span of~$f_1, \ldots, f_m$ is closed under the derivation~$d/dx$. 
This remark could also have a significance for the project of refining our qualitative linear independence results to quantitative measures of linear independence. 

At the same time, the proof of the Chudnovskys's theorem (namely: of Galo\v{c}kin's canceling factorials property and the global nilpotence of an integrable connection that admits at least one $G$-series formal solution with~$\C(x)$-linearly independent components) \emph{itself}  relies on a suitable qualitative Shidlovsky lemma~\cite[Prop.~VIII.2.3]{Dwork}, \cite[Theorem~3.1, Lemma~8.3]{ChudnovskyG}, \cite[\S~VI.2]{AndreG} in the dual form for simultaneous --- that is now \emph{type II Hermite--Pad\'e} --- functional rational approximants~$f_i \approx P_i/Q$, $1 \leq i \leq m$,  selected to have integer coefficients of controlled size as provided by the Thue--Siegel lemma.~\endofremark
\end{remark}

All these theorems also embed as very special cases into the broader subject of zero multiplicity estimates for functions
satisfying a possibly nonlinear algebraic differential system. This path was opened up by the groundbreaking works of Nesterenko~\cite{Nesterenko2} and
Brownawell--Masser~\cite{BrownawellMasser}. We refer to Binyamini~\cite{Binyamini}  for a survey, a modern treatment, and refinements of a large
portion of the literature on this rather vast topic; and to~\cite{NesterenkoThm,NesterenkoPhilippon} for applications to algebraic independence. In 
the special context of linear ODEs,
the separate streams opened up by Kolchin and Shidlovsky converged in the early
1980s with the resolution of Problem~\ref{kolchinsproblem}, independently by David and Gregory Chudnovsky~\cite{Chudnovsky4} and Osgood~\cite{Osgood4}. 

\begin{thm}[Chudnovsky, Osgood] \label{KolchinSolved}
Consider an  arbitrary set $\{f_1, \ldots, f_m\}$ 
 of holonomic functions in $\Q \llbracket x \rrbracket$. That is,  our only assumption now
is that each of the formal power series $f_i(x)$ separately satisfies some nonzero linear ODE $\LL_i(f_i) = 0$. For an arbitrary $\varepsilon > 0$, 
there exists a constant $C(\varepsilon) = C(\varepsilon; \LL_1, \ldots, \LL_m)$, effectively computable from the arguments in~\cite[\S~2]{Chudnovsky4}, 
such that for all $D \in \NwithoutzeroA$, the vanishing filtration jumps of the evaluation module $E_D = E_{D,[0,1)}^{\{1,\ldots,m\}}$ satisfy
$$
\mathcal{V}_{D,[0,1)}^{\{1,\ldots,m\}}  \subseteq \big\{  0, 1, 2, \ldots, (m+\varepsilon) D + C(\varepsilon)    \big\},   \qquad \#\mathcal{V}_{D,[0,1)}^{\{1,\ldots,m\}}  = mD.
$$
\end{thm}

In conjunction with Corollary~\ref{jumps Cartesian}, these theorems can be summarized into the following proposition.

\begin{lemma} \label{lem_Shidlovsky} Let $f_1, \ldots, f_m \in \C \llbracket x \rrbracket$ be $\C(x)$-linearly independent \emph{holonomic} power series: there
exist nonzero linear differential operators $\LL_i$ over $\Q(x)$ with $\LL_i(f_i) = 0$.  Then,
 for every $\varepsilon > 0$, there exists a constant $C(\varepsilon) \in \R$, in principle effectively computable from the datum $(\varepsilon; \LL_1, \ldots, \LL_m)$ alone, such that the following
is true.
 
 We consider an arbitrary positive integer $d \in \NwithoutzeroA$, and write
$$
\bx := (x_1, \ldots, x_d), \qquad f_{\bi}(\bx):=\prod_{s=1}^d f_{i_s}(x_s) \quad \textrm{ for } \bi := (i_1, \ldots, i_d) \in \{1,\ldots,m\}^d.  
$$  
Consider further an arbitrary positive integer $D \in \NwithoutzeroA$ and, over $\mathbf{i} \in \{1,\ldots,m\}^d$, an arbitrary set of polynomials
$$
Q_{\mathbf{i}}(\mathbf{x}) \in \C[x_1, \ldots, x_d] \ \textrm{ with } \ \deg_{x_j}Q_{\mathbf{i}} < D \textrm{ for all } j \in \{1,\ldots,d\} \textrm{ and } \bi \in \{1,\ldots,m\}^d. 
$$
Then, in the nonzero formal power series 
$$
F(\mathbf{x}) := \sum_{\bi \in \{1,\ldots,m\}^d} Q_{\bi}(\mathbf{x}) f_{\bi}(\mathbf{x}) \in \C \llbracket x_1, \ldots, x_d \rrbracket \setminus \{0\}, 
$$
every lowest-order nonzero monomial term $\beta \, \mathbf{x^n}$ in $F(\mathbf{x})$ has necessarily an exponent vector $\mathbf{n} = (n_1,\ldots,n_d)$, all of 
whose components satisfy
$$
n_j \leq (m+\varepsilon)D + C(\varepsilon).
$$
If moreover the $\Q(x)$-linear span of $f_1,\ldots,f_m$ is closed under the derivation $d/dx$, then $\varepsilon = 0$ could be taken. 
\end{lemma}

\begin{remark}
David and Gregory Chudnovsky conjecture~\cite[page~5161]{Chudnovsky4} that $\varepsilon = 0$ could be taken in Theorem~\ref{KolchinSolved},
and therefore --- as a consequence --- also in Lemma~\ref{lem_Shidlovsky}. However,  this conjecture remains unproved even for the case~\cite{JulieWang} of algebraic functions.~\endofremark
\end{remark}

\medskip

\emph{At this point, for the logic of the proofs, the reader may skip directly ahead to~\S~\ref{concentration of measure}. The remainder of~\S~\ref{functional transcendence} collects some examples, placed in their historical context, behind the theorems that we borrowed without proof in~\S~\ref{functional bad approximability}. }

\subsection{Some explicit constructions of Hermite--Pad\'e approximants}  \label{PerfectPadeExamples}
  For the rest of~\S~\ref{functional transcendence}, we collect a few simplest and most fundamental illustrating examples, aiming at a modest attempt at sketching the
 historical seeds of some of the basic ideas in the proofs of the theorems on functional bad approximability that we collected in~\S~\ref{functional bad approximability}, but also of the broader concept of holonomy bounds and the way we use them in our present paper. A quintessential illustrating example for the key point in the proofs of Shidlovsky type theorems on functional bad approximability can be taken as the explicit (in the simple case outlined here) determinantal identity~\eqref{Thue det} from the theory of the hypergeometric ODE. The number-theoretic relevance of such identities  was found by Thue when he created the subject of non-effective Diophantine approximation. Our approach here to Ap\'ery limits has perhaps some faint similarity to Thue's paradigm with its organic ineffectivity; the proofs that we have of the explicit holonomy bounds of Theorem~\ref{basic main} do not\footnote{Except in the case~$\bb = \mathbf{0}$ of integer coefficients. In that very special case, even a much more precise integral finiteness counterpart
  is contained in the work of Bost and Charles~\cite[\S~9.1]{BostCharles}, in an implicitly effective form. } contain, even in principle, any effective procedure for the far more elusive problem of outputting a set of $\Q(x)$-vector space generators for the finite-dimensional holonomic module
 $\mathcal{H}(b_1,\ldots,b_r;\varphi)$  attached to a given holomorphic mapping~$\varphi : (\Db,0) \to (\C,0)$ paired up to a given denominators type~$\prod_{i=1}^r [1,\ldots, b_i  n]$ subject to~$|\varphi'(0)| > e^{b_1 +\ldots + b_r}$. Yet, when favored by the presence of suitable anchors (such as we have in~\S~\ref{sec:purefunctions}) and levers  (such as we have in~\S~\ref{sec:YtoY0(2)}, \S~\ref{jointli} and~\S~\ref{sec:diffGA}), the Diophantine repellency principles can occasionally be turned around into true Diophantine inequalities and linear independence proofs.  With Thue's method, it took over seventy years until a fairly general-scope theory, on a scale comparable to the Gelfond--Baker method of linear forms in logarithms,  started to emerge at the hands of Bombieri and his coauthors~\cite{BombieriDyson,BombieriMueller,BombieriGm,BombieriCohenGm}. 

But we should probably begin this discussion by delving a bit into our subject's proper origin: the work of Hermite by which he proved the 
transcendence of~$e$.

\subsubsection{Hermite approximations}  \label{sec:Hermite approximations} The memoir~\cite{Hermite} on the exponential function was based on the explicit
Hermite--Pad\'e approximants to the functions $1$, $e^x$, $e^{2x}, \ldots$, $e^{rx}$ in order to prove the transcendence of~$e$ by specializing $x := 1$.
For $r = 1$ the formula is
\begin{equation}
\begin{aligned}  \label{exp Hermite}
 &  \quad \pFq{1}{1}{-m}{-m-n}{x}  -   e^{x} \cdot  \pFq{1}{1}{-n}{-m-n}{-x}    \\ 
 & = - \frac{e^x  \int_0^1 e^{-t x} t^m(t-1)^{n} \, dt  }{(m+n)!}   \, x^{m+n+1}  \\
 & = (-1)^{n-1} \frac{m!n!}{(m+n)!(m+n+1)!} \, x^{m+n+1} + O(x^{m+n+2}) \\
 \end{aligned}
 \end{equation}
for the unique (up to scalar multiple) combination $B(x) - e^x A(x)$ with $\deg{A} \leq n, \deg{B} \leq m$ 
that vanishes at $x = 0$ to order at least $m+n+1$. As we can see from the explicit formula, the vanishing order
is in fact exactly equal to $m+n+1$. The existence of such a regular array of formulas further proves that for \emph{any}
pair $A(x), B(x)$ of nonzero polynomials of degrees $n = \deg{A}$ and $m = \deg{B}$, the combination $B(x) - e^x A(x)$ 
has $x = 0$ vanishing order at most $m+n+1$, with equality if and only if the form $B(x) - e^x A(x)$ is a scalar
multiple of~\eqref{exp Hermite}. This means that the holonomic function $e^x$ is very badly approximable by rational functions. 
Hermite's philosophy, which was later taken up by Siegel who started his~1929 paper~\cite{Siegel1929SNS} in 
exactly the same way as Hermite~\cite{Hermite} 
---  outlining an analogy between numbers, to be approximated
in the archimedean absolute value, and functions, to be expanded in power series and approximated in terms of the $x=0$ vanishing order,  ---  was that 
the functional formulas could be specialized at algebraic arguments to yield a full set of small linear forms with integer coefficients 
in the numbers (the special values) of interest; which in turn can often be used to prove the $\Q$-linear independence of those numbers. The bad approximability property serves as the sieve for expressing and recognizing a full (linearly independent) set of linear forms, both in the holonomic functions and in their special values, once these are constructed to be reasonably small:
as in the functional formula~\eqref{exp Hermite} and its specializations at the algebraic arguments.  
In Hermite's method, the functional bad approximability of~$e^x$ (suitably generalized to include all the powers $1, e^x, e^{2x}, \ldots, e^{rx}$), {\it via} 
identities such as~\eqref{exp Hermite}, can be
specialized at $x := \alpha \in \Qbar^{\times}$ to derive the Hermite--Lindemann--Weierstrass theorem on the transcendence, and furthermore the bad approximability, of the special value~$e^{\alpha}$. 

The content of Shidlovsky's lemma~\ref{Shidlovsky} and the Chudnovsky--Osgood theorem~\ref{KolchinSolved} can approximately be described as the statement that a property (only very slightly relaxed) of bad approximability by rational functions is in place
for any set of holonomic functions. We illustrate this on the most classical cases of perfect systems. 

\subsubsection{The Hermite--Pad\'e approximants to $(1-x)^{\nu}$} \label{binomials discussion}
We have the hypergeometric polynomials identity of Jacobi~\cite[\S~8]{Jacobi} to describe explicitly the Pad\'e table
for the binomial function (cf. \cite[page~75]{Siegel1929SNS}, or Siegel's introductory
paper for Thue's Selected Works volume~\cite[\S~2]{ThueSelected}): 
\begin{equation}
\begin{aligned}
\label{hyper poly}
  & \quad \pFq{2}{1}{-\nu-n,,-m}{-m-n}{x}  -    (1-x)^{\nu} \cdot  \pFq{2}{1}{\nu-m,,-n}{-m-n}{x}   
 \\ 
 &  = (-1)^n \frac{\binom{m+\nu}{m+n+1}}{\binom{m+n}{n}}    \pFq{2}{1}{-\nu + n+1 ,m+1}{m+n+2}{x}  \, x^{m+n+1} \\ 
 & = (-1)^n \frac{\binom{m+\nu}{m+n+1}}{\binom{m+n}{n}} \, x^{m+n+1} + O(x^{m+n+2}),
 \end{aligned}
 \end{equation}
proved for instance by verifying that all three terms satisfy the second-order Gauss hypergeometric equation
with parameters $\alpha = -\nu-n, \beta = -m, \gamma = -m-n$, and are therefore $\C$-linearly dependent (the 
argument is also in~\cite[Hilfssatz~1]{SiegelThue}). 
In this identity, the hypergeometric series on the left-hand side~\eqref{hyper poly} terminate to polynomials of degrees~$m$ and~$n$,
and so they give precisely the $[m/n]$ Hermite--Pad\'e approximant $B_{m,n}(x) - (1-x)^{\nu} A_{m,n}(x)$ to the binomial function $(1-x)^{\nu}$, for any $\nu \in \C \setminus \Z$.  
As the $x^{m+n+1}$ coefficient in the formula~\eqref{hyper poly} is nonzero, we see here another explicit example of a bad approximability by rational functions. 

Now with $\nu \in \Q$, a standard game of Diophantine approximation, both in ineffective (the original
and simplest proof of Thue's theorem for the special case of $r$-th roots from rational numbers)
and, in favorable rare circumstances, effective works (Thue, Siegel, Baker, and Gregory Chudnovsky~\cite{ChudnovskyThueSiegel}), is to take the diagonal $[n/n]$ of the Pad\'e table, specialize $x$ 
to some rational number $\xi \in (0,1) \cap \Q$, and exploit the ensuing small linear forms whose generating function
\begin{equation}
\begin{aligned} \label{Pade gen}
& \qquad \sum_{n=0}^{\infty} \big(  B_{n,n}(\xi) - (1-\xi)^{\nu} A_{n,n}(\xi)  \big) \, z^n \\ 
& = \sum_{n=0}^{\infty} \left(   \pFq{2}{1}{-\nu-n,,-n}{-2n}{\xi}  - (1-\xi)^{\nu} \pFq{2}{1}{\nu-n,,-n}{-2n}{\xi}     \right) \, z^n  \\
 & \in \bigoplus_{n=0}^{\infty} \frac{z^n}{(\den(\nu)\den(\xi))^{2n} \binom{2n}{n}} \, \Z  + (1-\xi)^{\nu} \bigoplus_{n=0}^{\infty} \frac{z^n}{(\den(\nu)\den(\xi))^{2n} \binom{2n}{n}} 
\, \Z
\end{aligned}
\end{equation}
is holonomic on $\C \setminus \left\{ \left(  \frac{ 1 \pm \sqrt{1-\xi}}{2} \right)^{-2} \right\}$  and overconvergent 
at the smaller of these singularities $\left( \frac{1 + \sqrt{1-\xi}}{2} \right)^{-2}$; so the convergence disc of~\eqref{Pade gen} is 
$|z| < \left( \frac{1 - \sqrt{1-\xi}}{2} \right)^{-2}$, the distance to the next singularity.  
Baker~\cite{BakerCube} famously used the $\{2,3,\infty\}$-adic properties of~\eqref{Pade gen} with the choice $\nu := -1/3$  and $\xi := 3/128$   --- 
so $(1-\xi)^{\nu} = (8/5) \sqrt[3]{2}$, and the $3$-adic convergence radius is $1/\sqrt{3}$, thanks to $|\xi|_3  = |3/128|_3 = 1/3$, rather than the ``generic'' $1/3\sqrt{3}$,  --- 
to derive the explicit sub-Liouville inequality $|\sqrt[3]{2} - p/q| > 10^{-6} q^{-2.995}$. (Baker's analysis is synthesized by~\cite[
Theorem~3.5]{ChudnovskyThueSiegel}, following which Chudnovsky sets forth to compute the exact denominators asymptotic to refine the crude~$4^n$ from $\binom{2n}{n}$, and
thus improve Baker's effective irrationality measure to~$2.43$; see also~\cite{ChudnovskyHermite}.)

For the general cubic (or higher) root $\sqrt[r]{a/b}$, this analysis stands no
chance for a sub-Liouville effective irrationality measure (unless $b$ is much bigger than~$a$). But Thue~\cite[\S~9]{ThueSelected}, in his groundbreaking paper 
\emph{Bemerkungen \"uber gewisse N\"aherungsbr\"uche algebraischer Zahlen} written in~1907, proved the ineffective irrationality measure $1+r/2 + \epsilon$ by --- in effect --- observing that one excellent rational approximant $p/q \approx   \sqrt[r]{a/b} \in \Q^{\times} \cap (0,1)$ yields the infinite set of fair rational approximants 
$$
\frac{pB_{n,n}\left(1-aq^r/bp^r \right)}{qA_{n,n}\left(1-aq^r / bp^r \right)} \approx \sqrt[r]{\frac{a}{b}},
$$
and that these form a fairly dense net of fair approximants, thus precluding --- by the gap principle --- the existence of  a second excellent $p'/q' \approx \sqrt[r]{a/b}$. For that Thue used the  
 $x := 1-aq^r/bp^r$ specialization in the polynomial identity (cf.~\cite[Hilfssatz~2]{SiegelThue}, or~\cite[Lemma~2]{AlladiRobinson} for an axiomatization)
\begin{equation} \label{Thue det}
A_{n,n}(x) B_{n+1,n+1}(x) - A_{n+1,n+1}(x) B_{n,n}(x) =  (-1)^{n-1} \frac{(n!)^2}{(2n)!(2n+1)!} \, x^{2n+1},
\end{equation}
with the nonvanishing determinant proving at once the requisite non-equality 
$$
\frac{pB_{n+1,n+1}\left(1-aq^r/bp^r \right)}{qA_{n+1,n+1}\left(1-aq^r / bp^r \right)} \neq  \frac{pB_{n,n}\left(1-aq^r/bp^r \right)}{qA_{n,n}\left(1-aq^r / bp^r \right)}, \qquad \textrm{for all }    n = 0, 1, 2, \ldots. 
$$ 
For arbitrary algebraic targets~$\alpha \in \Qbar$ (other than $r$-th roots or cubic irrationalities), where Thue could not 
rely on the explicit Hermite--Pad\'e approximants to the binomial functions $(1-x)^{\nu}$, he 
instead employed~\cite[\S~11]{ThueSelected} the Dirichlet box principle, in a flash of insight in the~1908 paper  \emph{Om en generel i store hele tal ul{\o}sbar ligning},
to derive the existence of similar  (but vaguer) polynomial identities. His key discovery was that the inexplicit polynomial identities found 
nonconstructively by the Dirichlet box principle worked, {\it grosso modo}, in essentially the same way as in the explicit special case of~\eqref{hyper poly} and~\eqref{Thue det}.
In particular, Thue used a Wronskian determinant to replace the explicit determinant~\eqref{Thue det}, now evaluating
to some
nonzero degree-$2n+1$   polynomial of the form
$x^{(2-\eta)n}V(x)$, with small coefficients, for a suitably small parameter $\eta > 0$. As $\deg{V} \leq \eta n + 1$  ---  or alternatively, as Thue argued, since the coefficients of $V$ are small,  ---  the polynomial $V(x)$ has forcibly a low order of vanishing at the point $x  = 1-aq^r/bp^r$. 
Then Thue
runs the construction after taking the corresponding derivative of his auxiliary polynomials. (See also Zannier~\cite[\S~2]{ZannierLectureNotes} or Masser~\cite[\S~12]{MasserBook} for a detailed treatment and a discussion of 
nuances.)

Shidlovsky's lemma is a different generalization of Thue's Wronskian argument, whose proofs can still roughly be summarized by a (higher rank) determinantal identity akin to~\eqref{Thue det} (of which the latter is strictly speaking a particular and representative case) remaining ``almost in the monomial form.''  It includes Theorems~\ref{Shidlovsky} and~\ref{Chudnovsky Fuchsian form}, and their multiple variations such as~\cite[\S~3]{BombieriG} from the proof of Bombieri's $G$-functions theorem that we discuss in~\S~\ref{sec:G-function arithmetic}, and~\cite[Prop.~VIII.2.3]{Dwork},~\cite[Theorem~3.1, Lemma~8.3]{ChudnovskyG} from the proof of the Chudnovskys's fundamental theorem that we mentioned in Remark~\ref{remark on global nilpotence}. 

\subsubsection{The Hermite--Pad\'e approximants to $\log(1-x)$}  \label{log discussion} We have~\cite[\S~8]{Jacobi} (see also Feldman--Nesterenko~\cite[ch.~2, \S~3.2]{FeldmanNesterenko},  Jager~\cite{Jager}, and Chu~\cite{Chu} for various generalizations)
\begin{equation}
\begin{aligned}
\label{hermitepadelog}
& 2 \sum_{k=0}^{n} \binom{n}{k}^2 \left( H_{n-k}  - H_k \right) (1-x)^k + \log(1-x) \sum_{k=0}^{n} \binom{n}{k}^2 (1-x)^k  \\  
& =
x^{2n+1} \, \int_0^1 \frac{ t^n (t-1)^{n}}{(tx-1)^{n+1}} \, dt = - \frac{x^{2n+1}}{(2n+1) \binom{2n}{n}} + O(x^{2n+2}),
\end{aligned}
\end{equation}
where $H_r := \sum_{k=1}^{r} 1/k$ are the harmonic numbers.

 \begin{remark}
 In terms of the general Meijer $G$ function
 \begin{equation*}
 G_{p,q}^{m,n}  \left(    \begin{smallmatrix}   a_1 &  \cdots & a_p  \\ b_1 & \cdots & b_q \end{smallmatrix}   \big\vert \, z  \right)  = 
  \int_{\Re(s) = \sigma} \frac{\prod_{j=1}^{m} \Gamma(b_j-s)  \prod_{j=1}^{n} \Gamma(1-a_j+s)}{ \prod_{j=m+1}^{q} \Gamma(1-b_j+s) \prod_{j=n+1}^{p} \Gamma(a_j-s) }
  \, z^s \, \frac{ds}{2\pi i},
 \end{equation*}
  the remainder term in~\eqref{hermitepadelog} can be expressed also as $ - (n!)^2 \, G_{2,2}^{2,0}  \left(    \begin{smallmatrix}   n+1 & n+1 \\ 0 & 0 \end{smallmatrix}  \big\vert \, 1-x \right)$. This general definition as a Barnes integral is valid under the assumption that all poles of all $\Gamma(b_j-s)$ are on the right of the integration line $\Re(s) = \sigma$, while all poles of all $\Gamma(1-a_j+s)$ are on the left of that line.~\endofremark
\end{remark}

Here,
\begin{equation}
\begin{aligned}
 \label{Afun}
 \sum_{k=0}^{n} \binom{n}{k}^2 (1-x)^k & = 
  \pFq{2}{1}{-n,-n}{1}{1-x}  = x^n P_n\left(  \frac{2-x}{x} \right),  \\ 
  P_n(x) &  := \frac{1}{2^n \, n!} \left( \frac{d}{dx} \right)^n \, (x^2-1)^n
\end{aligned}
\end{equation}
in terms of the Legendre polynomials $P_n(x)$: the complete orthogonal system on $[-1,1]$ under the Lebesgue measure and the normalization $P_n(1) = 1$. Their 
generating series 
\begin{equation} \label{LegendreP}
\frac{1}{\sqrt{ 1 - 2xz + z^2  }} = \sum_{n=0}^{\infty} P_n(x) z^n
\end{equation}
is precisely the function whose integrality properties  ---  namely: that $P_n$ is integer-valued on the odd integers,
amounting to the $\Z[x]$ polynomials $(2x)^nP_n(1/x)$ in~\S~\ref{ARB} below  ---  we exploit in~\S~\ref{sec:logs}. 

If like in~\S~\ref{binomials discussion} we multiply~\eqref{hermitepadelog} by $z^n$ and sum the generating series over $n \in \NwithzeroA$, the resulting
$$
\Q[x] \llbracket z \rrbracket + \log(1-x) \, \Z[x]\llbracket  z \rrbracket
$$
function  is holonomic in~$z$ and has its $\Z[x]\llbracket  z \rrbracket$ and $\Q[x] \llbracket z \rrbracket $ components satisfy the 
homogeneous and inhomogeneous first-order ODEs
\begin{equation} \label{HermitePadelogODE}
\begin{aligned}
(-1 + 4z - 2xz - x^2z^2) Y'(z) + (2 - x - x^2 z) Y(z)  & = 0 \\ 
\textrm{and} \qquad & \\
(-1 + 4z - 2xz - x^2z^2) Y'(z) + (2 - x - x^2z) Y(z) &  = -x,
\end{aligned}
\end{equation}
respectively. These are holonomic functions on
$$
\C \setminus \left\{  p_-(x), p_+(x) \right\},   \qquad p_{\pm}(x) := \left( \frac{1 \pm \sqrt{1-x}}{x} \right)^2,
$$
where these singularities can be also directly obtained from~\eqref{Afun} and~\eqref{LegendreP}. 

Specializing $x= 1/n$ and $y = 1/m$, we have
\begin{equation*}
p_-(1/n) p_+(1/m)/|mn| = 1 + o_{|m/n| \to 1}(1). 
\end{equation*}
This asymptotic is related to the analyticity mechanism with Hadamard products in~\S~\ref{sec:Hadamard}, and could be also
used there as an alternative, but ultimately equivalent given~\S~\ref{ARB} just below, proof of Theorem~\ref{logsmain}.

\subsubsection{The Hermite--Pad\'e approximants to $\log\left( \frac{1-x}{1+x} \right)$}  \label{ARB} The change of variables $x \mapsto 2x/(1+x)$ in~\eqref{hermitepadelog} 
rewrites the formula thus:
\begin{equation*}
\begin{aligned}
& 2 \sum_{k=0}^{n} \binom{n}{k}^2 \left( H_{n-k}  - H_k \right) (1-x)^k(1+x)^{n-k} + \log\left(\frac{1-x}{1+x}\right) (2x)^n P_n(1/x) \\ 
& =
(x+x^2)^{2n+1}\, \int_0^1 \frac{ t^n (t-1)^{n}}{(2tx-1-x)^{n+1}} \, dt = - \frac{x^{2n+1}}{(2n+1) \binom{2n}{n}} + O(x^{2n+2}).
\end{aligned}
\end{equation*}
In~\S~\ref{sec:logs} we use the generating functions of these formulas specialized to $x := 1/a$ with $a$ a large odd integer.

\section{Concentration of measure} \label{concentration of measure}

If we randomly and independently sample  a very large number~$n-1$ of uniformly distributed points of the segment~$[0,1]$, the~$n$ spacings that remain
will be almost surely close to some ordering of the set $\left\{ \log(n/j)/n \, : \, 1 \leq j \leq n \right\}$, while the $n-1$
sample points themselves will be almost surely close to some ordering of the set $\{ j/n \, : \, 1 \leq j < n \}$. These facts are  simplest expressions of 
the Law of Large Numbers in statistics, with the precise quantitative decay rates being captured by
the \emph{concentration of measure phenomenon} of Dvoretzky and Milman~\cite{MilmanSchechtman,MilmanDvoretzky,Ledoux} for the high-dimensional $\ell^r$-ball,  
in the respective cases $r=1$ and $r = \infty$. A popular expression of the measure concentration principle, due to Gromov~\cite[\S~$3\frac{1}{2}.20$]{Gromov}, is to say that the \emph{observable
diameter} of the unit volume $\ell^r$-ball in the asymptotic of high dimension~$n$ is on the order of only $\frac{1}{\sqrt{n}} = o(1)$, in contrast to its
diameter as a metric space which is on the order of~$\sqrt{n}$. It is the observable and not the metric properties that are relevant to the various auxiliary polynomial constructions
 undertaken in Diophantine approximation. 
 
 These specific distributions (and the finer statistics) are best expressed by the fact~\cite[Theorem~1]{StochasticLp} (going back\footnote{The classic theorem relating the normal distribution to the Euclidean sphere is popularly ascribed to Poincar\'e in~1912, but see~\cite[\S~6]{DiaconisFreedman} for a scrupulous historical research, 
and a discussion of a broader context.} to \'Emile Borel~\cite[\S~V]{Borel} for the $r=2$ case of the Euclidean ball; see also~\cite[Lemma~1]{SchechtmanZinn} or~\cite[\S~3]{ApproximateIndependence}) that the normalized volume measure of the $n$-dimensional $\ell^r$ ball is generated stochastically by 
the random vector
$$
\left( \frac{X_1}{\left( |X_1|^r +\ldots + |X_n|^r + Z  \right)^{1/r}}, \ldots,  \frac{X_n}{\left( |X_1|^r +\ldots + |X_n|^r + Z  \right)^{1/r}}  \right),
$$
where $X_1, \ldots, X_n$ are independent and identically distributed random variables with probability density function 
$\frac{1}{2\Gamma(1+1/r)} e^{-|t|^r}$, and $Z$ is a jointly independent random variable with the exponential density function $e^{-t} \cdot \chi_{[0,\infty)}(t)$. Moreover,
 the concentration function is Gaussian. These features are general, 
while for our purposes here, only the simplest statement with the $\ell^{\infty}$-ball is used. In this $r \to \infty$ limiting case, one additional (but only technical) simplification is that the random vector components in the stochastic
generation of $\mu_{[-1,1]^n}$ are independent rather than merely asymptotically independent.

In the multivariable auxiliary Diophantine constructions~\S\S~\ref{fine section} and~\ref{slopes}, we will use measure concentration ideas as described in~\S~\ref{high dims}. One aspect of this is to constrict the component sets $\{k_j\}$ of the high-dimensional exponent vectors~$\bk$ in all the monomials~$\bx^\bk$ occurring in the 
auxiliary polynomial constructions from the evaluation modules that we introduced in~\S~\ref{filtrations} and studied in~\S~\ref{evaluation}.
 This type of application is among the most standard in Diophantine approximation, after the classic works of Roth~\cite[\S~6.3.5]{BombieriGubler} and Schmidt\footnote{
In the proof of Schmidt's Subspace theorem, the ``$d+1$'' exponent for the bad approximability in projective space $\P^d$ receives a probabilistic explanation as the reciprocal 
of the equal expectations of the individual coordinates of a point~$\xi$ taken at random from the surface boundary of the $d+1$-dimensional
standard simplex. The concentration property, used for the parameter count at the auxiliary polynomial construction in the Thue--Siegel lemma, states precisely that all the column sums in a tall $n \times (d+1)$ matrix made of  $n \to \infty$ such independent and identically distributed
random rows~$\xi$ converge in probability to the expectation $n/(d+1)$, at an asymptotic rate exponential in~$-n$. 
 See~\cite[\S~3, \emph{Example~1}]{FaltingsWustholz} for a broader context of Harder--Narasimhan filtrations on graded algebras of auxiliary functions. The work of Faltings and W\"ustholz made a deeper use of probability measures which, in combination with the Faltings product theorem for directly treating the nonvanishing of the auxiliary construction at the special point, ultimately eliminates the difficult geometry of numbers part from Schmidt's proof.}~\cite[\S~7.5.15]{BombieriGubler}, and especially Wirsing~\cite[\S~4.2]{Wirsing}  (see also~\cite[Theorem~7.2.1]{Stolarsky} for two alternative and more detailed treatments of the relevant material from Wirsing's argument). 
 Methodologically our high-dimensional Diophantine analysis in~\S~\ref{fine section} is rather similar to the multivariable auxiliary polynomial constructions that Wirsing used in the proof of his theorem on the bad approximability of a fixed algebraic target by algebraic approximants of a given degree. 
 
\subsection{The Erd\"os--Tur\'an bound}

Recall the definition of the box discrepancy function on the hypercube; cf~\cite[\S~2.5.3]{UDC}.

\begin{df}  \label{box discrepancy}
The (normalized, box) {\bf discrepancy function} $D : [0,1]^n  \to (0,1]$ is the supremum over all closed intervals
$I  = [a,b] \subset [0,1]$ of the defect between the length $\mu_{\mathrm{Lebesgue}}(I) = b-a$ of $I$ and the proportion of points falling inside $I$: 
$$
D(t_1, \ldots, t_n) := \sup_{I \subset [0,1) } \left| \mu_{\mathrm{Lebesgue}}(I) - \frac{1}{n}  \# \{ i \, : \, t_i \in I \}  \right|.
$$
\end{df}
With the identification $[0,1)^n   \iso \T^n$ induced from $e(t) := \exp(2\pi i t)$, harmonic analysis on the 
circle supplies a basic way to upper-bound the discrepancy function. The  Erd\"os--Tur\'an  inequality states~\cite[Theorem~1.14]{DrmotaTichy}
\begin{equation} \label{Erdos-Turan}
D(t_1, \ldots, t_n) \leq 3 \left(  \frac{1}{K+1} + \sum_{k=1}^K \frac{1}{k} \left| \frac{e(kt_1) +\ldots + e(kt_n)}{n} \right| \right) \quad \forall K \in \NwithzeroA,
\end{equation}
 in terms of the character sums on the group~$\T$.

\subsection{The large deviations bound}
The following estimate will be critical. 

\begin{thm} \label{thm_MeasureConcentration}
There exist two  absolute constants $C, c \in \R$ such that
for any $\varepsilon > 0$ and any $n \in \NwithzeroA$, 
the set
$$
B_{\varepsilon}^n :=
\left\{ \mathbf{t} \in [0,1]^n \,  : \, D(\mathbf{t}) \geq \varepsilon \right\}
$$
has $n$-dimensional Lebesgue measure smaller than $C e^{-  c \varepsilon^4 n}$. 
\end{thm}

For instance, the proof will show that we can take $c = 1/300$ and $C = 100$ in this theorem. 

\begin{remark} The existence of an exponential (in the negative of the dimension) asymptotic decay rate 
is a hallmark of basic concepts of entropy in the theory of large deviations~\cite{Ellis}. The specific rate estimate worked out in Theorem~\ref{thm_MeasureConcentration}
is of no consequence for our purposes, but its existence is used crucially in~\S\S~\ref{fine section},~\ref{slopes}. An alternative path
to Theorem~\ref{thm_MeasureConcentration}, not using harmonic analysis and the Erd\"os--Tur\'an bound (and with different, indeed better numerical constants~$c,C$), but instead taking for base the rudimentary Chebyshev estimate~\cite[Lemma~12]{Wirsing}, can be derived from 
the bound $\leq e^{- \pi r^2}$ on the concentration function~\cite[Prop.~2.8]{Ledoux} for the
uniform measure on~$[0,1]^n$. In Ledoux's book, the latter concentration inequality is obtained as a consequence~\cite[Cor.~2.6]{Ledoux} under a contraction of the sharp estimate
$\leq e^{-r^2/2}$ for the concentration function of the canonical Gaussian measure on $\R^n$. The latter, in turn, 
is traditionally a consequence of L\'evy's isoperimetric inequality on the Euclidean sphere~\cite[Theorem~2.3]{Ledoux}.~\endofremark
\end{remark}

Our proof of Theorem~\ref{thm_MeasureConcentration} will be based on the most standard form of Hoeffding's concentration inequality~\cite{Hoeffding}.
 For the sum of independent random variables $X_1, \ldots, X_n$ 
 taking values in the interval $[-1,1]$, Hoeffding's inequality~\cite[Theorem~2.8]{BoucheronLugosiMassart} bounds the large deviation tail probability
exponentially by
\begin{equation} \label{Hoeffding bound}
\P\left( \Big| X_1 +\ldots + X_n -  \mathbf{E} [X_1 +\ldots + X_n] \Big| \geq \varepsilon n \right)  \leq 2 e^{-  \varepsilon^2 n /2}. 
\end{equation}
On changing $\varepsilon$ to $\varepsilon/2$ and using the triangle inequality and the subadditivity of probability,
we can apply this to the real and imaginary parts of $\T$-valued independent random variables $Z_1, \ldots, Z_n$ to
get the following variant: 
\begin{lemma}[Hoeffding] 
The sum of independent random variables $Z_1, \ldots, Z_n$ taking values in the complex unit circle~$\T$ has the tail
probability large deviations bound
\begin{equation} \label{Hoeffding T bound}
\P\left( \Big| Z_1 +\ldots + Z_n -  \mathbf{E} [Z_1 +\ldots + Z_n] \Big| \geq \varepsilon n \right)  \leq 4 e^{-  \varepsilon^2 n /8}. 
\end{equation}
\end{lemma}
\begin{proof}[Proof of Theorem~\ref{thm_MeasureConcentration}]
In combination with the Erd\"os--Tur\'an bound, we derive a proof of the theorem, with the following explicit estimate.
Take $Z_1, \ldots, Z_n$ to be independent and uniformly distributed points of the circle~$\T$. Then $\mathbf{E}[Z_1^k 
+\ldots + Z_n^k] = 0$ for all $k = 1, 2, \ldots$, giving uniformly by Hoeffding's bound~\eqref{Hoeffding T bound}
\begin{equation}
\P\left( \frac{1}{k} \left| \frac{Z_1^k +\ldots + Z_n^k}{n} \right| \geq \varepsilon \right)  \leq 4 e^{-  \varepsilon^2 k^2 n /8}.
\end{equation}
It follows that for every $K \in \NwithzeroA$ and $\varepsilon > 0$ the probability 
\begin{equation*}
\P \left(
\frac{3}{K+1} +3 \sum_{k=1}^K \frac{1}{k} \left| \frac{Z_1^k +\ldots + Z_n^k}{n} \right| 
\geq \frac{3}{K+1} +3 K \varepsilon \right)  
\leq 4K e^{-\varepsilon^2 n /8}. 
\end{equation*}
If we firstly change $\varepsilon$ to $(\varepsilon/6)^2$ and then select $K := \lfloor 6/\varepsilon \rfloor$, we derive
\begin{equation} \label{eq:et}
\inf_{K \in \NwithzeroA} \left\{
\P \left(
\frac{3}{K+1} +3 \sum_{k=1}^K \frac{1}{k} \left| \frac{Z_1^k +\ldots + Z_n^k}{n} \right| 
\geq \varepsilon \right)   \right\}
\leq \min\left( 1, ( 24/\varepsilon) e^{-\varepsilon^4 n /288} \right). 
\end{equation}
By the Erd\"os--Tur\'an inequality~\eqref{Erdos-Turan}, the left-hand side of~\eqref{eq:et} is an upper bound on our
requisite $\vol(B_{\varepsilon}^n) =\vol(\{D(\bt) \geq \varepsilon\})$. The right-hand side of~\eqref{eq:et} 
is majorized by $100\,e^{-\varepsilon^4 n/300}$ for all $n \geq 1$ and all $\varepsilon > 0$. 
\end{proof}

\section{The cost of an integration}   \label{integration cost}

The basic idea of our solution~\cite{UDC} of the unbounded denominators conjecture is that
we can get useful holonomy  bounds on certain $\Q(x)$-linear spaces of 
algebraic functions that come from a supposed $\Z \llbracket q \rrbracket$ modular function $f(\tau)$  on a  noncongruence subgroup of $\SL_2(\Z)$,
expanded formally in the modular function $x = \lambda/16$ via the equality of rings $\Z \llbracket q \rrbracket = \Z \llbracket x \rrbracket$. In that setting, the key point was in getting an asymptotically tight holonomy bound which only 
runs into a contradiction upon successively including more and more functions with the transformation  $f(\tau) \rightsquigarrow f(p\tau)$ for a range of primes~$p$, 
and finding that the increase in the dimension of $\Z \llbracket q \rrbracket$ modular functions is more than the increase in the holonomy bound, unless $f(\tau)$ was congruence to begin with.

In our present paper, we have a somewhat analogous scheme where the role of the transformation $f(\tau) \rightsquigarrow f(p\tau)$ is taken up by
an integration $f(x) \rightsquigarrow \int (f(x) - f(0)) \frac{dx}{x}$. Here it is more of a gamble whether or not the increase in the dimension (which we compute in~\S~\ref{sec:lindep} and~\S~\ref{logsjointli} for
our main application to Theorems~\ref{mainA} and~\ref{logsmain}) 
turns out enough of a compensation for the increase in the bound (which comes entirely through the added denominators, and is handled in the present
section by a prime number theorem estimate). We find it astonishing that the integrations gamble succeeds as a crucial ingredient for both of our main applications
in the present paper: Theorem~\ref{mainA} on the $\Q$-linear independence of $1, \zeta(2)$, and $L(2,\chip)$, and Theorem~\ref{logsmain} on the irrationality of
certain products of two logarithms. 

 This section establishes some preparatory Lemmas
 which will be used to compute the added denominator cost for including such integrations into the
setup of Theorem~\ref{basic main}. The upshot will be the integration cost function of Definition~\ref{integrated lcm} of the next section, and our main result
Theorem~\ref{main:elementary form} where this function is used to define an added denominators term $\tau^{\sharp}$ to the $\tau(\mathbf{b})$ of Theorem~\ref{basic main}. 

The following lemma is a direct consequence of the prime number theorem.
\begin{lemma}\label{lcm lemma} \leavevmode
\begin{enumerate}
\item \label{lcmone}
If $k  \sim \gamma n$ for a fixed $\gamma \in (0,1]$, the lowest common multiple $L_{n,k}$ of the consecutive integers $n-k, n-k+1, \ldots, n$ is asymptotic under $n \to \infty$ to
$$
\exp \left(  \left(\sum_{h=1}^{\lfloor 1/\gamma \rfloor - 1} \frac{1}{h}\right) k + \frac{n}{\lfloor 1/\gamma \rfloor} + o(n)  \right).
$$
This bound is uniform for all $\gamma \geq \gamma_0$, where $\gamma_0>0$ is a fixed constant, in the following sense, for any $\epsilon>0$, there exists $N=N(\gamma_0, \epsilon)$ such that for all $n\geq N$ and $\gamma \geq \gamma_0$, the error term is at most $\epsilon n$.
\item \label{lcmtwo}
If $k=o(n)$, then the lowest common multiple $L_{n,k}$ of the consecutive integers $n-k, n-k+1, \ldots, n$ is  $\exp(o(n))$. Moreover, as $n\rightarrow \infty$, we have for all $0\leq k \leq n$, 
$$
L_{n,k} \leq \exp \left(  \left(\sum_{h=1}^{\lfloor 1/\gamma \rfloor - 1} \frac{1}{h}\right) k + \frac{n}{\lfloor 1/\gamma \rfloor} + o(n)  \right),
$$
where $\gamma=k/n$ and and if $k=0, \gamma=0$, the above formula is to be interpreted as $\exp( o(n))$. The error term $o(n)$ in this upper bound is uniform. 
\end{enumerate}
\end{lemma}

\begin{basicremark}  \label{remark:zudilingraph} If~$k \ge n/2$, then~$[1,2,\ldots,n] = [(n-k),\ldots,n]$,
because if~$p \le n$ then some multiple of~$p$ lies in~$[n/2,n]$. Hence~$L_{n,k}$
does not depend on~$k$ within this range. With~$\gamma = k/n$ and~$n \rightarrow \infty$,
the exponent in this inequality in terms of~$\gamma$ as a multiple of~$n$ is given in Figure~\ref{fig:zudilin}.
\addtocounter{subsubsection}{1}
 \begin{figure}[!h]  
\begin{center}
  \includegraphics[width=70mm]{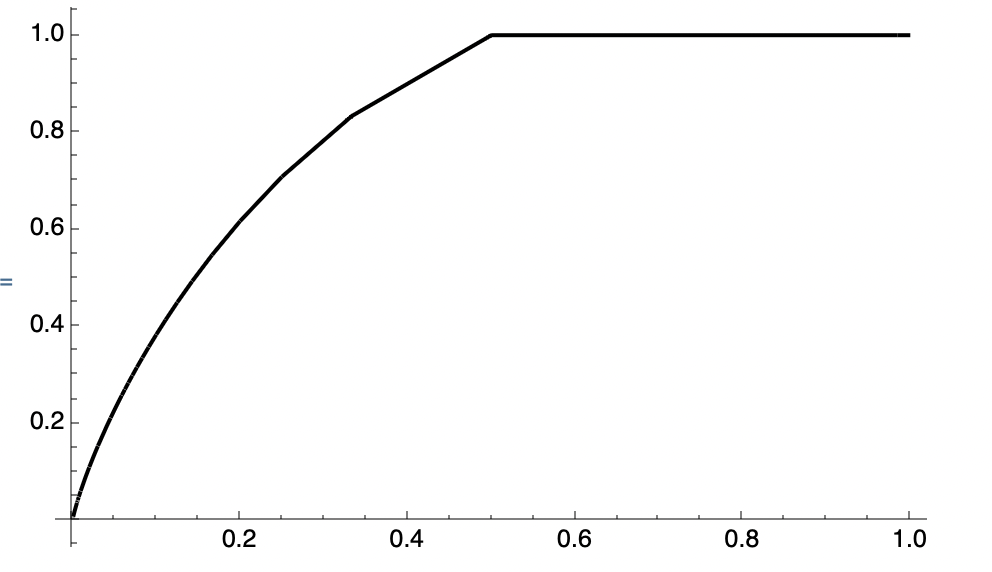}
\end{center}
\caption{The bound
for $\log(L_{n,k})/n$ as a function of~$\gamma = k/n$
with~$n \rightarrow \infty$.}
\label{fig:zudilin}
\end{figure}
\endofremark
\end{basicremark}

\begin{proof}[Proof of Lemma~\ref{lcm lemma}]
We begin with part~(\ref{lcmone}).
By the prime number theorem, the main term in $[1,\ldots,n]$ (after taking log) is given by $\sum_{p\leq n} \log p$ (i.e., we may just count primes without counting multiplicities). The error term here is independent of $\gamma$. Thus the exponential asymptotic rate of $[n-k,\ldots,n]$ is given by counting how many of the primes $p \leq n$ divide at least one among $n-k,\ldots, n$. The only primes $p\leq n$ 
not occurring in this count are those that admit an $a \in \NwithoutzeroA$ such that $ap < n-k$ and $(a+1)p > n$. Given such an $a$, the primes~$p$ in question are exactly the primes from the interval $(n/(a+1), (n-k)/a)$. This is a non-empty interval if and only if $a+1 < 1/\gamma$; in which case its length equals~$n((1-\gamma)/a - 1/(a+1))$.  Hence, $\log{[n-k,\ldots,n]}$ amounts to
 \[
 \begin{aligned}
 n\left(1-\sum_{a=1}^{\lfloor 1/\gamma \rfloor -1} ((1-\gamma)/a - 1/(a+1))\right)+o(n) \\
 =n\left(\gamma \left(\sum_{a=1}^{\lfloor 1/\gamma \rfloor -1} 1/a\right)+1/(\lfloor 1/\gamma \rfloor)\right)+o(n).
 \end{aligned}\]
Note that from our assumption $\gamma\geq \gamma_0$, the above sum is a finite sum with a uniform upper bound on its number of terms. Also, the error term from the prime number theorem is uniformly controlled as it is only being applied to intervals whose lengths and endpoints are controlled uniformly in~$n$. 

We now  consider case~(\ref{lcmtwo})
of Lemma~\ref{lcm lemma}.
The precise formulation of the first assertion is that for any $\epsilon>0$, there exist $N=N(\epsilon)$ and $\delta = \delta(\epsilon)$ such that for all $n > N$ and all $k < \delta n$, we have $L_{n,k} < \epsilon n$. This is a consequence of (1). More precisely, for $\delta < 1/2 $ by definition, for all $k < \delta n$, we have
\[
\begin{aligned}
L_{n,k} \leq L_{n, \delta n}  & =\exp \left(  \left(\sum_{h=1}^{\lfloor 1/\delta \rfloor - 1} \frac{1}{h}\right) \delta n + \frac{n}{\lfloor 1/\delta \rfloor} + o_\delta (n)  \right) \\
& \leq \exp \left(\left(\delta (1 + \log (1/\delta)) + (1/\delta -1)^{-1} \right)n + o_\delta (n) \right)
\end{aligned} .\]
Note that $\lim_{\delta \rightarrow 0} \delta (1 + \log (1/\delta)) + (1/\delta -1)^{-1} =0$; we pick a $\delta=\delta(\epsilon)$ with $ \delta (1 + \log (1/\delta)) + (1/\delta -1)^{-1} < \epsilon/2$. For this $\delta$, by (1), there exists $N=N(\delta, \epsilon)=N(\epsilon)$ such that for $n>N$, the error term $o_\delta (n)  < (\epsilon/2)n$. Then for all $n>N$ and $k<\delta n$, we have $L_{n,k} < \epsilon n$.

The precise formulation of the second assertion is that for any $\epsilon>0$, there exists $N=N(\epsilon)$ such that for all $n>N$ and all $0\leq k \leq n$, we have
\[\log L_{n,k} \leq \left(\sum_{h=1}^{\lfloor 1/\gamma \rfloor - 1} \frac{1}{h}\right) k + \frac{n}{\lfloor 1/\gamma \rfloor} + \epsilon n.\]
Note that from the proof of the first assertion above, there exists $\delta=\delta(\epsilon)$ such that the above inequality holds for all $n> N_1(\epsilon)$ and $k < \delta n$. Moreover, from part (1) with $\gamma_0=\delta$, the above inequality holds for all $n> N_2(\delta, \epsilon)=N_2(\epsilon)$. Thus the desired bound holds for all $n > \max\{ N_1(\epsilon), N_2(\epsilon)\}$.
\qedhere
\end{proof}

The following is a variant of the above lemma.

\begin{lemma} \label{lcm-2} \leavevmode
\begin{enumerate}
\item \label{lcmoneagain}
Fix $\gamma_0\leq \gamma_0' \in (0,1)$. For $k, l \leq n$ such that $\gamma_0 \leq k/n$ and $\gamma_0 \leq l/n\leq \gamma_0'$, the logarithm of the product $L_{n,k}^{\geq l}$ of the primes $p > l$ that have some multiple in the interval $[n-k,n]$ is asymptotic  ---  as $n \to \infty$ uniformly with respect to $k,l$  ---  to
$$
\left(
k \sum_{h=1}^{\lfloor (n-k)/\max(k,l) \rfloor} 1/h   \right) +\left(\frac{n}{\lfloor (n + (l-k)^+ )/\max(k,l) \rfloor} - l \right)^+  + o(n),
$$
where $\alpha^+ := \max(0,\alpha)$ and the convention being that $\sum_{h=a}^{b}$ is over all integers in the range $a \leq h \leq b$, and the empty sum is zero. 
\item \label{lcmtwoagain}
Moreover, as $n\rightarrow \infty$, for all $0\leq k, l \leq n$, we have
\[\log L_{n,k}^{\geq l} \leq \left(
k \sum_{h=1}^{\lfloor (n-k)/\max(k,l) \rfloor} 1/h   \right) +\left(\frac{n}{\lfloor (n + (l-k)^+ )/\max(k,l) \rfloor} - l \right)^+  + o(n),\]
where the error term is uniform with respect to all $k,l$. 

(If $k=l=0$, the right-hand side of the above bound is to be interpreted as~$o(n)$.)
\end{enumerate}
\end{lemma}

\begin{proof}
We begin with part~(\ref{lcmoneagain}).
As in the proof of Lemma~\ref{lcm lemma},
the  primes $p\leq n$ that do not have any multiples among $n-k,\ldots,n$ are the ones that lie in $\displaystyle \cup_{a=1}^{\lfloor n/k \rfloor -1} \left(n/(a+1), (n-k)/a\right)$. 
The new assumption here that $p>l$ implies that $a < (n-k)/l$, and hence that
\[
\begin{aligned}
& \left(\cup_{a=1}^{\lfloor (n-k)/k \rfloor} (n/(a+1), (n-k)/a)\right) \cap (l,n] \\
& =  \left(\cup_{a=1}^{h_0 -1} (n/(a+1), (n-k)/a)\right) \cup  (\max(n/(h_0+1),l), (n-k)/h_0),
\end{aligned}
\]
where $h_0= \lfloor (n-k)/\max(k,l) \rfloor$.

By the prime number theorem, our asymptotic is given by
  \begin{equation*}
  \begin{aligned}
  n - l - \left(\sum_{a=1}^{h_0} ((n-k)/a - n/(a+1)) - (\max(n/(h_0+1),l) - n/(h_0+1))\right)  +o(n) 
  \\ 
  =  k \left(\sum_{a=1}^{h_0} 1/a\right)+\max(n/(h_0+1),l)- l +o(n)
  =  k \left(\sum_{a=1}^{h_0} 1/a\right)+(n/(h_0+1)- l)^+ + o(n)
  \end{aligned}
  \end{equation*}
  and $\displaystyle n/(h_0+1)=\frac{n}{ \lfloor (n-k)/\max(k,l) \rfloor +1}= \frac{n}{\lfloor (n + (l-k)^+ )/\max(k,l) \rfloor} $.

The precise formulation of the second assertion~(\ref{lcmtwoagain}) is that for any $\epsilon>0$, there exists $N=N(\epsilon)$ such that for all $n>N$ and all $0\leq k,l \leq n$, we have
\[\log L_{n,k}^{\geq l} \leq \left(
k \sum_{h=1}^{\lfloor (n-k)/\max(k,l) \rfloor} 1/h   \right) +\left(\frac{n}{\lfloor (n + (l-k)^+ )/\max(k,l) \rfloor} - l \right)^+  + \epsilon n.\]
By Lemma~\ref{lcm lemma} (2), there exists $\delta=\delta(\epsilon)$ such that for all $n > N_1(\epsilon)$ and all $k\leq \delta n$, we have
\[\log L_{n,k}^{\geq l}  \leq \log L_{n,k} < \epsilon n.\]
Therefore, we now assume $k > \delta n$. For $l < \min\{\delta, \epsilon/2\} n < k$, by Lemma~\ref{lcm lemma}(1), we have that for $n> N_2(\delta, \epsilon/2) = N_2(\epsilon)$, 
\[
\begin{aligned}
\log L_{n,k}^{\geq l}  \leq \log L_{n,k} &  \leq \left(
k \sum_{h=1}^{\lfloor (n-k)/k \rfloor} 1/h   \right) +\frac{n}{\lfloor n  /k \rfloor} + (\epsilon/2)n \\
&  \leq \left(
k \sum_{h=1}^{\lfloor (n-k)/k \rfloor} 1/h   \right) +\left(\frac{n}{\lfloor n /k \rfloor} - l \right)^+  + \epsilon n,
\end{aligned}\]
which is the desired bound as $\max(k,l)=k$ in this case.
For $l> (1-\delta/2) n$, by definition, we have $\log L_{n,k}^{\geq l} \leq \log L_{n,k}^{\geq (1-\delta/2)n}$. Applying (1) to $\gamma_0=\delta, \gamma_0'= 1-\delta/2$, we have that there exists $N_3=N_3(\delta)=N_3(\epsilon)$ such that for all $n>N_3$, we have
 \[
 \begin{aligned}
&  \log L_{n,k}^{\geq (1-\delta/2)n}   \leq 
\left(
k \sum_{h=1}^{\lfloor (n-k)/\max(k,(1-\delta/2)n) \rfloor} 1/h   \right) \\
 & +\left(\frac{n}{\lfloor (n + ((1-\delta/2)n-k)^+ )/\max(k,(1-\delta/2)n) \rfloor} - (1-\delta/2)n \right)^+  + (\epsilon/2) n.
 \end{aligned}\]
Note that the first term is $0$ since $n-k < (1-\delta/2)n$ and the second term $\leq n -  ((1-\delta/2)n \leq (\delta/2)n$. For the above proof, we may shrink $\delta$ to make it $< \epsilon$ and then the above discussion shows that for all $l \geq (1-\delta/2)n$, we have the desired bound
\[\log L_{n,k}^{\geq l}  \leq \log L_{n,k} \leq \epsilon n.\]
Now we only remain to consider $k\geq \delta n$ and $\min\{\delta, \epsilon/2\} n  \leq l < (1-\delta/2)n$ and this case follows from (1).
\qedhere
\end{proof}

\section{The fine holonomy bound}   \label{fine section}

In this section, we arrive at our first main holonomy bound (Theorem~\ref{main:elementary form}), 
which we prove
 by revisiting the method in~\cite[\S~2.5]{UDC}
and enhancing it by the (standard) results of~\S~\ref{functional transcendence} and~\S~\ref{concentration of measure}. This 
elementary treatment of our bound suffices for the proof of Theorems~\ref{mainA} and~\ref{logsmain} and for all our other applications
in this paper. Later, in~\S~\ref{new slopes} and~\S~\ref{slopes}, we will prove other holonomy bounds, some of which involve a Bost--Charles double integral that is theoretically smaller than the rearrangement integral in~\eqref{basic rearrangement form}; however, we will find in Remark~\ref{nearly equal} and \S~\ref{sec_slope for rearrangement}  the difference to be negligibly small in practice. 
For our default treatment we have opted to highlight the increasing rearrangement feature which occurs, under a probabilistic interpretation, 
simultaneously in the top and bottom quantities in the holonomy quotient~\eqref{fine new bound elementary}.  

In order to state our holonomy bound, we first define (following Lemma~\ref{lcm-2}) a
 function
 to  measure up the additional contributions to the coefficient denominators in our multivariable evaluation module 
 under including the \emph{integrals}  of our original set of functions, as
 discussed in~\S~\ref{integration cost}. (In the statement of
 Theorem~\ref{main:elementary form}, these will be the functions~$f_i$ of the form~(\ref{den type int})
 with~$e_i > 0$.)
\begin{df} \label{integrated lcm}
For $0 \leq \max\{ u, 1\} \leq v$ and $w \leq v$, set
\begin{equation*}
\begin{aligned}
I_{u}^{v}(w) & := \int_{\min\{u,1\}}^1 \max\{t-w,0\} \, dt + \int_{\max\{u,1\}}^{v} \left\{   \sum_{h=1}^{\lfloor(t-1)/\max(1,w) \rfloor} 1/h \right\} \, dt  \\
& + 
\int_{\max\{u,1\}}^{v} \max \left\{ \frac{t}{ \lfloor  (t+\max(0,w-1))/\max(1,w) \rfloor} - w, 0 \right\} \, dt. \end{aligned}
\end{equation*}
\end{df}

We now have:

\begin{thm} \label{main:elementary form}
Consider two positive integers~$m,r \in \NwithoutzeroA$, a nonnegative integer vector $\mathbf{e} := (e_1, \ldots, e_m) \in \NwithzeroB^{m}$, and an $m \times r$ rectangular array of nonnegative real numbers 
\[\mathbf{b} := \big( b_{i,j} \big)_{\substack{ 1 \leq i \leq m, \,
1 \leq j \leq r }},\]
 all of whose columns have the form
\begin{equation} \label{column shape}
0 = b_{1,j} =\cdots = b_{u_j,j} < b_{u_j+1,j}= \cdots = b_{m,j}=: b_j,  \qquad \forall j = 1, \ldots, r,
\end{equation}
for some $u_j \in \{0, 1, \ldots, m\}$ depending on the column. 
Let
$$
\sigma_i := b_{i,1} +\ldots + b_{i,r}, \qquad i = 1, \ldots, m
$$
be the $i$-th row sum, and define
\begin{equation}\label{taub}
\tau^{\flat}(\mathbf{b}) := \frac{1}{m^2} \sum_{i=1}^{m} (2i-1) \sigma_i = \sigma_m -\frac{1}{m^2}\sum_{j=1}^r u_j^2 b_j \in [0, \sigma_m]. 
\end{equation}
and, with~\(I_{\xi}^{m}(\xi)\)  as in Definition~\ref{integrated lcm},
\begin{equation}\label{taue}
\tau^{\sharp}(\mathbf{e}) := (2/m^2) \min_{\xi \in [0,m]} \left\{  \xi \sum_{i=1}^{m} e_i 
   +
   \left(  \max_{1 \leq i \leq m}  e_i \right) I_{\xi}^{m}(\xi)
   \right\}.
\end{equation}
Define, finally, 
\begin{equation}   \label{tau sharp}
\tau(\mathbf{b;e}) := \tau^{\flat}(\mathbf{b}) + \tau^{\sharp}(\mathbf{e}).
\end{equation}
Consider a sequence of holomorphic mappings $\varphi_0, \ldots, \varphi_l : (\Db, 0) \to (\C,0)$ with derivatives (\emph{conformal sizes}) satisfying 
\begin{equation} \label{relaxedifholonomic}
|\varphi_0'(0)| <  |\varphi_1'(0)| < \cdots < |\varphi_l'(0)| \quad \textrm{ and } \quad  |\varphi_l'(0)| > e^{\max\left( \sigma_m, \tau(\mathbf{b;e}) \right)}.
\end{equation}
Accordingly, partition the segment $[0,m]$ by introducing the division point parameters 
\[0 = \gamma_0  < \cdots < \gamma_{l} < \gamma_{l+1} = m,\]
 and use these choices to define an $L^1$ function by piecewise patching the functions $\log{|\varphi_k|}$ on the circle~$\T$ according to the linear scaling of~$[0,1)$ to 
$\left[ \gamma_k/m, \gamma_{k+1}/m\right)$: 
\begin{equation}   \label{rearrangement fine}
\begin{aligned}
 g_{ \boldsymbol{\varphi},\boldsymbol{\gamma}}    & : [0,1) \to \R \cup \{-\infty\}, \\  
 g_{ \boldsymbol{\varphi},\boldsymbol{\gamma}}(t) & := \log{\left|\varphi_k\left( e^{2\pi i \frac{mt - \gamma_k}{\gamma_{k+1}-\gamma_k}} \right)\right|} \quad \textrm{on} \quad t \in \left[ \gamma_k/m, \gamma_{k+1}/m\right). 
\end{aligned}
\end{equation}
If there exists an $m$-tuple $f_1, \ldots, f_m \in \Q \llbracket x \rrbracket$ of $\Q(x)$-linearly independent formal functions with denominator types of the form
\begin{equation}   \label{den type int}
f_i(x) = a_{i,0} +  \sum_{n=1}^{\infty} a_{i,n} \frac{x^n}{n^{e_i} [ 1, \ldots, b_{i,1} \cdot n] \cdots [1,\ldots, b_{i,r} \cdot n]}, \qquad a_{i,n} \in \Z, 
\end{equation}
 such that $f_i(\varphi_k(z)) \in \C \llbracket z \rrbracket$ is the germ of a meromorphic function on $|z| < 1$ for all pairs $i = 1, \ldots, m$ 
 and $k = 0, \ldots, l$, 
  then 
\begin{equation}  \label{fine new bound elementary}
\begin{aligned}
m  \leq   \frac{
\displaystyle{ \int_0^1  2t \cdot g_{\boldsymbol{\varphi}, \boldsymbol{\gamma}}^*(t) \, dt + \frac{1}{m} \sum_{k=1}^{l} \gamma_k^2 \log{\frac{|\varphi_{k}'(0)|}{|\varphi_{k-1}'(0)|}} }}{  \log{|\varphi_l'(0)|} - \tau(\mathbf{b;e}) } \\
=  \frac{
\displaystyle{ \int_0^1 \int_0^1  \max\left( g_{\boldsymbol{\varphi},\boldsymbol{\gamma}}(s),g_{\boldsymbol{\varphi}, \boldsymbol{\gamma}}(t) \right)  \, ds \, dt + \frac{1}{m} \sum_{k=1}^{l} \gamma_k^2 \log{\frac{|\varphi_{k}'(0)|}{|\varphi_{k-1}'(0)|}}} }{  \log{|\varphi_l'(0)|} - \tau(\mathbf{b;e}) }. 
\end{aligned}
\end{equation}
If moreover all functions~$f_i$ are \emph{a priori} assumed to be holonomic, the 
assumption $ |\varphi_l'(0)| > e^{\max\left( \sigma_m, \tau(\mathbf{b;e}) \right)}$  on~$\varphi_l$ in equation~\eqref{relaxedifholonomic}
can be relaxed to  $ |\varphi_l'(0)| > e^{ \tau(\mathbf{b;e}) }$. 
\end{thm}

Here, in~$g^*$, we use the notation from~\eqref{incr} of the \emph{increasing rearrangement function} of $g$. This is why we will often refer to quantities like
$\int_0^1 2t \cdot g^*(t) \, dt = \iint_{[0,1]^2} \max\left( |g(s)|, |g(t)| \right) \, ds \, dt$ as to \emph{rearrangement integrals}. 

\begin{remark}[Musical Notation] \label{rem:musical}
In the notation of Theorem~\ref{basic main}, we have $\tau(\mathbf{b}) = \tau^{\flat}(\mathbf{b}) = \tau(\mathbf{b;0})$. Our reason
for the musical notation is to think of  $\tau = \tau^{\flat}(\mathbf{b})$ as the main reduction (flattening) of cruder values such as the value~$\tau = \sigma_m$
from~\cite{zeta5} when we remove the powers $n^{\mathbf{e}}$ from~\eqref{den type integrated}, and of $\tau^{\sharp}(\mathbf{e})$ as the extra term
from adding those integrations to the original list of functions.~\endofremark
\end{remark}

\begin{remark} \label{rem:overflow}
In Theorem~\ref{main:elementary form} (and all the other similar theorems that we prove), we may formally relax the denominator type~\eqref{den type int} to allow for the
looser form:
\begin{equation}   \label{rascal}
n^{e_i} [ 1, \ldots, b_{i,1} \cdot n + c_{i,1}] \cdots [1,\ldots, b_{i,r} \cdot n + c_{i,r}], 
\end{equation}
for any fixed set of integers~$c_{i,j}$. This follows upon applying the original statement of our theorem where all the nonzero~$b_{i,j}$
are changed to~$b_{i,j} + \varepsilon$, for some sufficiently small positive number $\varepsilon > 0$. This subsumes
the denominator type~(\ref{rascal}) for all but finitely many~$n$, and any finite initial string of 
coefficients
can  be made to have any given denominator type by scaling. 
Then one takes the limit~$\varepsilon \rightarrow 0$,
after noting that
 the bounds always depend continuously on the~$b_{i,j}$.~\endofremark
\end{remark}

As the special case~$l=0$ of a single analytic map~$\varphi$, we record the extension of the bound~\eqref{new bound 2}: 

\begin{cor}  \label{basic main corollary}
Assume the same conditions and notation as in Theorem~\ref{main:elementary form}, but consider more simply a single holomorphic
mapping $\varphi : (\Db,0) \to (\C,0)$ satisfying $|\varphi'(0)| > e^{\tau(\mathbf{b;e})}$
and such that the pullbacks~$\varphi^* f_i$
are  meromorphic functions on the open unit disc, that is, $\varphi^* f_i \in \mathcal{M}(\D)$.
  If either $|\varphi'(0)| > e^{\sigma_m}$, or if all functions~$f_i$ are holonomic, then 
\begin{equation} \label{basic rearrangement form}
m  \leq  \frac{\displaystyle{
 \iint_{\T^2} \log\left( \max (|\varphi(z)|, \varphi(w)| \right) \, \mv(z) \mv(w)} }{  \log{|\varphi'(0)|} - \tau(\mathbf{b;e}) }  
 =  \frac{\displaystyle{
 \int_0^1 2t \cdot ( \log{|\varphi(e^{2\pi i t})|} )^* \,  dt}}{  \log{|\varphi'(0)|} - \tau(\mathbf{b;e}) }.
 \end{equation}
\end{cor}

\begin{remark} \label{converse remark}
The \emph{a priori} holonomicity cannot be dropped if we only assume $|\varphi'(0)| > e^{\tau(\mathbf{b;e})}$. 
More precisely, for any given datum~$\left( \bb; \varphi \right)$ in the statement of Corollary~\ref{basic main corollary} (the 
case of a single map~$\varphi$ in Theorem~\ref{main:elementary form}), \emph{except} now with assuming
 the opposite inequality $|\varphi'(0)| \leq e^{\sigma_m}$, a simple inductive construction demonstrates the following. If there is at least one $m$-tuple of $\Q(x)$-linearly independent formal functions~$\{f_i\}$ obeying the arithmetic and analytic conditions of the datum~$\left( \bb; \varphi \right)$, then there are continuum-many such $m$-tuples;  this in particular implies non-holonomic such functions.

 To see the claim, upon keeping fixed~$f_1, \ldots, f_{m-1}$, it suffices to show that~$f_m \in \Q\llbracket x \rrbracket$ has continuum-many valid coefficient options $a_{m,\bullet} \in \Z$ in the form~\eqref{den type int} with~$e_m = 0$, under which   
 the pulled back power series~$\sum c_n z^n :=  f_m (\varphi(z)) \in \C \llbracket z \rrbracket$ has sub-exponentially small coefficients $|c_n| = \exp(o(n))$.  
 (Compare with~\cite[\S~6.4.2]{BostCharles}, \cite[\S~6]{PolyaOriginal}, or~\cite[\S~5]{robinson}.)  We show that each successive coefficient~$a_{m,n} \in \Z$ has at least two valid options after all the preceding coefficients~$a_{m,0}$, $\ldots$, $a_{m,n-1}$ have already been selected. This follows upon recursively expressing 
 $
 c_n = \mu_n a_{m,n} - P_n\left( a_{m,0}, \ldots, a_{m,n-1} \right)
 $
 with coefficient 
 $$
 \mu_n := \varphi'(0)^n \big/ \prod_{i=1}^{h}  [1,\ldots, b_{m,i} \cdot n]
 $$
  of sub-exponential growth by the prime number theorem,
and $P_n \in \C[x_0,\ldots,x_{n-1}]$ polynomials that depend on the map~$\varphi$. 
  This gives the two distinct valid options
 $a_{m,n} \in \left\{ \lfloor P_n(\boldsymbol{a}_{m,<n})/\mu_n \rfloor, \lfloor P_n(\boldsymbol{a}_{m,<n})/\mu_n \rfloor + 1 \right\}$, and altogether a construction of a set of~$f_m \in \Q\llbracket x \rrbracket$ with cardinality $2^{\#\NwithzeroA}  = \#\R$.~\endofremark
 \end{remark}

An essential technical feature in this section, and ultimately in the proofs of both Theorems~\ref{mainA} and~\ref{logsmain}, is the term~$\tau^{\sharp}$ accommodating added integrals to the principal denominators shape of Theorem~\ref{basic main}. We describe this feature on a few examples.

\begin{basicremark}   \label{lcm second half}
To revisit the simplest example from Basic Remark~\ref{basiclog}, let us compute the quantity $\tau(\mathbf{b;e})$ with
\begin{equation*}
\mathbf{b} = (0,0)^{\mathrm{t}}, \qquad \mathbf{e} = (0,1). 
\end{equation*}
Clearly $\tau^{\flat}(\mathbf{b}) = 0$, and we easily compute that 
\begin{equation*}
\tau^{\sharp}(\mathbf{e}) = \frac{1}{2}  \min_{\xi \in [0,2]} \left\{ \frac{3 + (\xi-1)^2}{2} \right\} = \frac{3}{4},
\end{equation*}
attained at the midpoint $\xi =1$. Hence we find the same value $\tau(\mathbf{b,e}) =  \tau^{\flat}(\mathbf{b}) + \tau^{\sharp}(\mathbf{e}) =  3/4$ as when we use the cruder
scheme 
\begin{equation*}
\mathbf{b} =(0,1)^{\mathrm{t}}, \qquad \mathbf{e} = (0,0), \qquad \tau^{\flat}(\mathbf{b}) = 3/4, \quad \tau^{\sharp}(\mathbf{e}) = 0,
\end{equation*}
and no improvement  is made over~\S~\ref{sec:holbasic} in this example. 

In a moment, we will revisit and refine the Diophantine approximation framework of~\cite[\S~2.1]{UDC}. In that framework on our running example, 
we have auxiliary functions (replicated to many variables) of the form $P(x) + Q(x) \log(1-x)$, where $P, Q \in \Z[x]$ are polynomials of degrees less than~$D$.  By the
discussion in~\S~\ref{log discussion}, the lowest order monomial $\beta x^n$ of any such function is necessarily in degree $n \leq 2D-1$, where the equality is
attained uniquely by the Hermite--Pad\'e approximants which are essentially given by  Legendre polynomials. Our proof scheme combines an analytic
upper bound on the coefficient~$\beta \in \Q^{\times}$ with the arithmetic lower bound $|\beta| \geq 1/\mathrm{den}(\beta)$ by the reciprocal of the denominator
of the nonzero rational number~$\beta$. We can directly see why in this case the finer denominators of $\log(1-x) = -\sum_{k=1}^{\infty} x^k/k$ 
do not give any extra help in the arithmetic lower bound $|\beta| \geq 1/\mathrm{den}(\beta)$. The denominator of the $x^n$ coefficient $\beta \in \Q$  is estimated
by the lowest common multiple of all integers from the interval $[n-D,n] \supset [n/2, n]$. As $[n/2,\ldots,n] = [1,\ldots, n]$ (for every integer  $k \in [1,n]$ has a unique $2$-power multiple fitting into $(n/2,n]$), this in the situation is equal to the lowest
common multiple $[1,\ldots,n]$ of the full initial string of integers: the estimate that we get from using $[1,\ldots,k]$ instead of~$k$ as the coefficients denominators
in the function $f(x) = \log(1-x)$. We also see that the finer denominators are expected to make a difference once we have at least $m \geq 3$ functions, as already $[2n/3,\ldots,n]$
is substantially smaller than $[1,\ldots,n]$. (See also
Basic Remark~\ref{remark:zudilingraph}.)  \endofremark
\end{basicremark}

\begin{remark} \label{concatenation}  Using a refined pair $\mathbf{b} \in M_{m \times r}(\NwithzeroB)$ with an $\mathbf{e} \in \NwithzeroB^m$ instead of a 
crude concatenation with $\mathbf{e \rightsquigarrow 0}$ may not always 
give an improvement in the estimate of Theorem~\ref{main:elementary form}. 
Consider the case  $m=3$ with the situation with the proof of Theorem~\ref{logcharacterization}, but with the finer types
\begin{equation*}
\mathbf{b} = (0,0,1)^{\mathrm{t}}, \qquad \mathbf{e} = (0,1,0),
\end{equation*}
giving a template of three functions with denominator types
\begin{equation}  \label{log char dens}
x^n, \qquad \frac{x^n}{n}, \qquad \frac{x^n}{[1,\ldots,n]}.
\end{equation}
This choice for the array $(\mathbf{b;e})$ has $\tau^{\flat}(\mathbf{b}) = \tau^{\sharp}(\mathbf{e}) = 5/9$, with the latter reaching the 
minimum over the whole interval $\xi \in [3/2,2]$. But the type~\eqref{log char dens} is also covered by the cruder choice
\begin{equation*}
\mathbf{b}_0 := (0,1,1)^{\mathrm{t}}, \qquad \mathbf{e}_0 = (0,0,0),
\end{equation*}
that we made for the proof of Theorem~\ref{logcharacterization}, and this basic choice gives in this case the better value $\tau(\mathbf{b}_0;\mathbf{e}_0) =
\tau(\mathbf{b}_0) = 8/9$ than~$\tau(\mathbf{b};\mathbf{e}) = 10/9$.
The reason for this is in how the denominators in the leading order coefficients of the auxiliary functions end up getting estimated 
in the proof of Theorem~\ref{main:elementary form}; the rate $\tau^{\flat}(\mathbf{b}) + \tau^{\sharp}(\mathbf{e})$
serves as an upper bound, and that upper bound estimation turns out to be strict and lossful in the example at hand.
We do not know whether or not the upper bound is an equality in the $m=14$ case that we ultimately devise for the 
proof of Theorem~\ref{mainA}, but we expect it to be a fairly sharp denominator estimate, and possibly an equality.~\endofremark
\end{remark}

\begin{example}   \label{ex:li11}
We give one final example, which we will use in~\S~\ref{sec:completion three elements} to complete the proof of Theorem~\ref{thm three elements}.  For the case
\begin{equation*}
\mathbf{b} = \left( \begin{array}{cc} 0 &  0 \\ 0 & 2 \\ 0 & 2 \\ 1 & 2  \end{array} \right), \qquad \mathbf{e} = (0, 0, 1, 0)
\end{equation*}
relevant to a set of four functions with denominator types
\begin{equation} \label{Li11 problem types}
x^n, \qquad \frac{x^n}{[1,\ldots,2n]}, \qquad \frac{x^n}{n[1,\ldots,2n]}, \qquad \frac{x^n}{[1,\ldots,n][1,\ldots,2n]}, 
\end{equation}
we have
\begin{equation*}
\tau^{\flat}(\mathbf{b}) =  \frac{37}{16}, \quad \tau^{\sharp}(\mathbf{e}) = \frac{7}{16}, \qquad \tau(\mathbf{b;e}) = \frac{21}{8} = 2.625, 
\end{equation*}
with the $7/16$ value being attained on the identical interval minimizer $\xi \in [2,3]$.

But even the intermediate crude choice
\begin{equation*}
\mathbf{b}_0 = \left( \begin{array}{cc} 0 &  0 \\ 0 & 2 \\ 1 & 2 \\ 1 & 2  \end{array} \right), \qquad \mathbf{e}_0 = (0, 0, 0, 0),
\end{equation*}
that minimally covers the types~\eqref{Li11 problem types} within the framework of Theorem~\ref{basic main} highlighted
for the introduction, already gives 
\begin{equation*}
\tau(\mathbf{b}_0; \mathbf{e}_0) = \tau(\mathbf{b}_0) = \frac{1}{16}\left( 1 \cdot 0 + 3 \cdot 2 + 5 \cdot 3 + 7 \cdot 3  \right) = \frac{21}{8} = 2.625. 
\end{equation*}
In this case the value is the same, similarly to the situation in Basic Remark~\ref{lcm second half}.  \endofremark
\end{example}

In contrast to the examples above, 
 exploiting a refined pair $(\mathbf{b;e})$ does often lead to strictly better results than are possible by capping  
up to some cruder $\mathbf{e=0}$ scheme. This is in particular true for the proof Theorem~\ref{mainA} laid out in~\S~\ref{sec:proofA}. 
There we have a local system of rank $m=14$ with added integrals, meaning $e_i = 1$ for six of the indices, and $e_i = 0$ for the remaining eight
indices. Such a vector has~$\tau^{\#}(\mathbf{e}) = 27/80$ after a simple computation~\eqref{tauraise}. Overall the fine $\tau(\mathbf{b;e})$ 
used computes to $191/49+27/80 = 16603/3920  = 4.235459\ldots$. But with $\mathbf{e=0}$ types in this example we do not get a better estimate
than the rather poor $865/196 = 4.413265\ldots$ of Remark~\ref{int is necessary}.

\subsection{Horizontal integration}  \label{horizontal integration}
Our proof scheme follows precisely the $d \to \infty$ asymptotic method 
that we originally devised for our first solution~\cite[route \S~2.5]{UDC} of the unbounded denominators conjecture. This was the idea that we dubbed a \emph{cross-variables 
integration}, where the given single-variable functions $f_i(x)$ were replicated in $d \to \infty$ splitting variables to form, using the Dirichlet box principle, a nonzero auxiliary function
\begin{equation} \label{function form}
F(x_1,\ldots,x_d) := \sum_{\substack{ 
\mathbf{i} \in \{1,\ldots,m\}^d \\
\mathbf{k} \in [0,D]^d \cap \Z^d }} c_{\mathbf{i,k}} \, x_1^{k_1} \cdots x_d^{k_d} \, f_{i_1}(x_1) \cdots f_{i_d}(x_d) \in \Q \llbracket x_1, \ldots, x_d \rrbracket
\setminus \{0\},
\end{equation}
with integer coefficients $c_{\mathbf{i,k}} \in \Z$ of sub-exponential  asymptotic size
 \[|c_{\mathbf{i,k}} |   = \exp\left( o_{d \to \infty} (dD) \right),\]
but yet with $F(x_1,\ldots,x_d)$
vanishing to almost the highest conceivable order $\alpha$ 
at $\mathbf{x = 0}$. 
 Curiously enough,  in this scheme the Nevanlinna characteristic
growth term  $\int_{\T} \log^+{|\varphi|} \, \mv$ arose not as a circle integral \emph{per se} (although the latter is also possible, 
by either~\cite[\S~2.3 or \S~2.4]{UDC}, or by our discussion based~\cite{BostCharles} and~\cite{zeta5} in~\S~\ref{new slopes} below); but rather  ---  by the standard Large Deviations bound mandating\footnote{This 
is just a restatement of Theorem~\ref{thm_MeasureConcentration} which we treated in detail in~\S~\ref{concentration of measure}. Indeed, the 
exponential function establishes an isomorphism of the measure spaces $([0,1]^d, \mu_{\mathrm{Lebesgue}})$ and $(\T^d, \mv)$, under
which the box discrepancy functions correspond.} that the low discrepancy set $D(\mathbf{z}) < \epsilon$ on the high-dimensional torus $\T^d$ has measure at least $1-Ce^{-c \epsilon^4 d}$ once $d \geq d_0(\epsilon)$  ---  from the preponderant growth rate of
the pulled-back monomials $\varphi(z_1)^{k_1} \cdots \varphi(z_d)^{k_d}$. This is why we would describe such an approach
as doing an integration in a cross-variables way, or ``horizontally'' if one pictures the dimension \emph{versus} degree \emph{versus} $\T$ (the complex analysis in 
any one fixed variable) aspects in the
Diophantine approximation construction.

 One of the main findings of the present paper is a certain combination of~\S~\ref{functional transcendence} and~\S~\ref{concentration 
 of measure} which allows to actually improve the meaning here 
of  the ``highest conceivable vanishing order'' $\alpha$ in~\eqref{function form}. In the more rudimentary treatment in~\cite[\S~2]{UDC}, 
we had only $\alpha = mdD/e - o_{d \to \infty}(dD)$ in the parameter count for the number of linear equations
to be solved in the unknown variables $c_{\mathbf{i,k}}$;  this owes to the~$e^d:1$ asymptotic  
volumes ratio for a standard high-dimensional simplex to its largest embedded subcube. 
We now explain how Corollary~\ref{jumps Cartesian}, on the commutativity with Cartesian products of the formation of the vanishing filtration jumps sets, 
 and the measure concentration material~\S~\ref{concentration of measure} work together
to improve this vanishing order to $\alpha = mdD/ 2 -o_{d \to \infty}(dD)$. For simplicity, since we will anyway need this later on for the general form of Theorem~\ref{main:elementary form}, we assume the strongest form of the Cartesian power structure (based on the holonomicity of the $f_i$): Lemma~\ref{lem_Shidlovsky}, stemming from the Chudnovsky--Osgood theorem coupled to Corollary~\ref{jumps Cartesian}.

\subsubsection{The Thue--Shidlovsky idea} \label{Siegel sketch}
This Lemma~\ref{lem_Shidlovsky} 
 ensures that the nonzero power series~\eqref{function form} has to 
posses a nonzero monomial $\beta \, \bx^\bn = \beta \, x_1^{n_1} \cdots x_d^{n_d}$ with $\bn = (n_1, \ldots, n_d) \in [0,(m+\delta)D]^d$, for any $\delta > 0$ and $d$, once 
$D \gg_{\delta,d} 1$. Therefore, in the linear system to solve for the total vanishing order in the Thue--Siegel lemma, we need not be concerned with the broad simplex region
$|\mathbf{n}| < \alpha$, but instead we
 can simply vanish the coefficients of $\mathbf{x^n}$ for all $\mathbf{n}$ ranging over
the hypercube $[0,(m-\delta)D]^d$ (clearly, this is the hypercube of the maximal conceivable size in the parameter count), 
as well as for all $\mathbf{n}$ outside of the low discrepancy part of the slightly bigger hypercube $[0,(m+\delta)D]^d$ (allowing us to also use the measure concentration for the component set of the vector $\mathbf{n}$). It is essential here for the application of Theorem~\ref{thm_MeasureConcentration} to take~$\delta > 0$ sufficiently small in terms of the discrepancy parameter~$\epsilon$; then the parameter count goes through. (This step is contained in Lemma~\ref{Siegel}.)
The upshot is that  in this construction, as $\epsilon \to 0$ (eventually: after $D \to \infty, d \to \infty, \delta \to 0$), all the lowest order  monomials $\beta \cdot x_1^{n_1} \cdots x_d^{n_d}$ in $F(x_1,\ldots,x_d)$ have their exponent vectors $(n_1,\ldots,n_d)$ asymptotically close to some
ordering of the set $\{ j m D/d \, : \, 0 \leq j < d \}$. In particular, the total vanishing order is indeed close to $mdD/2$.

 It is this improvement over the $mdD/e$ of~\cite[Lemma~2.1.2]{UDC} that recovers the $e \rightsquigarrow  2$ coefficient refinement under the elementary asymptotic framework of~\cite[\S~2]{UDC}.  
At this point, the Thue--Shidlovsky argument further supplies two
similar, and equally crucial, improvements on both the top and bottom of the fraction~\eqref{founding hol}. 

We next introduce these two improvements in turn, showcasing to the case of Corollary~\ref{basic main corollary} for simplicity, following which we start on the proof of the general Theorem~\ref{main:elementary form}.

\subsubsection{The archimedean sharpening} \label{sharp1}  
Firstly, this asymptotic scheme allows by elementary methods to improve the $\int_{\T} \log^+{|\varphi|} \, \mv$ integral to 
the strictly smaller quantity $\int_0^1 t \cdot ( \log{|\varphi(e^{2\pi i t})|} )^* \, dt$.  We have remarked already in~\cite[\S~2.3.3]{UDC} that some improvement in the holonomy bound can be made 
by exploiting that, by Theorem~\ref{thm_MeasureConcentration} again, the monomials exponents $\mathbf{k}$ in~\eqref{function form} can
be constricted to the low discrepancy part of the hypercube $[0,D]^d$. In concrete heuristic terms, this means that as $d \to \infty$, the exponents vectors $(k_1,\ldots,k_d)$ in~\eqref{function form} can be considered as being close to some ordering of the set $\{ j D / d \, : \, 0 \leq j < d \}$. 
In this way, upon noting that the largest value of 
$\left| \varphi(z_1)^{k_1} \cdots \varphi(z_d)^{k_d} \right|$ over all these orderings occurs when $(k_1,\ldots,k_d)$ is arranged in the same way as $\left( |\varphi(z_1)|, \ldots, |\varphi(z_d)| \right)$, we find the rearrangement integral
\begin{equation} \label{filtering rearrangement}
\int_0^1 t \, (\log{|\varphi|})^*(t) \, dt  = \frac{1}{2} \iint_{\T^2} \log \left(  \max(|\varphi(z)|, |\varphi(w)|) \right) \, \mv(z) \mv(w)
\end{equation} 
 precisely in the refined preponderant growth rate based on the uniform distribution of not only the torus points $\mathbf{z} \in \T^d$, but 
 also of the monomials exponent vector $\mathbf{k}$.
 An illustration of the saving thus made is in the explicit example of Figure~\ref{bivalc2}.
 
 \subsubsection{The arithmetic sharpening} \label{sharp2}
 To the uniform distribution of the 
 leading order jet exponents $\mathbf{n}$ of $F(x_1,\ldots,x_d)$, we can add yet another elementary application of the Law of Large Numbers: 
 without changing (broadly speaking)
 the asymptotic size
 of the parameter count, we can insist that the $m$ function
 species occur with equal frequency $1/m$ in all the split-variable products $f_{i_1}(x_1) \cdots f_{i_d}(x_d)$ in the
 make up of~\eqref{function form}. 
 In other words, we can assume that $d \equiv 0 \mod{m}$ and that the summation multi-index $\mathbf{i}$ in~\eqref{function form} is constricted to have 
  each index $i_0 \in \{1,\ldots,m\}$ arise $d/m$ times as a component of~$\mathbf{i}$. This permits us to integrate over the $m$ different
  denominator types of our $m$ function species, under our condition on the denominators cap array~$\mathbf{b}$, 
  and we find precisely~\eqref{finitary rearrangement} as the counterpart of~\eqref{maxintegral}.

\subsection{The auxiliary construction}\label{aux_eqdis}
We start here the proof of Theorem~\ref{main:elementary form}. 
In the case that all $f_i$ are \emph{a priori} holonomic functions, we have the following improvement on~\cite[Lemma~2.1.2]{UDC}.

For $d \in \NwithoutzeroA$ and $\epsilon \in (0,1]$, we denote by 
\begin{equation} \label{typical}
P_{\epsilon}^d := \left\{ \mathbf{t} \in [0,1]^d \, : \, D(\mathbf{t}) < \epsilon  \right\} \subset [0,1]^d
\end{equation}
the $\epsilon$-discrepancy part of the $d$-dimensional hypercube, with $D : [0,1]^d \to [0,1]$ being the 
normalized discrepancy function of Definition~\ref{box discrepancy}. The image of this set under the analytic isomorphism $\exp : [0,1)^d \to \T^d ,\, \mathbf{t} \mapsto e^{2 \pi i \mathbf{t}}$ will be denoted by $T_{\epsilon}^d \subset \T^d$. 
 In these notations, a cruder form of Theorem~\ref{thm_MeasureConcentration} can
be restated as the double limits
\begin{equation} \label{asymptotically full}
\lim_{\varepsilon \to 0} \lim_{d \to \infty} \mu_{\mathrm{Lebesgue}}(P_{\varepsilon}^d) = 1, \qquad
 \lim_{\varepsilon \to 0} \lim_{d \to \infty} \mv(T_{\varepsilon}^d) = 1.
\end{equation}
In the following, we fix an $m \in \NwithoutzeroA$ and restrict the asymptotic parameter~$d \in \NwithoutzeroA$ to the integers $\equiv 0\mod{m}$. For the plan outlined in~\S~\ref{sharp2}, we 
restrict the multi-index $\mathbf{i}$ in~\eqref{function form} to the equidistributed set
\begin{equation}
V_m^d := \left\{ \mathbf{i} \, : \, \forall i \in \{1,\ldots,m\}, \, \#\{ h \in \{1,\ldots,d\} \, : \, i_h = i \} = d/m  \right\}  \subset \{1,\ldots,m\}^d. 
\end{equation}

This set still has the asymptotically full size $m^{d-o(d)}$: 
\begin{lemma} \label{balanced indices}
Under our standing assumption $d \equiv 0 \mod{m}$, we have
\begin{equation}
\#V_m^d = \binom{d}{\frac{d}{m},\ldots,\frac{d}{m}} > m^d \Big/ \binom{d+m-1}{m-1}. 
\end{equation}
\end{lemma}

\begin{proof}
The $m$-fold expansion
$$
m^d = (1+\ldots+1)^d = \sum_{ \substack{ \mathbf{j} 
\in \NwithzeroB^m  \\ |\mathbf{j}| = m }} \binom{d}{j_1, \ldots, j_m}
$$
has $\binom{d+m-1}{m-1}$ terms, the maximal of which is the central multinomial coefficient with $j_1 = \cdots = j_m = d/m$. 
\end{proof}

Our Thue--Siegel construction (step~\ref{dirichletbox} of the general outline from~\S~\ref{sec:ideas outline}) is the following. 

\begin{lemma} \label{Siegel}
Suppose we have $m$ \emph{holonomic} power series $f_1,\ldots,f_m \in \Q\llbracket x \rrbracket$.  
 Assume that each $f_i(x)$ has denominators of the crude type
\begin{equation} \label{crude G cap}
f_i(x) = \sum_{n=0}^{\infty} a_{i,n} \frac{x^n}{A^{n+1} [1,\ldots,Bn]^{\sigma}}, \qquad \textrm{for some } A \in \NwithoutzeroA, \, B, \sigma \in \NwithzeroA
\end{equation}
and converges on a complex disc $|x| < \rho$, for some $\rho \in (0,1)$.

There exists a function
$$
d_0 : \NwithzeroB^3 \times (0,1) \times (0,1) \to \NwithzeroA
$$
such that the following holds.

For each  $\epsilon \in (0,1)$, there is a~$\delta = \delta(\epsilon) \in (0,\epsilon)$, such that 
for all $d \geq d_0(A,B,\sigma, \rho; \epsilon)$ with $d \equiv 0 \mod{m}$, 
there exists asymptotically
  for $D \to \infty$ a nonzero $d$-variate formal function $F(\mathbf{x})$ of the~\eqref{function form} form
\begin{equation}  \label{F construction}
 F(\mathbf{x}) =     \sum_{\substack{  \mathbf{i} \in V_m^d  \\ 
\mathbf{k}/D \in P_{\epsilon}^d \cap \Z^d }   } c_{\mathbf{i,k}} \,
\mathbf{x^k} \, \prod_{s=1}^d f_{i_s}(x_s) 
\in \Q\llbracket \mathbf{x} \rrbracket \setminus \{0\},
\end{equation}
   with $c_{ \mathbf{i,k} } \in \Z$ integers, all bounded in absolute value by
   $|c_{\mathbf{i,k}}| < e^{\epsilon d D}$, and such that the power series expansion $F(\mathbf{x}) = \sum b_{\mathbf{n}} \mathbf{x^n}$ of~\eqref{F construction}
   obeys the following main requirement:  
   \begin{equation}
  (\star) \qquad \qquad \label{eqstar}
  \begin{aligned}
       &  \text{All the minimal order monomials $\beta_{\mathbf{n}} \, \mathbf{x^n}$ in~\eqref{F construction}}  \\ 
             & \text{have an exponent vector $\mathbf{n}$ satisfying $\mathbf{n}/\left((m+\delta)D\right) \in P_{\epsilon}^d$. }
        \end{aligned}
        \end{equation}
   \end{lemma}

Note that ($\star$) implies in particular that $\beta_{\mathbf{n}} \neq 0$ for at least one such   $\mathbf{n}$
in the set $\left((m+\delta)D\right) P_{\epsilon}^d \cap \Z^d$. 

Before we proceed with the proof, we collect some consequences of our condition~$(\star)$. In the following, we will consider $\epsilon \in (0,1/4]$ which in the end will be let to approach zero. Throughout
this~\S~\ref{fine section}, we will write
\begin{equation} \label{alpha order}
\alpha := mdD/2.
\end{equation}
This notation reflects a related use of~$\alpha$  as
a vanishing parameter in~\cite[\S~2]{UDC}, see for 
example~\cite[Lemma~2.1.2(1)]{UDC}. In that previous paper, we
 (asymptotically) took~$\alpha$ to be (of order)
$mdD/e \sim  m (d!)^{1/d} D$, which here we  improve to~$mdD/2$.

\begin{cor} \label{F stats}
In Lemma~\ref{Siegel}, we have
\begin{equation} \label{vanishing order}
\mathrm{ord}_{\mathbf{x} = \mathbf{0}} F(\mathbf{x}) \in [(1-2\epsilon)\alpha, (1+2\epsilon)^2 \alpha]. 
\end{equation}
Every multi-index  $\mathbf{k} = (k_1,\ldots,k_d)$ in~\eqref{F construction} has $\max_{j=1}^d k_j \leq  \frac{2\alpha}{md}$ and admits a permutation $\psi = \psi_{\bk}$ of $\{1,\ldots,d\}$
such that $k_{\psi(1)} \leq \cdots \leq k_{\psi(d)}$ and 
\begin{equation} \label{F degrees}
 \frac{2\alpha j}{md^2}  - 2\epsilon \alpha/(md)  \leq k_{\psi(j)}  \leq (1+\epsilon)  \, \frac{2\alpha j}{md^2} + 4\epsilon \alpha/(md), \qquad \forall j \in \{1,\ldots,d\}. 
\end{equation}
Further, every exponent $\bn = (n_1,\ldots,n_d)$ of minimal total order
 \[n := |\bn| = \mathrm{ord}_{\mathbf{x} = \mathbf{0}} F(\mathbf{x})\]
 in the Taylor series of $F(\bx) \in \Q \llbracket \bx \rrbracket$ has a 
permutation $\pi$ of $\{1,\ldots,d\}$ such that $n_{\pi(1)} \leq \cdots \leq n_{\pi(d)}$ and
\begin{equation} \label{order n pi} 
2\alpha j/d^2 - 2\epsilon \alpha/d  \leq n_{\pi(j)}  \leq (1+\epsilon) 2\alpha j /d^2 + 4\epsilon \alpha/d, \qquad \forall j \in \{1,\ldots,d\},
\end{equation}
and for all $u,v \in [0,1]$ with $u \leq v$ it satisfies
\begin{equation} \label{uni F}
(1-2\epsilon) (v^2-u^2) \alpha \leq  \sum_{u d \leq j < v d} n_{\pi(j)} \leq (1+2\epsilon)^2 (v^2-u^2) \alpha. 
\end{equation}
\end{cor}

\begin{proof}
The partial degrees bound~\eqref{F degrees} is tautologically a rewriting of our definition~\eqref{alpha order}, but it is used
to organize the analysis around the leading asymptotic parameter~$\alpha$. The estimate~\eqref{vanishing order} follows from~\eqref{eqstar} upon
noting that the expected value $\mathbf{E}\left[t \in [0,1] \right] = \int_0^1 t \, dt = 1/2$, which by the Koksma--Hlawka inequality or an elementary bit of computation 
shows the implication
\begin{equation} \label{discrepancy t}
\mathbf{t} \in P_{\epsilon}^d \quad \Longrightarrow  \quad  \sum_{j=1}^d t_j \in  \left[d/2-\epsilon d, d/2+\epsilon d \right].
\end{equation}
In detail, the upper bound in~\eqref{vanishing order} is by the chain of trivial estimates $|\bn| \leq (m+\delta)D (d/2+\epsilon d) \leq (1+\epsilon)  mD  (d/2+\epsilon d) 
< (1+2\epsilon) mD  (d/2+\epsilon d)  = (1+2\epsilon) (2\alpha/d)  (d/2+\epsilon d) = (1+2\epsilon)^2 \alpha$ implied by~\eqref{discrepancy t} for all the nonzero monomials~$\beta_\bn \bx^\bn$ (which form a nonempty set!) in~$(\star)$. The
lower bound is similar with using
$|\bn| \geq (m+\delta)D (d/2 - \epsilon d) >  mD  (d/2 - \epsilon d) 
  =  (2\alpha/d)  (d/2-\epsilon d) = (1-2\epsilon) \alpha$ from~\eqref{discrepancy t} for all nonzero monomials~$\beta_\bn \bx^\bn$ in~$(\star)$, and the
  proof of~\eqref{uni F} is the same.  
  
  For an arbitrary $\bn$ from the leading order $|\bn| =  \mathrm{ord}_{\mathbf{x} = \mathbf{0}} F(\mathbf{x})$ jet~$(\star)$, take a 
  permutation~$\pi$  with $n_{\pi(1)} \leq n_{\pi(2)} \leq \cdots \leq n_{\pi(d)}$. 
  For the lower bound in~\eqref{order n pi}, remark that the interval $[0, n_{\pi(j)}]$ contains at least the~$j$ elements
  $n_{\pi(1)}, \ldots, n_{\pi(j)}$ of the component set $\{ n_j \}$. As $\bn \in \left((m+\delta)D\right)  P_{\epsilon}^d$, 
  the definition of discrepancy mandates that the interval $[0, n_{\pi(j)}] = ((m+\delta)D) \cdot  \left[0, n_{\pi(j)}/((m+\delta)D) \right)$ contains
  at most 
  $$
  n_{\pi(j)}d/((m+\delta)D) + \epsilon d < n_{\pi(j)}d/(mD) + \epsilon d = 
  n_{\pi(j)} d^2/(2\alpha) + \epsilon d
  $$
   of the components of~$\bn$. Hence
   $$
   j  \leq   n_{\pi(j)} d^2/(2\alpha) + \epsilon d,
   $$
   giving the claimed lower bound on~$n_{\pi(j)}$. 
   
   For the upper bound in~\eqref{order n pi}, the interval 
   $$
   \left[n_{\pi(j)},\left((m+\delta)D\right) \right]
   = \left((m+\delta)D\right) \cdot \left[ n_{\pi(j)}/\left((m+\delta)D\right), 1 \right]
   $$
   contains at least the $d-j+1$ elements $n_{\pi(j)}, \ldots, n_{\pi(d)}$ of the component set~$\{ n_j \}$ 
   of our $\bn \in \left((m+\delta)D\right)  P_{\epsilon}^d$. Thus, as before, 
   \begin{equation*}
   \begin{aligned}
d-j+1 \leq d \cdot \left( 1 - n_{\pi(j)}/\left((m+\delta)D\right) \right) + \epsilon d \\
\leq  d - n_{\pi(j)} d/\left((1+\epsilon)mD\right)  + \epsilon d  =  d - n_{\pi(j)} (1+\epsilon)^{-1} d^2/(2\alpha)  + \epsilon d 
   \end{aligned}
   \end{equation*}
   completing the proof of~\eqref{order n pi}. Finally, the bounds~\eqref{F degrees} follow by the same proof. 
\end{proof}

\subsection{The box principle step: proof of Lemma~\ref{Siegel}} Our proof combines
 the classical Thue--Siegel lemma~\cite[Lemma~2.9.1]{BombieriGubler}
 and Lemma~\ref{lem_Shidlovsky} on the functional bad approximability combined with 
 the product structure of the vanishing filtration jumps,
 following the plan we laid out in~\S~\ref{Siegel sketch}.  

\begin{proof}[Proof of Lemma~\ref{Siegel}] The first
parameter to consider is the $\epsilon$ that we use to measure the discrepancy from $\mu_{\mathrm{Lebesgue},[0,1)}$ or $\mv$; this will be the last parameter which
we  let  approach~$0$ in the proof. 
Then, in terms of the exponential decay rate function~$\kappa(\epsilon) := \epsilon^4/300  > 0$ in Theorem~\ref{thm_MeasureConcentration}, 
we take any $\delta \in (0,\epsilon)$ small enough to have 
\begin{equation} \label{take delta}
\frac{m - \delta}{m+\delta} >  e^{-\kappa(\epsilon)} = e^{-\epsilon^4/300}.
\end{equation}
We then
  set up a linear system of 
 $M \leq 2\left( (m-\delta)D \right)^d$
 linear equations in the $N = \left( m -o_{d \to \infty}(1))D \right)^d$
 unknown coefficients $c_{\mathbf{i,k}}$ in the function template form~\eqref{F construction}, by requiring that in the Taylor series expansion
 \begin{equation} \label{functional template}
 F(\mathbf{x}) =     \sum_{\substack{  \mathbf{i} \in V_m^d  \\ 
\mathbf{k}/D \in P_{\epsilon}^d \cap \Z^d }   } c_{\mathbf{i,k}} \,
\mathbf{x^k} \, \prod_{s=1}^d f_{i_s}(x_s)   =  \sum_{\mathbf{n} \in \NwithzeroB^d} \beta_{\mathbf{n}} \,\mathbf{x^n},
 \end{equation}
 all coefficients $\beta_{\mathbf{n}}$  vanish whenever either $\max_{j=1}^{d} n_j \leq (m-\delta) D$ or 
 $\mathbf{n} \notin (m+\delta) D  \cdot P_{\epsilon}^d$. Once the dimension $d \gg_{\epsilon,\delta} 1$ is sufficiently big, Theorem~\ref{thm_MeasureConcentration}
 and Lemma~\ref{balanced indices} 
 show that the number~$N$ of free parameters $c_\mathbf{{i,k}}$ in our linear system will exceed the quantity
 \begin{equation}
 N > \big( (m-\delta/2)D \big)^d >  2\big( (m-\delta)D \big)^d \geq M. 
 \end{equation}
 In the Siegel lemma, this gives a \emph{Dirichlet exponent}
 \begin{equation} \label{small Dirichlet exponent}
 \frac{M}{ N - M} < \frac{ 1}{\frac{1}{2} \cdot \left( \frac{m-\delta/2}{m-\delta} \right)^d -1} = o_{d \to \infty}(1).
 \end{equation}
 For the height of our linear system,
  a simple estimate based on the prime number theorem shows that the system can be expressed into the form $\mathbf{A} \cdot \mathbf{y} = 0$,
  to be solved nontrivially for an integer vector $\mathbf{y} \in \Z^N$ of a small height, with some $M \times N$
  integer matrix $\mathbf{A} \in \mathrm{M}_{M \times N}(\Z)$ whose entries are bounded in absolute
  value by $C_0(A,B,\sigma,\rho)^{\alpha}$. Here, $C_0(A,B,\sigma,\rho)$ is a simple computable function, immaterial to us, in the 
  parameters $A,B,\sigma,\rho$ that we assume for the form~\eqref{crude G cap} of the functions $f_i(x)$ (archimedeanly convergent on $|x| < \rho$). 
  At this point~\eqref{small Dirichlet exponent} and~\cite[Lemma~2.9.1]{BombieriGubler} prove that, 
  once $d \gg_{\epsilon,\delta} 1$ and then $D \gg_d 1$, there exists a nonzero formal function $F(\mathbf{x}) \in \Q \llbracket \mathbf{x} \rrbracket \setminus \{0\}$ 
  of the form~\eqref{functional template}, in which all coefficients $c_{\mathbf{i,j}} \in \Z$ on the left-hand side are rational integers smaller than $e^{\epsilon d D}$
  in absolute value, and having on the right-hand side the vanishing of all $\beta_{\mathbf{n}} = 0$ with  $\mathbf{n} 
  \notin (m+\delta)D \cdot P_{\epsilon}^d$, as well as for all $\mathbf{n} \notin \left[ 0, (m - \delta)D \right]^d$. 
  
  The desired property~$(\star)$ (see Equation~\ref{eqstar}) now follows by Lemma~\ref{lem_Shidlovsky},
  applied with~$\varepsilon := \delta/2$, after noting our assumption that $f_1,\ldots, f_m$ are holonomic functions.  
  \end{proof}

The condition in Lemma~\ref{Siegel} that all $f_i$ are holonomic functions is met by the hypotheses in Theorem~\ref{main:elementary form} currently under proof.  Indeed, an \emph{a priori} holonomicity is either directly an assumption, or else the stronger positivity condition~\eqref{relaxedifholonomic} is imposed. By Andr\'e's holonomicity criterion (Corollary~\ref{holonomic criterion}, also outlined in~\S~\ref{sec:ideas outline}, and completely proved in the self-contained~\S~\ref{app:PerelliZannier}), 
upon applying to~$f_i(x)$ the differential operator $\left( x \frac{d}{dx}\right)^{e_i}$ to remove the extra~$n^{e_i}$ terms from the denominators of~\eqref{den type int} (and observing that, thanks to the chain rule, the $d/dx$ derivation preserves the meromorphicity condition $\varphi_l^*f_i \in \mathcal{M}(\Db)$),  the condition $\log{|\varphi_l'(0)|} > \sigma_m \geq b_{i,1} +\ldots + b_{i,r}$ by itself forces~$f_i$ to be a holonomic function. 

This places us into a position to apply Lemma~\ref{Siegel}. In the following, we will fix an ``auxiliary function'' $F(\bx) \in \Q \llbracket \bx \rrbracket \setminus \{0\}$ supplied by that lemma, and write $\delta := \delta(\epsilon)$ for the~$\delta \in (0,\epsilon)$ under the thesis of the lemma. At the end of the proof we will let, in this order, $D \to \infty$, $d \to \infty$, and $\epsilon \to 0$, remembering that the latter also in particular makes~$\delta \to 0$. 

\subsection{Seeding}  \label{seeding} Consider now a nonzero minimal order monomial $\beta \, \bx^\bn$ in $F(\bx)$. Thus $\beta := \beta_\bn \in \Q^{\times}$ is a nonzero rational number of a certain denominator cap inherited from~\eqref{den type int}, that we will study in~\S~\ref{denominator arithmetic} below, and the exponent~$\bn = (n_1, \ldots, n_d) \in \left((m+\delta)D\right) P_{\epsilon}^d $ with $\delta<\epsilon$ has $n := |\bn| = n_1 +\ldots + n_d \in [(1-2\epsilon)\alpha, (1+2\epsilon)^2\alpha]$ by Corollary~\ref{F stats}. {\it Until the end of the proof, we fix this minimal order exponent~$\bn$, and then upon relabeling the variables $x_1, \ldots, x_d$, we may and will assume that $n_1 \leq n_2 \leq \cdots \leq n_d$.} 
 
We turn now to the piece of the argument  ---  which we omitted from the introductory sketch~\S~\ref{horizontal integration} (but we briefly described in~\S~\ref{high dimensional analysis} of our general introduction),  ---  needed to get the stronger bound~\eqref{rearrangement fine} in place of the more basic special case~\eqref{basic rearrangement form}. The idea is to partition the indexing set $\{1,\ldots,d\}$ into~$l+1$ groups so as to use the map~$\varphi_k$ in the analytic 
variable~$z_j$ for the case $\gamma_k/m \leq j/d < \gamma_{k+1}/m$, for $k = 0,\ldots,l$, with the understanding that $\gamma_{l+1}/m = 1$ and equality is meant on the right-hand side condition for $k = l$.  Let $\Phi : \Db^d \to \C^d$ be the diagonal map thus 
defined from using $\varphi_k(z_j)$ for its $j^{\mathrm{th}}$ coordinate function, where $k = k(j) \in \{0,\ldots,l\}$ is uniquely determined by $j \in \left\{
 \lceil \gamma_k d/m \rceil, \ldots, \lceil \gamma_{k+1} d/m \rceil - 1 \right\}$ (and $k=l$ for $j=d$). 
 This is a holomorphic mapping with $\Phi(\boldsymbol{0}) =\mathbf{0}$, and  ---  clearing a common holomorphic denominator for the meromorphic functions $f_i(\varphi_k(z))$, $1 \leq i \leq m$, $0 \leq k \leq l$,  ---  there is a holomorphic function $h \in \mathcal{O}(\Db)$ with $h(0) =1$ and
 \begin{equation}  \label{G pull}
 G(\bz) := h(z_1)\cdots h(z_d) \cdot (\Phi^* F) (\bz) \in \mathcal{O}\left(\Db^d\right)
 \end{equation}
 holomorphic on some neighborhood of the closed unit polydisc. 
 By construction, the $\bz^\bn$ coefficient of~$G(\bz)$ equals\footnote{Formally exponentiating the additive notation, choosing any branch for the logarithm.}
 \begin{equation} \label{G coeff}
[\bz^\bn] \left\{ G(\bz) \right\} = 
 \beta \cdot \exp \left( \sum_{k=0}^l  \left( \sum'_{ 
 \frac{d \gamma_k}{m} \leq j < \frac{d\gamma_{k+1}}{m} }   n_j \right)  \log{\varphi_k'(0)}  \right),
 \end{equation}
 where the dash in the inner summation is to remind us that for $k = l$ the term $j = d$ is supposed to also be included into the sum. 
By Lemma~\ref{F stats}, the inner sum over~$j$ satisfies
\begin{equation}  \label{seed sum est}
 (1-2\epsilon) (\gamma_{k+1}^2-\gamma_k^2) \frac{\alpha}{m^2}  \leq 
 \sum_{ 
 \frac{d \gamma_k}{m} \leq j < \frac{d\gamma_{k+1}}{m} }   n_j  \leq (1+2\epsilon)^2  (\gamma_{k+1}^2 - \gamma_k^2)\frac{\alpha}{m^2}.
\end{equation}

\subsection{Equidistribution} \label{equidistribution}  There are at least two ways~\cite[\S~2.4 or \S~2.5]{UDC} to handle the archimedean growth term
in a manner compatible with our finer analysis. The Vandermondian damping factors of~\cite[\S~2.5]{UDC} are based directly on the Cauchy formula,
and are more in line with the cross-variables integration technique that we exploit to carry out~\S~\ref{sharp1} and~\S~\ref{sharp2}. The Poisson--Jensen
method~\cite[\S~2.4]{UDC} is based on a lexicographical induction lemma~\cite[Lemma~2.4.1]{UDC} suggested to us by Andr\'e; this approach
is closer in spirit to our treatment in~\S~\ref{slopes}. The third proof~\cite[\S~2.3]{UDC} of our original holonomicity theorem for the solution of the
unbounded denominators conjecture does not seem to apply to the present refinement. 

Our choice hence will be to stick to the Vandermondians method for the details of the current section (nevertheless referring to~\cite[\S~2.5.1]{UDC}
for some of basic and well-known facts of potential theory).

\subsubsection{Vandermondians}
To set up our damping factor, we collect here some basic facts from the logarithmic potential theory in the complex plane. Given (a block of) variables~$\bz = (z_1,\ldots,z_d)$, we define 
\begin{equation} \label{Vandermonde}
V(\mathbf{z}) := \prod_{i< j}(z_i - z_j) = \det{ \begin{bmatrix} 1 & z_1 & z_1^2 & \cdots & z_1^{d-1} \\ 1 & z_2 & z_2^2 & \cdots & z_2^{d-1} \\
\vdots & \vdots & \vdots & \cdots & \vdots \\
1 & z_d & z_d^2 & \cdots & z_d^{d-1}\end{bmatrix} } \in \Z [z_1, \ldots, z_d] \setminus \{0\}.
\end{equation}
As in~\cite[~\S~2.5.1]{UDC}, we note: 
\begin{lemma}[Fekete] \label{Fekete}
The supremum of $|V(\mathbf{z})| = \prod_{1 \leq i < j \leq d} |z_i-z_j|$ over the unit polycircle $\mathbf{z} \in \T^d$ is equal to $d^{d/2}$, with equality if and only if the points $z_1, \ldots, z_d$ are the vertices of a regular $d$-gon. 
\end{lemma}

\begin{lemma}[Bilu] \label{non-equidistribution} 
There  are functions $c(\varepsilon) > 0$ and $d_0(\varepsilon) \in \R$  such that, for every $\varepsilon \in (0,1]$, if $d \geq d_0(\varepsilon)$  and $\mathbf{z} = (z_1, \ldots, z_d) \in \T^d$ is a $d$-tuple with discrepancy $D(\mathbf{z}) \geq \varepsilon$, then
\begin{equation} \label{damping grace}
|V(\mathbf{z})| = \prod_{1 \leq i < j \leq d} |z_i - z_j| < e^{-c(\varepsilon) d^2}. 
\end{equation}
\end{lemma}

\begin{proof} 
See the proof in~\cite[Lemma~2.5.8]{UDC}, which in turn is closely based on Bombieri and Gubler's treatment~\cite[page~103]{BombieriGubler}
of Bilu's equidistribution theorem for points of small canonical height on linear algebraic tori. 
\end{proof}

{\it At this point, we fix an $\epsilon > 0$ and a $\delta \in (0,\epsilon)$ until the end of the proof, and
 we assume~$d \geq d_0(\epsilon, \delta)$. } 

\subsubsection{Holomorphic dampener} We suppose now the $\{1,\ldots,d\}$ partitioning into~$l+1$ consecutive blocks from~\S~\ref{seeding}, and we
write $\bz = \left( \bz^{(0)}, \ldots, \bz^{(l)} \right)$ for the corresponding variable blocks. Explicitly, $\bz^{(k)}$ enlists, in increasing labels, the variables
$z_j$ where $\gamma_k d/m \leq j < \gamma_{k+1} d/m$, and the end term $ j = d$ is assumed to be included in the case $k = l$. Following~\cite[\S~2.5.14]{UDC}, we will dampen the integrand in the Cauchy integral formula for the coefficient~\eqref{G coeff} by using the following choice of multivariable holomorphic multiplier: 
\begin{equation}  \label{hol dampener}
W(\bz) := \prod_{k=0}^l V\left(\bz^{(k)}\right)^M \in \Z[\bz] \setminus \{0\},
\end{equation}
where~$M$ is a large integer parameter to be selected in the proof. 

\subsubsection{Cross-integration} \label{cross-integration}
The idea for the cross-integration is simple. Consider $\mathbf{z} \in \T^d$ on the
high-dimensional unit torus. If one of our $k \in \{0,\ldots,l\}$ blocks has discrepancy $D\left(\mathbf{z}^{(k)} \right) \geq \epsilon > 0$, Lemma~\ref{non-equidistribution} tells us
that the corresponding term $V\left(\mathbf{z}^{(k)}\right)$ in~\eqref{hol dampener} is uniformly\footnote{As a function of~$d \in \NwithoutzeroA$, but for the fixed~$\epsilon > 0$.} exponentially small in $-d^2$. This, in combination with Lemma~\ref{Fekete} used as a uniform upper bound on the other factors $V\left(\mathbf{z}^{(q)}\right)$ for $q \in \{0,\ldots,l\} \setminus \{k\}$, entails that the overall damping factor~$W(\bz)$ decays at the exponential rate~$-Md^2$, 
uniformly in~$d \in \NwithoutzeroA$ and $\left\{ \bz^{(q)} \, : \, q  \in \{0,\ldots,l\} \setminus \{k\} \right\}$. This proves that
\begin{equation} \label{W dampens}
\sup_{\substack{ \bz \in \T^d \\ \exists k : \, D(\bz^{(k)}) \geq \epsilon }} \left\{ |W(\bz) | \right\}  < e^{-c'(\epsilon) M d^2},
\end{equation}
with some function~$c'(\epsilon) > 0$ depending on~$\epsilon$ but not on~$d$. 

In addition, momentarily using 
$d = d_0 +\ldots + d_l$ to denote the partition of the variable slot cardinalities, Lemma~\ref{Fekete} also implies the uniform $\exp(o(Md^2))$ upper bound  
\begin{equation}   \label{W not too big}
\sup_{\mathbf{z} \in \T^d} |W(\mathbf{z})| \leq
\prod_{k=0}^l d_k^{Md_k/2} \leq  d^{Md/2}.  
\end{equation} 

The effect of using a multiplier with~\eqref{W dampens} and~\eqref{W not too big} is roughly the following.
Since the monomial exponent vectors $\bk$ in the make up of~$G$ via~\eqref{F construction} have asymptotically uniformly distributed components $\{k_j\} \subset [0,D]$, 
the rearrangement inequality
brings out the function~\eqref{rearrangement fine} and entails in the limit for the product $W(\bz) G(\bz)$ to sift out the mean growth rate 
\[ \exp \left( D \int_0^1 t \cdot (g_{\boldsymbol{\varphi},\boldsymbol{\gamma}})^*(t)   \, dt  \right)\]
as a uniform $\mathbf{z} \in \T^d$ supremum. We will make this into a precise argument below.

With this plan in mind, we turn to step~\ref{stepthree} of the general outline from~\S~\ref{sec:ideas outline}. We study analytically the coefficient $\beta \in \Q^{\times}$
of $\bx^\bn$ in $F(\bx)$. To approach it, we express~\eqref{G coeff} by a Cauchy integral: 
 \begin{equation} 
 \begin{aligned}
&  \beta \cdot \exp \left( \sum_{k=0}^l  \left( \sum'_{ 
 \frac{d \gamma_k}{m} \leq j < \frac{d\gamma_{k+1}}{m} }   n_j \right)  \log{\varphi_k'(0)}  \right)  \\
 & = [\bz^\bn] \left\{ G(\bz) \right\} \\
 & = \left[z_1^{n_1+M}  \cdots z_d^{n_d+dM}\right] \left\{ W(\bz) G(\bz) \right\}  \\
&  = \int_{\T^d}  \frac{W(\bz) G(\bz)}{z_1^{n_1+M} z_2^{n_2+2M} \cdots z_d^{n_d + dM}} \, \mv(\bz).
 \end{aligned}
 \end{equation}
Consequently, estimating the latter integrand pointwise by the supremum, we derive an analytic upper bound on the nonzero rational number~$\beta \in \Q^{\times}$: 
\begin{equation} \label{the analytic bound}
\begin{aligned}
 |\beta| & \leq  \exp \left( \sup_{\T^d} \left\{ \log{|WG|}  \right\}  - \sum_{k=0}^l  \left( \sum_{ 
 \frac{d \gamma_k}{m} \leq j < \frac{d\gamma_{k+1}}{m} }   n_j \right)  \log{|\varphi_k'(0)|}  \right) \\
 & \leq  \exp \left( \sup_{\T^d} \left\{ \log{|W(\bz) F(\Phi(\bz))|}  \right\}  - \sum_{k=0}^l  \left( \sum_{ 
 \frac{d \gamma_k}{m} \leq j < \frac{d\gamma_{k+1}}{m} }   n_j \right)  \log{|\varphi_k'(0)|}  + O_h(1)  \right). 
\end{aligned}
\end{equation}
Here, since $|\varphi_l'(0)| > 1$ and this is an upper bound, we have legitimately removed the dash proviso in the inner summation over~$j$. 
We further rework~\eqref{the analytic bound} using~\eqref{seed sum est} and apply an Abel summation to obtain 
(recalling for the boundary terms that $\gamma_{l+1} = m$ and $\gamma_0 = 0$),
the following bound on~$\log{|\beta|}$:
\begin{equation}  \label{Cauchy bound}
\begin{aligned}
 & \leq   \sup_{\T^d} \left\{ \log{|W(\bz) F(\Phi(\bz)))|}  \right\}  - \frac{\alpha}{m^2} \sum_{k=0}^l  (\gamma_{k+1}^2-\gamma_k^2) \log{|\varphi_k'(0)|}  + O(\epsilon \alpha) \\
& = \sup_{\T^d} \left\{ \log{|W(\bz)F(\Phi(\bz)))|}  \right\}  - \alpha \log{|\varphi_l'(0)|}  + \frac{\alpha}{m^2} \sum_{k=1}^l  \gamma_k^2 \log{\frac{|\varphi_{k}'(0)|}{|\varphi_{k-1}'(0)|}}  + O(\epsilon \alpha).
\end{aligned}
\end{equation}

At this point we follow~\cite[\S~2.5.14]{UDC} to upper-estimate the supremum term in~\eqref{Cauchy bound}. From~\eqref{F construction}, the
triangle inequality yields as a pointwise upper bound over $\bz \in \T^d$: 
\begin{equation}  \label{inner up}
\begin{aligned}
\log{|F(\Phi(\bz))|} \leq   \max_{ 
\bk/D \in P_{\epsilon}^d \cap \Z^d  } \left\{ \sum_{j=1}^{d} k_j \log{|\Phi_j(z_j)|} \right\} + O(\epsilon \alpha) + o(\alpha),
\end{aligned}
\end{equation}
where the splicing notation for the univariate components of the multivariable map (which we defined in~\S~\ref{seeding} above)
$$
\Phi(\bz) =: \left( \Phi_1(z_1), \ldots, \Phi_d(z_d) \right)
$$
uses 
$\Phi_j(z_j) 
:= \varphi_k(z_j)$ for the unique $k = k(j) \in \{0,\ldots,l\}$ determined by the~$\boldsymbol{\gamma}$ rule spelled out in~\S~\ref{seeding}.   

\subsubsection{The numerical integration}    \label{sec:num int}
Upon infinitesimally scaling down $z \mapsto (1-\varepsilon) z$ the coordinate of the unit disc $\D$, and taking the $\varepsilon \to 0$
limit at the very end, we may and do assume that none of the holomorphic functions $\varphi_0, \ldots, \varphi_l \in \mathcal{O}(\Db)$ have any zeros lying on 
the unit circle~$\T$. 

For ease of notation later, we define $T^d_{\boldsymbol{\gamma}, \epsilon}:=\{\bz\in \T^d : \forall k, D(\bz^{(k)}) <\epsilon\}$.
We denote $d^{(k)}:=\lceil \gamma_{k+1} d/m\rceil - \lceil \gamma_k d/m\rceil$  the length of the~$\bz^{(k)}$ variable block, and
we write~$\bz =: e^{2\pi i \bs}$ in block form~$\bz = (\bz^{(0)},\ldots, \bz^{(l)})$, so that~$\bz \in T_{\boldsymbol{\gamma}, \epsilon}^d$ is tantamount to
having~$D(\bs^{(k)} )<\epsilon$ for every~$k = 0, \ldots, l$.
Given a vector $\mathbf{w} = (w_1, \ldots, w_d) \in \R^d$, let us admit a slight abuse of notation and denote by $\mathbf{w}^* =: (w_1^*,  \ldots, w_d^*)$ the \emph{increasing\footnote{Or rather, nondecreasing.} rearrangement vector} of the component set of~$\mathbf{w}$.
By Corollary~\ref{F stats}, the running condition $\bk/D \in P_{\epsilon}^d$ implies
\begin{equation}  \label{tru 1}
k_j^* = D(j/d) + O(\epsilon D) = (2j/d^2)\frac{\alpha}{m} + O\left(\epsilon \alpha /d \right), \qquad j = 1,\ldots, d.
\end{equation}
Similarly, 
for $\bs^{(k)} \in [0,1)^{d^{(k)}}$ with $D(\bs^{(k)} )<\epsilon$, the increasing rearrangement $(\bs^{(k))})^*$ has components
\begin{equation} \label{tru 2}
(s^{(k)})_\ell^* = \ell /d^{(k)} + O(\epsilon), \qquad \ell = 1,\ldots, d^{(k)}.
\end{equation}

Since all the coordinate functions
$\log{|\Phi_j|} = \log{|\varphi_{k(j)}|}  : \T \to \R$ are of bounded variation,
Koksma's inequality (see, for instance, \cite[\S~2.5.1]{UDC} for
a discussion and further references) implies that for $j\in [\lceil \gamma_k d/m\rceil, \lceil \gamma_{k+1} d/m\rceil )$,
\[\log|\Phi_j(e^{2 \pi i (s^{(k)})_\ell^*})| = \log |\varphi_k(e^{2 \pi i \ell /d^{(k)}})| + O(\epsilon).\]
 Thus, by Koksma's inequality again, we arrive at the definition of the function~\eqref{rearrangement fine}:  
 \[g_{\boldsymbol{\varphi}, \boldsymbol{\gamma}}(j/d)= \log \left|\varphi_k\left(e^{2 \pi i \left(j - \sum_{h=0}^{k-1} d^{(h)}\right)/d^{(k)}}\right)\right| +O(1/d).\] 
  \silentcomment{The $O(1/d)$ comes from taking into account the nearest integer rounding for $j/d$ in the definition of each $\gamma_k$ interval.}
The increasing rearrangement notation then reads: 
\begin{equation}  \label{with inc rearr}
g_{\boldsymbol{\varphi},\boldsymbol{\gamma}}^*(j/d) =  \left( \log \left|\varphi_{k(j)}\left(e^{2 \pi i \left(j - \sum_{h=0}^{k(j)-1} d^{(h)}\right)/d^{(k(j))}}\right)\right|\right)_j^*+O(1/d),  
\end{equation}
where the index~$k = k(j) \in \{0,\ldots,l\}$ is determined by the rule of~\S~\ref{seeding}, which at these arguments reads: $j/d \in [\gamma_k/m, \gamma_{k+1}/m)$.

In summary, we have proved that
\[\left( \Phi_j(z_j) \right)^*_j = g_{\boldsymbol{\varphi},\boldsymbol{\gamma}}^*(j/d) + O(\epsilon) + O(1/d).\]
Koksma's inequality and the rearrangement inequality now
yield a numerical integration estimate: 
\begin{equation}  \label{numest}
\begin{aligned}
 & \bt \in P_{\epsilon}^d, \, t_1 \leq \cdots \leq t_d, \quad \bz \in T^d_{\boldsymbol{\gamma}, \epsilon} <\epsilon \Longrightarrow \\  
 & \sum_{j=1}^{d} 2t_j \cdot \log{|\Phi_j(z_j)|} \leq 
\sum_{j=1}^{d} (2j/d) g^*_{\boldsymbol{\varphi},\boldsymbol{\gamma}}(j/d) + O(\epsilon d) + O(1). 
\end{aligned}
\end{equation}

Therefore~\eqref{inner up} and~\eqref{numest} imply the following
upper estimate: 
\begin{equation}  \label{uni qlo 5}
\begin{aligned}
\sup_{\bz \in T_{\boldsymbol{\gamma}, \epsilon}^d} \left\{ 
\log{|F(\Phi(\bz))|} \right\} &  \leq   \frac{1}{d} \left( \sum_{j=1}^{d}  2(j/d) \cdot g_{\boldsymbol{\varphi},\boldsymbol{\gamma}}^*(j/d)   \right) \frac{\alpha}{m} + O(\epsilon \alpha) + O\left(\frac{\alpha}{d}\right)+ o(\alpha). 
\end{aligned}
\end{equation}
At this point, Koksma's inequality applies yet again to prove the following estimate,
  uniformly on the well-distributed part $T_{\boldsymbol{\gamma}, \epsilon}^d \subset \T^d$: 
  \begin{equation}  \label{uni qlo 2}
\begin{aligned}
\sup_{\bz \in T_{\boldsymbol{\gamma}, \epsilon}^d} \left\{ 
\log{|F(\Phi(\bz))|} \right\} \leq  \frac{\alpha}{m} \cdot  \int_0^1 2t \cdot g_{\boldsymbol{\varphi},\boldsymbol{\gamma}}^*(t) \, dt    + O(\epsilon \alpha) + O\left(\frac{\alpha}{d}\right) + o(\alpha). 
 \end{aligned}
 \end{equation}
  
  \subsubsection{Noise canceling}
  To handle the complementary (badly distributed) part of the integration torus~$\T^d$, we select the ``sufficiently big'' exponent~$M$
  of the damping Vandermondian: 
  \begin{equation}
  M := \left\lfloor  \frac{\sup_{\T} \log{|\Phi|}}{c'(\epsilon)} \frac{\alpha}{d^2}  \right\rfloor,
  \end{equation}
  where $c(\epsilon)$ is the function from~\eqref{W dampens}, and we recall that we have assumed the condition~$d \geq d_0(\epsilon)$ in that lemma. 
  On the poorly distributed part $\T^d \setminus T_{\gamma,\epsilon}^d$ we get: 
\begin{equation} \label{uni qlo 3}
\begin{aligned}
\sup_{\bz \in \T^d \setminus T_{\boldsymbol{\gamma}, \epsilon}^d} \left\{ 
\log{\left|W(\bz)F(\Phi(\bz))\right|} \right\}  = O(\epsilon\alpha) + o(\alpha). 
\end{aligned}
\end{equation}
Putting together \eqref{uni qlo 2},  \eqref{uni qlo 3}, and Lemma~\ref{Fekete}, we derive the uniform estimate 
\begin{equation} \label{uni qlo 4}
\begin{aligned}
\sup_{\bz \in \T} \left\{ 
\log{\left|W(\bz)F(\Phi(\bz))\right|} \right\}  \leq \frac{\alpha}{m} \cdot  \int_0^1 2t \cdot 
g_{\boldsymbol{\varphi},\boldsymbol{\gamma}}^*(t) \, dt  + O_{\epsilon}\left( \frac{\log{d}}{d} \alpha \right) + O(\epsilon\alpha) + o(\alpha). 
\end{aligned}
\end{equation}
\subsubsection{The Cauchy bound}
Our upper bound on the leading~$\bx^\bn$ coefficient~$\beta$ now follows as the combination of~\eqref{Cauchy bound} and~\eqref{uni qlo 4}
\begin{equation}  \label{Cauchys}
\begin{aligned}
\log{|\beta|}  \leq  \frac{\alpha}{m} \cdot  & \int_0^1 2t \cdot 
 g_{\boldsymbol{\varphi},\boldsymbol{\gamma}}^*(t)  \, dt  - \alpha \log{|\varphi_l'(0)|} + \frac{\alpha}{m^2} \sum_{k=1}^l  \gamma_k^2 \log{|\varphi_k'(0)|} \\  &  + O_{\epsilon}\left( \frac{\log{d}}{d} \alpha \right)   + O(\epsilon \alpha) + o(\alpha). 
\end{aligned}
\end{equation}

This asymptotic inequality, upon taking the limits in the order $\alpha \to \infty, d \to \infty$, and $\epsilon \to 0$, already proves the special case
  $\mathbf{b = 0}, \mathbf{e = 0}$ of the theorem, whereby $\beta \in \Z \setminus \{0\}$ is a nonzero rational
 integer and therefore at least one in magnitude. To complete the general case, it remains to estimate the denominator of the leading order coefficient~$\beta \in \Q^{\times}$
 in $F(\mathbf{x})$. 
 
 \subsection{Denominator arithmetic} \label{denominator arithmetic} This is a new aspect which we did not encounter in~\cite{UDC}. We consider all the possible 
 combinations~\eqref{F construction} with $c_{\mathbf{i,k}} \in \Z$, and in those, we estimate prime-by-prime the worst possible denominator
 that may arise in a leading order monomial coefficient~$\beta$, under the premises of Lemma~\ref{Siegel} and the denominator types~\eqref{den type int}. 
 We prove $\exp\left(\alpha \tau^{\flat}(\mathbf{b})  + o(\alpha) \right)$ as  the best-possible (exact) formula  in the $\mathbf{e = 0}$ case. 
 For the general case with added integrals, the exact denominator worst-case analysis seems subtle  ---  especially if
 in the actual $m =14$ case in~\S~\ref{sec:proofA} of our main application one tries to consider the finer denominators we indicate by Remark~\ref{central binomial};  ---  
 but we provide a handy upper estimate which turns out to be the
 quantity $\exp\left(\alpha \tau^{\flat}(\mathbf{b}) + \alpha \tau^{\sharp}(\mathbf{e})  + o(\alpha) \right) = \exp\left(\alpha \tau(\mathbf{b;e}) + o(\alpha) \right)$ of the 
 statement of Theorem~\ref{main:elementary form}. We suspect our estimate to be pretty sharp in the case that we use for the proof of Theorem~\ref{mainA}. 
 
 \subsubsection{A preview on $\tau^{\sharp}$} \label{sec:den arithmetic preview}
 It is plain from the way these growth rates $\tau^{\flat}(\mathbf{b})$ and $\tau^{\sharp}(\mathbf{e})$ are added up that we are separately estimating 
 the extra denominators that the factors $n^{e_i}$ introduce from~\eqref{den type int}. The  parameter $\xi \in [0,m]$ of the definition~\eqref{tau sharp} of $\tau^{\#}(\mathbf{e})$ is used for the cutoff in Lemma~\ref{lcm-2}
 to decide for which primes~$p$ to estimate the added power of~$p$ in $\den(\beta)$ based on the lemma, and for which primes to estimate it based, instead, directly on the remark that  every product $f_{i_1}(x_1) \cdots 
 f_{i_d}(x_d)$ has only a limited number of factors involving ``extra $n$ denominators'' in~\eqref{den type int}: namely, precisely $\left( 
 \sum_{i=1}^{m} e_i \right) d/m$ of the~$d$ factors contribute, if we count with multiplicities $n^{e_i}$.  For the primes $p \leq \xi D$, we use the latter trivial estimate;
  for the range $p > \xi D$, we estimate by Lemma~\ref{lcm-2} using that the leading $\prec$-order exponent vector $\mathbf{n}$ is, upon relabeling the variables, close
 to $(mD/d, 2mD/d, \ldots, dmD/d)$. 
 
 This means looking respectively into the left-hand side and the right-hand side of the
 identity $\sum_{\mathbf{i,k}} c_{\mathbf{i,k}} \mathbf{x^k} 
 f_{i_1}(x_1) \cdots f_{i_m}(x_m)  = \sum_{\bm} \beta_{\bm} \bx^\bm$. In our $d \to \infty$ asymptotic, our $\prec$-leading exponent vector $\mathbf{n} \approx \left( 
 mD/d, 2mD/d, \ldots, dmD/d 
\right)$ realizes once again the cross-variables 
 dimension, 
 with an integration variable  $t := jm/d \in [0,m]$ that 
 leads up to the function $\max(e_i) \cdot I_{\xi}^m(\xi)$ of Definition~\ref{integrated lcm}. More precisely, recalling that we will take~$\epsilon \to 0$ in the end, which will automatically force $\delta \rightarrow 0$, the latter emerges as the $
\int_0^{m+\delta}$ Riemann  integral of  the
 $$
\frac{\max_i (e_i)}{D} \, \log{ \frac{\left[ \max(1, tD - D), tD \right]}{ \mathrm{gcd} \left\{  [1,\ldots,\xi D], [\max(1, tD - D), tD]  \right\}  }  } \, dt
 $$
  estimates
 from Lemma~\ref{lcm-2}. 
 As this latter estimate goes ``across the variables,'' it 
 only ``sees''  the $n^{e_i}$ exponents through their common capping $\max_i (e_i)$, based on the remark that, separately in every variable $x_j$, all the terms $[x_j^{n}] \, f_i(x_j)$ 
 have added denominators multiplying at most by $n^{\max_i(e_i)}$; this has to be taken uniformly in $i \in \{1,\ldots,m\}$ (hence the 
 maximum over~$i$), as each given function species~$f_i$ will occur from some product of~\eqref{F construction} at every single coordinate~$x_j$. 
 To leverage that many of the $f_i$ could have smaller added denominators $n^{e_i}$ than the common
 $n^{\max_i(e_i)}$ capping of this cross-variable denominator estimation, we use a balancing parameter~$\xi$ and directly estimate the
 primes $p \leq \xi D$ from the multiplicity density $\left( 
 \sum_{i=1}^{m} e_i \right) d/m$ of affected factors in each product $f_{i_1}(x_1) \cdots f_{i_m}(x_d)$, in which an extra $p^{e_i}$
 denominator could possibly be hiding (this is a conservative estimate!).

We now execute both points $\tau^{\flat}(\mathbf{b})$ and $\tau^{\#}(\mathbf{e})$ of this $\den(\beta)$ majorization program. 
In these denominator estimates, the essential point is that our $\prec$-minimal exponent $\mathbf{n} \in (m+\delta)D \cdot P_{\epsilon}^d$ in $F(\bx)$ has uniformly distributed components, but 
the corresponding information on $\mathbf{k} \in D \cdot P_{\epsilon}^d$ is now ignored. (It is conceivable that the latter could be also exploited to give
a more precise bound; however, we were unable to do that in our applications at hand.) 

This is the exact opposite to
the archimedean growth estimate in~\S~\ref{equidistribution}. 

\subsubsection{The $\tau^{\flat}(\mathbf{b})$ piece}  \label{flat piece}  Consider 
any of the lowest order exponent vectors
$$
\mathbf{n} = (n_1, \ldots, n_d)  \in (m+\delta)D \cdot P_{\epsilon}^d,
$$
as given by Lemma~\ref{Siegel}. Recall that in~\S~\ref{seeding} we relabeled the coordinates to assume  ---  simply for a notational convenience  ---  that our $\bn$ has nondecreasing components: $n_1 \leq \cdots \leq n_d$. Then, by Corollary~\ref{F stats}, we have
\begin{equation} \label{n den}
\begin{aligned}
n_j & \leq (1+\epsilon) mD(j/d) + 2m\epsilon D \\
& \leq mD(j/d) + 3m\epsilon D, \qquad j =1, \ldots, d. 
\end{aligned}
\end{equation}
   We need to compute the lowest common denominator
of all the nonzero rational numbers $\beta \in \Q^{\times}$ that may arise
as the $\mathbf{x^n}$ coefficient in any product
\begin{equation} 
x_1^{k_1} \cdots x_d^{k_d} \cdot g_{\pi(1)}(x_1) \cdots g_{\pi(d)}(x_d) \in \Q \llbracket \mathbf{x} \rrbracket
\end{equation}
with some arbitrary $\mathbf{k} \in D \cdot P_{\epsilon}^d$, some arbitrary permutation~$\pi$ of~$\{1,\ldots,d\}$, and for each $i \in \{1,\ldots,m\}$, some arbitrary formal functions
\begin{equation} 
g_{1+(i-1)d/m }(x), \ldots, g_{id/m}(x)  \in \bigoplus_{n=0}^{\infty}  \frac{x^n}{ [1,\ldots, b_{i,1} \cdot n] 
\cdots [1, \ldots, b_{i,r} \cdot n] } \, \Z.
\end{equation} 
This means nothing more nor less than the lowest common multiple of all products
\begin{equation}   \label{lcm products}
\prod_{i=1}^{m} \prod_{s=1}^{d/m}  \prod_{h = 1}^{r}  \left[1,\ldots, b_{i,h} \cdot n_{\pi\left((i-1)d/m + s\right)}\right], \quad \pi \in S_{d},
\end{equation}
as~$\pi$ ranges over all permutations of~$\{1,\ldots,d\}$. 

We handle~\eqref{lcm products} with a prime-by-prime determination of the maximizing valuation. 
The following simple lemma is where the  special condition~\eqref{column shape} on the denominators shape matrix~$\bb$ is used,
in all our theorems in~\S\S~\ref{fine section},~\ref{new slopes}. 

\begin{lemma} \label{single step valuationwise}
For every prime~$p$, every vector~$(c_1,\ldots,c_m) \in \NwithzeroA^m$ of the form
\begin{equation} \label{rowc}
0 = c_1 = \cdots = c_u < c_{u+1} = \cdots = c_m =: c,
\end{equation}
and every nondecreasing sequence $n(1)  \leq \cdots \leq n(km)$ consisting of $km$~positive integers, 
the following equality holds: 
 \begin{equation} \label{shapelcm}
 \begin{aligned}
\max_{ \pi \in S_{km} } \quad  & \mathrm{val}_p \left\{ \prod_{i=1}^{m} \prod_{s=1}^{k} \left[ 1,\ldots, c_i \cdot n\left( \pi\left((i-1)k + s \right)\right) \right] \right\}  \\
= \, & \mathrm{val}_p \left\{ \prod_{i=1}^{m} \prod_{s=1}^{k} \left[ 1,\ldots, c_i \cdot n\left( (i-1)k + s\right) \right] \right\}. 
\end{aligned}
 \end{equation}
 In other words, as~$\pi \in S_{km}$
 ranges through all permutations of~$\{1,\ldots,km\}$, the identity permutation~$\pi = \mathrm{id}$ maximizes 
 the $p$-adic valuation in~\eqref{shapelcm}. 
\end{lemma}

\begin{proof} 
The condition~\eqref{rowc} simplifies the requisite product~\eqref{shapelcm} to
 \begin{equation}  \label{simpleshape}
 \prod_{i=u+1}^{m} \prod_{s=1}^{k} \left[ 1,\ldots, c \cdot  n\left( \pi\left((i-1)k + s \right)\right)  \right]. 
 \end{equation}
 The lowest common multiple~$[1,\ldots,N]$ of the first~$N$ positive integers has $p$-adic valuation 
 equal to~$\big\lfloor \frac{\log{N}}{\log{p}} \big\rfloor$, and so the quantity in~\eqref{shapelcm}
 under maximization is exactly equal to
 \begin{equation}  \label{simplesum}
 \sum_{i=u+1}^{m} \sum_{s=1}^{k} \Big\lfloor \frac{\log{c}}{\log{p}} + \frac{\log{n\left( \pi\left((i-1)k + s \right) \right)}}{\log{p}} \Big\rfloor. 
 \end{equation}
 From the $km$~positive integers~$\{n(1), \ldots, n(km)\}$, we have to pick~$k(m-u)$ with pairwise distinct indices 
 to maximize the sum~\eqref{simplesum}. Clearly this is maximized by picking the~$k(m-u)$ largest available numbers~$n(\bullet)$, so in 
 particular our monotonicity assumption on~$n(\bullet)$ gives that~\eqref{simplesum} is maximized by the identity permutation~$\pi = \mathrm{id}$, with maximum
 \[
 \sum_{j=ku+1}^{km} \Big\lfloor \frac{\log{c}}{\log{p}} + \frac{\log{n(j)}}{\log{p}} \Big\rfloor.  \qedhere
 \]
\end{proof}

Applying~\eqref{n den} on the nondecreasing sequence $\bn$, together with our condition~\eqref{column shape} on the $m \times r$ array,  
we find by Lemma~\ref{single step valuationwise} that as soon $D \gg_{\epsilon} 1$, all the lowest common multiple products~\eqref{lcm products} divide
\begin{equation} \label{lcm cap}
\prod_{i=1}^{m} \prod_{s=1}^{d/m}  \prod_{h = 1}^{r}  \left[1,\ldots, b_{i,h} \cdot \big(  (i-1)D + smD/d
\big) + \epsilon B D  \right],
\end{equation}
where the constant~$B := m \cdot \max_{i,h}\{b_{i,h}\}$. 

 By the prime number theorem, the lowest common multiple cap~\eqref{lcm cap} evaluates in the $D \to \infty$ asymptotic to
\begin{equation}
\begin{aligned}
& \exp \left(    \sum_{i=1}^{m} \sum_{s=1}^{d/m}  \sum_{h = 1}^{r} \left( b_{i,h} \cdot \big(  (i-1)D + smD/d
 \big)  + O(\epsilon D) \right)  \right)  \\  
& =  \exp \left(   mD  \sum_{i=1}^{m} \sum_{s=1}^{d/m} \sigma_i \cdot \big(  (i-1)/m + s/d
\big)  + O(\epsilon d D)   \right)  \\    \label{flat est}
& =  \exp \left(    \alpha  \sum_{i=1}^{m} \sigma_i \int_{(i-1)/m}^{i/m} 2t \, dt  + O(\epsilon \alpha) + o_{d \to \infty}(\alpha)  \right) \\ 
& = 
 \exp \left(   \alpha \tau^{\flat}(\mathbf{b})  +  O(\epsilon \alpha) + o_{d \to \infty, \epsilon \to 0}(\alpha)  \right) ,
\end{aligned}
\end{equation}
recalling our definition~\eqref{alpha order} of the vanishing order parameter~$\alpha  = mdD/2$. 

\begin{remark}  \label{genbcomplicated}
The statement of Lemma~\ref{single step valuationwise} ceases to be true if the condition~\eqref{rowc} is relaxed to an arbitrary monotonic~$0 \leq c_1 \leq \cdots
\leq c_m$.  Thus, with~$\tau^{\flat}(\bb) = \frac{1}{m^2} \sum_{i=1}^{m} (2i-1) \sigma_i$  as the definition in~\eqref{taub}, the proof of the theorem would
no longer hold if we relaxed the crude capping~\eqref{column shape} of our denominator types to an arbitrary
matrix~$\bb$ having columns with nondecreasing components.~\endofremark
\end{remark}

\begin{remark} \label{finedenominators}
Unlike for the archimedean growth estimate in~\S~\ref{sec:num int}, our computation here ignored the uniform distribution constraint $\mathbf{k} \in D \cdot P_{\epsilon}^d$ inside the trivial estimate $\bk \in [0,D]^d$. This was how the growth rate~$\tau^{\flat}$ was defined, not to take account of the distribution of the exponents~$\bk$ of the auxiliary polynomials; for this definition, it is an exact computation. 
 In contrast, it was crucial for the horizontal integration idea to exploit the uniformly distributed components of the $\prec$-leading $\bx = \mathbf{0}$ exponent $\mathbf{n} \in (m+\delta)D \cdot P_{\epsilon}^d$.

  In principle (or in practice), upon calculating a denominator rate still more involved than the term~$\QQQQ_N(\bb,\epsilon,d,\varepsilon)$ in Theorem~\ref{high dim BC convexity},  one could formulate a version of Theorem~\ref{main:elementary form} in which~$\tau(\bb;\be)$ is formally refined to a 
complicated limiting formula that does also takes account of the uniform components restriction on the exponent vectors~$\bk \in [0,D]^d$ in the make-up of the auxiliary function~\eqref{F construction}; and where denominator shapes still finer than our template form~\eqref{den type integrated} could be considered. 
The denominators in Remark~\ref{central binomial}, and similar forms involving products of binomial coefficients or products of primes from restricted intervals, are the typical example to have in mind for prospective applications; for deeper studies and more complicated examples, cf.~\cite{ViolaEulerIntegrals,RhinViolaZeta2,Sorokin,ZudilinHypergeometricTales,DaugetZudilin}. In the situation of our application to our main Theorem~\ref{mainA}, we 
did not succeed in making any (non-negligible) use of the uniformly distributed~$\bk$ for handling the more restricted denominator types of Remark~\ref{central binomial}.~\endofremark
\end{remark}

\subsubsection{The $\tau^{\#}(\mathbf{e})$ piece}  \label{sharp piece}  To estimate the denominator surplus from the
extra integration denominators $n^{\mathbf{e}}$, we will separately (as an upper bound) multiply the 
principal denominators cap 
\eqref{lcm cap} by the lowest common denominator of all the possible $\mathbf{x^n}$ coefficients $\beta \in \Q^{\times}$ of all  possible products
$$
x_1^{k_1} \cdots x_d^{k_d} \cdot g_{\pi(1)}(x_1) \cdots g_{\pi(d)}(x_d) \in \Q \llbracket \mathbf{x} \rrbracket,
$$
across all possible 
$\mathbf{k} \in D \cdot P_{\epsilon}^d$, some arbitrary permutation $\pi \in S_d$, and, for each $i \in \{1,\ldots,m\}$, arbitrary formal functions
\begin{equation} \label{generic h}
g_{1+(i-1)d/m}(x), \ldots, g_{id/m}(x)  \in \bigoplus_{n=0}^{\infty}  \frac{x^n}{ n^{e_i}  } \, \Z.
\end{equation} 
This multiplies separate local estimations of the highest possible power of a denominator at every prime~$p$. Consider 
$\xi \in [0,m]$ the parameter of the definition~\eqref{taue}. We estimate differently the 
cases $p \leq \xi D$ and $p > \xi D$.  Firstly we collect two basic standard facts, the first of which is a version of the 
prime number theorem, and the second, an immediate consequence: 

\begin{itemize}
\item[(a)] {\it The product of the primes $p \leq n$ is asymptotic to $\exp(n + o(n))$.}

\medskip 

\item[(b)] {\it  The product of the proper prime powers $p^a \leq n, a \geq 2$,
is bounded by $\exp\left( O(\sqrt{n} \right) =  \exp(o(n))$. }

\medskip

Together, they imply: 

\medskip

\item[(c)] {\it The lowest common multiple $[1,\ldots,n] = \exp( n+o(n))$.}
\end{itemize}

\medskip

These properties prove that  for the $\mathbf{x^n}$ coefficient denominator of the ``generic'' 
\begin{equation} \label{partial clearing}
\left[ 1, \ldots, \xi D \right]^{ \lfloor d \left( \sum_{i=1}^m e_i \right) /m \rfloor} \cdot x_1^{k_1} \cdots x_d^{k_d} \cdot g_{\pi(1)}(x_1) \cdots g_{\pi(d)}(x_d) \in \Q \llbracket \mathbf{x} \rrbracket,
\end{equation} 
we have:
\begin{enumerate}[label=(\roman*)]
\item
the totality of the primes $p \leq \xi D$ add only a negligible $\exp(o(\xi D)) = \exp(o(\alpha))$ factor to the denominators of~\eqref{partial clearing};
\item the clearing factor 
\[
\begin{aligned}
\left[ 1, \ldots, \xi D \right]^{ \lfloor d \left( \sum_{i=1}^m e_i \right) /m \rfloor} & = \exp\left(  \xi d D  \left( \sum_{i=1}^m e_i \right) /m + o(\alpha) \right) \\
& =  \exp\left(  \xi  \left( \sum_{i=1}^m e_i \right) \cdot 2\alpha/m^2 + o(\alpha) \right). 
\end{aligned}
\]
\end{enumerate}
 It is clear then that, up to an $\exp(o_{d \to \infty, \epsilon \to 0}(\alpha))$ factor, the lowest common denominator of all the $\mathbf{x^n}$ coefficients 
 of all the formal expressions~\eqref{partial clearing}, as the $h_j(x)$ range over~\eqref{generic h}, $\mathbf{n}$ ranges over $(m+\delta)D \cdot P_{\epsilon}^d$, and 
 $\mathbf{k}$ ranges over $D \cdot P_{\epsilon}^d$, is a divisor of the lowest common denominator
 of all formal expressions 
\begin{equation*}
\frac{
\left[ 1, \ldots, \xi D \right]^{\infty}}{ (n_1 - k_1)^{\max_i(e_i)} \cdots (n_d - k_d)^{\max_i(e_i)} },  \quad \mathbf{n} \in (m+\delta)D \cdot P_{\epsilon}^d, \quad \mathbf{k} 
\in D \cdot P_{\epsilon}^d;
\end{equation*}
and that, again up to an  $\exp(o_{d \to \infty, \epsilon \to 0}(\alpha))$ factor, this is a divisor\footnote{Even upon including all $\bk \in [0,D]^d$, that is once again ignoring the
$\bk \in P_{\epsilon}^d$ constraint.} of
\begin{equation}  \label{horizontal prod}
\prod_{j=1}^{d} \prod_{ \substack{ \textrm{primes } p > \xi D:  \\
p  \textrm{ divides some} \\
\textrm{positive integer in } \\ \left[mD(j/d) - D,  mD(j/d) \right] }}  p^{\max_i(e_i)}. 
\end{equation}
By Lemma~\ref{lcm-2}, if $m(j/d) >1$, the inner product in~\eqref{horizontal prod} is asymptotic to 
the exponential of
\begin{equation*}
\begin{aligned}
  &  \left( \max_i e_i \right) D  \left( \sum_{h=1}^{\lfloor(m(j/d)-1)/\max(1,\xi) \rfloor} 1/h  \right) \\
  + &  \left( \max_i e_i \right)  D
\max \left\{ \frac{m(j/d)}{ \lfloor  (m(j/d)  +\max(0, \xi -1))/\max(1,\xi) \rfloor} - \xi, 0 \right\}  + o(D).
\end{aligned}
\end{equation*}
In the case $m(j/d) \leq 1$,  the inner product in~\eqref{horizontal prod} is asymptotic to 
the exponential of
\[ \left( \max_i e_i \right) D \max\{0, m(j/d)-\xi\} + o(D).\]
Hence, recollecting our Definition~\ref{integrated lcm} of the integrated LCM cost function $I_u^v(w)$, 
we find that as $d \to \infty$, so that the discrete variable $t := m(j/d)$ converges
to  the  continuous Lebesgue measure of the segment $[0,m]$, the horizontal integration computes the asymptotic
full denominator product~\eqref{horizontal prod} to the following, up to an $\exp(o_{d \to \infty, \epsilon \to 0}(\alpha))$ factor: 
\begin{equation}
\begin{aligned}
\exp \left( (dD/m)  \left( \max_i e_i \right)   \cdot I_{\xi}^{m}(\xi) + o_{d \to \infty, \epsilon \to 0}(dD) \right) \\  = 
\exp \left( (2\alpha/m^2)  \left( \max_i e_i \right)   \cdot I_{\xi}^{m}(\xi) + o_{d \to \infty, \epsilon \to 0}(\alpha) \right). 
\end{aligned}
\end{equation}
All in all, we obtain for any $\xi \in [0,m]$ the upper estimate
\begin{equation}   \label{added estimate}
\exp \left(  (2\alpha/m^2)  \cdot \left(  \xi \sum_{i=1}^{m} e_i +   \left( \max_{1 \leq i \leq m} e_i \right)   \cdot I_{\xi}^{m}(\xi)  \right)   + o(\alpha) + o_{d \to \infty, \epsilon \to 0}(\alpha) \right)
\end{equation}
on the addition to the denominator from the $n^{\mathbf{e}}$ factors in~\eqref{den type int} in any 
leading order coefficient of our auxiliary function $F(\mathbf{x})$ in Lemma~\ref{Siegel}. 

{\it
Our definition of  the rate $\exp \left( \alpha \cdot \tau^{\sharp}(\mathbf{e}) + o(\alpha) \right)$ is as the minimum of
the total added denominators estimate~\eqref{added estimate} over the parameter $\xi \in [0,m]$. }

\subsection{Proof of Theorem~\ref{main:elementary form}}

 We combine the upper bound~\eqref{Cauchys} on the leading~$\bx^\bn$ coefficient~$\beta \in \Q^{\times}$ with the added up upper estimates that we computed in~\S~\ref{flat piece} and~\S~\ref{sharp piece} on the denominator $\den(\beta) \in \NwithoutzeroA$ of that coefficient. 
 The latter give: 
 \begin{equation}  \label{arithmetic lo}
 \log{|\beta|}  \geq - \alpha \cdot \left(  \tau^{\flat}(\mathbf{b}) + \tau^{\sharp}(\mathbf{e}) \right)   + o_{d \to \infty, \epsilon \to 0}(\alpha).
 \end{equation}
 The former simplifies to: 
 \begin{equation}  \label{analytic up}
\begin{aligned}
\log{|\beta|}  \leq  \frac{\alpha}{m} \cdot & \int_0^1 2t \cdot 
 g_{\boldsymbol{\varphi},\boldsymbol{\gamma}}^*(t) \, dt - \alpha \log{|\varphi_l'(0)|} \\
 & + \frac{\alpha}{m^2} \sum_{k=1}^l  \gamma_k^2 \log{\frac{|\varphi_k'(0)|}{|\varphi_{k-1}'(0)|}}   + o_{d \to \infty, \epsilon \to 0}(\alpha).
\end{aligned}
\end{equation}
 
 The combination of~\eqref{arithmetic lo} and~\eqref{analytic up} sifts out in the $\alpha \to \infty, d \to \infty, \epsilon \to 0$ limit to
 \begin{equation}
 m  \leq \frac{ \int_0^1 2t \cdot 
 g_{\boldsymbol{\varphi},\boldsymbol{\gamma}}^*(t) \, dt +  \frac{1}{m} \sum_{k=1}^l  \gamma_k^2 \log{\frac{|\varphi_k'(0)|}{|\varphi_{k-1}'(0)|}}  }{ \log{|\varphi_l'(0)|} - \tau^{\flat}(\mathbf{b}) - \tau^{\sharp}(\mathbf{e})  },
 \end{equation}
 which is precisely our claimed holonomy bound.      $\hfill{\square}$

 \bigskip 
 
\emph{At this point, a reader primarily interested in the proof of Theorems~\ref{mainA} and~\ref{logsmain} can skip directly ahead to~\S~\ref{sec:YtoY0(2)} on a first reading. }

 \subsection{Completion of the proof of Theorem~\ref{thm three elements}}   \label{sec:completion three elements}

  In~\S~\ref{bivalent app}, by a direct application of Theorem~\ref{basic main}, we already proved the
 property~$(\ast)$ in Theorem~\ref{thm three elements} towards the arithmetic characterization of the $\log^2(1-x)$ function.  This was the case  ---  the one of relevance to the sample application to $\Q$-linear independence
 proofs that we gave in~\S~\ref{mixed examples}  ---  that the minimal order differential operator~$\cL$ has an essential singularity at the ``fourth puncture'' $x=\delta$ from the statement of the theorem. With the feature of the multiple maps~$\boldsymbol{\varphi}$ in Theorem~\ref{main:elementary form}, we can now complete the proof of the full Theorem~\ref{thm three elements} by handling the case that~$x = \delta$ is at most an apparent singularity of~$\cL$.

 \begin{proof}[Proof of Theorem~\ref{thm three elements}]
 From the discussion in~\S~\ref{bivalent app}, as the setup of the theorem-under-proof implies that our type~$[1,\ldots,n][1,\ldots,n/2]$ formal function $f(x) \in \Q\llbracket x \rrbracket$
 has a meromorphic pullback under $\varphi(z) := \frac{8(z+z^3)}{(1+z)^4}$, where already $|\varphi'(0)| = 8 > \tau = 3/2$, we certainly get the existence of a minimal-order nonzero 
 linear differential operator~$\cL$ over $\Q(x)$ satisfying $\cL(f) = 0$. 
 It remained to cover the case that the linear ODE~$\cL(F) = 0$ has a full set of \emph{meromorphic} solutions in a neighborhood of~$x = \delta$. Upon multiplying by a nonzero polynomial~$Q \in \C[x] \setminus \{0\}$ to clear up the possible meromorphic poles, this is equivalent to the assumption that the holomorphic function germ~$Q(x)f(x) \in \C\llbracket x \rrbracket$ is analytically continuable as a \emph{holomorphic} function along all paths in~$\P^1 \setminus \{0,1,\infty\}$. By Proposition~\ref{overconvergence}, this condition in turn furnishes a meromorphic pullback~$\varphi^*f = \left(\varphi^*(Qf)\right)/\varphi^*Q \in \mathcal{M}(\Db)$ under all holomorphic mappings 
 $\varphi : \Db \to \C \setminus \{1\}$ subject to~$\varphi^{-1}(0) = \{0\}$.
 
 The reason that the previous argument breaks down in this case is that, in the absence of the fourth singularity $\delta \notin \{0,1,2,\infty\}$, there is no longer
 a reason for the~$\Q(x)$-linear independence of~$f\left( x/(x-1) \right)$ from the four other functions in~\eqref{li11 case 1}, and we only have~$m=4$ with the functions
 \begin{equation} \label{basic four}
 \begin{aligned}
 f_1(x) := 1, \quad f_2(x) := \log(1-x), \quad f_3(x) := \log^2(1-x), \quad f_4(x) := f(x), 
 \end{aligned}
 \end{equation}
 and the choice of denominators type given by $\frac{1}{2}\bb_0$ and~$\be = \mathbf{0}$ of Example~\ref{ex:li11}. 
  On the other hand, the absence of~$\delta$ spares us the trouble to have the analytic mapping~$\varphi$ to necessarily cover the value~$\delta$ only once, and we can take a completely arbitrary~$\varphi : \Db \to \C \setminus \{1\}$ subject only to~$\varphi^{-1}(0) = \{0\}$.
 
 Suppose for the contradiction that there exists a fourth $\Q(x)$-linearly independent function~$f(x)$ in~\eqref{basic four} still of the type~$[1,\ldots,n/2][1,\ldots,n]$, and therefore completing the combined type~$\frac{1}{2}\bb_0$ from Example~\ref{ex:li11}. 
 With~$S = T = \emptyset$ and~$\HH$ taken as the four-dimensional~$\Q(x)$-linear span of~\eqref{basic four}, we use the symmetrization dictionary~$\varphi_{Y(2)} \rightsquigarrow \varphi_{Y_0(2)}$ described in Basic Remark~\ref{equivalence} in the section~\S~\ref{sec:YtoY0(2)} on the~$y := x + x/(x-1) = x^2/(x-1)$ descent, with technical details given by Lemma~\ref{etalecover} (on the algebraic and arithmetic sides) and Lemma~\ref{analytic quotient} (on the analytic side).  Explicitly, using the involution~$w(x) := x/(x-1)$ and the symmetrization
 coordinate~$y := x+w(x) = xw(x) = x^2/(x-1)$, we have for $\HH^{w=1}$ the four-dimensional~$\Q(y)$-vector space of the denominator type denoted~$\bb_0$ in Example~\ref{ex:li11}, with~$\be = \mathbf{0}$ (no added integrations), and spanned by~$1$, $\sqrt{y(4-y)} \, \arcsin{\frac{\sqrt{y}}{2}}$, $\left( \arcsin{\frac{\sqrt{y}}{2}} \right)^2$, and the symmetrizations of~$f(x)$. 
In the notation of Basic Remark~\ref{equivalence}, where in particular~$h = \lambda^2/(\lambda-1) = -256q + \cdots$ denotes a hauptmodul~\eqref{defofh} of~$Y_0(2)$ written out in the coordinate~$q := e^{2\pi i \tau}$, we 
select for our ambiance the analytic mapping
\begin{equation} \label{app:Y02g}
\varphi_{Y_0(2)} := h \circ \HHH(1/2,10,3)  \in \mathcal{O}(\Db),
\end{equation}
where~$\HHH(1/2,10,3) : \Db \hookrightarrow \D$ is the domain described in~\S~\ref{gobble contours}. Lemma~\ref{confgobble} computes~$|\HHH'(1/2,10,3)(0)| = 198/505$ for the conformal radius of this domain, giving
\begin{equation}
|\varphi_{Y_0(2)}'(0)| = 256 \cdot \frac{198}{505} = \exp\left( 4.608886\ldots \right)
\end{equation}
for the conformal size of our ambient analytic map. It is usefully large in comparison to the denominators growth rate~$\tau(\bb_0) = 21/8 = 2.625$ from Example~\ref{ex:li11}. 

Here Corollary~\ref{stacky overconvergence}, 
as interpreted by Basic Remark~\ref{equivalence} and once again using a suitable polynomial multiplier~$Q \in \C[y] \setminus \{0\}$ to clear all the possible meromorphic poles, proves that our choice~\eqref{app:Y02g} resolves analytically 
the four-dimensional holonomic $\Q(y)$-module $\HH^{w=1}$:
 \begin{equation}
\dim_{\Q(y)} \HH^{w=1} = 4, \qquad \varphi_{Y_0(2)}^* \HH^{w=1} \subset \mathcal{M}(\Db). 
\end{equation} 
This is all conditional on the supposed falsity of the theorem under proof. It is to this $\Q(y)$ situation that we apply Theorem~\ref{main:elementary form}, with 
$m := 4$ and the following choices of the intermediate maps~$\varphi_k$ and division parameters~$\gamma_k$: 
\begin{equation}
\begin{aligned}
l  & := 2;  \gamma_1 := 3/5,  \gamma_2 := 2; \\
 \varphi_0(z) &  := \varphi_{Y_0(2)} \left( e^{-5} z \right), \\
  \varphi_1(z)  & := \varphi_{Y_0(2)} \left( e^{-1/2} z \right),  \\
 \varphi_2(z)  & := \varphi_{Y_0(2)}(z). 
\end{aligned}
\end{equation}
\texttt{Mathematica} then yields a holonomy quotient~\eqref{fine new bound elementary} of slightly under~$< 3.9$, which is the desired contradiction
to complete the proof of the~$\Q(x)$-linear independence of the four original functions~\eqref{basic four}.

Finally, the integral refinement over~$\Q\left[ x, \frac{1}{x}, \frac{1}{1-x} \right]$ follows at this point from the Hermite--Lindemann--Weierstrass and Mahler theorems in transcendence, exactly as in the proof of Theorem~\ref{logcharacterization}; and then the upgrade from$~\Q\left[ x, \frac{1}{x}, \frac{1}{1-x} \right]$
to~$\Q\left[ x, \frac{1}{1-x} \right]$ follows from our strict~$[1,\ldots,n][1,\ldots,n/2]$ denominators requirement exactly as
in point~\ref{no x0} in the proof of Theorem~\ref{logcharacterization}. 
\end{proof}

\begin{example}
With~$m=3$ and the true functions 
\[1, \sqrt{y(4-y)} \, \arcsin{\frac{\sqrt{y}}{2}}, \left( \arcsin{\frac{\sqrt{y}}{2}} \right)^2\]
 now using the type
\begin{equation*}
\mathbf{b} = \left( \begin{array}{cc} 0 &  0 \\ 0 & 2 \\ 1 & 2   \end{array} \right), \qquad \mathbf{e} = (0, 0,  0),  
\qquad \tau(\bb;\be) = \frac{7}{3},
\end{equation*}
a short numerical experimentation, which we have not attempted to make rigorous, suggests that for~$l = 2$ (two division points) and maps of the form 
$$
\varphi_{Y_0(2)}(z) := h \left( \HHH(r,e,f)(\delta z) \right), \quad \varphi_k(z) := \varphi_{Y_0(2)}\left( \gamma_k z \right), \quad k = 0, 1, 2, 
$$ 
the minimizing holonomy bound on the three functions should probably be attained at (for example) about the choice
\begin{equation*}
\begin{aligned}
(r,e,f) &  \approx (0.55,\infty,5); \quad \left( \gamma_1, \gamma_2 \right)   \approx (0.19,0.65),  \\
(r_0, r_1) &  \approx \left( e^{-4.3}, e^{-0.76} \right),  \quad \delta   \approx 0.77,
\end{aligned}
\end{equation*}
with holonomy quotient value~\eqref{fine new bound elementary} being at~$\approx 3.239$.  \endofremark
\end{example}

\section{Convexity in Bost's slopes method} \label{new slopes}

We begin with the following clean refinement of Theorem~\ref{basic main}, which finally we prove in this section by a single variable method based 
on~\cite[\S~5]{BostCharles}. This simple result by itself suffices for our application to Theorems~\ref{mainA} and~\ref{logsmain}, although (in the case of the former) only by the narrowest of margins. The tenor of this section, whose main results are stated in~\S~\ref{sec:convex} after a short introduction, is what we are calling
the \emph{convexity input} that leads up to sharpened holonomy bounds. The improvements are usually fairly small, but they are significant enough to comfortably pass the numerical margin in the requisite numeric in the proof of Theorem~\ref{mainA}, and thereby make fully convincing the Arakelov theory path to the irrationality proof of $L(2,\chip)$. 

\begin{thm} \label{main:BC form} 
With the same standing assumptions of Theorem~\ref{main:elementary form}, 
consider a holomorphic mapping $\varphi : (\Db, 0) \to (\C,0)$ with derivative ({\it conformal size}) satisfying the condition
\begin{equation} \label{stronger positivity condition}
 \log{|\varphi'(0)|} > \max\{\sigma_m, \tau(\mathbf{b;e})\}. 
 \end{equation}
Suppose there exists an $m$-tuple $f_1, \ldots, f_m \in \Q \llbracket x \rrbracket$ of $\Q(x)$-linearly independent formal functions with denominator types of the form
\begin{equation*}   \label{den type integrated}
f_i(x) = a_{i,0}+ \sum_{n=1}^{\infty} a_{i,n} \frac{x^n}{n^{e_i} [ 1, \ldots, b_{i,1} \cdot n] \cdots [1,\ldots, b_{i,r} \cdot n]}, \qquad a_{i,n} \in \Z, 
\end{equation*}
 such that $f_i(\varphi(z)) \in \C \llbracket z \rrbracket$ is the germ of a meromorphic function on $|z| < 1$ for all $i = 1, \ldots, m$. 
 Then 
\begin{equation}\label{BCbound} 
m  \leq  \frac{ \iint_{\T^2} \log|\varphi(z)-\varphi(w)| \, \mv(z) \mv(w) }{  \log{|\varphi'(0)|} - \tau(\mathbf{b;e}) }. 
\end{equation}
In particular, the formal functions $f_1, \ldots, f_m$ are holonomic. 

If moreover all functions $f_i$ are \emph{a priori} assumed to be holonomic, the condition~\eqref{stronger positivity condition} 
can be relaxed to $|\varphi'(0)| > e^{\tau(\mathbf{b;e})}$. 
\end{thm}

When~$\mathbf{e}=\mathbf{0}$, then (as previously
noted after
the statement
of Theorem~\ref{main:elementary form}) $\tau(\mathbf{b;e})$
coincides with the~$\tau(\mathbf{b})$
of Theorem~\ref{basic main}. Hence
Theorem~\ref{basic main}
is an immediate corollary
of Theorem~\ref{main:BC form}.

As discussed in~\S~\ref{sec:algebraic}, the $\mathbf{b = 0}$, $\mathbf{e=0}$ case of $\Z\llbracket x \rrbracket$ functions was established by Bost and Charles~\cite[Corollary~8.3.5]{BostCharles}. 
Charles explained to us that their method can be modified to take denominators into account and obtain the following weaker form\footnote{Conditional, as always, on the positivity of the denominator.} of~\eqref{BCbound}
for a starting bound: 
\begin{equation} \label{starting dens}
m  \leq  \frac{ \iint_{\T^2} \log|\varphi(z)-\varphi(w)| \, \mv(z) \mv(w) }{  \log{|\varphi'(0)|} -(\sigma_m +  \max_{1 \leq i \leq m}  e_i )}.
\end{equation}

The basic idea of Theorem~\ref{main:BC form} is to use the same archimedean estimate as in \cite{BostCharles}, but incorporate into it a closer nonarchimedean evaluation height counterpart that, sufficiently for our purposes in the paper, improves the denominator of the initial bound~\eqref{starting dens}. Remarkably, despite 
two seemingly 
rather different methods being used (the concentration of measure method exploiting multivariable approximations in~\S~\ref{fine section}, and a single variable evaluation heights scheme in the present section), the final denominator
term is exactly the same in both Theorems~\ref{main:elementary form} and all the results in the present section. The measure concentration method is nevertheless theoretically more precise in the general denominators aspect, even though no difference is made to any case of relevance to this paper. In the next section~\S~\ref{slopes}, we formulate the most precise of our holonomy bounds by blending together the measure concentration feature of~\S~\ref{fine section} with the Bost--Charles feature of the present~\S~\ref{new slopes}. 

 For the present section and the next, we use Bost's slopes inequality in Arakelov theory. A practically equivalent framework would be the elementary dynamic box principle of \S~\ref{app:PerelliZannier}; we stick to the former choice for variety in our paper, and because the archimedean evaluation height estimate requires in any event the appeal to some relatively deep theorems in Arakelov theory. We hope that the elementary evaluation heights arrangement in \S~\ref{app:PerelliZannier} could nevertheless be helpful to some readers as an introduction to the tenor of the more elaborate method  (due to Bost) that we take up here.

We  also note that one can more intrinsically formulate the proof of~\eqref{BCbound} in terms of Bost's theta invariants $h_{\theta}^0$ and~$h_{\theta}^1$ as in~\cite{BostBook, BostCharles}. The latter pursue the concept of absolute dimension for the ``space of global sections'' of a (countably) infinite-dimensional Hermitian vector bundle over an arithmetic curve, under the traditional analogy between number fields and algebraic curves over finite ground fields. We do not pursue this approach here as the subsequent ``convexity improvements'' of the archimedean growth term in~\eqref{BCbound} seem to be more of an analytic than a geometric nature.

\subsection{Improvements from convexity}  \label{sec:convex} 
In \S~\ref{iii coeff}, the dynamic pigeonholing involves a certain upper bound~\eqref{Cauchy est} on the interval of possibilities for each Taylor coefficient of a power series $F(x) = \sum_{i=1}^{m} Q_i(x) f_i(x) \in \Q\llbracket x\rrbracket$ in the range of the auxiliary evaluation map, given all previous coefficients of
that power series. This upper bound is the result of estimating a Cauchy contour integrand~\eqref{Cauchy int} over~$\T$ to express the $z^n$ coefficient of the $x = \varphi(z)$ pullback of an element of $\psi_D\left(E_D^{(n)} \right) \subset x^n \cdot \Q\llbracket x \rrbracket$. As in Theorem~\ref{main:elementary form}, there is no particular reason to stick to a single fixed analytic map~$\varphi$ for each filtration layer~$n$ in this set of analytic estimates. Notably, depending on~$n/D$, the ambient map~$\varphi$, and the choice of Euclidean metric structure in the evaluation module~$E_D$ (cf. \S~\ref{slopes review} below), there is a certain optimal choice of an intermediate radius~$r = r(n) \in (0,1]$ for estimating this \emph{$n$-th archimedean evaluation height} analytically via the pullback by $x = \varphi(r(n)z)$; the only difference is that now, in the single-variable analysis, the integration procedure over~$n/D$ is ``vertical'' along the vanishing order, rather than ``horizontal'' across variables. The choice of~$r(n)$ has a geometric significance with convex hulls; incidentally giving another nuance to the name \emph{slopes method}. It aligns with the well-known fact that the Nevanlinna growth characteristic, and various cognate quantities, are convex increasing functions of the logarithm of the radius. 

In this section, we compute these optimal choices~$r(n)$ for two types of Euclidean metrics in the evaluation module: the \emph{Bost--Charles metric} from~\cite{BostCharles}, and an explicit family of \emph{binomial metric weights} $\lambda t^r + \mu t$ depending on three real parameters~$(r,\lambda,\mu)$, which are better amenable to numerical computation and still tend in practice to return close-to-optimal bounds for the best triple~$(r,\lambda,\mu)$. With either of these variations, the results of this section alone (which are independent of~\S~\ref{concentration of measure} and~\S~\ref{fine section}, see also~\S~\ref{sec:leitfaden}) lead to a proof of Theorems~\ref{mainA} and~\ref{logsmain}. We spell them out in~\S~\ref{sec:BC convexity} and~\S~\ref{sec: easy convexity} next, and prove them in~\S~\ref{sec_BCconvexity} and~\S~\ref{sec: proof easy convexity} after preparations in~\S~\ref{slopes review}. Along the way, the proof of the more basic Theorem~\ref{main:BC form} appears in~\S~\ref{one variable slopes}. 
After that, in~\S~\ref{sec_BCfull}, we discuss a further improvement that lines up with --- and theoretically\footnote{At least in all the cases that we encountered in practice; see Remark~\ref{}} strengthens --- Theorem~\ref{main:elementary form} of the previous section. 

\subsubsection{The Bost--Charles characteristic} \label{sec:BC convexity}
Inspired by~\cite[\S~5]{BostCharles}, we introduce a (doubled) Nevanlinna-type growth characteristic, sticking for simplicity to the
case of relevance here of \emph{holomorphic} (rather than meromorphic) disc maps $\Db \to \C$. 
Crucially for this approach, there turns out to be an interpretation of this growth characteristic as an arithmetic intersection number in the Bost--Charles theory. 

\begin{df} \label{BC characteristic}
For a nonconstant holomorphic\footnote{In the general meromorphic case, which we will not consider here, a suitable polar counting term would have to be added. } mapping $\varphi : \Db \to \C$, define the \emph{Bost--Charles characteristic function}
\begin{equation} \label{BCchar}
\hT (\cdot, \varphi) : (0, 1] \to \R, \qquad
\hT(r,\varphi) := \iint_{\T^2} \log|\varphi(z)-\varphi(rw)| \, \mv(z) \mv(w). 
\end{equation} 
\end{df}

As with the usual Nevanlinna and Ahlfors--Shimizu characteristics, see~\cite[\S~3.3.5]{Nevanlinna} or~\cite[Remark~13.3.8]{BombieriGubler}, we have: 

\begin{lemma}  \label{BC convex}
The Bost--Charles characteristic~$\widehat{T}$ is a convex increasing function of~$\log{r}$. 
\end{lemma}

\begin{proof}
The Poisson--Jensen formula allows the following rewriting of the double continuous integral~$\hT(r,\varphi)$ as a single continuous integral of a discrete sum: 
\begin{equation} \label{Tdef}
\hT(r,\varphi) = \int_{\T} \left\{  \log{|\varphi(z)|} + \sum_{ \substack{ u \in \D \\ \varphi(u) = \varphi(z)  } } \log^+{ \frac{r}{u} } \right\} \, \mv(z).
\end{equation}
We have substituted here $u := rw$. The 
lemma now follows upon remarking that the finite sum in the curly brackets is itself a convex increasing function of~$\log{r}$
 for each given $z \in \T$. 
 \end{proof}

We give two essentially equivalent formulations for the main theorem of this section. 

\begin{thm}\label{main:BC conv discrete}
Assume the same conditions and notation as in Theorem~\ref{main:BC form}.
Let
\[1 = r_l > r_{l-1} > \cdots > r_0>0\]
 be a sequence of subradii,  and consider the \emph{slopes}
\begin{equation} \label{slopes alpha}
\alpha_k := \frac{\hT(r_k, \varphi)- \hT(r_{k-1},\varphi)}{\log{r_k} - \log{r_{k-1}}}.
\end{equation} 
Assume that $\alpha_l \in [0, m]$.
 Then we have the following improvement  over the bound~\eqref{BCbound}: 
\begin{equation} \label{disc bounds}
m  \leq  \frac{ \iint_{\T^2} \log|\varphi(z)-\varphi(w)| \, \mv(z) \mv(w)  - \frac{1}{m} \sum_{k=1}^l \frac{ \left( \hT(r_k, \varphi)- \hT(r_{k-1},\varphi)\right)^2 }{\log{r_k} - \log{r_{k-1}}}}{  \log{|\varphi'(0)|} - \tau(\mathbf{b;e}) }.
\end{equation}
\end{thm}

By Lemma~\ref{BC convex}, the bound of Theorem~\ref{main:BC conv discrete} is only
improved if one refines the sequence of subradii, subject to the inequality~$\alpha_l \le m$ on the slopes.
Thus, in the limit, we obtain a continuous version of this theorem (see  Theorem~\ref{main:BC convexity} below).
In our experience, 
the extra numerical saving obtained in the limit is negligible  once one chooses just a few
 division points. Moreover,  it seems more practical from a computational standpoint to compute bounds on specific
values of~$\hT(r,\varphi)$ rather than integrals in terms of this function.

In order to formulate the continuous version of Theorem~\ref{main:BC conv discrete},
we firstly introduce the following
 positive increasing 
 function\footnote{This is, in fact, a continuous function. We will not use this, and we do not give a proof of the~$C^1$ property of~$\hT(r,\varphi)$. For the abstract purpose (not used elsewhere in our paper, neither) of making an {\it almost everywhere} sense of the~$d/dr$ derivative in~\eqref{A function} and the~$r \in [0,1]$ Riemann integral in~\eqref{BCbound-cov}, it suffices to appeal to Lebesgue's theorem that a monotone function~$[0,1] \to \R$ is differentiable almost everywhere.} of $r \in (0,1]$: 
\begin{equation}  \label{A function}
\hA(r,\varphi) := r  \frac{d}{dr} \hT(r,\varphi),
\end{equation}
whose notation mirrors the  traditional 
covering spherical area function 
$$
\oA(r,\varphi) :=  \frac{1}{\pi} \iint_{D(0,r)} \frac{|\varphi'|^2}{(1+|\varphi|)^2} \, dx dy = \iint_{D(0,r)} \varphi^* \omega_{\mathrm{FS}}
=: r \frac{d}{dr} \oT(r,\varphi)
$$
of the Ahlfors--Shimizu theory.

Now by interpreting the numerator of equation~(\ref{disc bounds}) as a Riemann
sum, in the limit (under the assumption~$\hA(r,\varphi) \le \hA(1,\varphi) \leq m$ for~$r \le 1$)
we obtain the following:

\begin{thm}\label{main:BC convexity}
Assume the same conditions and notation as in Theorem~\ref{main:BC form}.
Assume that $\hA(1,\varphi) \leq m$. 
Then
\begin{equation}\label{BCbound-cov} 
\begin{aligned}
m  &  \leq  \frac{ \iint_{\T^2} \log|\varphi(z)-\varphi(w)| \, \mv(z) \mv(w) - \frac{1}{m}  \int_0^1  \hA(r,\varphi)^2 \, \frac{dr}{r} }{  \log{|\varphi'(0)|} - \tau(\mathbf{b;e}) } \\
& =  \frac{ \int_{\T} \log{|\varphi|} \, \mv + \int_0^1 \hA(r,\varphi) \, \frac{dr}{r} - \frac{1}{m}  \int_0^1  \hA(r,\varphi)^2 \, \frac{dr}{r} }{  \log{|\varphi'(0)|} - \tau(\mathbf{b;e}) }. 
\end{aligned}
\end{equation}
\end{thm}

A further improvement is given in Theorem~\ref{main:BC fullconv} by using a (heuristically speaking) optimal choice of the Euclidean metric in the evaluation module $E_D$.

\subsubsection{Binomial metrics}  \label{sec: easy convexity}
As in Theorem~\ref{main:elementary form}, the set of analytic maps $\varphi_n$ used to estimate the $n^{\mathrm{th}}$ archimedean evaluation height
does not need to be of the particular form $\varphi_n(z) = \varphi(r(n) z)$. We include here one more elementary and fully explicit bound using~$l+1=2$ knots with $\varphi_0(z) = R z$ and $\varphi_1(z) = \varphi(z)$, but  ---  unlike with~\S~\ref{sec:BC convexity}  ---   taking a family of metrics independent of the map~$\varphi$. On the space of real auxiliary polynomials~$E_D \otimes_{\Z} \R = \R[x]_{< D}$, the metric can be described as diagonalizing the monomials basis $\left\{ x^k \right\}_{k=0}^{D-1}$ and giving the weights $\| x^k \| := \exp \left(  \lambda  D (k/D)^r + \mu k \right)$. The bound works out to the following explicit form, in which the triple of \emph{binomial metric weights} $\{r, \lambda, \mu\}$ is to be optimized. Unlike for all our other holonomy bounds in this section~\S~\ref{new slopes} and the next~\S~\ref{slopes}, the proof of this theorem does not require any of the results from~\cite{BostCharles}. 

\begin{thm}\label{main: easy convexity}
Assumptions and notation as in Theorem~\ref{main:BC form}. We further assume that $f_1,\ldots, f_m\in \Q\llbracket x \rrbracket$ are all convergent on the complex disc $|x|<R$. 
For $r\in \R_{>1}, \lambda\in \R_{>0}, \mu \in \R$, set   
\begin{equation*}
\begin{aligned}
\Gamma(x; r, \lambda, \mu) & :=\min\left\{(r-1)\frac{\big(\max\{0,x-\mu\}\big)^{r/(r-1)}}{(r^r\lambda)^{1/(r-1)}}, \max\{(r-1)\lambda, x-\lambda-\mu\}\right\}, \\
  T(\varphi;r,\lambda, \mu)  & := \int_{\T} \Gamma(\log|\varphi(z)|; r, \lambda, \mu)\, \mv(z), \\
  T_{r,\lambda,\mu}(\varphi) & : = \frac{\lambda}{r+1} + \frac{\mu}{2} + T(\varphi;r,\lambda,\mu),  \\
 \chi_0 & :=\min\left\{1, \left(\frac{\max\{0, \log{R} - \mu\} }{\lambda r}\right)^{1/(r-1)}\right\}. 
\end{aligned}
\end{equation*}
Suppose that $\mu \leq \log{R} < \log{|\varphi'(0)|}$ and 
$$
\chi_0 \leq \chi_1 := \frac{ T(\varphi;r,\lambda,\mu) - \Gamma(\log{R}; r, \lambda, \mu) }{\log |\varphi'(0)|- \log{R}} \leq m, \qquad \chi_0 < 1.
$$

 Then
\begin{equation}\label{EasyConvexityBound}
m \leq  \frac{2T_{r,\lambda,\mu}(\varphi) - \frac{2}{m} \left(\frac{1}{2}\chi_1^2 \log{\frac{|\varphi'(0)|}{R}}+\chi_0\Gamma(\log R; r, \lambda,\mu) - \chi_0^2 (\log R - \mu)\left(\frac{1}{2} - \frac{1}{r(r+1)}\right)\right)}{\log |\varphi'(0)| - \tau(\bb;\be)}.
\end{equation}
\end{thm}

 The extra assumption in this theorem concerning the inequalities among $\log R$, $\mu$, 
 $\log |\varphi'(0)|$, $m$, $\chi_1$, $\chi_0$ can be bypassed in the proof that follows to get, on a
 case-by-case basis, \emph{some} bound on~$m$ in every case. We pick these particular conditions as they are satisfied in our applications. See Example~\ref{Ex-easyconv}.

The function $\Gamma(x;r,\lambda,\mu)$ emerges as the \emph{Legendre transform}~\cite[\S~VI]{Ellis} of the binomial metric weight function~$\lambda t^r + \mu t$. This basic explicit computation is the content of our next lemma.  

\begin{lemma}\label{T-formula}
For arbitrary $r\in \R_{>1}, \lambda\in \R_{>0}, \mu \in \R$, and $x \in \R$, we have 
\[\max_{0\leq t \leq 1} \left\{ tx - \lambda t^r - \mu t \right\} = \Gamma(x; r, \lambda, \mu). \]
Therefore,
\[T(\varphi;r,\lambda, \mu)=\int_{\T}\max_{0\leq t \leq 1} \left\{ t \log |\varphi(z)| - \lambda t^r - \mu t \right\}  \, \mv(z).\]
\end{lemma}

\begin{proof}
The proof is a direct computation. Set $F(t) := xt - \lambda t^r - \mu t$, regarding $x, \lambda$, and $\mu$ as fixed. Then $F'(t)= x - r\lambda t^{r-1} - \mu$, and so $F$ admits a critical point in $t\in \R_{\geq 0}$ if and only if $x -\mu \geq 0$, in which case the unique such critical point is 
\[t_0 := \left(\frac{x - \mu}{r \lambda}\right)^{1/(r-1)}.\]
Therefore 
\[\max_{0\leq t \leq 1} F(t) =
	\begin{cases}
               F(0)= 0&\text{ if } x -\mu \leq 0\\
               F(t_0)= (r-1)\frac{(x-\mu)^{r/(r-1)}}{(r^r\lambda)^{1/(r-1)}}&\text{ if }0\leq x -\mu \leq \lambda r\\
               F(1)= x - \lambda  - \mu &\text{ if } x -\mu \geq \lambda r
            \end{cases}           .\]
This is why we defined $\Gamma(x; r, \lambda, \mu)$ the way we did in the statement of Theorem~\ref{main: easy convexity}. 
\end{proof}

\subsection{A brief review of Bost's slopes method}  \label{slopes review}

We review Bost's slopes method and related background material from Arakelov theory. The main references are 
\cite[\S\S 4.1, 4.2]{BostFoliations} and \cite[Chapter 1]{BostBook}. For simplicity, we only recall the theory over $\Q$ as that is sufficient for our applications. Everything in this section holds verbatim for any number field, see \cite{zeta5} and Remark~\ref{BCboundK}.

\subsubsection{Hermitian vector bundles on $\spec{\Z}$}  \label{sec:Hermitian vector bundles}

\begin{df}
A \emph{Euclidean lattice} is a pair $\overline{E} = (E, \| \cdot \|)$ made of a finite rank free $\Z$-module $E$ and 
a Euclidean norm $\| \cdot \|$ on the vector space $E_\R$. In other words: $\| \cdot \|^2$ is a positive definite quadratic form on $E_{\R} := E \otimes_{\Z} \R$. 
\end{df}

In Arakelov geometry, this coincides with the notion of a Hermitian vector bundle on $\spec{\Z}$. We thus use the notion of the \emph{arithmetic degree}
defined as the negative of the logarithm of the covolume of the Euclidean lattice:
\begin{equation} \label{arithmetic deg} 
\ardeg{\overline{E}} := - \log{ \mathrm{covol} (E, \| \cdot \|)} = - \frac{1}{2}  \log{ \big|  \det \big( \langle e_i, e_j  \rangle \big)_{i,j = 1}^{r} \big| }.
\end{equation}
Here, $r := \rank{E}$, 
$$
\langle e, f \rangle := \frac{ 1}{2}  \| e+f \|^2 - \|e\|^2 - \|f\|^2  
$$
is the associated inner product giving the quadratic form $\| e \| = \sqrt{ \langle e, e \rangle }$, and $e_1, \ldots, e_r$ is any $\Z$-module basis of 
$E$.\\

Let $M_{\Q}$ denote the equivalence classes (\emph{places}) of absolute values $\Q \to [0,\infty)$.
At the place $v \in M_{\Q}$, we select the representative absolute value $|\cdot|_v$ with the usual 
normalizations: $|\cdot |_{\infty}$ is the usual absolute value for the archimedean place $\infty \in M_{\Q}$, and $|p|_p = 1/p$ for the $p$-adic place $p \in M_{\Q}$. Thus $\prod_{v \in M_{\Q}} |x|_v = 1$ for all $x \in \Q^{\times}$. Let $M_{\Q}^{\mathrm{fin}} := M_{\Q} \setminus \{\infty\}$ denote the set of all finite places of $\Q$, which is identified with the set of rational primes.

Along with the quadratic form $\| \cdot \|$ on $E_\R$, it is convenient 
to consider the $p$-adic norms $\| \cdot \|_p$ defined on $E_{\Q_p}$ by
\[
 \Big\|  \sum_{i=1}^{r} c_i e_i \Big\|_p := \max_{1 \leq i \leq r}  |c_i|_p. 
\]
 Note that $\| \cdot \|_p$ is independent of the choice of the basis $e_1,\ldots,e_r$ of~$E$: more intrinsically, we have $p^{\Z}$ for 
 the value group of $\| \cdot \|_p$, with $\| w \|_p=p^{-n}$ if and only if $w\in E\otimes_\Z p^{n}\Z_p$ and $w\notin E\otimes_\Z p^{n+1}\Z_p$.
Thus the Euclidean structure combines with the integral lattice structure $E$ of the $\Q$-vector space $E_\Q := E \otimes_{\Z} \Q$ to define an adelic metric $(E_{\Q}, (\| \cdot \|_v)_{v \in M_{\Q}})$. Conversely, we can recover the lattice $E \subset E_{\Q} \subset E_{\R}$ as the $w \in E_{\Q}$ defined by
the simultaneous conditions $\|w\|_p \leq 1$ for all primes~$p$.
 
In these notations, the arithmetic degree formula~\eqref{arithmetic deg} takes the following adelic form: 
 \begin{equation*}
 \ardeg{\overline{E}}  =    - \frac{1}{2}  \log{ \left|  \det \big( \langle v_i, v_j  \rangle \big)_{i,j = 1}^{r} \right| }
  - \sum_{p \in M_{\Q}^{\mathrm{fin}}} \sum_{i=1}^{r} \log{ \|  v_i \|_p}, 
 \end{equation*}
 where $\{v_1, \ldots, v_r\}$ is any $\Q$-basis of $E_\Q$.

Given two Euclidean lattices $\overline{E}, \overline{F}$, let $\overline{E}\oplus \overline{F}$ denote $E\oplus F$ equipped with the norm given by the orthogonal direct sum of the norms on the subspaces $E_\R$ and $F_\R$. By definition, we have
\begin{equation}\label{degsum}
\ardeg(\overline{E} \oplus \overline{F})=\ardeg(\overline{E}) + \ardeg(\overline{F}).
\end{equation}

\subsubsection{Slopes of Euclidean lattices and heights of morphisms}\label{sec_def_slopes}

\begin{df}
The \emph{slope} $\widehat{\mu}(\overline{E})$ of a Euclidean lattice $\overline{E} = (E, \| \cdot \| )$ is defined as 
$$
\widehat{\mu}(\overline{E}) := \frac{\ardeg{\overline{E}}}{\rank{E}} \in \R.
$$
The \emph{maximal slope} of $\overline{E}$ is defined as
$$
\widehat{\mu}_{\max}(\overline{E}) := \sup_{0 \subsetneq F \subseteq E} \widehat{\mu}(\overline{F}),
$$ 
where $F$ runs through all nonzero $\Z$-submodules of $E$ and $\overline{F}$ denotes the induced Euclidean lattice of $F$ equipped with the quadratic form obtained from restricting $\| \cdot \|$ to $F_\R$. 
\end{df}

The following lemma follows from the definition.
\begin{lemma}\label{slope-tensor}
Let $\overline{E}, \overline{F}$ be two Euclidean lattices. Let $\overline{E}\otimes \overline{F}$ denote $E\otimes_{\Z} F$ equipped with the tensor norm. Then we have
\[\slope(\overline{E} \otimes \overline{F})=\slope(\overline{E}) + \slope(\overline{F}).\]
\end{lemma}

\begin{proof}
See, for instance, \cite[Lemma~2.3]{Chen}, which we briefly summarize here. By definition of the arithmetic degree, we have $\ardeg (\overline{E})=\ardeg (\wedge ^{\rank E} \overline{E})$, and for any two rank $1$ Euclidean lattices $\overline{L}_1, \overline{L}_2$, we have $\ardeg (\overline{L}_1 \otimes \overline{L}_2) = \ardeg (\overline{L}_1) + \ardeg( \overline{L}_2) $. Note that \[\wedge ^{\rank E\otimes F} (\overline{E} \otimes \overline{F})\cong (\wedge ^{\rank E} \overline{E})^{\otimes \rank F} \otimes (\wedge ^{\rank F} \overline{F})^{\otimes \rank E}.\]
Thus 
\[\ardeg(\overline{E}\otimes \overline{F})= (\rank F) \, \ardeg(\overline{E}) + ( \rank E) \, \ardeg(\overline{F}),\]
completing the proof of the lemma.
\end{proof}

Consider two Euclidean lattices $\overline{E} = (E, \| \cdot\|_E)$ and $\overline{F}  = (F, \|\cdot \|_F)$ and an \emph{injective} homomorphism
$\psi : E_\Q \hookrightarrow F_\Q$. If $\psi$ sends the Euclidean lattice $E$  isometrically into a sublattice of $F$, then $\widehat{\mu}(\overline{E}) \leq \widehat{\mu}_{\max}(\overline{F})$ by the definition of the maximal slope.  In general,  the slope of the source lattice $\overline{E}$ can be upper-estimated in terms of the maximal slope of 
the range lattice $\overline{F}$ and the \emph{height} of the homomorphism $\psi$. 

\begin{df} 
The \emph{local $v$-adic height} (at a place $v \in M_{\Q}$) of the monomorphism $\psi$
is defined as the logarithm of the norm of the induced monomorphism 
$$
\big( E_{\Q_v}, \| \cdot \|_{E, v} \big) \hookrightarrow \big( F_{\Q_v}, \| \cdot \|_{F,v} \big)
$$
of normed $\Q_v$-vector spaces: 
\[
h_v(\psi)  :=\sup_{e \in E_{\Q_v} \setminus \{0\}}  \log  \,  \frac{\| \psi(e) \|_{F, v} }{ \| e\|_{E,v}  }=  \sup_{e \in E \setminus \{0\}}  \log  \,  \frac{\| \psi(e) \|_{F, v} }{ \| e\|_{E,v}  }.
\]
The \emph{global height} of  $\psi$ is the sum of the local $v$-adic heights over all places $v \in M_{\Q}$:
$$
h(\psi) := \sum_{v \in M_{\Q}} h_v(\psi).
$$
\end{df}

The tautological inequality $\widehat{\mu}(\overline{E}) \leq \widehat{\mu}_{\max}(\overline{F})$ for the isometric injections $E \hookrightarrow F$ then generalizes to 
arbitrary monomorphisms $\psi : E_\Q \hookrightarrow F_\Q$, in the  following way.

\begin{lemma}[\cite{BostFoliations}, Prop.~4.5]
For every monomorphism 
$$
\psi : E_\Q \hookrightarrow F_\Q
$$
of the induced $\Q$-vector spaces of the Euclidean lattices $\overline{E}$ and $\overline{F}$, we have
\begin{equation}  \label{basic slopes}
\widehat{\mu}(  \overline{E} )    \leq  \widehat{\mu}_{\max}(\overline{F}) + h(\psi).
\end{equation}
\end{lemma}

\subsubsection{Bost's slopes inequality}\label{sec_slope_ineq}   
For filtered Euclidean lattices, the slopes inequality~\eqref{basic slopes} generalizes as follows. Let $F$ be a free $\Z$-module, which we no longer require  to be of finite rank.
We suppose that there is a filtration on $F_{\Q}$ 
$$
 F^{\bullet}_{\Q}  \, : \, F_{\Q} = F^{(0)}_{\Q} \supseteq F^{(1)}_{\Q} \supseteq F^{(2)}_{\Q} \supseteq \cdots
$$
with finite-dimensional graded quotients $\Gr_n(F^{\bullet}_{\Q}) := F_{\Q}^{(n)} / F_{\Q}^{(n+1)}$ and such that $\cap_{n=0}^\infty F^{(i)}_{\Q}=\{0\}$.

We consider $\overline{E}, \psi$ as in \S~\ref{sec_def_slopes}. In particular,
the linear monomorphism $\psi : E_\Q \hookrightarrow F_{\Q}$ induces a filtration on $E$:
$$
E^\bullet\, : \,E = E^{(0)} \supseteq E^{(1)} \supseteq \cdots, \quad \text{where } E^{(n)}   =  E \cap \psi^{-1}(F_{\Q}^{(n)}).
$$
Note that since $E$ is a finite rank free $\Z$-module and $\psi$ is injective, we have that $\Gr_n(E^\bullet)$ are finite rank free $\Z$-modules and the above filtration stabilizes to $\{0\}$ after finitely many steps. Moreover, the restriction of $\|\cdot \|_E$ to $E^{(n)}$ gives $E^{(n)}$ a Euclidean lattice structure and the corresponding quotient metric on $E^{(n)} / E^{(n+1)}$ equipped it with a Euclidean lattice structure $\overline{E^{(n)} / E^{(n+1)}}$.

We also assume that each graded quotient piece $\Gr_n(F_{\Q}^{\bullet})$ is endowed with a Euclidean
lattice structure. More precisely, for each $n$, we have a Euclidean lattice 
$$
\overline{G^{(n)}} =  \big(  G^{(n)},  \| \cdot \|_{G^{(n)}} \big),
$$
where $G^{(n)} \subset \Gr_n(F_{\Q}^{(n)}) $  a $\Z$-submodule such that
$G^{(n)}_\Q = \Gr_n(F_{\Q}^{(n)})$.

The map $\psi$ induces a linear monomorphism between the graded quotients: 
\begin{equation} \label{gr map}
\psi_D^{(n)}  : \, \Gr_n(E^{\bullet})_\Q \hookrightarrow \Gr_n(F^{\bullet}_\Q) 
\end{equation}
and its height $h( \psi_D^{(n)})$ is defined using the above-mentioned Euclidean lattices structures $\overline{  E^{(n)} / E^{(n+1)} }$ and $\overline{ G^{(n)} }$.
\begin{lemma}[\cite{BostFoliations}, Prop.~4.6] 
In this situation, 
\begin{equation}\label{slope-inequality} 
\ardeg{(\overline{E})}   \leq \sum_{n=0}^{\infty} \rank(E^{(n)} / E^{(n+1)}) \Big[ \widehat{\mu}_{\max}(\overline{ G^{(n)} } ) + 
h( \psi_D^{(n)} )  \Big]. 
\end{equation}
Note that the above sum is a finite sum since $E^{(N)}=0$ for $N\gg 1$.
\end{lemma}

\subsection{The Bost--Charles bound}  \label{one variable slopes}  
We follow \cite{BostBook, BostCharles} with a slight modification to take denominators into account.

We recall the setting of their work for our application. Consider $\cX=\P^1_\Z$ and the line bundle $\cL:=\OL(1)$ on $\cX$. We denote by $x := X_1/X_0$  the coordinate of (an affine line in) $\P^1_\Z = \mathrm{Proj} \,\Z[X_0,X_1]$, and then for $D\in \Z_{>0}$ we follow the usual identifications $\cL=\OL([0])$ (here $[0]$ denotes the divisor of the point $x=0$) and $ \Z[1/x]_{\leq D}= \Gamma(\cX, \cL^{\otimes D})$.

\subsubsection{The Bost--Charles metric} \label{BC metric}
Following the ideas in \cite[\S\S 8.2, 8.3]{BostCharles}, using $\varphi:(\Db, 0) \rightarrow (\P^1(\C), 0)$, we endow $\cL = \OL(1)$ with the Hermitian metric $\|\cdot\|_{\overline{\cL}}$ 
 defined by \[\|\mathbf{1}(y)\|_{\overline{\cL}} := \exp\left(-\sum_{z\in \varphi^{-1}(y)}\log^+{\frac{1}{|z|}} \right)=\prod_{z\in \Db, \, \varphi(z)=y} |z|,\]
where $\mathbf{1} = \mathbf{1}_{[0]}$ is the canonical section (``constant function'') of $\cL = \OL([0])$ corresponding to the divisor~$[0]$, and $y\in \P^1(\C)$. This Hermitian metric has $\cC^{\mathrm{b}\Delta}$ regularity in the sense of Bost--Charles (see \cite[\S\S 4.1.1, 4.2.1.2]{BostCharles}). We shall denote this Hermitian line bundle by $\overline{\OL(1)}$, or by~$\overline{\OL(1)}_{\varphi}$ if we wish to indicate the dependence on~$\varphi$. We work in the framework of arithmetic intersection theory using such Hermitian line bundles and Arakelov divisors with $\cC^{\mathrm{b}\Delta}$ regularity in the sense of Bost--Charles~\cite[\S 4.5]{BostCharles}.

\begin{remark}\label{rmk_BCmetric}
In~\cite{BostCharles}, the Hermitian metric of~\S~\ref{BC metric} is given by the Arakelov divisor $([0], \varphi_* \left(\log^+|z|^{-1})\right)$. More precisely, in the sense of \cite[\S 6.2.1 using Example~4.3.1]{BostCharles}, 
we have the compactly supported Arakelov divisor $([0], \log^+|z|^{-1})$ on the smooth
formal-analytic (hereafter, f.-a.) ~arithmetic surface 
\[\wV(\varphi) := (\spf \Z \llbracket x\rrbracket, (\Db, 0), i_{\varphi})\]
 over $\Z$, where $x\mapsto \varphi$ defines an isomorphism $\C \llbracket x \rrbracket \iso \C \llbracket z \rrbracket$ (here $z$ denotes the coordinate on $\Db$, and we use the assumption that $\varphi'(0)\neq 0$), and thus its compositional inverse induces an isomorphism $i_\varphi: \spf \C \llbracket x \rrbracket \cong \widehat{\D}_0$. See~\cite[\S 10.6.1]{BostBook} or~\cite[\S 6.1.1]{BostCharles} for the general definition of smooth f.-a.~arithmetic surface over a number field, and \S~6.4.1.1 in \emph{loc.~cit.} for this construction~$\wV(\varphi)$, which also comes with a distinguished nonconstant regular function $(\iota,\varphi) : \wV(\varphi) \to \A_{\Z}^1$ on the f.-a.~arithmetic surface; cf~\cite[\S 7.1.1.1]{BostCharles}. 
 Here, $\iota: \spf \Z \llbracket x \rrbracket \hookrightarrow \spec \Z[x] = \A_{\Z}^1$ is the natural formal immersion. 
 In the setting of \cite[\S 7.2.1]{BostCharles}, the Arakelov divisor  $([0], \varphi_* (\log^+|z|^{-1}))$ is the direct image of $([0], \log^+|z|^{-1})$ by the morphism $(\iota, \varphi) : \wV(\varphi) \to \A_{\Z}^1$, where the pushforward map $\varphi_*$ on Green functions is defined in \cite[\S 3.4.2.1]{BostCharles}. By \cite[Corollary 4.4.2(ii)]{BostCharles}, this pushforward map preserves $\cC^{\mathrm{b}\Delta}$ regularity of Green functions, which is essential for having a well-behaved arithmetic intersection theory. The explicit formula of the Hermitian metric associated to $\varphi_* (\log^+|z|^{-1})$ is given in \cite[\S 5.1.2]{BostCharles}.~\endofremark
\end{remark}

\subsubsection{Direct images and arithmetic Hilbert--Samuel} \label{direct images}
Fix a smooth probability measure $\nu$ on $\P^1(\C)$, for instance the Fubini--Study form $\omega_{\mathrm{FS}} = \frac{\sqrt{-1}}{2\pi} \frac{dz \wedge d\bar{z}}{(1+|z|^2)^2}$; the choice of $\nu$ is immaterial to the proof.  As in \cite[\S~6]{BostCharles}, the Hermitian metric on $\overline{\cL}$ combines with fiberwise integration over the manifold $\cX(\C)$ to define a Euclidean lattice structure on the $\Z$-module $\Gamma(\cX, \cL^{\otimes D})$. 
Explicitly, we norm $s\in \Gamma(\cX, \cL^{\otimes D})$ by
\[\|s\|:= \sqrt{\int_{\cX(\C)} \|s\|_{\overline{\cL}}^2 \, \nu}.\]
Following~\cite[\S~6.1.2.2]{BostCharles}, we denote by $\Gamma_{L^2}\left( \cX, \nu;  \overline{\cL}^{\otimes D} \right)$ this Euclidean lattice. 
Up to the integration metric weight~$\nu$, in a $D \to \infty$ asymptotic sense, this is essentially the zeroth direct image of~$\overline{\cL} = \overline{\OL}(1)_{\varphi}$ under the structure morphism $\wV(\varphi) \to \spec{\Z}$. 
As in \cite[\S 8, Theorem 8.2.5]{BostCharles}, we can express the arithmetic Hilbert--Samuel formula on the arithmetic surface
$\cX = \P_{\Z}^1$ into the form 

\begin{equation}\label{P1aHS}
 \ardeg \Gamma_{L^2}\left( \cX, \nu;  \overline{\cL}^{\otimes D} \right) = \frac{1}{2}(\overline{\cL} \cdot \overline{\cL}) D^2+ o(D^2).
  \end{equation}
When the Hermitian metric in $\cL = \OL(1)$ is smooth, this formula is due to Zhang \cite[Theorem~1.4]{Zhang} with an additional input by Bost in comparing two Hermitian metrics in the proof of~\cite[Theorem~10.3.2]{BostBook}. 
Zhang's theorem is a refinement to (non-pointwise-strict) semipositive curvature (Chern form) $c_1(\cL, \| \cdot \|) \geq 0$ of the work of Gillet--Soul\'e \cite{GilletSoule} and Bismut--Vasserot \cite{BismutVasserot}; see also Abbes--Bouche~\cite{AbbesBouche} for an outline of a more direct approach. 
Following the idea in \cite[\S 5]{BostL21} and \cite[\S\S 3--4]{BostCharles} to separate the Green function into a smooth Green function and a $\cC^{\mathrm{b}\Delta}$ function, the same arithmetic Hilbert--Samuel formula holds for ample line bundles with $\cC^{\mathrm{b}\Delta}$ Hermitian metrics of pointwise non-negative Chern form (as defined in~\cite[\S~4.2.1.2]{BostCharles}), and so the formula is also valid in our setting. 

In terms of $\varphi$, Bost and Charles~\cite[Theorem 5.4.1 and Proposition 5.4.2]{BostCharles} provide the following formula for the self-intersection number: 
\begin{equation}\label{BCintersection}
(\overline{\cL}\cdot \overline{\cL}) = \left(\overline{\OL(1)}_{\varphi} \cdot \overline{\OL(1)}_{\varphi} \right) = \iint_{\T^2} \log{|\varphi(z) - \varphi(w)|} \, \mv(z) \, \mv(w)
= \hT(1,\varphi). 
\end{equation}

We will review their computation in our mild generalization in Lemma~\ref{BCconv-intersection} further down. 

\begin{proof}[Proof of Theorem~\ref{main:BC form}]
Note that the choice of $\bb$ is not unique; we may permute the columns of $\bb$ without changing the form of the~$f_i$. Therefore, after a suitable permutation of the columns, we may assume $u_1 \leq u_2 \leq \ldots \leq u_r$. {\it We keep this convention for all the proofs in \S~\ref{new slopes}.}

For $D\in \NwithzeroA$, we take for our evaluation module the following free $\Z$-module of rank $m(D+1)$: 
\begin{equation} \label{eval factors}
E_D:= \bigoplus_{h=0}^r \bigoplus_{i=u_h+1}^{u_{h+1}} \frac{[1,\ldots, \xi D]^{e_i}}{[ 1, \ldots, y_{h+1} D] \cdots [1,\ldots, y_r D]} \, f_i \cdot \Z[1/x]_{\leq D},
\end{equation} 
where $u_0:=0, u_{r+1}:=m$, $\xi \in [0,m]$ and  $y_h\in [0, b_h m] \subset \R$ are auxiliary parameters to be optimized in the proof.
Note that the indexing  of~$E_D$ in equation~(\ref{eval factors})
differs slightly from the notation of~\S~\ref{filtrations}; the difference amounts to considering polynomials
of degree~$<D$ rather than~$\le D$. We use this normalization  ---  which asymptotically makes no
difference and so is ultimately an aesthetic choice  ---  so that we can
talk below about sections of~$\cL^{\otimes D}$ rather than~$\cL^{\otimes (D-1)}$.

To endow $E_D$ with a Euclidean norm, we take the orthogonal direct sum~\eqref{eval factors} of the lattices
\[ \frac{[1,\ldots, \xi D]^{e_i}}{[ 1, \ldots, y_{h+1} D] \cdots [1,\ldots, y_r D]} \Z[1/x]_{\leq D} \subset \Gamma\left(\cX, \cL^{\otimes D}\right)_\R,\] with each of these summands inheriting the norm induced from $\overline{\OL(1)}_{\varphi} = \overline{\cL}$. 
We shall denote this Euclidean lattice by~$\overline{E}_D$. 

By \eqref{P1aHS} and~\eqref{BCintersection}, we have
\begin{equation}\label{BCardegED}
\begin{aligned}
\ardeg \ovE_D & = \left(\frac{m}{2}\left(\overline{\OL(1)}_{\varphi} . \overline{\OL(1)}_{\varphi} \right) +\sum_{h=1}^r u_h y_h - \xi \left( \sum_{i=1}^{m} e_i \right) \right)D^{2} +o(D^{2}) \\ 
& = \left(\frac{m}{2} \, \hT(1,\varphi) +\sum_{h=1}^r u_h y_h - \xi \Big( \sum_{i=1}^{m} e_i \Big) \right)D^{2} +o(D^{2}). 
\end{aligned}
\end{equation}

Let $X$ denote $\cX_{\Q} = \P_{\Q}^1$. We  identify $\Spf \Q \llbracket x \rrbracket =\widehat{X}_{0}$ as the formal completion of~$X$ at its closed subscheme~$0$. This designates $f_i (x) \in \Gamma(\widehat{X}_{0}, \OL_{\widehat{X}_{0}})$, for $i = 1,\ldots,m$. Let $\Gamma(\widehat{X}_{0}, \cL^{\otimes D})$ denote the global sections of $\cL^{\otimes D}|_{\widehat{X}_{0}}$; these get identified with 
\[\Gamma(\widehat{X}_{0}, \cL^{\otimes D})=x^{-D} \Q \llbracket \bx \rrbracket=: F_{\Q}.\] (Here, $x^{-D} \Q \llbracket \bx \rrbracket$ denotes the $\Q$-vector space generated by $x^k$, where $k\geq -D$.) Thus $f_i \Gamma(\cX, \cL^{\otimes D}) \subset  \Gamma(\widehat{X}_{0}, \cL^{\otimes D})$, and we have the evaluation map 
\begin{equation} \label{eval D}
\psi_D : E_{D} \otimes_{\Z} \Q \hookrightarrow F_{\Q}, \qquad (Q_i)_{1\leq i\leq m} \mapsto \sum_{i=1}^m f_iQ_i,
\end{equation}
where $Q_i \in \Gamma(\cX, \cL^{\otimes D})_\Q$. It is an injective homomorphism due to our assumed $\Q(x)$-linear independence of the formal
functions~$f_i(x)$. 

We filter $F$ by the $x=0$ vanishing order of the formal sections of $\cL^{\otimes D}|_{\widehat{X}_{0}}$:
\[F_{\Q}=F_{\Q}^{(0)}\supseteq F_{\Q}^{(1)} \supseteq \cdots \supseteq F_{\Q}^{(n)} \supseteq \cdots.\]
Concretely, $F_{\Q}^{(n)}= \Cspan_{\Q}\{x^{k-D} \mid k\geq n\}$. The graded piece $F_\Q^{(n)}/F_\Q^{(n+1)}$ is a one dimensional $\Q$-vector space generated by the image of $x^{n-D}$ under the quotient map. We take the Euclidean lattice structure on $F_\Q^{(n)}/F_\Q^{(n+1)}$ given by the free rank one $\Z$-module generated by the image of $x^{n-D}$ and the Euclidean norm with $\| x^{n-D} \|=1$. Note that these Euclidean lattice structures on graded piece are all induced from the free $\Z$-module $F=x^{-D} \Z\llbracket x \rrbracket$ and the Euclidean norm on $x^{-D}\R[x]$ that has $\{x^n\}_{ n\in \Z_{\geq -D}}$ for an orthonormal basis.

As in \S~\ref{sec_slope_ineq}, we use ${E}^{(n)}_{D} := \psi_D^{-1} \left( F^{(n)}_\Q \right) \cap E_D$ to denote the preimage of $F^{(n)}_\Q$ in $E_D$ under $\psi_D$. For each $n \in \NwithzeroB$, the evaluation map~\eqref{eval D} induces an injective homomorphism 
\[\psi_D^{(n)}: {E}^{(n)}_{D}/ {E}^{(n+1)}_{D}\hookrightarrow F_\Q^{(n)}/F_\Q^{(n+1)}.\] 
In particular, as in~\eqref{van fj}, we have $\rank \left( {E}^{(n)}_{D}/{E}^{(n+1)}_{D} \right) \in \{0,1\}$. Let 
\[\mathcal{V}_{D} := \left\{n\in \NwithzeroB \mid \rank \left( {E}^{(n)}_{D}/{E}^{(n+1)}_{D} \right)=1 \right\}.\] We have $\# \mathcal{V}_{D} = \rank E_D=m(D+1)$.

We now provide upper bounds on the evaluation heights $h_\infty(\psi_D^{(n)})$ and $h_{\textrm{fin}}(\psi_D^{(n)})$, where $h_{\textrm{fin}}(\psi_D^{(\bn)}):=\sum_{v\in M_\Q, v\nmid \infty} h_v(\psi_D^{(\bn)})$.

The archimedean evaluation height bound stems from the work of Bost and Bost--Charles: 
\begin{equation}\label{BC-arch-ht-bound}
h_{\infty}(\psi_D^{(n)}) \leq  - n \log | \varphi'(0)| +   \left(\overline{\OL(1)}_{\varphi} \cdot \overline{\OL(1)}_{\varphi}\right) D + o(D).
\end{equation}
The proof details for our specific setting are in either~ \S~\ref{sec_BCconvexity} (Lemma~\ref{arch-ht-BCconv}, specializing to $r=1$) or in~\S~\ref{sec:slopes archimedean}) (specializing to $d=1$), where respectively we will need a refinement of this estimate to incorporate convexity and handle the high dimensional setup.
For the original source we refer to the two paragraphs following Theorem 8.2.2 on page~127 in \cite{BostCharles}, which in turn summarize the relevant sections of~\cite{BostBook}. The specific bound is essentially \cite[\S 10.5.5, Theorem 10.5.3, Corollary 10.5.4]{BostBook}. 

Next we estimate $h_{\mathrm{fin}}(\psi_D^{(n)} )$. For each prime $p$, by the definition of $h_p$, our task is to consider an arbitrary element $(Q_i)_{1\leq i\leq m} \in {E}^{(n)}_{D} \setminus {E}^{(n+1)}_{D}$, and to provide an upper bound on $\log | c_{n} |_p$, where $c_{n}$ denotes the (leading order~$n$) coefficient of  $x^{n}$ in $\sum_{i=1}^m f_iQ_i = c_n x^n +\ldots$. 

Recall the notation of the indices cutoffs~$u_j \in \{0,1,\ldots,m\}$, for $1 \leq j \leq r$, from the statement of Theorem~\ref{main:elementary form}. 
Let $h_i$ be the index in~$\{0,1,\ldots,r\}$ defined by $u_{h_i} < i \leq u_{h_i +1}$. The ultrametric triangle inequality for $|\cdot |_p$ directly gives

\begin{equation}  \label{vals p}
\begin{aligned}
&  \frac{\log |c_n|_p}{\log p} \\
& \leq  \max_{ \substack{ 1\leq i \leq m \\ 0\leq k \leq \min\{n-1, D\} }} \left\{ \mathrm{val}_p \left(\prod_{j=1}^{h_i} [1,\ldots, b_j (n-k)] \cdot \prod_{j=h_i+1}^r [1,\ldots, y_j D] \right) + \mathrm{val}_p\left( \frac{(n-k)^{e_i}}{[1,\ldots, \xi D]^{e_i}} \right) \right\} \\
& \leq  \left(\sum_{h=1}^r \mathrm{val}_p\left([1,\ldots, b_h \max\{n, (y_h/b_h) D\}] \right) \right) + \left(\max_{1\leq i \leq m} e_i\right) \mathrm{val}_p\left( \frac{[\max\{n-D,1\}, \ldots, n]}{[1,\ldots, \xi D]} \right).
\end{aligned}
\end{equation}
Note that for $n^{1/2} < p\leq \xi D$ we have $\mathrm{val}_p\left([\max\{n-D,1\}, \ldots, n]/[1,\ldots, \xi D] \right) \leq 0$. Since all terms under the $p$-adic valuation in~\eqref{vals p}  
are independent of $p$, the prime number theorem gives
\begin{equation*}
\begin{aligned}
h_{\mathrm{fin}}(\psi_D^{(n)})   \leq & \sum_{h=1}^r \log ([1,\ldots, b_h \max\{n, (y_h/b_h) D\}]) \\  
&+ \left(\max_{1\leq i \leq m} e_i\right) \sum_{p > \max\{n^{1/2}, \xi D\} } \mathrm{val}_p([\max\{n-D,1\}, \ldots, n]) \log p + o(n)  \\
\leq & \sum_{h=1}^r  b_h \max\{n, (y_h/b_h) D\}  \\
& +  \left(\max_{1\leq i \leq m} e_i\right) \sum_{ \substack{  p > \max\{n^{1/2}, \xi D\}, \\ p \mid [\max\{n-D,1\}, \ldots, n] } } \log p  +o(n+D). 
\end{aligned}
\end{equation*}

By Lemma~\ref{lcm-2}, we have for $n\geq \max\{\xi,1\} D$ 
\begin{equation} \label{finite>xi}
\begin{aligned} 
h_{\mathrm{fin}}(\psi_D^{(n)})  & \leq  \left(\max_{1\leq i \leq m} e_i\right)\left(\left(
D \sum_{j=1}^{\lfloor (n/D - 1)/\max(1,\xi ) \rfloor} 1/j   \right)  \right. \\
 &  \left. +  \left(\frac{n}{\lfloor (n/D + (\xi -1)^+ )/\max(1,\xi ) \rfloor} - \xi D \right)^+\right) \\
& + \sum_{h=1}^r  b_h \max\{n, (y_h/b_h) D\}  +o(n+D);  
\end{aligned}
\end{equation}
Again by Lemma~\ref{lcm-2} (taking $k=n-1$), we have for $\min\{\xi,1\} D \leq n <D$
\begin{equation}\label{finite-mid}
\begin{aligned}
h_{\mathrm{fin}}(\psi_D^{(n)}) & \leq \left(\max_{1\leq i \leq m} e_i\right)(n-\xi D)^+\\
& + \sum_{h=1}^r  b_h \max\{n, (y_h/b_h) D\}  +o(n+D);  
\end{aligned}
\end{equation}
for $n<\xi D$ the estimate is just
\begin{equation} \label{finite<xi}
h_{\mathrm{fin}}(\psi_D^{(n)}) \leq \sum_{h=1}^r  b_h \max\{n, (y_h/b_h) D\}  +o(n)+o(D).
\end{equation}

Andr\'e's Corollary~\ref{holonomic criterion}, proved in Appendix~\ref{app:PerelliZannier} but also in the self-contained Lemma~\ref{f_holonomic} below, 
permits us to apply the Chudnovsky--Osgood Theorem~\ref{KolchinSolved} on functional bad approximability. In the $D \to \infty$ asymptotic, by the continuity of  $v \mapsto I_u^v(w)$, this entails the total evaluation height upper bounds
\begin{equation}\label{BC-sum-arch-ht}
\sum_{n \in \mathcal{V}_D} h_\infty(\psi_D^{(n)}) \leq \left( - \frac{m^2}{2} \log | \varphi'(0)| +  m (\overline{\OL(1)} \cdot \overline{\OL(1)}) \right) D^2 + o(D^2).
\end{equation}
and
\begin{equation}
\begin{aligned}
\sum_{n \in \mathcal{V}_D} h_{\mathrm{fin}}(\psi_D^{(n)})&  \leq  \left( \left(\max_{1\leq i \leq m} e_i\right) I_{\xi}^m(\xi) + \sum_{h=1}^r  b_h \int_0^m \max\{s, y_h/b_h\} \, ds \right)D^2 + o(D^2) 
 \\
&=   \left(\left(\max_{1\leq i \leq m} e_i\right) I_{\xi}^m(\xi) + \frac{1}{2}\left(\sigma_m m^2 + \sum_{h=1}^r   y_h^2/b_h \right) \right)D^2 + o(D^2). \label{BC-sum-finite-ht}
\end{aligned}
\end{equation}

Let us give more details on how to obtain \eqref{BC-sum-finite-ht} using Theorem~\ref{KolchinSolved}; the verbatim reasoning applies also to \eqref{BC-sum-arch-ht} and to other similar evaluation height estimates in the rest of the section. 
Firstly, Lemma~\ref{f_holonomic} and our standing assumptions in Theorem~\ref{main:BC form} imply that $f_1,\ldots, f_m$ are holonomic functions. Let~$\varepsilon$ and $C(\varepsilon)$ be as in the statement of Theorem~\ref{KolchinSolved}. 
Throughout this section, the evaluation heights $h_\infty(\psi_D^{(n)})$ and $h_{\mathrm{fin}}(\psi_D^{(n)})$ get
asymptotically upper-estimated with $o(D+n)$ implicit error terms, and with certain explicit nonnegative main terms that are, in all cases, certainly $\geq - C' D$,  uniformly in $D\gg 1$ and~$\varepsilon$; in these estimates, $C'$ as well as the decay rates in the $o(D+n)$ of the error terms only depend on $m, \{f_i\}$, and $\varphi$. 
(For the $h_{\mathrm{fin}}(\psi_D^{(n)})$ situation under current highlight, we may of course simply take $C'=0$;  we keep $C'$ to illustrate how the argument 
in the other situations.) For any $\varepsilon>0$, Theorem~\ref{KolchinSolved} gives $\mathcal{V}_D \subset [0, (m+\varepsilon)(D+1) + C(\varepsilon)]$ with
$\#\mathcal{V}_D = m(D+1)$. Thus the total evaluation height $\sum_{n\in \mathcal{V}_D} h_{\mathrm{fin}}(\psi_D^{(n)})$ is majorized by the $0 \leq n \leq (m+\varepsilon)(D+1) + C(\varepsilon)$ sum of \eqref{finite>xi}, resp.~\eqref{finite-mid},~\eqref{finite<xi}, minus the overcount of at most $ \varepsilon(D+1) + C(\varepsilon)$ terms, to all of which we apply the $\geq - C' D$ lower bound to compensate. In the situation at hand, we get asymptotically

\begin{equation*}
\begin{aligned}
& \sum_{n \in \mathcal{V}_D} h_{\mathrm{fin}}(\psi_D^{(n)})  \\
& \leq   \left( \left(\max_{1\leq i \leq m} e_i\right) I_{\xi}^{m+\varepsilon +C(\varepsilon)/D}(\xi) + \sum_{h=1}^r  b_h \int_0^{m+\varepsilon +C(\varepsilon)/D} \max\{s, y_h/b_h\} \, ds \right)D^2 \\
& +  C'D(\varepsilon(D+1) + C(\varepsilon)) + o(D^2).
\end{aligned}
\end{equation*}
By the continuity of $I_u^v(w)$ in $v$, we derive with an arbitrary $\varepsilon > 0$ the upper estimate
\[\lim_{D\rightarrow \infty} \frac{\sum_{n \in \mathcal{V}_D} h_{\mathrm{fin}}(\psi_D^{(n)})}{D^2} \leq \left( \left(\max_{1\leq i \leq m} e_i\right) I_{\xi}^{m+\varepsilon}(\xi) + \sum_{h=1}^r  b_h \int_0^{m+\varepsilon} \max\{s, y_h/b_h\} \, ds \right)+ C' \varepsilon.\]
As the left-hand side 
is independent of the choice of $\varepsilon$ in Theorem~\ref{KolchinSolved}, we can let $\varepsilon \rightarrow 0$ and  obtain
\[\lim_{D\rightarrow \infty} \frac{\sum_{n \in \mathcal{V}_D} h_{\mathrm{fin}}(\psi_D^{(n)})}{D^2} \leq  \left(\max_{1\leq i \leq m} e_i\right) I_{\xi}^{m}(\xi) + \sum_{h=1}^r  b_h \int_0^{m} \max\{s, y_h/b_h\} \, ds.\]
This is exactly \eqref{BC-sum-finite-ht}.

Now Bost's slopes inequality~\eqref{slope-inequality} reads, in our situation: 
\begin{equation}
\ardeg \ovE_D \leq  \sum_{n \in \mathcal{V}_D} h_{\infty}(\psi_D^{(n)})  + \sum_{n \in \mathcal{V}_D} h_{\mathrm{fin}}(\psi_D^{(n)}). 
\end{equation}

Combining the upper bounds \eqref{BCardegED}, \eqref{BC-sum-arch-ht}, and \eqref{BC-sum-finite-ht}, and picking out the coefficients of  the leading $D^2$ of the $D \to \infty$ asymptotic, we derive by Bost's inequality~\eqref{slope-inequality} the upper bound
\begin{equation*}
\begin{aligned}
(\log|\varphi'(0)| - \sigma_m)m^2 &  \leq m \left(\overline{\OL(1)}_{\varphi} \cdot \overline{\OL(1)}_{\varphi} \right)
 + 2\left(\xi \left(\sum_{i=1}^m e_i\right) +\left(\max_{1\leq i \leq m} e_i\right) I_{\xi}^m(\xi) \right)  \\ 
 & \qquad +\left(\sum_{h=1}^r   y_h^2/b_h - 2\sum_{h=1}^r u_h y_h\right). 
\end{aligned}
\end{equation*}
The quadratic form $\sum_{h=1}^r   y_h^2/b_h - 2\sum_{h=1}^r u_h y_h$ reaches its minimum when $y_h=u_hb_h$ for all $1\leq h \leq r$, and thus
with~\eqref{BCintersection} 
 we obtained the desired bound.
\end{proof}

\begin{lemma}\label{f_holonomic}
If $\log |\varphi'(0)| > \sigma_m$, then all the $f_i$ are holonomic.
\end{lemma}

(See also Corollary~\ref{holonomic criterion} and its proof in Appendix~\ref{app:PerelliZannier}.)

\begin{proof}
Applying the differential operator $\left(x\frac{d}{dx}\right)^{e_i}$ to remove the terms $n^{e_i}$ from the denominators of the coefficients of $f_i$, we may and do assume --- with no loss of generality for the goal of proving the present lemma --- that $\mathbf{e=0}$. 
The following then is the familiar calculation as in Appendix~\ref{app:PerelliZannier}, which in effect gives another proof of~\eqref{basic basic}, now in the framework of Bost's slopes inequality. 
For our concrete purposes here, we need for every~$i$ to construct a $\Q(x)$-linear dependency among the derivatives $f_i, f_i', f_i'', \ldots$. Suppose to the
contrary that all those derivatives are $\Q(x)$-linearly independent. 
With an arbitrary $m' \in \NwithzeroA$, we apply a stripped down form of the main argument of the present section, now to the rank-$(m'+1)D$ evaluation module
\[E_D= \bigoplus_{j=0}^{m'} f_i^{(j)}\Z[x]_{<D},\]
equipped with the Euclidean norm in which $\{ f_i^{(j)}x^k\}$ is an orthonormal basis of $E_{D} \otimes_{\Z} \R$.

Then $\ardeg(\overline{E_D})=0$, and for any $n\in \mathcal{V}_D$, we have by the Poisson--Jensen formula (see for instance\footnote{We take $d=1$ and $p(x)=x$ in \cite{UDC}. There is a minor difference with the assumptions in \cite{UDC}:  we supposed there  $f_i(\varphi(z))$ to be holomorphic on $|z|\leq 1$, whereas here we only assume that $f_i(\varphi(z))$ are meromorphic on $|z|< 1$. For all our estimates on archimedean heights in \S~\ref{new slopes}, this difference is insignificant up to the error term of~$o(n+D)$. The point is that for any $\varepsilon>0$, we have $f_i(\varphi((1-\varepsilon)z))$ meromorphic on $|z|\leq 1$. There are at most finitely many meromorphic poles of all $f_i(\varphi((1-\varepsilon)z))$ in some neighborhood of~$\Db$, and we may take a polynomial $h(z)\in \C[z]$ with $h(0) = 1$ and such that $h(z)f_i(\varphi((1-\varepsilon)z))$ are all holomorphic on a neighborhood of $|z|\leq 1$. We observe that replacing all $f_i(\varphi(z))$ by $h(z)f_i(\varphi((1-\varepsilon)z))$ yields the same archimedean estimate once we let $\varepsilon\rightarrow 0$ at the end. See the proof of Lemma~\ref{arch-ht-BCconv} for details.} \cite[\S~2.4]{UDC}), 
\[h_\infty (\psi^{(n)}_D)\leq  - n \log |\varphi'(0)| +D  \int_\T \log^+|\varphi(z)| \, \mv(z) +o(n+D)\]
and by the prime number theorem,
\[h_{\mathrm{fin}} (\psi_D^{(n)})\leq \sigma_m n +o(n).\]

Fix an $\epsilon > 0$ with $\log{|\varphi'(0)|} > \sigma_m + \epsilon$. Then there is an 
$N_0 \in \NwithzeroA$ such that,
for all $n \geq N_0$, we have
\[h_\infty (\psi^{(n)}_D)\leq  - n \log |\varphi'(0)| +D  \int_\T \log^+|\varphi(z)| \, \mv(z) + (\epsilon/2) n + o(D);\]
\[h_{\mathrm{fin}}(\psi_D^{(n)})\leq (\sigma_m + \epsilon/2) n.\]
Hence, for all $n \in \NwithzeroB$, 
\[h_\infty (\psi^{(n)}_D)\leq  - n \log |\varphi'(0)| +D  \int_\T \log^+|\varphi(z)| \, \mv(z) + (\epsilon/2) n + o(D);\]
\[h_{\mathrm{fin}}(\psi_D^{(n)})\leq (\sigma_m + \epsilon/2)  n +o(D).\]
By the slopes inequality~\eqref{slope-inequality} and $\log{|\varphi'(0)|} > \sigma_m + \epsilon$, it ensues that
\begin{equation}
\begin{aligned}
0&= \ardeg(\overline{E_D}) \leq  \sum_{n =0}^{\infty} \mathrm{rank}(E_D^{(n)}/E_D^{(n+1)}) \cdot h(\psi_D^{(n)})  = \sum_{n\in \mathcal{V}_D} h(\psi_D^{(n)})  \\ 
&\leq - \left(\sum_{n \in \mathcal{V}_D}  n\right) \left(\log |\varphi'(0)| - \sigma_m -\epsilon\right) + (m'+1)D^2 \int_\T \log^+|\varphi(z)| \mv(z) +o(D^2)  \\
&\leq - \left(\sum_{n =0}^{(m'+1)D-1}  n\right) \left(\log |\varphi'(0)| - \sigma_m -\epsilon\right) + (m'+1)D^2 \int_\T \log^+|\varphi(z)| \mv(z) +o(D^2) \\
&= - \binom{(m'+1)D}{2} \left(\log |\varphi'(0)| - \sigma_m -\epsilon\right) + (m'+1)D^2 \int_\T \log^+|\varphi(z)| \mv(z) +o(D^2)
\end{aligned}
\end{equation}
Comparing the leading asymptotic order~$D^2$ coefficients and then letting $\epsilon\rightarrow 0$, we have
\[m'+1\leq \frac{2\int_\T \log^+|\varphi(z)| \mv(z) }{\log |\varphi'(0)| - \sigma_m} < \infty,\]
contrary to our assumption that~$m'$ could be arbitrarily large. 
\end{proof}

\begin{example}
For the case (see~\S~\ref{sec:proofA}) 
$$
\mathbf{b} := \left(  \begin{array}{llllllllllllll}   0 & 2 & 2 & 2 & 2 & 2 & 2 & 2 & 2 & 2 & 2 & 2 & 2 & 2 \\
0 & 0 & 0 & 2 & 2 & 2 & 2 & 2 & 2 & 2 & 2 & 2 & 2 & 2  \end{array} \right)^{\mathrm{t}}
$$
of relevance to the proof of Theorem~\ref{mainA}, we compute
$$
\tau^{\flat}(\mathbf{b}) = (2+2) - \frac{2\cdot 1^2 + 2 \cdot 3^2}{14^2}=\frac{191}{49} . 
$$
\end{example}

\begin{example}
For the case (see~\S~\ref{sec:proofC}) 
$$
\mathbf{b} := \left(  \begin{array}{lllllllllllllllll}   0 & 2 & 2 & 2 & 2 & 2 & 2 & 2 & 2 & 2 & 2 & 2 & 2 & 2 & 2 & 2 & 2 \\
0 & 0 & 0 & 2 & 2 & 2 & 2 & 2 & 2 & 2 & 2 & 2 & 2 & 2 & 2 & 2 & 2  \end{array} \right)^{\mathrm{t}}
$$
of relevance to the proof of Theorem~\ref{logsmain}, we also compute
$$
\tau^{\flat}(\mathbf{b}) = (2+2) - \frac{2\cdot 1^2 + 2 \cdot 3^2}{17^2}=\frac{1136}{289}. 
$$
\end{example}

\subsection{The convexity enhancement of the Bost--Charles bound}\label{sec_BCconvexity}

We follow the same outline as in~\S~\ref{one variable slopes}, working with the same evaluation module $\ovE_D$  and using the same non-archimedean 
evaluation heights estimate. The improvement is from the optimal use of the dilated maps $\varphi_r(z) := \varphi(rz)$ at the step of the archimedean evaluation height estimate.

In order to estimate $\psi_D^{(n)}$ at the vanishing filtration jumps $n\in \mathcal{V}_D$, we consider $Q_i \in \Z[x^{-1}]_{\leq D} = \Gamma(\P^1_\Z, \OL(D))$ such that $s:= \sum_{i=1}^m f_i Q_i = c_n x^n +\ldots$ has exact vanishing order~$n$ at $x=0$. We can view~$s$ as a formal section of $\OL(D)$, and then $s(x)\cdot x^D$ is canonically a formal function. 
Recall that we have endowed $E_D$ with a Euclidean norm induced by the Hermitian line bundle $\overline{\cL}=\overline{\OL(1)}$, on which the Hermitian metric is induced by $\varphi_* \log^+ |z|^{-1}$. By extension, we define $\overline{\cL}_r$ to be the line bundle $\OL(1)$ equipped with the Hermitian metric induced from $\left( \varphi_{r}\right)_* \log^+ |z|^{-1}$. Explicitly: 
\[\|\mathbf{1}(y)\|_{\ovcL_r} := \exp\left(-\sum_{z\in \varphi_r^{-1}(y)}\log^+{\frac{1}{|z|}} \right)=\prod_{z\in D(0,r), \, \varphi(z)=y} |z/r|. \]

As a generalization of~\eqref{BC-arch-ht-bound}, we have the following archimedean evaluation height estimate in which we can take the optimal radius
parameter $r = r(n)$: 
 
\begin{lemma}\label{arch-ht-BCconv}
For any $0<r\leq 1$, we have  
\[h_\infty(\psi_D^{(n)}) \leq -n \log|\varphi'_r(0)| + D(\ovcL \cdot \ovcL_r) +o(D).\]
\end{lemma}
\begin{proof}
We assumed the functions $\varphi^* f_i(z)=f_i(\varphi(z))$ to be meromorphic on $|z|<1$. Let us firstly remark that we can reduce the proof to the stronger assumption that $\varphi^* f_i \in \mathcal{M}(\Db)$, namely that $f_i(\varphi(z))$ is meromorphic on an open neighborhood of $|z|\leq 1$, for all $1\leq i \leq m$. Indeed, in the following proof, for $r<1$, we only use the assumption that the $\varphi^* f_i$ are meromorphic on an open neighborhood of $|z|\leq r$; therefore we only need to discuss the reduction step for $r=1$. Of course this particular $r=1$ case is indeed the estimate by Bost and Charles recalled in \eqref{BC-arch-ht-bound}. Nevertheless, we spell out a limit argument for deducing the $r=1$ case from the $r<1$ case, for the same reduction can be applied in the proofs in \S\S~\ref{new slopes}--\ref{slopes} to allow us to assume $\varphi^* f_i \in \mathcal{M}(\Db)$. Note that
\[\lim_{r\rightarrow 1} \log|\varphi'_r(0)| = \log|\varphi'(0)|, \qquad \lim_{r\rightarrow 1} \ovcL \cdot \ovcL_r = \ovcL \cdot \ovcL.\]
Therefore, the $r \to 1^-$ limit of the inequality on $h_\infty(\psi_D^{(n)}) $ gives
\[h_\infty(\psi_D^{(n)}) \leq -n \log|\varphi'(0)| + D(\ovcL \cdot \ovcL) +o(n+D).\]
This gives the desired inequality with $r=1$ since $o(n+D)=o(D)$ by Theorem~\ref{KolchinSolved} and Lemma~\ref{f_holonomic}.

And so we start from the meromorphy $\varphi^* f_i \in \mathcal{M}(\Db)$ of all pullbacks. Choose and
fix a holomorphic function $h \in \mathcal{O}(\Db)$ such that $h(0) = 1$ and all $h \cdot \varphi^* f_i \in \mathcal{O}(\Db)$
are holomorphic. We follow the notation $h_r(z) := h(rz)$ for $r \in (0,1]$ and $z \in \Db$. Then, for any $s = \sum_{i=1}^m f_iQ_i = c_nx^n +\ldots$ 
 as above, $z^{-n} h_r(z) \cdot \varphi_r^* \left(s(x)\cdot x^D \right) \in \mathcal{O}(\Db)$ is a holomorphic function whose $z = 0$ value equals $c_n\varphi'_r(0)^n \neq 0$. Therefore $\log{\left|z^{-n} h_r(z) \cdot \varphi_r^* \left(s(x)\cdot x^D \right)\right|}$ is a subharmonic function on $\Db$. 

We modify the computation in \cite[\S 10.5.5]{BostBook}. Instead of using $\varphi$ as in~{\it loc. cit.}, we apply the Poisson--Jensen formula  ---  or the subharmonic property  ---  to 
this subharmonic function $\log |z^{-n} h_r\varphi_r^* (s(x)\cdot x^D)|$. This gets us the upper bound
\begin{equation} \label{PJ est}
\begin{aligned}
\log |c_n|  & \leq -n \log |\varphi'_r(0)| + \int_\T  \log |h_r\varphi_r^* (s(x)\cdot x^D)| \, \mv \\
& = -n \log |\varphi'_r(0)| + \int_\T  \log |\varphi_r^* (s(x)\cdot x^D)| \, \mv + O(1).
\end{aligned}
\end{equation} 
Now we claim the identity
\begin{equation}  \label{intersection number}
D (\overline{\cL} \cdot \overline{\cL}_r)  =  - \int_\T \log \|\varphi_r^* x^{-D}\|_{\varphi_r^* \overline{\cL}^{\otimes D}}\mv. 
\end{equation}

To prove it, 
we start from the Poincar\'e--Lelong formula that gives:
\[
\begin{aligned}
\frac{i}{\pi} \partial \overline{\partial} \log \|\varphi_r^*(x^{-D})\|_{\varphi_r^*\overline{\cL}^{\otimes D}}
& = - c_1\left(\varphi_r^*\overline{\cL}^{\otimes D}\right), \\
 \frac{i}{\pi} \partial \overline{\partial} \log^+ |z|^{-1} &  = -\delta_0 + \mv. \end{aligned}
 \]
Therefore, by the Green--Stokes formula, we find
\begin{equation*}
\begin{aligned}
& - \int_\T \log \|\varphi_r^* x^{-D}\|_{\varphi_r^* \overline{\cL}^{\otimes D}}\mv  \\
&=  \int_{\Db} -\log \|\varphi_r^* x^{-D}\|_{\varphi_r^* \overline{\cL}^{\otimes D}} \frac{i}{\pi} \partial \overline{\partial} \log^+ |z|^{-1} 
  -\log \|\varphi_r^* x^{-D}\|_{\varphi_r^* \overline{\cL}^{\otimes D}}|_{z=0}  \\
&=  \int_{\Db} - \frac{i}{\pi} \partial  \overline{\partial}\log \|\varphi_r^* x^{-D}\|_{\varphi_r^* \overline{\cL}^{\otimes D}}  \log^+ |z|^{-1} + \ardeg \left(\varphi_r^*\overline{\cL}^{\otimes D}|_{z=0} \right)\\
&=  \int_{\Db} \log^+ |z|^{-1} c_1(\varphi_r^*\overline{\cL}^{\otimes D}) + \ardeg \left( \overline{\cL}^{\otimes D}|_{x=0}  \right)\\
&= D \cdot \left(  (\iota, \varphi_r)^*\overline{\cL} \cdot  ([0], \log^+|z|^{-1}) \right) \\
&=  D \cdot \left(  \overline{\cL} \cdot  (\iota, \varphi_r)_* ([0], \log^+|z|^{-1}) \right) \\
&=  D \cdot \left(  \overline{\cL} \cdot  \overline{\cL}_r \right).
\end{aligned}
\end{equation*}
proving~\eqref{intersection number}. 

At this point, by~\eqref{PJ est} and the pointwise decomposition 
\begin{equation}  \label{dec integrand}
 \log |\varphi_r^* (s(x)\cdot x^D)|  =   \log \|\varphi_r^* s\|_{\varphi_r^* \overline{\cL}^{\otimes D}} - \log \|\varphi_r^* x^{-D}\|_{\varphi_r^* \overline{\cL}^{\otimes D}}
\end{equation}
inside the $\T$ integrands, the lemma follows if we prove an $o(D)$ bound on the $\T$ integral of the first term on the right-hand side of~\eqref{dec integrand}
under the assumption $\| s\| \leq 1$. 

If in place of the integration measure $\mv$ we had a continuous measure~$\mu$ on~$\Db$, we would have had a constant $C$ such that $C \varphi^* \nu \geq \mu$, and then we would have had
\begin{equation*}
\begin{aligned}
0  & \geq \log \|s \|  \geq  \max_{1\leq i \leq m} \log \| Q_i\| = \frac{1}{2} \max_{1\leq i \leq m} \log{ \int_{\cX(\C)} \|Q_i\|_{\ovcL}^2 \, \nu } \\
 & \geq  C' +  \frac{1}{2} \log{ \int_{\Db} \|\varphi^*s\|_{\varphi^*\ovcL}^2 \, \mu } \geq C' + \int_{\T} \log \| \varphi^*s\|_{\varphi^*\overline{\cL}} \, \mu, 
\end{aligned}
\end{equation*}
where $C'$ is a constant depending only on $m, f_i, C$, and is independent of $D$. 
We obtain the desired bound with $\mv$ in place of~$\mu$ upon approximating $\log^+{\frac{1}{|z|}}$ (the Green function of~$\mv$) by smooth Green functions. See \S~\ref{sec:slopes archimedean} for details.
\end{proof}

Next, we explicitly calculate the arithmetic intersection number~$\left( \overline{\cL} \cdot  \overline{\cL}_r \right)$: 
\begin{lemma}\label{BCconv-intersection}
For $0<r\leq 1$, we have
\[ \left( \overline{\cL} \cdot  \overline{\cL}_r \right) = \int_{\T^2} \log |\varphi(z) - \varphi(rw)| \, \mv(z) \mv(w).\]
\end{lemma}

\begin{proof}
Recall from the discussion above that we have, straight from the definition, 
\[ \left( \overline{\cL} \cdot  \overline{\cL}_r \right) = \int_{\Db} \log^+ |z|^{-1} c_1(\varphi_r^*\overline{\cL}) +  \ardeg(\overline{\cL}|_{x=0}).
\]
From the definition of $\overline{\cL}$ and the functorial behavior of the Chern form under pushforward \cite[Proposition~3.4.5(2)]{BostCharles}, we have
\[c_1(\varphi_r^*\overline{\cL})= \varphi_r^*\varphi_* \mv, \qquad  \left( \overline{\cL} \cdot  \overline{\cL}_r \right) =\int_{\Db} \log^+ |z|^{-1} \varphi_r^*\varphi_* \mv + \ardeg(\overline{\cL}|_{x=0}).\]

Recall again from the discussion above, and from $\| \mathbf{1}(x) \|_{\overline{\cL}}=\prod_{z\in \Db,\varphi(z)=x} |z|$ and using the Poisson--Jensen formula, that 
\begin{equation*}
\begin{aligned}
\ardeg(\overline{\cL}|_{x=0}) & = - \log\| x^{-1} \|_{\ovcL}  \, \big|_{x=0}=-\log \left((\varphi(z))^{-1}\prod_{z\in \Db, \, \varphi(z)=x} |z| \right) \Big|_{x=0} \\
&=  \log |\varphi'(0)| + \sum_{0\neq z \in \Db, \, \varphi(z) =0} \log |z|^{-1}\\
&=  \int_{\T}   \log |\varphi(z) | \mv(z),
\end{aligned}
\end{equation*}
where both the product and sum 
count with multiplicities. 

We follow the same computation as in \cite[\S 5.4 and Example 5.3.2.1]{BostCharles}. Note that on $\C^2$ (with coordinates $x,y$), we have
\[ \frac{i}{\pi} \partial \overline{\partial} \log | x - y |^{-1} = - \delta_{\Delta(\C)},\]
where $\Delta(\C)$ denotes the diagonal divisor on $\C^2$.

For every $z\in \T \subset \Db$, we have (here we also view $\varphi(z)$ as the constant function that maps all points on $\Db$ to $\varphi(z)$)
\[
\begin{aligned}
\varphi_r^* \varphi_* \delta_z  & = \varphi_r^* \delta_{\varphi(z)} = (\varphi_r, \varphi(z))^*\delta_{\Delta(\C)}  \\
& = (\varphi_r, \varphi(z))^* \frac{i}{\pi} \partial \overline{\partial} \log | x - y | =  \frac{i}{\pi} \partial \overline{\partial} \log | \varphi_r(w) - \varphi(z) |. \end{aligned} \]

Therefore for a fixed $z$, by using the Green--Stokes formula (it is alright here even though the Green functions are not smooth like in \cite{BostCharles}), we have
\begin{equation*}
\begin{aligned}
\int_{\Db} \log^+|w|^{-1} \varphi_r^* \varphi_* \delta_z & =  \int_{\Db} \log^+|w|^{-1} \frac{i}{\pi} \partial \overline{\partial} \log | \varphi_r(w) - \varphi(z) | \\
&= \int_{\Db} \left( \frac{i}{\pi} \partial \overline{\partial} \log^+ |w|^{-1} \right) \cdot  \log | \varphi_r(w) - \varphi(z) | \\
&=  \int_{\Db} \left(\mv(w) - \delta_{w=0} \right) \cdot  \log | \varphi_r(w) - \varphi(z) | \\
&= \int_\T \log | \varphi_r(w) - \varphi(z) | \mv(w) - \log |\varphi(z) |.
\end{aligned}
\end{equation*}

We now integrate over $z$ on $\T$ and then we have
\[\begin{aligned}
\int_{\Db} \log^+|w|^{-1} \varphi_r^* \varphi_* \mv  = &  \int_{\T^2} \log | \varphi_r(w)   - \varphi(z) | \mv(z) \mv(w)\\
&  - \int_{\T}   \log |\varphi(z) | \mv(z). \end{aligned} \]

We arrive at the claimed formula
\[
\begin{aligned}
\left( \overline{\cL} \cdot  \overline{\cL}_r \right)
 & = \int_{\Db} \log^+ |z|^{-1} \varphi_r^*\varphi_* \mv  +  \int_{\T}   \log |\varphi(z) | \mv(z) \\
  & =\int_{\T^2} \log | \varphi_r(w) - \varphi(z) | \mv(z) \mv(w). \end{aligned} \]
\end{proof}

\begin{proof}[Proof of Theorems~\ref{main:BC conv discrete} and~\ref{main:BC convexity}]
By Lemma~\ref{BC convex}, we have that $\widehat{T}(r, \varphi)$ as a function in $\log r$ is nondecreasing and convex. Therefore it suffices to prove  Theorem~\ref{main:BC conv discrete}. 
 By Lemmas~\ref{arch-ht-BCconv} and \ref{BCconv-intersection}, using the slopes notation~\eqref{slopes alpha}, we have that for $n/D\in [\alpha_k, \alpha_{k+1}]$, the optimal bound for $h_\infty(\psi_D^{(n)})$ is obtained by using $r_k$ among $0\leq k \leq l$; here we set $\alpha_0=0, \alpha_{l+1}=m$. Note, once again, that by Corollary~\ref{holonomic criterion} and Theorem~\ref{KolchinSolved} we have  (for any $\varepsilon > 0$ and $D \gg_{\varepsilon } 1$) the containment $\mathcal{V}_D \subset [0, (m+\varepsilon)(D+1) + C(\varepsilon)]$. For $n \in [mD, (m+\varepsilon)(D+1) + C(\varepsilon)]$, we use the bound for $h_\infty(\psi_D^{(n)})$ obtain with the full radius $r_l=1$. Letting $\varepsilon \to 0$, by a similar argument to the proof of Theorem~\ref{main:BC form} (see the proof of \eqref{BC-sum-finite-ht}), we deduce from Lemmas~\ref{arch-ht-BCconv} and \ref{BCconv-intersection} and a straightforward computation that
\begin{equation*}
\begin{aligned}
& \sum_{n\in \mathcal{V}_D} h_{\mathrm{fin}}(\psi_D^{(n)}) \\
 & \leq  \left( - \frac{m^2}{2} \log |\varphi'(0)|+ \sum_{k=0}^l (\alpha_{k+1} - \alpha_k)\widehat{T}(r_k, \varphi) - \frac{1}{2}(\alpha_{k+1}^2 - \alpha_k^2) \log r_k \right)D^2 +o(D^2)\\
&= \left( - \frac{m^2}{2} \log |\varphi'(0)|+ m \widehat{T}(1, \varphi) - \frac{1}{2}\sum_{k=1}^l \alpha_k^2 ( \log r_k - \log r_{k-1}) \right)D^2 +o(D^2). 
\end{aligned}
\end{equation*}

Then the desired bound follows this estimate combined with~\eqref{BC-sum-finite-ht} and~\eqref{BCardegED}. 
\end{proof}

\begin{example}\label{BCconv-pfA}
In the proof of Theorem~\ref{mainA}, we use $\varphi$ as in \S~\ref{numerical integration}, where we have $\left( \overline{\cL} \cdot \overline{\cL}  
\right) = 11.844\ldots$ and $m=14$. 

We first apply Theorem~\ref{main:BC conv discrete} with $l=1$ and $r_0=e^{-1/2}, r_1=1$. We compute 
\[\left( \overline{\cL}_{e^{-1/2}} \cdot \overline{\cL} \right) = 10.5739\ldots,\]
 and the slope 
 \[\alpha_1=\frac{ \left( \overline{\cL} \cdot \overline{\cL}  \right) - \left( \overline{\cL}_{e^{-1/2}} \cdot \overline{\cL} \right) }{\log r_1 - \log r_0}= 2.5410\ldots.\]
Therefore, the convexity saving is 
\[\frac{  \frac{1}{m} \alpha_1^2 (\log r_1 - \log r_0)}{\log |\varphi'(0)| - \tau(\bb; \be)}=0.27243\ldots.\]
In other words, we obtain the proof by contradiction with 
\begin{equation}
\label{thirdbest}
\frac{ \left( \ovcL \cdot \ovcL\right) -  \frac{1}{m} \alpha_1^2 (\log r_1 - \log r_0)}{\log |\varphi'(0)| - \tau(\bb; \be)}=13.99303\ldots - 0.27243\ldots = 13.7206\ldots < 14.
\end{equation}

We refine this further by taking more radii: 
$$
r_0=e^{-1}, \quad r_1=e^{-1/2},  \quad r_2=e^{-1/4},\quad  r_3=1.
$$
 At these radii, we compute the corresponding Bost--Charles characteristic integrals: 
\[
\begin{aligned}
\widehat{T}(r_3, \varphi)  & = \left( \ovcL, \ovcL \right) =11.844\ldots, \\
 \widehat{T}(r_2, \varphi) & =  \left( \ovcL_{e^{-1/4}}\cdot \ovcL \right) = 11.049\ldots, \\
\widehat{T}(r_1, \varphi)  & =   \left( \ovcL_{e^{-1/2}}\cdot \ovcL \right) =10.573\ldots,  \\
  \widehat{T}(r_0, \varphi) & =  \left( \ovcL_{e^{-1}} \cdot \ovcL \right) =9.8766\ldots; \end{aligned}\]
and the corresponding slopes:
\[\alpha_1= 1.3943\ldots, \alpha_2=1.9018\ldots, \alpha_3=3.1802\ldots.\]
We thus derive the following convexity saving in the holonomy bound: 
\[\frac{  \frac{1}{m} \sum_{k=1}^3 \alpha_k^2 (\log r_k - \log r_{k-1})}{\log |\varphi'(0)| - \tau(\bb; \be)}=0.37171\ldots;\]
In other words, the refined holonomy bound is 
\begin{equation}
\label{best}
13.99303\ldots - 0.37171\ldots = 13.621\ldots
\end{equation}
\end{example}

\begin{example}\label{Ex_BCconv9}
In order to prove Theorem~\ref{mainA}, as discussed in Remark~\ref{improvement}, we could try to only use the $9$ functions (without integrations). Then Theorem~\ref{main:BC conv discrete} with $l=3$ in the above example gives
\[\frac{11.844\ldots - \frac{1}{9}\sum_{k=1}^3 \alpha_k^2 (\log r_k - \log r_{k-1})}{\log \left( 256  \cdot
\frac{5448339453535586608000000000}{8658833407565631122430056127}
\right)
- 2 \cdot 157/81}=9.4203\ldots <10,\]
which comes nearer to the~$9$ threshold but it remains insufficient to draw a contradiction with $9$ functions. This is why 
we need the integrations idea. 
\end{example}

\subsection{Binomial metrics: proof of Theorem~\ref{main: easy convexity}}  \label{sec: proof easy convexity}
We recall our assumption on the denominator types of $\{f_i\}$. Set $u_0 :=0$ and $u_{r+1} :=m$. For $0\leq h \leq r$, if $u_{h} < i \leq u_{h+1}$, then 
\[f_i(x) = a_{i,0}+ \sum_{n=1}^{\infty} a_{i,n} \frac{x^n}{n^{e_i} [ 1, \ldots, b_{1} \cdot n] \cdots [1,\ldots, b_{h} \cdot n]}, \qquad a_{i,n} \in \Z. \]
We take our evaluation module to be the following free $\Z$-module of rank $mD$: 
\[E_D=\bigoplus_{h=0}^r \bigoplus_{i=u_h+1}^{u_{h+1}} \frac{[1,\ldots, \xi D]^{e_i}}{[ 1, \ldots, u_{h+1}b_{h+1} D] \cdots [1,\ldots, u_rb_{r} D]} f_i \, \Z[x]_{<D}.\]
We endow $E_D$ with the Euclidean norm that has $\{f_i x^k\}_{1\leq i \leq m, 0 \leq k <D}$ as an orthogonal basis with vector lengths 
$\|x^k\|=e^{D(\lambda t^r + \mu t)}$, where $t=k/D$. Recall from our assumption that $\lambda>0$ and $r> 1$. We use $\overline{E}_D$ to denote this Euclidean lattice.

Applying the defining formula~\eqref{arithmetic deg} of $\ardeg{\overline{E}_D}$ to the $\Q$-basis $\{f_i x^k\}_{1\leq i \leq m, 0 \leq k <D}$ of $E_{D} \otimes \Q$, we compute, as $D\rightarrow \infty$: 
\begin{equation}
\begin{aligned}
\ardeg{\overline{E}_D} & =  -m\left(\int_{0}^1 \lambda t^r + \mu t  \, dt \right)D^2 +\left(\sum_{h=0}^{r-1} (u_{h+1}-u_h)\sum_{j=h+1}^r u_jb_j\right)D^2 \\
& \quad  - \left(\xi \sum_{i=1}^m e_i \right)D^2 +o(D^2) \\
&=    - m\left(\frac{\lambda}{r+1} + \frac{\mu}{2}\right)D^2 + \left(\sum_{h=1}^r u_h^2 b_h - \xi \sum_{i=1}^m e_i \right)D^2 +o(D^2). \label{cvx-ardeg}
\end{aligned}
\end{equation}

Again we let $F_\Q:=\Q\llbracket x\rrbracket$, and we filter it by the $x=0$ vanishing order: 
\[F_\Q=F_\Q^{(0)}\supseteq F_\Q^{(1)} \supseteq \cdots \supseteq F_\Q^{(n)} \supseteq \cdots,\]
where 
\[F_\Q^{(n)}:= \Cspan_{\Q}\{x^{k} \,  : \, k\geq n\}.\]
 The graded piece $F_\Q^{(n)}/F_\Q^{(n+1)}$ is a one dimensional $\Q$-vector space generated by the image of~$x^{n}$ under the quotient map. The Euclidean lattice structure on $F_\Q^{(n)}/F_\Q^{(n+1)}$ is given by the free rank one $\Z$-module generated by the image of~$x^{n}$ and the Euclidean norm with $\| x^{n} \|=1$. This is the same as in \S~\ref{one variable slopes} up to a shift by~$-D$ in the power of~$x$.

As in \S~\ref{one variable slopes}, we the have natural injective \emph{evaluation map} $\psi_D : E_D \rightarrow F_\Q$, inducing injections on the graded pieces $\psi_D^{(n)}: {E}^{(n)}_{D}/ {E}^{(n+1)}_{D}\rightarrow F_\Q^{(n)}/F_\Q^{(n+1)}$.
We still have $\rank {E}^{(n)}_{D}/{E}^{(n+1)}_{D}\in \{0,1\}$, and the cardinality $\# \mathcal{V}_{D} = \rank E_D=mD$ of the vanishing filtration jumps set $\mathcal{V}_{D} = \left\{n\in \NwithzeroB \, : \, \rank {E}^{(n)}_{D}/{E}^{(n+1)}_{D}=1 \right\}$.\\

We now provide upper bounds on the evaluation heights $h_\infty(\psi_D^{(n)})$ and $h_{\mathrm{fin}}(\psi_D^{(n)})$. 
 For $v\in M_\Q$, by the definition of the local evaluation height $h_v$, we consider an arbitrary $(Q_i)_{1\leq i\leq m} \in {E}^{(n)}_{D} \setminus {E}^{(n+1)}_{D}$, and our task is to provide an upper bound on $\log | c_{n} |_v - \log \| ( Q_i)_{1\leq i\leq m}\|_{E_D,v}$, where $c_{n}$ denotes the coefficient of  $x^{n}$ in $\sum_{i=1}^m f_i(x)Q_i(x)$,  and $|\cdot |_v$ is the usual $v$-adic norm on $\Q$. 

For $v=\infty$, we use the equivalent interpretation of $h_\infty(\psi_D^{(n)})$ by considering any $(Q_i)_{1\leq i\leq m} \in {E}^{(n)}_{D,\R} \setminus {E}^{(n+1)}_{D,\R}$ with $\| (Q_i)_{1\leq i\leq m}\|_{E_D,\infty} \leq 1$, and providing an upper bound on $\log | c_{n} |_\infty$, which is then our upper bound for $h_\infty(\psi_D^{(n)})$. By definition of our binomial metric, upon writing momentarily $t:=k/D$, the unit ball condition $\| ( Q_i)_{1\leq i\leq m}\|_{E_D,\infty} \leq 1$ implies the bounds $|\alpha_{i,k}|_\infty \leq e^{-D(\lambda t^r + \mu t)}$ on the coefficients of $Q_i(x)=\sum_{k=0}^{D-1} \alpha_{i,k} x^k$,
 for all $1\leq i \leq m$. 

For simplicity of notation, we write $f_i(x)=\sum_{n=0}^\infty a'_{i,n} x^n$. By assumption, 
all $f_i$ converge on the closed disc $\overline{D(0,\rho)}$ for all $\rho <R$. We use this information to derive an upper bound on $h_\infty(\psi_D^{(n)})$ which is useful on a certain range of~$n/D$. The analyticity on $\Db_{\rho}$ 
means that $|a'_{i,k}|_\infty =O_\rho(\rho^{-k})$, where the implicit constant only depends on $\rho, f_1,\ldots, f_m$, but not on~$k$. Hence, for arbitrary $n\in \NwithzeroB$, we derive from 
\[c_n = \sum_{i=1}^m \sum_{k=0}^{\min\left( n,D-1\right)} \alpha_{i,k} a'_{i, n-k}, \text{ thus } |c_n|_\infty \leq mD \max_{ 
\substack{ 1\leq i \leq m, \\ 0 \leq k \leq \min{(n,D-1)} }} |\alpha_{i,k}|_\infty \cdot  |a'_{i, n-k}|_\infty.\]
the following archimedean evaluation heights bound: 
\begin{equation}  \label{pre bound}
h_\infty(\psi_D^{(n)}) \leq \left(\max_{0\leq t \leq \min\{1,n/D\}} \left\{ -\lambda t^r - \mu t - (n/D - t) \log \rho \right\} \right)D +o_\rho(D).
\end{equation}

The function of~$t \in [0, \min\{1,n/D\}]$ under the maximum is concave, and in particular unimodal. 
From here it is easy to justify the $\rho \rightarrow R^-$ limit:
\begin{equation}  \label{R bound}
h_\infty(\psi_D^{(n)})\leq \left( \max_{0\leq t \leq \min\{1,n/D\}} \left\{ -\lambda t^r - \mu t - (n/D - t) \log R \right\} \right)D +o(n+D).
\end{equation}
We include the details of this limiting argument as it also reveals the limit point $\rho = R^-$ to indeed be the optimal choice to make in~\eqref{pre bound} across $\rho \in (0,R)$.
Let us denote the maximizers of the unimodal functions under the curly brackets in~\eqref{R bound} and~\eqref{pre bound} to be at $t := t_R$ and $t :=t_{\rho}$, respectively. 
We have, noting that by definition $0 \leq t_R, t_{\rho} \leq \min\{1,n/D\}$:  
\begin{equation*}
\begin{aligned}
& \max_{0\leq t \leq \min\{1,n/D\}} \left\{ -\lambda t^r - \mu t - (n/D - t) \log \rho \right\}  \geq  -\lambda t_R^r - \mu t_R - (n/D - t_R) \log \rho \\
 & =  (\log R - \log \rho)(n/D -t_R)  +  \max_{0\leq t \leq \min\{1,n/D\}} \left\{ -\lambda t^r - \mu t - (n/D - t) \log R \right\}\\
&  \geq   \max_{0\leq t \leq \min\{1,n/D\}} \left\{ -\lambda t^r - \mu t - (n/D - t) \log R)\right\},
\end{aligned}
\end{equation*}
and similarly, 
\begin{equation*}
\begin{aligned}
& \max_{0\leq t \leq \min\{1,n/D\}} \left\{ -\lambda t^r - \mu t - (n/D - t) \log R\right\}  \geq  -\lambda t_{\rho}^r - \mu t_{\rho} - (n/D - t_{\rho}) \log{R}  \\
&  =   - (\log R - \log \rho)(n/D -t_{\rho}) + \max_{0\leq t \leq \min\{1,n/D\}} \left\{ -\lambda t^r - \mu t - (n/D - t) \log \rho \right\} \\
&  \geq  - (\log R - \log \rho)(n/D) + \max_{0\leq t \leq \min\{1,n/D\}} \left\{ -\lambda t^r - \mu t - (n/D - t) \log \rho \right\}.
\end{aligned}
\end{equation*}

This proves~\eqref{R bound}, and also the optimality of taking the limit $\rho \to R^-$ in~\eqref{pre bound}. 

Continuing with the proof, we set $s:=n/D$, and recall our notation 
$$
\chi_0:=\min\left\{1, \left(\frac{\max\{0, \log R - \mu\} }{\lambda r}\right)^{1/(r-1)}\right\}
$$
from the statement of the theorem under proof. Its meaning is the following. 
 By the same computation as in the proof of Lemma~\ref{T-formula}, the following archimedean evaluation height bound is 
 valid in the range $s\geq \chi_0$: 
\begin{equation}\label{eq:arch_id}
h_\infty(\psi_D^{(n)})\leq \big(\Gamma(\log R, r, \lambda, \mu) - s \log R\big)D +o(n+D);
\end{equation}
whereas in the range $0\leq s\leq \chi_0$, the following improvement holds: 
\begin{equation}\label{eq:arch_id2}
h_\infty(\psi_D^{(n)})\leq (-\lambda s^r - \mu s)D +o(n)+o(D).
\end{equation}
Saying that $s \in [0,\chi_0]$ is the domain of improvement of~\eqref{eq:arch_id2} over~\eqref{eq:arch_id} is exactly the definition of~$\chi_0$.

For the range $s > \chi_0$, we instead use the Poisson--Jensen formula applied to the logarithm of the holomorphic function
$$
h(z) \cdot \varphi^*\left(\sum_{i=1}^m f_iQ_i \right) \cdot z^{-n} \in \mathcal{O}(\Db),
$$
where an arbitrary $h \in \mathcal{O}(\Db)$ is fixed subject to $h(0) = 1$ and $h \cdot \varphi^* f_i \in \mathcal{O}(\Db)$ for 
all~$i = 1, \ldots, m$. 
 We derive the usual bound (here we also use $|\cdot |_\infty$ to denote the usual absolute value on $\C$): 
\begin{equation*}
\begin{aligned}
& \log |c_n|_\infty \\
 & \leq  - n \log |\varphi'(0)| +\int_{\T} \left|h \cdot \varphi^*\left(\sum_{i=1}^m f_iQ_i\right)\right|_\infty \,\mv\\
&\leq  - n \log |\varphi'(0)| + \int_{\T}\max_{1\leq i \leq m, 0\leq k \leq \min\{n, D-1\}} (\log |\alpha_{i,k}|_\infty + k\log|\varphi(z)|_\infty ) \,\mv+o(D)\\
&\leq  - n \log |\varphi'(0)| + \int_{\T}\max_{0\leq k \leq \min\{n, D-1\}} (D(-\lambda(k/D)^r-\mu(k/D)) + k\log|\varphi(z)|_\infty ) \,\mv+o(D).
\end{aligned}
\end{equation*}
 Therefore, by the definition of $T(\varphi;r,\lambda, \mu)$ via the Legendre transform~$\Gamma(x;r,\lambda,\mu)$ of the binomial metric weight function $\lambda t^r+\mu t$, we derive the following for our upper bound on all the archimedean evaluation heights: 
\begin{equation}\label{arch_phi}
h_\infty(\psi_D^{(n)}) \leq - n \log |\varphi'(0)| +D  T(\varphi;r,\lambda, \mu) +o(D).
\end{equation}
\silentcomment{As above, for $n/D\leq 1$, one may improve the bound by observing that we only need to consider $0\leq t \leq n/D$; however, in our application, we will always use this bound for $n/D\geq 1$, so we will not work out 
the exact formula the improved bound for the $n/D\leq 1$.}
Note that bound \eqref{eq:arch_id} is better than \eqref{arch_phi} if and only if
\[\frac{n}{D} \leq \chi_1:= \frac{T(\varphi;r,\lambda, \mu)- \Gamma(\log R; r, \lambda, \mu) }{\log |\varphi'(0)|- \log R}.\]
Hence, by Theorem~\ref{KolchinSolved} (letting  $\varepsilon\rightarrow 0$ right after taking $D \to \infty$), we have 
\begin{equation}
\begin{aligned}
& \limsup_{D \to \infty} \left\{ 
D^{-2} \sum_{n\in \mathcal{V}_D} h_\infty(\psi_D^{(n)})  \right \} \\ 
& \leq   \int_{0}^{\chi_0} (-\lambda s^r - \mu s)\,ds  +\int_{\chi_0}^{\chi_1}\big( \Gamma(\log R; r, \lambda, \mu) - s \log R \big)\,  ds \\
& \quad  +  \  \int_{\chi_1}^m \big(T(\varphi;r,\lambda, \mu) - s \log |\varphi'(0)| \big) \,ds 
\\
 & =   m T(\varphi;r,\lambda, \mu) - \frac{m^2}{2} \log |\varphi'(0)| - \frac{(T(\varphi;r,\lambda, \mu) - \Gamma(\log R; r, \lambda, \mu))^2}{2(\log |\varphi'(0)|  - \log R)} 
 \\
& \quad  -    \chi_0\Gamma(\log R; r, \lambda,\mu) + \chi_0^2 (\log R - \mu)\left(\frac{1}{2} - \frac{1}{r(r+1)}\right). \label{cvx-sum-arch-ht}
\end{aligned}
\end{equation}

Next we turn to estimating $h_{\mathrm{fin}}(\psi_D^{(n)})$. Considering an arbitrary $( Q_i)_{1\leq i\leq m} \in {E}^{(n)}_{D} \setminus {E}^{(n+1)}_{D}$, our task is for each prime $p$ to provide an upper bound on $\log | c_{n} |_p$. Since the $\Z$-lattice $E_D$ here has essentially the same structure as the one in \S~\ref{one variable slopes}, the argument there yields for the total finite evaluation height the upper bound: 
\begin{equation}\label{cvx-sum-finite-ht}
\limsup_{D \to \infty} \left\{ D^{-2} \sum_{n\in \mathcal{V}_D} h_{\mathrm{fin}}(\psi_D^{(n)})  \right\} 
\leq  \frac{1}{2} \left(\sigma_m m^2 + \sum_{h=1}^r  u_h^2b_h \right) + \left(\max_{1\leq i \leq m} e_i\right) I_{\xi}^m(\xi). 
\end{equation}

We plug \eqref{cvx-ardeg}, \eqref{cvx-sum-arch-ht}, and \eqref{cvx-sum-finite-ht} into \eqref{slope-inequality} and derive: 
\begin{equation*}
\begin{aligned}
& \,  \frac{m^2}{2}\left(\log|\varphi'(0)| - \sigma_m +\frac{1}{m^2}\left(\sum_{h=1}^r u_h^2 b_h  - \frac{2}{m^2} (\xi \sum_{i=1}^m e_i + \left(\max_{1\leq i \leq m} e_i\right)I_{\xi}^m(\xi)) \right) \right)\\
  \leq & m\left( T(\varphi;r,\lambda, \mu)+ \frac{\lambda}{r+1} + \frac{\mu}{2} \right)  - \frac{(T(\varphi;r,\lambda, \mu) - \Gamma(\log R; r, \lambda, \mu))^2}{2(\log |\varphi'(0)|  - \log R)} \\ & \, 
-\chi_0\Gamma(\log R; r, \lambda,\mu) + \chi_0^2 (\log R - \mu)\left(\frac{1}{2} - \frac{1}{r(r+1)}\right),
\end{aligned}
\end{equation*}
which rearranges into the claimed bound on $m$.   $\hfill{\square}$

\begin{example}\label{Ex-easyconv}
In the proof of Theorem~\ref{mainA}, we use $\varphi$ as in \S~\ref{numerical integration}, and recall that
\[ \begin{aligned}
m & = 14, \\
\log |\varphi'(0)|  & =\log \left( 256  \cdot
\frac{5448339453535586608000000000}{8658833407565631122430056127}
\right), \\
 \tau(\bb;\be) & = \frac{27}{80}  + \frac{191}{49}, \end{aligned}\]
and that all $f_i$ have convergence radii at least $R:=4$.

Select the following for the binomial metric weight parameters: 
\[r=4.7, \lambda=10, \mu= -4.5.\]
A numerical computation gives
\begin{equation*}
\begin{aligned}
T(\varphi;4.7,10, -4.5) & =6.5316\ldots, \\
\Gamma(\log 4; 4.7,10, -4.5) & =2.6429\ldots,  
\end{aligned}
\end{equation*}
with
$$
\chi_1= 1.0522\ldots, \quad \chi_0=0.57035\ldots,
$$
meeting the special assumptions that we made in Theorem~\ref{main: easy convexity}, and supplying the holonomy bound $m \leq 13.8527\ldots <14$. The contradiction supplies a proof of Theorem~\ref{mainA} (see~\S~\ref{sec:proofA}), with a better numeric than when we use Theorem~\ref{main:BC form} alone like in~\S~\ref{numerical integration}, prior to the convexity enhancement by Theorem~\ref{main:BC convexity}.

If we only work with $9$ functions as in Remark~\ref{improvement}, with the same parameters above replacing $m=9, \tau(\bb';0)=2\cdot \frac{157}{81}$, we have the bound in Theorem~\ref{main: easy convexity} is $9.5234\ldots <10$, but not enough to draw a contradiction to deduce Theorem~\ref{mainA}.  \endofremark
\end{example}

\subsection{A further improvement} \label{sec_BCfull} 
The setup is similar to Theorem~\ref{main:BC conv discrete}: fix a set of subradii $1 = r_l > r_{l-1} > \cdots > r_0>0$. The following is the counterpart --- and ultimate sharpening --- of the archimedean term in the refined bound from~\S~\ref{fine section}.  
In place 
of using the Hermitian line bundle $\ovcL=(\iota, \varphi)_* ([0], \log^+ |z|^{-1})$,  
we introduce weights $s_0, \ldots, s_l \in [0,1]$ with total mass $\sum_{h=0}^l s_h=1$, and use the $\mathbf{s}$-weighted average of
the Hermitian line bundles defined by the restricted maps $\varphi(r_k z)$: 
$$
\ovcL' := \prod_{h=0}^l \ovcL_{r_h}^{\otimes s_h}.
$$
Here, by a mild abuse of $\R$-line bundle notation, 
this is the line bundle $\cL = \OL(1)$ over~$\cX$ with Hermitian metric defined by 
$$
\| \cdot \|_{\ovcL'} :=  \prod_{h=0}^l \| \cdot \|_{\ovcL_{r_h}}^{s_h}. 
$$
As in~\S~\ref{sec_BCconvexity}, for $1\leq k \leq l$, set 
\begin{equation} \label{slopy}
\beta_k:=\frac{\ovcL' \cdot \ovcL_{r_k} - \ovcL' \cdot \ovcL_{r_{k-1}} }{ \log r_k - \log r_{k-1} }= \frac{\sum_{h=0}^l s_h(\ovcL_{r_h}\cdot  \ovcL_{r_k} - \ovcL_{r_h} \cdot \ovcL_{r_{k-1}} )}{ \log r_k - \log r_{k-1} }.
\end{equation}

We note that by the same argument\footnote{These are actually the same statement upon changing~$\varphi(z)$ to~$\varphi(r_k z)$ and~$r$ to~$r_h/r_k$, 
if $h \leq k$.} as in Lemma~\ref{BCconv-intersection} we have
\[\ovcL_{r_h}\cdot  \ovcL_{r_k}=\int_{\T^2} \log |\varphi(r_hz) - \varphi(r_kw)| \, \mv(z) \mv(w).\]

Therefore, since all $s_h \geq 0$, Lemma~\ref{BC convex} on convexity shows that 
\[0 \leq \beta_1 \leq \beta_2 \leq \cdots \leq \beta_l.\]

As in Theorem~\ref{main:BC conv discrete}, we assume $\beta_l \leq m$. (If this condition fails, it usually serves as a stronger bound on $m$ anyhow.) We extend the notation by setting $\beta_0 := 0$ and $\beta_{l+1} := m$. For $n/D\in [\beta_k, \beta_{k+1})$, we estimate the archimedean evaluation height in terms of $\varphi_{r_k}$. Namely, the proof of Lemma~\ref{arch-ht-BCconv} with $\ovcL$ replaced by $\ovcL'$ gives
\[h_\infty(\psi_D^{(n)}) \leq -n \log|\varphi'(0)| - n \log r_k + D(\ovcL' \cdot \ovcL_{r_k}) +o(D).\]
Similarly to the proof of Theorem~\ref{main:BC conv discrete}, using Theorem~\ref{KolchinSolved}, we have that $\mathcal{V}_D \subset \left[0, (m+\epsilon)D\right]$ once $D \gg_{\epsilon} 1$, and for $n \in \left[mD, (m+\epsilon)D\right]$, we continue to use the bound on $h_\infty(\psi_D^{(n)})$ from taking the full radius $r_l=1$. As $D \to \infty$ and then $\epsilon \rightarrow 0$, we obtain
\[ 
\begin{aligned}
\limsup_{D \to \infty} \left\{ D^{-2} \sum_{n\in \mathcal{V}_D} h_\infty(\psi_D^{(n)}) \right\}
&  \leq   - \frac{m^2}{2} \log |\varphi'(0)| \\
& \quad + \sum_{k=0}^l (\beta_{k+1} - \beta_k)\ovcL' \cdot \ovcL_{r_k}- \frac{1}{2}(\beta_{k+1}^2 - \beta_k^2) \log r_k. \end{aligned} \] 
We have the same estimate on $h_{\mathrm{fin}}(\psi_D^{(n)})$ as in \eqref{BC-sum-finite-ht}. Finally,  
\[\ardeg \ovE_D = \left(\frac{m}{2}(\ovcL' . \ovcL') +\sum_{h=1}^r u_h y_h - \xi \Big( \sum_{i=1}^{m} e_i \Big) \right)D^{2} +o(D^{2}).\]
Hence, by the slopes inequality~\eqref{slope-inequality} as before, we get
\begin{equation} \label{BCfull-1}
\begin{aligned}
& m^2 \left( \log |\varphi'(0)| - \tau(\bb;\be) \right) \\
&  \leq     -m \ovcL' \cdot \ovcL' + \sum_{k=0}^l \left(-(\beta_{k+1}^2 - \beta_k^2) \log r_k + 2 (\beta_{k+1} - \beta_k)\ovcL' \cdot \ovcL_{r_k} \right)  
\\
&  =  2m \ovcL' \cdot \ovcL_{1} -m \ovcL' \cdot \ovcL'  - \sum_{k=1}^l \left(-\beta_k^2(\log r_k - \log r_{k-1}) +2\beta_k (\ovcL' \cdot \ovcL_{r_k} - \ovcL' \cdot \ovcL_{r_{k-1}} ) \right) 
 \\
& =  2m \ovcL' \cdot \ovcL_{1} -m \ovcL' \cdot \ovcL' -  \sum_{k=1}^l \frac{(\ovcL' \cdot \ovcL_{r_k} - \ovcL' \cdot \ovcL_{r_{k-1}} )^2}{ \log r_k - \log r_{k-1} } 
\\
& = 2m \ovcL' \cdot \ovcL_{1} -m \ovcL' \cdot \ovcL' -  \sum_{k=1}^l \beta_k (\ovcL' \cdot \ovcL_{r_k} - \ovcL' \cdot \ovcL_{r_{k-1}} )
\\
& =  m \ovcL' \cdot \ovcL_{1} -m \ovcL' \cdot \ovcL' +  \sum_{k=0}^l (\beta_{k+1} - \beta_k) \ovcL' \cdot \ovcL_{r_k}.
\end{aligned}
\end{equation}
Note that \eqref{BCfull-1} gives a bound on $m$ for every choice of partition $\mathbf{s} = \{s_h\}_{h=0}^l$, and it recovers the convexity saving of Theorem~\ref{main:BC conv discrete} as the special case $\mathbf{s} = (0,0,\ldots,0;1)$. 

We propose the following choice of $\{s_h\}_{h=0}^l$. Let us postulate the following system of~$l+1$ inhomogeneous linear equations in the~$l+1$ unknowns~$s_h$: 
\silentcomment{(we would expect this condition always by using some properties of the intersection pairing. This condition holds in our particular application; see the Example below)}
\begin{equation} \label{equi s}
s_h = \frac{1}{m} (\beta_{h+1} - \beta_h), \quad 0 \leq h \leq l, \quad \beta_0 = 0, \, \beta_{l+1} = m.
\end{equation}
This choice is explained in Remark~\ref{stationary choice} below.
We suppose the $(l+1)\times (l+1)$ coefficient matrix of this linear system to have a nonzero determinant, and furthermore that the unique solution $\{s_h^*\}$
has nonnegative components. 
This solution clearly has $\sum_{h=0}^l s_h^* =1$,
and we can set our $s_h :=s_h^*$ in \eqref{BCfull-1}. The inequality~\eqref{BCfull-1}  then reads: 
\[m^2(\log |\varphi'(0)| - \tau(\bb;\be) ) \leq m \sum_{h=0}^l s_h^* \ovcL_{r_k}\cdot \ovcL_{1}. \]
In this situation we derive the following refined holonomy bound: 
\[m \leq \frac{\sum_{h=0}^l s_h^* \cdot \ovcL_{r_k}\cdot \ovcL_{1}}{\log |\varphi'(0)| - \tau(\bb;\be)}.\]

We summarize our findings into a theorem: 
\begin{theorem}\label{main:BC fullconv}
Assume the same conditions and notation as in Theorem~\ref{main:BC form}. Fix a sequence of subradii $1 = r_l > r_{l-1} > \cdots > r_0>0$. Assume that the following system of~$l+1$ linear inhomogeneous equations in the~$l+1$ unknowns $\{s_h\}_{h=0}^l$
\begin{equation} \label{see ls}
\begin{aligned}
m \sum_{h=0}^{k-1} s_h & = \frac{\sum_{h=0}^l s_h(\ovcL_{r_h}\cdot  \ovcL_{r_k} - \ovcL_{r_h} \cdot \ovcL_{r_{k-1}} )}{ \log r_k - \log r_{k-1} }, 
\qquad  k  =  1, \ldots, l, \\
\sum_{h=0}^l s_h & =1
\end{aligned}
\end{equation}
has a unique solution $\{s_h^*\} \in [0,1]^{l+1}$.

 Then, 
\begin{equation} \label{equilibrium bound}
m \leq \frac{\sum_{h=0}^l s_h^* \ovcL_{r_k}\cdot \ovcL_{1}}{\log |\varphi'(0)| - \tau(\bb;\be)} = \frac{\ovcL_1 \cdot \ovcL_1 - \sum_{h=0}^{l-1} s_h^* (\ovcL_1 \cdot \ovcL_1 - \ovcL_{r_h}\cdot \ovcL_{1})}{\log |\varphi'(0)| - \tau(\bb;\be)}.
\end{equation} 
\end{theorem}

\begin{remark}  \label{stationary choice}
The special assumptions about the linear system~\eqref{see ls} having a unique solution with nonnegative components appears to hold in practice. We do not know if it is a general feature. 
The heuristic behind this particular choice of $s_h=s_h^*$ is to emulate the Euler--Lagrange stationary action principle on our upper bound~\eqref{BCfull-1}. Namely, we compute the $d/ds_h$ derivatives of that upper bound,
\[2m \ovcL' \cdot \ovcL_{1} -m \ovcL' \cdot \ovcL' -  \sum_{k=1}^l \frac{(\ovcL' \cdot \ovcL_{r_k} - \ovcL' \cdot \ovcL_{r_{k-1}} )^2}{ \log r_k - \log r_{k-1} }\]
and set these derivatives to~$0$.~\endofremark
\end{remark}

\begin{example}\label{Ex_BCfull}
Let us revisit now the first case in Example~\ref{BCconv-pfA}: $l=1, r_0=e^{-1/2}$, and $14$ putative functions for Theorem~\ref{mainA}. We have
\[
\begin{aligned}
C:=\ovcL_1\cdot \ovcL_1  & = 11.844\ldots \\
B:=\ovcL_1 \cdot \ovcL_{e^{-1/2}} & =10.573\ldots \\
A:=\ovcL_{e^{-1/2}}\cdot \ovcL_{e^{-1/2}} & =8.3717\ldots.\end{aligned}
\]
Now the point of the improvement over the previous bound is that
\[2 \ovcL_1 \cdot \ovcL_{e^{-1/2}} > \ovcL_1\cdot \ovcL_1 + \ovcL_{e^{-1/2}}\cdot \ovcL_{e^{-1/2}}.\]
We calculate
\[\beta_1= \frac{s_0(B-A) + s_1(C-B)}{\log r_1 - \log r_0}=2(s_0(B-A) + s_1(C-B)) > 2(C-B) = \alpha_1,\]
whence
\[s_0^*= \frac{1}{m} \beta_1(s_0^*, s_1^*) > \alpha_1 /m.\]
The new convexity saving is by $\frac{s_0^*(C-B)}{\log |\varphi'(0)| - \tau(\bb;\be)} > \frac{\alpha_1(C-B)/m}{\log |\varphi'(0)| - \tau(\bb;\be)}$, where the latter was the previous convexity saving in Example~\ref{BCconv-pfA}. 
Explicitly, we find $s_0^*$ and $s_1^*$ as the solution of the two linear equations
\[s_0=\frac{2}{m}(s_0(B-A) + s_1(C-B)), \quad s_1= 1-\frac{2}{m}(s_0(B-A) + s_1(C-B)) .\]
The solution is 
\[s_0^*=\frac{C-B}{m/2 + (A+C - 2B)}=0.20936\ldots\]
resulting in the following improvement convexity saving over
 Example~\ref{BCconv-pfA} (compare with equation~\eqref{thirdbest}): 
\[\frac{s_0^*(C-B)}{\log |\varphi'(0)| - \tau(\bb;\be)} = 0.31426\ldots (> 0.27243\ldots).\]

In other words, the refined holonomy bound on~$m$ is here 
\begin{equation} \label{secondbest}
13.99303\ldots - 0.31426\ldots = 13.678\ldots,\end{equation}
giving a still more comfortable numerical margin 
 for the ultimate contradiction to~$m=14$.   
 (A similar but more complicated computation with four
 radii instead of two would also yield a slight improvement
 on equation~(\ref{best}) in Example~\ref{BCconv-pfA}.)
 \endofremark
\end{example}

{\it  At this point, a reader primarily interested in the proof of Theorems~\ref{mainA} and~\ref{logsmain} can skip directly ahead to~\S~\ref{sec:YtoY0(2)} on a first reading. }

\subsection{On bypassing the Kolchin--Shidlovsky type theorems from \S~\ref{functional bad approximability}} \label{sec_withoutShid}  At least for our specific and qualitative applications in the present paper, it is technically possible to avoid all recourse to the --- fairly technical --- zero estimates we collected in~\S~\ref{functional bad approximability}. 
 As the general case of our abstract holonomy bounds seems somewhat awkward to approach in its full generality\footnote{We do not have such a proof.} without using the functional bad approximability theorems, while on the other hand the Shidlovsky (or Chudnovsky--Osgood) type of input is a golden standard in the subject which --- furthermore and more importantly --- turns out indispensable for all quantitative refinements in our method to deriving actual Diophantine inequalities on the bad approximability of a period vector by an integer vector, we limit ourselves here to only a few brief indications on how one could technically avoid the appeal to~\S~\ref{functional bad approximability} 
 or the purpose of proving certain relaxed versions of our holonomy bounds, still sufficient for all our present applications in this paper. 

We recall that for all the proofs in this section, we have made the assumption that $0 = u_0 \leq u_1 \leq \cdots \leq u_r  < u_{r+1}= m$ in the denominators form~\eqref{column shape}, as permutation on the columns of $\bb$ does not change the assumption on $f_i$.
 
\subsubsection{Discussion for Theorem~\ref{main:BC form}}\label{noShidBC}
We sketch a proof of~\eqref{BCbound} that bypasses~\S~\ref{functional bad approximability} under assuming the stronger positivity condition $\log |\varphi'(0)| > \sigma_m + \max (e_i)$ in place of~\eqref{stronger positivity condition}.  In our application to Theorem~\ref{mainA}, and to at least some weaker form (i.e., with $10^{-6}$ replaced by a smaller explicit positive number) of Theorem~\ref{logsmain}, this condition is satisfied since $\log |\varphi'(0)| = \log \left( 256  \cdot
\frac{5448339453535586608000000000}{8658833407565631122430056127}
\right)
> 5.08 > 4+1 = \sigma_m + \max_{i=1}^m (e_i)$ in~\S~\ref{numerical integration}. 

Let $\overline{h}_{\infty}(\psi_D^{(n)}), \overline{h}_{\mathrm{fin}}(\psi_D^{(n)})$ denote the main terms in the bounds on the archimedean, resp.~finite evaluation heights $h_{\infty}(\psi_D^{(n)}), h_{\mathrm{fin}}(\psi_D^{(n)})$ that we proved in~\eqref{BC-arch-ht-bound}, resp.~\eqref{finite>xi}, \eqref{finite-mid}, and~\eqref{finite<xi}; in~$\overline{h}_{\mathrm{fin}}$, we take the optimal parameters choices $y_h :=b_h u_h$ that we used at the end of our proof.  Thus the global evaluation height has an upper bound with main term 
$\overline{h}(\psi_D^{(n)}) =\overline{h}_{\infty}(\psi_D^{(n)}) + \overline{h}_{\mathrm{fin}}(\psi_D^{(n)})$ given by
\begin{equation}
\begin{aligned}
& - n \log | \varphi'(0)| +  D \left(\overline{\OL(1)} \cdot  \overline{\OL(1)} \right) + \left( \sum_{h=1}^r  b_h \max\{n, u_h D\} \right)  \\
& + n \chi_{[0,\xi]}(n/D) \left(\max_{1\leq i \leq m} e_i\right)J_\xi(n/D), \end{aligned}
\end{equation}
where $\xi \in [0,m]$ is our cutoff parameter from the definition of $\tau^{\sharp}$ in our estimates, and we set for $s\geq 1$
\[J_\xi (s):= \left(
\frac{1}{s} \sum_{j=1}^{\lfloor (s - 1)/\max(1,\xi ) \rfloor} 1/j   \right) +\left(\frac{1}{\lfloor (s + (\xi -1)^+ )/\max(1,\xi ) \rfloor} - \frac{\xi}{s} \right)^+;\]
and for $s<1$,
\[J_\xi (s):= \left(1- \frac{\xi}{s}\right)^+.\]

By our proof, we only need to show that 
\[\sum_{n\in \mathcal{V}_D} \overline{h}(\psi_D^{(n)}) \leq \sum_{n=0}^{mD-1} \overline{h}(\psi_D^{(n)}) + o(D^2).\]
 To this end, it is sufficient to show that for $n\geq mD$, we have 
\begin{equation}\label{withoutShid}
\overline{h}(\psi_D^{(n)}) \leq \min_{0 \leq n' < mD} \overline{h}(\psi_D^{(n')}) + o(D).
\end{equation}

We observe $J_\xi(s)$ is a continuous function in $s$ which is piecewise smooth on the intervals of the form $(k\xi +1, (k+1)\xi), ((k+1)\xi, (k+1)\xi +1)$, where $k \in \NwithzeroB$. Moreover
\[J'_\xi(s)= -\frac{1}{s^2} \sum_{j=1}^k 1/j \leq 0\]
 on $s \in (k\xi +1, (k+1)\xi)$, and 
 \[J'_\xi(s)= -\frac{1}{s^2} \left(- \xi + \sum_{j=1}^k 1/j \right) \leq \xi/s^2\]
  on $s\in ((k+1)\xi, (k+1)\xi +1)$. 
  \silentcomment{ Since  \[J'_\xi(s) = \frac{\xi}{s^2}\] for  $\xi<1$ and $s\in (\xi, 1)$,
    this assertion even holds for $k=0$ and $\xi <1$.}
  Moreover, $J_\xi(s)=0$ for $s\in (0, \min\{1,\xi\}]$. Hence $J_\xi(s)+\xi/s$ is a decreasing function of~$s \in \R_{>0}$, and in particular, its $[\xi,\infty)$ maximum is taken at~$s = \xi$, with value~$1$. 
We analyze the function 
\[F(s) := -s \left(\big(\log|\varphi'(0)| -\sigma_m - (\max_{i=1}^m e_i)J_\xi(s)\big)\right).\]
 It is continuous on~$s \in \R_{> 0}$ and piecewise smooth on the intervals of the form $(k\xi +1, (k+1)\xi), ((k+1)\xi, (k+1)\xi +1)$. For $s\geq \xi$, in each of those intervals,
\begin{equation*}
\begin{aligned}
F'(s) & = - \left(\log|\varphi'(0)| -\sigma_m - \left(\max_{1\leq i \leq m} e_i \right)J_\xi(s)\right) + \left(\max_{1 \leq i \leq m} (e_i)\right)sJ'(s) \\
& \leq  -  \left(\log|\varphi'(0)| -\sigma_m - \left(\max_{1\leq i \leq m} e_i \right)J_\xi(s)\right)  + \left(\max_{1\leq i \leq m} e_i\right)\xi /s\\
& \leq  - \log|\varphi'(0)|  + \sigma_m + \left(\max_{1\leq i \leq m} e_i \right) <0.
\end{aligned}
\end{equation*}

Note that since $u_h < m$ and $\xi \leq m$, we have for $n \geq mD$: 
\[\overline{h}(\psi_D^{(n)}) = - n \log | \varphi'(0)| +  D \left(\overline{\OL(1)} \cdot \overline{\OL(1)}\right) + \sigma_m n+ n \left(\max_{1\leq i \leq m} e_i\right)J_\xi(n/D).\]

Our discussion then shows that $\overline{h}(\psi_D^{(n)})$ is a decreasing function in $n$ on the requisite range $n \geq mD$. 
In fact, the same argument applies to see that  $\overline{h}(\psi_D^{(n)})$ is a decreasing function on the whole $n \in \NwithzeroA$. 
More precisely, for $n/D \geq \xi$ and $n/D \in [u_h, u_{h+1}]$, we have
\begin{equation}
\begin{aligned}
\overline{h}(\psi_D^{(n)}) = & - n \log | \varphi'(0)| +  D \left(\overline{\OL(1)} \cdot \overline{\OL(1)} + \sum_{k=h+1}^r u_kb_k\right) + \left(\sum_{k=1}^h b_k\right) n\\
&+ n \left(\max_{1\leq i \leq m} e_i\right)J_\xi(n/D).
\end{aligned}
\end{equation}
Since $\sum_{k=1}^h b_k\leq \sigma_m$, the same argument as above shows that 
\[F_h(s):=-s \left(\big(\log|\varphi'(0)| - \sum_{k=1}^h b_k - (\max_{i=1}^m e_i)J_\xi(s)\big)\right)\] 
has negative derivative and hence $\overline{h}(\psi_D^{(n)})$ is a decreasing function in $n$.
For $n/D < \xi$ and $n/D \in [u_h, u_{h+1}]$, we have 
\[\overline{h}(\psi_D^{(n)}) = - n \log | \varphi'(0)| +  D \left(\overline{\OL(1)} \cdot \overline{\OL(1)} + \sum_{k=h+1}^r u_k b_k\right) + \left(\sum_{k=1}^h b_k\right) n.\]
Since $\sum_{k=1}^h b_k\leq \sigma_m < \log | \varphi'(0)| $, we conclude that $\overline{h}(\psi_D^{(n)})$ is a decreasing function. 

Thus we obtain \eqref{withoutShid}, giving a proof of Theorem~\ref{main:BC form} free of appeal to the Shidlovsky type theorems from~\S~\ref{functional bad approximability}, but under the stronger assumption that $\log |\varphi'(0)| > \sigma_m + \max_{i=1}^m(e_i)$. In particular, in a manner
free of any of the references in~\S~\ref{functional transcendence}, these remarks already suffice for proving Theorem~\ref{basic main} except for the clause that $|\varphi'(0)| > e^{\max(\sigma_m, \tau(\bb))}$ can be relaxed to $|\varphi'(0)| > e^{\sigma_m}$ when the $f_i$ are \emph{a priori} supposed holonomic.

\subsubsection{Discussion for the $\be = \mathbf{0}$ case of Theorem~\ref{main:BC conv discrete}}
As the behavior of the function $J_\xi$ is the main obstacle to devising a clean proof of our general holonomy bounds not relying on functional bad approximability theorems for holonomic functions, and since we do not logically need an alternative proof for any of our applications, we are content  
(still in the context of explaining how to bypass~\S~\ref{functional bad approximability})
for Theorem~\ref{main:BC conv discrete} with demonstrating how to handle the case $\be=\mathbf{0}$ and under the seemingly mild extra condition that
\begin{equation}\label{condition on m}
m \geq \max_{0 \leq k \leq l}  \left\{  \frac{ \left( \ovcL \cdot \ovcL \right) -  \left( \ovcL_{r_k} \cdot \ovcL  \right) + \alpha_k \left(\log |\varphi'_{r_k}(0)| - \sum_{j=1}^{h(k)} b_j\right) - \sum_{j=h(k)+1}^r u_jb_j }{\log |\varphi'(0)| - \sigma_m} \right\},
\end{equation}
where, for a given $k$, we pick the $ h(k)$ with $\alpha_k \in [u_{h(k)}, u_{h(k)+1})$, and we recall the convention $\alpha_0 = 0$.
In other words, we will show without appealing to~\S~\ref{functional bad approximability} that when $\be=\mathbf{0}$ and under the conditions of Theorem~\ref{main:BC conv discrete},  at least one of the dimension bounds \eqref{disc bounds} or 
\begin{equation} \label{contrapositive of condition on m} m \leq \max_{0 \leq k \leq l}  \left\{  \frac{ \left( \ovcL \cdot \ovcL \right) -  \left( \ovcL_{r_k} \cdot \ovcL  \right) + \alpha_k \left(\log |\varphi'_{r_k}(0)| - \sum_{j=1}^{h(k)} b_j\right) - \sum_{j=h(k)+1}^r u_jb_j }{\log |\varphi'(0)| - \sigma_m} \right\}
\end{equation}
is in place.  
 In practice, we have always found that the inequality~\eqref{contrapositive of condition on m} is already implied by the contrapositive of the condition inequality~\eqref{disc bounds}, in which case
 the conclusion~\eqref{disc bounds} certainly follows.

For $n/D \in [\alpha_k, \alpha_{k+1}] \cap [u_{h}, u_{h+1}]$ (here $h$ does not need to be $h(k)$ defined above), the proof of Theorem~\ref{main:BC conv discrete} shows that
\begin{equation}\label{ovh}
\overline{h}(\psi_D^{(n)})= -n \log|\varphi'(0)| - n \log r_k + D(\ovcL \cdot \ovcL_r) + n \sum_{j=1}^{h} b_j + D\sum_{j=h+1}^r u_j b_j ,
\end{equation}
which, if viewed as a piecewise linear function in $s=n/D$, is continuous on $s\in [0, \infty)$. The local minima 
can only occur at points of the form $s=\alpha_k$, or over a line segment of slope~$0$. 
But for $n/D \geq m$ we use as archimedean height evaluation bound
\[\overline{h}(\psi_D^{(n)})= -n \log|\varphi'(0)| + D(\ovcL \cdot \ovcL) + n \sigma_m,\]
which decreases monotonically in $s = n/D$. Now to show~\eqref{withoutShid}, which as in~\eqref{noShidBC} suffices to bypass the appeal to~\S~\ref{functional bad approximability} in our analysis in~\S~\ref{sec_BCconvexity}, it is enough to check that
\[\overline{h}(\psi_D^{(mD)}) \leq \min_{0\leq k \leq l} \overline{h}(\psi_D^{(\alpha_k D)});\]
here as $\alpha_k D$ might not be an integer, by a slight abuse of notation, $\overline{h}(\psi_D^{(\alpha_k D)})$ means replacing $n$ in~\eqref{ovh} by $\alpha_k D$.
This unfolds to the definition of the condition~\eqref{condition on m}.

\subsubsection{Discussion for Theorem~\ref{main: easy convexity}}
Under the conditions
\begin{equation} \label{cfc}
\begin{aligned}
\xi > u_1 \geq 1,  \quad b_1 & > \log R \geq 0,  \quad  \log |\varphi'(0)| > \sigma_m + \max_{1 \leq i \leq m} e_i, \\
 \Gamma(\log R, r, \lambda, \mu) & \geq u_1 \log R,
\end{aligned}
\end{equation}
we show independently of~\S~\ref{functional bad approximability} that at least one of the dimension bounds~\eqref{EasyConvexityBound} or 
\begin{equation}\label{m condition in easyconvexity}
m \leq \frac{T(\varphi; r,\lambda, \mu) - \sum_{j=1}^r u_jb_j }{\log|\varphi'(0)| -\sigma_m - (\max_{1\leq i \leq m} e_i)J_\xi(m)}
\end{equation}
is in place. 
In the practical situations of many applications, usually~\eqref{m condition in easyconvexity} is expected to be a smaller bound than~\eqref{EasyConvexityBound}. This applies for example to our proof of Theorem~\ref{mainA} via Example~\ref{Ex-easyconv}, where the conditions~\eqref{cfc} are met, and the right-hand side of~\eqref{m condition in easyconvexity} is negative, whereas the holonomy bound of~\eqref{EasyConvexityBound} is at $\sim 13.8527$. 
\silentcomment{For instance, in our proof for Theorem~\ref{mainA}, the assumption is satisfied since we take $\xi=2, R=4, \sigma_m=4,  (\max_{1\leq i \leq m} e_i)=1$ and 
\[\log |\varphi'(0)|=\log \left( 256  \cdot
\frac{5448339453535586608000000000}{8658833407565631122430056127}
\right).\]
 Recall from Example~\ref{Ex-easyconv}, we have $m=14, T(4.7,10, -4.5)=6.5316\ldots$. We compute 
\[J_2(14)=\left(
\frac{1}{14} \sum_{j=1}^{\lfloor (14 - 1)/\max(1,2 ) \rfloor} 1/j   \right) +\left(\frac{1}{\lfloor (14 + (2 -1)^+ )/\max(1,2 ) \rfloor} - \frac{2}{14} \right)^+=\frac{7}{40}.\]
We then have the bound in \eqref{m condition in easyconvexity} would give
\[m \leq \frac{6.5316\ldots -1*2 - 3*2}{\log \left( 256  \cdot
\frac{5448339453535586608000000000}{8658833407565631122430056127}
\right) - 4 - \frac{7}{40}} <0\]}

Let $\overline{h}_{\infty}(\psi_D^{(n)})$ and $\overline{h}_{\mathrm{fin}}(\psi_D^{(n)})$ denote the respective main terms of our bounds on $h(\psi_D^{(n)})$, 
given in~\eqref{eq:arch_id2}, \eqref{eq:arch_id}, and \eqref{arch_phi} (for the archimedean estimates), and~\eqref{finite<xi} and \eqref{finite>xi} (for the finite estimates), according to the various cases in dependence on~$n/D$.  
In these notations we have $\overline{h}(\psi_D^{(n)})$ continuous in $s :=n/D$, and hence again we only need to ensure \eqref{withoutShid}. 

In the case at hand, our assumptions imply $m \geq \max\{ \chi_1, \xi, u_1, \ldots, u_r\}$. For $n \geq mD$, we derive for the left-hand side of~\eqref{withoutShid}: 
\[\overline{h}(\psi_D^{(n)}) = -n \left(\log|\varphi'(0)| -\sigma_m - \left(\max_{1\leq i \leq m} e_i\right)J_\xi(n/D)\right) + DT(r,\lambda, \mu).\]
By the same analysis as in \S~\ref{noShidBC}, we derive that $\overline{h}(\psi_D^{(n)}) $ is a decreasing function of~$n$ in the range $n\geq \max\{\xi, \chi_1\} D$; therefore for $n \geq mD$, we have
\[\overline{h}(\psi_D^{(n)}) \leq \min_{\max\{\xi, \chi_1\} D \leq n' < mD} \overline{h}(\psi_D^{(n)}).\]

It remains to consider the range $n' < \max\{\xi, \chi_1\} D$.
If $\xi > \chi_1$, then for $s':=n'/D\in [\chi_1, \xi] \cap [u_h, u_{h+1}]$, we have
\[\overline{h}(\psi_D^{(n)}) = -n' \left(\log|\varphi'(0)| -\sum_{j=1}^h b_j \right) + D\left(T(r,\lambda, \mu) +\sum_{j=h+1}^r u_j b_j \right),\]
and therefore $\overline{h}(\psi_D^{(n)})$ is a decreasing function in $s'$ in the range $[\chi_1, \xi]$; 
continuity the gives 
\[\overline{h}(\psi_D^{(n)}) \leq \min_{\chi_1 D \leq n' < mD} \overline{h}(\psi_D^{(n)}).\]

Next we take up the range $\xi \leq s' <\chi_1$. Here we have for $s' \in [u_h, u_{h+1}]$,
\begin{equation*}
\begin{aligned}
\overline{h}(\psi_D^{(n')}) = & n' \left(-\log R +\sum_{j=1}^h b_j + \left(\max_{1\leq i \leq m} e_i\right)J_\xi(n/D)\right) \\
&+ D\left(\Gamma(\log R, r,\lambda, \mu)+ \sum_{j=h+1}^r u_jb_j \right).
\end{aligned}
\end{equation*}
Since $u_1 < \xi$, we have $h\geq 1$ for these $n'$ and then $\overline{h}(\psi_D^{(n')}) $ is an increasing function of $s'$ due to our assumption $\log R < b_1$ (while, by definition, $J_\xi \geq 0$). 

We continue: for $u_1\leq  s' <\min\{\xi, \chi_1\}$ (if $u_1> \chi_1$, then this set is empty and move on to the next case below), we have (the $h$ in the formula may vary depending on $s'$ as above)
\[\overline{h}(\psi_D^{(n')}) = n' (-\log R +\sum_{j=1}^h b_j) + D\left(\Gamma(\log R, r,\lambda, \mu)+\sum_{j=h+1}^r u_jb_j \right),\]
also an increasing function of $s'$.
 
 For $\chi_0 \leq s' < \min\{u_1, \chi_1\}$, 
 \[\overline{h}(\psi_D^{(n')}) = - n' \cdot \log R + D\left(\Gamma(\log R, r,\lambda, \mu)+\sum_{j=1}^r u_jb_j\right),\]
 is a decreasing function of $s'$ due to our assumption $\log R \geq 0$.

Finally, for $0 \leq s' <\chi_0$, we have 
\[\overline{h}(\psi_D^{(n')}) = D\left(-\lambda (s')^r - \mu s' +\sum_{j=1}^r u_jb_j \right).\] 
Recall that $\lambda>0$. If $\mu >0$, then $\overline{h}(\psi_D^{(n')})$ is a decreasing function. If $\mu \leq 0$, 
the critical point $s_0$ satisfies that
\[s_0=\left(\frac{ - \mu}{r\lambda }\right)^{1/(r-1)}\leq \left(\frac{\log R - \mu}{r\lambda }\right)^{1/(r-1)} = \chi_0,\]
under our assumed conditions;
therefore $\overline{h}(\psi_D^{(n')})$ is an increasing function on $[0, s_0]$ and a decreasing function on $[s_0, \chi_0]$.

The above discussion shows that
\[\min_{0 \leq n' \leq mD} \overline{h}(\psi_D^{(n)})= \min\{\overline{h}(\psi_D^{(0)}), \overline{h}(\psi_D^{(u_1D)}), \overline{h}(\psi_D^{(mD)})\}.\]
Since
\[\overline{h}(\psi_D^{(0)})=D \sum_{j=1}^r u_jb_j, \,\,  \overline{h}(\psi_D^{(u_1D)})= -u_1D \log R + D\left(\Gamma(\log R, r,\lambda, \mu)+\sum_{j=1}^r u_jb_j\right),\]
we have $\overline{h}(\psi_D^{(0)}) \leq \overline{h}(\psi_D^{(u_1D)})$ due to the assumption $\Gamma(\log R, r, \lambda, \mu) \geq u_1 \log R$.

In upshot, the requisite inequality 
\[\overline{h}(\psi_D^{(mD)}) \leq \min_{0 \leq n' < m D} \overline{h}(\psi_D^{(n')}) \]
boils down to securing that $\overline{h}(\psi_D^{(mD)}) \leq D \sum_{j=1}^r u_jb_j$, which is equivalent to 
\[m \geq \frac{T(r,\lambda, \mu) - \sum_{j=1}^r u_jb_j}{\log|\varphi'(0)| -\sigma_m - (\max_{1\leq i \leq m} e_i)J_\xi(m)}.\]
In other words, we have either proved  \eqref{withoutShid} and hence \eqref{EasyConvexityBound} holds, or else \eqref{m condition in easyconvexity} holds.

\section{The finer holonomy bound with the Bost--Charles integral} \label{slopes} 
This section combines the measure concentration input of~\S~\ref{fine section} with the Bost--Charles refinements of~\S~\ref{new slopes} to give the most accurate general holonomy bound of all the theorems worked out in our paper. The added strengthening turns out to be zero for the particular applications to Theorems~\ref{mainA} and~\ref{logsmain}, and the theorem becomes somewhat complicated to state, nevertheless we hope that the principle of the abstract refinement could be useful in future applications of our holonomy bounds. 

The theoretical improvement from adding high-dimensional methods to the Bost--Charles calculus is, as far as we were able to tell in the framework of~\S~\ref{seriousintro}, reflected only in the denominator term~$\tau(\mathbf{b;e})$. It is inevitable to ask which of the growth integrals in~\S~\ref{fine section} versus~\S~\ref{new slopes} is the smaller. 
In~\S~\ref{rmk_Nazarov}, we present a proof by Fedor Nazarov that the Bost--Charles integral is strictly better than the rearrangement integral.
This, in particular, implies that Theorem~\ref{main:BC fullconv} is more precise than Theorem~\ref{main:elementary form}, granting our heuristic Remark~\ref{stationary choice}, and at least as far as the stated denominator is concerned in the latter theorem. (Remark~\ref{finedenominators} indicates that the multidimensional proof in~\S~\ref{fine section} can go further than~\S~\ref{new slopes} in the general denominator aspect, and at least as far as our choice of treatment in the present section~\S~\ref{slopes}.) Remark~\ref{nearly equal} further down in this section suggests that the difference in the archimedean growth terms is usually very small, implying that little is to be lost from working with the rearrangement integrals.  For purposes of sampling and testing the holonomy bounds with different maps~$\varphi$, the latter integrals have the practical advantage to allow for faster and more reliable numerical computations.

    We now proceed to formulating our unifying theorem. 
The following sums up the sharpest\footnote{Possibly up to considering other $\varphi$ in addition to $\varphi_r$, as in Theorems~\ref{main:elementary form} or~\ref{main: easy convexity}. Such variations are straightforward to incorporate as well, but we refrain from doing this here.} of all the holonomy bounds we prove in this paper. 

\begin{thm}\label{high dim BC convexity}
We relax all assumptions on the denominators in Theorem~\ref{main:elementary form}, and allow for an
arbitrary denominators matrix $\mathbf{b} \in M_{m \times r} \left( \R_{\geq 0} \right)$ with nonnegative coefficients. 
Define
\begin{equation} \label{tau bar}
\begin{aligned}
\ovtau(\bb)  & :=\limsup_{\epsilon\rightarrow 0, \varepsilon \rightarrow 0}\limsup_{d\rightarrow \infty}  \left\{ \frac{1}{d}  \limsup_{N\rightarrow \infty} 
\QQQQ_N(\bb,\epsilon,d,\varepsilon) \right\}, \\
\QQQQ_N(\bb,\epsilon,d,\varepsilon) & :=
 \frac{2}{N} \sum_{k \in \NwithoutzeroA}  \max_{\bn \in P_{\varepsilon}^d(N) , \bi \in V_m^d(\epsilon)} \#\{(j,h) 
\, : \, k \leq b_{i_j,h} \cdot n_j \}, \end{aligned} 
\end{equation}
where $(j,h)$ runs through $\{1,\ldots,d\} \times \{1,\ldots,r\}$, 
the set $V_m^d(\epsilon) \subset \{1,\ldots, m\}^d$ is defined as
\[
\begin{aligned}
V_m^d(\epsilon)  := \big\{ &  \bi \in \{ 1,\ldots, m \}^d \, : \, \forall i_0 \in \{1,\ldots,m\},  \\
&  d/m - \epsilon d <  \#\{1\leq j\leq d \mid i_j=i_0\} < d/m +\epsilon d \big\},
\end{aligned}
\]
and
$ P_{\varepsilon}^d(N)   \subset [0, N]^d \cap \Z^d$ denotes the subset  of those $\bn$ for which the normalized $\left( [0,1), \mu_{\mathrm{Lebesgue}} \right)$ discrepancy of $\{n_i/N\}_{i=1}^{d}$ is $\leq \varepsilon$.

Assume either that 
\(\displaystyle{\log |\varphi'(0)| > \max \left\{\sum_{h=1}^r \max_{1\leq i \leq m} b_{i,h} ,\ovtau(\bb)+ \tau^\sharp(\be)\right\}}\), or that $\log{\varphi'(0)} > ,\ovtau(\bb)+ \tau^\sharp(\be)$ and all~$f_i$ are holonomic. 

 Then we have 
\begin{equation}\label{BCbound-highdim} 
m  \leq  \frac{ \iint_{\T^2} \log|\varphi(z)-\varphi(w)| \, \mv(z) \mv(w) }{  \log{|\varphi'(0)|} - (\ovtau(\bb)+ \tau^\sharp(\be)) }. 
\end{equation}

Further, for any subradii sequences $1 = r_l > r_{l-1} > \cdots > r_0>0$ as in Theorem~\ref{main:BC conv discrete} and
using the Bost--Charles characteristic~$\hT$ of Definition~\ref{BC characteristic}, we have 
\begin{equation}\label{BCbound-highdim-cov} 
m  \leq  \frac{ \iint_{\T^2} \log|\varphi(z)-\varphi(w)| \, \mv(z) \mv(w)  - \frac{1}{m} \sum_{k=1}^l \frac{ (\hT(r_k,\varphi)- \hT(r_{k-1},\varphi))^2 }{\log{r_k} - \log{r_{k-1}}}}{  \log{|\varphi'(0)|} -(\ovtau(\bb)+ \tau^\sharp(\be)) }.
\end{equation}

Moreover, if the $s_h^*$ for the given $\{r_k\}$ are defined as in Theorem~\ref{main:BC fullconv} 
then we have
\begin{equation}\label{BCbound-fullconv}
m \leq  \frac{  \iint_{\T^2} \log|\varphi(z)-\varphi(w)| \, \mv(z) \mv(w) - \sum_{h=0}^{l-1} s_h^* \cdot (\hT(1,\varphi) - \hT(r_h, \varphi))}{\log |\varphi'(0)|-(\ovtau(\bb)+ \tau^\sharp(\be)) }.
\end{equation}
\end{thm}

\begin{remark} \label{equalityoftaus}
We observe that $\ovtau(\bb) \in [0,\infty)$ by definition. More precisely, we have the trivial bound $\ovtau(\bb) \leq \sum_{h=1}^r \max_{i=1}^m  \{ b_{i,h} \}$.

For the denominator matrices $\bb$  of the form considered throughout~\S\S~\ref{fine section}--\ref{new slopes}, we have $\ovtau(\bb) = \tau^\flat(\bb)$. Indeed, writing~$t$ for the continuous limit of the discrete variable~$k/N$, 
we have in the setup of Theorem~\ref{main:elementary form}~$ \ovtau(\bb)$ bounded above by
\begin{equation*}
\begin{aligned}
&    \sum_{h=1}^r \limsup_{\epsilon, \varepsilon \rightarrow 0} \limsup_{d\rightarrow \infty} 
\left\{ \frac{1}{d}  \limsup_{N\rightarrow \infty} \left\{ \frac{2}{N} \sum_{
k \in \NwithoutzeroA}  \max_{\bn \in P_{\varepsilon}^d(N) , \bi \in V_m^d(\epsilon)} \#\{(j,h) \, : \,  k \leq b_{i_j, h} \cdot n_j \} \right\} \right\} \\
& \leq  2 \sum_{h=1}^r \int_0^{b_h} \min\{ 1 - u_h/m, 1- t/b_h\}   \, dt = \sum_{h=1}^r \left(b_h - \frac{b_h u_h^2}{m^2} \right) = \tau^\flat(\bb),
\end{aligned}
\end{equation*}
where the inequality stems from the observation that the restriction to the \emph{balanced}~$\bi$ (meaning: each $i_0 \in \{1,\ldots,m\}$ occurs with the same asymptotic frequency~$1/m$) supplies the upper bound constraint $1 - u_h/m$ in the integrand on the second line, while the restriction to the \emph{balanced} $\bn$ (meaning: the components set~$\{n_j\}$ takes asymptotically the uniform distribution on~$[0,N]$)  supplies the upper bound constraint $1- t/b_h$ in that integrand. Moreover, equalities are reached in the case where both~$\bn$ and~$\bi$ are arranged in non-decreasing order: $n_1 \leq \cdots 
\leq n_d$ and $i_1 \leq \cdots \leq i_d$.  This
proves that $\ovtau(\bb)=\tau^\flat(\bb)$ under the standing assumptions throughout~\S~\ref{fine section} and~\S~\ref{new slopes}. (See also Lemma~\ref{single step valuationwise}
and Remark~\ref{genbcomplicated}.)~\endofremark
\end{remark}

\begin{example} \label{nonequalityoftaus}
Here is a simple example to illustrate that, for a given~$\bb$,
if one lets~$\bb'$ range over all arrays which dominate~$\bb' \ge \bb$ coefficient-wise
and which additionally meet the constraints of
Theorem~\ref{main:elementary form}, the inequality in $\ovtau(\bb) \leq \tau^{\flat}(\bb')$ can be strict. 

Consider $\bb=[0,1,2]^{\mathrm{t}}$. In order to use $\tau^\flat$, the optimal choice of~$\bb'$ is to take $\bb'=[0,2,2]^{\mathrm{t}}$ and we have $\tau^\flat(\bb')=\frac{16}{9}$. On the other hand, we have $\ovtau(\bb)=\frac{5}{3} = \frac{15}{9}$. Indeed, writing $t :=k/N$, we have: 
\[\ovtau(\bb)= 2 \left(\int_0^{2/3} \frac{2}{3}\, dt + \int_{2/3}^1  \left( 1-t +\frac{1}{3} \right) \, dt + \int_1^{4/3} \frac{1}{3} \, dt + \int_{4/3}^2 \left( 1-\frac{t}{2} \right) \, dt   \right).\]
This is to be compared with
\[\tau^\flat(\bb') = 2 \int_0^2 \min\left\{\frac{2}{3}, 1- \frac{t}{2} \right\} \, dt.\]
To explain the difference between two formulas, we notice that in the range $t\in [2/3, 4/3]$, for every $\bn$ to be considered in the definition of $\ovtau(\bb)$, there are at most $(\max\{0, 1-t\} + o(1)) d$ among the $n_j$ with $i_j=2$ (corresponding to $b_{2,1}=1$) to contribute to $\{(j,1) \mid k \leq n\}$, and there are at most  $(1/3+ o(1)) d$ among the $n_j$ with $i_j=3$ (corresponding to $b_{3,1}=2$) to contribute to $\{(j,1) \mid k \leq 2 n\}$ and both bounds can be reached with suitable choice of $\bn$.~\endofremark
\end{example}

\subsection{Comparison of the Bost--Charles and the Rearrangement integrals} \label{rmk_Nazarov}
This section, due entirely to Fedor Nazarov,  treats the clean comparison of the two integrals  ---  beneath the empirical observation that they are practically the same in the situations
we encounter in~\S~\ref{contour choiceA}, as well as in practice for most of the multivalent cases. This theorem strictly speaking is not used for any of our proofs in the paper, 
and hence it can be omitted on a first (and on a second) reading.

\begin{basicremark} \label{Bochner remark}
As the considerations that follow rely on the potential theory in the plane, consider first the easier situation on the simplest of all the Lie groups: the circle $\T$. 
The integrable function $G(z) := \log{\frac{1}{|1-z|}}$, $G : \T \to \R \cup \{\infty\}$ has the nonnegative Fourier coefficients
\begin{equation}
\widehat{G}(n) := \int_{\T} G(z) z^{-n} \, \mv(z) =  \int_{\T} z^{-n} \log{\frac{1}{|1-z|}} \, \mv(z)
= \left\{  \begin{array}{cc}  0, &  \textrm{ if } n = 0;  \\
\displaystyle{\frac{1}{2|n|}}, &  \textrm{ if } n \neq 0. \end{array} \right.
\end{equation}
By the general Bochner theorem, the positivity of the Fourier transform implies that $G(z)$ is 
a positive-definite function on the locally compact abelian group $\T$: that is, 
$\iint_{\T^2} G(zw^{-1}) \, \nu(z)\nu(w) \geq 0$ for all \emph{reasonable} signed measures $\nu$ on $\T$. \emph{Reasonable}
 here may be taken to mean $\nu = \nu^+ - \nu^-$ with finite positive measures $\nu^{\pm}$ satisfying $\int_{\T} G \, \nu^{\pm} < \infty$.
 For any such signed measure~$\nu$, this computation shows more precisely that
 \begin{equation} \label{circle energy}
 \begin{aligned}
 I(\nu)  & := \iint_{\T^2} \log{\frac{1}{|z-w|}} \, \nu(z)\nu(w) \\
 &  = \iint_{\T^2} \log{\frac{1}{|1-zw^{-1}|}} \, \nu(z)\nu(w) 
 = \sum_{n \in \Z \setminus \{0\}} \frac{|\widehat{\nu}(n)|^2}{2|n|}  \geq 0,
 \end{aligned}
 \end{equation}
 where manifestly the equality holds if and only if the Fourier transform $\widehat{\nu}$ is a scalar multiple of the Dirac mass 
 at $0 \in \widehat{\T} = \Z$, and that in turn is the case if and only if $\nu$ is a scalar
 multiple of the Haar measure $\mv$. 
This reflects the basic potential theory on the circle. 
 
 On the circle $|z| = R$ of arbitrary radius, the left-hand side of~\eqref{circle energy} scales by the additive summand $(\log R)  \left( \int \nu \right)^2$, and so whereas for $R > 1$ the inequality~\eqref{circle energy} is false for arbitrary measures $\nu$ supported by that circle, it continues to be in place for the measures that are \emph{balanced} in the
 sense that $\int \nu = \nu(\C) = 0$. As the logarithmic kernel is also invariant under additive translations, the latter remark continues to hold for balanced
 measures carried by any circle in $\C$. \endofremark
\end{basicremark}

As we review next,  the positivity of the energy integral is a completely general fact about balanced measures on~$\C$. 

\subsubsection{The energy principle} In the Newtonian  gravitational field created by a point mass at the origin $\mathbf{0} \in \R^n$, the potential energy function $U : \R^n  \setminus \{ \mathbf{0} \} \to \R$ is determined by the  distributional Laplace equation
$$
\Delta U := \sum_{i=1}^n \frac{\partial^2 U}{ \partial x_i^2} =- \delta_{\mathbf{0}},
$$
namely as the fundamental solution $U := U_n^{2}$ of that equation, where more generally it is useful to consider the \emph{Riesz potential}~\cite{Riesz}  
defined by
\begin{equation}  \label{Riesz potential}
U_n^{\alpha}(\mathbf{x}) := \left\{ \begin{array}{ll} 
\displaystyle{ \frac{\Gamma\left( \frac{n}{2}\right)}{2\pi^{\frac{n}{2}}} \log{\frac{1}{\|\mathbf{x}\|}}}, & \, {\textrm{ for }} n = \alpha; \\
\displaystyle{ \frac{\Gamma \left( \frac{n-\alpha}{2} \right)}{2^{\alpha} \Gamma\left( \frac{\alpha}{2} \right)\pi^{\frac{n}{2}}} \frac{1}{\| \mathbf{x} \|^{n-\alpha}}}, & \, \textrm{ for } n > \alpha > 0. \end{array} \right.
 \end{equation}
 Here, $\| \mathbf{x} \| := \sqrt{ \sum_{i=1}^{n} x_i^2 }$ is the Euclidean distance function in $\R^n$. The gravitational potential $U = U_n^2$ is rotationally invariant
 and defines a  kernel function $k(\mathbf{x,y}) := U_n^2(\mathbf{x-y})$, which by definition is furthermore translationally invariant. 
 For $n > 2$ this kernel is also positive-definite, by the fundamental formula of Frostman and  Marcel Riesz~\cite[\S~I.3]{Riesz}, which generalizes
 the Dirichlet energy integral $\frac{1}{2}\iint \| \nabla F \|^2 \, \mathrm{d vol}$, and is tantamount
 to the computation~\cite{Deny,Schwartz,Landkof,NikishinSorokin} of the distributional Fourier transform of~$U_n^2$: 
 \begin{equation}
 \begin{aligned}
 I(\nu) &  := \iint_{\R^n \times \R^n} k(\mathbf{x,y})  \, \nu(\mathbf{x})  \nu(\mathbf{y})  = 
\iint_{\R^n \times \R^n} U_n^2(\mathbf{x-y})  \, \nu(\mathbf{x})  \nu(\mathbf{y})  \\  \label{energy formula}
& =  \int_{\R^n} \left(  U_n^1 \ast \nu \right)^2 \, \mu_{\mathrm{Lebesgue}}.
\end{aligned}
\end{equation}
This is the spatial analogue of the energy formula~\eqref{circle energy} for the circle. 
To be more precise, this formula gives the strict 
positivity $I(\nu) > 0$ 
 of  the  energy 
of any nonzero compactly supported signed measure $\nu = \nu^+ - \nu^-$ on $\R^n$ expressible as the difference of two finite positive Borel measures
$\nu^+, \nu^-$ of finite energies $I(\nu^+), I(\nu^-) < \infty$. 
For the logarithmic kernel (the case $n = 2$, whose proper physical interpretation is rather in electrostatics on a plate), Riesz observed~\cite[\S~I.4]{Riesz}
that the analogy becomes almost perfect upon additionally requiring the signed measure $\nu$ to be \emph{balanced}: $\nu(\R^2) = 0$. 
The \emph{energy principle} states that for balanced signed measures subject to the above regularity conditions (with the balancing condition being
only required in the case $n = 2$), the energy $I(\nu) \geq 0$ is nonnegative, and equality holds if and only if $\nu = 0$. This refines the uniqueness theorem for 
the equilibrium probability measure of a compact. 
The $n=2$ case, which is the one of relevance to us, is treated in detail in~\cite[Theorem~16.4.2]{Hille}. 

\subsubsection{Fuglede's inequality} For completeness, we summarize the more general situation due to Fuglede~\cite{Fuglede} for non-balanced measures. This is 
a different generalization of the $n=2$ case, this time to the logarithmic kernel 
$$
k(\mathbf{x,y}) := U_n^n(\mathbf{x-y}) =  \frac{\Gamma(n/2)}{2 \pi^{n/2}} \,  \log\frac{1}{\| \mathbf{x-y} \| }
$$
 on $\R^n$. For this kernel, the analog of the energy formula~\eqref{energy formula} for the case of the balanced measures is also due to Riesz~\cite[\S~I.4]{Riesz}:
 $$
\int_{\R^n}\nu = 0 \quad \Longrightarrow \quad  I(\nu) =  \int_{\R^n} \left(  U_n^{n/2} \ast \nu \right)^2 \, \mu_{\mathrm{Lebesgue}} \geq 0, 
 $$
 and more generally,
 Fuglede~\cite[\S~4]{Fuglede} proves the sharp energy lower bound
 \begin{equation} \label{Fugledes}
 I(\nu) \geq  \log{ \left( \frac{a_n}{R} \right) }  \cdot \left( \int_{\R^n} \nu \right)^2
 \end{equation} 
 for all signed measures $\nu$ on $\R^n$ expressible as $\nu = \nu^+ - \nu^-$ with $\nu^{\pm}$ finite positive measures of convergent energy
 integrals $I(\nu^{\pm}) < \infty$ and having $\mathrm{supp}(\nu) \subseteq \{ \| \mathbf{x} \| \leq R \}$, and with the optimal constant $a_n$ being precisely
 $$
 a_n := \left\{  \begin{array}{ll}
 \displaystyle{ \exp\left( \frac{1}{2}  +\ldots + \frac{1}{n-4} + \frac{1}{n-2} \right)}, &  \textrm{ for $n$ even};  \\  
\displaystyle{ \exp\left( \frac{1}{1}  +\ldots + \frac{1}{n-4} + \frac{1}{n-2} - \log{2} \right)}, &  \textrm{ for $n$ odd}. \end{array}  \right.
 $$
 The cases $a_1 = 1/2$ and $a_2 = 0$ (for the case $n=2$ of relevance to us) are already in de la Vallee-Poussin~\cite[\S~47]{dVP}. In any case, \eqref{Fugledes} 
 certainly implies the requisite positivity $I(\nu) \geq 0$ for all (reasonable) \emph{balanced} measures, and more generally, for all measures supported by a sufficiently small ball.

 \subsubsection{Michelli's criterion for positive-definite kernels} \label{Michelli criterion} A different generalization, for which we refer to~\cite{Mattner} and the references there, admits an arbitrary kernel of the form $k(\bx,\mathbf{y}) = U\left(\| \mathbf{x}-\mathbf{y} \| \right)$, where $U(t) \in C^{\infty}\left(\R_{>0}\right)$ obeys $(-1)^n\left( \frac{d}{dt} \right)^n U(t) \geq 0$ for all~$n \geq n_0$, and now the additional (``balancing'') constraints on the compactly supported signed measure $\nu = \nu^+ - \nu^-$ with $I(\nu^{\pm}) < \infty$ on~$\R^d$ being $\int_{\R^d} \bx^\bm \, \nu(\bx) = 0$ for all $\bm \in \NwithzeroA^d$ with $|\bm| < n_0$. The condition on~$U(t)$ is equivalent to the existence of an integral representation 
 for $\left( \frac{d}{dt} \right)^{n_0} U(t) = \int_0^{\infty} e^{- tu} \, \mathrm{d}\alpha(t)$ as a Laplace--Stieltjes transform of a \emph{positive} Borel measure $\mathrm{d}\alpha$ on~$\R_{> 0}$. 

\subsubsection{Intersection pairing and signature} Another way of informally summarizing this discussion\footnote{This is taking~$n_0 := 1$ in~\S~\ref{Michelli criterion} if we are to include the more general kernels~$k(x,y)$ there. For arithmetic geometry, the case of relevance is~$n=2$ and the kernel~$k(z,w) = -\log{|z-w|}$.} is to say that
 the infinite-dimensional quadratic form 
 $$
 \langle \mu, \nu \rangle := \iint_{\R^n \times \R^n} k(\mathbf{x,y}) \, \mu(\mathbf{x})\nu(\mathbf{y})
 $$
  has one `$-$' sign  on the space of \emph{reasonable} (non-balanced) signed measures on $\R^n$. For the case $n=2$ of relevance
  to the rest of~\S~\ref{rmk_Nazarov}, it could be interesting to know if a more precise connection could be drawn to the arithmetic Hodge index 
  formula in Arakelov theory and the computations in~\cite[\S~5]{BostCharles} that led to the Bost--Charles double integral. Is there a proof of Nazarov's inequality (Proposition~\ref{nazarovbound} below) that works directly into the arithmetic intersection theory framework of~\cite{BostCharles}? 
  A basic remark in the algebraic model is that for any two line bundles $L,M$  on a polarized normal projective algebraic surface, if
   $(L.L) = (M.M)$  and $\deg{L} = \deg{M}$, then $(L.L) \leq (L.M)$ following from the Hodge index theorem for
   the line bundle $L \otimes M^{-1}$.

\subsubsection{Nazarov's inequality} Consider now the case $n = 2$ as the complex plane $\R^2 \cong \C$. Then the rotation group
is realized by the unitary transformations $t_a \, : \, z \mapsto az$, indexed by the circle points $a \in \T$, 
and our rotationally-invariant, positive-definite kernel is given by $k(z,w) = \frac{1}{2\pi} \log{\frac{1}{|z-w|}}$ on the $\C$-linear space
of  balanced  signed measures on $\C$ with the regularity conditions we described.  We consider a non-constant continuous
function $\varphi : \T \to \C$, and for any $a \in \T$ we apply the energy principle $I(\nu) \geq 0$ to the balanced signed measure 
$\nu = \nu_a := \varphi_*(\mv) - t_a^*\varphi_*(\mv)$. As the kernel $k(z,w)$ is symmetric and rotationally invariant, the resulting inequality rewrites as
\begin{equation} \label{sliced}
\iint_{\T^2} \log{|\varphi(z) - \varphi(w)|} \, \mv(z)\mv(w) \leq \iint_{\T^2} \log{|a \varphi(z) - \varphi(w)|} \, \mv(z)\mv(w),
\end{equation}
for any $a \in \T$, and with equality holding if and only if $\nu_a = 0$. Integrating over $a \in \T$ we get: 

\begin{proposition}[Nazarov] \label{nazarovbound}
For any continuous function $\varphi : \T \to \C$, 
\begin{equation*}
\begin{aligned}
\iint_{\T^2} \log{|\varphi(z) - \varphi(w)|} \, \mv(z)\mv(w) & \leq \iint_{\T^2}  \log \max \left( |\varphi(z)|,
|\varphi(w)| \right)  \, \mv(z)\mv(w) \\
& = \int_0^1 2t \cdot ( \log{|\varphi(e^{2\pi i t})|} )^* \, dt, 
\end{aligned}
\end{equation*}
where $g^*$  denotes the increasing rearrangement~\eqref{incr} of a continuous function $(0,1) \to \R$. 
Furthermore, equality holds if and only if $\varphi(z) = cz^m$ for some $c \in \C$ and $m \in \Z$. 
\end{proposition}

\begin{proof} 
Using the Poisson formula
\begin{equation*}
\int_{\T} \log{| a x - y |} \, \mv(a) = \log \max{(|x|,|y|)}
\end{equation*}
in the termwise integration $\int_{\T} I(\nu_a) \, \mv(a) \geq 0$  of~\eqref{sliced}. The equality requires $\nu_a = 0$ for almost all $a \in \T$, hence that
$\varphi_*(\mv)$ is rotationally-invariant, and hence  that the supporting loop $\varphi(\T) \subset \C$ is rotationally-invariant and therefore 
a centered circle, and that $\varphi^*(\mv)$ is a scalar multiple of the Haar measure of that circle. 
\end{proof}

In particular, we get a clean proof that the Bost--Charles integral is always  strictly majorized by the doubled Nevanlinna characteristic that we
have in~\cite{zeta5} (and not merely by the slightly larger doubled Ahlfors--Shimizu characteristic 
\[\int_{\T} \log{\sqrt{1+|\varphi|^2}} \, \mv\]
 noted in~\cite[Prop.~5.4.5]{BostCharles},  which is a more basic estimate following simply by the trivial pointwise inequality $|x-y|^2 \leq (1+|x|^2)(1+|y|^2)$): 

\begin{cor} \label{crude indeed}
Every continuous function $\varphi : \T \to \C$ satisfies 
\begin{equation*}
\iint_{\T^2} \log{|\varphi(z) - \varphi(w)|} \, \mv(z)\mv(w) \leq  
2 \int_{\T} \log^+{|\varphi|} \, \mv.
\end{equation*}
\end{cor}

\addtocounter{subsubsection}{1}
 \begin{figure}[!h]  
\begin{center}
  \includegraphics[width=100mm]{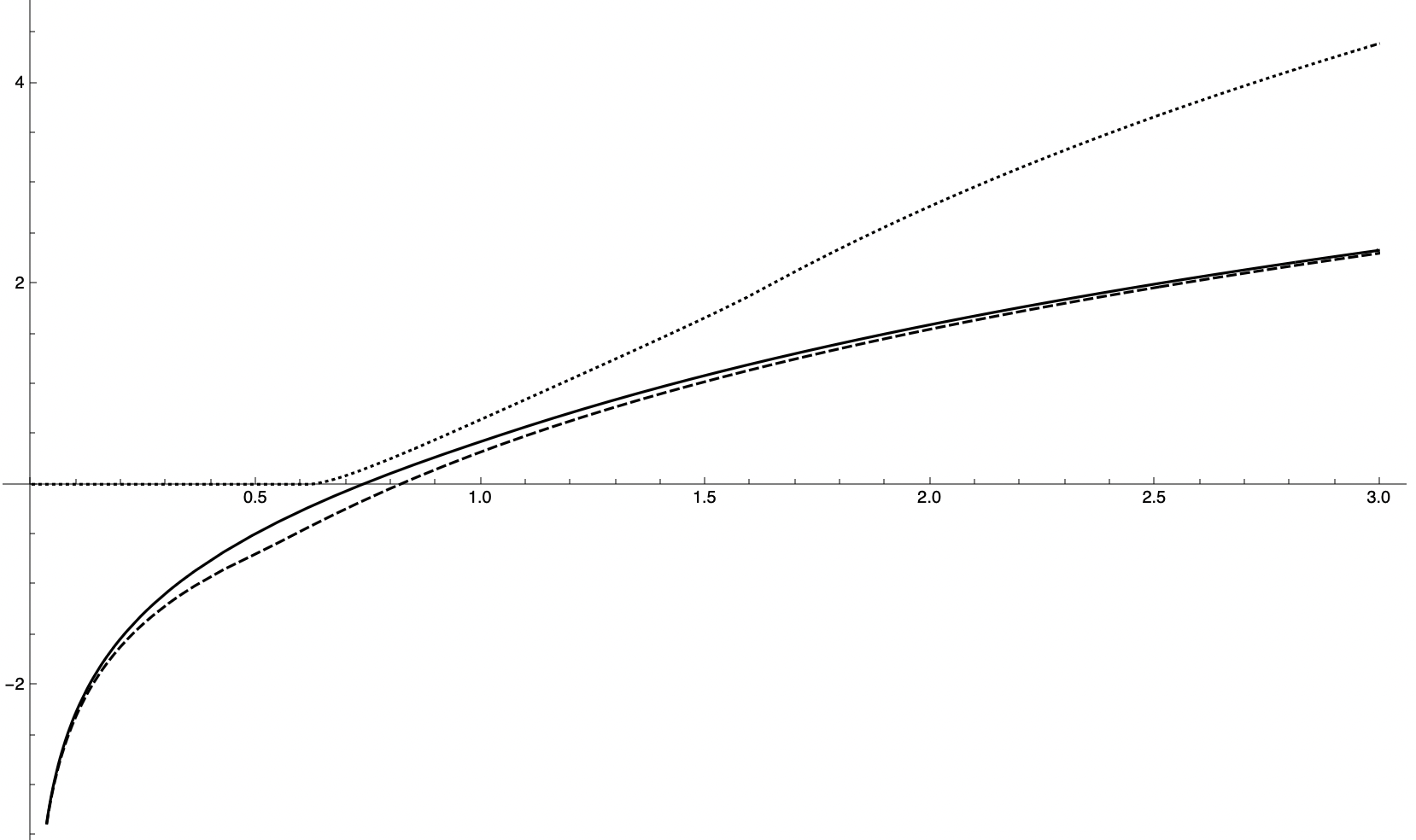}
\end{center}
\caption{The plots of the Rearrangement (in solid) and the Bost--Charles integrals (dashed) for the bivalent function $\varphi_r(z) = (rz) - (rz)^2$. Above them (dotted), 
the plot of $2\int_{\T} \log^+{|\varphi_r|} \, \mv$. }
   \label{bivalc2}
\end{figure}

\begin{example} \label{examplecompare}
Consider the function $\varphi_r(z) := rz - (rz)^2$ with the varying radius $r$.
For $r \geq 1$, 
the Bost--Charles integral amounts to
\begin{equation*}
\begin{aligned}
 \iint_{\T^2} \log{|\varphi_r(z)-\varphi_r(w)|} \, \mv(z) \mv(w) &
  = 2\log{r} +  \int_{\T} \log^+{\left| \frac{r}{1-rz} \right|} \, \mv(z) \\
 &  = 2 \log{r} + \frac{1}{\pi r} + \frac{1}{72 \pi r^2} +  \ldots. \end{aligned}
\end{equation*}
  In comparison, a computation reveals the 
rearrangement integral as the explicit function
\begin{equation*}
\begin{aligned}
\int_0^1 2t \cdot ( \log{|\varphi_r(e^{2\pi i t})|})^* \, dt  = 2\log{r} + & \frac{1}{2\pi^2} \left(  8 \, \li_3 \left( 1/r \right)  
- \li_3\left( 1/r^2 \right) \right) \\
&  = 2 \log{r} + \frac{4}{\pi^2 r} + \frac{4}{27 \pi r^2} +  \ldots.  \end{aligned}
\end{equation*} 
The comparison is pretty tight for most values of~$r$, as illustrated by Figure~\ref{bivalc2}.
 In contrast, the Nevanlinna characteristic upper bound by 
\begin{equation*}
\begin{aligned}
2 \int_{\T} \log^+{|\varphi_r|} \, \mv & = 2 \int_{|z| = r} \log^+{|z-z^2|} \, \mv(z) \\
& =   4 \log r \ \text{when} \ r \ge \frac{\sqrt{5}+1}{2} \\ 
\end{aligned}
\end{equation*}
is quite
crude.~\endofremark 
\silentcomment{
Note that these integrals can all be computed explicitly in \texttt{mathematica}
though the answers are very messy in two of the three cases. For example, specializing at~$r=1$ gives the inequalities
$$\frac{3 \sqrt{3} L(2,\chi_{-3})}{2 \pi}
> \frac{7 \zeta(3)}{2 \pi^2} > \frac{3 \sqrt{3} L(2,\chi_{-3})}{4 \pi}.
$$}
\end{example}

\begin{remark}  \label{nearly equal} 
The asymptotic equivalence of the two growth characteristic integrals at a ``big'' radius $|z| = r$
(as observed in
Example~\ref{examplecompare} and  Figure~\ref{bivalc2})  seems to be a fairly general feature
that reflects the near-rotational invariance of the $\varphi_*$ pushforward of the uniform measure $\mv$ of the expanding circle $|z| = r$, considering 
the proof of Proposition~\ref{nazarovbound} via the slice-by-slice inequality~\eqref{sliced}. An example ``in action'' is in our proof
of Theorem~\ref{logsmain} in~\S~\ref{sec:proofC}.   In this situation, the  Bost--Charles integral~(\ref{logBC}) compares
with the rearrangement integral~(\ref{logrearrange}) as follows: 
\begin{equation*}
 9.963 \sim \begin{aligned}
\iint_{\T^2} \log{|\varphi(z) - \varphi(w)|} \, \mv(z,w) \le  \int_0^1 2t \cdot ( \log{|\varphi(e^{2\pi i t})|} )^* \, dt \,
\sim 9.972,
\end{aligned}
\end{equation*}
an improvement on the order of merely a tenth of a percent. (In contrast, the  Nevanlinna characteristic quantity~$2T(\varphi)
= 2\int_{\T} \log^+{|\varphi|} \, \mv$ from~\eqref{basic basic} is in this case as big as~$\sim 14.08$.)

  Another example 
of this kind can be seen by comparing
the upper bounds occurring in the two proofs of Theorem~\ref{mainA}  
coming from Theorem~\ref{main:elementary form} and Theorem~\ref{main:BC fullconv}
respectively. The proof via Theorem~\ref{main:elementary form} (with two division points,
see~\S~\ref{sec:proofA})
gives a bound in terms of rearrangement
integrals of the form~$m < 13.731$.
 On the other hand, 
Theorem~\ref{main:BC fullconv}
    gives
   a bound in terms of arithmetic intersection numbers
   which are expressible as integrals which are slight generalizations
   of the Bost--Charles integrals (see  Lemma~\ref{BCconv-intersection}).
   With an analogous choice of division points, this leads to the bound~$m < 13.679$
    (see  Example~\ref{Ex_BCfull}, in particular Equation~\ref{secondbest}),
   an improvement of less than half a percent.
 This comparison is not literally an example
 of Proposition~\ref{nazarovbound} because of the slightly modified forms of  both
 the integrals and the bounds
 arising from the convexity argument. However,  it accurately reflects
the amount of improvement   between the 
 rearrangement
integrals and  Bost--Charles integrals that
we observed numerically
without convexity.~\endofremark
\end{remark}

 \subsection{Proof of Theorem~\ref{high dim BC convexity}}\label{highdimBC}
We firstly prove~\eqref{BCbound-highdim-cov}, which contains~\eqref{BCbound-highdim} as a special case. We will discuss at the end how to modify the proof to get \eqref{BCbound-fullconv}. 

For the following, we fix an $\epsilon>0$ and the number of variables $d$. 
Only at the end of the proof we will let, firstly, $d\rightarrow \infty$, followed by $\epsilon \rightarrow 0$.
To compare to \S~\ref{fine section} and \eqref{fine new bound elementary}, we remark that the bounds in Theorem~\ref{high dim BC convexity} only have $\ovtau(\bb)$ in the denominator using a high dimensional equidistribution feature, while the other terms are the same as those in the one dimensional bounds in \S~\ref{new slopes}. Therefore, to prove Theorem~\ref{high dim BC convexity}, we only need to incorporate the feature of a balanced index $\bi=(i_1,\ldots, i_d)$ in $\prod_{j=1}^d f_{i_j}(x_j)$, while we do not need an equidistributed degree $\bk$ in $\bx^\bk$ as in \S~\ref{aux_eqdis}. This motivates the following choice
for the Euclidean lattice to underlie our auxiliary evaluation module. 

\subsubsection{The Euclidean lattice}  \label{sec:hermitian lattice}
Consider the free $\Z$-module 
\begin{equation} \label{hd ev module}
E_D := \bigoplus_{\bi \in V_m^d(\epsilon)}   f_{\bi}(\bx) \,   \Z[1/x_1, \ldots, 1/x_d]_{\leq D},
\end{equation}
where we recall that $V_m^d(\epsilon) \subset \{1,\ldots, m\}^d$ is defined (depending on the fixed small positive constant $\epsilon$) as
\[
\begin{aligned}
V_m^d(\epsilon)=\big\{ & \bi \in \{ 1,\ldots, m \}^d \, : \, \forall i_0 \in \{1,\ldots,m\}, \\
&  d/m - \epsilon d <  \#\{1\leq j\leq d \mid i_j=i_0\} < d/m +\epsilon d \big\}; \end{aligned}
\]
and, as usual, $f_\bi (\bx) :=\prod_{j=1}^d f_{i_j}(x_j)$. Here, $\Z[1/x_1, \ldots, 1/x_d]_{\leq D}$ denotes the free $\Z$-module consisting of integer polynomials in $1/x_1, \ldots, 1/x_d$, all of whose partial degrees  ---  with respect to each $1/x_j$  ---  are at most $D$. 

By construction, $\rank E_D=(D+1)^d \cdot \# V_m^d(\epsilon)$. By Theorem~\ref{thm_MeasureConcentration} (actually, the weak law of large numbers suffices here), 
we have 
\[\lim_{d\rightarrow \infty} \frac{\# V_m^d(\epsilon)}{m^d}=1, \text{ thus } \lim_{d\rightarrow \infty}\lim_{D\rightarrow \infty} \frac{\rank E_D}{m^dD^d}=1.\]

In order to endow $E_D$ with the suitable norm, we consider the smooth projective arithmetic scheme $\cX :=(\mathbb{P}_{\Z}^1)^d$ and the natural very ample line bundle $\cL
:= \otimes_{j=1}^{d} \mathrm{pr}_j^* \, \OL(1)$ on $\cX$, where $\mathrm{pr}_j: \cX\rightarrow \P^1_\Z$ denotes the projection onto the $j$-th component. 
Then we can identify $\Gamma(\cX, \cL^{\otimes D})$ with $\Z[1/x_1, \ldots, 1/x_d]_{\leq D}$, where $x_j := Y_j/Z_j$ is an affine coordinate of the $j$-th $\P^1_\Z
= \mathrm{Proj} \, \Z[Y_j,Z_j]$.  

Recall that we are given a holomorphic map $\varphi:(\Db, 0) \rightarrow (\P^1(\C), 0)$, where ``$0$'' can mean $x_j=0$ for each copy of $\P^1_{\C}$ in $\cX_{\C}$.
The Bost--Charles metric from~\S~\ref{BC metric} is thus defined using~$\varphi$ on every factor $\mathrm{pr}_j^* \, \OL(1)$.
This induces a Hermitian line bundle structure $\ovcL = \left( \cL, \| \cdot \|_{\ovcL} \right)$ on $\cL$, and in turn, as in~\S~\ref{direct images}, a Euclidean lattice
$\Gamma_{L^2} \left(\cX, \nu; \overline{\cL}^{\otimes D} \right)$ after we fix a smooth probability measure~$\nu$ on $(\P^1)^d(\C)$. The choice of~$\nu$
is immaterial to the proof, since $D \to \infty$ for the fixed~$\nu$ and we only need to study the asymptotic leading order term given by the arithmetic
Hilbert--Samuel formula for $\ardeg \, \Gamma_{L^2} \left(\cX, \nu; \overline{\cL}^{\otimes D} \right)$. For concreteness, we pick~$\nu$ to be the smooth measure 
$$
\nu := \bigwedge_{j=1}^{d} \mathrm{pr}_j^* \, \omega_{\mathrm{FS}} =  \left( \frac{\sqrt{-1}}{2\pi}  \right)^d \frac{dz_1 \wedge \cdots 
\wedge dz_d \wedge d\bar{z_1} \wedge \cdots \wedge d\bar{z_d}}{\prod_{j=1}^d (1+|z_j|^2)^2}
$$
 with $\omega_{\mathrm{FS}} = \frac{\sqrt{-1}}{2\pi} \frac{dz \wedge d\bar{z}}{(1+|z|^2)^2}$
being the Fubini--Study form on~$\P^1(\C)$. 
Then, as in~\S~\ref{direct images}, the Euclidean norm  $\|\cdot \|$ on $\Gamma\left(\cX, \overline{\cL}^{\otimes D}\right)$ is defined by
\begin{equation} \label{norme globale}
\|s\|:= \sqrt{\int_{\cX(\C)} \|s\|_{\ovcL}^2 \, \nu}. 
\end{equation}

Just as in~\S~\ref{direct images}, we take the orthogonal direct sum~\eqref{hd ev module} of the Euclidean lattices
$$
 \Z[1/x_1, \ldots, 1/x_d]_{\leq D} \cong  \Gamma(\cX, \cL^{\otimes D}) \subset \Gamma\left( \cX, \cL^{\otimes D} \right)_{\R}
$$
induced from the above norm~\eqref{norme globale} on $\Gamma\left(\cX, \overline{\cL}^{\otimes D} \right)$. We use $\overline{E}_D = \left( E_D, \| \cdot \| \right)$ to denote this Euclidean lattice.

\subsubsection{Arithmetic degree}
The asymptotic calculation of the arithmetic degrees of direct images to $\spec{\Z}$ is the subject of the arithmetic Hilbert--Samuel formula: 
\begin{lemma}\label{ardegasymp}
As $D\rightarrow \infty$, we have the following asymptotics of arithmetic degrees: 
\begin{equation}
\begin{aligned}
\ardeg \Gamma_{L^2} \left(\cX, \nu; \overline{\cL}^{\otimes D} \right) & =\frac{d}{2} \left(\overline{\OL(1)} \cdot \overline{\OL(1)} \right) D^{d+1} +o(D^{d+1}), \\
 \ardeg \overline{E}_D  & = \frac{d\#V_m^d(\epsilon)}{2} \left(\overline{\OL(1)} \cdot \overline{\OL(1)}\right) D^{d+1} +o(D^{d+1}).
 \end{aligned}
 \end{equation}
\end{lemma}
\begin{proof}
As Euclidean lattices, $\Gamma(\cX, \overline{\cL}^{\otimes D}) \cong \otimes_{j=1}^d \Gamma\left(\P^1_\Z, \overline{\OL(1)}^{\otimes D}\right)$, where the norm on each factor
$\Gamma\left(\P^1_\Z, \overline{\OL(1)}^{\otimes D}\right)$ is the $L^2$-norm using the Hermitian line bundle $\overline{\OL(1)}$ and the Fubini--Study form $
\omega_{\textrm{FS}}$ on $\P^1(\C)$. Therefore, by Lemma~\ref{slope-tensor}, we have

\[\ardeg \Gamma\left(\cX, \overline{\cL}^{\otimes D} \right)=d (D+1)^{d-1} \ardeg \Gamma\left(\P^1_\Z, \overline{\OL(1)}^{\otimes D}\right).\]
Now the first assertion follows from the arithmetic Hilbert--Samuel formula~\eqref{P1aHS}. 
The second assertion follows from the first one by \eqref{degsum}).
\end{proof}

\begin{remark}
A proof of this calculation is also a consequence of the general arithmetic Hilbert--Samuel formula 
in Krull dimension~$d+1$. In
\cite[Theorem~1.4]{Zhang}, the Hilbert--Samuel formula is proven for an arithmetic variety of any dimension and an ample Hermitian line bundle with smooth metric of pointwise non-negative Chern form. Using the idea in \cite[\S 5]{BostL21} and \cite[\S\S 3--4]{BostCharles}, the same formula continues to hold for an ample line bundle with a $\cC^{\mathrm{b}\Delta}$ Hermitian metric of pointwise non-negative Chern form over an arithmetic variety of any dimension. 
Thus we may also deduce Lemma~\ref{ardegasymp} directly from the arithmetic Hilbert--Samuel formula for $(\cX, \overline{\cL})$: 
\[\ardeg \Gamma_{L^2}\left(\cX, \nu;  \ovcL^{\otimes D} \right)=\frac{\overline{\cL}^{d+1}}{(d+1)!} D^{d+1}+o(D^{d+1}).\]
Here, $\overline{\cL}^{d+1}$ denotes the arithmetic self-intersection number $\widehat{c}_1(\ovcL)^{d+1}.[\cX]$. In our situation, we write $\overline{\cL}=\otimes_{j=1}^d \mathrm{pr}_j^* \, \overline{\OL(1)}$, and expand the self-intersection number by multilinearity. The only nonzero terms come from crossing all but one of the  $\mathrm{pr}_j^* \, \overline{\OL(1)}$ factors once and the remaining $\mathrm{pr}_j^* \, \overline{\OL(1)}$ factor twice. 
By the projection formula, that gives us
\[\overline{\cL}^{d+1}=d \binom{d+1}{2} (d-1)! \left(\overline{\OL(1)} \cdot \overline{\OL(1)} \right),\]
and we recover Lemma~\ref{ardegasymp} by this perspective also.~\endofremark
\end{remark}

\subsubsection{Evaluation filtration}
Writing $X := \cX_{\Q}$, we can identify $\Spf \Q \llbracket \bx \rrbracket =\widehat{X}_{\mathbf{0}}$, giving in particular elements 
$$
f_\bi = f_\bi (\bx) \in \Gamma(\widehat{X}_{\mathbf{0}}, \OL_{\widehat{X}_{\mathbf{0}}}).
$$
 The space $\Gamma\left(\widehat{X}_{\mathbf{0}}, \cL^{\otimes D}\right)$ of global sections of $\cL^{\otimes D}|_{\widehat{X}_{\mathbf{0}}}$ is then naturally identified with
\[ \Gamma\left(\widehat{X}_{\mathbf{0}}, \cL^{\otimes D}\right) =\bx^{-D} \Q \llbracket \bx \rrbracket=: F_{\Q},\] where $\bx^{-D} \Q \llbracket \bx \rrbracket$ denotes the $\Q$-vector space generated by the $\bx^\bk$ with $k_j \geq -D$ for all $1\leq j \leq d$. Thus $f_\bi \, \Gamma(\cX, \cL^{\otimes D}) \subset  \Gamma(\widehat{X}_{\mathbf{0}}, \cL^{\otimes D})$, and we have the injective evaluation map 
\[\psi_D: E_D \hookrightarrow F_\Q, \quad  \left(Q_\bi\right)_{\bi\in V_m^d(\epsilon)} \mapsto \sum_{\bi\in V_m^d(\epsilon)} f_\bi Q_\bi,\]
where $Q_\bi \in \Gamma(\cX, \cL^{\otimes D})$ and $( Q_\bi)_{\bi\in V_m^d(\epsilon)}\in E_D$.

Similarly to~\S~\ref{vanfiltration}, we filter $F_\Q$ using the total vanishing order and then the lexicographical ordering within every jet space: 
\[F_\Q=F_\Q^{(\bzo)}\supseteq \cdots \supseteq F_\Q^{(\bn)} \supseteq \cdots,\]
where $\bn \in \NwithzeroB^d$, and $F_\Q^{(\bn)}:=\Cspan_{\Q}\{\bx^\bm \, : \, \bn \prec \bm + D \text{ or } \bn=\bm+D \}$. 
Here,  the total order $\prec$ on $\NwithzeroB^d$ is defined in \S~\ref{vanfiltration}, and for $\bm=(m_j)_{j=1}^d$, we define $\bm+D := (m_j+D)_{j=1}^d$. The ordering $x_1,\ldots,x_d$ of the variables used to define the lexicographical ordering is immaterial to the proof. 
 In this $\prec$-filtering notation, since $\prod_{j=1}^d x_j^{-D} \in \Gamma(\cX, \cL^{\otimes D})$ is a generator of $\cL^{\otimes D}_{\mathbf{0}}$, where we use $\cL_{\mathbf{0}}$ to denote the restriction of $\cL$ to the $\Z$-point $\bx=\mathbf{0}$, we observe that the $\mathbf{x=0}$ vanishing order
 $\ord_{\mathbf{0}}(g(\bx))$ of a $g(\bx) \in F_\Q=\Gamma(\widehat{X}_{\mathbf{0}}, \cL^{\otimes D})$ (as a regular section of $\cL^{\otimes D}|_{\widehat{X}_{\mathbf{0}}}$, not as a Laurent series in $\bx^{-D} \Q \llbracket \bx \rrbracket$) is at least $n$ if and only if $g(\bx) \in F_\Q^{(0,\ldots, 0, n)}$.
We also observe that $g(\bx) \in F_\Q^{(\bn)}$ if and only if either $\ord_0(g(\bx)) > |\bn|$ or else  $\ord_0(g(\bx))=|\bn|$ and the lowest lexicographical order term in the homogenous degree $|\bn|$ part in $g$ has an exponent vector~$\bm$ such that $\bn \prec \bm$. 
(If here one prefers to think of $g$ as a Laurent series in $\bx^{-D} \Q \llbracket \bx \rrbracket$ rather than as a section of
 $\cL^{\otimes D}|_{\widehat{X}_{\mathbf{0}}}$, one would have to shift all exponent vectors $\bn$ to $\bn-D$ in these statements.)
 \silentcomment{To compare to the convention used in Bost's paper, we endow the \emph{tangent space} --- the free $\Z$-module spanned by $d/dx_i$ --- with the standard metric declaring  $\{d/dx_i\}$ for an orthonormal basis.  Write $F^{(n)}:=F_\Q^{(0,\ldots, 0, n)}$, which is our notation for vanishing to order at least~$n$.  Instead of endowing $F^{(n)} / F^{(n+1)}$ with the tensor product of the $(n+1)^{\mathrm{st}}$ symmetric power of the dual of tangent space (this part has arithmetic degree and $\mu_{\max}$ to be both $0$) and the Hermitian line $\cL_0$ (the fiber of $\cL$ at $0$, with the restricted metric), we opt to skip the latter factor. The upshot of our change of convention is that $\mu_{\max} (\cL_0)$ inters into the height estimate. The two choices are mathematically and  trivially equivalent.}

As in \S~\ref{vanfiltration}, we use $\bn^+$ to denote the successor of $\bn$ under the total order $\prec$. The graded piece $F_\Q^{(\bn)}/F_\Q^{(\bn^+)}$ is a one dimensional $\Q$-vector space generated by the image of $\bx^{\bn-D}$ under the quotient map. The Euclidean lattice structure on $F_\Q^{(\bn)}/F_\Q^{(\bn^+)}$ is given by the free rank one $\Z$-module generated by the image of $\bx^{\bn-D}$ and the Euclidean norm with $\| \bx^{\bn-D} \|=1$. Note that these Euclidean lattice structures on graded piece are all induced from the free $\Z$-module $F=\bx^{-D} \, \Z\llbracket \bx \rrbracket$ and with the Euclidean norm that has $\{\bx^\bm\}_{\bm \in \Z_{\geq -D}^d}$ for an orthonormal basis.

Let ${E}^{(n)}_{D} := \psi_D^{-1} \left( F^{(n)}_\Q \right) \cap E_D$  denote the preimage of $F^{(\bn)}_\Q$ in $E_D$ under $\psi_D$.
Then $\psi_D$ induces injective maps 
\[\psi_D^{(\bn)}: {E}^{(\bn)}_{D}/ {E}^{(\bn^+)}_{D} \hookrightarrow F_\Q^{(\bn)}/F_\Q^{(\bn^+)}.\] 
into the one-dimensional $\Q$-vector space $F_\Q^{(\bn)}/F_\Q^{(\bn^+)}$. 
Therefore $\rank {E}^{(\bn)}_{D}/{E}^{(\bn^+)}_{D}\in \{0,1\}$ for all $\bn \in \NwithzeroB^d$. Let 
\[\mathcal{V}_{D}^d :=\{\bn\in \NwithzeroB^d \mid \rank {E}^{(\bn)}_{D}/{E}^{(\bn^+)}_{D}=1\}\]
be the vanishing filtration jumps.  We have $\# \mathcal{V}_D^d = \rank E_D$.

In the next two subsections~\S~\ref{sec:slopes archimedean} and \S~\ref{sec_BCfinite}, we provide an upper bound on  the local evaluation heights $h_v(\psi_D^{(\bn)})$ at all $v\in M_\Q$.
From the definition of the local evaluation height, we need to consider an arbitrary $(Q_\bi)_{\bi\in V_m^d(\epsilon)} \in {E}^{(\bn)}_{D} \setminus {E}^{(\bn^+)}_{D}$ and then  provide an upper bound on $\log | c_{\bn} |_v - \log \| ( Q_\bi)_{\bi\in V_m^d(\epsilon)}\|_{E_D,v}$, where $c_{\bn}$ denotes the coefficient of  $\bx^{\bn-D}$ in $s:=\sum_{\bi\in V_m^d(\epsilon)} f_\bi Q_\bi$. Here, $|\cdot |_v$ denotes the usual $v$-adic norm on $\Q$. 

\begin{df}
For a formal power series $F(t_1,\ldots,t_d) \in k \llbracket t_1, \ldots, t_d \rrbracket = \sum_{\bn} a_{\bn} \mathbf{t^n}$ over a field~$k$, the \emph{$n^{\mathrm{th}}$ order jet} of $F$ 
is the degree-$n$ homogeneous polynomial $J_n(F) \in k[t_1,\ldots,t_d]_{(n)}$ given by the sum of all the degree-$n$ terms: 
$$
J_n(F)(t_1,\ldots,t_d) := \sum_{|\bn| = n} a_{\bn} \mathbf{t^n} \in k[\mathbf{t}]_{(n)}. 
$$
\end{df}

\begin{remark} \label{split filtration remark}
In earlier work~\cite{BostFoliations,BostGerms} on the algebraization of higher dimensional formal-analytic arithmetic varieties,
it sufficed to filter the auxiliary evaluation module using only the total vanishing order at the point~$\mathbf{0}$. That we use the finer filtration with one-dimensional quotients has two advantages: 
\begin{enumerate} 
\item \label{split var comment} For the estimate of $h_v$ at an archimedean place $v$, using the product structure of $\Db^d\rightarrow \cX(\C)$, we have an easy ``variable by variable'' subharmonic estimate similar to \cite[Lemma~2.4.1]{UDC}. 
This obviates the blowing-up method in \cite[\S 4.3.2]{BostFoliations}, which gives an upper bound on 
the Mahler measure of the leading order jet polynomial $J_n(F) \in \Q[\bx]_{(n)}$, whereas the quantities
that need to be estimated are the individual coefficients~$c_{\bn}$ in that leading order jet. The discrepancy in these quantities is too sensitive in the 
dimension $d = \dim{\cX_{\Q}}$, which is fixed in~\cite{BostFoliations} whereas we want to have $d \to \infty$ at the end. 
\item The advantage for the $h_v$ estimate at a non-archimedean place $v$ is crucial. Among all $\bn$ with $|\bn|=n$, due to our specific construction of $E_D$, we will have a much better estimate of $h_v\left(\psi_D^{(\bn)} \right)$ under the condition that $\{n_j\}_{j=1}^d$ has asymptotically equidistributed components. Our complete filtration with including the lexicographical ordering allows us to take stock of this improvement.~\endofremark
\end{enumerate}
\end{remark}

\subsubsection{Archimedean estimate} \label{sec:slopes archimedean}
Recall that by the same reduction argument as in the beginning of the proof of Lemma~\ref{arch-ht-BCconv}, we may assume that $f_i(\varphi(z))$ is meromorphic on an open neighborhood of $|z|\leq 1$, for all $1\leq i \leq m$.

For ease of notation, we use $|\cdot |$ to denote the usual absolute value $|\cdot|_{\infty}$. The broader outline of the proof is similar to \cite[\S 4.3.3]{BostFoliations}, with caveat the modification we described in Remark~\ref{split filtration remark}\eqref{split var comment}. Given $s \in \psi_D^{(\bn)}(E^{\bn}_D \setminus E^{\bn^+}_D)$, we will follow~\cite[\S~2.4]{UDC} in studying the $n^{\textrm{th}}$ (leading) order jet $J_n\left(\varphi^*s\right)(\bz) \in 
\C[\bz]_{(n)}$
at the point $\bz = \mathbf{0} \in \Db^d$. Here and in the following, we use the notation $n := |\bn|$, and by a slight abuse of notation, we continue to denote $\varphi: \Db^d \rightarrow \cX(\C)$ for the analytic morphism $\bx \mapsto \left(\varphi(z_1), \ldots, \varphi(z_d) \right)$ given by $\varphi:\Db \rightarrow \C$ 
diagonally on each factor. By extension of that notation, and in a manner unifying~\S~\ref{fine section} and~\S~\ref{sec_BCconvexity}, we consider for every $\br \in (0,1]^d$ the 
analytic morphism 
$$
\varphi_{\br} : \Db^d \to \C^d \hookrightarrow \cX(\C), \qquad  \bz \mapsto \left(\varphi(r_1z_1), \ldots, \varphi(r_dz_d) \right).
$$
  We will use the Poisson--Jensen formula to bound $\log|c_\bn|$ in terms of the jet function, then relate the resulting bound to the Chern form of~$\ovcL$ by means of the Poincar\'e--Lelong formula. We follow the notations in~\S~\ref{sec_BCconvexity}, 
  and we borrow from~\eqref{slopes alpha} the notation for the slopes~$\alpha_k$. In addition, for each $\bn\in \NwithzeroB^d$, we define $r(\bn) := \left(
  r(n_1), \ldots, r(n_d) \right)$, where
  \begin{equation} \label{radii spec}
  r(t) := \left\{  \begin{array}{ll} r_k, & \quad \textrm{if } t/D \in [\alpha_{k},\alpha_{k+1}), \\ 
  1, & \quad \textrm{if } t/D > m.  \end{array} \right.
  \end{equation} 
  
Throughout this section, $\bz=(z_1,\ldots, z_d)$ denotes the coordinate on $\Db^d$. We use the trivialization $\cL^{\otimes D}|_{\C^d}\iso  \OL_{\C^d}$ of $\cL^{\otimes D}$
over $\A_{\C}^d$ given by 
$\bx^{-D} \mapsto \mathbf{1}$. Under this identification, $s/\bx^{-D}=s \cdot \bx^D =:G(\bx)$ is naturally in $\Gamma(\widehat{X}_{\mathbf{0}}, \OL_{\widehat{X}_{\mathbf{0}}})$ with vanishing order $n$. 
Since our analytic $\varphi_{r(\bn)}$-pullback is defined by $x_j=\varphi_{r(n_j)}(z_j)=(\varphi'(0) \cdot r(n_j)) z_j + O(z_j^2)$, we have by
construction that $\varphi^*_{r(\bn)} G$ has $\bz = \mathbf{0}$ vanishing order $n$, with $c_{\bn}\varphi'(0)^n \prod_{j=1}^d r(n_j) \, \bz^\bn$ for the 
lexicographically minimal term in~$J_n(\varphi_{r(\bn)}^*G)$ (as well as the overall $\prec$-minimal term in~$\varphi_{r(\bn)}^*G$). 
However, since we only assume in this theorem the \emph{meromorphy} (as opposed to the holomorphy) of
the pullbacks: $\varphi^*f_i \in \mathcal{M}(\Db)$, the analytic germ $\varphi^*_{r(\bn)} G \in \C \llbracket \bz \rrbracket$ only extends meromorphically, rather than holomorphically through~$\Db$. But if we choose $h \in \mathcal{O}(\Db)$ a holomorphic function such that $h(0) = 1$ 
and all $h \cdot \varphi^* f_i \in \mathcal{O}(\Db)$ are holomorphic, then $h(z_1) \cdots h(z_d) \cdot \left( \varphi^*_{r(\bn)}\right) G (\bz) 
\in \mathcal{O}\left(\Db^d\right)$ is holomorphic throughout~$\Db^d$, and has the same leading order jet $J_n\left(  \varphi^*_{r(\bn)} G \right)$ as~$\varphi^*_{r(\bn)} G$. 

This puts us in a position to use~\cite[Lemma 2.4.1]{UDC} for upper-bounding the requisite coefficient $|c_{\bn}|$ in terms of the Mahler measure of the 
(homogeneous) polynomial $J_n\left( \varphi^*_{r(\bn)} G \right) = J_n\left( h(z_1) \cdots h(z_d) \cdot  \varphi^*_{r(\bn)} G \right)$, in which certainly the overall lexicographically lowest term is the monomial $c_{\bn} \bz^\bn$ in the multidegree~$\bn$:
 \begin{equation} \label{Mahler jet}
\log|c_{\bn}| \leq - n \log |\varphi'(0)| - \sum_{j=1}^d n_j \log r(n_j) + \int_{\T^d} \log \left|J_n\left( h(z_1) \cdots h(z_d) \cdot  \varphi^*_{r(\bn)} G \right)\right| \, \mu_{\rm{Haar}}.
\end{equation}
We can connect this to Bost's blowing up argument in \cite[Equation (4.28) in Lemma 4.13]{BostFoliations} which shows
\begin{equation} \label{Mahler Bost}
\begin{aligned}
\int_{\T^d} \log \left|J_n\left( h(z_1) \cdots h(z_d) \cdot  \varphi^*_{r(\bn)} G \right)\right| \, \mv  &  \leq \int_{\T^d} \log \left| h(z_1) \cdots h(z_d) \cdot  \varphi^*_{r(\bn)} G \right|\, \mv \\
& =  \int_{\T^d} \log \left|   \varphi^*_{r(\bn)} G \right|\, \mv + O_h(1). 
\end{aligned}
\end{equation}
To recall Bost's argument, write $h(z_1) \cdots h(z_d) \cdot G(\varphi^*_{r(\bn)}(\bz)) =: \sum_{\bk\in \NwithzeroB^d}c'_{\bk} \bz^{\bk}$, and define 
\begin{equation} \label{homotope}
U_t(\mathbf{z}) := \sum_{\bk\in \NwithzeroB^d}c'_{\bk}t^{|\bk|-n}\bz^{\bk}, \quad \textrm{ for $t \in \C$ with } |t| \leq 1.
\end{equation}
Then, by the $\bz \mapsto t \bz$ substitution, 
 \begin{equation}  \label{rotationally}
 \int_{\T^d} \log |U_t(\bz)| \, \mu_{\rm{Haar}}(\bz)= \int_{|t|\T^d} \log \left| 
 h(z_1) \cdots h(z_d) \cdot \varphi_{r(\bn)}^* G\right| \, \mu_{{\rm Haar}}(\bz) - n \log|t|.
 \end{equation}
 In particular,  $\int_{\T^d} \log |U_t(\bz)| \, \mv(\bz)$, 
 which by~\eqref{homotope} is clearly a subharmonic
function in the single complex variable~$t$, only depends on that complex variable~$t$ through $|t|$ by means of the right-hand side of the identity~\eqref{rotationally}. Hence the value of that function at $t = 1$, which equals $\int_{\T^d} \log \left| h(z_1) \cdots h(z_d) \cdot  \varphi^*_{r(\bn)} G \right| \, \mv$  
 but is also the $t \in \T$ integral of $\int_{\T^d} \log |U_t(\bz)| \, \mv(\bz)$, is no less than its value at $t = 0$, which is the left-hand side of the requisite bound~\eqref{Mahler Bost}. 
 
Next we follow the proof idea of \cite[Theorem~10.5.3]{BostBook}, along with the discussion in \cite[\S\S4.2--4.3]{BostCharles} in order
to pass from smooth Green functions to $\cC^{\mathrm{b}\Delta}$ ones.
By the product structure of $\overline{\cL}$ and the computation from the proof of Lemma~\ref{arch-ht-BCconv},
we have
\begin{equation}  \label{Chern form d}
\begin{aligned}
& \int_{\T^d} \log\| \varphi^*_{r(\bn)} \bx^{-D}\|_{(\varphi^*_{r(\bn)}\overline{\cL})} \, \mv  =  \sum_{j=1}^d \int_{\T} \log \| \varphi^*_{r(n_j)} (x_j ^{-D}) \|_{\overline{\OL(D)}} \, \mv\\ 
&  =  \sum_{j=1}^d \left(- \int_{\Db} \log^+|z|^{-1} \, c_1\left(\varphi^*_{r(n_j)} \overline{\OL(D)} \right) + \| \varphi^*_{r(n_j)}(x ^{-D}) \|_{\varphi^*_{r(n_j)}\overline{\OL(D)}} \mid_{z=0}\right).
\end{aligned}
\end{equation}
Since $x^{-D}\mid_{x=0}$ is a $\Z$-generator of the free $\Z$-module $\OL(D)_0$, we have 
\[ \| \varphi^*_{r(n_j)}(x^{-D}) \|_{\overline{\OL(D)}} \mid_{z=0}=-\ardeg \overline{\OL(D)}_0 = - D \, \ardeg  \overline{\OL(1)}_0.\]
Putting together~\eqref{Mahler jet}, \eqref{Mahler Bost},  and~\eqref{Chern form d},  we have
\begin{equation} \label{upper d}
\begin{aligned}
&\log |c_\bn| - \int_{\T^d} \log \| \varphi^*_{r(\bn)} s\|_{(\varphi^*_{r(\bn)} \overline{\cL})} \, \mv \\
& \leq  - n \log | \varphi'(0)| - \sum_{j=1}^d n_j \log r(n_j)\\
& \quad + \int_{\T^d} \left( \log |\varphi^*_{r(\bn)} G (\mathbf{z})| -  \log \| \varphi^*_{r(\bn)} s\|_{\varphi^*_{r(\bn)}\overline{\cL}} \right)  \mv +O_h(1)\\
&= - n \log | \varphi'(0)| - \sum_{j=1}^d n_j \log r(n_j) - \int_{\T^d} \log\| \varphi^*_{r(\bn)} \bx^{-D}\|_{(\varphi^*_{r(\bn)} \overline{\cL})} \, \mv +O_h(1)\\
& =  - n \log | \varphi'(0)| - \sum_{j=1}^d n_j \log r(n_j)\\
& \qquad \qquad \qquad 
 + D\left(\sum_{j=1}^d \ardeg (\overline{\OL(1)}_0) + \int_{\Db} \log^+|z|^{-1} c_1(\varphi^*_{r(n_j)} \overline{\OL(1)}) \right) +O_h(1)\\
& =  - n \log | \varphi'(0)| - \sum_{j=1}^d n_j \log r(n_j) + D \sum_{j=1}^d \hT(r(n_j),\varphi) +O_h(1),
\end{aligned}
\end{equation}
where $\hT(r(n_j),\varphi))$ is the Bost--Charles characteristic we defined in~\ref{BC characteristic},  and the final equality derives from the projection formula \cite[Proposition~7.2.2]{BostCharles} applied to the morphism $(\iota, \varphi_{r(n_j)})$ in Remark~\ref{rmk_BCmetric} and Lemma~\ref{BCconv-intersection}. We spell out that last step in more detail. Following the definitions of pullback of Hermitian vector bundles in\cite[\S 7.1.1.1]{BostCharles} and the arithmetic intersection number in \cite[Equation~6.2.4]{BostCharles}, we have 
\[(\iota, \varphi_{r(n_j)})^*\overline{\OL(1)} \cdot  ([0], \log^+|z|^{-1})= \ardeg \overline{\OL(1)}_0 + \int_{\Db} \log^+|z|^{-1} c_1(\varphi^*(\overline{\OL(1)})).\]
By the projection formula \cite[Proposition~7.2.2]{BostCharles} together with Lemma~\ref{BCconv-intersection}, we derive the requisite equality
\[(\iota, \varphi_{r(n_j)})^*\overline{\OL(1)} \cdot  ([0], \log^+|z|^{-1}) =\hT(r(n_j),\varphi)),\]
and the final line on~\eqref{upper d} follows.

We claim that~\eqref{upper d} implies for an arbitrary $(Q_\bi)_{\bi \in V_m^d(\epsilon)}\in E^{(\bn)}_D \setminus E^{(\bn^+)}_D$ a uniform upper bound: 
\begin{equation}\label{arch-ht-bound}
\log |c_\bn| - \|(Q_\bi)_{\bi \in V_m^d(\epsilon)}\|  \leq  - n \log | \varphi'(0)| - \sum_{j=1}^d n_j \log r(n_j) + D \sum_{j=1}^d \hT(r(n_j),\varphi)) + o(D). 
\end{equation}

As $\mv$ is not a continuous measure on $\Db^d$, we begin by approximating $\log^+|z|^{-1}$ by a sequence  $(g_k)_{k \in \NwithzeroB}$ of smooth 
rotationally symmetric Green functions on~$\Db$ for the divisor~$[0]$. Precisely, by~\cite[page~268]{BostBook}, 
we choose $g_k \in C^{\infty}(\Db)$ with $\mathrm{supp}(g_k) \subset \D$, such that $g_k(z) = g_k(|z|)$ and 
 $ g_k - \log^+|z|^{-1} \to 0$ uniformly on~$\Db$. 
 Following the same argument as in the proof of Lemma~\ref{arch-ht-BCconv} (see, for instance, \cite[pages.~270--271]{BostBook} for the one-dimensional case), we write $\frac{i}{\pi} \partial \overline{\partial} g_k = -\delta_0 + \mu_k$, where $\mu_k$ is a smooth probability measure on $\Db$, and then, denoting by~$\mu_k^d$ the product measure induced from~$\mu_k$ on $\Db^d$: 
\[
\begin{aligned}
\int_{\Db^d} \log\| \varphi^*_{r(\bn)} \bx^{-D}\|_{\varphi^*_{(r(\bn)}\overline{\cL})} \, \mu_k^d & = -D \sum_{j=1}^d \left( \int_{\Db} g_k c_1(\varphi^*_{r(n_j)} \overline{\OL(1)}) +  \ardeg \overline{\OL(1)}_0 \right)
\\ & =-D \sum_{j=1}^d (\iota, \varphi_{r(n_j)})^*\overline{\OL(1)} \cdot  ([0], g_k).
\end{aligned}
\]
Moreover, since $g_k$ and $\mu_k$ are rotationally invariant, Poisson--Jensen gives:
\begin{equation*}
\begin{aligned}
\log|c_{\bn}| + n \log |\varphi'(0)| + \sum_{j=1}^d n_j \log r(n_j) & \leq  \int_{\Db^d} 
\left( \log |J_n(\varphi^*_{r(\bn)} G)(\mathbf{z})| - \log |\bz|^{\bn}  \right)  \mu_k^d\\
& = \int_{\Db^d} \log |J_n(\varphi^*_{r(\bn)} G)| \, \mu_k^d - n \int_{\Db} \log |z| \mu_k;
\end{aligned}
\end{equation*}
and our previous argument gives 
\[\int_{\Db^d} \log |J_n(\varphi^*_{r(\bn)} G)| \, \mu_k^d \leq \int_{\Db^d} \log |\varphi^*_{r(\bn)} G|  \, \mu_k^d.\]
Since $\mu_k$ is smooth on $\Db^d$ and $\Db^d$ is compact, there exists a constant $C_k>0$ depending only on 
$d, g_k, \varphi, \mathbf{r}$ (but independent of $n,D$, and~$s$), such that $C_k\varphi^*_{r(\bn)} \nu = C_k  \bigwedge_{j=1}^{d} \mathrm{pr}_j^* \varphi_{r(n_j)}^* \, \omega_{\mathrm{FS}}  \geq \mu_k^d$ as a pointwise inequality for smooth $(d,d)$-forms on $\Db^d$. Therefore, we have
\begin{equation*}
\begin{aligned}
& \log \|(Q_\bi)_{\bi \in V_m^d(\epsilon)}\|  
\\ & \geq  \max_{\bi\in V_m^d(\epsilon)} \log \|Q_\bi \| = \frac{1}{2} \max_{\bi\in V_m^d(\epsilon)} \log  \int_{\cX(\C)} \|Q_\bi\|_{\ovcL}^2 \, \nu  \\
& \geq - \frac{1}{2} \log C_k  + \frac{1}{2} \max_{\bi\in V_m^d(\epsilon)} \log \int_{\Db^d } \|\varphi^*_{r(\bn)} Q_\bi\|_{(\varphi^*_{r(\bn)}\ovcL)}^2 \, \mu_k^d \\
& \geq  - \frac{1}{2} \log C_k - \log \left\{m^d(D+1)^d \max_{\substack{ 1\leq i \leq m, \\ z\in \Db}} |f_i(z)|^d \right\} + \frac{1}{2} \log  \int_{\D^d} \|\varphi^*_{r(\bn)} s\|_{(\varphi^*_{r(\bn)} \ovcL)}^2 \, \mu_k^d \\
& >  - C_k' - d \log D + \frac{1}{2} \log  \int_{\D^d} \|\varphi^*_{r(\bn)} s\|_{(\varphi^*_{r(\bn)} \ovcL)}^2 \, \mu_k^d  \\
& \geq   - C_k' -d \log D +  \int_{\D^d} \log \| \varphi^*_{r(\bn)} s\|_{(\varphi^*_{r(\bn)} \overline{\cL})}  \, \mu_k^d,
\end{aligned}
\end{equation*}
where $C_k'>0$ is a constant only depending on $d,m, \{f_i\}, g_k, \varphi, \mathbf{r}$ (but independent of $n,D$, and~$s$), and the last inequality follows from 
the quadratic mean  ---  geometric mean inequality since $\mu_k^d$ is a probability measure on $\Db^d$ (see for instance \cite[(1.4.10)]{BostGilletSoule}).

We get: 
\begin{equation}
\begin{aligned}
\log|c_\bn| - \log \|(Q_\bi)_{\bi \in V_m^d(\epsilon)}\|  \leq  & \  -n \log |\varphi'(0)| - \sum_{j=1}^d n_j\log r(n_j) - n \int_{\Db} \log |z| \, \mu_k 
 \\
 & \ +  D \sum_{j=1}^d (\iota, \varphi_{r(n_j)})^*\overline{\OL(1)} \cdot  ([0], g_k) + o(D), \label{BCconv-ht-arch-pre}
\end{aligned}
\end{equation}
giving for any fixed $\gamma \geq 0$ independent of~$k$ and~$D$, and for all large enough $k \gg 1$: 
\begin{equation}\label{BCconv-ht-arch-precise}
\begin{aligned}
\limsup_{\substack{ D\rightarrow \infty \\ |\bn| \geq \gamma D}}' \frac{h_{\infty}(\psi_D^{(\bn)})+ \sum_{j=1}^d n_j \log r(n_j)}{D}  & \leq  \sum_{j=1}^d (\iota, \varphi_{r(n_j)})^*\overline{\OL(1)} \cdot  ([0], g_k)  \\
& - \gamma \left(\log |\varphi'(0)| + \int_{\Db} \log |z|  \, \mu_k \right),
\end{aligned}
\end{equation}
where the dash over the limit supremum indicates that we consider all 
 $\bn \in \mathcal{V}_D^d$ with $|\bn| \geq \gamma D$ and sharing some fixed $r(\bn)$ (note that there are only $(l+1)^d$ many possibilities of $r(\bn)$). 
 Here, remarking that $\log{|\varphi'(0)|} > 0$ while the uniform limit
 $g_k \to \log^+|z|^{-1}$ on $\Db$ implies
 $$
 \int_{\Db} \log |z| \, \mu_k \to 0,
 $$
 the meaning of ``large enough $k$'' is specified by the positivity of the term  $\log |\varphi'(0)| + \int_{\Db} \log |z|  \, \mu_k$.
At this point, the requisite estimate~\eqref{arch-ht-bound} follows\footnote{To be fully rigorous, the proof of~\eqref{arch-ht-bound} is completed by the paragraph below. We firstly note a mild sloppiness in our formulation due to the involvement of $\gamma$ and the requirement to only work with the $\bn$ constrained by $|\bn|/D \geq \gamma$. In practical terms, we use a large deviations bound to show that $|\bn|/D = d (m/2 +o(1))$ for most $\bn\in \mathcal{V}_D^d$, and take $\gamma = d(m/2 - \epsilon_0)$ with $\epsilon_0\rightarrow 0$ in the end. Then \eqref{BCconv-ht-arch-precise} is used for these $\bn$, while a trivial estimate applies to the leftover meagre set of~$\bn$. 

In any case, here is a rigorous completion of the proof of the requisite bound~\eqref{arch-ht-bound}. For any $\epsilon>0$, we can pick $k\gg 1$ (depending on $d$) such that $ \int_{\Db} \log |z| \, \mu_k < \epsilon$, and $\left|(\iota, \varphi_{r(\bn)})^* \left( \overline{\OL(1)} \cdot  ([0], g_k)  \right) - \sum_{j=1}^d \hT(r(n_j),\varphi))\right|< \epsilon$. For this specific $k$, we consider $D\gg1$ such that the $o(D)$ in \eqref{BCconv-ht-arch-pre} is $< \epsilon D/2$. From~\eqref{BCconv-ht-arch-precise},  
\[h_{\infty}(\psi_D^{(\bn)}) \leq  - n \log | \varphi'(0)| - \sum_{j=1}^d n_j \log r(n_j) + D \sum_{j=1}^d \hT(r(n_j),\varphi)) + \epsilon(n+D).\]
At this point~\eqref{arch-ht-bound} follows by Lemma~\ref{f_holonomic} and Shidlovsky's lemma, which give $n = |\bn| =O(dD)$ for $\bn\in \mathcal{V}_D^d$.} form the convergence
\[
\begin{aligned}
\lim_{k\rightarrow \infty} \left( (\iota, \varphi_{r(\bn)})^*\overline{\OL(1)}  \cdot  ([0], g_k) \right) 
& = \left(\overline{\OL(1)} \cdot (\iota, \varphi_{r(\bn)})_* ([0], \log^+ |z|^{-1})\right) \\
& =\sum_{j=1}^d \hT(r(n_j),\varphi). \end{aligned} \]

Armed with~\eqref{arch-ht-bound}, and using inputs from the functional bad approximability theorems in~\S~\ref{functional bad approximability} (which are in place since, in all cases, the $f_i$ are holonomic functions; noting that the proof of Lemma~\ref{f_holonomic} entails automatic holonomicity from the condition $\log |\varphi'(0)| > \sum_{h=1}^r \max_{1\leq i \leq m} b_{i,h}$), we now estimate the total contribution 
$$
\sum_{\bn \in \NwithzeroB^d} \rank\left(E^{(\bn)} / E^{(\bn^+)}\right) \cdot  h_{\infty}(\psi_D^{(\bn)})= \sum_{\bn \in \mathcal{V}_D^d} h_{\infty}(\psi_D^{(\bn)})
$$
 of the archimedean height showing in the right-hand side of Bost's slopes inequality \eqref{slope-inequality}. 

For any $\epsilon'>0$, Lemma~\ref{lem_Shidlovsky} shows that all $D\gg_{\epsilon', \{f_i\}} 1$ satisfy 
\[\mathcal{V}_D^d \subset \left[0,(m+ \epsilon')D\right]^d.\]
 (In the Lemma, we may pick $\varepsilon := \epsilon'/2$ and we consider $D\gg_{\epsilon', \{f_i\}} 1$ such that $\epsilon' D/2 > C(\varepsilon)$.)
By \eqref{arch-ht-bound}, which by definition is an upper bound on the $\bn^{\mathrm{th}}$ archimedean evaluation height $h_{\infty}(\psi_D^{\bn})$, we have:
\begin{equation}  \label{arch up up}
\begin{aligned}
&\sum_{\bn \in \mathcal{V}_D^d} h_{\infty}(\psi_D^{(\bn)}) \\
&\leq  -  \log |\varphi'(0)| \left(\sum_{\bn \in \mathcal{V}_D^d} |\bn| \right) +D \sum_{\bn \in \mathcal{V}_D^d} \sum_{j=1}^d\left(-\frac{n_j}{D} \log r(n_j) + \hT(r(n_j),\varphi)) \right) +o(D^{d+1}).
\end{aligned}
\end{equation}

We estimate  the quantities
\[\sum_{\bn \in \mathcal{V}_D^d} |\bn|,  \sum_{\bn \in \mathcal{V}_D^d} \sum_{j=1}^d (n_j/D)\log r(n_j), 
\sum_{\bn \in \mathcal{V}_D^d} \sum_{j=1}^d \hT(r(n_j),\varphi))\]
 using the following consequence of Theorem~\ref{thm_MeasureConcentration}, to the effect that most $\bn \in \mathcal{V}_D^d$ have uniformly distributed components. Similarly to~\S~\ref{aux_eqdis}, recall from the statement of Theorem~\ref{high dim BC convexity}, we use $ P_{\varepsilon}^d(N)   \subset [0, N]^d \cap \Z^d$ to denote the subset  of those $\bn$ for which the normalized $\left( [0,1), \mu_{\mathrm{Lebesgue}} \right)$ discrepancy of $\{n_i/N\}_{i=1}^{d}$ is $\leq \varepsilon$, and $B_d^{\varepsilon}(N)$ to denote the complement of $P_{\varepsilon}^d(N)$ in~$ [0, N]^d \cap \Z^d$.

\begin{lemma}\label{lem_n_eqdistr}
There is a function $c : (0,1) \to (0,1)$, such that the following holds. 

Consider an $\epsilon'' \in (0,1)$. 
Then, for all $\epsilon' > 0$ small enough with respect to $\epsilon''$, 
\[\lim_{D\rightarrow \infty} \frac{\# \left\{ \mathcal{V}_D^d \cap B_d^{\epsilon''}\left((m+\epsilon')D\right) \right\}}{\# \mathcal{V}_D^d}=O\left(e^{-c(\epsilon'') d} \right),\]
where 
 the implicit coefficient is absolute. 
 \end{lemma}
\begin{proof}
Note that $\# \left\{ \mathcal{V}_D^d \cap B_d^{\epsilon''}\left((m+\epsilon')D\right) \right\} \leq \# B_d^{\epsilon''}\left((m+\epsilon')D\right) $, and recall that
\[\lim_{d\rightarrow \infty}\lim_{D\rightarrow \infty}\frac{\# \mathcal{V}_D^d}{m^d D^d}=\frac{\rank E_D}{m^d D^d}=1.\]
Hence, it suffices to show that $\lim_{D\rightarrow \infty}\frac{\# B_d^{\epsilon''}\left((m+\epsilon')D\right) }{m^d D^d}=O(e^{-c d})$ for some $c = c(\epsilon'') > 0$ and a small enough $\epsilon'$. 

By Theorem~\ref{thm_MeasureConcentration}, we have
\[\lim_{D\rightarrow \infty}\# B_d^{\epsilon''}\left((m+\epsilon')D\right) / ((m+\epsilon')D)^d = O(e^{-c_0(\epsilon'') d})),\]
 with a certain $c_0(\epsilon'') > 0$ and an absolute implicit coefficient. 
For $\epsilon'$ sufficiently small in terms of $c_0(\epsilon'')$, we will have $\lim_{D\rightarrow \infty}((m+\epsilon')D)^d/ (mD)^d = O(e^{c_0(\varepsilon'') d/2})$. We obtain the desired bound with $c :=c_0/2$.
\end{proof}

Now, for an arbitrary $\epsilon''>0$, we pick the sufficiently small $\epsilon'>0$ as guaranteed by Lemma~\ref{lem_n_eqdistr}, 
and apply Lemma~\ref{lem_Shidlovsky} as discussed above to obtain that $\mathcal{V}_D^d \subset \left[0,(m+ \epsilon')D\right]^d$ for $D\gg 1$. We obtain: 
\begin{equation*}
\begin{aligned}
\lim_{d\rightarrow \infty}  & \lim_{D\rightarrow \infty}  \left\{  \frac{\sum_{\bn \in \mathcal{V}_D^d} |\bn|}{d m^d D^{d+1}}  \right\} \geq  \  \lim_{d\rightarrow \infty} \lim_{D\rightarrow \infty} \frac{\sum_{\bn \in \mathcal{V}_D^d\cap P_d^{\epsilon''}\left((m+\epsilon')D\right)} |\bn|}{d D \#\mathcal{V}_D^d}\\
\geq & \  \lim_{d\rightarrow \infty} \lim_{D\rightarrow \infty} \frac{(\#\mathcal{V}_D^d\cap P_d^{\epsilon''}\left((m+\epsilon')D\right))(m+\epsilon')(1-2\sqrt{\epsilon''})/2}{ \#\mathcal{V}_D^d}\\
=  & \ (m+\epsilon')(1-2\sqrt{\epsilon''})/2.
\end{aligned}
\end{equation*}
Here, the second inequality follows upon remarking that the definition of the discrepancy function implies  $|\bn| \geq dD(m+\epsilon')(1-2\sqrt{\epsilon''})/2$ for all $\bn \in P_d^{\epsilon''}\left((m+\epsilon')D\right)$;  
\silentcomment{here we could more directly appeal to Hoeffding's inequality in Lemma~\ref{lem_n_eqdistr}, however, we choose to state the lemma in discrepancy notation to be consistent with the needs of next subsubsection; to deduce the expectation of~$|\bn|$ from the discrepancy estimate, we can take intervals of length $1/N$ with $1/N= \sqrt{\epsilon''}$;}
and the last equality follows from Lemma~\ref{lem_n_eqdistr} which implies that, for fixed $\epsilon'', \epsilon'$,
\[\lim_{d\rightarrow \infty} \lim_{D\rightarrow \infty} \frac{\#\mathcal{V}_D^d\cap P_d^{\epsilon''}\left((m+\epsilon')D\right)}{ \#\mathcal{V}_D^d}=\lim_{d\rightarrow \infty} \left\{ 1-O(e^{-c(\epsilon'')d}) \right\} =1.\]
Now we let $\epsilon''\rightarrow 0$ (this will force $\epsilon'\rightarrow 0$). We get: 
\begin{equation}  \label{ord van eval}
 \lim_{d\rightarrow \infty} \lim_{D\rightarrow \infty} \left\{ \frac{\sum_{\bn \in \mathcal{V}_D^d} |\bn|}{d m^d D^{d+1}}  \right\} \geq \frac{m}{2}.
\end{equation}
Similarly, still directly from the definition of the discrepancy function, we have the following evaluation for all $\bn \in \mathcal{V}_D^d \cap P_d^{\epsilon''}\left((m+\epsilon')D\right)$: 
\begin{equation}
\begin{aligned}
&\sum_{j=1}^d\left\{-(n_j/D) \log r(n_j) + \hT\left(r(n_j),\varphi)\right) \right\}  \\
&= \frac{d}{m} \sum_{k=0}^{l}\left\{ (\alpha_{k+1} - \alpha_k)\hT\left(r_k,\varphi\right) - \frac{1}{2}(\alpha_{k+1}^2 - \alpha_k^2) \log r_k\right\} + O\left(d(\epsilon' + \sqrt{\epsilon''})\right),
\end{aligned}
\end{equation}
recalling the notations~\eqref{slopes alpha} and~\eqref{radii spec}. (The implicit coefficient here depends linearly on $m|\log{r_0}| + \hT(1,\varphi)$.)
A partial summation now gives 
\begin{equation}  \label{main piece eval}
\begin{aligned}
\lim_{d\rightarrow \infty}  & \lim_{D\rightarrow \infty} \left\{ \frac{\sum_{\bn \in \mathcal{V}_D^d \cap P_d^{\epsilon''}\left((m+\epsilon')D\right)} \sum_{j=1}^d\left(-(n_j/D) \log r(n_j) + \hT(r(n_j),\varphi)) \right)}{dm^d D^{d}} \right\} \\
& =\hT(1,\varphi)- \frac{1}{2m} \sum_{k=1}^l \frac{ (\hT(r_k,\varphi)- \hT(r_{k-1},\varphi))^2 }{\log{r_k} - \log{r_{k-1}}} + O((\epsilon' + \sqrt{\epsilon''})),
\end{aligned}
\end{equation}
and the last error term goes to $0$ once we take $\epsilon'' \rightarrow 0$ (which implies $\epsilon' \to 0$).

Further, for all $\bn \in \mathcal{V}_D^d$, we certainly have the following trivial bound
\[\sum_{j=1}^d\left(-(n_j/D) \log r(n_j) + \hT\left(r(n_j),\varphi)\right) \right) \leq d(m+\epsilon') |\log r_0| + d \hT(1,\varphi),\]
and thus
\begin{equation} \label{supplemental piece eval}
\lim_{d\rightarrow \infty} \lim_{D\rightarrow \infty} \frac{\sum_{\bn \in \mathcal{V}_D^d \cap B_d^{\epsilon''}\left((m+\epsilon')D\right)} \sum_{j=1}^d\left(-(n_j/D) \log r(n_j) + \hT(r(n_j),\varphi)) \right)}{dm^d D^{d}} =0.
\end{equation}

Combining~\eqref{arch up up}, \eqref{ord van eval}, \eqref{main piece eval}, and~\eqref{supplemental piece eval}, we arrive at our total archimedean evaluation height bound: 
\begin{equation}\label{BChighdim-conv-arch}
\begin{aligned}
\lim_{d\rightarrow \infty} & \lim_{D\rightarrow \infty}   \left\{  \frac{\sum_{\bn \in \mathcal{V}_D^d} h_{\infty}(\psi_D^{(\bn)})}{dm^d D^{d+1}} \right\} \\
&  \leq   - \frac{m}{2}\log |\varphi'(0)| + \hT(1,\varphi) -  \frac{1}{2m} \sum_{k=1}^l \frac{ (\hT(r_k,\varphi)- \hT(r_{k-1},\varphi))^2 }{\log{r_k} - \log{r_{k-1}}}. 
\end{aligned}
\end{equation}

\subsubsection{Non-archimedean estimate}\label{sec_BCfinite}
The main idea here is similar to~\S~\ref{denominator arithmetic}. 

Let $h_{\mathrm{fin}}(\psi_D^{(\bn)})$ denote $\sum_{v\nmid \infty} h_v(\psi_D^{(\bn)})$. For $\bn\in \mathcal{V}_D^d$ and a prime $p$, by definition, $h_p(\psi_D^{(\bn)})$ is $\log p$ times the maximal $p$-adic valuation $v_{p,\bn}$ of the denominators of the ($\bn-D$)-th coefficient of $\sum_{\bi\in V_m^d(\epsilon)} f_\bi Q_\bi$ across all $(Q_\bi)_{\bi\in V_m^d(\epsilon)}\in E_D$. Since all $Q_\bi(\bx)$ are $\Z$-linear combinations of monomials $\bx^{\bk}$ with $\bk \in [-D,0]^d$, it follows that $v_{p,\bn}$ is at most the maximum of the $p$-adic valuations of the denominators of the $\bx^{\bm}$ coefficients of all $f_\bi(\bx)$ for all $\bm$ with $(n_j - D)^+ \leq m_j \leq n_j$ for all $1\leq j \leq d$; here for once we write $(n_j-D)^+ := \max(n_j-D,0)$. 
We consider separately the $\prod_{h=1}^{r} [1,\ldots,b_{i,h}n]$ and the $n^{e_i}$ pieces of the $p$-denominators of the coefficients of $f_i(x)$: 
\[v_{p,\bn}^\flat:=\max_{\bm \preceq \bn, \, \bi\in V_m^d(\epsilon)} \left\{\sum_{j=1}^d \mathrm{val}_p([ 1, \ldots, b_{i_j,1} \cdot m_j] \cdots [1,\ldots, b_{i_j,r} \cdot m_j])\right\},\]
\[v_{p,\bn}^\sharp:=\max_{ \substack{  \bi\in V_m^d(\epsilon)  \\
(n_j -D)^+ \leq m_j \leq n_j, \forall 1\leq j \leq d} } \left\{\sum_{j=1}^d e_{i_j}\mathrm{val}_p(\max\{m_j,1\})\right\},\]
where $\mathrm{val}_p$ denotes the usual $p$-adic valuation with $\mathrm{val}_p(p)=1$. 
Here we use the convention that for $m_j=0$, we set $[1,\ldots, b_{i_j,h} \cdot m_j]=1$.
By definition, we have
\begin{equation}
v_{p,\bn} \leq v_{p,\bn}^\flat+ v_{p,\bn}^\sharp.
\end{equation}

We continue with the notations from~\S~\ref{sec:slopes archimedean}; in particular, $\mathcal{V}_D^d \subset [0, (m+\epsilon')D]^d$ by Lemma~\ref{lem_Shidlovsky}.  We firstly discuss $v_{p,\bn}^\flat$.
For the case $\bn \in \mathcal{V}_D^d \cap  B_d^{\epsilon''}\left((m+\epsilon')D\right)$, we stick to the trivial bound.  Observe that $\prod_{h=1}^r [1,\ldots, (\max_{1\leq i \leq m} b_{i,h}) \cdot n]$ is a multiple of the denominators of the $x^n$-coefficients of all $f_1,\ldots, f_m$. 
By the prime number theorem, it follows that
\begin{equation}   \label{triv fin 1}
\sum_p v_{p,\bn}^\flat \log p \leq \left(\sum_{h=1}^r \max_{1\leq i \leq m} b_{i,h}\right) |\bn| +o(|\bn|).
\end{equation}
Summing over all $\bn \in \mathcal{V}_D^d \cap  B_d^{\epsilon''}\left((m+\epsilon')D\right)$, so that in particular $|\bn| \leq d(m+\epsilon')D$, we have, as $d\rightarrow \infty$:
\begin{equation}  \label{B part}
\begin{aligned}
&\limsup_{D\rightarrow \infty} \left\{ \frac{\sum_{\bn \in \mathcal{V}_D^d\cap B_d^{\epsilon''}\left((m+\epsilon')D\right)}  \sum_{p} v_{p,\bn}^\flat \log p }{dm^d D^{d+1}} \right\}
\\ 
&\leq \limsup_{D\rightarrow \infty} \left\{ \frac{\left(\sum_{h=1}^r \max_{1\leq i \leq m} b_{i,h}\right)(m+\epsilon')Dd}{dD} \frac{\#\mathcal{V}_D^d\cap B_d^{\epsilon''}\left((m+\epsilon')D\right)}{m^dD^d}  
\right\}\\
& = O\left(\left(\sum_{h=1}^r \max_{1\leq i \leq m} b_{i,h}\right) (m+\epsilon') e^{-c(\varepsilon'')d}  \right) = o_{d \to \infty}(1). 
\end{aligned}
\end{equation}
For the case $\bn \in \mathcal{V}_D^d \cap P_d^{\epsilon''}\left((m+\epsilon')D\right)$: at a fixed $p$, our assumption on the denominator types of the~$f_i$ says that the $p$-adic valuation of the denominators of the coefficients of all the monomials with exponent vectors $\preceq \bn$ in $f_{\bi}$ is at most
\[\sum_{j=1}^d \mathrm{val}_p\left([ 1, \ldots, b_{i_j,1} \cdot n_j] \cdots [1,\ldots, b_{i_j,r} \cdot n_j]\right);\]
 for $p > \sqrt{(\max_{i,h} b_{i,h} )(m+\epsilon') D}$, this equals $\# \{(j,h) \, : \, p \leq b_{i_j,h} n_j\}$.

Therefore, by the prime number theorem and the definition of $\ovtau(\bb)$, we have
\begin{equation}  \label{P part}
\lim_{d\rightarrow \infty}\lim_{D\rightarrow \infty} \left\{ \frac{\sum_{\bn \in \mathcal{V}_D^d\cap P_d^{\epsilon''}\left((m+\epsilon')D\right)} \sum_{p} v_{p,\bn}^\flat \log p }{dm^d D^{d+1}} \right\} \leq \frac{1}{2}(m+\epsilon')\ovtau(\bb) +o_{\epsilon'' \to 0}(1).
\end{equation}
Combining~\eqref{B part} and~\eqref{P part}, we get (noting that $\epsilon'' \to 0$ has, by implication, $\epsilon' \to 0$): 
\begin{equation}\label{BChighdim-conv-fin-1}
 \lim_{d\rightarrow \infty}\lim_{D\rightarrow \infty} \left\{ \frac{\sum_{\bn \in \mathcal{V}_D^d} \sum_{p} v_{p,\bn}^\flat \log p }{dm^d D^{d+1}}  \right\}
 \leq \frac{1}{2}m \ovtau(\bb). 
\end{equation}

Finally, we turn to $v_{p,\bn}^\sharp$. 
For any $\bn\in \mathcal{V}_D^d \subset [0, (m+\epsilon')D]^d$, at a given $p$, we defined $v_{p,\bn}^\sharp$ as the maximal valuation of the denominators of the coefficients of all monomials $\prod_{j=1}^{d} x_j^{m_j}/m_j^{e_{i_j}}$ ranging over all $\bi \in V_m^d(\epsilon)$ and all exponent vectors $\bm$ 
such that $(n_j-D)^+ \leq m_j \leq n_j$ for all $1\leq j \leq d$.
It satisfies 
\begin{equation}
\begin{aligned}
v_{p,\bn}^\sharp
&  \leq \max_{\bi \in V_m^d(\epsilon)} \left\{ \sum_{j=1}^d e_{i_j}\mathrm{val}_p([\max\{n_j-D,1\}, \ldots, n_j])  
\right\} \\
& \leq \left(\max_{1\leq i \leq m} e_i \right)  \sum_{j=1}^d \mathrm{val}_p\left([\max\{n_j-D,1\}, \ldots, n_j]\right).
\end{aligned}
\end{equation}
It is here that we use the cutoff parameter~$\xi$ of the defining formula~\eqref{taue}. 
Summing over all $p\geq \xi D$, we derive the estimate
\begin{equation}
\sum_{p \geq \xi D}   v_{p,\bn}^\sharp \log p \leq \left(\max_{1\leq i \leq m} e_i\right)  \sum_{j=1}^d \sum_{p\geq \xi D} \mathrm{val}_p([\max\{n_j-D,1\}, \ldots, n_j])\log p.
\end{equation}
Since the $\bn \in P_d^{\epsilon''}\left((m+\epsilon')D\right)$ have uniformly distributed components up to normalized discrepancy $\leq \epsilon''$, 
Lemma~\ref{lcm-2} with the prime number theorem 
yields 
\begin{equation}
\begin{aligned}
& \sum_{\bn\in \mathcal{V}_D^d \cap P_d^{\epsilon''}\left((m+\epsilon')D\right)}   \left\{ \sum_{p \geq \xi D}  v_{p,\bn}^\sharp \log p \right\} \\ &  \leq   d \left(\max_{1\leq i \leq m} e_i\right) (m+\epsilon')^{d-1}D^{d+1} I^{m+\epsilon'}_{\xi}(\xi) \left(1+ O(\sqrt{\epsilon''})+o(D^{d+1}) \right).
\end{aligned}
\end{equation}
The complementary meagre set of~$\bn$ is once again handled by the trivial estimate: 
\begin{equation}
\begin{aligned} 
\sum _{\bn\in \mathcal{V}_D^d \cap B_d^{\epsilon''}\left((m+\epsilon')D\right)}
\left\{ \frac{\sum_{p \geq \xi D}  v_{p,\bn}^\sharp \log p}{dm^dD^{d+1}} \right\}  
& = O\left(  \left(\max_{1\leq i \leq m} e_i\right)(m+\epsilon') e^{-c(\epsilon'')d} \right) \\
&  = o_{d \to \infty}(1). \end{aligned}
\end{equation}

Therefore, taking $d \to \infty$ and noting again that $\epsilon'' \to 0$ entails $\epsilon' \to 0$, we arrive at the limit majorization
\begin{equation}
\lim_{d\rightarrow \infty} \lim_{D\rightarrow \infty} \left\{ \frac{\sum_{\bn\in \mathcal{V}_D^d}\sum_{p \geq \xi D}  v_{p,\bn}^\sharp \log p }{dm^d D^{d+1}}
 \right\} \leq \frac{\left(\max_{1\leq i \leq m} e_i\right) I^{m}_{\xi}(\xi)}{m}.
\end{equation}

It remains to estimate the $p\leq \xi D$ contribution. Here we use the fact that the multi-index $\bi\in V_m^d(\epsilon)$ is $\epsilon$-balanced. Therefore, once again by the prime number theorem (indicating that we only count $\mathrm{val}_p=1$), we have asymptotically as $D\rightarrow \infty$: 
\begin{equation*}
\begin{aligned}
\sum_{p\leq \xi D}  v_{p,\bn}^\sharp \log p & \leq  \sum_{p\leq \xi D} \left\{ \max_{\bi  \in V_m^d(\epsilon)} \sum_{j=1}^d e_{i_j} \mathrm{val}_p([1,\ldots, n_j]) \right\}
 \log p  \\ 
 &\leq \sum_{p\leq \xi D}  \left \{ \max_{\bi\in V_m^d(\epsilon)} \sum_{j=1}^d e_{i_j} \right\} \log p + o(D) \\
& =  d \xi D \left( \frac{\sum_{i=1}^m e_i}{m} + O(\epsilon) \right) +o(D).
\end{aligned}
\end{equation*}

Therefore, once we let $D \to \infty$ followed by $d \to \infty$ and then $\epsilon'' \rightarrow 0$ (entailing $\epsilon' \to 0$), we derive
\begin{equation}\label{BChighdim-conv-fin-2}
\lim_{d\rightarrow \infty} \lim_{D\rightarrow \infty}\frac{\sum_{\bn\in \mathcal{V}_D^d} v_{p,\bn}^\sharp \log p }{dm^d D^{d+1}} \leq \frac{\xi (\sum_{i=1}^m e_i) + O(\epsilon) + (\max_{1\leq i \leq m} e_i) I^{m}_{\xi}(\xi)}{m}.
\end{equation}

\begin{proof}[Conclusion of the proof of Theorem~\ref{high dim BC convexity}]  \label{sec_slope_conclusion}
We derive the desired bound from Bost's slopes inequality~\eqref{slope-inequality}, upon collecting Lemma~\ref{ardegasymp} and the estimates~\eqref{BChighdim-conv-arch}, \eqref{BChighdim-conv-fin-1}, and \eqref{BChighdim-conv-fin-2}, and finally by letting $\epsilon \to 0$ in the end. 

The proof of~\eqref{BCbound-fullconv} differs only in the archimedean evaluation height estimate, upon replacing $\ovcL$ by $\ovcL' =\prod_{h=0}^l \ovcL_{r_h}^{\otimes s_h}$ and inputting the bound from~\S~\ref{sec_BCfull}. 
\end{proof}

\begin{remark}
In the case (such as we have in all our applications in this paper) that the $m$-dimensional $\Q(x)$-vector space $\Cspan_{\Q(x)}\{f_1, 
\ldots, f_m\}$ is closed under differentiation, we can apply the more elementary Shidlovsky lemma (Theorem~\ref{Shidlovsky}, which is much easier to prove and effectivize than Theorem~\ref{KolchinSolved}); in this situation, the statement of Lemma~\ref{lem_Shidlovsky} applies even with $\epsilon=0$. In this situation, the large deviations bound quoted from Theorem~\ref{thm_MeasureConcentration} in the proof can be replaced by the most rudimentary weak law of large numbers.~\endofremark
\end{remark}

\begin{remark}\label{BCboundK}
The proof immediately gives the following formal generalization to a number field~$K$. 
For each $\sigma: K\hookrightarrow \C$, we consider a holomorphic mapping $\varphi_\sigma: (\Db, 0) \to (\C,0)$  with $\varphi'_\sigma(0)\neq 0$. 
Assume there exists an $m$-tuple $f_1, \ldots, f_m \in K \llbracket x \rrbracket$ of $K(x)$-linearly independent formal functions with denominator types of the form
\[
f_i(x) = a_{i,0} + \sum_{n=1}^{\infty} a_{i,n} \frac{x^n}{n^{e_i} [ 1, \ldots, b_{i,1} \cdot n] \cdots [1,\ldots, b_{i,r} \cdot n]}, \qquad a_{i,n} \in \OL_K, 
\]
where $e_i, b_{i,j}$ are the same as in Theorem~\ref{main:elementary form}, 
and such that for all $i \in \{1,\ldots,m\}$ and $\sigma:K\hookrightarrow \C$,we have $f_i(\varphi_\sigma(z)) \in \C \llbracket z \rrbracket$ convergent on $|z| < 1$.
 If 
\[
\frac{1}{[K:\Q]}\sum_{\sigma:K\hookrightarrow \C} \log |\varphi'_\sigma(0)| > \sigma_m ,\] 
then all $f_i$ are holonomic functions, and 
\[
m  \leq  \frac{\sum_{\sigma:K\hookrightarrow \C}  \iint_{\T^2} \log|\varphi_{\sigma}(z)-\varphi_{\sigma}(w)| \, \mv(z) \mv(w) }{ (\sum_{\sigma:K\hookrightarrow \C}  \log{|\varphi_{\sigma}'(0)|}) - [K:\Q](\ovtau(\bb)+ \tau^\sharp(\be)) }.
\]
The convexity improvements also extend in the obvious way.~\endofremark
\end{remark}

\subsection{The slopes method in Theorem~\ref{main:elementary form}}\label{sec_slope for rearrangement}
We showed in~\S~\ref{rmk_Nazarov} that Theorem~\ref{main:BC form} formally implies the corresponding particular case Corollary~\ref{basic main corollary} of Theorem~\ref{main:elementary form}. In this brief section,
we comment how Theorem~\ref{main:elementary form} can be more directly recovered in the framework of the preceding proof.

In constructing the Euclidean lattice $E_D$, in addition to only considering the split-variable products~$f_{\bi} $ with $\bi\in V_m^d(\epsilon)$, we may also  ---  as in the Thue--Siegel lemma construction in~\S~\ref{aux_eqdis}  ---  constrict the monomials~$\bx^\bk$ to have exponent vectors~$\bk$ with uniformly distributed components~$\{k_i\}$. More precisely,
define the free $\Z$-module: 
\[E_D := \bigoplus_{\bi \in V_m^d(\epsilon),  \, \bk/D \in P^d_\epsilon}   f_\bi \,   \Z\bx^{\bk},
\]
where $P_\epsilon^d$ was defined in \eqref{typical}. Then indeed we have  the requisite double limit
\[\lim_{d\rightarrow \infty} \lim_{D\rightarrow \infty} \left\{  \frac{\rk E_D}{m^d D^d} \right\}=1.\]
We equip $E_D$ with the Euclidean metric that makes $\{f_\bi \bx^\bk\}_{\bi \in V_m^d(\epsilon), \bk/D \in P^d_\epsilon}$ an orthonormal basis.

In Theorem~\ref{main:elementary form}, we are given 
a set of $l+1$~holomorphic mappings~\eqref{relaxedifholonomic}, and corresponding division point parameters
 $0=\gamma_0 < \gamma_1 < \ldots < \gamma_l <\gamma_{l+1}:=m$. 
  Like in~\S~\ref{seeding}, the meaning of these numbers is that we will use the Poisson--Jensen formula for $\varphi_k(z_j)$ for the unique $k = k(j)$ determined by $n_j/D \in [\gamma_k, \gamma_{k+1})$; we use $\varphi_l(z_j)$ for $n_j/D \in [m, m+\epsilon')$. We continue to use $\mathcal{V}_D^d$ to denote the set of $\bn$ such that $\rank E_D^{(\bn)}/ E_D^{(\bn^+)} =1$; recall that for $D\gg 1$, we have $\mathcal{V}_D^d \subset [0, (m+\epsilon')D]^d$. For all $\bn \in \mathcal{V}_D^d$, the ensuing evaluation height estimate is
 \[
 \begin{aligned}
 h_{\infty}(\psi_D^{(\bn)})  &   \leq D \int_{\T^d} \max_{\bt \in P_\epsilon^d} \left\{ 
\sum_{j=1}^d t_j \log |\varphi_{k(j)}(z_j)|  \right\} \mv \\
&  - |\bn| \log{|\varphi'_l(0)|} - \sum_{j=1}^d n_j \log{|\varphi_{k(j)}'(0)/\varphi'_l(0)|} + o(D). \end{aligned} \]
As in any case we have the trivial bound 
\[\max_{\bt \in P_\epsilon^d} \left\{ \sum_{j=1}^d t_j \log |\varphi_{k(j)}(z_j)| \right\} \leq d \max_{k, \T} \log |\varphi_k|,\]
an argument similar to~\S~\ref{sec:slopes archimedean} shows 
\begin{equation*}
\begin{aligned}
\lim_{d\rightarrow \infty} \lim_{D\rightarrow \infty} & \left\{ \frac{\sum_{\bn \in \mathcal{V}_D^d} h_{\infty}(\psi_D^{(\bn)})}{dm^d D^{d+1}} \right\}
 \leq  \lim_{d\rightarrow \infty} \lim_{D\rightarrow \infty} \left\{ \frac{\sum_{\bn \in \mathcal{V}_D^d\cap P_d^{\epsilon''}\left((m+\epsilon')D\right)} h_{\infty}(\psi_D^{(\bn)})}{dm^d D^{d+1}}
 \right\} \\
&  \leq  \lim_{d\rightarrow \infty}   \left\{ \frac{1}{d} \int_{\T^d} \max_{  \substack{ \bt \in P_\epsilon^d, \\ \bn \in P_d^{\epsilon''}\left((m+\epsilon')D\right)}} \left\{ \sum_{j=1}^d t_j \log |\varphi_{k(j)}( z_j)| \right\} \, \mv  \right\}  \\ 
& \quad  -  \frac{m}{2}\log |\varphi_l'(0)|
 - \frac{1}{2} \sum_{k=0}^l (\gamma_{k+1}^2 - \gamma_k^2) \log{\frac{|\varphi_k'(0)|}{|\varphi'_l(0)|}}   +O(\epsilon' + \sqrt{\epsilon''})
\end{aligned}
\end{equation*}

Given $\bn \in P_d^{\epsilon''}\left((m+\epsilon')D\right)$, asymptotically as $d\rightarrow \infty$, almost all $\bz \in \T^d$ have the property that for every $k \in \{0, 1,\ldots, l\}$, the set $\{z_j\}_{k(j) = k}$ is equidistributed in the uniform measure~$\mv$ of $\T$. Mirroring~\S~\ref{seeding}, we thus define a function $\Phi_{\boldsymbol{\varphi},\boldsymbol{\gamma}}$ on $\T$ by the piecewise splicing rule\footnote{Note that this function is different from the multivariable $\Phi$ in \S~\ref{seeding}.} 
\[\Phi(e^{2 \pi i t})_{\boldsymbol{\varphi},\boldsymbol{\gamma}}:= \varphi_k \left( e^{2 \pi i \frac{mt - \gamma_k}{\gamma_{k+1} - \gamma_k}}\right),
\ \text{for} \  t\in [\gamma_k/m, \gamma_{k+1}/m).\]
In~\eqref{rearrangement fine}, we have $g_{\boldsymbol{\varphi},\boldsymbol{\gamma}}(t) = \log{|\Phi_{\boldsymbol{\varphi},\boldsymbol{\gamma}}(e^{2\pi i t})|}$.  
Then, as in~\S~\ref{sec:num int}, we have
\begin{equation*}
\begin{aligned}
\lim_{\epsilon \rightarrow 0} \lim_{\epsilon', \epsilon'' \rightarrow 0} \lim_{d\rightarrow \infty} 
\left\{ \frac{1}{d}\int_{\T^d} \max_{ \substack{ \bt \in P_\epsilon^d,  \\ \bn \in P_d^{\epsilon''}\left((m+\epsilon')D\right) }} \left\{ \sum_{j=1}^d t_j \log |\varphi_{k(j)}(z_j)|  \right\} \,  \mv \right\}\\
 = \int_0^1 t \cdot g_{\boldsymbol{\varphi},\boldsymbol{\gamma}}^*(t)  \, dt .
\end{aligned}
\end{equation*}

The argument for the non-archimedean estimate is the same as \S~\ref{sec_BCfinite}, and we recover the thesis of Theorem~\ref{main:elementary form}: 
\begin{equation}\label{rearrangment-convexity}
m \leq \frac{\int_0^1 2t \cdot g_{\boldsymbol{\varphi},\boldsymbol{\gamma}}^*(t)  \, dt   + \frac{1}{m} \sum_{k=1}^l \left\{ \gamma_k^2 \log{\frac{|\varphi_k'(0)|}{|\varphi'_{k-1}(0)|}}  \right\} }{\log |\varphi_l'(0)|-(\ovtau(\bb)+ \tau^\sharp(\be)) }.
\end{equation}

Let us for concreteness now specialize to the setup $\varphi_k(z) := \varphi(r_k z)$ of~\S~\ref{sec_BCconvexity}; the argument there can equally be adapted to the general situation.
If now we select our division parameters $\boldsymbol{\gamma}$ to be the slopes $\gamma_k := \beta_k(s_h^*)$ in Theorem~\ref{main:BC fullconv} (assume the linear algebra condition of that theorem to be satisfied), then by \S~\ref{rmk_Nazarov}, we have
\begin{equation*}
\begin{aligned}
 \int_0^1 2t \cdot g_{\boldsymbol{\varphi},\boldsymbol{\gamma}}^*(t)  \, dt   & +  \frac{1}{m}\sum_{k=1}^l \gamma_k^2 \log{\frac{r_k}{r_{k-1}}} \\ 
& = \int_0^1 2t (\log|\Phi|)^*(e^{2\pi i t}) \, dt -  \frac{1}{m}\sum_{k=0}^l (\gamma_{k+1}^2 - \gamma_k^2) \log r_k 
\\ 
& \geq  \ovcL' \cdot \ovcL' -  \frac{1}{m}\sum_{k=0}^l (\beta_{k+1}^2 - \beta_k^2) \log r_k \\
& =  \ovcL' \cdot \ovcL'  + \frac{1}{m} \sum_{k=1}^l \beta_k^2 (\log r_k - \log r_{k-1})\\
&=  \ovcL' \cdot \ovcL' + \frac{1}{m} \sum_{k=1}^l \beta_k (\ovcL_{r_k} \cdot \ovcL' - \ovcL_{ r_{k-1}} \cdot \ovcL')\\
&= \ovcL' \cdot \ovcL' +  \ovcL' \cdot \ovcL_1 -  \frac{1}{m}\sum_{k=0}^l (\beta_{k+1} -\beta_k ) \ovcL_{r_k} \cdot \ovcL'\\
&= \ovcL' \cdot \ovcL_1,
\end{aligned}
\end{equation*}
which is the numerator of the bound that we obtained in \S~\ref{highdimBC}. In practice, \S~\ref{examplecompare} suggests these two bounds to be pretty close. In particular, numerically speaking, we do not need $\gamma_k$ or $\beta_k$ to be exact the heuristically optimal choice. The bound \eqref{rearrangment-convexity} holds for any choice of $\boldsymbol{\gamma} = \{\gamma_k\}$.

\section{The relationship between~\texorpdfstring{$Y(2)$}{Y(2)} 
and~\texorpdfstring{$Y_0(2)$}{Y0(2)}}
\label{sec:YtoY0(2)}
There is a natural identification of~$Y(2)$ with~$\mathbf{P}^1 \setminus \{0,1,\infty\}$ given by the coordinate~$\lambda$ with equation~(\ref{lambdadef}). If we let~$Y_0(2)$ denote the modular
curve of level~$\Gamma_0(2)$, then~$Y_0(2)$ is also rational with a hauptmodul
\begin{equation}
\label{defofh}
h:=  \lambda + \frac{\lambda}{\lambda - 1} = - 256 q \prod_{n=1}^{\infty} (1 + q^{n})^{24} = -256 \cdot \frac{\Delta(2 \tau)}{\Delta(\tau)}, \quad q = e^{2 \pi i \tau},
\end{equation}
where this time we write~$q =e^{2 \pi i \tau}$ in comparison
to~$q = e^{\pi i \tau}$ in  equation~(\ref{lambdadef}).
 Just as with~$\lambda$, we also view~$h$ by abuse of notation as a function of the parameter~$q \in \D$. 
The parameter~$h$ gives an identification
of~$Y_0(2)$ with~$\mathbf{P}^1 \setminus \{0,\infty\}$. The map~$Y(2) \rightarrow Y_0(2)$
of modular curves is smooth as a map of algebraic stacks, but not of the underlying
coarse moduli spaces. To properly account for this, it is better to remember that~$Y_0(2)$ has an elliptic point of order~$2$ at~$h = 4$,  which is the branch point of the double covering $\lambda \mapsto h =\lambda^2/(\lambda-1)$, with~$\lambda = 2$ for its unique preimage: the ramification divisor of the 
branched covering $Y(2) \to Y_0(2)$ as algebraic curves.  
On the stacks level, $Y(2) \rightarrow Y_0(2)$ is an \'etale map 
which is a Galois covering of degree~$2$, and so there is a natural relation between invariant functions on~$Y(2)$ and functions on~$Y_0(2)$.

However, as we shall see below, this relationship also respects some arithmetic properties of the corresponding power series expansions.
First we remark that
the transformation
\begin{equation}
\label{involutionfirst}
w: x \mapsto \frac{x}{x-1} \in \Z \llbracket x \rrbracket 
\end{equation}
is an involution of~$\mathbf{P}^1 \setminus \{0,1,\infty\}$ that preserves~$0$
and swaps~$1$ and~$\infty$. The integrality properties
of this map (and its inverse) means that if~$f(x) \in \Q \llbracket x \rrbracket$ has a certain denominator type then so does~$f(w(x))$. 
Second,
we note that the map~$w$ of equation~(\ref{involutionfirst}) is precisely the non-trivial Galois automorphism of (the function field of)~$Y(2)$ over~$Y_0(2)$.

\begin{lemma} \label{etalecover} 
 Let~$S \subset \mathbf{P}^1 \setminus \{0,1,\infty\}$ be a finite set invariant under the involution~$w$ of equation~\eqref{involutionfirst}, and define~$T \subset \P^1 \setminus \{0,\infty\}$ to be the image of~$S$ under the map $x \mapsto y$, where
\begin{equation} \label{y in x}
y := x + w(x) = x + \frac{x}{x-1} = \frac{x^2}{x-1}. 
\end{equation}
Consider $c_1, \ldots, c_r \in [0,\infty)$, and let~$f(x) \in \Q \llbracket x \rrbracket$ be a 
power series of the form
\begin{equation}
\label{nice}
f(x) = \sum_{n=0}^{\infty} a_n \, \frac{x^n}{\prod_{i=1}^{r}[1,\ldots,c_in]} \in \Q\llbracket x \rrbracket, \qquad a_n \in \Z \quad \forall n \in \NwithzeroA,
\end{equation}
which converges on a neighborhood of~$x=0$ and continues analytically as a holomorphic function along all paths
in~$\mathbf{P}^1 \setminus \{0,1,S,\infty\}$. 
Then:
\begin{enumerate}
\item \label{itemnice} The function~$f \left(w(x) \right) \in \Q \llbracket x \rrbracket$ is also of the form~(\ref{nice}).
\item \label{itemnicer} If~$f(x) = f(w(x))$, 
then we may use~\eqref{y in x} to formally write~$f(x)$ as a power series~$f(x) = F(y) \in \Q \llbracket y \rrbracket$ that satisfies
\begin{equation}
\label{nicer}
2 F(y) = \sum_{n=0}^{\infty} b_n \, \frac{y^n}{{\prod_{i=1}^{r}[1,\ldots,2c_in]}} \in \Q\llbracket y \rrbracket, \qquad b_n \in \Z \quad \forall n \in \NwithzeroA,
\end{equation}
converges in a neighborhood of~$y=0$, continues analytically as a holomorphic function along all paths in~$\mathbf{P}^1 \setminus \{0,4,\infty,T\}$,
and has finite local monodromy\footnote{In this generality, by a ``finite local monodromy of order dividing~$2$'' we simply mean that if~$\gamma : (0, 1] \to \P^1 \setminus \{0,4,\infty,T\}$ is any path with origin $\lim_{t \to 0^+} \gamma(t) =0$, and $\pi$ is a simple loop around~$y=4$ in~$\P^1 \setminus \{0,4,\infty,T\}$ based
at the endpoint~$\gamma(1)$, then the analytic continuations of~$f(y)$ at the ends of the concatenated paths~$\gamma$ and~$\pi^2 \cdot\gamma$ are equal.
This, of course, agrees with the usual notion in the special (finite-dimensional) case of a local system on~$\mathbf{P}^1 \setminus \{0,4,\infty,T\}$.} 
of order dividing~$2$ around the point~$y=4$.  
\end{enumerate}
Conversely, if~$2F(y)$ has the form~(\ref{nicer}), then~$2 \displaystyle{F \left(x + \frac{x}{x-1}\right)}$ has the form~(\ref{nice}), and if~$F(y)$ 
has the analyticity properties on~$\mathbf{P}^1 \setminus \{0,4,\infty,T\}$ spelled out in~\eqref{itemnicer}, then~$\displaystyle{F \left(x + \frac{x}{x-1}\right)}$ has
the analytic continuation property on~$\mathbf{P}^1 \setminus \{0,1,S,\infty\}$. 
\end{lemma}

\begin{proof} The function-theoretic claims follow directly
from Galois theory and the fact that~$x \mapsto w(x)$
is the automorphism of~$Y(2)$ over~$Y_0(2)$.
Hence it suffices to establish the claims concerning integrality. Property~\eqref{itemnice} is clear from the integral
coefficients in the expansion $w(x)^n = (-1)^nx^n(1-x)^{-n} \in x^n  \Z\llbracket x \rrbracket$ together
with the remark that the denominator type~$\prod_{i=1}^{r}[1,\ldots,c_in]$ is nested by division under~$n \mapsto n+1$. 
 Let~$x$ and~$y$ be related by the identity~\eqref{y in x}.
Since both the elementary symmetric functions in~$x$ and~$w(x) = x/(x-1)$ are equal to~$y = x+w(x) = xw(x)$, we may define polynomials~$P_n(y)$ by the rule
\begin{equation}  \label{polysym}
P_n(y) :=  x^n + (x/(x-1))^n.
\end{equation}
Then~$P_0(y)=2$, $P_1(y) = y$, and there is the 
elementary recurrence 
$$P_{n}(y) =   y P_{n-1}(y) -  y P_{n-2}(y).$$
We find  that~$P_n(y)$ has degree~$n$ and vanishes at~$y=0$ to order~$\lceil n/2 \rceil$.
Let us suppose now that we have a function~$f(x) = \sum A_n x^n$ 
whose coefficients~$A_n$ are rational numbers with~$a_n := A_n \prod_{i=1}^r [1,\ldots,c_in] \in \Z$. Then, in property~\eqref{itemnicer} under proof, we exploit the assumption $f(w(x)) = f(x)$ to write
$$f(x) + f\left( \frac{x}{x-1} \right)= 2F(y) =  \sum A_k P_k(y) 
=: \sum B_n y^n.$$
The middle equality defines a legitimate~$\Q\llbracket y \rrbracket$ series since~$P_k(y)$ is divisible by~$y^{\lceil k/2\rceil}$, and it can be
taken as a definition. 
To be more precise, all the nonzero coefficients of the polynomial~$P_k(y)$ occur in the degree range~$[k/2,k]$, and they are integers. Thus~$P_k(y)$ contributes 
to the~$y^n$ term only for~$k \in [n,2n]$, and $b_n := B_n  \prod_{i=1}^r [1,2,\ldots,2c_i n] \in \Z$.
Conversely, since
$\displaystyle{x + \frac{x}{x-1} = \frac{x^2}{x-1}},$ if we write
$$\sum A_n x^n = \sum B_k \left(x + \frac{x}{x-1}\right)^k = \sum B_k \frac{x^{2k}}{(x-1)^k},$$
then the terms on the right-hand side contributing to~$A_n$ occur only for~$k \le n/2$.
\end{proof}

Motivated by this lemma, we have the following:
\begin{df} Let~$F(x) \in \Q \llbracket x \rrbracket$.
We define the plus and minus \emph{symmetrization} functions~$F^{+}(y)$
and~$F^{-}(y)$ to be the elements of~$\Q \llbracket y \rrbracket$ such that
\begin{equation} \label{plusminus}
\begin{aligned}
F^{+}(y) = & \  F(x) + F \left(\frac{x}{x-1}\right), \\
F^{-}(y) =  & \ \left(x - \frac{x}{x-1} \right) \left(F(x) - F \left(\frac{x}{x-1}\right) \right), \end{aligned}
\end{equation}
where~$y := x + \displaystyle{\frac{x}{x-1}} = \frac{x^2}{x-1}$.  \endofremark
\end{df}

We connect these symmetrizations to the analytic resolvents~$\varphi^*f \in \mathcal{O}(\D)$ in the context of the arithmetic holonomy bounds.
We firstly introduce an {\it ad hoc} definition (which will only be used
in Lemma~\ref{analytic quotient}):

\begin{df}  A \emph{holonomic descent datum} 
is a tuple $$\cR_f = \left( U_{Y(2)}, \Sigma_{Y(2)}^0, \Sigma_{Y(2)}^1,  f \right)$$ consisting of: 
\begin{enumerate}
\item A contractible open neighborhood~$0 \in U_{Y(2)} \subset \C \setminus \{2\}$ 
which is
invariant under the involution~$w$.
\item  Finite subsets $\Sigma_{Y(2)}^0 \subset U_{Y(2)}$ and $\Sigma_{Y(2)}^1 \subset Y(2) = \C \setminus \{0,1\}$, both  invariant under the involution~$w$.
\item A holomorphic function~$f \in \mathcal{O}\left(U_{Y(2)} \right)$  
which is $w$-invariant ($f(w(x)) = f(x)$) and analytically continuable as a holomorphic function along all paths in~$\P^1 \setminus \{0, 1, \Sigma_{Y(2)}^0, \Sigma_{Y(2)}^1, \infty\}$. 
\end{enumerate}
To every such datum~$\cR_f$, we attach a \emph{quotient datum} $$\cQ_F =  \left( U_{Y_0(2)}, \Sigma_{Y_0(2)}^0, \Sigma_{Y_0(2)}^1,  F \right)$$ as follows.
By expressing the $w$-invariant power series expansion~$f(x) \in \C \llbracket x \rrbracket$ formally into $y := x+w(x)$, we have attached as\footnote{Except that, now in the analytic context, we can assume~$f \in \C\llbracket x \rrbracket$ rather than~$f \in \Q \llbracket x \rrbracket$; the proof, of course, is the same. } in Lemma~\ref{etalecover} a unique formal power series~$F(y) = F_f(y) \in \C \llbracket y \rrbracket$ 
 such that
 \begin{equation} \label{symf}
 F(y) = F(x+w(x)) = F\left(x + \frac{x}{x-1}\right) = f(x).
 \end{equation}
In a similar manner, we
define~$\Sigma_{Y_0(2)}^0$, $\Sigma_{Y_0(2)}^1$, and~$U_{Y_0(2)}$ to be the images of~$\Sigma_{Y(2)}^0$, $\Sigma_{Y(2)}^1$, and~$U_{Y(2)}$, under the map 
$$
x \mapsto y := x + w(x) = x + \frac{x}{x-1}
$$
of~\eqref{y in x}.   \endofremark
\end{df}

For brevity we also adopt the following definition, which formalizes the idea of the univalent leaves from Proposition~\ref{overconvergence}. 
\begin{df}  \label{def:univalentleaf}
Consider two pointed Riemann surfaces~$(D,O)$ and~$(X,P)$ and an open neighborhood~$P \in U \subset X$.
A holomorphic mapping~$\varphi : (D,O) \to (X,P)$ 
 \emph{has a univalent leaf over~$U$ at~$O$} if~$\varphi$ maps the 
connected component of~$\varphi^{-1}(U)$ containing~$O$ 
conformally isomorphically onto~$U$: 
$$
\varphi \, : \, \left(  \varphi^{-1}(U) \right)_O \iso U. 
$$
We refer to~$\left(  \varphi^{-1}(U) \right)_O \subset D$ itself as the univalent leaf (at~$O$ over~$U$).  \endofremark
\end{df}

\begin{lemma}  \label{analytic quotient}
Let~$\cR_f$ be  a 
holonomic descent datum with quotient~$\cQ_F$. 
Let~$\varphi_{Y(2)}$ be a holomorphic mapping that obeys
\[\varphi_{Y(2)} : \Db \to \C \setminus \{1, \Sigma_{Y(2)}^1\}, \quad \varphi_{Y(2)}^{-1}(0) = \{0\},\]
and which has a univalent leaf~$\left( \varphi_{Y(2)}^{-1}\left( U_{Y(2)} \right) \right)_0$ over~$U_{Y(2)}$ at~$0 \in \D$ containing all the pre-images of~$\Sigma_{Y(2)}^0$ under~$\varphi_{Y(2)}$.  

Suppose that 
 \begin{equation} \label{phi_invariance} 
 w(\varphi_{Y(2)}(z)) = \varphi_{Y(2)}(-z). 
 \end{equation}
The pullback~$\varphi_{Y(2)}^*f$ is holomorphic on~$\Db$.
If~$\varphi_{Y_0(2)}$ is the holomorphic map
 \begin{equation} \label{passage quot}
 \varphi_{Y_0(2)}(z) :=  \varphi_{Y(2)}\left(\sqrt{z}\right) + \varphi_{Y(2)}\left(-\sqrt{z}\right)
\in \mathcal{O}(\Db),
\end{equation}
then:
\begin{enumerate}
\item \label{0 prepped Y02} $\varphi_{Y_0(2)}^{-1}(0) = \{0\}$.
\item \label{omitted set Y02} The range of~$\varphi_{Y_0(2)}$ omits~$\Sigma_{Y_0(2)}^1$.
\item \label{ramification Y02} The ramification indices of~$\varphi_{Y_0(2)}$ are even at all points of the fiber~$\varphi_{Y_0(2)}^{-1}(4)$.
\item \label{Y02 leaf} The neighborhood $U_{Y_0(2)} \ni 0$ is a contractible domain, and~$\varphi_{Y_0(2)}$ has a univalent leaf over~$U_{Y_0(2)}$
at~$0 \in \D$, which furthermore contains all the pre-images of~$\Sigma_{Y_0(2)}$ 
under~$\varphi_{Y_0(2)}$. 
\item
$F|_{U_{Y_0(2)}} \in \mathcal{O}\left(U_{Y_0(2)}\right)$ is holomorphic,
and the following relation holds:
\begin{equation} \label{descent an}
\begin{aligned}
F\left( \varphi_{Y_0(2)}(z) \right)  & = f\left( \varphi_{Y(2)}(\sqrt{z}) \right) =  f\left(w\left( \varphi_{Y(2)}(\sqrt{z})  \right) \right)  \\ 
& = \frac{f\left( \varphi_{Y(2)}(\sqrt{z}) \right) +  f\left(\left( \varphi_{Y(2)}(-\sqrt{z})  \right) \right) }{2} \in \mathcal{O}(\Db).
\end{aligned}
\end{equation}
In particular, the pullback of~$F$ by~$\varphi_{Y_0(2)}$ is holomorphic on~$\Db$.
\end{enumerate}

Conversely, every holomorphic mapping~$\varphi_{Y_0(2)} \in \mathcal{O}(\Db)$ obeying the conditions~\eqref{0 prepped Y02} 
through~\eqref{Y02 leaf}
 determines through~\eqref{passage quot} a unique pair~$\left\{ \varphi_{Y(2)}, w \circ \varphi_{Y(2)}  \right\}$ of holomorphic mappings
 \[\varphi : \Db \to \C \setminus \{1, \Sigma_{Y(2)}^1 \}, \quad \varphi^{-1}(0) = \{0\} \]
  subject 
  to~$w\left( \varphi(z) \right) = \varphi(-z)$.
\end{lemma}

\begin{proof}
The holomorphy of~$\varphi_{Y(2)}^*f$ on~$\Db$
 follows directly from Proposition~\ref{overconvergence} with~$\Omega := \left( \varphi_{Y(2)}^{-1}(U_{Y(2)}) \right)_0$, as~$f \in \mathcal{O}\left( U_{Y(2)} \right)$. 
 We observe that~$U_{Y_0(2)}$ 
 \silentcomment{One way to see the claim is to take $u=x/(2-x)$, then on the $u$-plane, the $w$ map is nothing but $u \mapsto -u$ and then we see anything stable under $-1$ around $0$, we can take half of it (say $\Re \ge 0$) part such that the interior is biholomorphic to the image; the boundary corresponding to y-axis glue together under $u \mapsto u^2$ map and thus the image of  domain is simply connected.}
  --- a domain, by the open mapping theorem --- is also a topological disc, as~$U_{Y(2)} \subset \C \setminus \{2\}$ while the map~$y := x^2/(x-1)$ has~$x \in \{0, 2\}$ for its only ramification points, with branching values~$y \in \{0, 4\}$.

Assume now the symmetries~$f(x) = f(w(x))$ and~$w(\varphi_{Y(2)}(z)) = \varphi_{Y(2)}(-z)$, and define the manifestly holomorphic map~$\varphi_{Y_0(2)} \in \mathcal{O}(\Db)$ by~\eqref{passage quot}. It is the $w(\varphi_{Y(2)}(z)) = \varphi_{Y(2)}(-z)$ symmetry that allows to descend the 
analytic data to the~$Y_0(2)$ picture, as the plus-symmetrization of~$\varphi_{Y(2)}$: 
\begin{equation}  \label{an sym}
\begin{aligned}
\varphi_{Y_0(2)}(z) := & \ \varphi_{Y(2)}\left(\sqrt{z}\right) + \varphi_{Y(2)}\left(-\sqrt{z}\right)
\in \mathcal{O}(\Db) \\
= & \ \varphi_{Y(2)}\left(\sqrt{z}\right) + w \left(\varphi_{Y(2)}\left(\sqrt{z}\right) \right) \\
= & \ \frac{\varphi_{Y(2)}(\sqrt{z})^2}{\varphi_{Y(2)}(\sqrt{z})-1}.
\end{aligned}
\end{equation}
The second line --- together with the definitional fact that~$\Sigma_{Y(2)}^1$ is the full inverse image of~$\Sigma_{Y_0(2)}^1$ under the double covering map~$y = x+w(x)$ --- shows that the range of~$\varphi_{Y_0(2)}$ omits the set~$\Sigma_{Y_0(2)}^1$, as~$\varphi_{Y(2)}$ omits the corresponding set~$\Sigma_{Y(2)}^1$. The third line shows that~$\varphi_{Y_0(2)}$ satisfies~$\varphi^{-1}(0) = \{0\}$ and has even ramification indices at all points in the fiber~$\varphi_{Y_0(2)}^{-1}(4) = \varphi_{Y(2)}^{-1}(2)$.  Applying~$F(x+w(x)) = f(x) = f(w(x))$ for the~$x = \varphi_{Y(2)}(\sqrt{z})$ given on the second line in~\eqref{an sym}, we get~\eqref{descent an}, and in particular the holomorphy of~$\varphi_{Y_0(2)}^*F \in \mathcal{O}(\Db)$. Lastly, the holomorphy~$F|_{U_{Y_0(2)}} \in \mathcal{O}\left( U_{Y_0(2)} \right)$ follows directly from the corresponding holomorphy~$f|_{U_{Y(2)}} \in \mathcal{O}\left( U_{Y(2)} \right)$ thanks to the defining equation~$F(x+w(x)) = f(x)$ and the definition of~$U_{Y_0(2)}$ as the~$x+w(x)$ image of~$U_{Y(2)}$. 

For the converse,
 we get by the formal binomial expansion --- choosing any branch for the square root signs --- a power series~$\varphi_{Y(2)} \in \C\llbracket z \rrbracket$
from resolving the quadratic relation on the third line in~\eqref{an sym} with~$z$ changed to~$z^2$: 
\begin{equation} \label{pulling back}
\varphi_{Y(2)}(z) := \frac{ \varphi_{Y_0(2)}(z^2) + \sqrt{\varphi_{Y_0(2)}(z^2)} \cdot \sqrt{\varphi_{Y_0(2)}(z^2) - 4}  }{2}.
\end{equation}
The conditions on~$\varphi_{Y_0(2)}^{-1}(0) = \{0\}$ and on even ramification indices for~$\varphi_{Y_0}(2)$ along~$\varphi_{Y_0(2)}^{-1}(4)$  show that the formal function~\eqref{pulling back} is in fact
holomorphic on a neighborhood of~$\Db$, 
and satisfies~$\varphi_{Y(2)}^{-1}(0) = \{0\}$ and~$w\left( \varphi_{Y(2)}(z) \right) = \varphi_{Y(2)}(-z)$. The other choice of the square roots sign in~\eqref{pulling back} leads to the argument sign swap~$\varphi_{Y(2)}(-z)$, and the pair~$\left\{ \varphi_{Y(2)}, w \circ \varphi_{Y(2)}  \right\}$ is uniquely determined from~$\varphi_{Y_0(2)}$ and satisfies~\eqref{an sym}, whence the range property~$\varphi_{Y(2)} : \Db \to \C \setminus \{ 1, \Sigma_{Y(2)}^1 \}$ is also inherited. 
\end{proof} 

We spell out as a separate corollary the case that we will use of analytic pullbacks of the hauptmodul map~\eqref{defofh}. This should be regarded as a stacky version for~$Y_0(2)$
of Proposition~\ref{overconvergence} on overconvergence. 

\begin{cor}  \label{stacky overconvergence}
Consider an arbitrary power series~$F \in \C\llbracket y \rrbracket$ that defines a holomorphic function on a contractible open neighborhood~$0 \in U_{Y_0(2)} \subset \C \setminus \{4\}$. Suppose~$\Sigma_{Y_0(2)}^0 \subset U_{Y_0(2)}$ and~$\Sigma_{Y_0(2)}^1 \subset \C$ are finite subsets such that~$F(y)$ continues analytically as a holomorphic function along all paths in~$y \in \P^1 \setminus \{0, 4, \Sigma_{Y_0(2)}^0, \Sigma_{Y_{0(2)}}^1 \}$ and has around~$y=4$ a finite local monodromy of order dividing~$2$. Let~$h : \D \to \C$ be the map~\eqref{defofh}.

Then, 
 under any holomorphic mapping~$\varphi_{Y_0(2)} : \Db \to \C \setminus \Sigma_{Y_0(2)}^1$ that	has a univalent leaf over~$U_{Y_0(2)}$
 at~$0 \in \D$ containing~$\varphi^{-1}\left( \Sigma_{Y_0(2)}^0 \right)$, and which
  factors as a composition~$\varphi_{Y_0(2)} = h \circ \psi_{Y_0(2)}$ for some holomorphic~$\psi_{Y_0(2)} : \Db \to \D$ with~$\psi_{Y_0(2)}^{-1}(0) = \{0\}$, 
the pullback of~$F$ is holomorphic: 
 $\varphi_{Y_0(2)}^*F \in \mathcal{O}(\Db)$. 
\end{cor}

\begin{proof}
Define~$U_{Y(2)} := y^{-1}\left(U_{Y_0(2)}\right)$ as the full inverse image under the map~$y := x+w(x) = x^2/(x-1)$. Since~$4 \notin U_{Y_0(2)}$, 
 this neighborhood~$U_{Y(2)} \ni 0$ is also contractible. Setting also~$f(x) := F(y) = F(x+w(x))$ and~$\Sigma_{Y(2)}^0 := y^{-1} \left( \Sigma_{Y(2)}^0 \right)$, $\Sigma_{Y(2)}^1 := y^{-1} \left( \Sigma_{Y(2)}^1 \right)$, we have thus constructed a holonomic descent datum~$\cR_f$ with quotient~$\cQ_F = \left( U_{Y_0(2)}, \Sigma_{Y_0(2)}^0,  \Sigma_{Y_0(2)}^1, F \right)$. 
 
 Since (with our assumptions on~$\psi_{Y_0(2)}$) the maps of the form~$\varphi_{Y_0(2)} = h \circ \psi_{\psi_{Y_0(2)}}$ satisfy the conditions~\eqref{0 prepped Y02} through~\eqref{Y02 leaf} in Lemma~\ref{analytic quotient},
the converse direction of the lemma then constructs a holomorphic mapping~$\varphi_{Y(2)} : \Db \to \C \setminus \{1, \Sigma_{Y(2)}^1\}$ with~$\varphi_{Y(2)}^{-1}(0) = \{0\}$ and $w\left( \varphi_{Y(2)} (z) \right) = \varphi_{Y(2)}(-z)$, and inducing a conformal isomorphism~$\varphi_{Y(2)}^{-1}\left( U_{Y(2)} \right)_0 \iso U_{Y(2)}$: a univalent leaf over~$U_{Y(2)}$ at~$0 \in \D$. The forward direction of Lemma~\ref{analytic quotient} now proves the holomorphy~$\varphi_{Y(2)}^*f \in \mathcal{O}(\Db)$  together with the symmetrization relation~\eqref{descent an}, which in particular  manifests the holomorphy~$\varphi_{Y_0(2)}^* F \in \mathcal{O}(\Db)$.
\end{proof}

\begin{basicremark} \label{equivalence}
We combine and interpret Lemmas~\ref{etalecover} and~\ref{analytic quotient} as follows. For the remainder of our paper, we will consistently reserve the letter~$y$ to denote the covering~$y := x^2/(x-1)$. Suppose given a~$\Q(x)$-vector space~$\HH$ generated by $\Q\llbracket x \rrbracket$ power series of  the arithmetic type~\eqref{nice}, holomorphic on some neighborhood~$U_{Y(2)} \ni 0$, and analytically continuing as holomorphic functions along all paths in~$\mathbf{P}^1 \setminus \{0,1,S,\infty\}$. 
Lemma~\ref{etalecover} then constructs a corresponding~$\Q(y)$-vector
space~$\HH^{w = 1}$ over~$\Q(y)$ of functions on~$\mathbf{P}^1 \setminus \{0,4,T,\infty\}$, with at most $\Z/2$ local monodromy around~$y=4$, and satisfying the arithmetic condition~(\ref{nicer}). Moreover,
from basic Galois theory, we have
\begin{equation}
\label{galoisdescent}
\dim_{\Q(x)} \HH = \dim_{\Q(y)} \HH^{w = 1}.
\end{equation}
Explicitly, we have~$\dim_{\Q(y)}(\Q(x)) = 2$ with a basis given by~$1$ and~$y^{-}  =  x  -\displaystyle{\frac{x}{x-1}}$.
There is an isomorphism of~$\Q(y)$-vector spaces~$\HH = \HH^{w=1} \oplus \HH^{w=-1}$ given 
by
$$F(x) \mapsto \left(F^{+}(y), F^{-}(y)/y^-\right) = \left(F(x) +F\left( \frac{x}{x-1} \right), F(x) - F\left( \frac{x}{x-1} \right) \right),$$
and an isomorphism~$\HH^{w=1} \rightarrow \HH^{w=-1}$ given by multiplication by~$y^{-}$. 
Hence
$$\dim_{\Q(y)} \HH^{w = 1} = \frac{1}{2} \dim_{\Q(y)} \HH = \dim_{\Q(x)} \HH.$$
One can now ask what happens (for example) to a holonomy bound of the form:
\begin{equation} \label{BCagain}
\dim_{\Q(x)} \mathcal{H}  \le \frac{ \iint_{\T^2}  \log{|\varphi(z) - \varphi(w)  |} \, \mv(z) \mv(w)  }{\log{|\varphi'(0)|} - \sigma}, 
\end{equation}
when translated from the~$Y(2)$ or~$\Q(x)$ domain into the~$Y_0(2)$ or~$\Q(y)$ domain?

The answer to this question is that the corresponding bounds~\eqref{BCagain} are, like the 
dimensions~\eqref{galoisdescent} themselves, also equivalent in the framework of Theorem~\ref{basic main}. Firstly, we need to make precise what we 
mean by the $Y(2)$ versus the~$Y_0(2)$ domain in the context of formal-analytic arithmetic surfaces
and arithmetic holonomy bounds. This is the content and purpose of Lemma~\ref{analytic quotient}.
Coming from the setting of~\S~\ref{the lambda template}, ``using the~$Y(2)$ domain'' refers to the holomorphic mappings
$\varphi = \varphi_{Y(2)}: \Db \to \C \setminus \{1\} = Y(2) \cup \{0\}$ with~$\varphi^{-1}(0) = \{0\}$, and therefore factorizing as
$$\varphi_{Y(2)} = \lambda \circ \GGG_{Y(2)},$$
where~$\GGG_{Y(2)}: \Db \to \D$ is a holomorphic 
map still having~$\GGG_{Y(2)}^{-1}(0) = \{0\}$. 
The proof of this factorization (cf.~\cite[\S\S~4.11, 4.12]{Caratheodory} for the details) reduces to the fact that $\tau \mapsto \lambda(\tau)$, $\tau : \H \to  \P^1 \setminus \{0,1,\infty\}$ 
is a universal covering map at~$\tau(i) = 1/2$. On the other hand, the basic properties of the modular lambda map also include 
$\lambda(\tau+1) = \lambda(\tau)/(\lambda(\tau)-1)$ in the~$\tau \in \H$ domain, that is $\lambda(-q) = w(\lambda(q))$ in the~$q 
= e^{\pi i \tau} \in \D$ domain. Therefore, if we impose the condition $\GGG_{Y(2)}(-z) = - \GGG_{Y(2)}(z)$ on the map~$\GGG_{Y(2)}$, then
the involution~$w$ acts as~$w(\varphi_{Y(2)}(z)) = \varphi_{Y(2)}(-z)$. In the special case that~$\GGG_{Y(2)} : (\D,0) \iso (\Psi,0)$ is 
the Riemann map of a contractible domain~$\Psi$ with~$0 \in \Psi \subset \overline{\Psi} \subset \D$, this condition simply amounts
to asking for the domain~$\Psi$ to be symmetric across the origin. We further assume that the open neighborhood~$U_{Y(2)}$ of the origin meets the conditions of Lemma~\ref{analytic quotient}: namely, $w \left( U_{Y(2)} \right) = U_{Y(2)}$, and $\varphi_{Y(2)}^{-1}\left( U_{Y(2)}  \right)_0 \iso U_{Y(2)}$
is a univalent leaf of~$\varphi_{Y(2)}$ at~$0 \in \D$ containing all pre-images of $\Sigma^0_{Y(2)}$.  
We note that this property implies, but is stronger than, the corresponding property for the inner map~$\psi_{Y(2)}$ in the factorization
$\varphi_{Y(2)} = \lambda \circ \psi_{Y(2)}$. 

We are now in the realm of Lemma~\ref{analytic quotient}, where we may look for a decomposition~$S = \Sigma_{Y(2)}^0 \sqcup 
\Sigma_{Y(2)}^1$ with~$\Sigma_{Y(2)}^0 \subset U_{Y(2)}$ and
$\varphi_{Y(2)} : \Db \to \C \setminus \{ 1, \Sigma_{Y(2)}^1 \}$. Under these conditions, 
we obtain from~$\varphi_{Y(2)}$ a map~$\varphi_{Y_0(2)}$
such that~$\varphi_{Y(2)}^* f$ 
and~$\varphi_{Y_0(2)}^*F$ are
meromorphic on~$\Db$ for any~$f \in \HH$ or~$F \in \HH^{w=1}$ respectively.

It is with these choices~$\varphi_{Y(2)}$ and~$\varphi_{Y_0(2)}$ for the analytic mapping~$\varphi$, and with correspondingly the terms~$\sigma := \tau(\bb)$, resp.~$\sigma := \tau(2\bb) = 2\tau(\bb)$ under
the formulation of Theorem~\ref{basic main}, that we are comparing the holonomy quotients~\eqref{BCagain} under the dictionary supplied by Lemma~\ref{etalecover}. 

We now substantiate our claim that these two quotients~\eqref{BCagain} are exactly equal. 
The preceding analysis relies on the fact that the map~$Y(2) \rightarrow Y_0(2)$ of algebraic stacks is \'{e}tale. On the other hand, as the branched double covering of rational algebraic curves $X(2) \rightarrow X_0(2)$ is totally ramified over the center~$h=0$ (the cusp~$\tau = i\infty$) of our formal function expansions,  it follows by the projection formula in Lemma~\ref{BCconv-intersection}
that \emph{both} the corresponding integral and conformal radius terms on the~$Y_0(2)$ version of the holonomy quotient~\eqref{BCagain} are exactly scaled
by the degree of that 
 covering (which in our case is equal to two).
In our basic situation, we can see this in a very direct and explicit way as follows.
Let us write~$H(\tau)$ for the hauptmodul~$h$ evaluated at~$e^{2 \pi i \tau}$, and~$L(\tau)$ for~$\lambda$ evaluated at~$e^{\pi i \tau}$,
both with~$\tau$ in the upper half plane~$\H$.
If~$\GGG_{Y(2)}(e^{i \theta}) = e^{2 \pi i \tau}$ with~$\tau \in \H$, then, 
making some (consistent) choice of square roots, we have
$$\varphi_{Y_0(2)}(e^{i \theta}) = h(e^{2 \pi i \tau}) = H(\tau),$$
whereas
$$\varphi_{Y(2)}(e^{i \theta/2}) = \lambda(e^{\pi i \tau}) = L(\tau), 
\quad
\varphi_{Y(2)}(-e^{i \theta/2}) = L(\tau + 1) = \frac{L(\tau)}{L(\tau) - 1}.$$
In particular, the~$Y(2)$ integral involving~$\log |L(\tau) - L(\sigma)|$ becomes, after dividing the integral up into four pieces, an integral 
of 
$$\begin{aligned}
 \log  \left| L(\tau) - L(\sigma) \right|  & + \log \left|   \frac{L(\tau)}{L(\tau) -  1} - L(\sigma) \right| \\
+    \log \left|  L(\tau) - \frac{L(\sigma)}{L(\sigma) -  1}  \right| 
& +\log \left|  \frac{L(\tau)}{L(\tau) -  1}  - \frac{L(\sigma)}{L(\sigma) -  1}  \right|. \end{aligned}$$
But now (using only  the multiplicativity property of the logarithm and some elementary algebra) this is exactly
$$2 \log\left| L(\tau) +  \frac{L(\tau)}{L(\tau) -  1}  - L(\sigma) -  \frac{L(\sigma)}{L(\sigma) -  1} \right|
= 2 \log | H(\tau) - H(\sigma)|.$$
Taking into account the factors of~$2$ coming from the various scalings, 
this means that the integral in the~$Y_0(2)$ domain is \emph{precisely double} the integral in the~$Y(2)$ domain.
On the other hand, the conformal radius is also squared (this is clear for the factors of~$\GGG_{Y(2)}$
and~$\GGG_{Y'(2)}$ together with the equality~$|h'(0)| = 256 = |\lambda'(0)|^2$),
and so the logarithm of the conformal radius is also doubled.
At the same time,
in the context of Lemma~\ref{etalecover}, the invariant~$\sigma = \sum_{i=1}^r c_i$ is \emph{also} doubled, and so is 
the invariant~$\tau(\bb)$ in the context of Theorem~\ref{basic main}.

In summary,
the bound~(\ref{BCagain}) applied to ~$\dim_{\Q(x)} \HH$ and~$\dim_{\Q(y)} \HH^{w=1}$ (which are equal
by equation (\ref{galoisdescent}))
 gives the same result in both cases. This therefore gives a (rough) equivalence between these two problems
on both the arithmetic and  the analytic sides.
However, it  is also important to note is that this crisp equivalence of the bounds only applies to the framework of the crude denominator types as stated in Theorem~\ref{basic main} or Lemma~\ref{etalecover}, and that there \emph{is} still a difference once we start to consider refined denominators data, such as with  the~$\tau^{\sharp}$ from the added integrals in~\S~\ref{fine section}. We shall see in~\S~\ref{mpoly} that it can then be advantageous to perform the~$\varphi_{Y(2)} \rightsquigarrow \varphi_{Y_0(2)}$ analytic descent passage that we detailed in this Basic Remark. 
\endofremark
\end{basicremark}

\section{Pure functions on~
\texorpdfstring{$\mathbf{P}^1 \setminus \{0,1,\infty\}$}{P1m0,1,oo}
and on~
\texorpdfstring{$\mathbf{P}^1 \setminus \{0,4,\infty\}$}{P1m0,4,oo}
}
\label{sec:purefunctions}
The goal of this section is to write down a number of~$G$-functions
with nice integrality properties on~$\mathbf{P}^1 \setminus \{0,1,\infty\}$,
and then, using the translation discussion
in~\S~\ref{sec:YtoY0(2)}, on~$\mathbf{P}^1 \setminus \{0,4,\infty\}$ as well,
where the point~$y = 4$ should be though of as an elliptic point of order~$2$. In terms of 
local systems on an orbifold, the proper way to think of these domains is as the modular curves 
$Y(2)$ in the coordinate $x = \lambda$, and respectively, $Y_0(2)$ in
the coordinate~$y = x^2/(x-1) = \lambda^2/(\lambda-1) = h$. 

\subsection{Five functions of type~\texorpdfstring{$[1,\ldots,n]^2$}{[1,2,..,n]2} on~\texorpdfstring{$\mathbf{P}^1 \setminus \{0,1,\infty\}$}{P1m0,1,oo}} \label{sec:pureX2}
There are four obvious~$\Q(x)$-linearly independent~$G$-functions we can write down on~$\mathbf{P}^1 \setminus \{0,1,\infty\}$
with denominator type~$\tau = [1,2,3,\ldots,n]^2$. Namely:
\begin{equation}
\begin{aligned}
A_1(x)  =  & \ 1, \\
A_2(x) = -\log(1-x)  =  & \ \sum_{n=1}^{\infty} \frac{x^n}{n}, \\
A_3(x) = \log^2(1-x) = & \ \left(   \sum_{n=1}^{\infty} \frac{x^n}{n} \right)^2, \\
A_4(x) = \Li_2(x) = & \ \sum_{n=1}^{\infty} \frac{x^n}{n^2}.
\end{aligned}
\end{equation}
Clearly~$A_2(x)$ additionally has type~$\tau = [1,2,3,\ldots,n]$, and~$A_3(x)$ has 
denominator type $[1,2,\ldots,n][1,2,\ldots,n/2]$. 
These functions are linearly independent over~$\Q(x)$. 
Using symmetrizations, we also obtain~$4$ linearly independent functions over~$\Q(y)$. These can be given explicitly as follows: 
\begin{equation}
\label{defB}
\begin{aligned}
B_1(y) = & \ 1, \\
B_2(y) = & \ 
\sum_{n=2}^{\infty} 2 y^n \cdot \frac{(n-2)! n!}{(2n)!}
= 2y -  2\sqrt{y(4-y)} \arcsin\left(\frac{\sqrt{y}}{2}\right) \\ 
B_3(y) = & \
\sum_{n=1}^{\infty} y^n \cdot \frac{(n-1)!^2}{(2n)!}
= 2 \arcsin(\sqrt{y}/2)^2 \\
B_4(y) = & \  \mathrm{Sym}^{-} \Li_2(y)
= \left(x - \frac{x}{x-1}\right) \left(\Li_2(x) - \Li_2 \left(\frac{x}{x-1}\right)\right)
\\
= &  \  -2 \sqrt{y(4-y)} \int \frac{
\displaystyle{\arcsin\left(\frac{\sqrt{y}}{2}\right)}
}{y} \\
= & \ 4 \cdot \sum_{n=0}^{\infty} \frac{y^{n+1}}{16^n} 
\left(
\sum_{k=0}^{n}
\binom{2k}{k} \binom{2n-2k}{n-k}
\frac{1}{(2k-1)(2n-2k+1)^2} \right) \\
= & -4 y + \frac{4 y^2}{9} + \frac{31 y^3}{900} + \frac{389 y^4}{88200} + \ldots.  \end{aligned}
\end{equation}

These functions all have denominator type subsumed by~$[1,2,\ldots,2n]^2$
(for a more precise description, see Lemma~\ref{bdenominators} and Remark~\ref{central binomial}).

We spent a possibly embarrassing period of time believing that the four functions~$A_1(x), \ldots, A_4(x)$ spanned
the~$\Q(x)$-vector space of functions on~$\mathbf{P}^1 \setminus \{0,1,\infty\}$
with denominator type~$\tau = [1,2,\ldots,n]^2$. However, there is also a fifth function one can write down.
It arises more naturally in the~$\Q(y)$-domain, namely as
\begin{equation}
\label{defB5}
B_5(y) = \sum_{n=1}^{\infty} y^n \cdot \frac{(n-1)!^2}{(2n-1)! \cdot (2n-1)} =  y
\cdot  \pFq{3}{2}{1/2,1,1}{3/2,3/2}{\frac{y}{4}}.
\end{equation}
The function~$B_5(y)$ arises in Nesterenko's approximations~\cite{Nesterenko} to Catalan's constant~$G = L(2,\chi_{-4})$~\cite{Catalan}
in association with the equality
\begin{equation}
\label{Nielsen}
B_5(4)= 8 G
\end{equation}
due to Nielsen~\cite[page 166]{Nielsen}.
Here~$B_5(y)$ is of type~$[1,2,3,\ldots, 2n]^2$ (and even somewhat
better than this, see Lemma~\ref{bdenominators}).
One can easily define a corresponding function
\begin{equation}
\label{newfunction}
A_5(x) = J(x) = x \cdot  \pFq{3}{2}{1/2,1,1}{3/2,3/2}{
\frac{1}{4} \left(x + \frac{x}{x-1} \right)}
\end{equation}
on~$\mathbf{P}^1 \setminus \{0,1,\infty\}$ with denominator type~$\tau = [1,2,3,\ldots,n]^2$,
which we originally missed!
If
$$\LL = 2x(1-x)^2 \frac{d^2}{dx^2} + (2-x)(1-x) \frac{d}{dx} + 1,$$
then~$\LL J(x) = 2 - 2 x$.
We also see that
\begin{equation}
\label{odeJ}
2x(x-1) \frac{d J(x)}{dx} - x J(x) = 2 (1-x) \log(1-x).
\end{equation}

From the differential equation we see that~$J(x)$
is defined on~$\mathbf{P}^1 \setminus \{0,1,\infty\}$.
The solutions to the homogenous differential 
equation~$\LL(F)=0$ are given by
$$A(x) = \sqrt{1-x},$$
$$B(x) = \sqrt{1-x} \cdot \mathrm{arctanh}(\sqrt{1-x})
= \frac{1}{2} \cdot \sqrt{1-x} \cdot \log \left( \frac{1 + \sqrt{1-x}}
{1 - \sqrt{1-x}}\right).$$

Using the method of variation of parameters,
an explicit solution to the ODE is given by~$H(x)$ defined as follows:
$$2 \sqrt{1-x} \cdot\mathrm{arctanh}(\sqrt{1-x}) \log(-1 + x)
- 2 \sqrt{1-x} \left(- \Li_2 \left(- \sqrt{1-x} \right) + \Li_2 \left(\sqrt{1-x}\right)
\right),$$
and, having made suitable choices for the various analytic continuations of these terms, one can write
$$J(x) = H(x) - 2 \pi i B(x) + \frac{\pi^2}{2} A(x).$$
In retrospect,  an easier (if equivalent) way to write a new~$\Q(x)$-linearly
independent function (that gives the same span as~$J(x)$) while
remaining entirely on~$\mathbf{P}^1 \setminus \{0,1,\infty\}$ is to consider the integral:
\begin{equation}
\label{twomystery}
\frac{1}{\sqrt{1-x}} \int_0^x \frac{\log(1-t)}{t \sqrt{1-t}} \, dt.
\end{equation}
What is surprising in this formulation is the unexpected lack of extra powers of~$2$ in the denominators
of the Taylor series expansion of~(\ref{twomystery}). While the individual factors $1/\sqrt{1-x}$ and $ \int \frac{\log(1-x)}{x \sqrt{1-x}} dx$
have $2$-adic convergence discs $|x|_2 < 1/4$ at $x = 0$, their product overconverges to the full unit open $2$-adic disc $|x|_2 < 1$. 

\begin{remark} \label{mystery three}
One can prove that the $k \in \NwithoutzeroA$ for which the Taylor series of 
\begin{equation}
\frac{1}{\sqrt{1-x}} \int \frac{\log^{k-1}(1-x)}{x \sqrt{1-x}} dx
\end{equation}
converges on the $2$-adic unit disc $|x|_2 < 1$ are exactly the positive even integers, and that for these~$k$, the Taylor expansion belongs to 
$$
J_k(x) := \frac{1}{\sqrt{1-x}} \int \frac{\log^{k-1}(1-x)}{x \sqrt{1-x}} dx \in  \sum_{n=1}^{\infty} \frac{x^n}{[1,\ldots, n]^k} \, \Z. 
$$
This defines a sequence $J_2, J_4, J_6, \ldots$ (with~$J_2 = J$) of $G$-functions holonomic on $\mathbf{P}^1 \setminus \{0,1,\infty\}$, of denominator 
types $x^n/[1,\ldots,n]^{\ast}$, and independent over the multiple polylogarithm ring~\S~\ref{mpoly}.~\endofremark
\end{remark}

\subsection{Added integrations}
We define two more functions~$B_6(y)$ and~$B_7(y)$ as follows:

\begin{equation} 
\label{defB67}
\begin{aligned}
B_6(y) = & \ \int \frac{B_3(y)}{y} \, dy = \int \frac{2 \arcsin(\sqrt{y}/2)^2}{y}  \, dy = \sum_{n=1}^{\infty} y^n \cdot \frac{(n-1)!^2}{n(2n)!},  \\
B_7(y) = & \ \int \frac{B_4(y)}{y} \, dy. \\
\end{aligned}
\end{equation}

We have:
\begin{lemma} \label{bdenominators}The denominator types of~$B_i(y)$ for~$i = 1,\ldots 7$, as defined in equations~(\ref{defB}), (\ref{defB5}), and~(\ref{defB67})
are as follows:
\begin{enumerate}
\item $B_1(y)$ has trivial denominator type.
\item $B_2(y)$ has denominator type~$[1,2,,\ldots,2n]$.
\item $B_3(y)$ has denominator type~$[1,2,\ldots,2n]n$.
\item $B_4(y)$ has denominator type~$[1,2,,\ldots,2n]^2$.
\item $B_5(y)$ has denominator type~$[1,2,,\ldots,2n](2n-1)$, and thus in particular of denominator type~$[1,2,\ldots,2n]^2$. 
\item $B_6(y)$ has denominator type~$[1,2,\ldots,2n]n^2$, and thus in particular 
of denominator type~$[1,\ldots,n][1,\ldots,2n] n$, and \emph{a fortiori}~$[1,2,\ldots,2n]^2n$. 
\item $B_7(y)$ has denominator type~$[1,2,,\ldots,2n]^2 n$. 
\end{enumerate}
\end{lemma}

\begin{proof}
 This follows in the case of~$B_4(y)$  from Lemma~\ref{etalecover},
 and in the case of~$B_7(y)$ from direct integration from the~$n=4$ case.
 For the remainder, it follows by direct computation since there is an explicit
 expression in terms of factorials for the general coefficient.
\end{proof}

\begin{remark} \label{central binomial}
In fact the denominators of these functions have a somewhat better type, namely the
$[1,\ldots,2n]$ can be relaxed to $n(n-1) \binom{2n}{n}$. In practice this means that the prime product 
$\prod_{2n/3 < p  < n} p$ is absent from the $[1,\ldots,2n]$ part of these denominators. 
 This remark seems to not make any improvement in the setup for Theorem~\ref{mainA}, but the possibility
 of canceling prime products could be useful to exploit in other contexts.~\endofremark  
\end{remark}

We shall prove in \S~\ref{sec:lindep} that the seven functions~$B_i(y)$ are linearly independent over~$\Q(y)$.

\subsection{The multiple polylogarithm ring} \label{mpoly}
(This section  is more of an extended aside and can be omitted on first reading.) 
Over~\texorpdfstring{$\mathbf{P}^1 \setminus \{0,1,\infty\}$}{P1m0,1,oo}, a basic construction of
$G$-functions of type $[1,\ldots,n]^{\bullet}$ is supplied by the single variable multiple polylogarithm functions
\begin{equation} \label{polym}
\Li_{k_1,\ldots,k_d}(x) =  \sum_{n_1>n_2> \ldots> n_d}
\frac{x^{n_1}}{n^{k_1}_1 n^{k_2}_2 \cdots n^{k_d}_d}.
\end{equation}
Of these, the following eight functions form a maximal $\Q(x)$-linearly independent
set with type $n[1,\ldots,n]^2$:
\begin{equation} \label{mpoly2}
1, \quad \Li_1, \quad \Li_{1,1},  \quad  \Li_2, \quad  \Li_{1,1,1}, \quad \Li_1 \cdot \Li_2, \quad \Li_{1,2}, \quad \Li_3. 
\end{equation}
In Remark~\ref{mystery three}, we found that the multiple polylogarithms do not exhaust all the~\texorpdfstring{$\mathbf{P}^1 \setminus \{0,1,\infty\}$}{P1m0,1,oo} 
functions of the type $[1,\ldots,n]^{\bullet}$, and in particular, that we can add to~\eqref{mpoly2} a ninth independent function~\eqref{twomystery}
of the $[1,\ldots,n]^2$ type. By symmetrization, these nine $\Q(x)$-linearly independent functions go to nine $\Q(y)$-linearly independent functions on~\texorpdfstring{$\mathbf{P}^1 \setminus \{0,4,\infty\}$}{P1m0,4,oo} with $\Z/2$ local monodromies at the elliptic point $y =4$. 
However, whereas in Lemma~\ref{etalecover} we proved that the two symmetrization operations $F(x) \rightsquigarrow F^{\pm}(y)$ 
take the type $[1,\ldots,n]^{\sigma}$ to the type $[1,\ldots,2n]^{\sigma}$, an examination of the polynomials~\eqref{polysym} of the proof
reveals that the plus symmetrization $F^+$ takes the integrated type $n[1,\ldots,n]^{\sigma}$ to the 
integrated type $n[1,\ldots,2n]^{\sigma}$, but the minus symmetrization $F^-$ spoils the integrated type $n[1,\ldots,n]^{\sigma}$
into $[1,\ldots,2n]^{\sigma+1}$. This is why, as it turns out, the nine $\Q(x)$-independent functions of type $n[1,\ldots,n]^2$ on~\texorpdfstring{$\mathbf{P}^1 \setminus \{0,1,\infty\}$}{P1m0,1,oo} 
go to only seven $\Q(y)$-independent functions of type $n[1,\ldots,2n]^2$: the above $B_i(y)$. A similar remark shall
apply to~\S~\ref{jointli}, where under a supposed $\Q$-linear dependency among $1,\zeta(2), L(2,\chip)$ we would get
as many as~$17$  independent  functions over $x \in \P^1 \setminus \{0, 1/9, -1/8, 1, \infty\}$ with the 
integrated type  $n[1,\ldots,n]^2$ and
holomorphic at $\{0, 1/9, -1/8\}$ (only the above-listed nine of which really exist), but only $14$ of the symmetrizations (seven of them genuine) have the corresponding
integrated type  $n[1,\ldots,2n]^2$. 

One of the key ideas in our paper is that while  ---  as explained in Basic Remark~\ref{equivalence}  ---  our holonomy bounds are equivalent for the 
data  $(\varphi; \prod [1,\ldots, \mathbf{b}n] ) := \left( \lambda^2/(\lambda-1), [1,\ldots,2n]^2 \right)$ and $(\varphi; \prod [1,\ldots, \mathbf{b}n] ) := \left( \lambda, [1,\ldots,n]^2 \right)$, 
 the integrated type  
 \[(\varphi; \prod n^{\mathbf{e}} [1,\ldots, \mathbf{b}n] ) := \left( \lambda^2/(\lambda-1), n[1,\ldots,2n]^2 \right)\]
 yields significantly better bounds than the integrated type 
 \[(\varphi; \prod n^{\mathbf{e}} [1,\ldots, \mathbf{b}n] ) := \left( \lambda,  n[1,\ldots,n]^2 \right);\]
 so much so that the~$14$ functions in the former type turn out to be a far stronger constraint than the~$17$ functions in the latter type.
  We elaborate on this comparison in our next remark.

 \begin{remark} (This remark is best appreciated after reading the entire proof of Theorem~\ref{mainA}, although
 it still makes the most sense to place it in this section.)
 The above~$17$ functions of denominator type $n[1,\ldots,n]^2$ fit into the following refined
 denominators scheme in Theorem~\ref{main:elementary form}: 
 $$
\mathbf{b} := \left(  \begin{array}{lllllllllllllllll}   0 & 0 & 1 & 1 & 1 & 1 & 1 & 1 & 1 & 1 & 1 & 1 & 1 & 1 & 1 & 1 & 1 \\
0 & 0 & 0 & 0 & 1 & 1 & 1 & 1 & 1 & 1 & 1 & 1 & 1 & 1 & 1 & 1 & 1  \end{array} \right)^{\mathrm{t}}
$$
and the integrations vector  
$$
\mathbf{e} :=  (0, 1, 1, 1, 1, 1, 1, 1; 0;  0,  0, 0, 0, 1,1,1,1). 
$$
The first eight entries here are indexed by the row~\eqref{mpoly2}, in precisely this order, where, in view 
of the term $\max_i(e_i)$ in the definition~\eqref{taue} of $\tau^{\sharp}$, 
we opt to subsume $\Li_2$ into the type $n[1,\ldots,n]$ and $\Li_3$ into the type $n [1,\ldots, n]^2$. 
The ninth entry is the function~\eqref{twomystery}. Finally, writing $H(x) \in \Q \llbracket x \rrbracket$ for the 
function in Proposition~\ref{functionsH} below, the eight last entries are the fictive functions
\begin{equation*}
\begin{aligned}
& H(x), \, H'(x),  \, H\left( \frac{x}{x-1} \right), \, H'\left( \frac{x}{x-1} \right), \\
&  \int \frac{H(x) - H(0)}{x} \, dx, \,  \int \frac{H\left( \frac{x}{x-1} \right) - H(0)}{x} \, dx, \\
&   \int \frac{H(x) - H(0)}{x-1} \, dx, \,  \int \frac{H\left( \frac{x}{x-1} \right) - H(0)}{x-1} \, dx,
 \end{aligned}
 \end{equation*}
which turn out to be $\Q(x)$-linearly independent and holonomic on $\P^1 \setminus \{0, 1/9, -1/8, 1, \infty \}$. 

Recall that~$\tau(\mathbf{b;e}) = \tau^{\flat}(\mathbf{b}) + \tau^{\sharp}(\mathbf{e})$ is built out of two pieces.
For these denominator types, we calculate  
\begin{equation*}
\begin{aligned}
\tau^{\flat}(\mathbf{b}) & = \frac{ (1+3) \cdot 0 + (5+7) \cdot 1 + (9 + 11 + 13 +\ldots  + 33 ) \cdot 2 }{17^2}  \\
& = \frac{558}{289} = 1.93079584\ldots
\end{aligned}
\end{equation*}
which improves over the crude main denominator cap $\sigma = 2$.
The value  for~$\tau^{\flat}(\mathbf{b})$ we obtain here is even better than
the corresponding value~$191/49 = 2 \cdot 1.948979\ldots$ that we will use in~\S~\ref{sec:proofA}; for here
we can further exploit the special integrated type~$n$ of $\Li_1 = \int dx/(x-1)$. 
But for the other piece $\tau^{\sharp}(\mathbf{e})$ of $\tau(\mathbf{b;e})$ we get 
\[\tau^{\sharp}(\mathbf{e}) = 83711/242760 = 0.34483\ldots,\]
 with 
the optimal $\xi$ in~\eqref{taue} being a certain short interval containing the choice~$\xi = 57/40$. In total here,
\[\tau(\mathbf{b;e}) =558/289 +  83711/242760 = 552431/242760 =2.275626\ldots.\]
 This is \emph{very} much 
inferior to the value $16603/3920 = 2 \cdot 2.1173\ldots$ in~\eqref{taufine}, and 
the three additional functions are not nearly enough to compensate, as we now explain. 

To look into the numerics of the holonomy quotients, we can choose the map~$\varphi$ as the
optimal map of the form
$$
\varphi(z) := \lambda(G(z)), \qquad G :  (\Db,0) \to (\D,0), \qquad \varphi'(0) = 16 G'(0),
$$ 
where concretely $G$ can be (for example) the Riemann mapping for the topological disc inside $\D$ 
constrained by any simple closed contour that encircles the origin,
precisely like the contours we study in~\S~\ref{contour choiceA}. To be admissible, the contour must not enclose any 
of the non-real fiber points in $\lambda^{-1}(1/9)$ and $\lambda^{-1}(-1/8)$, but let us even ignore this point since
it will only make the numerics worse. Then the holonomy bound, which would have to compare to~$m = 17$, is by the quotient
\begin{equation} \label{X2quot}
\frac{ \iint_{\T^2} \log{|\varphi(z) - \varphi(w)|} \, \mv(z) \mv(w) }{\log{16} + \log{|G'(0)|} - \frac{552431}{242760}}.
\end{equation}
Using a (lightly) optimized choice~$\HHH(0.92, 110,23)$ from the gobble contours defined in~\S~\ref{gobble contours},
we find that $|G'(0)| = 0.9163768\ldots$ 
and the quotient~\eqref{X2quot} comes out to approximately~$22.7527$, a rather long distance  from
the requisite threshold of~$17$.~\endofremark
 \end{remark}

\subsubsection{Perspective on Theorems~\ref{mainA} and~\ref{logsmain}}
  This is why we shall henceforth stick with the type $n[1,\ldots,2n]^2$ functions $f_i(y)$ (such as the above $B_i(y)$), holonomic with singularities at $y  = \infty$, at $y = 4$ with a $\Z/2$ local monodromy, 
 and with all other singularities being overconvergent for the $f_i(y)$, and close enough to~$0$. In the application to Theorems~\ref{mainA} and~\ref{logsmain}, these latter ``overconvergent''
 singularities turn out to be $\{0, -1/72\}$, as we find out in the next section. Armed with Theorem~\ref{main:elementary form}, we will find in~\S~\ref{sec:proofA} that these 
 singularities $\{0,-1/72\}$ are indeed close enough to~$0$, and that a holonomy bound smaller than~$14$ can fortuitously be reached:  proving that such~$14$ independent functions cannot simultaneously exist. 
 Ultimately, this contradicts the supposed $\Q$-linear dependency among $1,\zeta(2)$, and $L(2,\chip)$, where as many as~$7$ of the $14$ functions arise from 
any such linear relationship via Lemma~\ref{14functions}. 
  
\section{Zagier's sequences {\bf A\rm} and {\bf C\rm}}  \label{Zagier local system}

\subsection{Definitions and basic properties}  \label{X06}
 
In this section, we construct a number of holonomic functions converging on the unit disc and extending to holomorphic
functions on the universal cover of~$\mathbf{P}^1 \setminus \{0,1/9,-1/8,1,\infty\}$, and also
on the universal cover of (the orbifold/stack)~$\mathbf{P}^1 \setminus \{0,-1/72,4,\infty\}$ where~$y=4$ is an elliptic point of order~$2$. 
Under the hypothesis that there is a~$\Q$-linear
relation between~$1$, $\zeta(2)$, and~$L(2,\chip)$, these functions would have rational coefficients and bounded denominator growth.
These constructions all come   ---   in a form very close to what is presented here   ---   from a paper of Zagier~\cite{Zagier}, but the sequences themselves
were certainly considered before then in similar contexts, including in particular in~\cite{Stienstra}, and the arguments
required to prove the required identities were first observed by Beukers~\cite{Beukers}. They arise more or less as solutions to Picard--Fuchs
equations associated to modular curves with precisely four cusps.
 The observation that certain linear combinations of solutions in~$\Q \llbracket x \rrbracket$ whose coefficients are interesting periods are overconvergent 
 (that is, extend analytically across the singular point of the ODE closest to~$x = 0$) was exploited by Beukers~\cite{Beukers} to give a reinterpretation
 of Ap\'{e}ry's original proof that~$\zeta(2)$ and~$\zeta(3)$ are irrational. As Zagier notes, however, the particular functions we consider
 (associated to the sequences~{\bf A\rm} and~{\bf C\rm} in the notation of~\cite{Zagier})   ---  
 while giving sequences which converge to both~$\zeta(2)$ and~$L(2,\chip)$   ---   
 ``do not converge quickly enough to yield the irrationality of the limit'' (\cite[p.~360]{Zagier}). To be precise, 
 these simultaneous approximations $u_n/q_n \to L(2,\chi_{-3}), 
 v_n/q_n \to \zeta(2)$ converge at the rate $q_n^{ - c }$ where $c = \log{9} \big/ \big( 2 + \log{9} \big) = 0.52349\ldots$, whereas,
 in the classical scheme for irrationality proofs, an exponent $c > 1$ would be required. 
    And yet, these functions are precisely  the required input  in our method to prove
 the desired irrationality results, by exploiting not simply the convergence properties of these functions on the unit disc $|x| < 1$, but also their analytic continuations
beyond the boundary point $x = 1$. 
 
The following is standard, and can also be read off from~\cite[Table~3, p.~357]{Zagier}:
 
\begin{lemma}
The function
\begin{equation}
\label{haupt}
x =  q \prod_{n=1}^{\infty} \frac{(1 - q^n)^4 (1 - q^{6n})^8}{(1 - q^{2n})^8 (1 - q^{3n})^4} = q - 4 q^2 + 10 q^3 + \ldots 
\end{equation}
with~$q = e^{2 \pi i \tau}$  defines a uniformization map
$$
x \, : \, Y_0(6) = \H/\Gamma_0(6) \rightarrow \mathbf{P}^1 \setminus \{0,1/9,1,\infty\}
$$
taking the $\Gamma_0(6)$ cusps $\tau = i\infty, 0, 1/3, 1/2$ to the respective cusps $x = 0, 1/9, 1, \infty$.
\end{lemma}

Note that one can formally invert this power series and write
\begin{equation}
\label{invert}
q = x + 4 x^2 + 22 x^3 + \ldots \in \Z \llbracket x \rrbracket.
\end{equation}
It follows that any power series in~$\Z \llbracket q \rrbracket$ can be written formally as a power
series in~$\Z \llbracket x \rrbracket$, and any power series in~$\Q \llbracket q \rrbracket$ can be written formally as a power
series in~$\Q \llbracket x \rrbracket$.

Let
$$\chi_{-3}(n) =  \left( \frac{-3}{n} \right) $$
be the unique primitive character of conductor~$3$. Consider the theta function of the Eisenstein
lattice $\Z[\zeta_3]$:
$$
\theta_{-3}(\tau)  := \sum_{m,n \in \Z} q^{m^2 + mn + n^2}.
$$
This is a weight one modular form of level~$\Gamma_0(3)$, and incidentally also an Eisenstein series 
$$
\theta_{-3}(\tau) = 1 + 6 \sum_{n=1}^{\infty} \left( \sum_{d \mid n} \chip(n) \right) q^n  \in M_1(\Gamma_0(3),\chip).  
$$
On $\Gamma_0(6)$, we get the weight one Eisenstein series 
\begin{equation*}
\begin{aligned}
A := \frac{\theta_{-3}(\tau) + \theta_{-3}(2\tau)}{2}  & =  1 + 3  \sum_{n=1}^{\infty}  \frac{\chip(n) q^n}{1 - q^n}  + 3  \sum_{n=1}^{\infty}  \frac{\chip(n) q^{2n}}{1 - q^{2n}}
\\ & =  1 + 3 q + 3 q^2 + 3 q^3 + \ldots \in M_1(\Gamma_0(6), \chip).
\end{aligned}
\end{equation*}
Further we have these weight three
 Eisenstein series in~$M_3(\Gamma_0(6),\chip)$: 
 $$
 \sum_{n=1}^{\infty} \left(\sum_{d|n} (-1)^{d-1} \chip(n/d)  d^2  \right) q^n
$$
and
$$
 \sum_{n=1}^{\infty} \left(\sum_{d|n} \chip(d)  d^2  \right) q^n -  \sum_{n=1}^{\infty} \left(\sum_{d|n} \chip(d)  d^2  \right) q^{2n}.
$$
Let us write them respectively as~$\theta^2 B$ and $\theta^2 C$, 
where~$\theta = (2\pi i)^{-1} d/d \tau = q d/dq$, and the \emph{Eichler integrals} $B$ and $C$ compute to the following: 
\begin{equation*}
\begin{aligned}
B & =   \sum_{n=1}^{\infty} \left(\sum_{d|n} (-1)^{d-1} \chip(n/d)  d^2  \right) \frac{q^n}{n^2}  = \sum_{n=1}^{\infty}  \frac{\chip(n) q^n}{n^2(1 - q^n)} - 2  \sum_{n=1}^{\infty}  \frac{\chip(n) q^{2n}}{n^2(1 - q^{2n})} \\ & = \sum_{n=1}^{\infty} \chip(n) n^{-2} \frac{q^n}{1+q^n} = q - \frac{5 q^2}{4} + q^3 -  \frac{11 q^4}{16} + \frac{44 q^5}{25} +  \ldots, \\
C & = \sum_{n=1}^{\infty} \left(\sum_{d|n} \chip(d)  d^2  \right) \frac{q^n}{n^2} - \frac{1}{4} \sum_{n=1}^{\infty} \left(\sum_{d|n} \chip(d)  d^2  \right) \frac{q^{2n}}{n^2} \\
& = \frac{1}{4} \sum_{n=1}^{\infty} \chip(n) \left( 4 \Li_2(q^n) - \Li_2(q^{2n}) \right)  =  
q -  q^2 + \frac{q^3}{9} +  q^4 -  \frac{24 q^5}{25} + \ldots 
\end{aligned}
\end{equation*}
These formulas make it plain\footnote{We have $\Li_2(1) = \zeta(2)$, and the Ces\`aro regularization $\sum_{n=1}^{\infty} \chip(n) := 1/3$ out of the average of the 
partial sums $1, 0, 0, 1, 0,0, 1,0,0,\ldots$ is the relevant interpretation in this context. This calculation and heuristic are readily made rigorous after an Abel summation.} that $\lim_{q \to 1} B = \frac{1}{2} L(2,\chip)$ and $\lim_{q \to 1} C = \frac{1}{4}\zeta(2)$. Moreover, by canceling modularity factors
in the opposite weights~$1$ and~$-1$ (the latter coming
from basic properties~\cite{WeilHeckeLemma} of Eichler integrals), the forms
$A(B -  \frac{1}{2} L(2,\chip))$ and $A(C -  \frac{1}{4}\zeta(2))$ have a weight zero symmetry around the $\Gamma_0(6)$ cusp $\tau = 0$. Expressing these two products in the Hauptmodul
coordinate $x$ leads to holonomic functions on~$Y_0(6) \cong \PP^1 \setminus \{0,1/9, 1, \infty\}$ (see, for example, \cite[\S~2.3]{KontsevichZagier}) which are overconvergent at~$x = 1/9$,
and such that the  coefficients of the factors $AB$ and~$AC$  in~$\Q \llbracket x \rrbracket$ give rise to
 simultaneous Ap\'ery limits $ \frac{1}{2} L(2,\chip)$ and $\frac{1}{4}\zeta(2)$ when compared to the coefficients of~$A$ (also
 considered as a function of~$x$). This is Beukers's framework~\cite{Beukers} for irrationality proofs. Our next lemma collects these remarks with indications on how to read them off from~\cite{Zagier}.

\begin{lemma} \label{defabc} Define power series~$H_A(x)$, $H_B(x)$,  and~$H_C(x)$
in terms of the following formulas:
$$H_A(x) = A(q), \qquad \frac{H_B(x)}{H_A(x)}  =  B(q), \qquad \frac{H_C(x)}{H_A(x)}  =  C(q),$$
where~$x = x(q)$ is as in Equation~\ref{haupt}, so
\begin{equation*}
\begin{aligned}
H_A(x) & = 1 + 3 x + 15 x^2 + 93 x^3 + \ldots = \sum a_n x^n,\\
H_B(x) &  = x + \frac{23 x^2}{4} + \frac{145 x^3}{4} + \frac{3993 x^4}{16}  + \ldots = \sum b_n x^n,\\
H_C(x) & = x + 6 x^2 + \frac{343 x^3}{9} + \frac{788 x^4}{3} + \ldots = \sum c_n x^n.
\end{aligned}
\end{equation*}
Then:
\begin{enumerate}
\item The functions~$H_A(x)$, $H_B(x)$, $H_C(x)$ are multivalued  holonomic functions
on~$\mathbf{P}^1 \setminus \{0,1/9,1,\infty\}$; that is, they extend to holomorphic functions
on the universal cover.
\item We have~$a_n \in \Z$ and~$[1,2,\ldots,n]^2 b_n, [1,2,\ldots,n]^2 c_n \in \Z$.
\item The radius of convergence of~$H_A(x)$, $H_B(x)$, and~$H_C(x)$ is~$R = 1/9$. However,
any linear combination of the  following two functions:
$$H_B(x) - \frac{L(2,\chip)}{2} H_A(x),    \qquad H_C(x) - \frac{\zeta(2)}{4} H_A(x)$$
 has radius of convergence~$R = 1$, where
 $$ L(2,\chip) = \sum_{n=1}^{\infty} \frac{\chip(n)}{n^2}, \qquad
 \zeta(2) = \sum_{n=1}^{\infty} \frac{1}{n^2} = \frac{\pi^2}{6}.$$
 \end{enumerate}
\end{lemma}

\begin{proof}
This result follows from~\cite[Table~3]{Zagier} and~\cite[Table~5]{Zagier} (using an argument previously used by Beukers~\cite{Beukers}).
To orient the reader, note that the sequence~$a_n$ is none other than Zagier's sequence~{\bf C\rm} from~\cite{Zagier}.
Namely, there is an equality
$$a_n = \sum_{k=0}^{n} \binom{n}{k}^2 \binom{2k}{k},$$
and~$a_n$
 satisfies the recurrence
\begin{equation}
\label{zagr}
(n+1)^2 a_{n+1} - An(n+1) a_{n} + B n^2 a_{n-1} = \lambda a_n
\end{equation}
for all~$n$ (\cite[Equation~(3)]{Zagier}). Moreover, $b_n$ satisfies the same recurrence~(\ref{zagr}) for all~$n \ne 0$.
In particular, the facts above concerning~$H_A(x)$ and~$H_B(x)$ are explained in~\S6 of~\cite{Zagier}.
On the other hand, $c_n$ satisfies the recurrence
\begin{equation}
\label{zagrt}
(n+1)^2 c_{n+1} - An(n+1) c_{n} + B n^2 c_{n-1} =1 +  \lambda c_n
\end{equation}
for all~$n \ge 0$. While this sequence is not explicitly in~\cite{Zagier}, it is a disguised form of Zagier's sequence {\bf A\rm}.
More precisely, if one defines the functions
\begin{equation} \label{Gs}
\begin{aligned}
G_A(x) = & \ \frac{1}{1+x} \cdot H_A \left(\frac{x}{x+1}\right)  = 1 + 2 x + 10 x^2 + 56 x^3 + \ldots \\
G_B(x) = & \ \frac{1}{1+x} \cdot H_B \left(\frac{x}{x+1}\right) =  x + \frac{15 x^2}{4} + 22 x^3+ \ldots \\
G_C(x) = & \ \frac{1}{1+x} \cdot H_C \left(\frac{x}{x+1}\right)=  x + 4 x^2 + \frac{208 x^3}{9} + \ldots 
\end{aligned}
\end{equation}
then the coefficients of~$G_A(x)$ are exactly Zagier's sequence~{\bf A\rm}, that is,
they satisfy equation~(\ref{zagr}) except now for the values~$(A,B,\lambda) =  (7,-8,2)$, and the coefficients of~$G_C(x)$ now satisfy
the same  recurrence for~$n > 0$. This can be easily proved  by showing that both functions satisfy the same ODE and then checking
that the first few coefficients are in agreement.
\end{proof}

Using this, we deduce the following:

\begin{proposition} \label{functionsH} Suppose there exists a~$\Q$-linear relationship between~$1$, $\zeta(2)$, and~$L(2,\chip)$,
namely, suppose that
$$a + b \cdot L(2,\chip)/2 + c \cdot \zeta(2)/4  = 0$$
for rational numbers~$a$, $b$, and~$c$. Let
\begin{equation}
\label{defH}
\begin{aligned}
H(x)  & := a H_A(x) + b H_B(x) + c H_C(x)  \\
& = b 
\left(H_B(x) - \frac{L(2,\chip)}{2} H_A(x)\right) + c \left(H_C(x) - \frac{\zeta(2)}{4} H_A(x)\right).
\end{aligned}
\end{equation}
Then~$H(x) \in \Q\llbracket x \rrbracket$ with denominators of shape~$[1,2,\ldots,n]^2$,
and~$H(x)$ satisfies the ODE
\begin{equation}
x(1-x)(1 - 9x)y'' + (1 - 20x +27x^2)y' + 3(-1 + 3x)y = b + \frac{c}{1-x}
\end{equation}
\end{proposition}

\begin{proof} The rationality claims were established above.
Either from equation~(\ref{zagr} and~(\ref{zagrt}) or more directly using the definitions
in terms of Eichler integrals following~\cite{Beukers}, one verifies that~$y = H(x)$ satisfies the given differential equation.
\end{proof}

\begin{remark}  \label{euler rep}
In~\cite{BeukersIntegral}, Beukers gave alternate proofs
of the irrationality of~$\zeta(2)$ and~$\zeta(3)$ in terms of multiple
integrals.  For example, Ap\'{e}ry's  approximations to~$\zeta(2)$
were seen to be 
 coming directly from the integral
\begin{equation} \label{Z2integral}
\iint_{[0,1]^2} \frac{t^n (1-t)^n s^n (1-s)^n}{(1 - s t)^{n+1}} \, ds dt,
\end{equation}
when evaluated as~$a_n - b_n \zeta(2)$. 
We note that the approximations to~$L(2,\chi_{-3})$ and~$\zeta(2)$
considered here (and in~\cite{Zagier}) can also be viewed in the same way. In particular, 
one can easily verify (either by hand using a little effort or  by using~\cite{AZ} without any effort)
the identities:
\begin{equation} \label{Cintegral}
L(2,\chi_{-3}) H_A(x) - 2 H_B(x) 
 =   \sum_{n=0}^{\infty}  x^n \iint_{[0,1]^2} \frac{9^n s^n t^n (1-s^3)^n (1 - t^3)^n}{(1 + s t + s^2 t^2)^{2n+1} } ds dt,
 \end{equation}
 and
 \begin{equation} \label{Aintegral}
 \zeta(2) G_A(x) - 4 G_C(x)
 =   \sum_{n=0}^{\infty} x^n (-1)^n \iint_{[0,1]^2} \frac{(1-s^2)^n (1-t^2)^{n}}{(1 - s t)^{n+1}}  ds dt.
 \end{equation}
One easily finds that
\begin{equation} \label{boundary singularity as a max}
\max_{[0,1]^2} \left| \frac{9 s t (1-s^3)(1 - t^3)}{(1 + s t + s^2 t^2)^2}\right|  = 1,
\qquad
\max_{[0,1]^2} \left| \frac{(1-s^2)(1-t^2)}{1 - s t}\right|  = 1,
\end{equation}
the maxima being obtained at~$s=t=1/2$ in the first case and~$s=t=0$ in the second. It
follows that the
integrals are all  bounded  by~$1$ which gives a transparent proof that
 the functions~$H(x)$ considered in Proposition~\ref{functionsH} overconverge beyond the singularity at~$x = 1/9$ to the entire unit disc.
 
We can also evaluate the geometric series to express the integral formula~\eqref{Cintegral} as
\begin{equation} \label{Erep}
L(2,\chi_{-3}) H_A(x) - 2 H_B(x)  = \iint_{[0,1]^2} \frac{1 + s t + s^2 t^2}{(1 + s t + s^2 t^2)^2 - 9st(1-t^3)(1-s^3) x} \, ds  dt. 
\end{equation}
By~\eqref{boundary singularity as a max}, this formula represents the function $L(2,\chi_{-3}) H_A(x) - 2 H_B(x)$ as a 
continuous integral of a family of (rational, as it happens) functions $f_{s,t}(x) \in \mathcal{O}\left(\C \setminus [1,\infty) \right)$ \emph{holomorphic}
on $\C \setminus [1,\infty)$. The property of being a holomorphic function over a complex domain is 
inherited by any continuous integration over a parameter $(s,t) \in [0,1]^2$, and so the integral 
representation makes equally transparent the analyticity of~\eqref{Erep} on $\C \setminus [1,\infty)$.  This is Zudilin's point of view in~\cite{ZudilinDet}.~\endofremark
\end{remark}

\begin{remark}
These integral representations, and especially~\eqref{Erep}, may also be compared to Zudilin's~\cite{ZudilinCatalan,RivoalCatalan,Nesterenko}
\begin{equation*}
\begin{aligned}
L(2,\chi_{-4})U(x) - V(x)  & := \iint_{[0,1]^2} \frac{ds \, dt}{\sqrt{(s-s^2)(t-t^2)} \cdot \big( 1 - st - (s-s^2)(t-t^2)x \big) } \\
 & = - \sum_{n=0}^{\infty }   x^n \iint_{[0,1]^2}  \frac{(s-s^2)^{n-1/2} (t-t^2)^{n-1/2}}{(1-st)^{n+1}}  \, ds dt
 \end{aligned}
 \end{equation*}
 giving rational approximants to Catalan's constant $G = L(2,\chi_{-4})$. In this instance, 
 the integrand peak rate is 
 $$
 \max_{[0,1]^2} \left| \frac{(s-s^2)(t-t^2)}{1-st}\right|  = \left(  \frac{1+\sqrt{5}}{2}  \right)^{-5},
$$
and indeed the singularities of the linear ODE are at $0, \left(  \frac{1 \pm \sqrt{5}}{2} \right)^5$, and
$\infty$: precisely the same as for Ap\'ery's approximants to $\zeta(2)$. Unfortunately, due to the 
half-integral exponents in this integral representation, the denominators in these rational approximations
are as big as $16^n [1,\ldots,2n]^2$. 

A different holonomic sequence of rational approximants to~$G$ was given by case~{\bf E} in~\cite{Zagier}, where the ODE
singularities are $\{0,1/8; 1/4, \infty \}$ (the first two of which are overconvergent), and the denominator types
are $[1,\ldots,n]^2$. 

Since  $e^2 > 16 \cdot (1/4)$ and $e^4 >  \left(  \frac{1+\sqrt{5}}{2}  \right)^{5}$ (by a wide margin!), this 
definitely precludes an approach to the irrationality of the Catalan constant by our method using either
of these particular families of rational approximants,
 unless
some completely new idea is discovered.~\endofremark
\end{remark}

\subsection{The symmetrization of~\texorpdfstring{$H(x)$}{H(x)}}  \label{SympH}
Let~$a$, $b$, and~$c$ be complex numbers such that
$$
a + b \cdot L(2,\chip)/2 + c \cdot \zeta(2)/4  = 0.
$$
Then we may define~$H(x) \in \C \llbracket x \rrbracket$ as in equation~(\ref{defH}).
If we 
additionally assume that~$1$, $\pi^2$, and~$L(2,\chi_{-3})$ are linearly dependent over~$\Q$, then we can choose~$a$, $b$, and~$c$
to be rational, although the arguments of this section will not require this hypothesis.

 We now let~$G(y)$ be the symmetrization of~$H(x)$:

\begin{df} \label{defG} Let~$G(y) = \Sym^{+} H(x)$ as defined in equation~(\ref{plusminus}), so
$$G(y) = H(x) + H\left( \frac{x}{x-1} \right) \in \C \llbracket y \rrbracket,$$
and let~$G_A(x) = \Sym^{+} H_A(x)$, so
$$G_A(y) = H_A(x) + H_A\left( \frac{x}{x-1} \right) \in \Z \llbracket y \rrbracket.$$
\end{df}
Note that~$G_A(y) = 2 - 27 y + 1014 y^2 - 49536 y^3 + \ldots $;
the function~$G_A(y)$ satisfies an order~$4$ ODE
$$\sum_{i=0}^{3} c_i(y) G^{(i)}_A(y) = 0$$
which we give explicitly later in equation~(\ref{explicitode}).
The span of~$G_A(y)$ and its derivatives generates, over~$\Q(x)$, the  space spanned by~$H_A(x)$ and~$H_A(x/(x-1))$  and their derivatives (which
are both vector spaces of dimension~$4$).
By  Lemma~\ref{etalecover} (\ref{itemnicer}), we immediately  have the following:

\begin{lemma} \label{Gdenominatortype} If~$a,b,c \in \Q$, then~$G(y)$
has denominator type~$[1,2,\ldots,2n]^2$.
\end{lemma}

\section{Functional Linear Independence}

\label{sec:lindep}

Let~$A_i(x) \in \Q \llbracket x \rrbracket $ be a collection of holomorphic functions functions on~$\mathbf{P}^1 \setminus S$ for some finite set~$S$. (In our situation, they will
all be Siegel~$G$-functions.)
Suppose we wish to prove that the~$A_i(x)$ are linearly independent over~$\Q(x)$ or~$\C(x)$. One strategy is as follows. Let~$\gamma$ be a path in~$\C$; for example,
take a path
starting at~$x=0$, avoiding other points in~$S$, and then returning to~$0$. The functions~$A_i(x)$ can be analytically continued along~$\gamma$, and as we
return to~$x=0$ we obtain a sequence of functions~$\AA_i(x)$ which may now have singularities at~$x=0$.
Certainly any polynomial relationship between the~$A_i(x)$ extends to (the same) polynomial relationship between the~$\AA_i(x)$,
and hence also to a polynomial relationship
between the~$\AA_i(x) - A_i(x)$, which can sometimes be useful.
But we can alternatively consider
any identity between the~$\AA_i(x)$ \emph{modulo} functions which are holomorphic at~$0$.
What may (and often does) happen in principle is that this reduces a linear relationship between a  large number of functions to a smaller number of functions,
and one can hope to employ some form of inductive strategy to establish full linear independence.
 A typical example is as follows: Suppose that the path~$\gamma$ starts at~$0 \in S$ and is a simple loop around a single point~$1 \in S$.
Then, if a proper subset of the functions~$A_i(x)$ are actually holomorphic at~$x=1$, the corresponding~$\AA_i(x)$ vanish modulo holomorphic functions, and 
we obtain a corresponding linear relationship between the~$\AA_i(x)$ with fewer terms.
A basic example of this is as follows. Suppose that
$$A_1(x) = 1, \quad A_2(x) = \log(1-x), \quad A_3(x) = \Li_2(x).$$
After a  suitably oriented loop around zero, we have
$$\AA_1(x) = 1, \quad \AA_2(x) = \log(1-x) + 2 \pi i, \quad \AA_3(x) = \Li_2(x) + 2 \pi i \log(x).$$
Now, modulo holomorphic functions at zero, we obtain a linear relationship between the three functions~$0$, $0$, and~$2 \pi i \log(x)$. Clearly this forces the coefficient of~$\AA_3(x)$ to be zero,
and reduces us to showing that~$A_1(x)$ and~$A_2(x)$ are linearly independent because~$\log(1-x)$ is not a rational function.
We will use this strategy a number of times below.
Note that another argument in this case would
be to consider the functions~$\AA_i(x) - A_i(x)$ which reduces the problem
to the~$\C(x)$-linear independence of~$1$ and~$\log x$.

\begin{lemma} \label{seven} The  seven functions~$B_i(y)$ for~$i=1,\ldots,7$ defined
in Section~\ref{sec:purefunctions} in equations~(\ref{defB}), (\ref{defB5}), and~(\ref{defB67})
respectively are linearly independent over~$\C(y)$.
\end{lemma}

\begin{figure}
\begin{center}
\begin{tikzpicture}[scale = 0.7]
  \draw[->] (-3,0) -- (6,0) node[below]{};
  \draw[->] (0,-3) -- (0,3) node[left] {};
  \fill (0,0) circle (2pt) node[below left] {0};
  \fill (3,0) circle (2pt) node[below] {1};
  \node at (0.6,-0.3) {$\alpha$};
 \draw[line width=1.5pt,postaction={decorate, decoration={markings, mark=at position 0.5 with {\arrow{>}}}}] (0,0) -- (1.5,0);
  \draw[line width=1.5pt, postaction={decorate, decoration={markings, mark=at position 0.20 with {\arrow{<}}}}]
  (4.5,0) arc (0:360:1.5);
  \node at (3.4,1) {$\beta$};
    \draw[line width=1.5pt, postaction={decorate, decoration={markings, mark=at position 0.20 with {\arrow{<}}}}]
  (1.5,0) arc (0:360:1.5);
  \node at (0.4,1) {$\delta$};
\end{tikzpicture}
\caption{The path~$\gamma = \alpha^{-1} \beta^{-1} \delta \beta \alpha$}
\label{figone}
\end{center}
\end{figure}
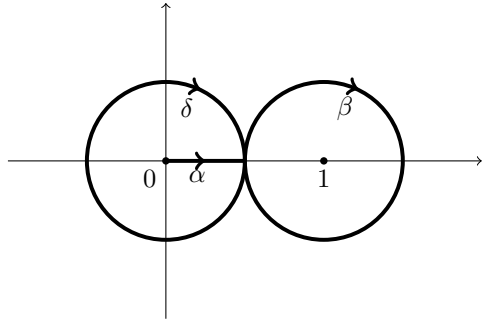

\begin{proof} We first of all note that the~$B_i(y)$ are all elements of~$\Q  \llbracket y \rrbracket$.
Therefore any linear dependency over~$\C(y)$ upgrades to one over~$\Q(y)$, and so it suffices to prove the result
over~$\Q(y)$.

We begin by proving the linear independence of the~$B_i(y)$ for~$i=1,\ldots,5$.
By Lemma~\ref{etalecover}, it suffices to prove the linear independence
of the~$5$-functions~$1$, $\log(1-x)$,  $\log^2(1-x)$, $\Li_2(x)$,
and
and~$J(x)$ of equation~(\ref{newfunction}) over~$\Q(x)$.
Certainly~$1$, $\log(1-x)$, and~$\log^2(1-x)$
are independent since~$\log(1-x)$
is transcendental over~$\Q(x)$.
These three functions are also defined
on~$\mathbf{P}^1 \setminus \{1,\infty\}$
which distinguishes then from~$\Li_2(x)$   ---  
take a path~$\gamma$ from~$0$
which winds around~$x=1$, then winds around~$x=0$,
then winds (in the opposite way) around~$x=1$, and
returns to zero (as in Figure~\ref{figone}). The first three functions will be invariant,
but~$\Li_2(x)$ has non-trivial monodromy on this path. So~$\Li_2(x)$
is independent of these previous functions.
Now
suppose that~$J(x)$ was a~$\Q(x)$-linear
combination of~$1, \log(1-x), \log^2(1-x), \Li_2(x)$.
If the coefficient of~$\Li_2(x)$ was non-trivial,
then, after scaling, we may assume that it is~$1$.
But now by differentiation, and using
the ODE~\ref{odeJ} for~$J(x)$,
we obtain a new relation over~$\Q(x)$ with~$1, \log(1-x),
\log^2(1-x)$, and~$J(x)$ only (with a non-trivial coefficient of~$J(x)$
because of the ODE). But~$J(x)$ also
has non-trivial monodromy over the path~$\gamma$,
so we conclude as for~$\Li_2(x)$ above.
Now let us return to the~$\mathbf{P}^1 \setminus \{0,4,
\infty\}$ domain and consider the functions~$B_6(y)$
and~$B_7(y)$.
Recall that:
$$B_6(y) = \int \frac{B_3(y)}{y} \, dy, \quad
B_7(y) = \int \frac{B_4(y)}{y} \, dy.$$
Using the derivation formula
\begin{equation}
\label{simpledif}
\frac{d}{dy} \left\{ A(y) \int F(y) \, dy \right\} = A'(y) \int F(y) \,  dy + A(y) F(y),
\end{equation}
we can firstly assume that the  coefficients of our linear relation in~$\Q(y)$
are polynomials,
and then by differentiation reduce to an equality of the form
\begin{equation}
\label{lasttwo}
a_0  \int \frac{B_3(y)}{y} \, dy + a_1  \int \frac{B_4(y)}{y} \, dy
= \sum_{i=1}^{5} b_i(y) B_i(y),
\end{equation}
where now~$a_0$ and~$a_1$ are constants which
are not both zero and~$b_i(y) \in \Q(y)$.
We note the 
\begin{equation}
\label{diffformulas}
\begin{aligned}
y(4-y) \frac{d}{dy} B_2(y)  = & \  (2-y) B_2(y)  + y^2, \\
y(4-y) \frac{d}{dy} B_3(y) = & \ - B_2(y) + 2 y, \\
y(4-y) \frac{d}{dy} B_4(y) = & \ (2-y) B_4(y) + (4-y) B_2(y) - 2 y(4-y), \\
2 y(4 - y) \frac{d}{dy} B_5(y) =  & \  (4-y) B_5(y) - 2 B_2(y) + 4 y), \\
\end{aligned}
\end{equation}
But now let us differentiate~(\ref{lasttwo}) to get
\begin{equation}
\label{lasttwotwo}
a_0   \frac{B_3(y)}{y} + a_1  \frac{B_4(y)}{y}
= \sum_{i=1}^{5} b_i(y) B'_i(y) + b'_i(y) B_i(y).
\end{equation}
We see from equation~(\ref{diffformulas}) and equating the coefficients
of~$B_4(y)$ and~$B_3(y)$ respectively (and using the linear independence
of~$B_i(y)$ for~$i=1,\ldots,5$ that
\begin{equation}
\begin{aligned}
\frac{a_0}{y} = & \ b'_3(y), \\
\frac{a_1}{y} = & \ b'_4(y) + \frac{(2-y)}{y(4-y)} b_4(y).
\end{aligned}
\end{equation}
From the former equation we get (up to constant)
$$b_3(y) = a_0 \log(y),$$
which, since~$b_3(y) \in \Q(y)$, can only happen if~$a_0 = 0$.
We may also write the latter equation as
$$\frac{a_1 \sqrt{y(4-y)}}{y}= \frac{d}{dy} \left(  b_3(y) \sqrt{y(4-y)} \right)$$
But the integral of the left-hand side is not algebraic (if~$a_1 \ne 0$) but
the integral of the right-hand side is, so once more this can only happen
when~$a_1 = 0$, and the linear independence is established.
\end{proof}

\addtocounter{subsubsection}{1}

\subsection{Linear Independence of pure functions
and Zagier functions}   \label{jointli}

Recall from Definition~\ref{defG} that 
$$G(x)= H(x) + H \left( \frac{x}{x-1}\right).$$
We also let
$$G_A(x) = H_A(x) + H_A \left( \frac{x}{x-1}\right) \in \Z \llbracket x \rrbracket,$$
which is a homogenous solution to a degree~$4$ ODE (given explicitly in equation~(\ref{explicitode}) below).
Our final functional linear independence result is as follows:

\begin{lemma}[14 functions] \label{14functions}
The seven functions
$$\begin{aligned}
\int & G(y) \, dy,   \int \frac{G(y) - G(0)}{y} \, dy,
\int \frac{G(y) - G(0) -  G'(0)y}{y^2} dy, \\
& G(y),G'(y),G''(y),G'''(y), \end{aligned}$$
together with the seven functions~$B_i(y)$ for~$i=1,\ldots,7$,
are linearly independent over~$\C(y)$.
\end{lemma}

\begin{remark} \label{whyisthistrue} It is easy enough to discover Lemma~\ref{14functions}
experimentally.
A  collection of power series~$A_i(x)$ which satisfy a linear relation
over~$\C(y)$ also satisfy a linear relation with coefficients in~$\C[y]$,
and thus with coefficients which are polynomials of degree~$\le D$
for some~$D \in \NwithoutzeroA$.
But the  question as to whether there exists such a relation
for any given~$D$  is equivalent to the vanishing of  the
determinant of an explicitly computable matrix. Once one establishes
that there are no such linear relations for~$D$ of moderate size (say~$D=20$),
one is sufficiently convinced the result is true and then one writes down a proof.
We admit that this is how we arrived at both Lemma~\ref{14functions}
and Lemma~\ref{17functions}, even though there is most likely a higher level
proof which better explains the precise numerology.
See also Remark~\ref{not15}.~\endofremark
\end{remark}

\begin{proof} Since the~$B_i(y)$ are linearly independent by Lemma~\ref{seven},
any dependence must include at least one of the terms above with a non-zero coefficient.
Let~$\gamma$ denote a path which first traverses~$4$,
then~$-1/72$, then~$4$ in the opposite direction, and then back to~$0$.
The function~$G(y)$ is replaced by~$\hatG(y)$,
which is a solution to the same non-homogenous differential equation at~$G(y)$.
On the other hand, the functions~$\BB_n(y) = B_n(y)$ remain invariant
Since~$G(y) \ne \hatG(y)$, we obtain an equivalent relation
between the functions
$$\begin{aligned}
\int  & \hatG(y) \, dy, \int \frac{\hatG(y) - G(0)}{y} \, dy,
\int \frac{\hatG(y) - G(0) -  G'(0)y}{y^2} \, dy, \\
& \hatG(y),\hatG'(y),\hatG''(y),\hatG'''(y) \end{aligned} $$
and the~$B_i(y)$ with the same coefficients.
Hence,  with~$\Delta = \hatG(y) - G(y)$, we obtain a non-zero~$\C(y)$-linear relationship between
the seven functions
$$\int \Delta(y) \, dy, \int \frac{\Delta(y)}{y} \, dy,
\int \frac{\Delta(y)}{y^2} \, dy,
\Delta(y),\Delta'(y),\Delta''(y),\Delta'''(y)$$
But~$\Delta(y)$ is now a homogenous solution to the corresponding degree~$4$ ODE
which is irreducible, and so by replacing~$\Delta(y)$ by its translates
under elements of the monodromy group~$\pi_1(\mathbf{P}^1 \setminus \{0,-1/72,4,\infty\})$,
we deduce that the corresponding linear relation must hold for
any such~$\Delta(y)$, including in particular the holomorphic solution~$\Delta(y) = G_A(y)$.

Assume such a linear relation exists. After scaling,
we may assume that the coefficients lie in~$\C(y)$
Using~(\ref{simpledif}) again:
\begin{equation}
\label{simpleagain}
\frac{d}{dy} \left\{ A(y) \int F(y) \, dy \right\} = A'(y) \int F(y) \, dy + A(y) F(y ),
\end{equation}
after repeated differentiation we may assume that the coefficients
of
the three integral terms are all constants, and that at least one is non-zero.
Hence there exists a relation
\begin{equation}
\label{todiff}
a_0 \int G_A(y) \, dy + a_{-1} \int \frac{G_A(y)}{y} \, dy + a_{-2}
\int \frac{G_A(y)}{y^2} \, dy
= \sum_{i=0}^{3} b_i(y)  G_A^{(i)}(y).
\end{equation}
Note that we cannot insist that the~$b_i(y) \in \C[y]$, for two reasons.
First is that the derivative terms from the integrals involve~$G_A(y)$ divided by powers of~$y$.
But also when differentiating~$G^{(3)}_A(y)$ we obtain~$G^{(4)}_A(y)$, and to write this in terms
of lower order derivatives in~$G^{(i)}_A(y)$ we need to divide by the leading term in the differential equation.
In fact,
$G_A(x)$
 satisfies the ODE
$$\sum_{i=0}^{4} c_i(y) G_A^{(i)}(y) = 0,$$
where~$c_i(y)$ are defined as follows:
\begin{equation} \label{explicitode}
\begin{aligned}
c_0(y) = & \   
-18 (3 + 126 y - 712 y^2 + 360 y^3), \\
c_1(y) = & \    2 (-2 - 2761 y + 141632 y^2 - 
   280328 y^3 + 176412 y^4 - 95616 y^5 + 20736 y^6), \\
 c_2(y) = & \     2 y (-34 - 
   6353 y + 690355 y^2 - 1065613 y^3 + 867876 y^4 - 438336 y^5 + 
   72576 y^6),  \\
   c_3(y) = & \  2 (-4 + y) y^2 (10 + 204 y - 118195 y^2 + 146946 y^3 - 
   142848 y^4 + 41472 y^5), \\
   c_4(y) = & \  (-4 + y)^2 y^3 (1 + 72 y) (-1 + 118 y - 
   122 y^2 + 144 y^3). \end{aligned}
   \end{equation}
Thus we can assume that~$b_i(y) \in \C[y,c_4(y)^{-1}]$.
Differentiating equation~(\ref{todiff}) one more time gives an identity
\begin{equation} \label{remaining}
a_0 G_A(y) + a_{-1} \frac{G_A(y)}{y}  + a_{-2}
\frac{G_A(y)}{y^2} 
= \sum_{i=0}^{3} b'_i(y)  G_A^{(i)}(y) +  b_i(y) G_A^{(i+1)}(y).
\end{equation}
We rewrite~(\ref{remaining}) as
\begin{equation}
\begin{aligned}
& a G_A(y) + b \frac{G_A(y)}{y}  + c
\frac{G_A(y)}{y^2}  \\
& =   \sum_{i=0}^{3} b'_i(y)  G_A^{(i)}(y) +  
\sum_{i=1}^{3} b_{i-1}(y) G_A^{(i)}(y)
+ b_3(y) G^{(4)}(y), \\
& =   \sum_{i=0}^{3} b'_i(y)  G_A^{(i)}(y) +  
\sum_{i=1}^{3} b_{i-1}(y) G_A^{(i)}(y)
- \sum_{i=0}^{3} \frac{c_i(y)}{c_4(y)}  b_3(y) G^{(i)}(y),  \end{aligned}
\end{equation}
and thus we deduce the simultaneous equations
\begin{equation}
\begin{aligned}
b'_3(y) + b_2(y) - \frac{c_3(y)}{c_4(y)} b_3(y) = & \ 0, \\
b'_2(y) + b_1(y) - \frac{c_2(y)}{c_4(y)} b_3(y) = & \ 0, \\
b'_1(y) + b_0(y) - \frac{c_1(y)}{c_4(y)} b_3(y) = & \ 0, \\
b'_0(y)  - \frac{c_0(y)}{c_4(y)} b_3(y) = & \ a_0 + \frac{a_{-1}}{y} + \frac{a_{-2}}{y^2}, 
\end{aligned}
\end{equation}
Recall that~$b_3(y) \in \C[y,c_4(y)^{-1}]$. Moreover, given~$b_3(y)$, one can inductively solve for~$b_i$ for~$i \in \{2,1,0\}$
from the equation
$$b_i(y) =  \frac{c_{i+1}(y)}{c_4(y)} b_3(y) - b'_{i+1}(y).$$
Our strategy is as follows. Since~$b_3(y) \in \C(y)$, we can consider the power series expansion of~$b_3(y)$
around~$\infty$ and around any point~$\alpha \in \C$. Then, by considering the final equation, we obtain
an explicit bound on the order of any pole of~$b_3(y)$ at~$\alpha$ (note that~$b_3(y)$
will be holomorphic unless~$\alpha = \infty$ or is a root of~$c_4(y)$). But that confines~$b_3(y)$
to be a rational function such that divisor~$(b_3(y)) + D \ge 0$ for an explicit divisor~$D$ supported at the roots of~$c_4(y)$.
This is a finite dimensional (explicitly computable) vector space, and then we can solve for all possible~$b_3(y)$
using linear algebra.
Another way to view this is to think of this system as a (non-homogenous) ODE in~$b_3(y)$, and we are computing
the (possible) local expansions around any point using the Frobenius method.
Explicitly, with~$\alpha \ne \infty$, and
$$b_3(y) = \sum_{i=N}^{\infty} r_{i}(y -\alpha)^i,$$
(with~$N = N_{\alpha}$ and~$r_i -= r_{i,\alpha}$, and suppressing the subscript below)
then:
\begin{enumerate}
\item If~$\alpha = 0$, 
 the last equality becomes:
$$a_0 + \frac{a_{-1}}{y} + \frac{a_{-2}}{y^2} =  \frac{-1}{4} (3-N)^2 (5-2N)^2 r_{N} y^{N-4} + \ldots $$
\item If~$\alpha = 4$, 
the last equality becomes:
$$a_0 + \frac{a_{-1}}{y} + \frac{a_{-2}}{y^2} =  \frac{-1}{4} (3-N) (2-N)(5-2N)(3-2N) r_{N} (y-4)^{N-4} + \ldots $$
\item If~$\alpha = -1/72$,
the last equality becomes:
$$a_0 + \frac{a_{-1}}{y} + \frac{a_{-2}}{y^2} =   (3-N)^2 (2-N)(1-N) r_{N} (y + 1/72)^{N-4} + \ldots $$
\item If~$\alpha$ is a root~$\beta$ of~$144 y^3 - 122 y^2 + 118 y - 1 = 0$, then
$$a_0 + \frac{a_{-1}}{y} + \frac{a_{-2}}{y^2} =   (3-N) (2-N)(1-N)(1+N) r_{N} (y - \beta)^{N-4} + \ldots $$
\item At~$\alpha \rightarrow \infty$, with
$$b_3(y) =  y^N \sum_{i=N}^{\infty} r_{i} y^{-i},$$
we have
$$a_0 + \frac{a_{-1}}{y} + \frac{a_{-2}}{y^2} = \frac{-1}{4} (3-N)^2 (5-2N)^2 r_{N} y^{N-4} + \ldots $$
\end{enumerate}
From this we deduce that:
\begin{equation}
\begin{aligned}
N_0 \ge & \ 2 \\
N_4 \ge & \ 2 \\
N_{-1/72} \ge & \  1 \\
N_{\beta} \ge & \ -1 \\
N_{\infty} \le & \ 4.
\end{aligned}
\end{equation}
From this, it follows that
\begin{equation}
\label{b3eq}
b_3(y) = \frac{y^2(y-4)^2 (y + 1/72)}{(144 y^3 - 122 y^2 + 118 y - 1)} Q(y),
\end{equation}
where~$Q(y)$ is a polynomial of degree at most~$2$.
However, if we write
$$Q(y) =q_0 + q_1 y + q_2 y^2,$$
then we find that
\begin{equation} \label{degree12}
a_0 + \frac{a_{-1}}{y} + \frac{a_{-2}}{y^2} =  \frac{
52542464 y^{12} q_2 + \ldots }
{36 y^2 (-1 + 118 y - 122 y^2 + 144 y^3)^4} 
\end{equation}
where the numerator on the right-hand side is a degree~$12$ polynomial with coefficients linear in~$\Z q_0 \oplus \Z q_1 \oplus \Z q_2$.
But now by linear algebra one can directly check that there are no
choices
of the parameters~$q_i$ to even make the numerator vanish to order (at least) one at a non-zero root of the denominator.
Hence no such~$b_3(y)$ exists, and we are done.
\end{proof}

\begin{remark} \label{not15} Suppose instead we had tried to prove the (false!) linear independence of
the seven functions~$B_i(y)$ together with
$$\begin{aligned}&  \int G(y) \, dy, \int \frac{G(y) - G(0)}{y} \, dy,
\int \frac{G(y) - G(0) -  G'(0)y}{y^2} \, dy, \\
 & \int \frac{G(y) - G(0) -  G'(0)y - G''(0) y^2}{y^3} , dy,
G(y),G'(y),G''(y),G'''(y), \end{aligned}$$
that is, adding another integral.
Then the argument would have proceeded exactly as above except now we could only deduce that~$N_{0} \ge 1$ rather than~$N_0 \ge 2$. Then, writing
$$Q(y) = \frac{q_{-1}}{y} + q_0 + q_1 y + q_2 y^2,$$
just as in equation~(\ref{degree12}), we would have found that
$$a_0 + \frac{a_{-1}}{y} + \frac{a_{-2}}{y^2} + \frac{a_{-3}}{y^3} =  
\frac{
52542464 y^{13} q_2 + \ldots }
{36 y^3 (-1 + 118 y - 122 y^2 + 144 y^3)^4} $$
 There is now a unique choice of the parameters~$q_i$ up to
scalar which allows us to remove a single factor of the numerator, namely (up to scalar)
$$Q(y) = \frac{1}{y}  - 278  - 844 y - 4644 y^2.$$
With this choice of~$b_3(y)$ as coming from equation~(\ref{b3eq}), the powers of~$(144 y^3 - 122 y^2 + 118 y - 1)$
 disappear completely, and we  arrive at the equality
$$a_0 + \frac{a_{-1}}{y} + \frac{a_{-2}}{y^2} + \frac{a_{-3}}{y^3} =  \frac{676}{9 y^2} +  \frac{2}{y^3}.$$
This, of course, does now have solutions.
This reflects that there \emph{is} a linear dependence between these functions. In fact, there is already a~$\C(y)$-linear
dependence between the functions
$$\int \frac{G_A(y)}{y^2} \, dy,  \, \int \frac{G_A(y)}{y^3}  \, dy, \, 
G_A(y), \, G'_A(y), \, G''_A(y), \, G'''_A(y), \textrm{ and } 1.$$
Note as another consistency check with the solution~$a_{-2} = 676/9$ and~$a_{-3} = 2$, we have
(in analogy with equation~\ref{remaining}) that
$$\begin{aligned}
a_0 G_A(y) + a_{-1} \frac{G_A(y)}{y}  + a_{-2}
\frac{G_A(y)}{y^2} + a_{-3} \frac{G_A(y)}{y^3}
& = \frac{676}{9} \frac{G_A(y)}{y^2} + 2  \frac{G_A(y)}{y^3} \\
& =
\frac{4}{y^3} + \frac{866}{9 y^2} - \frac{68728}{3} + \ldots  \end{aligned}$$
with no~$1/y$ term, consistent with it being a derivative of a meromorphic function at~$y=0$.~\endofremark
\end{remark}

\section{Proof of the linear independence of 
\texorpdfstring{$1, \zeta(2)$}{1,zeta(2)}, 
and 
\texorpdfstring{$L(2,\chi_{-3})$}{L(2chi3)}
}  \label{sec:proofA}

In this section, we complete the proof of Theorem~\ref{mainA} using the results of Appendix~\ref{contour choiceA}. 
The argument is by a contradiction, by proving that a certain $G$-function cannot exist. Suppose for the contradiction
that there exists a $\Q$-linear relation among the periods $1, \zeta(2), L(2,\chip)$, which 
we could write as
\begin{equation} \label{absurd}
a + b \cdot L(2,\chip)/2 + c \cdot \zeta(2)/4  = 0
\end{equation}
with some rational integers $a,b, c \in \Z$, not all zero. Proposition~\ref{functionsH} then constructs a certain $G$-function
$H(x) \in \Q \llbracket x \rrbracket$ with denominator type $[1,\ldots,n]^2$ and continuing holonomically on~$\P^1 \setminus \{0,1/9,1,\infty\}$.

Now~\S~\ref{SympH} converts the $G$-function~$H(x)$ to a $G$-function
$G(y) := \Sym^{+} H(x) \in \Q \llbracket y \rrbracket$ in the symmetrization coordinate 
$$
y := x + x/(x-1) = x^2/(x-1).
$$
 Lemma~\ref{Gdenominatortype}
shows that the denominator type of $G(y) \in \Q\llbracket y \rrbracket$ is $[1,\ldots,2n]^2$. 
On the other hand, $G(y)$ is holonomic on $y \in \P^1 \setminus \{0, 4, \infty, -1/72\}$, holomorphic
on~$y \in \C \setminus [4,\infty)$, and with $\Z/2$ local monodromy around~$y=4$.

We apply Theorem~\ref{main:elementary form} with the $14 \times 2$ denominators type array 
\begin{equation}  \label{typesA}
\mathbf{b} := \left(  \begin{array}{llllllllllllll}   0 & 2 & 2 & 2 & 2 & 2 & 2 & 2 & 2 & 2 & 2 & 2 & 2 & 2 \\
0 & 0 & 0 & 2 & 2 & 2 & 2 & 2 & 2 & 2 & 2 & 2 & 2 & 2  \end{array} \right)^{\mathrm{t}}
\end{equation}
 and the integrations vector 
$$
\mathbf{e} :=  (0, 0, 1; 0, 0, 0, 0, 0, 0; 1, 1, 1, 1, 1),
$$
taking over in~\eqref{den type int} after replacing the letter~$x$ there by the symmetrization letter 
$$
y := x+x/(x-1) = x^2/(x-1); 
$$
and taking the following ordered list of functions $\{f_i\}_{i=1}^{14}$, see~\eqref{defB}, \eqref{Nielsen}, and \eqref{defB67}: 
\begin{equation*}
\begin{gathered}
B_1(y), \, B_2(y), \, B_3(y);  \, B_4(y), \, B_5(y), \, G(y), \, G'(y), \, G''(y), \, G'''(y), \\ B_6(y), \, B_7(y), \, \int G(y) \, dy, \int \frac{G(y) - G(0)}{y} \, dy,
\int \frac{G(y) - G(0) -  G'(0)y}{y^2}  \, dy.
\end{gathered}
\end{equation*}
The~$\Q(y)$-linear independence of these~$14$ functions was proved in Lemma~\ref{14functions}. 
The denominator types were computed in Lemma~\ref{bdenominators}.
For the  integrals, we note that the shift in indexing caused by dividing by powers of~$y$ means
that these functions are not literally of denominator type~$n[1,2,\ldots,2n]^2$ but rather of 
type~$n[1,2,\ldots,2n+3]^2$;  this is not an issue by Remark~\ref{rem:overflow}.

Let us denote by~$\HH_{Y_0(2)}$ the~$14$-dimensional $\Q(y)$-linear span of these functions. 
For the analytic maps~$\boldsymbol{\varphi}$ figuring in the various holonomy bounds we have developed, 
we take restrictions~$\varphi(rz)$ of the holomorphic mapping~$\varphi \in \mathcal{O}(\Db)$ of 
Lemma~\ref{choicephiA}. From Corollary~\ref{stacky overconvergence} 
used with~$\Sigma_{Y_0(2)}^0 := \{-1/72\}$, $\Sigma_{Y_0(2)}^1 := \emptyset$, 
$\varphi_{Y_0(2)} := \varphi$, and~$U_{Y_0(2)}$ a sufficiently small open neighborhood of the line segment~$[-1/72, 0]$, 
we have the analyticity~$\varphi^*\HH_{Y_0(2)} \subset \mathcal{M}(\Db)$.

\addtocounter{subsubsection}{1}
 \begin{figure}[!h]  
\begin{center}
  \includegraphics[width=75mm]{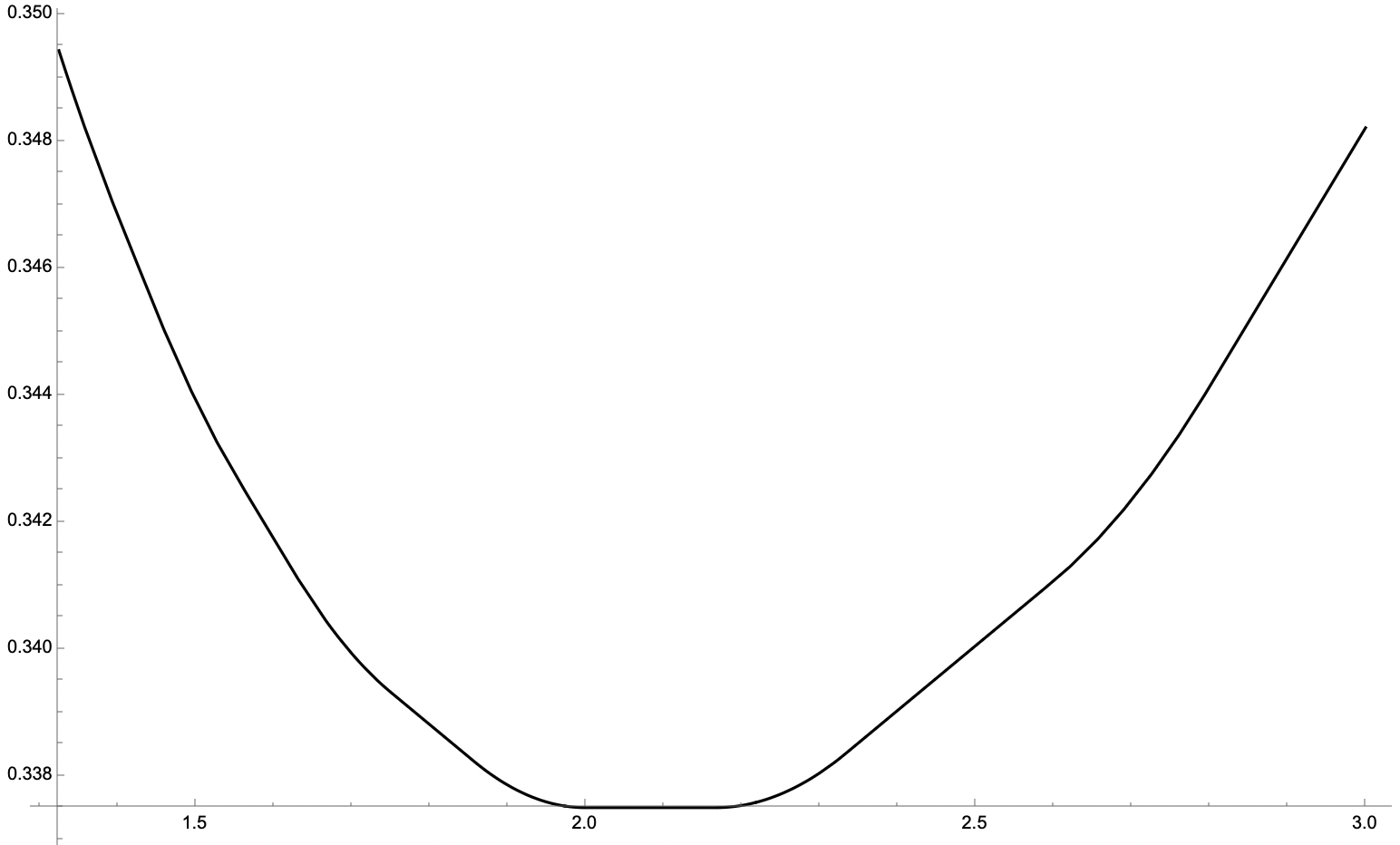}
\end{center}
\caption{The $\xi \in [1.325, 3]$ fragment of the graph of $(6\xi  
   + I_{\xi}^{14}(\xi))/98 $, displaying the interval $\xi \in [2,13/6]$ as the identical minimizer.}
   \label{xigraph}
\end{figure}

For the denominator rates, we calculate
\begin{equation}
\tau^{\flat}(\mathbf{b}) = \frac{ 1 \cdot 0 + (3+5) \cdot 2  + (7 + 9 + 11 + 13 +\ldots  + 27) \cdot 4 }{14^2}  
= \frac{191}{49}
\end{equation}
and, from Figure~\ref{xigraph} which reveals $\xi \in [2,13/6]$ to be the identical minimizer, 
\begin{equation}
\begin{aligned}
\tau^{\sharp}(\mathbf{e}) & =  \frac{2}{m^2}  \min_{\xi \in [0,m]} \left\{  \xi \sum_{i=1}^{m} e_i 
   +
   \left(  \max_{1 \leq i \leq m}  e_i \right) I_{\xi}^{m}(\xi)
   \right\} \\  \label{tauraise}
   & = \min_{\xi \in [0,14]} \left\{  \frac{6\xi  
   + I_{\xi}^{14}(\xi)}{98} \right\} =  \frac{ 12 + I_{2}^{14}(2) }{98}  = \frac{27}{80}. 
   \end{aligned}
   \end{equation}
   
   We obtain
   \begin{equation}  \label{taufine}
   \tau(\mathbf{b;e}) = \tau^{\flat}(\mathbf{b}) + \tau^{\sharp}(\mathbf{e}) = \frac{191}{49} +  \frac{27}{80}
   = \frac{16603}{3920}  = 4.235459\ldots,
   \end{equation}
   arriving at the number $ \frac{191}{49} +  \frac{27}{80} = \frac{16603}{3920}$ in~\eqref{lessthan14}.

We can now connect to our holonomy bounds to prove Theorem~\ref{mainA}. 
By Proposition~\ref{functionsH} and Lemma~\ref{14functions}, we have a set of~$m=14$ (holonomic) functions linearly independent over~$\Q(y)$ 
that are in $\Q\llbracket y \rrbracket$ with the denominator types~\eqref{typesA},
 contingent upon the  $\Q$-linear dependency~\eqref{absurd}. Hence it suffices to prove that
  (any one of) our holonomy bounds yields~$m < 14$: this will refute~\eqref{absurd}. 
    
  For example:

\begin{proof}[Proof via Theorem~\ref{main:BC form}]
Applying Theorem~\ref{main:BC form}, we obtain the upper bound $m \leq 13.9938\ldots$,
as computed in~\eqref{lessthan14}.
\end{proof}

\begin{proof}[Proof via Theorem~\ref{main:elementary form}]
Apply Theorem~\ref{main:elementary form} with $l=1, r_0=e^{-1/2}, \gamma_1=14\cdot 0.209=2.926$;
we pick this particular parameter based on the numerics in Example~\ref{Ex_BCfull}. In this case, we have
\[ \int_0^1  2t \cdot g_{\boldsymbol{\varphi, \gamma}}^*(t) \, dt = 11.316,\ldots\]
and thus the holonomy bound reads
\[
m \leq \frac{
\displaystyle{\raisebox{0.5ex}{$\displaystyle{11.316\ldots + \frac{1}{14} \cdot 2.926^2 \cdot \frac{1}{2}}$}\vphantom{\frac{1}{2}}}
}
{\displaystyle{\log \left( 256  \cdot
\frac{5448339453535586608000000000}{8658833407565631122430056127}
\right)
- \left(\frac{27}{80}  + \frac{191}{49}\right)}} = 13.730\ldots < 14. \qedhere
\]
\end{proof}

\begin{proof}[Proof via Theorem~\ref{main:BC conv discrete}]
With the choice of parameters as in Example~\ref{BCconv-pfA} with~$r_0 = e^{-1/2}$ and~$r_1= 1$,
we obtain the bound (see equation~\ref{thirdbest})
\[m \le 13.7206 \ldots < 14.  \qedhere\] 
\end{proof}

 Of course we may also apply other holonomicity bounds in \S~\ref{new slopes}. See Example~\ref{BCconv-pfA}
 with four parameters~$r_i$ rather than two,  and Examples~\ref{Ex-easyconv} and~\ref{Ex_BCfull}. 

\begin{remark}  \label{int is necessary} 
Had we stayed in the cruder framework $\mathbf{e=0}$ of Theorem~\ref{basic main} without added integrals, we would have had
to augment $\mathbf{b}$ to the array
$$
\mathbf{b}' := \left(  \begin{array}{llllllllllllll}   0 & 2 & 2 & 2 & 2 & 2 & 2 & 2 & 2 & 2 & 2 & 2 & 2 & 2 \\
0 & 0 & 0 & 2 & 2 & 2 & 2 & 2 & 2 & 2 & 2 & 2 & 2 & 2 \\
0 & 0 & 1 & 0 & 0 & 0 & 0 & 0 & 0 & 1 & 1 & 1 & 1 & 1    \end{array} \right)^{\mathrm{t}}.
$$ 
Note that although this $\bb'$ is not of the particular form in Theorem~\ref{main:elementary form}, we could apply the more general denominator formula~\eqref{tau bar} in Theorem~\ref{high dim BC convexity}. However, for the sake of simplicity, we give a good enough estimate of $\ovtau(\bb')$ which is sufficient to illustrate the necessity of working with nonzero $\be$.
On one hand, the upper bound argument in Remark~\ref{equalityoftaus} applies to those $\bb'$ such that every column in $\bb'$ has two values including one of which is $0$.
Therefore we have $\ovtau(\mathbf{b}') \leq (2+2+1) - \frac{1} {14^2} (1^2 \cdot 2 + 3^2 \cdot 2 + 8^2 \cdot 1)=\frac{32}{7}= 4.571\ldots$.

On the other hand, by the definition of~\eqref{tau bar}, we have an easy lower bound by only considering $\bn$ satisfying that $n_{j_1} > n_{j_2}$ implies $i_{j_1} \geq i_{j_2}$ and then we have~$\ovtau(\mathbf{b}')$
is at least
{\small
\begin{equation*}
\begin{aligned}
&  \frac{1\cdot 0+ 3\cdot 2 + 5\cdot (2+1) + (7 + 9 + \ldots + 17)\cdot (2+2) + (19+ 21 + 23 + 25 + 27)\cdot (2+2+1)}{14^2}\\
& \quad =\frac{884}{196}= 4.510\ldots,
\end{aligned}
\end{equation*}
}
a significantly worse value than~\eqref{taufine}.~\endofremark
\end{remark}

\section{Products of two logarithms} \label{sec:logs}

In this section, we apply our methods to certain products of logarithms.
Baker's theorem~\cite{Baker} gives a definitive result for linear forms in logarithms, even over~$\Qbar$, but we still do not
know how to show that~$\log 2 \cdot \log 3$ or~$\pi  \cdot \log 2$ is irrational.
While our methods cannot (as yet!) handle those cases either, we do prove Theorem~\ref{logsmain}, which we recall again here:

\begin{thm} \label{logsmainagain}
Let~$m, n \in \Z \setminus \{-1,0\}$ be  integers such that
$\displaystyle{\left| \frac{m}{n} - 1 \right| < \frac{1}{10^6}}$.
Then
\begin{equation} \label{logresult}
\log \left(1 + \frac{1}{m} \right) \log \left( 1 + \frac{1}{n} \right)
\end{equation}
is irrational.
Moreover,
for~$m \ne n$, the following are linearly independent over~$\Q$:
\begin{equation}
\label{fourlogs}
1, \quad \log \left(1 + \frac{1}{m} \right), \quad  \log \left( 1 + \frac{1}{n} \right), \quad
\log \left(1 + \frac{1}{m} \right) \log \left( 1 + \frac{1}{n} \right). 
\end{equation}
\end{thm}

\begin{remark} We could certainly improve the constant~$10^{-6}$ by
our methods, but some computation suggests that it is unlikely one could do
better than (say)~$10^{-4}$, and most likely not even that far; we make
this choice of constant for its relative simplicity.~\endofremark
\end{remark}

The degenerate case of~$m=n$ is a trivial consequence of the transcendence of~$\log r$ for~$r > 0$ in~$\Q \setminus \{1\}$, and so we shall assume that~$m \ne n$.
We begin by recalling a proof of the irrationality of
\begin{equation}
\label{singlelog}
\log \left(1 + \frac{1}{m} \right)
\end{equation}
for~$m \ge 1$ from~\cite{AlladiRobinsonTwo,AlladiRobinson,ChudnovskyHermite,Apery,vdP2}), based on the method of Ap\'ery limits. It is closely related to the construction we recounted in~Basic Remark~\ref{log 3 example}, and also to the Hermite--Pad\'e construction in~\S~\ref{log discussion} for the logarithm function. 

 Let~$a>1$ be an integer. The function
\begin{equation}
\label{beukersA}
A(a,x) := \frac{1}{\sqrt{1 - 2 a x + x^2}}  = \sum_{n=0}^{\infty} u_n(a) x^n 
\end{equation}
lies in~$\Z \llbracket x \rrbracket$ if~$a$ is odd and in~$\Z \llbracket x/2 \rrbracket$ otherwise,
and satisfies the first order ODE
\begin{equation}
\label{odeone}
(1 - 2 a x + x^2) y' + (x - a) y = 0.
\end{equation}
There is a unique solution to the non-homogenous ODE
$$(1 - 2 a x + x^2) y' + (x - a) y = 1$$
with coefficients in~$\Q$ that  is holomorphic and vanishes at~$0$; it is given by
\begin{equation}
\label{beukersB}
\begin{aligned}
H(a,x) := & \  \frac{1}{\sqrt{1 - 2 a x + x^2}} \int_{0}^{x} \frac{dt}{\sqrt{1 - 2 a t + t^2}}  \\ 
= & \ \frac{1}{\sqrt{1 - 2 a x + x^2}} \left( \log\left(a - x - \sqrt{1 - 2 a x + x^2} \right) - \log(a - 1) \right) \\
= & \ x +  \frac{a x^2}{2} + \frac{(3 a^2 - 1) x^3}{6} + \ldots = \sum_{n=0}^{\infty} v_n(a) x^n \in \Q \llbracket x \rrbracket, \end{aligned}
\end{equation}
and moreover the coefficients~$v_n(a)$ satisfy~$[1,2,\ldots,n] v_n(a) \in \Z$ if~$a$ is odd and
satisfy $[1,2,\ldots,n] 2^n v_n(a) \in \Z$ otherwise. 
By~\eqref{beukersB}, we have the formula
$$H(a,x) - \frac{1}{2} \log \left( \frac{a+1}{a-1} \right)  A(a,x) = \frac{1}{\sqrt{1 - 2 a x + x^2}} \int_{a-\sqrt{a^2-1}}^{x} \frac{dt}{\sqrt{1 - 2 a t + t^2}},$$
whose right-hand side
 overconverges at the singularity~$x = a - \sqrt{a^2 - 1}$ due to multiplying~$(-1)$ monodromies of both factors after an analytic continuation along a simple loop enclosing that singularity. This is the same mechanism for overconvergence as in~\S~\ref{mixed examples}, as well as in~\S~\ref{X06} with the canceling automorphy weights in the Eisenstein series~$A$ and the Eichler integral~$B - \frac{1}{2}L(2,\chip)$. The case at hand is readily seen to be equivalent, upon notational changes, to the respective ODEs~\eqref{HermitePadelogODE} and formulas~\eqref{hermitepadelog} arising from the diagonal Hermite--Pad\'e table for the logarithm function, which 
we recounted in~\S~\ref{log discussion} and~\S~\ref{ARB}.

 It follows that
\begin{equation}
\label{beukerslimit}
\lim_{n \rightarrow \infty} \frac{v_n(a)}{u_n(a)} \rightarrow \frac{1}{2} \log \left( \frac{a+1}{a-1} \right) 
\end{equation}
sufficiently quickly to prove the irrationality of this quantity for any odd~$a \ge 3$ or any even~$a \ge 4$ in light of the inequalities
$$5.828 \ldots = 3 + 2 \sqrt{2} > e = 2.718\ldots,$$
$$7.872 \ldots = 4 + \sqrt{15} > 2 \cdot e = 5.43656 \ldots $$
If we let~$a = 1+2m$, then 
$$ \log \left( \frac{a+1}{a-1} \right) = \log \left(1 + \frac{1}{m} \right),$$
giving the irrationality of~(\ref{singlelog}), as promised (with a pretty decent irrationality measure, improved further by Chudnovsky~\cite{ChudnovskyHermite,ChudnovskyThueSiegel} by a closer study of this argument).

Now let us consider the  arithmetic of the quantities
\begin{equation}
\label{prodlogs}
 \log \left( \frac{a+1}{a-1} \right)   \log \left( \frac{b+1}{b-1} \right) 
 \end{equation}
for pairs of integers~$a \ne b$. 
From~(\ref{beukerslimit}), it is obvious that
\begin{equation}
\lim_{n \rightarrow \infty} \frac{v_n(a)}{u_n(a)} \cdot  \frac{v_n(b)}{u_n(b)} 
\rightarrow \frac{1}{4} \log \left( \frac{a+1}{a-1} \right) \log \left( \frac{b+1}{b-1} \right).
\end{equation}
However, this certainly does not converge fast enough to prove irrationality of the right-hand side through any elementary analysis.
We shall nevertheless see that as long as~$a/b$ is sufficiently close to~$1$, this quantity is approachable via our new methods based on the function-theoretic properties of the generating series themselves, and basic properties of the Hadamard product operation which allows to construct new~$G$-functions with the desired Ap\'ery limit.

For each pair of integers~$a,b \in \Z \setminus \{-1,0,1\}$, let us write
$$\begin{aligned}
 \eta_{a} & := \frac{1}{2} \log \left( \frac{a+1}{a-1} \right), \quad
\eta_{b} := \frac{1}{2} \log \left( \frac{b+1}{b-1} \right), \\
\eta_{a,b} & := \eta_{a} \eta_{b} := \frac{1}{4}  \log \left( \frac{a+1}{a-1} \right)  \log \left( \frac{b+1}{b-1} \right).
\end{aligned}$$
We shall assume that~$a \ne \pm b$, as the irrationality and, indeed, the transcendence of~$\eta_{a,a} = -\eta_{a,-a} = \eta_a^2$
is already known. Our approach to the arithmetic properties of the product of the Ap\'ery limits~$\eta_{a,b} = \eta_a\eta_b$ is via the 
Hadamard product of the underlying~$G$-functions. 

\subsection{Hadamard products and Ap\'ery limits}  \label{sec:Hadamard}
Let~$\star$ denote the Hadamard product operation on power series: $\left( \sum a_n x^n \right) \star \left( \sum b_n x^n \right) := \sum a_nb_n x^n$.
The function
\begin{equation}
\label{defPA}
P_A(x):=A(a,x) \star A(b,x) = \sum u_n(a) u_n(b) x^n 
\end{equation}
satisfies the following ODE~$\MMM(P_A(x)) = 0$:
\begin{equation}
\label{ODE}
\begin{aligned}
& \quad  (-1+x)x(1+x)(1 -4 a b x - 2 x^2 + 4 a^2 x^2 + 4 b^2 x^2 - 4 a b x^3 + x^4) y'' \\
&  + (-1 + 8a b x + 5 x^2 - 12 a^2 x^2 - 12 b^2 x^2 + 16 a b x^3 - 7 x^4 + 4 a^2 x^4 + 4 b^2 x^4
-8 a b x^5 + 3 x^6)y' \\
 & + (a b - (-1 + 3 a^2 + 3 b^2) x + 8 a b x^2 - (2 + a^2 + b^2) x^3 - a b x^4 + x^5) y = 0. \end{aligned}
\end{equation}
The points~$x = 1$ and~$x = -1$ are 
only apparent singularities, as long as~$(a-b) \ne 0$ and~$(a+b) \ne 0$ respectively.
This follows both by general properties~\cite{HadamardProd} of the Hadamard product but can also be verified directly by computing
the indicial equation (which is~$R(R-2) = 0$), and then verifying that there are two linearly
independent power series solutions. For example, for a putative 
solution~$\sum c_n (x+1)^n$, the coefficients~$c_n$ satisfy a recurrence of the form
$$\begin{aligned} (a+b)^2 n(n-2) c_n = & \  12(a+b)^2 c_{n-1} - 4 (a+b)^2(n-2) c_{n-1}(7n-9)(n-2)\\
  & \  +  6(a+b)^2 c_{n-2} + (n-2) c_{n-2}(\ldots) + \ldots, \end{aligned}$$
which implies that there is at least one solution of the form~$x^2 + \ldots$, but there is another of the form:
$$1 + \frac{(x+1)}{2}  - \frac{(x+1)^3}{4} + \ldots \in \Q \llbracket a,b,x+1 \rrbracket,$$
and the case of~$x = 1$ is similar.
The four roots of the quartic are exactly the products of the singularities of the order one ODEs,
namely
$$(a \pm  \sqrt{a^2 - 1}) \times (b \pm \sqrt{b^2 - 1}).$$

\begin{basicremark}
If~$a$ and~$b$ are large and~$a/b$ is very close to one, then the four nonzero finite singularities of~(\ref{ODE})
are grouped as follows:
\begin{enumerate}
\item One singularity~$\alpha$ very close to~$0$.
\item Two singularities very close to~$1$.
\item One very large singularity (``close to~$\infty$'').
\end{enumerate}
In addition to these,~$0$ and~$\infty$ themselves are also (essential) singularities. 
We shall construct (assuming a linear relationship over~$\Q$ between the quantities~(\ref{fourlogs})) a function~$H \in \Q\llbracket x \rrbracket$
with denominator type~$\tau = [1,2,3,\ldots,n]^2$, satisfying a non-homogenous version of~(\ref{ODE}), and overconvergent beyond~$\alpha$.
When considering functions on~$\mathbf{P}^1 \setminus \{0,1,\infty\}$ with~$\tau = [1,2,\ldots,n]^2$, we are required by~\S~\ref{the lambda template} to choose an auxiliary function~$\varphi: \D \rightarrow \C \setminus \{1\}$
with~$\varphi(0) = 0$ and~$\varphi^{-1}(0) = 0$. If we restrict~$\varphi$ to the disc~$D(0,1 - \varepsilon)$ for any~$\epsilon \in (0,1/2]$, then the image of~$\varphi$ will
avoid a small open ball containing~$1$ and an open ball containing~$\infty$, and satisfy~$\# \varphi^{-1}(\alpha) = 1$. This is the type of setting where our holonomy bounds can be applied,
for we can include not only the (presumably non-existent!) functions~$H(x)$ and their derivatives, but also the pure functions on~$\mathbf{P}^1 \setminus \{0,1,\infty\}$
the we devised in~\S~\ref{sec:purefunctions}.
As a practical matter, our maps~$\varphi$ are of the form~$\lambda \circ \GGG$ or, in the equivalent~$\mathbf{P}^1 \setminus \{0,4,\infty\}$ setting,~$h \circ \GGG$ for some map~$\GGG : \Db \to \D$, which we take as the Riemann map of a suitably chosen domain in~$\D$.
But even though the function~$\lambda$ for example avoids~$1$ on~$\D$, to avoid the values within~$\varepsilon$ of~$1$ requires taking~$\GGG$ to have significantly smaller
conformal radius unless~$\varepsilon$ is extremely small. For example, if the image of~$\GGG(\D)$ inside~$\D$ included the point~$3/4$, then~$\varphi = \lambda \circ \GGG$ on~$\D$
would already include the value
$$\lambda(3/4) = 0.9999999999999798332 \ldots $$
This numerology is ultimately what forces the hypothesis that~$|m/n - 1|$ is very small. 
\endofremark
\end{basicremark}

\subsubsection{The overconvergent space}  \label{sec:ovc}
We now consider solutions to a non-homogenous version of the product ODE~(\ref{ODE}).
Denoting the corresponding differential operator of~(\ref{ODE}) by~$\MMM$, then we also have the following identities:
$$\begin{aligned}
\MMM(H(a,x) \star  A(b,x)) =  & \  - b + 3 a x - 3 b x^2 + a x^3, \\
\MMM(A(a,x) \star H(b,x)) = & \  -a + 3 b x - 3 a x^2 + b x^3, \\
\MMM (H(a,x) \star H(b,x)) =  & \  -1 + 2 x^2 - x^4. \end{aligned}
$$
We further claim that  the functions:
\begin{equation} \label{Pabdef}
\begin{aligned}
P_{a}  :=  & \ \left( H(a,x) - \eta_{a}A(a,x) \right) \star A(b,x)  =  \sum (v_n(a) -  \eta_{a}  u_n(a)) u_n(b) x^n, \\
P_{b}  := & \ A(a,x) \star \left( H(b,x) - \eta_{b}    A(b,x) \right)   =  \sum u_n(a) (v_n(b) - \eta_{b} a_n(b))  x^n, \\
P_{ab}  := & \ H(a,x) \star H(b,x) - \eta_{a,b}  A(a,x) \star  A(b,x)  \\ 
  = & \  \sum \left( v_n(a) v_n(b) - \eta_{a,b} u_n(a) u_n(b) \right)  x^n \\
     = & \  \sum \big\{ \eta_a u_n(a) \left(v_n(b) - \eta_b u_n(b) \right)    + \eta_b  \left(v_n(a) - \eta_a u_n(a) \right) u_n(b)  \big\}  x^n \\
 + & \  \sum  \left(v_n(a)   - \eta_a u_n(a) \right) (v_n(b) - \eta_b u_n(b)) x^n   
 \end{aligned}
 \end{equation} 
are overconvergent beyond the smallest cusp~$(a-\sqrt{a^2-1})(b-\sqrt{b^2-1})$. This follows from the bounds
$$|v_n(a)  - \eta_a u_n(a)| = O((a - \sqrt{a^2-1})^{n(1 - \epsilon)}),$$
$$|v_n(b)  - \eta_b u_n(b)| = O((b - \sqrt{b^2-1})^{n(1 - \epsilon)}),$$
together with
$$\begin{aligned} |v_n(a)|, |u_n(a)| & = O((a + \sqrt{a^2-1})^{n(1 +\epsilon)}), \\
|v_n(b)|, |u_n(b)| & = O((b + \sqrt{b^2-1})^{n(1 +\epsilon)}). \end{aligned}$$
Another way to express the power series decomposition of~$P_{ab}$ in the last two lines of~\eqref{Pabdef} 
is as follows:
$$
\begin{aligned}
P_{ab} & =  \eta_a  A(a,x) \star \left( H(b,x) - \eta_b A(b,x) \right) + \eta_b \left( H(a,x) - \eta_a A(a,x) \right) \ast A(b,x) \\
& + 
\left( H(a,x) - \eta_a A(a,x) \right) \star \left( H(b,x) - \eta_b A(b,x) \right) 
\end{aligned}
$$
The general fact~\cite{HadamardProd} that we exploited here is the overconvergence of the Hadamard product of any set of holonomic power series, at least one among which is an overconvergent branch in the sense of~\S~\ref{local univalent leaves}. This was essentially combined with the Jacobson  identity~$1-xy = (1-x) + (1-y) - (1-x)(1-y)$, familiar for example from the proof of the nilpotence of the augmentation ideal of the~$\F_p$-group ring of
a finite~$p$-group. 

\subsubsection{Construction of the unlikely~$G$-function}
Assume now and until the end of the proof of Theorem~\ref{logsmainagain} that there exist integers~$r_0$, $r_a$, $r_b$, and~$r_{ab}$ not all zero such that
\begin{equation}
\label{rabc}
r_a  \eta_{a} + r_b \eta_{b} + r_{ab} \eta_{a,b} = r_0.
\end{equation}
Then the linear combination
\begin{equation}
\label{defP}
\begin{aligned}
P :=  & \ r_a P_a + r_b P_b + r_{ab} P_{ab} \\
= & \   r_a H(a,x) \star A(b,x) + r_b A(a,x) \star H(b,x)  \\
& \quad + r_{ab} H(a,x) \star H(b,x) + r_0 A(a,x) \star A(b,x) \\
= & \ \sum c_n(a,b) x^n
\in \Q \llbracket x\rrbracket
 \end{aligned}
\end{equation}
is also from the overconvergent space~\S~\ref{sec:ovc}, but now it has rational coefficients. 
This is the~$G$-function, contingent upon our absurd hypothesis of a linear dependency~\eqref{rabc}, that will ultimately be rejected by our holonomy bounds. 
The analytic properties of this unlikely function follow from the overconvergence in~\S~\ref{sec:ovc}; we now collect the arithmetic properties. 
If both~$a$ and~$b$ are odd, then~$u_n(a), u_n(b) \in \Z$ and moreover both  $[1,2,\ldots,n] v_n(a)$
 and $[1,2,\ldots,n] v_n(b)$ lie in~$\Z$. 
Therefore, in the Hadamard products construction, 
$$[1,2,\ldots,n]^2 c_n(a,b) \in \Z.$$
Moreover, with~$a=1+2m$ and~$b=1+2n$, the non-existence of a linear
relationship~(\ref{rabc}) is exactly the thesis of Theorem~\ref{logsmain}.
(The conditions $a,b \in \Z \setminus \{-1,0,1\}$ become $m,n \in \Z \setminus \{-1,0\}$.)
Thus to prove Theorem~\ref{logsmain} it suffices to assume the existence
of a relationship~(\ref{rabc}) and a function~$P(x)$ as in~(\ref{defP}), and
establish a contradiction.

\begin{df}  \label{PtoG} With~$r_a$, $r_b$, and~$r_{ab}$ satisfying~(\ref{rabc}), let~$P(x)$
be defined as in equation~(\ref{defP}), and let
$$G(y) := P(x) + P \left( \frac{x}{x-1} \right) \in \Q \llbracket y \rrbracket.$$
With~$P_A(x) = A(a,x) \star A(b,x)$ as in~(\ref{defPA}), let
$$G_A(y) := P_A(x) + P_A \left(\frac{x}{x-1} \right),$$
hence~$G_A(y) \in \Z \llbracket y \rrbracket$ if~$a,b$ are odd, and~$G_A(y) \in \Z[1/2] \llbracket y \rrbracket$ otherwise.  \endofremark
\end{df}

As in the proof of Theorem~\ref{mainA}, we will work with the~$Y_0(2)$ picture in the dictionary of~\S~\ref{sec:YtoY0(2)}, 
and refute the existence of a~$G$-function~$G \in \Q\llbracket y \rrbracket$ (contingent upon the existence of a relation~\eqref{rabc}), of the denominators type~$[1,\ldots,2n]^2$ and ``close to'' the
$\P^1 \setminus \{0, 4, \infty\}$ type that we studied in~\S~\ref{sec:purefunctions}.
We record the following properties of~$G$ in Proposition~\ref{fictive G} below, after the following definitions:

\begin{df} \label{yplus}
Let~$y_{a^{\pm},b^{\pm}}$ denote the
\(\displaystyle{y := \pi(x) = x + \frac{x}{x-1} }\)
images of
\[x = \left(a \pm \sqrt{a^2-1} \right) \left( b \pm \sqrt{b^2-1} \right),\]
 where all four pairs of signs are being considered. 
 Let~$\LL$ denote the pushforward of~$\MMM$ under~$\pi$, so that~$\LL(G_A(y)) = 0$.
  \endofremark
\end{df}

Recall that our present discussion is conditional on supposing a~$\Z$-linear relation~\eqref{rabc}. 
At this point, we make the additional assumption that the integers~$a,b \in \Z \setminus \{\pm 1\}$ are odd. 

\begin{proposition} \label{fictive G}
With~$a, b \in \Z \setminus \{ \pm 1\}$ odd,
the functions $G(y) \in  \Q\llbracket y \rrbracket$ and~$G_A(y)  \in \Q\llbracket y \rrbracket$  of~\S~\ref{PtoG} have
 denominator types~$[1,\ldots,2n]^2$ and~$1$, respectively.
 Moreover,
 \(\LL(G_A(y)) = 0\) and \(\LL(G(y)) \in \Q[y]\) for some non-zero linear differential operator~$\LL$ over~$\Q(y)$ satisfying: 
\begin{enumerate}
\item $\LL$ has no singularities
besides
$y \in \left\{0, 4, y_{a^{\pm},b^{\pm}}, \infty\right\}$.
\item $\LL$ has~$\Z/2$ local monodromy around the singularity~$y=4$.
\end{enumerate} 
\end{proposition}

We shall write down~$\LL$ explicitly in~\S~\ref{sec:diffGA} below;  the exact form
of the polynomial~$\LL(G(y)) \in \Q[a,b,y]$ can be computed but will not be important.

\subsection{The differential equation  \texorpdfstring{$\LL(G_A) = 0$}{LL}}
\label{sec:diffGA}
Before giving the statement and proof of Lemma~\ref{17functions} (the analog of 
of Lemma~\ref{14functions}), we shall  examine the ODE~$\LL(G_A)=0$ in more detail.
Because of the length of this computation, it makes more sense to present it separately
rather than interweave it with the proof of Lemma~\ref{17functions}, compared to the corresponding
facts concerning the Zagier functions which are proved during the proof of Lemma~\ref{14functions}.
However, the reader may well want to look ahead to the statement of Lemma~\ref{17functions} to
see where we are going.
One  can compute from~(\ref{ODE}) the following explicit form of~$\LL$:
 \[\displaystyle{\LL(G_A) = \sum_{i=0}^{4} c_i(y) G_A^{(i)}(y) = 0},\]
where~$c_i(y)$ are certain polynomials with
\begin{equation} \label{coefficientsofL}
\begin{aligned}
c_0(y) = & \ R_{12}(y) \\
c_1(y) = & \  R_{16,A}(y) \\
c_2(y) = & \ y R_{16,B}(y) \\
c_3(y) = & \ y^2 (y-4) R_{15}(y) \\
c_4(y) =  & \ (y-4)^2 y^3  R_4(y) R_{10}(y).
\end{aligned}
\end{equation}
Here~$R_d(y)$ denotes an irreducible polynomial of degree~$d$ (with respect to $y$) in~$\Q(a,b,y)$,
and the subscripts~$A$ and~$B$ denote that~$R_{16,A}$ and~$R_{16,B}$ are distinct.
The polynomial
$$(1-x)^4 R_4 \left( x + \frac{x}{x-1} \right)$$
has~$(a \pm \sqrt{a^2-1})(b \pm \sqrt{b^2-1})$ as~$4$ of
its~$8$ roots, together with the
images of these roots under the involution~$w(x) = x/(x-1)$.
In particular, the roots of~$R_4(y)$ are  the singularities
$y_{a^{\pm},b^{\pm}}$ of~$\LL$. 
The polynomial~$R_4(y)$
is given explicitly by
$$
\begin{aligned}
R_{4}(y) = & \ 1+4 y - 8 a^2 y - 4 a b y - 8 b^2 y + 16 a^2 b^2 y + 4 y^2 - 
  12 a^2 y^2 + 16 a^4 y^2  \\
  & \ + 20 a b y^2 - 16 a^3 b y^2 - 12 b^2 y^2 - 
  16 a b^3 y^2    + 16 b^4 y^2 - 8 a^2 y^3 + 16 a b y^3  \\
  & - 16 a^3 b y^3 - 
  8 b^2 y^3   + 32 a^2 b^2 y^3 - 16 a b^3 y^3 + 4 a^2 y^4 - 8 a b y^4 + 
  4 b^2 y^4, \end{aligned}
  $$
  Unlike with~$R_4(y)$ or the other accompanying powers of~$y$ and~$(y-4)$
  appearing in~$c_4(y)$ of~\eqref{coefficientsofL}, the roots of~$R_{10}(y)$
  are not genuine singularities of~$\LL$.  More precisely,  this is true if the roots
  of~$R_{10}(y)$ are distinct from those of~$R_4(y) y(y-4)$,
  and this will hold under our assumptions by
  Lemma~\ref{det} and Lemma~\ref{res} proved below.
The polynomial~$R_{10}(y)$ is given explicitly by
  {\small
  $$
  \begin{aligned}
 & \
 4 a^2 b^2 + 9 y - 27 a^2 y - 102 a b y + 144 a^3 b y 
 - 
 27 b^2 y + 207 a^2 b^2 y - 128 a^4 b^2 y \\
  & \  + 144 a b^3 y - 
 160 a^3 b^3 y - 128 a^2 b^4 y + 64 y^2 - 93 a^2 y^2 + 
 228 a^4 y^2 - 284 a b y^2 + 740 a^3 b y^2 \\
  & \  - 
 992 a^5 b y^2 - 93 b^2 y^2 - 818 a^2 b^2 y^2 - 
 16 a^4 b^2 y^2 + 740 a b^3 y^2 + 1208 a^3 b^3 y^2 + 
 228 b^4 y^2 - 16 a^2 b^4 y^2 \\
  & \  + 64 a^4 b^4 y^2 - 
 992 a b^5 y^2 + 164 y^3 - 132 a^2 y^3 - 804 a^4 y^3 + 
 432 a^6 y^3 - 100 a b y^3 + 1218 a^3 b y^3 \\
  & \   + 
 2360 a^5 b y^3 - 1536 a^7 b y^3 - 132 b^2 y^3 - 
 1108 a^2 b^2 y^3 - 856 a^4 b^2 y^3 + 1856 a^6 b^2 y^3 + 
 1218 a b^3 y^3 \\
 & \   - 3120 a^3 b^3 y^3 
   - 448 a^5 b^3 y^3 - 
 804 b^4 y^3 - 856 a^2 b^4 y^3 - 16 a^4 b^4 y^3 + 
 2360 a b^5 y^3 - 448 a^3 b^5 y^3 + 432 b^6 y^3 \\
  & \ + 
 1856 a^2 b^6 y^3
- 1536 a b^7 y^3 
+ 168 y^4 - 
 18 a^2 y^4 - 1476 a^4 y^4 + 432 a^6 y^4 - 108 a b y^4 - 
 840 a^3 b y^4 \\
  & \  + 688 a^5 b y^4 + 384 a^7 b y^4 - 
 18 b^2 y^4 + 4384 a^2 b^2 y^4 - 8864 a^4 b^2 y^4 + 
 4384 a^6 b^2 y^4 - 840 a b^3 y^4 \\
  & \  + 15744 a^3 b^3 y^4 - 
 11744 a^5 b^3 y^4 - 1476 b^4 y^4 - 8864 a^2 b^4 y^4 + 
 13920 a^4 b^4 y^4 + 688 a b^5 y^4 - 11744 a^3 b^5 y^4  \\
  & \  + 
 432 b^6 y^4 + 4384 a^2 b^6 y^4 + 384 a b^7 y^4 + 32 y^5 + 
 180 a^2 y^5 + 2467 a^4 y^5 + 792 a^6 y^5 
  - 720 a^8 y^5 \\ 
    & \  - 
 424 a b y^5 - 4228 a^3 b y^5 + 440 a^5 b y^5 + 
 1824 a^7 b y^5 + 180 b^2 y^5 
 + 3554 a^2 b^2 y^5 \\
   & \   - 
 15928 a^4 b^2 y^5 
 + 5104 a^6 b^2 y^5 - 4228 a b^3 y^5 + 
 29392 a^3 b^3 y^5 - 25088 a^5 b^3 y^5 + 2467 b^4 y^5 \\
    & \  - 
 15928 a^2 b^4 y^5 + 37760 a^4 b^4 y^5 
 + 440 a b^5 y^5 - 
 25088 a^3 b^5 y^5 + 792 b^6 y^5 + 5104 a^2 b^6 y^5 + 
 1824 a b^7 y^5   \\ 
   & \ 
 - 720 b^8 y^5 - 32 y^6 
  - 336 a^2 y^6 + 
 258 a^4 y^6 - 3840 a^6 y^6 + 480 a^8 y^6 + 736 a b y^6 + 
 1064 a^3 b y^6 \\
     & \  + 6680 a^5 b y^6 - 4320 a^7 b y^6 - 
 336 b^2 y^6 - 2676 a^2 b^2 y^6 + 6976 a^4 b^2 y^6 + 
 10688 a^6 b^2 y^6 + 1064 a b^3 y^6 \\
     & \  - 19632 a^3 b^3 y^6 - 
 12256 a^5 b^3 y^6 + 258 b^4 y^6 + 6976 a^2 b^4 y^6 + 
 10816 a^4 b^4 y^6 + 6680 a b^5 y^6 - 12256 a^3 b^5 y^6\\
     & \  - 
 3840 b^6 y^6 + 10688 a^2 b^6 y^6 - 4320 a b^7 y^6 + 
 480 b^8 y^6 - 576 a^2 y^7 - 1392 a^4 y^7 + 4368 a^6 y^7 + 
 1152 a b y^7 \\
     & \  + 5312 a^3 b y^7 - 12848 a^5 b y^7 + 
 576 a^7 b y^7 - 576 b^2 y^7 - 7840 a^2 b^2 y^7 + 
 12912 a^4 b^2 y^7 - 1104 a^6 b^2 y^7 \\
     & \  + 5312 a b^3 y^7 - 
 8864 a^3 b^3 y^7 - 768 a^5 b^3 y^7 - 1392 b^4 y^7 + 
 12912 a^2 b^4 y^7 + 2592 a^4 b^4 y^7 - 12848 a b^5 y^7 \\
     & \  - 
 768 a^3 b^5 y^7 + 4368 b^6 y^7 - 1104 a^2 b^6 y^7 + 
 576 a b^7 y^7 + 192 a^2 y^8 + 168 a^4 y^8 - 960 a^6 y^8 \\
     & \ - 
 384 a b y^8 - 864 a^3 b y^8 + 1568 a^5 b y^8 + 
 192 b^2 y^8 + 1392 a^2 b^2 y^8 + 2368 a^4 b^2 y^8 - 
 288 a^6 b^2 y^8 \\
     & \  - 864 a b^3 y^8 - 5952 a^3 b^3 y^8 + 
 1152 a^5 b^3 y^8 + 168 b^4 y^8 + 2368 a^2 b^4 y^8 - 
 1728 a^4 b^4 y^8 + 1568 a b^5 y^8 \\
     & \ + 1152 a^3 b^5 y^8 - 
 960 b^6 y^8 - 288 a^2 b^6 y^8 + 160 a^4 y^9 - 
 640 a^3 b y^9 + 448 a^5 b y^9 + 960 a^2 b^2 y^9 \\
     & \  - 
 1792 a^4 b^2 y^9 - 640 a b^3 y^9 + 2688 a^3 b^3 y^9 + 
 160 b^4 y^9 - 1792 a^2 b^4 y^9 + 448 a b^5 y^9 - 
 32 a^4 y^{10} \\
     & \  + 128 a^3 b y^{10} - 192 a^2 b^2 y^{10} + 
 128 a b^3 y^{10} - 32 b^4 y^{10}. 
  \end{aligned}
 $$
 }
We also find that
$$\frac{c_3(y)}{c_4(y)} 
= \frac{d}{dy} 
\log \left( \frac{ (y-4)^3 y^5 R_4(y)^3}
{R_{10}(y)} \right).$$
We compute that the discriminant
of~$R_{10}(y)$ has the form (up to an element of~$\Q^{\times}$):
\begin{equation}
\label{discriminant}
\begin{aligned}
\Delta_y(R_{10}(y)) =  & \ 
(a-b)^{12} (a+b)^6 
(1 + 4 a^2 - 4 a b)(1 + 4 b^2 - 2 a b)
\\
 & \  \times
 (-3 + 4 a^2 - 4 a b + 4 b^2)  \Phi_{14}(a,b) \Phi_{79}(a,b),
 \end{aligned}
\end{equation}
where~$\Phi_{d}(a,b)\in \Q(a,b)$ is irreducible,
satisfies~$\Phi_d(a,b) = \Phi_d(b,a)$,
and is of degree~$d$ when considered as a univariate
polynomial in either~$a$ or~$b$.
We also compute the resultant~$\mathrm{Res}_y(R_4(y),R_{10}(y))$
to be, up to a non-zero rational scalar; equal to
\begin{equation}
\label{resultant}
\begin{aligned}
 & \
(a-b)^8 (a+b)^6
(1 + 4 a^2 - 4 a b)(1 + 4 b^2 - 2 a b) \\
\ & \times (-3 + 4 a^2 - 4 a b + 4 b^2)^2 
(9 + 16 a^2 - 40 a b + 16 b^2)
  \Phi_{26}(a,b).
\end{aligned}
\end{equation}
It is easy to verify 
(reduce modulo~$2$)
that none of the quadratic factors vanish for
integer~$a,b \in \Z$ and
any of the~$\Phi$ above. One strongly suspects
that there are no other integral solutions to~$\Phi_d(a,b)=0$
for the other~$d$
except for certain degenerate solutions for some of
these polynomials when~$a=b$ or~$a=-b$. The general
Siegel theorem~\cite[\S~II.I]{Siegel1929SNS}, see also~\cite[Thm.~7.3.9]{BombieriGubler}, certainly guarantees that every irreducible nonrational affine
algebraic curve has at most a finite number of integral points; and here, as each of the polynomials $\Phi_{14}, \Phi_{26}$, and $\Phi_{79}$
turns out to have its highest degree homogeneous piece divisible by $ab(a-b)$ (the degrees of these polynomials 
are, respectively, $18, 34$, and~$92$), and hence is not proportional to a power of an irreducible polynomial over~$\Q$, Runge's method~\cite[\S~9.6.5]{BombieriGubler}
(see also~\cite[\S~4]{MasserBook} for a gentle and practical introduction)
provides in principle an exhaustive algorithm to enumerate all the integer solutions of these equations. 
For our purposes here, since we are studying the pairs $(a,b)$ with $|a| \asymp |b|$,   we shall exploit
this hypothesis in the sequel as it spares us the routine but grueling task of carrying out
these computations.

  \begin{lemma} \label{res}
  Assume that~$a,b \in \Z \setminus \{1,0,-1\}$ with~$a \ne \pm b$
  satisfy one of the following inequalities:
  $$\left| \frac{a}{b} - 1 \right| < \frac{1}{2}, \qquad
  \left| \frac{a}{b} + 1 \right| < \frac{1}{2}.$$
  Then:
  \begin{enumerate}
  \item $R_4(y)$ is irreducible.
  \item $R_{10}(y)$ is co-prime to~$R_4(y)$.
  In particular, the resultant~(\ref{resultant}) is non-vanishing.
  \end{enumerate}
  \end{lemma}

  \begin{proof}  The roots 
  of
  $\displaystyle{R_4 \left(x + \frac{x}{x-1} \right)}$
include~$(a - \sqrt{a^2 - 1})(b - \sqrt{b^2 - 1})$
  as a root. Hence, if we show that~$\Q(\sqrt{a^2-1})$ and~$\Q(\sqrt{b^2-1})$
  are distinct non-trivial real quadratic fields, then~$R_4(y)$ is absolutely
  irreducible since it has at least one root of degree~$4$.
  The assumptions on~$a$ and~$b$ certainly  imply that~$a^2 - 1$ and~$b^2 - 1$
  are not squares, so~$\Q(\sqrt{a^2 - 1})$ and~$\Q(\sqrt{b^2 - 1})$
  are quadratic fields. 
  If they define the same field, then there exist 
 integers $D, X,Y \in \Z$
  with~$D$ squarefree such that~$(a^2-1) = X^2 D$ and~$(b^2 - 1) = Y^2 D$, and so~$(a,X)$ and~$(b,Y)$
  are solutions to the Pell equation~$u^2 - D v^2 = 1$. 
  Let us consider the case when~$a>b>0$, the proof
  applies in the other cases \emph{mutatis mutandis}.
    After checking the small cases explicitly, we may assume that~$b  \ge 8$ 
  (note that bounding~$b$ also bounds~$a$).
We deduce that,
for positive  algebraic integer unit~$\varepsilon  > 1$ in~$\Q(\sqrt{D})$,
there is an equality
 $$ (a + X \sqrt{D}) = \varepsilon (b + Y \sqrt{D}).$$
 The left hand side lies in the interval~$[2a-1,2a]$. The right hand side
 lies in the interval~$\varepsilon [2b-1,2b]$. Hence
 \begin{equation}
 \label{tocontra}
 1 < \varepsilon < \frac{2a}{2b-1} = \frac{a/b}{1 - \frac{1}{2b}} 
 \le \frac{3/2}{1 - 1/16} = \frac{16}{10}.
 \end{equation}
 On the other hand, any  unit~$\varepsilon > 1$ of a real quadratic field satisfies
 $$\varepsilon \ge  \frac{\sqrt{5}+1}{2} > \frac{16}{10},$$
 contradicting equation~(\ref{tocontra}). 
    The first claim follows.
    Now if~$R_{10}(y)$ has a common factor with~$R_4(y)$, it must be divisible
    by~$R_{4}(y)$. However, we may now synthetically divide one polynomial
    by the other and the four coefficients of the remaining polynomial
    of degree~$\le 3$ must all be zero. 
  But there are no such solutions in~$a$ and~$b$ to these four equations   ---  
  already taking the resultant of any two of them gives an explicit polynomial
  in~$a$ 
  with no integers roots in~$ \Z \setminus \{-1,0,1\}$.
  \end{proof}
  
  We also have the following slightly unpleasant calculus exercise:
  \begin{lemma} \label{det}
  Assume that~$a,b \in \Z \setminus \{1,0,-1\}$  
  satisfy:
  $$
  0 < \left| \frac{a}{b} - 1 \right| < \frac{1}{10^3}.
  $$
  Then~$R_{10}(y)$ is separable. Moreover,
  $ \Res_{y}(R_{10}(y),y(4-y))$ which equals
  $$\begin{aligned}
   48 a^2 b^2 (9 + 16 a^2 - 40 a b + 16 b^2) (-45 + 80 a^2 - 128 a b + 
   80 b^2) \Phi_4(a,b) \end{aligned}
   $$
   is non-vanishing.
   \end{lemma}

\begin{proof} 
Let~$\varepsilon = 1/1000$.
First let us consider~$\Delta_y(R_{10}(y))$ as a  polynomial where
the coefficient~$a$ varies while~$b \in \Z$ is fixed.
 From~(\ref{discriminant}), the only factors which could
possibly vanish for~$a \in \Z$ are~$\Phi_{14}(a,b)$ and~$\Phi_{79}(a,b)$.
We examine each of these cases in turn. 
Consider the case of~$\Phi_{14}(a,b)$, and with~$a$ and~$b$ of the same sign.
Let~$b = a(1+x)$, so~$|x| \le \varepsilon$.
Then
\begin{equation}
\label{14limit}
\begin{aligned}
\frac{\Phi_{14}(a,a(1+x))}{\Phi_{14}(a,a) x^2 a^4}
= & \  Q_4(x) + a^{-2} Q_2(x) + a^{-2} Q_0(x) (a x)^{-2} \\
 & +  \frac{ \Psi_{-1,12}(a)}{a^3 (a x) \Phi_{14}(a,a)}
 +  \sum_{i=0}^{12} x^i \frac{\Psi_{i,12}(a)}{a^4 \Phi_{14}(a,a))},
    \end{aligned}
\end{equation}
where~$\Psi_{i,12}$ for~$i=-1,\ldots,12$ are explicit polynomials in~$a$ of degree at most~$12$,
and~$Q_i(x)$ is an explicit polynomial in~$x$ with~$Q_i(0) \ne 0$.
Moreover, $Q_4(0) = 2$ and  is bounded below on the 
interval~$x \in [-1/1000,1/1000]$ by something only very slightly less than~$2$.
Now we exploit the fact that~$b \in \Z$ is an integer to deduce that~$ax \in \Z$,
and so~$|ax| \ge 1$. But assuming~$|ax| \ge 1$ and~$|x| \le 1/1000$,
 all the other terms in~(\ref{14limit}) are clearly of order~$O(a^{-2})$ with explicitly computable constants,
 and so with the na\"{\i}ve triangle inequality bound, the left-hand side does not vanish as soon
 as~$a$ is large enough.
 To be completely explicit, we find that, for~$|a| \ge 1000$,
$$\begin{aligned}
|Q_4(x)| \ge & \  1.984, \\
|Q_0(x)|, |Q_2(x)| \le & \  10^{-5}, \\
\left| \frac{\Psi_{-1,12}(a)}{a^3 \Phi_{14}(a,a)} \right| \le & \ 10^{-10}, \\
\left| \frac{\Psi_{i,12}(a)}{a^4 \Phi_{14}(a,a)} \right| \le & \ 10^{-10}, \ i = 0,\ldots,12 \\
\end{aligned}
$$
from which the non-vanishing of~$\Phi_{14}(a,b)$ comfortably follows
 from equation~(\ref{14limit}).
For~$|a| < 1000$, note that there are no integers~$b$ satisfying the assumed
inequalities on~$a$ and~$b$. Alternatively, for any integer~$|a| \le 1000$, one can
check that~$\Phi(a,b) = 0$ has no integer roots except for~$(a,b) = (1,1)$ and~$(-1,-1)$.
The argument for~$\Phi_{79}(a,b)$ is entirely similar. The analogue of~(\ref{14limit}) in this case is
\begin{equation}
\label{79limit}
\begin{aligned}
\frac{\Phi_{79}(a,a(1+x))}
{
a^{24} x^{12} \Phi_{79}(a,a)
}  & =   \ 
 Q_{24}(x) \\
 +  & \sum_{i = 0}^{5} a^{-2-2i} Q_{22-4i}(x) (a x)^{-2i}
  + a^{-2-2i} Q_{20-4i}(x) (a x)^{-2i-2} \\
   + &  \sum_{i=1}^{11} \frac{ \Psi_{-i,68}(a)  }{a^{24 -i} \Phi_{79}(a,a) (a x)^i} 
  + \sum_{i=0}^{67} \frac{\Psi_{i,68}(a)}{a^{24} \Phi_{79}(a,a)} x^i
  \end{aligned}
\end{equation}
Where~$\Psi_{i,68}(a)$ has degree at most~$68$,
and~$\Phi_{79}(a,a)$ has degree~$70$.
Precisely the same argument as above holds (for~$|a| \ge 1000$),
again with a (very) comfortable margin, namely, 
$$\begin{aligned}
|Q_{24}(x)| \ge & \ 173210, \\
a^{-2} |Q_{k}(x)|  \le & \ 30, \  1 \le k \le 12, \\
\left| \frac{a \Psi_{-i,68}(a)}{a^{24 -i} \Phi_{79}(a,a)} \right| 
\le & \  10^{-20},  \ i = 1, \ldots, 11\\
\left| \frac{\Psi_{i,68}(a)}{a^{24} \Phi_{14}(a,a)} \right| \le & \ 10^{-20}, \ i = 0,\ldots, 67. \\
\end{aligned}
$$
\end{proof}

\subsection{Linear Independence of pure functions
and functions arising from~\texorpdfstring{$G$}{G}} \label{logsjointli}
Now, in a manner similar to~\S~\ref{jointli} (and with corresponding notation!), we have the following analogue of Lemma~\ref{14functions}.
(Remark~\ref{whyisthistrue} concerning
Lemma~\ref{14functions} is equally relevant in this case.)

\begin{lemma}[17 functions, logarithmic version] \label{17functions}
Assume that~$a$ and~$b$ satisfy the assumptions of Lemma~\ref{det}. Then
the ten functions
$$\int y G(y) \, dy, \int G(y) \, dy, \int \frac{G(y) - G(0)}{y} \, dy,
\int \frac{G(y) - G(0) - G'(0)y}{y^2} \, dy, $$
$$\int \frac{G(y) -  G(0) - G'(0)y- G''(0)\frac{y^2}{2}}{y^3} \, dy, 
\int \frac{G(y) - G(0) - G'(0)y- G''(0)\frac{y^2}{2} -G'''(0)\frac{y^3}{6}}{y^4} \, dy, 
$$
$$G(y),G'(y),G''(y),G'''(y),$$
together with the seven functions~$B_i(y)$ for~$i=1,\ldots,7$,
are linearly independent over~$\C(y)$.
\end{lemma}

\begin{proof} We proceed exactly as in the proof of Lemma~\ref{14functions}.
Namely, using a monodromy argument we replace~$G(y)$ by~$\hatG(y)$
and then with~$\Delta = \hatG(y) - G(y)$ we reduce to having to
show that a certain combination of derivatives and integrals
of~$\Delta$ only are linearly independent. However, $\Delta$
will now be a homogenous solution to the ODE $\LL = 0$, and
so it suffices to consider the case $\Delta(y) = G_A(y)$.
As in  the proof of Lemma~\ref{14functions}, we are reduced
 to an equation of the form
\begin{equation}
\label{todifflog}
\sum_{i=-4}^{1} a_i \int G_A(y) y^i \, dy 
= \sum_{i=0}^{3} b_i(y)  G_A^{(i)}(y).
\end{equation}
which we analyze by considering the local expansions
at the singular points of~$\LL$ described in Proposition~\ref{fictive G}.

Here the roots of~$R_4(y)$ are genuine singularities of the ODE,
whereas the roots of~$R_{10}(y)$ are not.
  Now, writing
$$b_3(y) = \sum_{i=N}^{\infty} r_{i}(y -\alpha)^i,$$
with~$N = N_{\alpha}$ and~$r_i = r_{i,\alpha}$,  just as in the proof of Lemma~\ref{14functions}, we have:
\begin{equation}
\label{onceagain}
\begin{aligned}
b'_3(y) + b_2(y) - \frac{c_3(y)}{c_4(y)} b_3(y) = & \ 0, \\
b'_2(y) + b_1(y) - \frac{c_2(y)}{c_4(y)} b_3(y) = & \ 0, \\
b'_1(y) + b_0(y) - \frac{c_1(y)}{c_4(y)} b_3(y) = & \ 0, \\
b'_0(y)  - \frac{c_0(y)}{c_4(y)} b_3(y) = & \  \sum_{i=-4}^{1} a_i  y^i,
\end{aligned}
\end{equation}
and inductively solving for~$b_i(y)$ the last equality in equation~(\ref{onceagain})
around various~$\alpha$ is as follows:
\begin{enumerate}
\item If~$\alpha = 0$, 
 the last equality becomes:  
$$\sum_{i=-4}^{1} a_i  y^i  
 =  \frac{-1}{4} (3-N)^2 (5-2N)^2 r_{N} y^{N-4} + \ldots $$
\item If~$\alpha = 4$, 
it becomes:
$$\sum_{i=-4}^{1} a_i  y^i 
=  \frac{-1}{4} (3-N) (2-N)(5-2N)(3-2N) r_{N} (y-4)^{N-4} + \ldots $$
\item If~$\alpha$ is a root $\beta$ of~$R_4(y)$, 
the last equality becomes:
$$\sum_{i=-4}^{1} a_i  y^i 
 =  - (3-N)^2 (2-N)(1-N) r_{N} (zy-\alpha)^{N-4} + \ldots $$
\item If~$\alpha$ is a root $\gamma$ of~$R_{10}(y)$, 
$$\sum_{i=-4}^{1} a_i  y^i  
=   (3-N) (2-N)(1-N)(1+N) r_{N} (y - \alpha)^{N-4} + \ldots $$
\item At~$\alpha \rightarrow \infty$, with~$b_3(y) =  y^N 
\displaystyle{\sum_{i=N}^{\infty} r_{i} y^{-i}}$,
we have
$$\sum_{i=-4}^{1} a_i  y^i 
 =-(5 - N) (4 - N)^2 (3 - N) r_{N} z^{N-4} r_{N} y^{N-4} + \ldots $$
\end{enumerate}
From these we deduce that:
\begin{equation}
\begin{aligned}
N_0 \ge & \ 0 \\
N_4 \ge & \ 2 \\
N_{\beta} \ge & \  1 \\
N_{\gamma} \ge & \ -1 \\
N_{\infty} \le & \ 5.
\end{aligned}
\end{equation}
This allows us to write:
\begin{equation}
\label{b3eqlog}
b_3(y) = \frac{(y-4)^2 R_4(y)}{R_{10}(y)} Q(y),
\end{equation}
Hence we may write
$$Q(y) =
q_0 + q_1 y + q_2 y^2 + \ldots q_9 y^9.$$
We find that, as a ratio of polynomials in~$y$, we have
\begin{equation} \label{degree12log}
 \sum_{i=-4}^{1} a_i  y^i=  \frac{
S_{42}(y)}
{y^4 R^4_{10}(y)} 
\end{equation}
for a polynomial~$S_{42}(y) \in \Q(a,b,y)$ of degree~$42$.
Note for degree reasons, this already implies that~$a_{1}=a_0=a_{-1}=0$),
Now solving for~$q_{0},\ldots,q_{9}$ in order to account
for a single factor of~$R_{10}(y)$,
we obtain a system of~$10$ linear equations in~$10$ unknowns.
If we take the corresponding determinant of the matrix, we obtain
a (symmetric) polynomial in~$a$ and~$b$ of the form:
$$(a-b)^{98} a^4 b^4 (a+b)^{12} (1 + 4 a^2- 4 a b)^2
(1 + 4 b^2 - 4 a b)(-3 + 4 a^2 - 4 a b + 4 b^2)^2$$
$$(9 + 16 a^2 - 40 a b + 16 b^2)^2
(-45 + 80 a^2 - 128 a b + 80 b^2)^2$$
$$\Phi_4(a,b) \Phi_6(a,b)  \Phi_{14}(a,b) \Phi_{26}(a,b) \Phi_{79}(a,b).$$
But each of these irreducibles factor is also a factor of
$$\Delta_{y} R_{10}(y)
\Res_y (R_{10}(y)  R_4(y))
\Res_{y} (R_{10}(y),y(y-4)),$$
which under our assumptions 
do not vanish
by Lemmas~\ref{res} and~\ref{det} respectively. Hence the determinant is non-zero,
which means that the~$q_i=0$, but then all the~$a_i$ are zero,
and there are no linear relationships, as claimed.
\end{proof}

\subsection{Location of the singularities}
\label{sec:location}
As noted in Proposition~\ref{fictive G},
the singularities of~$\LL$ in the~$Y_0(2)$ domain away from~$0,4,\infty$
are located at the points~$y_{a^{\pm},b^{\pm}}$ 
of Definition~\ref{yplus}, given by
the~$y := x^2/(x-1) = x + x/(x-1)$ images of
\[x = \left(a \pm \sqrt{a^2-1} \right) \left( b \pm \sqrt{b^2-1} \right).\]
Let us assume that~$\varepsilon < 10^{-6}$, and that
$$\left| \frac{m}{n} - 1 \right| < \varepsilon.$$
If~$m$ and~$n$ are distinct integers, then~$1 \le |m - n| < |n| \varepsilon$,
so~$|m|,|n| \ge \varepsilon^{-1}$.
Let~$a=2m+1$ and~$b=2n+1$. 
An elementary computation shows that, if~$\eta = (a + \sqrt{a^2-1})(b - \sqrt{b^2 - 1})$, that
$$\left|\eta + \frac{\eta}{\eta - 1} \right| > \frac{1}{\varepsilon}.$$
On the other hand, if~$\xi = (a - \sqrt{a^2-1})(b - \sqrt{b^2 - 1})$, and if~$a \ge b$, then~$\xi <  \varepsilon^2/4$, and
$$ \left| \xi + \frac{\xi}{\xi - 1} \right| \le \frac{\varepsilon^4}{16}.$$
Moreover, $\xi^{-1} > 4/\varepsilon^2$, and
$$ \left| \xi^{-1} + \frac{\xi^{-1}}{\xi^{-1} - 1} \right| > \frac{16}{\varepsilon^2} > \frac{1}{\varepsilon}.$$
It follows that the singularities of~\(\LL\) are all contained either
within the disc
$$D(0,\varepsilon^4/16) = D(0,10^{-12} 2^{-4}),$$
or outside the disc
$$D(0,\varepsilon^{-1}) = D(0,10^6).$$

\subsection{The proof of Theorem~\ref{logsmain}} \label{sec:proofC}

The overall argument is entirely similar to~\S~\ref{sec:proofA}. A putative~$\Z$-linear dependency~\eqref{rabc} with odd integers
$a,b \in \Z \setminus \{\pm 1\}$ reduces to the general~$\Z$-linear dependency
$$
2r_a\log{\left( 1+ \frac{1}{m} \right)}  + 2r_b\log{\left( 1+ \frac{1}{m} \right)} + r_{ab}\log{\left( 1+ \frac{1}{m} \right)} \log{\left( 1+ \frac{1}{n} \right)} = 4r_{0}, 
$$
writing $a = 2m+1, b = 2n+1$. We want to prove that there is no such relation if~$0 < |1-m/n| < 10^{-6}$, and 
we argue for the contradiction. 
By Proposition~\ref{fictive G}, the supposed relation produces a $G$-function
with unlikely analytic and arithmetic properties, including denominator type
$[1,\ldots,2n]^2$, which Lemma~\ref{17functions} promotes to some further associated functions, giving with~\S~\ref{sec:purefunctions} a totality of~$17$ functions of type~$n[1,\ldots,2n]^2$ and linearly independent over~$\Q(y)$. 
We are now in a position to reject by this~$G$-function an application of either one among Theorems~\ref{main:elementary form}, \ref{main:BC form},
or \ref{main: easy convexity}. 

All these theorems are to be used after changing the letter~$x$ of their respective statements to the symmetrization letter~$y := x^2/(x-1)$, and with the following ordered list~$\{f_i\}_{i=1}^{17}$ of~$17$ functions in Lemma~\ref{17functions}: 
\begin{equation*}
\begin{gathered}
B_1(y), \, B_2(y), \, B_3(y);  \, B_4(y), \, B_5(y), \, G(y), \, G'(y), \, G''(y), \, G'''(y); \\ B_6(y), \, B_7(y), \, \int y G(y) \, dy, \, \int G(y) \, dy, \int \frac{G(y) - G(0)}{y} \, dy,
 \\
\int \frac{G(y) - G(0) - G'(0)y}{y^2} \, dy, \, \int \frac{G(y) -  G(0) - G'(0)y- G''(0)\frac{y^2}{2}}{y^3} \, dy, \\
\int \frac{G(y) - G(0) - G'(0)y- G''(0)\frac{y^2}{2} -G'''(0)\frac{y^3}{6}}{y^4} \, dy. 
\end{gathered}
\end{equation*}
See~\eqref{defB}, \eqref{Nielsen}, and \eqref{defB67} for the functions~$B_1,\ldots,B_7$, and Definition~\ref{PtoG} for the function~$G$ (which in the end will not exist). 
The principal denominator types for this ordered list of functions forms the $17 \times 2$ array 
$$
\mathbf{b} := \left(  \begin{array}{lllllllllllllllll}   0 & 2 & 2 & 2 & 2 & 2 & 2 & 2 & 2 & 2 & 2 & 2 & 2 & 2 & 2 & 2 & 2 \\
0 & 0 & 0 & 2 & 2 & 2 & 2 & 2 & 2 & 2 & 2 & 2 & 2 & 2 & 2 & 2 & 2  \end{array} \right)^{\mathrm{t}}, 
$$
and the added integrations vector is
$$
\mathbf{e} :=  (0, 0, 1; 0, 0, 0, 0, 0, 0; 1, 1, 1, 1, 1,1,1,1).
$$

\addtocounter{subsubsection}{1}
 \begin{figure}[!h]  
\begin{center}
  \includegraphics[width=75mm]{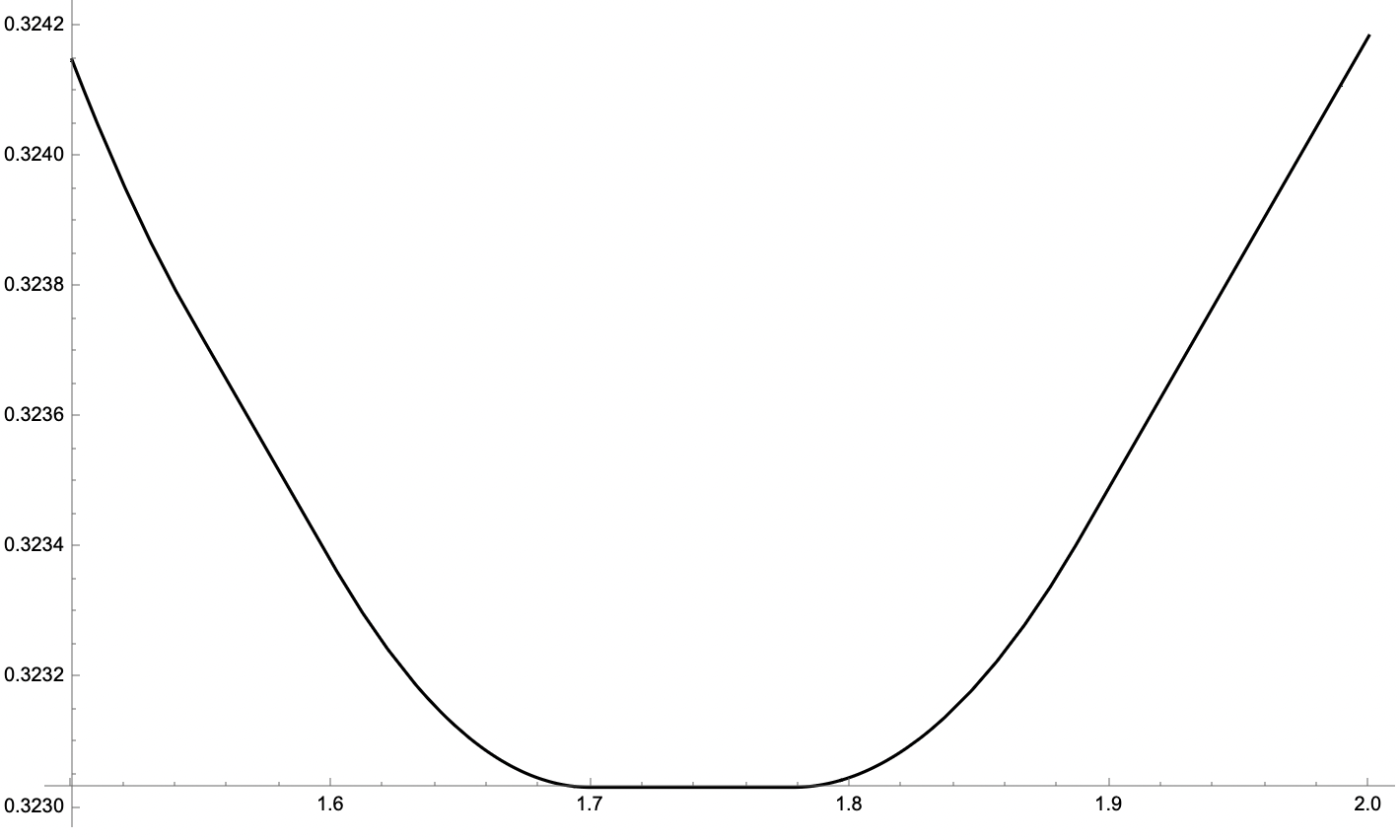}
\end{center}
\caption{The $\xi \in [1.5, 2]$ fragment of the graph of $(2/17^2)(9\xi  
   + I_{\xi}^{17}(\xi)) $, displaying the interval $\xi \in [1.7,1.78]$ being contained in the range of minimizer.}
   \label{xigraphlogs}
\end{figure}

For the ambient analytic map~$\varphi \in \mathcal{O}(\Db)$ we now select~$\varphi := h \circ \psi$, where --- yet again --- $h$ is the~$Y_0(2)$
hauptmodul written in the~$\tau = i\infty$ cusp-filling coordinate~$q = e^{2\pi i \tau}$ on the disc by the power series 
formula~\eqref{defofh}, and~$\psi : \Db \to \D$ is the holomorphic mapping from~\S~\ref{contour choiceC}. Corollary~\ref{stacky overconvergence}
now applies with Proposition~\ref{fictive G}, taking~$\Sigma_{Y_0(2)}^0 := \{ y_{a^-, b^-} \}$ to be the~$y := x+w(x) = x^2/(x-1)$ image of
\begin{equation*}
\begin{aligned}
\Sigma_{Y(2)}^0 := & \  \left\{ (a - \sqrt{a^2 - 1})(b - \sqrt{b^2 - 1}) \right\} \\
=  & \ \left\{ \frac{1}{(2n+1 + 2\sqrt{n^2+n})(2m+1+2\sqrt{m^2+m})}  \right\}; 
\end{aligned}
\end{equation*}
and~$\Sigma_{Y_0(2)}^1 := \emptyset$, $U_{Y_0(2)} := D(0,1/100)$, and of course,~$\varphi_{Y_0(2)}
:= \varphi = h \circ \psi$. Thus the analyticity conditions for our holonomy bounds are satisfied. 

For the denominator rates, the previous calculation now modifies to
\begin{equation}
\tau^{\flat}(\mathbf{b}) = \frac{ 1 \cdot 0 + (3+5) \cdot 2 + (7 + 9 + 11 + 13 +\ldots  + 33) \cdot 4 }{17^2}  
= \frac{1136}{289}, 
\end{equation}
and, from Figure~\ref{xigraphlogs} which reveals $\xi \in [1.7,1.78]$ to be contained by the minimizing interval, 
\begin{equation}
\begin{aligned}
\tau^{\sharp}(\mathbf{e}) & =  \frac{2}{m^2}  \min_{\xi \in [0,m]} \left\{  \xi \sum_{i=1}^{m} e_i 
   +
   \left(  \max_{1 \leq i \leq m}  e_i \right) I_{\xi}^{m}(\xi)
   \right\} \\  \label{tauraise2}
  & = \min_{\xi \in [0,17]} \left\{  (2/17^2) \left( 9\xi  
   + I_{\xi}^{17}(\xi)\right) \right\}  \\ & =   (2/17^2) \left( 9 \cdot 7/4  
   + I_{7/4}^{17}(7/4)\right)   = \frac{78419}{242760}.
   \end{aligned}
   \end{equation}
   
  Hence this time we obtain
   \begin{equation} \begin{aligned}
   \tau(\mathbf{b;e})  & = \tau^{\flat}(\mathbf{b}) + \tau^{\sharp}(\mathbf{e}) \\
   &  = \frac{1136}{289} + \frac{78419}{242760}
   = \frac{1032659}{242760}  = 4.2538\ldots, \end{aligned}
   \end{equation}
   arriving at the number $\frac{1032659}{242760}$ in~\eqref{logproductnumber}.

We can once again conclude the proof by deriving a contradiction of the form $m < 17$.
Just as in the proof of Theorem~\ref{mainA}, this can be done in a number of ways.
For example: 
applying Theorem~\ref{main:BC form}, we obtain the upper bound $m \leq 16.2$
as computed in~\eqref{logproductnumber}.  \hfill $\qed$

\section{Complements and further questions}  \label{sec:further questions}
We close our paper with a discussion on some further open problems naturally posed by our method.
But first, we discuss more closely the relationship between some of our results and more established methods.

\subsection{Comparison to the Siegel--Bombieri--Chudnovsky theory}  \label{sec:G-function arithmetic}
Theorem~\ref{logsmain} is an irrationality result in two parameters (subject to certain archimedean constraints). In this section, we compare 
Theorem~\ref{logsmain} to results previously available through the general arithmetic theory of special
values of~$G$-functions. 
We begin by recalling Siegel's  definition of a~$G$-function, fixing for this purpose a field embedding~$\Qbar \hookrightarrow \C$: 
\begin{df}[$G$-function] \label{gfunctiondef}
A power series~$f(x) = \sum_{n=0}^{\infty} a_n x^n  \in \Qbar \llbracket x \rrbracket$
 is a \emph{$G$-function} if it satisfies the following two properties:
\begin{enumerate}
\item \label{coefficientsinK} $f(x)$ is holonomic: it satisfies an ODE with coefficients in~$\Q[x]$. 
\item Both~$a_n$ and the denominators of the~$a_n$ have moderate growth, namely, the common denominator of~$a_0, \ldots, a_n$
grows at most exponentially in~$n$, and the largest Galois conjugate~$\ho{a_n}$ of~$a_n \in \Qbar \hookrightarrow \C$ grows as most exponentially in~$n$.
\end{enumerate}
\end{df}
From condition~(\ref{coefficientsinK}),  it follows that~$f(x) \in K \llbracket x \rrbracket$ for some number field~$K/\Q$.
As recalled in~\S~\ref{sec:introdenominators},
one expects (\cite{FischlerRivoal}) that for such an~$f(x)$, 
 there should exist $A \in \NwithoutzeroA$, $b \in \Q_{>0}$, and~$\sigma \in \NwithzeroA$ such that
\begin{equation}
\label{integralitytwo}
 a_n A^{n+1} [1, \ldots, b n]^{\sigma} \in \OL_K \qquad \forall n \in \NwithzeroA,
\end{equation}
and moreover  this is known unconditionally under the (conjecturally unnecessary) additional
assumption that~$f(x)$ arises from geometry~\cite[\S~V app.]{AndreG}. In any case, all the holonomic functions in our paper
do have denominators subsumed by~\eqref{integralitytwo}; they are manifestly~$G$-functions.

The basic paradigm of the arithmetic theory of~$G$-functions is captured by the following theorem:

\begin{thm}  \label{BombieriChudnovsky} 
For any $\Q(x)$-linearly independent set~$f_1, 
\ldots, f_h  \in \Q\llbracket x \rrbracket$ of $G$-functions with \emph{rational} coefficients,  
there is a constant
~$N_0 = N_0(f)$, effectively computable from the minimal ODEs of all the~$f_i$, 
such that the set
\begin{equation}  \label{Siegel like}
\begin{aligned}
&
\left\{ n \in \Z \, : \,  \emph{the $f_i(1/n)$ are $\Q$-linearly dependent or contain a divergent value} \right\}  \\ 
&  \subset [-N_0, N_0]. 
\end{aligned}
\end{equation}
\end{thm}

This result was  envisioned in Siegel's~1929 paper~\cite[\S~VII]{Siegel1929SNS} 
and proved, in the degree of abstraction that we state here, by
David and Gregory Chudnovsky~\cite{ChudnovskyG}, after the groundbreaking works of 
Galo\v{c}kin~\cite{Galochkin3} (who had to assume the `factorials canceling property' that reflects in the global nilpotence
of the integrable connection; a difficulty already noted by Siegel himself), and Bombieri~\cite{BombieriG} (who
proved a general adelic theorem under the similar and ultimately equivalent condition --- but by far easier to check than Galo\v{c}kin's --- 
 that the linear differential system 
is `Fuchsian of arithmetic type'.) The Chudnovskys' main result~\cite[Theorem~III]{ChudnovskyG}, \cite[Theorem~VIII.1.5]{Dwork},
\cite[\S~VI]{AndreG}, \cite{diVizio} was precisely the proof of the
global nilpotence property for all irreducible ODEs that possess at least
one $G$-series formal solution.

\begin{remark}  \label{general arithmetic theory}
In line with the discussion in~\S~\ref{binomials discussion},
the quantitative results on Siegel's program are, of course, stronger and more general than this quintessential form extracted from the works of Galo\v{c}kin, Bombieri, and the Chudnovskys. See Bombieri's Main Theorem\footnote{Noting Andr\'e's remark~\cite[page~79]{AndreG} that a scalar coefficient~'$2$' should be 
added in front of the summation over~$\zeta \in \mathrm{sing}_0(L)$ in the term~$c_{24}$ in Bombieri's Main Theorem.} in~\cite[page~49]{BombieriG}, and~\cite[Theorems~I and~II]{ChudnovskyG},  \cite[\S~1.2 Th\'eor\`eme Principal]{DebesG}, \cite[\S~VII]{AndreG} for other treatments with closely related results. A great picture of the pre-1997 state of the subject is in~\cite[ch.~5 \S~7]{FeldmanNesterenko}. For a more recent survey we refer to~\cite[\S~5.6]{RivoalSur}, as well as to~\cite{FischlerRivoalIneq} for further developments. Many of the standard relaxations, such as the admission of the apparently more general special arguments~$x = a/n \in \Q$ with~$|a/n| < c_1 \exp \left( - c_2 \sqrt{\log{n}  \cdot \log{\log{n}}} \right)$, can 
be subsumed into the form~\eqref{Siegel like} upon making explicit the dependence of~$N_0$ on the differential operator
following~\cite{BombieriG}.
  But~\eqref{Siegel like} is also a form that directly connects  to integral points on affine algebraic curves, and also to our framework in particular cases~\S~\ref{integral bounds}. We comment on the former connection in our next paragraph.

While 
Bombieri's general inequality is given in an adelic form over an arbitrary number field, 
the condition in Theorem~\ref{BombieriChudnovsky} on rational coefficients is of a fundamentally arithmetic nature, and it is crucial for the effectivity clause on~$N_0$. 
If for instance in the~$\{f_1,\ldots,f_h\} := \{1, f\}$ case one wants to handle algebraic number coefficients~$f \in \Qbar\llbracket x \rrbracket$ like in Remark~\ref{BCboundK}, the Siegel--Shidlovsky-style proof  logic in~\cite[\S~7]{ChudnovskyG} based on symmetric powers mandates that the hypothesis~$f(x) \notin \Qbar(x)$ (non-rational functions) would need to be strengthened to~$f(x) \notin \overline{\Q(x)}$ (transcendental functions).  
And indeed, Siegel's finiteness theorem on the integral points of non-rational affine algebraic curves has, to this day, not been resolved with an effective upper bound on the heights of the solutions\footnote{This is exactly the content of Hilbert's Tenth problem for the case of Diophantine equations in two variables.}, but it can be shown to be equivalent to the~$f(x) \in \overline{\Q(x)}$, $\{f_1, \ldots, f_h\} = \{1,f\}$ case of the statement~\eqref{Siegel like} \emph{with the assumption $f \in \Q\llbracket x \rrbracket$ \emph{(of rational coefficients)} relaxed to~$f \in \overline{\Q}\llbracket x \rrbracket$ \emph{(coefficients from a number field)}}. Hence, a statement
such as~\eqref{Siegel like} is a wide open question for the case of algebraic power series with coefficients from a number field other than~$\Q$ or an imaginary quadratic field.  The rational coefficients case handled by Theorem~\ref{BombieriChudnovsky}  reduces, in the algebraic case of $f(x) \in \overline{\Q(x)} \cap \Q\llbracket x \rrbracket$, to Bombieri's extension~\cite[Theorem~9.6.6]{BombieriGubler}, \cite[\S~IV]{BombieriWeildec}, \cite{DebesBombieri} of the classical \emph{Runge theorem}: an effective resolution in~$(x,y) \in \Z \times \Q$ of an irreducible bivariate Diophantine equation~$F(x,y) = 0$ over~$\Q$
when  the highest-order homogeneous part of~$F(x,y)$ is not proportional to a power of an irreducible polynomial over~$\Q$. (More intrinsically, under the \emph{Runge splitting condition}: the ``divisor at infinity'' used to give meaning to the integral points problem does not consist of a single Galois orbit of algebraic points on the algebraic curve. As 
is apparent from the explicit form of Bombieri's inequality~\cite[page 49]{BombieriG},
  the condition is arithmetic in nature and cannot be attained by extending to a number field; see~\cite[Equation~(9.26)]{BombieriGubler} for the general form of Runge's condition over the ring~$O_{K,S}$ of $S$-integers of a number field~$K$.)~\endofremark
\end{remark}

In these optics, our Theorem~\ref{logsmain} may be considered as an~$(x,y) = (1/n,1/m)$ special values analog for the particular \emph{bivariate} $G$-function
\begin{equation} \label{bi log}
F(x,y) := \log(1-x)\log(1-y).
\end{equation}
At least some basic results~\cite{Galochkin3,Galochkin1,Galochkin2,HataLogsOpp,Lysov} of this type are, of course, contained by the single variable Theorem~\ref{BombieriChudnovsky}, for instance one can evaluate the univariate $G$-function $f(x) = \log(1-x)\log(1+x)$ at the point~$x=1/n$. Already on an example as simple as this, the threshold term~$N_0$ arising from the general theory is extremely big;  it is estimated in~\cite{HataLogsOpp} to be on the order of~$e^{170}$ in this example. As far as we are aware, the record-lowest threshold on which the irrationality~$\log(1-1/n)\log(1+1/n) \notin \Q$ has been proved is Lysov's $n \geq 33$ in~\cite{Lysov}, by explicit (special!) Hermite--Pad\'e constructions. In this section, we investigate the scope of the general $G$-function  methods on our Theorem~\ref{logsmain}.

\subsubsection{The scope of the single variable theory}
To apply the single variable theory, we should treat~$k := m-n$ as a parameter, and consider the $G$-function
\begin{equation} \label{single app}
f(x)  := \log(1-x) \log\left( 1 - \frac{x}{1+kx} \right) \in \Q\llbracket x \rrbracket; \ 
 k := m-n \in \Z,
\end{equation}
whose value at~$x=1/n$ gives the desired product of two logarithms:
\[f(1/n)  = \log(1-1/n) \log(1-1/m) = F(1/n,1/m).\]
For notational simplicity alone, we shall only be concerned here with the irrationality of the product~$\log(1-1/n)\log(1-1/m)$, and not with its linear independence from the individual factors; this, of course, suffices for demonstrating the limitations of the general  arithmetic theory of special values of $G$-functions.   
Thus we apply Theorem~\ref{BombieriChudnovsky} with~$\{f_1, \ldots, f_h \} := \{1,f\}$. Then we need to quantify the $N_0 = N_0(k) = N_0(m-n)$ in Theorem~\ref{BombieriChudnovsky} as a function of~$k$, that is essentially of the height of the linear ODE. 

We claim that~$\log{N_0(k)} \asymp \log{|k|}$ 
for the minimal ODE of the function~$f_k$, by any of methods from the references that we listed in Remark~\ref{general arithmetic theory} on the general arithmetic $G$-function theory.\footnote{We omit the details, but the enterprising reader can find them
in the latex source code for this paper available on the arXiv.}
 Given this claim, the condition that guarantees~$f(1/n) \notin \Q$ from~\eqref{Siegel like} becomes~$|n|> |k|^c$ for some absolute constant~$c \in \R_{>0}$, that is the condition
$$
|1-m/n| < |k/n| < |n|^{1-1/c}. 
$$ 
This means that the cases of Theorem~\ref{logsmain} that were implicitly known through the general $G$-functions theory are all under a condition of the form
$$
0 < |1-m/n| \ll |n|^{-\kappa}, \qquad \textrm{for some } \kappa \in (0,1), 
$$
a condition necessary and sufficient for these general methods to apply; 
but a condition  significantly stronger than our  $0 < |1-m/n| < 10^{-6}$. 

\begin{quietcomment}
\subsubsection{Discussion through Bombieri's theorem} \label{bomb1}
To justify our claim on~$N_0(k)$,  we shall spell out Bombieri's ``Main Theorem'' from~\cite[page~49]{BombieriG}; see also~\cite[\S~IV]{BombieriWeildec} 
for a handy statement summary. (We note that the hypothesis of~$K(x)$-linear independence of the components of~$\mathbf{Y}$ is erroneously missing 
from the summary statement in~\cite[\S~IV]{BombieriWeildec}, and that, according to Andr\'e~\cite[page~79]{AndreG}, the~$\sum_{\zeta}$ term in both should be multiplied by the scalar factor~$2$.) In that setting, 
we have a number field~$K$ 
 and a 
first-order linear differential system~$\nabla = d/dx - \Gamma$, where~$\Gamma \in M_r(K(x))$ is an~$r \times r$ matrix
of rational functions over~$K$ and~$r$ is called the \emph{rank} of the system. Let $\mathrm{sing}(\nabla) \subset \overline{K}$, the set of \emph{singularities}
of~$\nabla$ (including apparent singularities, but not counting the possible singularity at~$\infty$),  be the (finite and~$\mathrm{Gal}(\overline{K}/K)$-stable)  set of~$\zeta \in \overline{K}$ that occur as a pole of some rational function component of the matrix~$\Gamma$. We may then view~$\nabla$ like an integrable connection over the  affine algebraic curve~$\A_K^1 \setminus \mathrm{sing}(\nabla)$ over~$K$. 

 We normalize the absolute value~$|\cdot |_v$ at each place~$v \in M_K$ in the usual way with the absolute conventions of~\cite[\S~1.4]{BombieriGubler}: 
 \[|\alpha|_v := |N_{K_v/\Q_p}(\alpha)|_p^{1/[K:\Q]}\]
  if~$v$ divides the rational place~$p \in M_{\Q}$, where, as in~\S~\ref{sec:Hermitian vector bundles}, the places on~$M_{\Q}$ are normalized in the usual way with $|p|_p = 1/p$ and~$|2|_{\infty} =2$. For a non-archimedean place~$v \in M_K^{\mathrm{fin}}$, let~$\C_v$ be the completion of some fixed choice of algebraic 
closure of~$K_v$, and let~$t_v \in \C_v$ (``a generic element'') be some unit (i.e., $|t_v|_v=1$) whose image in the residue field of~$\C_v$ is transcendental over the (finite) residue 
field of~$K_v$. The \emph{generic $v$-adic solvability radius} $R_v = R_v(t_v)$ of~$\LL$ is defined as the largest~$R \in [0,\infty]$ such that~$\ker(\LL)$ has  a full (i.e., rank~$r$) set of~$\C_v \llbracket x - t_v \rrbracket$ formal power series solutions all convergent on the open $v$-adic disc~$|x-t_v|_v < R$ in~$\C_v$. It 
is independent of the choice of~$t_v$, and so it can be denoted~$R_v(\nabla)$. Bombieri's condition (\emph{Fuchsian differential operators of arithmetic type}) is
to require that
\begin{equation}
\label{Bombieri condition}
\rho(\nabla) := \sum_{v \in M_K^{\mathrm{fin}}} \log^+{\frac{1}{R_v(\nabla)}} < \infty. 
\end{equation}
Similarly, on the level of a vector $\mathbf{Y} = \left( Y_1(x), \ldots, Y_r(x) \right)$ with power series components
$Y_j(x) = \sum_{n=0}^{\infty} a_{j,n} x^n \in K \llbracket x \rrbracket$ (not all of which are identically zero), we let $r_v(\mathbf{Y}) \in [0,\infty]$ 
to denote the smallest $v$-adic convergence radius of any of the power series~$Y_1,\ldots,Y_r$.
\emph{Bombieri's height} is defined~\cite[\S~7]{BombieriG} as 
\begin{equation}
\begin{aligned}
h(\mathbf{Y}) := \limsup_{N \to \infty} \frac{1}{N} & \left\{ \sum_{v \in M_K} 
\max_{\substack{ j \leq r \\ n \leq N }} \left\{  -n \log^+{\frac{1}{r_v(\mathbf{Y})}}  +  \log{|a_{j,n}|_v}   \right\} \right\} \\
&  + \sum_{v \in M_K} \log^+{\frac{1}{r_v(\mathbf{Y})}}, 
\end{aligned}
\end{equation}
and in this notation, the condition to be a \emph{formal $G$-series vector} simply reads: $h(\mathbf{Y}) < \infty$. 

Now assume that~$\nabla(\mathbf{Y}) = \mathbf{0}$ as formal power series regular in some open neighborhood of the origin~$x=0$, at every local completion~$K_v$ of~$K$. We can note in passing 
here that $h(\mathbf{Y}) < \infty$ if the Bombieri condition~$\rho(\nabla) < \infty$ is satisfied (this is a theorem of Bombieri and  Andr\'e~\cite[Theorems~V.5 and~V.6.2]{AndreG}, see
also~\cite[Theorems~VII.3.3, VII.3.4, and~VII.4.2]{Dwork}, or~\cite{DworkDiff} for a refinement); and that conversely, if~$\mathbf{Y}$ has $K(x)$-linearly 
independent components, the Chudnovskys's theorem states that~$\rho(\mathbf{Y}) < \infty$ implies~$\rho(\nabla) < \infty$ (more precisely,  
\cite[page~299]{Dwork} following \cite[Theorem~VI.4]{AndreG} yields --- in particular --- the explicit bound 
\[\rho(\nabla) \leq 4.95 \cdot r^2t\cdot h(\mathbf{Y}), \ \text{where}  \ t := \max \left\{ \deg{Q}, 1+ \deg{(\Gamma_{i,j}Q} )\right\}\]
 for a common polynomial denominator~$Q$ of the connection matrix entries~$\Gamma_{i,j} \in K(x)$).  
 Also, the Bombieri condition~\eqref{Bombieri condition} 
implies the global nilpotence property of the integrable connection~$\nabla$ and, by the Honda--Katz theorem, that the connection is of the Fuchsian class 
(regular singular points) and with rational local exponents. 
We do not need any of these theorems for stating Bombieri's inequality. Instead: 
simply consider a point~$\xi \in K \setminus \left\{ 0, \mathrm{sing}(\nabla) \right\}$ 
and a finite set of places~$S \subset M_K$ such that $|\xi|_v < \min(1, r_v(\mathbf{Y} ))$ (so that the $v$-adic numbers $Y_j(\xi) \in K_v$ can 
be defined at each~$v \in M_K$ as convergent power series values), and~$l \in \{1, \ldots, r\}$ linearly independent relations 
$$
\sum_{j=1}^{r} \lambda_{i,j} Y_j(\xi) = 0, \qquad \lambda_{i,j} \in K, \qquad i = 1, \ldots, l
$$
with $K$-rational coefficients, and valid simultaneously in every completion~$K_v$ at the places~$v \in S$. 
Let $h : \overline{K} \to [0,\infty)$ denote the absolute Weil canonical height function~\cite[\S~1.5.7]{BombieriGubler}. 
 For $F/K$ any finite extension and~$u \in M_F$ any place of~$F$, if we denote by~$v(u) := u|_K$
the place of~$M_K$ lying below~$u \in M_F$, then for each~$\alpha \in F$ we thus have
$$
h(\alpha) := \frac{1}{[F:K]} \sum_{u \in M_F} \log^+{| N_{F_u / K_{v(u)}} (\alpha) |_{v(u)}},
$$
a formula independent of the choices of \emph{both} number fields~$K$ and~$F/K$, and defining~$h : \overline{\Q} \to [0,\infty)$
with $h(\alpha) = 0$ if and only if~$\alpha \in \mu_{\infty} \cup \{0\}$ is a root of unity or zero. We have only used this particular presentation of Weil's canonical~$\mathbb{G}_m$ height as a shorthand convenience in~\eqref{Bombieri's inequality} below. 

In this setting and notation, we spell out the explicit content of the~$g=1$ case of~\cite[Main Theorem, page~49]{BombieriG}, with Andr\'e's remark accounted from~\cite[page~79]{AndreG}: 
\begin{BombieriInequality}
Assume that the~$G$-function vector solution~$\nabla(\mathbf{Y}) = \mathbf{0}$ has $K(x)$-linearly independent components. 
Then, for every $\varepsilon \in (0,1/2]$, we have
\begin{equation}  \label{Bombieri's inequality}
\begin{aligned}  
(r-2\varepsilon)  & \sum_{v \in S} \log{|\xi|_v}  + \left( r - l +  2\varepsilon \cdot \#\mathrm{sing}(\nabla)  \right) \sum_{v \in M_K} \log^+{|\xi|_v} \\
& \geq  -(r-l) \Big\{ \frac{r^2}{2\varepsilon} h(\mathbf{Y})   + \varepsilon \rho(\nabla) +2 \varepsilon \sum_{\zeta \in \mathrm{sing}(\nabla)} h(\zeta)
\\ & \ \  + \varepsilon \cdot \#\mathrm{sing}(\nabla)\log{2}    +\varepsilon (r-1) + (1+\varepsilon) \log{2} \Big\}. 
\end{aligned}
\end{equation}
Furthermore, the condition~$\varepsilon \in (0,1/2]$ can be replaced by~$\varepsilon \in (0,(r-1)/2]$ in the case that every non-zero element of~$\ker{\nabla}$ has~$K(x)$-linearly independent components.
\end{BombieriInequality}

Essentially equivalent inequalities (only with different constants and treatments) are~\cite[\S~1.2, Th\'eor\`eme Principale]{DebesG} and~\cite[\S~VI.3.6 Proposition on page~133]{AndreG}.

\begin{proof}[Derivation of Theorem~\ref{BombieriChudnovsky}]
We start with the hypothesis of~$\Q(x)$-linearly independent
$G$-functions~$f_1, \ldots, f_h$. Define~$r_1 \in \NwithzeroA$ 
so that~$f_1, f_1', \ldots, f_1^{(r_1)}; f_2, \ldots,f_h$ are still $\Q(x)$-linearly 
independent, but~$f_1^{(r_1+1)}$ becomes~$\Q(x)$-linearly dependent over these. 
Then define~$r_2 \in \NwithzeroA$ so that~$f_1, f_1', \ldots, f_1^{(r_1)}; f_2, f_2',\ldots, f_2^{(r_2)}; f_3,\ldots   \ldots,f_h$
are still~$\Q(x)$-linearly independent, but~$f_2^{(r_2+1)}$ becomes~$\Q(x)$-linearly dependent over these; and so on.
In this way, with~$r := h + r_1 + \cdots + r_m$ and upon ordering the list of functions
to beging with~$f_1, \ldots, f_h$, we obtain an~$r$-tuple
$$
(f_1, \ldots, f_h; g_1, \ldots, g_{r-h}) \in \Q\llbracket x \rrbracket^{r}
$$
 of $G$-series in~$\Q\llbracket x \rrbracket$  
whose~$\Q(x)$-linear span is a~$d/dx$-closed $\Q(x)$-vector space of dimension~$r$. Expressing
the linear operator~$d/dx$ in a matrix~$\Gamma \in M_{r \times r} (\Q(x))$ under this~$\Q(x)$-basis, 
we obtain a connection $\nabla = d/dx - \Gamma$ for applying Bombieri's inequality.  
David and Gregory Chudnovsky's fundamental theorem~\cite[\S~VI]{AndreG}, an explicit form of which (see \cite[page~299]{Dwork}) we spelled out in the discussion above, proves that the Bombieri condition~\eqref{Bombieri condition} 
is satisfied on~$\nabla$. 
 Bombieri's inequality then contains Theorem~\ref{BombieriChudnovsky} 
as the special case 
\[K = \Q,\,  l = 1, \,  S = \{\infty\},  \, \xi = 1/n, \, \mathbf{Y} =\left(f_1, \ldots, f_h; g_1, \ldots, g_{r-h}\right)\]
 with~$\lambda_{1,j} = 0$ for
$j \notin \{1, \ldots, h\}$, upon remarking that the two sums on the top row in~\eqref{Bombieri's inequality} now reduce to just~$-\log{|n|}$ and~$+\log{|n|}$, respectively.  We thus simply need to fix an~$\varepsilon \in (0,1/2]$ small enough that the coefficient~$r-2\varepsilon$ of the first sum exceeds the coefficient~$r-1 +2\varepsilon \cdot \#\mathrm{sing}(\nabla)$ of the second sum. 
An explicit~$N_0 = N_0(\nabla)$ then ensues in~\eqref{Siegel like} in terms of the invariants~$\rho(\nabla)$ and~$\sum_{\zeta \in \mathrm{sing}(\nabla)} h(\zeta)$ of the integrable connection~$\nabla$. 
\end{proof}

\subsubsection{The scope on Bombieri's inequality on Theorem~\ref{logsmain}}
We  apply Bombieri's inequality to the~$\{f_1, \ldots, f_h\} := \{1, f\}$ situation at hand with the $G$-function~\eqref{single app}. In this case, the finite singularities of the 
integrable connection~$\nabla$ constructed in the preceding derivation are at the points~$\zeta = 1$ and~$\zeta/(1+k\zeta) = 1$, to wit, $\zeta = -1/(k-1)$. 
These give a total height term~$h(1) + h(-1/(k-1)) = \log{|k-1|}$ in the sum over~$\zeta$ on the middle line in~\eqref{Bombieri's inequality}, 
and one easily also computes that the other two quantities~$h(\mathbf{Y}) \asymp \log{|k|}$ and~$\rho(\nabla) \asymp \log{|k|}$. 
Bombieri's inequality then reads: 
\begin{equation*}
\begin{aligned}
-(2-2\varepsilon) \log{|n|} + \left(  1+ 4\varepsilon  \right) \log{|n|} 
\gg - \left( \varepsilon + \varepsilon^{-1} \right) \left( 2+ \log{|k|} \right),
\end{aligned}
\end{equation*}
where the implicit coefficient is an absolute constant of no importance to us. 

In other words, for~$|k| \geq 2$, we derive
\[\left( 1 - 2\varepsilon \right) \log{|n|} \ll (\varepsilon + \varepsilon^{-1}) \cdot \log{|k|}\]
 under the assumption of rationality~$f(1/n) = \log(1-1/n) \log(1-1/m) \in \Q$. This 
gives exactly (and no better than this) the claimed bound~$|n| \leq |k|^c$, with some absolute constant~$c \in \R_{> 0}$, upon
choosing say~$\varepsilon = 1/16$. 
\end{quietcomment}

\subsubsection{Multivariable $G$-function theory}  \label{bomb2}
The arithmetic theory of multivariable~$G$-functions is still in its infancy; we refer to~\cite{AndreBaldassari} for the geometric foundations, and
to~\cite{Nagata2d} (for dimension two) and~\cite{diVizio} (for arbitrary dimension) for the generalizations of the  Chudnovskys' fundamental theorem. 
 It seems likely that a two-dimensional version of Siegel's approximating forms scheme~\S~\ref{sec:Hermite approximations}, as carried out by Bombieri in~\cite{BombieriG}, might combine  
 with standard nonvanishing methods~\cite{Dyson,BombieriDyson}  for Diophantine auxiliary constructions at a special point, to give the irrationality~$f(1/n,1/m) \notin \Q$ of the specializations of a bivariate $G$-function $f(x,y) \in \Q\llbracket x,y\rrbracket \setminus \Q(x,y)$ at arguments of the form~$x = 1/n, y = 1/m$, where~$n,m \in \Z \setminus \{0\}$ with~$\log{|n|} \gg_f 1$ and~$\log{|m|}/\log{|n|} \gg_f 1$; to our knowledge, this kind of program has not as yet been worked out in the literature. The range~$|1-m/n| < 10^{-6}$  that we obtained for~\eqref{bi log} from the Ap\'ery limits method is entirely orthogonal to this!

\subsubsection{Avenues from Hermite--Pad\'e constructions} As explained in~\S\S~\ref{log discussion},~\ref{ARB}, our proof in~\S~\ref{sec:logs} of Theorem~\ref{logsmain} can be conceptually linked 
to the Hermite--Pad\'e approximants to the logarithm function, used with the Hadamard product construction~\S~\ref{sec:Hadamard}. 
While products of three or more logarithms appear unreachable by our method here (by the numerology~$e^3 > 16$), 
it could be worthwhile to attempt linear independence of more than a single pair of products of two logarithms, starting from
the simultaneous Hermite--Pad\'e approximation theory with several logarithms worked out explicitly in~\cite{RhinToffin} and~\cite{DavidHirataKohnoKawashima}. 

A more general scheme, such as we indicated on the most basic examples in~\S~\ref{binomials discussion} and~\S~\ref{log discussion}, could be sought with forming the generating function of the special linear forms obtained 
from evaluating a regular sequence of functional Hermite--Pad\'e approximants to a basic function. Beukers~\cite{BeukersPade,BeukersPolylogs}, using polylogarithms, and Pr\'evost~\cite{Prevost}, having $\zeta(3,1+1/y)$ for the basic function to be evaluated at the points of the form~$y=1/n$, 
 each were able to interpret the Ap\'ery sequences inside such a scheme. The former type was vastly generalized by Fischler and Rivoal~\cite{FischlerRivoalPade}. It could be interesting to find a similar interpretation with simultaneous Hermite--Pad\'e approximants for the simultaneous linear forms in~$\zeta(2)$ and~$L(2,\chip)$ that we exploited in~\S~\ref{Zagier local system}. 
We note however that such generating function procedures far from always give rise to~$G$-functions, 
even if the starting function for the Hermite--Pad\'e approximation is algebraic~\cite{BombieriCohenPade}; for the Pr\'evost type, some non-examples 
related to zeta values are in~\cite[\S~8]{PrevostRivoalD}. 

Two other subjects that we have omitted here (in part, for reasons of space) are applications to non-rational algebraic arguments for the two logarithms, as well as $p$-adic logarithms. Another reason for omitting the
later application is that, in this paper, we have emphasized the archimedean place as special when it comes to  overconvergence.
In~\cite{zeta5}, we plan to write our holonomicity bound in a more general Arakelov adelic form over a global field.

 Finally, speaking more broadly of holonomic explicit constructions by any method, we remark that in many of the more intricate ones in the literature --- 
 such as in Zudilin's work~\cite{ZudilinHypergeometricTales} on simultaneous approximation to~$\zeta(2)$ and~$\zeta(3)$, and in Brown and Zudilin's work~\cite{BrownZudilin} on~$\zeta(5)$ ---  progress towards a not-yet-attained irrationality goal is measured by setting up a complex set of parameters to maximize the {\it worthiness exponent} $\limsup\left\{ -\frac{\log{|\eta -p/q|}}{\log{q}} \right\}$. The latter is very far from faithfully measuring worthiness as a potential for applying in our framework of rational holonomy bounds.

\subsection{Integral holonomic modules} \label{integral bounds} 
Although~$N_0$ in Theorem~\ref{BombieriChudnovsky} is effective, it is a wide open\footnote{Once again, the exception is the algebraic case~$f_i(x) \in \overline{\Q(x)} \cap \Q\llbracket x \rrbracket$, in which case a finite computer search is at least in theory enough to  finish off this problem in finite computational time. } problem to precisely (in principle) determine the left-hand side set in~\eqref{Siegel like}. Moreover, as we have seen with a case as simple as~\cite{HataLogsOpp}, 
the bound on~$N_0$ in the general property~\eqref{Siegel like} is, in practice, very big. Our findings with the $\Q \left[ x,  \frac{1}{1-x}  \right]$-integrality refinements in~\S\S~\ref{sec:logarithmchar},~\ref{sec:beyondlog}, see especially Remarks~\ref{reverse} and~\ref{contains in itself}, seem to  
point towards a completely different approach to those of the cases of~\eqref{Siegel like} whose holonomicity is recognized by Andr\'e's arithmetic  criterion (Corollary~\ref{holonomic criterion}). Our inspiration for hoping to reverse the proof logic in Remark~\ref{reverse} for the purpose of applying to similar other potential linear independence setups --- many of them unproved conjectures --- 
  stems from the ideas of B\'ezivin and Robba~\cite{BR} which they used for reproving the Hermite--Lindemann--Weierstrass theorem as an application of Bertrandias's arithmetic rationality criterion~\cite[Th\'eor\`eme~5.4.6]{Amice}  (see also~\cite{BezivinRobbaAMM} for a historical dissection of that proof); and, ultimately, the refinement of those ideas at the hands of  Andr\'e~\cite{AndreGevreyI,AndreGevreyII} and Beukers~\cite{BeukersFreeness}, using the Chudnovskys's theorem and the Fourier--Laplace duality 
between $E$- and $G$-functions, to reprove (and further refine) the qualitative Siegel--Shidlovsky theorem on special values of~$E$-functions.

Let~$\partial := x \cdot (d/dx)$ be the multiplicatively invariant derivation. 
Consider (for simplicity here) an \emph{\'etale}\footnote{In other words: the derivative~$\varphi'$ is nowhere vanishing on~$\D$.} holomorphic mapping~$\varphi : \D \to \C$ taking~$\varphi(0) = 0$, and a vector~$\bb := (b_1,\ldots,b_r) \in [0,\infty)^r$ with~$|\varphi'(0)| > e^{b_1+\cdots+b_r}$. 
Consider the set~$\mathcal{D}$ comprised of the formal power series of 
the shape 
\begin{equation} \label{type dens again}
f(x) = \sum_{n=0}^{\infty} a_n \frac{x^n}{\prod_{h=1}^r [1,\ldots, b_i \cdot n]}, \qquad a_n \in \Z \quad \forall n \in \NwithzeroA
\end{equation}
such that~$\varphi^*f \in \mathcal{O}(\D)$ is a holomorphic function on the disc. This is a module over the noncommutative
ring~$\Z[x,\partial]$: the derivation~$\partial$ acts on the monomials by~$\partial(x^n) = nx^n$, which preserves the integrality type in~\eqref{type dens again}, 
while the chain rule with the \'etaleness of~$\varphi$ show that if~$f(\varphi(z))$ is holomorphic, so also is 
\[\frac{\varphi(z)}{\varphi'(z)} \left( f(\varphi(z)) \right)' 
= \varphi(z) f'(\varphi(z))= (\partial f)(\varphi(z)).\] 
By construction, the~$\Z[x,\partial]$-module $\mathcal{D}$ is embedded as a submodule of the ring~$\Q \llbracket x \rrbracket$. 
  Within this ambient ring,
       $\mathcal{D}$ contains the subring
comprised of the~$\alpha \in \Z\llbracket x \rrbracket$ with~$\varphi^*\alpha \in \mathcal{O}(\D)$. 
Let us denote this ring by~$\mathcal{O}(\wV)$, for reasons related to~\cite{BostCharles} and our Remark~\ref{rmk_BCmetric}, 
in which this ring is the ring of regular functions on the formal-analytic\footnote{If~$\varphi$ extends to a holomorphic function on some open neighborhood of the closed disc~$\Db$, to match the convention in~\cite{BostCharles}. For applying the finiteness theorem~\cite[Theorem 9.1.1]{BostCharles}, this is not a restriction upon considering~$\widetilde{\varphi}(z) := \varphi\left( (1-\varepsilon) z\right)$ with an~$\varepsilon > 0$ small enough to still have~$|\widetilde{\varphi}'(0)| > e^{b_1+\cdots+b_r} \geq 1$.}  arithmetic surface we denoted~$\wV := \wV(\varphi)$. 
Then~$\mathcal{D}$ is a module over the ring~$\mathcal{O}(\wV)$. 
By Corollary~\ref{holonomic criterion}, the field of fractions~$\mathrm{Frac} ( \mathcal{O}(\wV) )$ is a finite field 
extension of~$\Q(x)$, and $\mathcal{D} \otimes_{\mathcal{O}(\wV)} \mathrm{Frac} ( \mathcal{O}(\wV) )$ is a finite-dimensional
vector space over that field. Since furthermore~$\mathcal{D}$ is preserved by the derivation~$\partial$, 
there is a finite and~$\Gal(\Qbar/\Q)$-stable set of complex algebraic points~$\Sigma \subset \Qbar \hookrightarrow \C$ such that all elements
of~$\mathcal{D} \otimes_{\mathcal{O}(\wV)} \mathrm{Frac} ( \mathcal{O}(\wV) )$  (and, \emph{a fortiori}, all elements of~$\mathrm{Frac} ( \mathcal{O}(\wV) )$) continue analytically as \emph{meromorphic} functions along all paths in~$\C \setminus \Sigma$.

However, Bost and Charles proved a deeper finiteness theorem \cite[Theorem~9.1.1]{BostCharles}:  the ring~$\mathcal{O}(\wV)$ is a finitely generated~$\Z$-algebra. Moreover, their proof leads in principle to an effective algorithm for listing a finite set of generating elements for this~$\Z$-algebra. On the other hand, as shown by Remark~\ref{integral module log} on the example~$\bb = (1)$ and~$\varphi(z) = 4z/(1+z)^2$, where~$\OL(\wV) = \Z[x,1/(1-x)]$, 
the $\OL(\wV)$-module~$\mathcal{D}$ is in general infinite. 
We hence tensor it with~$\Q$ and consider~$\mathcal{D}_{\Q} := \mathcal{D} \otimes_{\Z}\Q$,  
which is a torsion-free module over~$\OL(\wV_{\Q}) := \OL(\wV) \otimes_{\Z} \Q$. The latter\footnote{This is a definition in our \emph{ad hoc} notation here, which does not occur in~\cite{BostBook} or~\cite{BostCharles}.} ring is a finitely generated $\Q$-algebra but also, being 
of Krull dimension one and integrally closed in its fraction field, 
it has the added simplicity of being a Dedekind domain. A module over a Dedekind domain is finite and torsion free if and only if it is locally free of finite rank, if and only if it is projective and {\it generically finite} (where the latter means the finite-dimensionality of the induced vector space over the field of fractions). The content of Theorem~\ref{logcharacterization} is that, in the previous example of~$\bb = (1)$ and~$\varphi(z) = 4z/(1+z)^2$, 
the $\OL(\wV_{\Q})$-module~$\mathcal{D}_{\Q}$ is free of rank~$2$ with basis~$\{1, \log(1-x)\}$. However, the content of Remark~\ref{x=0 special role}
is that, in the slightly modified example~$\bb = (1 + 1/100)$ and~$\varphi(z) = 4z/(1+z)^2$ still having~$\OL(\wV) = \Z[x,1/(1-x)]$
and~$\OL(\wV_{\Q}) = \Q[x,1/(1-x)]$, the $\OL(\wV_{\Q})$-module $\mathcal{D}_{\Q}$ is infinite while the $\OL(\wV_{\Q})[1/x]$-module
 $\mathcal{D}_{\Q}[1/x]$ is once again free of rank~$2$ with basis~$\{1, \log(1-x)\}$. 
Similar remarks apply to Theorem~\ref{thm three elements}. In line with these we could ask:

\begin{question}  \label{freeness query}
Can one effectively construct an~$h \in \mathcal{O}(\wV) \setminus \{0\}$ so that:
\[ \label{latestar}
(\star) \  \textrm{The module }  
\mathcal{D}_{\Q}[1/h] \textrm{ is locally free over the ring } 
 \mathcal{O}(\wV_{\Q})[1/h]?
\]
Are there natural  verifiable conditions  such that this holds even with~$h=x$?
\end{question}
As remarked above, the rings~$\mathcal{O}(\wV_{\Q})$ and their localizations are Dedekind domains, 
and hence, since the~$\mathcal{O}(\wV_{\Q})[1/h]$-module~$\mathcal{D}_{\Q}[1/h]$ is torsion-free and generically finite,  
the local freeness in~$(\star)$ is equivalent to the module being finite, and also to the module being projective. 
   We will see why the insistence on effectivity is the important point for deriving irrationality proofs on~$x=1/n$ special values of certain functions from the holonomic module~$\mathcal{D}$. 
The reason for 
such a connection is the same as with the  Andr\'e--Beukers (qualitative) refinement~\cite{AndreGevreyI,AndreGevreyII,BeukersFreeness} of the Siegel--Shidlovsky theorem  being ultimately derived from a commutative algebra statement formally
similar to~$(\star)$: 
\begin{fact} \label{E is free}
 The ring of~$\mathbf{E} \subset \Q\llbracket x \rrbracket$ of~$E$-functions with rational coefficients generates over the Laurent polynomial
 ring~$\Q[x,x^{-1}]$ an infinite free $\Q[x,x^{-1}]$-module $\mathbf{E}[1/x] = \mathbf{E} \otimes_{\Q[x]} \Q[x,x^{-1}]$. 
  \end{fact}
  
  See~\cite[Theorem~1.5]{BeukersFreeness} for the statement\footnote{A theorem of Kaplansky, see~\cite[\S~4.1.2]{BostBook}, states that over a Dedekind domain any module which is projective and countably generated, but is not finitely generated, is a free module. Hence the~$\Q[x,x^{-1}]$-freeness of~$\mathbf{E}[1/x]$
  reduces to the freeness of all its finitely generated $\Q[x,x^{-1}]$-submodules, that is to~\cite[Theorem~1.5]{BeukersFreeness}.} and~\cite[proof of Cor.~2.2]{BeukersFreeness} for the mechanism.
  The inversion of~$x$ here is also necessary, just as we saw with~$(\star)$ on the example~$\bb = (1 + 1/100)$ and~$\varphi(z) = 4z/(1+z)^2$.
  In the  setting of $E$-functions, take 
  \[\{ (d/dx)^j \left\{ (e^x-1)/x \right\} \, : \, j \in \NwithzeroA \}\]
  as an example. This set  generates an infinite $\Q[x]$-module
  which localizes to a rank~$2$ free~$\Q[x,x^{-1}]$-module with basis~$\{1,e^x\}$. 
The special role of~$x=0$ in the Andr\'e--Beukers theory reflects the presence
of transcendental $E$-functions like~$f(x) = e^x$ whose  minimal ODE does not have~$0$ for singularity,
but yet the special value~$f(0) = 1 \in \Q$ is rational. 
  
  There is a similar formal mechanism to~\cite{BR,AndreGevreyII,BeukersFreeness} for deriving linear independence proofs if Question~\ref{freeness query} has a positive answer. Suppose~$h \in \mathcal{O}(\wV) \setminus \{0\}$ is such that the localized $\mathcal{O}(\wV_{\Q})[1/h]$-module~$\mathcal{D}_{\Q}[1/h]$ is locally free. 
Assume now additionally that, as in~\S~\ref{local univalent leaves},  there is a contractible open neighborhood $0 \in \Omega \subset \D$ to  which~$\varphi$ restricts as a univalent map, and such that~$\varphi^{-1}(\Sigma) \subset \Omega$. Consider  $f(x) \in \mathcal{D}_{\Q}$ and a finite,~$\Gal(\Qbar/\Q)$-stable set of complex algebraic points~$S_f \subset \Qbar \hookrightarrow \C$ containing~$\Sigma$, such that~$f(x)$ continues~\emph{holomorphically}  along all paths in~$\C \setminus S_f$. Assume furthermore that
those various analytic continuations end up taking on at least two distinct values at the point~$x = 1/n$. 

Now suppose~$n \in \Z \setminus \{0\}$ obeys the following restrictions: 
\begin{enumerate}
\item \label{local univalence at the point} $ \varphi^{-1}(1/n) \subset \Omega$;
 \item  \label{here nons} $1/n \notin S_f$;
\item \label{R regular} All analytic continuations of the element~$h$ in~$\C^1 \setminus \Sigma$ take nonzero values at the point~$1/n$. 
\end{enumerate}
{\it Then the special\footnote{This is meant as~$\Q\llbracket x \rrbracket$ series evaluated at~$x = 1/n$; hence
the usual dichotomy with ``either irrational or divergent.''  We can, however, say here the more precise conclusion of irrationality of the value~$f(1/n) \in \C$, well-defined by analytic continuation from~$x=0$ staying within the univalent leaf~$\Omega \supset \{0,1/n\}$. An inclusion into a Ball--Rivoal framework of special values beyond the disc of convergence for certain $G$-functions has been recently achieved by Fischler and Rivoal~\cite{FischlerRivoalO}.} value~$f(1/n)$ satisfies an irrationality property as in~\eqref{Siegel like}: either~$f(1/n) \notin \Q$, or else~$f(1/n)$ is divergent. }

The point is formally the same as in~\cite{BeukersFreeness}. There is a non-zero polynomial~$Q_f \in \Z[x] \setminus \{0\}$ such that~$\{Q_f = 0\} = S_f$ and~$Q_f(x) f(x)$ continues~\emph{holomorphically} along all paths in~$\C \setminus \Sigma$.  
If~$f(1/n) = p/q$ were rational, the local univalence property~\eqref{local univalence at the point} in our setup
from~\S~\ref{local univalent leaves} would apply to the function
$$
\widetilde{f}(x) :=  Q_f(x) \frac{ f(x) - f(1/n) }{1-nx}  \in  \Q\llbracket x \rrbracket,
$$ 
with $\Sigma^1 := \{ s \in \Sigma \, : \, \varphi^{-1}(s) = \emptyset \}$ and $\Sigma^0 := \{1/n\} \cup (\Sigma \setminus \Sigma^1)$ and~$\Sigma$ of Proposition~\ref{overconvergence} augmented by~$\Sigma \cup \{1/n\}$, to derive that 
\begin{equation} \label{all derivatives}
\partial^j \left\{  \widetilde{f}(x) \right\} \in \mathcal{D}_{\Q}, \qquad \textrm{ for all }  j \in \NwithzeroA. 
\end{equation}
The assumptions~$(\star)$,~\eqref{here nons}, and~\eqref{R regular} imply that the functions~\eqref{all derivatives} generate a finite~$\C(x)$-module,
 as well as a finite~$\C\llbracket \frac{1}{1-nx}\rrbracket$-module. Hence~$\LL(\widetilde{f}) = 0$ for some nonzero linear differential operator~$\LL$ over~$\C(x)$ which is
non-singular at the point~$x=1/n$. But this conflicts with our condition that~$\widetilde{f}$  has some analytic continuation~$\widetilde{F}$, necessarily also a solution of the ODE $\LL(\widetilde{F}) = 0$, such that~$x=1/n$ is a meromorphic pole of~$\widetilde{F}(x)$. 

This puts some prize on constructing an element~$h \in \OL(\wV) \setminus \{0\}$ for the finiteness property~$(\star)$
in Question~\ref{freeness query}. As a simple example, if~$h=1$ or even just~$h=x$ (analogously to the Andr\'e--Beukers Fact~\ref{E is free} and the discussion preceding it) is admissible 
for the bivalent map~$\varphi(z) = 8(z+z^3)/(1+z)^4$ of Basic Remark~\ref{bivalent clean} and the~$[1,\ldots,n]^2$
denominator type ($\bb = (1,1)$), that by itself would suffice --- in lieu of the full Conjecture~\ref{conjecture five elements} in that context, which
could be more difficult --- to embed Remark~\ref{contains in itself} into the above discussion, and conclude
at one stroke the irrationalities~$\Li_2(1/n) \notin \Q$ for the remaining values~$n \in  \{-4, -3, -2, 2, 3, 4, 5\}$. 
Indeed, on this example, an easy computation  (see, for example, \cite{robinson}) gives $\mathcal{O}(\wV) 
 = \Z\left[ x, 1/(1-x) \right]$ with~$\Sigma = \{0, 1, \infty\}$, and for $\Omega$ we can take any open neighborhood
 of~$(-1,1)$ in~$\D$ that is small enough to have~$\varphi^{-1}(\varphi(\Omega)) = \Omega$.

\subsection{Quantitative aspects of linear independence}
As is usual\footnote{\emph{Except} for Andr\'e's \emph{transcendance sans transcendance~\cite{AndreGevreyI,AndreGevreyII}} applying the arithmetic theory of~$G$-functions to recover the Siegel--Shidlovsky theorem on~$E$-functions.} in transcendental number theory, our linear independence proofs in this paper can in principle by promoted to quantitative lower bounds on the linear forms in the relevant periods. In this case, however, the transition is not straightforward and requires a substantial amount of added work that we decided to not engage with in the present paper. The discussion in~\S~\ref{binomials discussion} points to a first methodological clue for making such a transition. We plan to turn to this in a future work.

\subsection{The structure ring} It would be interesting to clarify the scope of  arithmetic characterization theorems of the kind of Theorem~\ref{logcharacterization} (on~$\log{x}$) and Theorem~\ref{thm three elements} (on~$\log^2{x}$), and of the more precise Conjecture~\ref{conjecture five elements} on the~$[1,\ldots,n]^2$ layer~$G$-functions on~$\P^1 \setminus \{0, 1, \infty\}$. One could for example ask how much of the multiple polylogarithm ring~\S~\ref{mpoly} may be captured in arithmetic algebraization terms.  The following question falls short of our methodology:
\begin{question} \label{tauthree}
Consider~$\HH$ to be the~$\Q(x)$-vector space generated by functions~$f(x)$
of the form
\begin{equation}
f(x) = \sum_{n=0}^{\infty} a_n \, \frac{x^n}{[1,\ldots,n]^{3}} \in \Q\llbracket x \rrbracket, \qquad a_n \in \Z \quad \forall n \in 
\NwithzeroA
\end{equation}
arising from~$G$-functions of geometric origin on~$\mathbf{P}^1 \setminus \{0,1,\infty\}$.
Is~$\HH$ finite dimensional?
\end{question}

The universal map~$\varphi : \D \to \C \setminus \{1\}$ taking $\varphi(0)=  0$ 
and satisfying~$\varphi^{-1}(0) = \{0\}$ is~$\lambda$.
Since~$16 < e^{\tau} = e^3$, our methods have nothing directly
to say about Question~\ref{tauthree};
we do not even
have a guess as to what the answer might be. (To contrast, 
for the denominator types --- for example --- $[1,\ldots,n]^2[1,\ldots,n/2]$ or $\prod_{k=1}^{8} [1,\ldots,n/k]$, Corollary~\ref{holonomic criterion} proves
that the corresponding~$G$-functions on~$\mathbf{P}^1 \setminus \{0,1,\infty\}$ form a finite-dimensional space, although it is probably quite difficult
to determine these spaces, even conjecturally.)

We emphasize that the problem does not necessarily become any
easier even when~$\tau = 0$; one can ask for which~$\alpha \in \Q$ the~$\Q(x)$-vector
space generated by \emph{algebraic} power series~$f(x)$ in~$\Z \llbracket x \rrbracket$ 
on~$\mathbf{P}^1 \setminus \{0,\alpha,\infty\}$ is infinite.
This space is finite when~$\alpha > 1/16$,
and infinite when~$\alpha = 1/16$, where one can construct
such functions by writing modular functions with integer coefficients 
(on congruence subgroups)
in terms of~$x = \lambda/16$.
The main result of~\cite{UDC}  was to show that all such~$f(x)$ arise in this way.
But as soon as~$\alpha < 1/16$, we have no methods to understand this 
problem, or even to determine whether there exists a single non-rational
function in~$\HH$. 
One can, however, leverage the~$\alpha = 1/16$ example to show that
 for certain templates with~$\tau > 0$, the space of~$G$-functions is infinite dimensional.

\begin{proposition}
Let~$\HH$ denote the $\Q(x)$-vector space generated by functions~$f(x)$
of the form
\begin{equation}
f(x) = \sum_{n=0}^{\infty} a_n \, \frac{x^n}{[1,\ldots,2n]^{2}} \in \Q\llbracket x \rrbracket, \qquad a_n \in \Z \quad \forall n \in \NwithzeroA
\end{equation}
arising from~$G$-functions of geometric origin on~$\mathbf{P}^1 \setminus \{0,1,\infty\}$. 
Then~$\HH$ is infinite dimensional.
\end{proposition}

\begin{proof} 
Let~$g(q) \in \Z\llbracket q \rrbracket$ be a modular function on~$X_0(N)$ which is holomorphic
away from the cusps. Then, writing~$f(q)$
as a function of~$x = \lambda/16$, we find that
$$g(x) \in \Z \llbracket x \rrbracket$$
is an algebraic function on~$\mathbf{P}^1 \setminus \{0,1/16,\infty\}$,
and the space of such~$g(x)$ is infinite dimensional over~$\Q(x)$ (by
taking larger and larger~$N$). Now let
$$h(x) =  \pFq{3}{2}{1,1,1}{1/2,1/2}{\frac{x}{16}} 
= \sum \frac{x^n}{\binom{2n}{n}^2}.$$
The function~$h(x)$ has denominator type~$\tau = [1,2,3,\ldots,2n]^2$
and is a~$G$-function of geometric origin over~$\mathbf{P}^1 \setminus \{0,16,\infty\}$.
Now the Hadamard products
$$f(x) = g(x) \star h(x)$$
lie in~$\HH$ and also generate an infinite dimensional space over~$\Q(x)$.
\end{proof}
Note that the numerology in this case corresponds to~$e^4 > 16$. 

\subsection{Algorithmic Questions}
Given an ODE of geometric origin, one can generally always give a bound on the denominator type. However, determining the precise growth of the denominators
appears to be difficult in general. Two enlightening examples can be given as follows.
In~\cite{Cooper}, the following example is considered.
Let
$$x = q \prod_{n=1}^{\infty} \left(\frac{(1 - q^{7n})}{(1 - q^n)} \right)^4,$$
which is a Hauptmodul for~$X_0(7)$. There is a corresponding uniformizer for~$X_0(7)^{+}$ as follows:
$$y = \frac{x}{1 + 13 x + 49 x^2}.$$
If we now take the weight~$2$ Eisenstein series
$$E = \frac{7 E_2(\tau) - E_{2}(7 \tau)}{1 - 7} = 
1 + 4 q + 12 q^2 + 16 q^3 + 28 q^4 + \ldots $$
and then write it in terms of~$y$, we get an order
three ODE without singularities
on~$\mathbf{P}^1 \setminus \{0,1/27,-1,\infty\}$,
given explicitly by~$\LL H_A(x) = 0$ with
\[
\begin{aligned}
 \LL &  = x^2 (1+x)(-1+27x) \frac{d^3}{dx^3}
+ 3x(-1 + 39 x + 54 x^2) \frac{d^2}{dx^2} \\
& \ \ + (-1 + 86 x + 186 x^2) \frac{d}{dx} 
+ 4(1+6x). \end{aligned} \]
If one considers the non-homogenous version~$\LL H(x) = -1$,
then there exists a holomorphic solution~$H(x)$ which is overconvergent beyond the cusp~$1/27$.
This solution is of modular origin, namely, there exists a weight~$4$ modular form whose triple (Eichler) integral
gives rise to~$H(x)/H_A(x)$  written in terms of the parameter~$q$. What is unusual about this example
is that the weight four form is meromorphic rather than holomorphic; it has the form~$h/E$ for
the unique Hecke eigenform~$h \in S_6(\Gamma_0(7),\Q)$,
and so~$h/E$ has poles away from the cusps. But perhaps more surprisingly, the form~$h/E$
is \emph{magnetic} in the sense of~\cite{magnetic}; that is, if~$h/E = \sum a_n q^n$ then~$n | a_n$ for all~$n$.
As a consequence, the holomorphic solution~$H(x)$, which from all appearances (and 
in light of the three integrations required to construct~$H(x)$) one should
expect to have denominator type~$\tau = [1,2,\ldots,n]^3$, actually has denominator type~$\tau = [1,2,\ldots,n]^2$.
This is not at all apparent from the ODE, and it is not clear whether there is an algorithm to compute this \emph{a priori}.
As a curious consequence, it also means that the irrationality of  the Ap\'{e}ry limit associated to Cooper's sequence is amenable
to our methods; however, since the corresponding constant appears to be~$\pi^2/42$, we have not pursued this!

The second example which highlights the difficulty in computing denominator
types is as follows.  Associated to Ramanujan's modular form~$\Delta = \sum \tau(n) q^n$, one can write down an ODE with a non-homogeneous solution
corresponding to the Eichler integral~$\sum \tau(n) n^{-11} q^n$. One expects
that the denominator type of the resulting function will be~$[1,2,\ldots,n]^{11}$. If one can prove that it is
\emph{not} of the form~$A^n [1,2,\ldots,n]^{10}$,  however, then one would have proven that there exist infinitely many
ordinary primes for~$\Delta$, a somewhat notorious open problem.

\subsection{The Gelfond--Schnirelman topic} \label{Chebyshev ex} 
Recall the numerology of the very basic special case of  Theorem~\ref{cap rationality}
on which we based our arithmetic characterization of the logarithm: the slit plane domain $\Omega := \C \setminus [1, \infty)$ has conformal mapping radius $\rho(\Omega,0) = 4$ with Riemann map
$\varphi(z) = 4z/(1+z)^2$, and it admits the transcendental analytic function $\log(1-x) = - \sum_{n=1}^{\infty} x^n/n$ whose denominator 
type can be expressed into the form~\eqref{s den type} with $r=1$ and $b_1 = 1$. We have $(2/3) \log{4} = 0.924196\ldots$ 
for the right-hand side 
of~\eqref{den gap} in this example.

This broaches another popular topic that was considered~\cite[\S~II]{ChudnovskyExt}, \cite[pages~493--494]{SelbergIntegral},
presumably for its methodological relevance, by several of the creators of the arithmetic theory of $G$-functions that we described in~\S~\ref{sec:G-function arithmetic}. This is the old idea of Gelfond and Schnirelman who observed in~1936 that some prime counting lower bound~$\pi(X) > (\log{2}) X/\log{X}$ for all~$X \gg 1$ follows at one stroke just by
remarking upon the pointwise~$\leq 4^{-n}$ integrand in $[1, \ldots, 2n+1] \int_0^1 (t-t^2)^n \, dt \in \NwithoutzeroA$. 
Using the functional bad approximability
property~\S~\ref{log discussion} of $\log(1-x)$ (the normality of the Hermite--Pad\'e table), the use of the prime number theorem in the 
proof of Theorem~\ref{logcharacterization} can be turned around to devise a Gelfond--Schnirelman style elementary proof of the integrated
prime counting function estimate $\int_1^X \psi(t) \, dt > (4/3)\log{2} \cdot X^2/2$, for all sufficiently big~$X$. The coefficient here is
slightly better than Chebyshev's $\log\left( 2^{1/2}3^{1/3}5^{1/5}30^{-1/30} \right) = 0.921292\ldots$, and now the point is that this proof is
not really new:  it is an isomorphic argument to  Bombieri, Nair, and Chudnovsky's opening estimates which they got from using the discriminant\footnote{In 
Bombieri's case, this inquiry led to the re-discovery of the Selberg integral and,  bearing with this for the true mathematical fruit, the historic proof of 
the Dyson--Mehta conjecture along with cases of the Macdonald conjectures. This is the story recounted in~\cite[\S\S~1.2, 1.3]{SelbergIntegral}.} polynomial in the multivariable method; see~\cite[page~94]{ChudnovskyExt}, and note that the Cauchy determinant from the proof of~\cite[Theorem~1]{Nair}
is none other than the Hankel determinant of the function~$\log(1-x)$ as applied to our Remark~\ref{remark on the slit plane}. As the logic of the prime number theorem can be reversed in every single one of our
arithmetic holonomy bounds, one cannot help but be curious about the denominator arithmetic in ``asymptotic near-misses'' of our holonomy bounds.

\subsection{A historical note and acknowledgments}

We originally conceived of our new approach to irrationality in 2020, starting with an easy proof of the irrationality of the $2$-adic avatar of~$\zeta(5)$ (which now we finally exposit in a companion paper~\cite{zeta5}), and even during that year we realized with the help of~\cite{Zagier} that the method could  apply in principle to $L(2,\chi_{-3})$. However, at the time, the holonomy bounds
that we could prove following~\cite[\S~VIII]{AndreG} were  totally insufficient, as explained in~\S~\ref{seriousintro}.
Our first serious improvements --- such as~\eqref{founding hol} --- over Andr\'e's holonomy bound~$\left( \sup_{\T} \log{|\varphi|} \right) / \left( \log{|\varphi'(0)|} - \tau \right)$  were still insufficient for this rather (as it seemed back then!) elusive application, but we
found them nonetheless to carry a certain asymptotic precision which was the key to the proof~\cite{UDC} 
of the unbounded denominators conjecture (the case~$\tau=0$ of algebraic functions). 
The paper~\cite{BostCharles} cites~\cite{UDC} as a significant influence. In turn, \cite{BostCharles}, which implicitly already has the bound~\eqref{BC integral 2}, has clearly been a crucial inspiration for our present paper, and (in part) it was by trying to synthesize our
ideas with those of~\cite{BostCharles} that lead to the optimal holonomy bounds here. 

The whole~\S~\ref{rmk_Nazarov} is due to Fedja Nazarov. We are grateful to him for explaining to us the precise analytic comparison between the Bost--Charles integral and the rearrangement integral. Remark~\ref{converse remark} is based on a discussion with Samuel Goodman. 

In addition, we would like to thank a number of people for conversations throughout the past four years
on ideas related to this paper, including Yves Andr\'e, Jean-Beno\^{\i}t  Bost, Alin Bostan, Fran\c{c}ois Charles, David and Gregory  Chudnovsky,
Tom Hutchcroft, 
Javier Fr\'esan,
Lars K\"uhne, Peter Sarnak, Umberto Zannier, Wadim Zudilin.

\appendix  

\section{Choosing a contour}  \label{contour choiceA}

\begin{figure}[!ht]
\begin{center}
  \includegraphics[width=85mm]{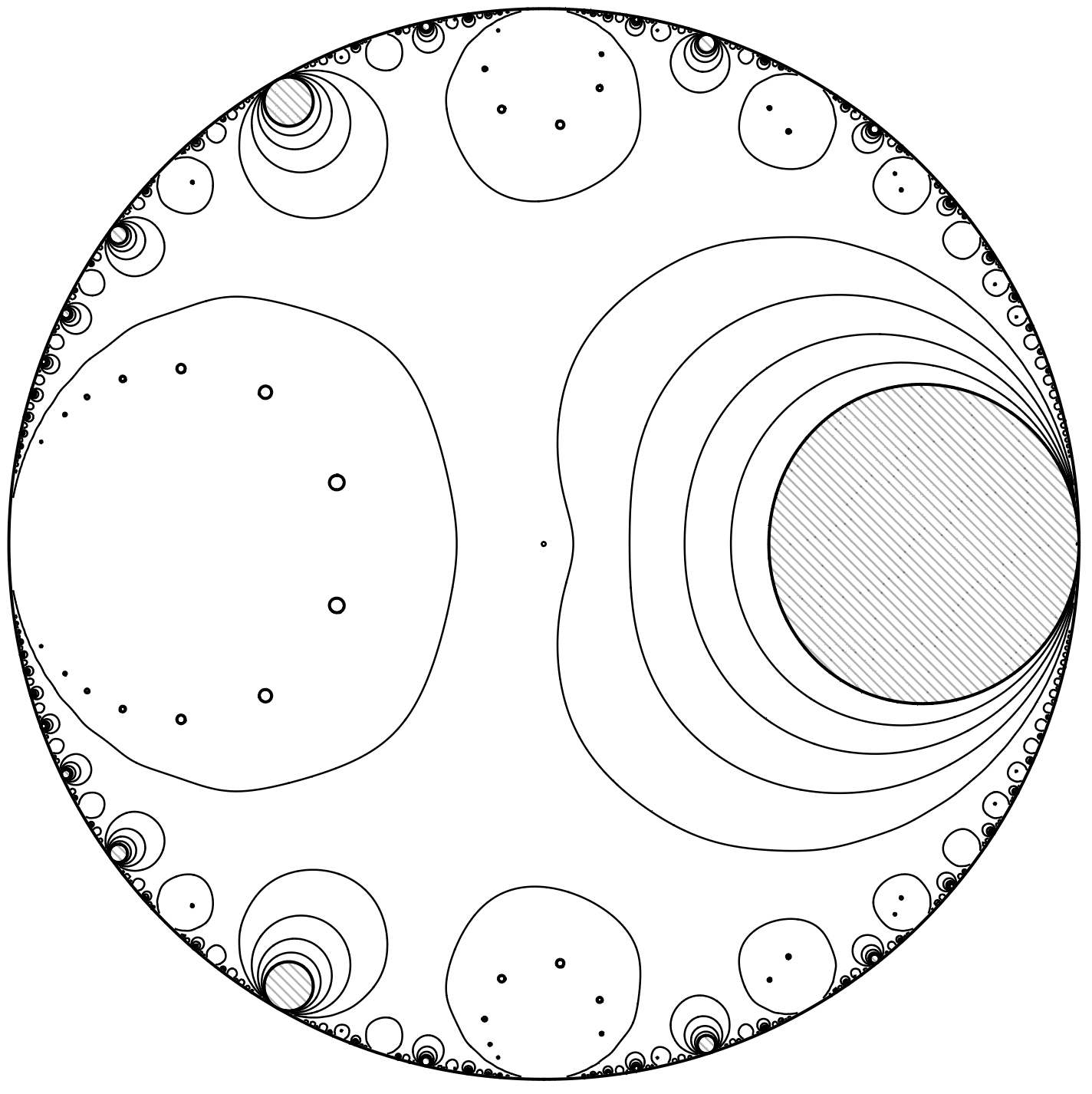}
  \caption{The shaded region consists of~$|z| \le 1$ with~$|h(z)| \ge e^{20}$ near
  the cusps~$z = e^{2\pi i p/q}$
  with~$q$ odd. The other curves are the level sets~$|h(z)| \in \{1,e^4,e^8,e^{12} e^{16}\}$,
  as well as the level sets~$|h(z)+1/72|=1/200$ around the
  preimages of~$-1/72$. The
  connected components of the 
  former level sets are distinguished from the latter by containing
 subsequences converging to the boundary.    }
  \label{top}
\end{center}
\end{figure}

Recall~(\ref{defofh}) the function~$h$ defined as follows:
\begin{equation}
h:=  \lambda + \frac{\lambda}{\lambda - 1} = - 256 q \prod_{n=1}^{\infty} (1 + q^{n})^{24} = -256 \cdot \frac{\Delta(2 \tau)}{\Delta(\tau)}, \quad q = e^{2 \pi i \tau}. 
\end{equation}
For any biholomorphic map~$\GGG: \D \rightarrow \Omega \subset \D$ with~$\GGG(0) = 0$,
let~$\varphi = h(\psi(x))$.
Our task is to choose a function~$\GGG$ for which:
\begin{enumerate}
\item The image of~$\GGG$ inside~$\D$ avoids all preimages of the point~$-1/72$
under~$h$ except for the one preimage~$0.0000541\ldots \in \R$. (This is the only preimage on the real line.)
\item 
The quantity (compare 
equation~(\ref{BCbound}))
\begin{equation}
\label{tomax}
 \frac{
 \displaystyle{ \iint_{\T^2}  \log{|\varphi(z) - \varphi(w)  |} \, \mv(z) \mv(w)  }
}
{ \log{ |256 \psi'(0)|} - \tau(\mathbf{b;e})}
\end{equation}
is as small as possible, for a  certain explicit constant
\[ \tau(\mathbf{b;e})   = \frac{16603}{3920}
= \frac{27}{80}  + \frac{191}{49}.\]
   \end{enumerate}
Here~$\log |\GGG'(0)|$ is the conformal radius of~$\Omega$.  The friction
here is that we want the denominator of~(\ref{tomax}) to be large, and so~$\Omega \subset \D$
to be large; at the same time,
the function~$h(z)$ has asymptotic growth 
\[\log |h(z)| \sim \frac{1}{2 \pi^2 q^2 (1 - |z|)}\]
as~$z$ varies  in a straight line
from~$0$ to the cusp~$e^{2 \pi i \alpha}$ where~$\alpha = p/q$
and~$q$ is odd. (All of this follows easily from the fact that~$h$ is a modular function of level~$\Gamma_0(2)$.)
In order to choose~$\Omega$, it is instructive first to examine
the topography of~$h$. 
The shaded region in
Figure~\ref{top} indicates the~$|z| \le 1$ for which~$|h(z)| \ge e^{20}$, the level sets~$|h(z)| \in \{1, e^4, e^8, e^{12}, e^{16}\}$,
and then finally level sets \( \left| h(z) + 1/72 \right| =1/200\)
around the preimages of~$-1/72$. (The two
types of level sets can be distinguished by whether any connected component
has a subsequence tending towards the boundary or not.)
\addtocounter{subsubsection}{1}
The basic idea for constructing a~$\Omega$ is to choose a circle
centered at the origin with
radius avoiding the preimages of~$-1/72$ in the vicinity of~$z = \pm i$, and then to (approximately) remove from~$\Omega$ the following:
\begin{enumerate}
\item The intersection of this circle with a horoball near~$z=1$.
\item Slits from this circle to the remaining preimages
 of~$-1/72$ along the (approximate) horoball in the vicinity of~$z=-1$.
\end{enumerate}
The Riemann mapping theorem guarantees
the existence of a~$\GGG(x)$ for any such region~$\Omega$.
However, we additionally want to choose~$\GGG(x)$
in an explicit form as follows
in order to be able to rigorously
estimate~(\ref{tomax}).
Thus in practice we choose simple explicit functions which approximate this
region. 
The construction of~$\Omega$ is very much bespoke, and it is completely unclear
(to us!) how to actually minimize~(\ref{tomax}) over all~$\Omega$, except
to say from our experience that we believe our construction is not a long way
from being optimal. 
Our ultimate choice of~$\Omega$ is displayed 
in Figure~\ref{level} (which also has a more detailed topographic map of~$\log |h(z)|$).

\subsection{Preliminaries on Lunes}
Let~$D(c,R)$ denote the disc of radius~$R$ centered at~$c$.
Fix~$c > 1$, and consider the map
$$f(z,c)
= z \cdot  \frac{(c^2 + 1) + (c^2 - 1) z}
{ (c^2 - 1) + (c^2 + 1) z}.$$
This map has the following property; it
is a conformal map from the lune~$L(c)$
 consisting of
 $$L(c) :=\D \setminus \D \cap D\left( - \frac{c^2 + 1}{c^2 - 1}, \frac{2c}{c^2 - 1}\right)$$
 to the unit disc sending~$z=0$ to~$0$.  That is, the unit disc minus the intersection of two discs which intersect at~$|z| = 1$ at right angles.  (This guarantees the existence of an explicit  and elementary conformal map).
 The ``innermost'' point of the right circle is the point
 $$ -\frac{c^2 + 1}{c^2 - 1} + \frac{2c}{c^2 - 1} =  - \frac{c-1}{c+1}.$$
In the limit~$c \rightarrow \infty$, this point tends to~$-1$, and the region~$L(c)$
tends to the entire disc, and~$f(z)$ tends to~$z$.
The function~$f(z)$ has an explicit inverse map as follows:
\begin{equation}
\label{firstgobble}
h(z,c) = \frac{ z(1 + c^2) -1 - c^2 +  \sqrt{(1 + c^2)^2 (1+z)^2 -16 c^2 z}}
{2 (c^2 - 1)},
\end{equation}
and hence the conformal radius of~$L(c)$ is
\begin{equation}
\label{firstgobblesize}
\frac{c^2-1}{c^2+1}.
\end{equation}

\subsection{Gobbles} \label{gobble contours}
We do not use the contours of this section in 
the proofs of Theorems~\ref{mainA} or~\ref{logsmain},
 having replaced them  by a combination of lunes with the slits considered in~\S~\ref{sec:slits} below.
 However, they are used in the proof 
 of Theorem~\ref{noextrassmall} in~\S~\ref{sec:completion three elements}.
 Moreover,
preliminary
versions of our argument did employ them, and they do provide convenient
contours on which to provide benchmarks for other examples,
and are also more flexible than our
 somewhat custom use of slits.

Suppose we wish to remove \emph{two} symmetrically opposite discs, not necessarily of equal sizes. There is no easily expressible simple conformal map in this case.
However, as a first approximation, we can first remove one disc, and then
remove the other. We define
$$\HHH(z,e,f) := h(-h(z,f),e).$$
Using equation~(\ref{firstgobble}), it has a somewhat
messy but completely explicit form.

\begin{df} The \emph{gobble} inside~$\D$ with parameters~$r \in (0,1]$, $e \in (1,\infty]$, and $f \in (1,\infty]$  is the image of~$\D = D(0,1)$
under the map
\[\HHH(r,e,f) : \D \to \C, \, z \mapsto r \cdot \HHH(z,e,f).\]
\end{df}

One key property is that, for a wide range of parameters, the gobble is visually indistinguishable from
the complement in~$D(0,1)$ of two discs, while at the same time being much more explicit and
thus easier to compute with.
From the explicit formula, we easily obtain: 
\begin{lemma} \label{confgobble}
The conformal radius of~$\HHH(r,e,f)(\D)$ centered at~$z=0$ is equal to
$$ r \cdot \frac{(e^2-1)(f^2-1)}{(e^2+1)(f^2+1)}.$$
\end{lemma}

\subsection{Slits} \label{sec:slits}
For a real number~$r \in (0,1)$, a conformal isomorphism
$$\left( \D,0 \right) \iso  \left( \D \setminus (-1,-r], 0\right)$$
is given by the function~$\Slit(z,r)$ defined by the following formula: 
{\small
\begin{equation}
\label{slit}
\begin{aligned}
& \frac{(r+1)^2 - 2(r-1)^2 z + (r+1)^2 z^2 
+ (1+r)(-1+z) \sqrt{(1+r)^2 - 2(1- 6 r + r^2) z + (1+r)^2 z^2}}{8 r z} \\
& = \frac{4r}{(1+r)^2} z + \frac{8r(1-r)^2}{(1+r)^4}z^2 +  \frac{4 (1-r)^2 r (3 -14r + 3 r^2)}{(1 + r)^6}z^3 + \ldots.
\end{aligned}
\end{equation}
}
 In particular, the conformal radius of~$\D \setminus (-1,-r]$ at the origin is equal to
 \begin{equation}
 \label{slitsize}
|\Slit'(0,r)| = \frac{4r}{(1+r)^2}.
 \end{equation} 
We include a sketch of the derivation of the map~\eqref{slit}. The starting
point is to remark that the rational transformation 
\begin{equation*}
\begin{aligned}
z  & \mapsto  z / (1+z)^2, \\
0 \mapsto 0; \quad -1 \mapsto \infty; \quad & 1 \mapsto 1/4; \quad -r \mapsto -r/(1-r)^2, \\ 
\textrm{with inverse } \,  z & \mapsto \frac{ 1-2z + \sqrt{1-4z} }{2z},
\end{aligned}
\end{equation*}
takes our slit disc $\D \setminus (-1,-r]$ conformally isomorphically onto the~$z \mapsto 1/z$ image of $ \P^1 \setminus \left[-(1-r)^2 / r, 4 \right]$.
Now a line segment $[A,B] \subset \R$ has transfinite diameter (\emph{$[\infty]$-capacitance}) equal to a quarter of its length.
This already proves the formula~\eqref{slitsize} on the conformal size, for the complement in~$\P^1 = \C \cup \{\infty\}$
of a contractible compact~$K \subset \C$ is a topological disc whose Riemann mapping radius from~$\infty$ is equal to the reciprocal 
of the transfinite diameter of~$K$. For the actual Riemann map~\eqref{slit}, we continue further
by observing that the  inverse Riemann map from~$\infty$ for a segment complement $\P^1 \setminus [A,B]$
in the Riemann sphere is given by a square root function:
\begin{equation*}
\begin{aligned}
\P^1 \setminus [A,B]    \iso \H,  \quad
z & \mapsto i \sqrt{(z-A) / (z-B)}, \quad   \infty \mapsto i, \\
 \textrm{with inverse }  \, z  & \mapsto \frac{A+Bz^2}{1+z^2}, \quad i \mapsto \infty. 
\end{aligned}
\end{equation*}
Following this through by the Cayley transform 
$
z \mapsto (z-i)/(z+i),   \, \H    \iso \D,
$
whose inversion is
$
z \mapsto i \, (1+z)/(1-z),
$
we arrive at the composed Riemann map~\eqref{slit}.

\subsection{Combining multiple slits and lunes} \label{mixing} 
Suppose we wish to remove \emph{four} slits and, additionally, a lune.
 There is no easily expressible simple conformal map in this case.
 However,  both the conformal maps in~(\ref{firstgobble}) and~(\ref{slit}) have the property that
 they are well--approximated by the identity map~$z \mapsto z$ for ``most'' points in the circle (namely, the points away from
 the lune and the slits, respectively).
 Thus one very primitive way to construct such maps is simply to \emph{compose} these
 maps in succession. 
 To this end, we consider the following function:
 {\small
\[
 G(z) := -R  \cdot h \left(- e^{2 \pi i \theta_1} \cdot \Slit\left( e^{2 \pi i \theta_2}  \cdot  \Slit\left( e^{2 \pi i \theta_3}
 \Slit\left( e^{2 \pi i \theta_4}   \cdot \Slit\left( z,r_1\right),r_2\right),r_3\right),r_4\right),c \right).
\]
 }
 We fix the first parameter~$R = 77/100$ to ensure that the initial circle only contains preimages of~$-1/72$ in the horoball
 around~$-1$, and indeed that there are only~$4$ such preimages that we need to exclude.
 We also fix the lune parameter~$c = 75/10$ which
 measures the (approximate)
  horoball we remove near~$z=1$.
 The angle parameters~$\theta_i$ allow us to ``line up'' the slits so that they include these preimages,
 and the length parameters~$r_i$ allow us to minimize the lengths of these slits so they do not go
 beyond the preimages we wish to exclude. 
 The final choice of parameters is as follows:
\begin{equation}
\begin{aligned}
& & R &= \frac{77}{100}, & \quad c &= \frac{75}{10},   \\
r_1 &= \frac{91}{100}, & \quad r_2 &= \frac{6188}{10000}, & \quad r_3 &= \frac{55515}{100000}, & \quad r_4 &= \frac{772}{1000}, \\
\theta_1 &= \frac{7977}{100000}, & \quad \theta_2 &= \frac{11543}{100000}, & \quad \theta_3 &= \frac{3525}{100000}, & \quad \theta_4 &= - \frac{783}{10000}.
\end{aligned}
\end{equation}
 These parameters are chosen from an \emph{ad hoc} computation making the ends of the slits
 as close to the four parameters as possible. Since it is not possible (numerically) to choose these
 parameters so that the bad preimages lie exactly on these slits, we finally define
 \begin{equation}
 \label{defGgob} \GGG(z) = G \left( \frac{995}{1000} \cdot z \right).
 \end{equation}
 By restricting to this open disc, we are removing not simply (curved) slits but open regions, which enables
 one to easily prove that the bad preimages are excluded.
 It is simple enough to compute that the conformal radius of~$\GGG(z)$ (using~(\ref{firstgobblesize}) and~(\ref{slitsize})) and
 it is equal to
 \begin{equation} \label{exactconformal}
\begin{aligned}
  |\GGG'(0)| & = \frac{995}{1000} \cdot  
  R \cdot \frac{c^2-1}{c^2 + 1} \cdot \prod_{i=1}^{4} \frac{4 r_i}{(1 + r_i)^2} \\
& = \frac{5448339453535586608000000000}{8658833407565631122430056127}
 = 0.6292232680 \ldots
 \end{aligned}
 \end{equation}
Recall that~$h: \D \rightarrow \C$ is given by~$-256 q \prod_{n=1}^{\infty} (1 + q^n)^{24}$ with~$q \in \D$.
The parameters~$r_i$ and~$\theta_i$
are chosen above to ensure that the following holds:

\begin{lemma} \label{choicephiA}
 The map~$\varphi := h \circ \GGG: \Db \rightarrow \C$ has a unique preimage of~$-1/72$.
\end{lemma}

\begin{proof} There is a preimage with~$x =0.0000541829\ldots$. But one can determine (to any precision) the other preimages by passing back
to~$\H$ and then the preimages are obtained by the action of~$\Gamma_0(2)$. 
The preimages in the region (approximating a horoball) near~$z=\pm i$ have absolute value at least
$$0.782767 \ldots > R  = \frac{77}{100}.$$
The closest other preimages lie near the horoball at~$z = -1$; but
the precise choice of the parameters~$r_i$ and~$\theta_i$ ensure that they lie outside
 image of~$\GGG$ as can be confirmed
by a simple numerical computation. 
\end{proof}

\begin{figure}[!h]
\begin{center}
  \includegraphics[width=85mm]{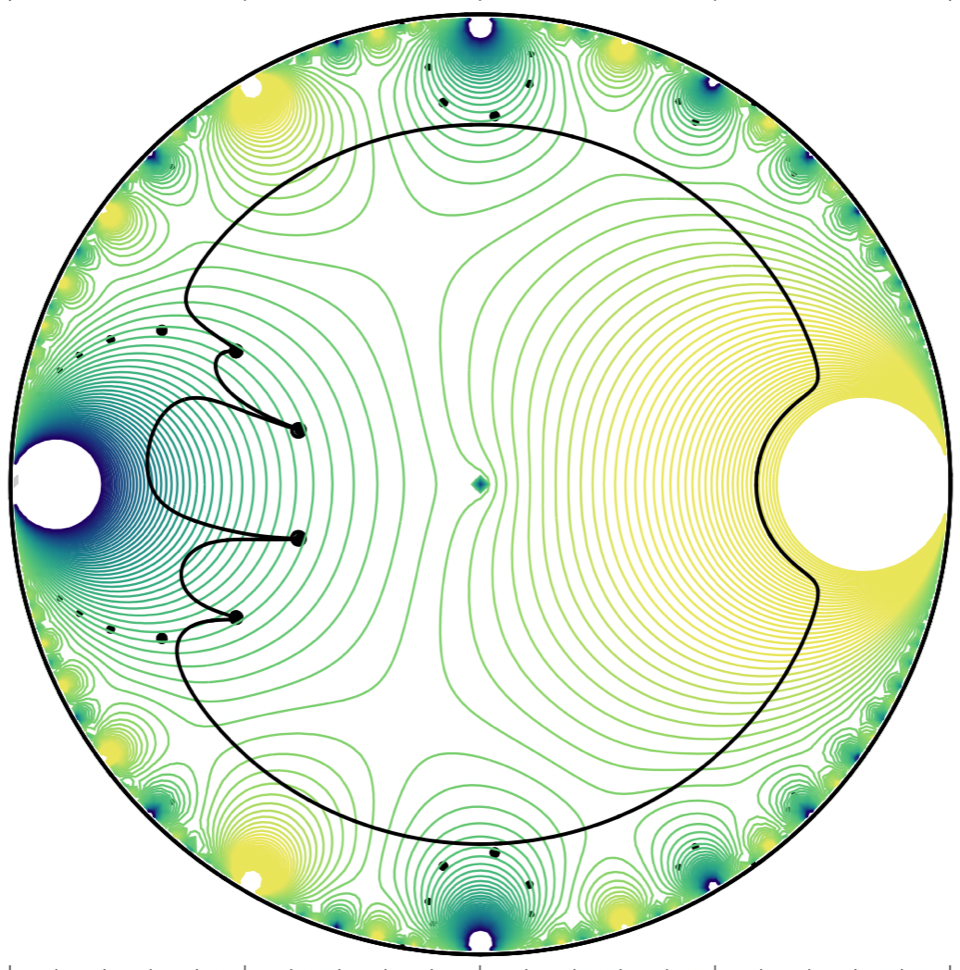}
\caption{The image of~$|z|=1$ under~$\GGG(z)$ together with preimages of~$-1/72$ under~$h: \D \rightarrow \C$, together with level sets for~$\log |h|$ at values in an arithmetic progression;
 the color scheme transitions between yellow for large positive values of  $\log |h(z)|$
  and blue for large negative values.}
\label{level}
\end{center}
\end{figure}

The contour~$\psi(\T)$ is drawn in Figure~\ref{level}. The asymmetry is due to our using four successive 
compositions of single slit maps, rather than having all four slits taken out at once. 

\addtocounter{subsubsection}{1}  

\begin{remark} \label{choiceofconstants}
Our choice of constants reflects merely the principle of finding an example ``which works''
rather than is the most aesthetically pleasing. 
There is no doubt some scope for improvement but since it is not necessary we
have not tried to optimize these choices   ---  we expect improvements in either respect would anyway be quite modest.~\endofremark
\end{remark}

\subsection{A contour for the \texorpdfstring{$L(2,\chi_{-3})$}{L2chi3} problem}   \label{numerical integration}
With our specific choice of~$\GGG(z)$ as in Definition~\ref{defGgob}, we now define
$$\varphi(z) = h(\GGG(z)),$$
which is uniformly continuous on~$\D$, and which is explicit enough as to be amenable
to rigorous numerical estimates.
We then finally obtain the estimate
$$
 \displaystyle{ \iint_{\T^2}  \log{|\varphi(z) - \varphi(w)  |} \, \mv(z) \mv(w)  }
 = 11.844\ldots $$ 
 and thus~(\ref{tomax}) is bounded above by
\begin{equation}
\label{lessthan14}
\frac{
11.845
}
{
\displaystyle{\log \left( 256  \cdot
\frac{5448339453535586608000000000}{8658833407565631122430056127}
\right)
- \left(\frac{27}{80}  + \frac{191}{49}\right)
}
}
=  
13.9938\ldots. 
< 14.
\end{equation}

\begin{remark}[Bounds
without integrations] \label{improvement}
Even with all our improvements, the best bound we could achieve
before integrations, for either Theorems~\ref{mainA} or~\ref{logsmain}, 
 was also above~$9$;
we give some of the numerics now.
Consider the following basic application of Theorem~\ref{basic main}.
We consider the functions~$A_i(x)$ for~$i=1,\ldots,9$ 
in the~$\mathbf{P}^1 \setminus \{0,1,\infty\}$ domain. (Here the first five functions
are given explicitly
in~\S~\ref{sec:purefunctions}, and the four functions~$A_6(x), \ldots, A_9(x)$
correspond to~$B_6(z), \ldots, B_9(z)$ via the transformations of that section.
Note that these last four functions only exist if there is a~$\Q$-linear
relation between our three periods.) Consider Theorem~\ref{basic main} for
$$
\mathbf{b} := \left( \begin{array}{lllllllll} 0 &  1 & 1 & 1  & 1 & 1 & 1 & 1 & 1 \\ 0 & 0 & 1 & 1  & 1 & 1 & 1 & 1 & 1     \end{array} \right)^{\mathrm{t}},
$$
hence 
$$
\sigma_1 =0,  \, \sigma_2 = 1; \quad \sigma_3 = \cdots = \sigma_9 = 2,
$$
and
$$
\tau(\mathbf{b}) = \frac{1 \cdot 0 + 3 \cdot 1 + (5+7+9+11+13+15+17) \cdot 2 }{81}  =  \frac{157}{81}.
$$
As explained in~Remark~\ref{equivalence},
if we use the same contour as given
in Definition~\ref{defGgob} except  pulled back to the~$X(2)$ domain,
both the integral and the conformal radius terms are doubled. Equivalently,
they remain the same and the~$\tau$ term is halved.
Hence the corresponding bound we obtain in this case is: 
\begin{equation}
\label{biggerthannine}
\frac{11.844\ldots }
{ 5.081\ldots-  2 \cdot 157/81}
 = 9.833  \ldots  < 10.
\end{equation}
which comes close but is not a contradiction because this term is not less than~$9$. While
this can be refined slightly (using the Bost--Charles integral and modifying the contour),
it seems unlikely that one may reach a direct contradiction by our methods
without involving added integrations; see Examples~\ref{Ex_BCconv9} and \ref{Ex-easyconv}.~\endofremark
\end{remark}

\subsection{A contour for the logarithm problem}  \label{contour choiceC}
We could literally use the same contour as above
to complete the proof of Theorem~\ref{logsmain}, except with a somewhat worse
constant. Following the arguments
of Section~\ref{sec:location}, it would suffice
to find the~$\varepsilon$ such that the image of~$\varphi$ above
excludes the regions where~$z$ is not too small and~$h(z)$ lies in~$D(0,\varepsilon^2/16)$
and also where~$h(z)$ lies outside~$D(0,\varepsilon^{-1})$. This leads to a choice
of~$\varepsilon$ somewhere between~$10^7$ and~$10^8$.
However, a compromise between optimizing over various conformal maps
and using the same map as above is just to write down simple lunes:
If we take
$$\GGG(z) =  -\frac{3}{4} h\left(z, \frac{23}{10}\right)$$
of conformal radius~$1287/2516$,
then the image of~$\GGG(z)$ avoids all the required discs
as well as the regions where~$|h(z)| \ge 10^6$ and~$|h(z)| \le 10^{-12} 2^{-4}$
(except for the preimages near~$z=0$). In this case, we obtain
the bound
\begin{equation} \label{logBC}
 \displaystyle{ \iint_{\T^2}  \log{|\varphi(z) - \varphi(w)  |} \, \mv(z) \mv(w)  }
  \sim 9.963 \ldots
< 10,
\end{equation}
and we have the very easy bound
\begin{equation}  \label{logproductnumber}
\frac{
10
}
{
\displaystyle{
\log \left( 256  \cdot \frac{1287}{2516}  \right)}
- \frac{1032659}{242760}
} =  16.103\ldots  < 17
\end{equation}
which we use for the proof of Theorem~\ref{logsmain} in~\S~\ref{sec:proofC}.
In comparison, we may also estimate the rearrangement integral
\begin{equation}
 \label{logrearrange}
  \int_0^1 2t \cdot ( \log{|\varphi(e^{2\pi i t})|} )^* \, dt \sim 9.972 \ldots 
  \end{equation}
which also suffices to prove Theorem~\ref{logsmain}, this time via
Theorem~\ref{main:elementary form} with only the trivial partition of~$[0,m]$.

\begin{figure}[!h]
\begin{center}
  \includegraphics[width=85mm]{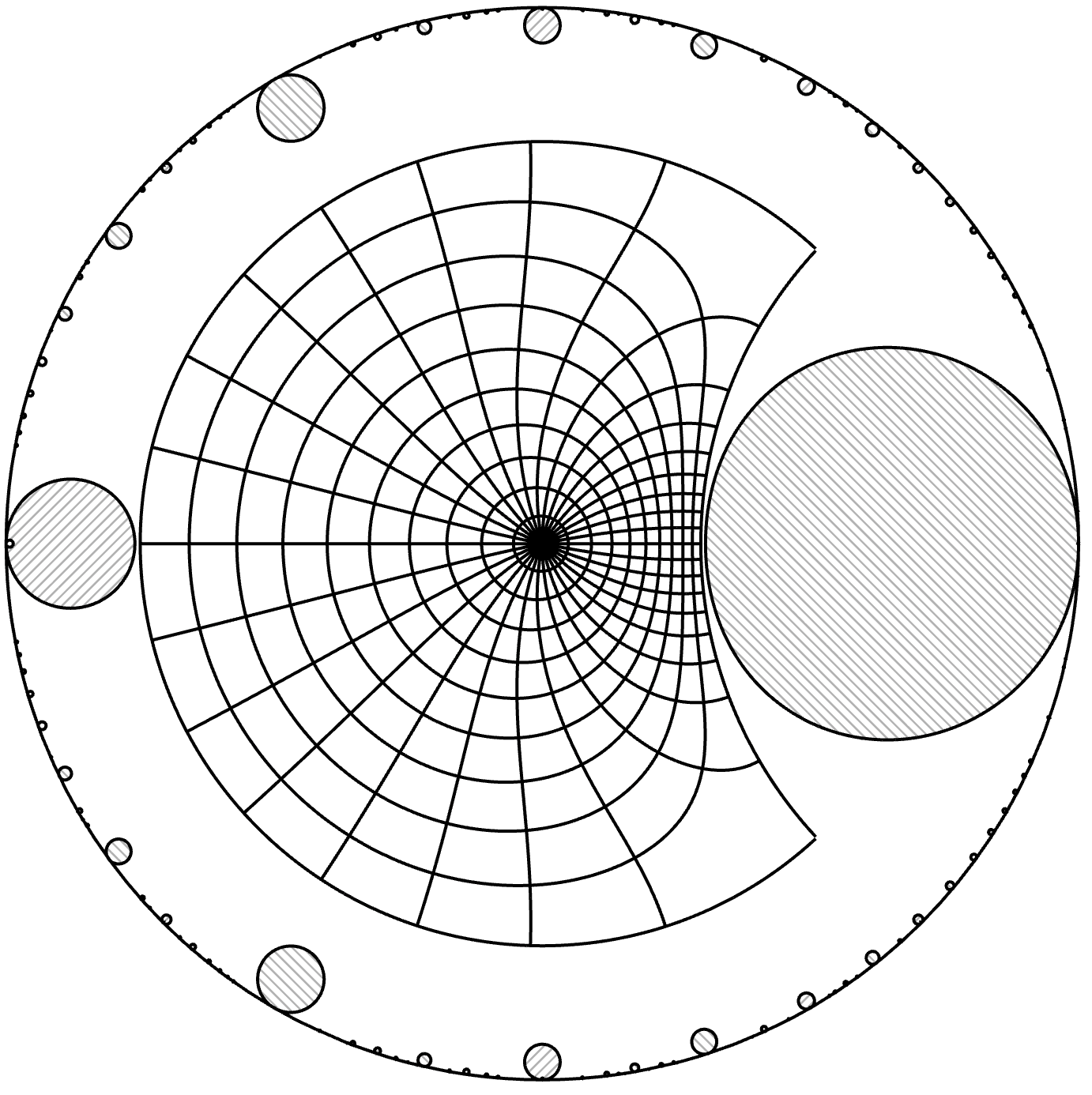}
\caption{The image of~$|z|=1$ together with the images
of the level sets~$|z|=k/10$ and~$\arg(z)=2 \pi k/32$ for integers~$k$,
together with the
regions with~$|h| \ge 10^6$ and~$|h| < 10^{-12} 2^{-4}$ (distinguished
by the angle of the shading)}
\label{logfig}
\end{center}
\end{figure}

A graph of the image of~$\GGG(z)$ together with the regions
where~$|h| \ge 10^6$ and~$|h| \le 10^{-12} 2^{-4}$ is given in
Figure~\ref{logfig}

\addtocounter{subsubsection}{1}

\section{A dynamic box principle} \label{app:PerelliZannier}  In this appendix we give a short new proof of the basic holonomy bound
\begin{equation} \label{basic basic}
m  \leq \frac{2 T(\varphi)}{ \log{|\varphi'(0)|} - b_1 - \cdots - b_r},
\end{equation}
under the condition of the positive denominator,
for a $\Q(x)$-linearly independent set of formal functions $f_1, \ldots, f_m \in \Q \llbracket x \rrbracket$ of the types~\eqref{gen den form} and
such that $\varphi; \varphi^* f_1, \ldots, \varphi^* f_m \in \mathcal{M}(\Db)$ are simultaneously meromorphic on a neighborhood of the closed unit disc~$\Db$. 
Here, 
\begin{equation}  \label{nchar}
T(\varphi) := \int_{\T} \log^+{|\varphi|} \, \mv  + \sum_{ \substack{ \rho \in \D \\ \textrm{poles of $\varphi$} }} \log\frac{1}{\rho} 
\end{equation}
is the \emph{Nevanlinna characteristic} of the meromorphic mapping~$\varphi$ (the meromorphic poles being taken with their multiplicities). 

This is based on the idea of Perelli and Zannier~\cite{PerelliZannier} with a dynamic box principle such as they formulate with their Lemma~1 of {\it loc.cit}. 
It may be considered as a more elementary form of Bost's technique in~\S~\ref{new slopes}, to which it is both an introduction and an alternative, and to our companion
paper~\cite{zeta5}, where these ideas are pursued further. We divide the proof into three steps according to the dissection in~\S~\ref{sec:ideas outline}. 

\subsection{Evaluation module}  Suppose we have a $\Q(x)$-linearly independent set
of functions $f_1(x), \ldots, f_m(x) \in \Q\llbracket x \rrbracket$
of the type~\eqref{gen den form} such that $\varphi(z) \in \C \llbracket z \rrbracket$, as well as each power series $f_i(\varphi(z)) \in \C \llbracket z \rrbracket$, are
germs of meromorphic functions on a neighborhood of the closed unit complex disc
$|z| \leq 1$. (Having this slightly bigger disc is no loss of generality upon replacing $\varphi(z)$ by $\varphi(\rho z)$ for 
some $\rho < 1$ still having $\rho |\varphi'(0)| > e^{b_1+\ldots + b_r}$.)  We introduce two positive integer parameters $D$ and $T$, and we consider the collection
\begin{equation}
\mathcal{I}_D(T) :=  \left\{ (Q_1, \ldots, Q_m) \in \Z[x] \, : \,  Q_i(x) = \sum_{j=0}^{D-1} c_{i,j} x^j, \, c_{i,j} \in [0,T) \cap \Z   \right\},
\end{equation}
of cardinality
\begin{equation} \label{input size}
\#\mathcal{I}_D = T^{mD}. 
\end{equation}
By the assumed $\Q(x)$-linear independence of the~$m$ formal power series $f_i(x) \in \Q\llbracket x \rrbracket$, the $\Z$-module \emph{evaluation map}
\begin{equation} \label{eval}
\psi_D \, : \, \Z[x]_{\deg < D}^{\oplus m} \hookrightarrow \Q \llbracket x \rrbracket, \qquad  (Q_1, \ldots, Q_m) \mapsto \sum_{i=1}^{m} Q_i(x) f_i(x) \in \Q\llbracket x \rrbracket
\end{equation}
is \emph{injective}. Hence, the image 
$$
\mathcal{O}_D := \psi_D(\mathcal{I}_D) \subset \Q \llbracket x \rrbracket
$$
under this map also has cardinality
\begin{equation} \label{output cardinality}
\# \mathcal{O}_D = \#\mathcal{I}_D = T^{mD}.
\end{equation} 
 In the free $\Z$-module~\eqref{eval} of rank~$mD$, we define the \emph{vanishing filtration jumps} 
\begin{equation} \label{rank jumps}
r_D^{(n)} := \dim_{\Q} \left\{ \left( \psi_D^{-1}\left( x^n \Q \llbracket x \rrbracket \right) \otimes \Q \right) \big/
\left( \psi_D^{-1}\left( x^{n+1} \Q \llbracket x \rrbracket \right) \otimes \Q \right)   \right\}.
\end{equation}
They are in $\{0,1\}$, because the linear injective map $\psi_D$ induces a linear injection
$$
 \left( \psi_D^{-1}\left( x^n \Q \llbracket x \rrbracket \right) \otimes \Q \right) \big/
\left( \psi_D^{-1}\left( x^{n+1} \Q \llbracket x \rrbracket \right) \otimes \Q \right)  \hookrightarrow x^n \Q \llbracket x \rrbracket \big/ x^{n+1} \Q
\llbracket x \rrbracket \cong \Q \cdot x^n
$$
into a one-dimensional $\Q$-vector space. On the other hand, we have 
\begin{equation} \label{total rank}
\sum_{n=0}^{\infty} r_D^{(n)} = \dim_{\Q} \psi_D^{-1} \Q \llbracket x \rrbracket = mD.
\end{equation}
Hence there is a size-$mD$ set of possible $x=0$ vanishing orders 
\begin{equation}
\begin{aligned}
&  \left\{ n \in \NwithzeroB \, : \,   \exists (Q_1,\ldots, Q_m) \in \Z[x]_{\deg < D}^{\oplus m}, \, 
\mathrm{ord}_{x=0} \left(  \sum_{i=1}^{m} Q_i(x) f_i(x) \right) = n \right\} \\  \label{filtration jumps}
& = \big\{  0 \leq  u(1) < u(2) < \cdots <  u(mD)  \big\}
\end{aligned}
\end{equation}
for our auxiliary functions. They depend only on the module~\eqref{eval}  ---  in other words,
on $f_1, \ldots, f_m$ and the parameter~$D$,  ---  but not on the parameter~$T$, which 
remains free to select in the following. (The parameter~$T$ will be taken to be any 
sufficiently big integer in dependence of the filtration jumps~\eqref{filtration jumps}.)

\medskip

This fulfills step~(i) of~\S~\ref{sec:ideas outline}. 

\subsection{Box principle}  \label{dynamic box}  For step~\ref{dirichletbox}, we measure up the tendency of the Taylor series of the auxiliary
function~$F(x)$ to depend recursively on its string of initial coefficients under
 the critical condition $|\varphi'(0)| > e^{b_1 +\ldots + b_r}$.
 
  We can upper-estimate the output
cardinality $\#\mathcal{O}_D$ by a product $\prod_{p=1}^{mD} \gamma_p$, where $\gamma_p$
is an upper estimate on the largest possible number of distinct $x^{u(p)}$ coefficients $\beta_{u(p)} \in \Q$ in any set of output functions $F(x) = \sum_{k=0}^{\infty} 
\beta_k x^k$ that share a common string $(\beta_0, \beta_1, \ldots, \beta_{u(p)-1})$ for their preceding coefficients:
\begin{equation*}
\begin{aligned}
&  \qquad \forall (\beta_0, \ldots, \beta_{u(p)-1}) \in \mathbb{Q}^{u(p)},  \\
& \# \left\{  \beta \in \mathbb{Q} \, : \, \exists  (Q_1, \ldots, Q_m) \in \mathcal{I}_D, \, \sum_{i=1}^{m} Q_i(x) f_i(x) - \sum_{k=0}^{u(p)-1} \beta_k x^k  \right. \\
&  \qquad \left. \vphantom{\sum_{i=1}^{m} Q_i(x) f_i(x)} = \beta x^{u(p)} + O\big(x^{u(p)+1}\big)  \right\} \leq \gamma_{u(p)}
\end{aligned}
\end{equation*}
At this point, the integrality condition~\eqref{gen den form} is used to remark that all such rational numbers $\beta \in \Q$ belong in fact 
to a $\Z$-module given by the requisite denominators type: 
$$
\beta \in \frac{1}{[1,\ldots, b_1 u(p)] \cdots [1,\ldots, b_r u(p)]} \, \Z. 
$$
Hence, if $A_p \in \R^{> 0}$ is such that any two such coefficients $\beta$ differ by some real number in $[-A_p, A_p]$, then we can take
\begin{equation} \label{conditioning}
\gamma_p := 1 + 2A_p
 \cdot [1,\ldots, b_1 u(p)] \cdots [1,\ldots, b_r u(p)]  = 1 + A_p \cdot  e^{(b_1 +\ldots + b_r)u(p) + o(u(p)) + O(1)}  \end{equation}
 as a total bound on the output possibilities of $\beta_{u(p)}$ given $(\beta_0, \ldots, \beta_{u(p)-1})$. This step 
 emulates the dynamic box principle of Perelli and Zannier~\cite[\S~2 Lemma~1]{PerelliZannier}. 
 
 \subsection{Diophantine analysis of the lowest order coefficient}  \label{iii coeff} The bound for~$A_p$ comes analytically from using the simultaneous
 meromorphic uniformization map~$\varphi$. Consider $\varphi = v/u$ any representation of $\varphi$ as the quotient of two
 convergent power series $u$ and $v$ on $\Db$ such that $u(0) = 1$. Let $h$ be a convergent power series on $\Db$ such that $h(0) = 1$ and
 $hf_i$ is holomorphic for each $i = 1, \ldots, m$. Let us write $n := u(p)$ for the following. Any two output functions~$F_1(x)$ and $F_2(x)$ as above whose 
 $x=0$ Taylor series coincide up to $O(x^n)$, and whose respective $x^n$ coefficients are $\beta^{(1)}$ and $\beta^{(2)}$, will have 
 $$
 V(z) :=h(z) u(z)^D \cdot \left( F_1(\varphi(z)) - F_2(\varphi(z)) \right)
 $$
  holomorphic (convergent) on
 a neighborhood of the closed unit disc $\Db$, and with leading order term $\varphi'(0)^n (\beta^{(1)} - \beta^{(2)}) z^n + O(z^{n+1})$ expressible as a 
 $|z| = 1$ contour integral by Cauchy's formula: 
 \begin{equation} \label{Cauchy int}
 \varphi'(0)^n (\beta^{(1)} - \beta^{(2)}) = \int_{\T} \frac{V(z)}{z^{n+1}} \, \mv. 
 \end{equation}
Estimating by the supremum of the integrand, we can take for our~$A_p$ the upper bound on the bottom row of
 \begin{equation}   \label{Cauchy est}
 \begin{aligned}
 \left| \beta^{(1)} - \beta^{(2)} \right|  & \leq 
 |\varphi'(0)|^{-n} \cdot \sup_{\T} |V| \\  
 &  \leq  |\varphi'(0)|^{-u(p)} \cdot  T \, \left( \sup_{\T}  \max \left( |u|, |v| \right) \right)^D \cdot m D \cdot \sup_{\T} |h \cdot \varphi^* f_i | =: A_p, 
  \end{aligned}
 \end{equation} 
 used with $n := u(p)$.
 
We  get for the $\prod_{p=1}^{mD} \gamma_p$ output possibilities the upper estimate
$T^{mD} = \#\mathcal{O}_D$ 
\begin{equation*}
\begin{aligned}
\leq \prod_{p=1}^{mD} & \Bigg\{ 1 + T \exp \Bigg( - \big( \log{|\varphi'(0)|} - 
 b_1 - \cdots - b_r  +o(1) \big) u(p) \Bigg. \\
&\quad \Bigg. \vphantom{\exp \Bigg( - \big( \log{|\varphi'(0)|} - 
 b_1 - \cdots - b_r  +o(1) \big) u(p) \Bigg)} + 
D \, \sup_{\mathbb{T}} \log{\max(|u|,|v|)}
 + \log D + O_{m,h}(1) \Bigg) \Bigg\}.
\end{aligned}
\end{equation*}
 At this point, we look at the last inequality asymptotically in~$T \to \infty$, or more concretely, we select a~$T$ so big that all~$mD$ factors of the product are $\geq 2$. 
 Using the trivial inequality $1+ x \leq 2x$ for $x \geq 1$ and canceling the common ensuing $T^{mD}$ from both sides, we 
 get (after taking the logarithm)
 \begin{equation*}
 \begin{aligned}
- \big(   \log{|\varphi'(0)|} - 
 b_1 - \cdots - b_r \big) (1-o(1))&  \sum_{p=1}^{mD}  u(p)  + 
mD^2 \, \sup_{\T} \log{\max(|u|,|v|)} + O(D\log D) \\  & \geq - mD \log{2} - O_{m,h}(D). 
\end{aligned}
 \end{equation*}
 As $0 \leq u(p) < u(2) < \cdots < u(mD)$ are a strictly increasing sequence of nonnegative integers, we have 
\begin{equation}
\sum_{p=1}^{mD}  u(p)  \geq \sum_{n=0}^{mD-1} n = \binom{mD}{2}. 
\end{equation}
 We derive 
 \begin{equation}
(1-o(1)) \binom{mD}{2}  \big(   \log{|\varphi'(0)|} - 
 b_1 - \cdots - b_r \big)   \leq mD^2 \, \sup_{\T} \log{\max(|u|,|v|)}  + O_{m,h,1}(D), 
 \end{equation}
 which in the $D \to \infty$ asymptotic filters down to the \emph{arithmetic holonomy bound}
 \begin{equation}
 m  \leq \frac{2   \sup_{\T} \log{\max(|u|,|v|)}  }{ \log{|\varphi'(0)|} - 
 b_1 - \cdots - b_r }. 
 \end{equation}
 This is true for any meromorphic quotient representation $\varphi = v/u$ with $u(0) = 1$. 
 A well-known lemma of Nevanlinna (cf.~\cite[\S~VII.1.4]{Nevanlinna} or~\cite[\S~VII.5]{Goluzin}), 
 based on the canonical Blaschke products and the canonical decomposition $\log = \log^+ - \log^-$
in the Poisson--Jensen formula, constructs \emph{on the open disc $\D$} a quotient representation $\varphi 
= v/u$ with $u(0) = 1$ and with both $\sup_{\D}{|u|}$ and $\sup_{\D}{|v|}$ bounded by $\exp \left( T(\varphi)  \right)$. 
Dilating the radius a little bit, we get for any $\varepsilon > 0$ a quotient representation $\varphi = v/u$, now
on some neighborhood of the closed disc $\Db$ as required in the above analysis, with
$\sup_{\T}{|u|}$ and $\sup_{\T}{|v|}$ bounded by $\exp \left( T(\varphi)  + \varepsilon \right)$.
This concludes the proof of the bound~\eqref{basic basic}.      \hfill{$\square$}

\hbadness=10000
\vbadness=10000

\nocite{HermiteOeuvres,DworkDiff}
\bibliographystyle{amsalpha}
\bibliography{L2chi}

\renewcommand{\MR}[1]{}
\providecommand{\bysame}{\leavevmode\hbox to3em{\hrulefill}\thinspace}
\providecommand{\MR}{\relax\ifhmode\unskip\space\fi MR }
\providecommand{\MRhref}[2]{%
  \href{http://www.ams.org/mathscinet-getitem?mr=#1}{#2}
}
\providecommand{\href}[2]{#2}
\begin{thebibliography}{BGMN05}

\bibitem[AB95]{AbbesBouche}
Ahmad Abbes and Thierry Bouche, \emph{Th\'{e}or\`eme de {H}ilbert-{S}amuel
  ``arithm\'{e}tique''}, Ann. Inst. Fourier (Grenoble) \textbf{45} (1995),
  no.~2, 375--401. \MR{1343555}

\bibitem[AB97]{AndreBaldassari}
Yves Andr\'e and Francesco Baldassarri, \emph{Geometric theory of
  {$G$}-functions}, Arithmetic geometry ({C}ortona, 1994), Sympos. Math., vol.
  XXXVII, Cambridge Univ. Press, Cambridge, 1997, pp.~1--22. \MR{1472489}

\bibitem[Ada87]{Adams}
Colin~C. Adams, \emph{The noncompact hyperbolic {$3$}-manifold of minimal
  volume}, Proc. Amer. Math. Soc. \textbf{100} (1987), no.~4, 601--606.
  \MR{894423}

\bibitem[Ami75]{Amice}
Yvette Amice, \emph{Les nombres {$p$}-adiques}, Collection SUP: ``Le
  Math\'{e}maticien'', vol.~14, Presses Universitaires de France, Paris, 1975,
  Pr\'{e}face de Ch. Pisot. \MR{0447195}

\bibitem[And89]{AndreG}
Yves Andr\'{e}, \emph{{$G$}-{F}unctions and {G}eometry}, Aspects of
  Mathematics, no. E13, Friedr. Vieweg Sohn, Braunschweig, 1989. \MR{0990016}

\bibitem[And96]{AndreGtranscendence}
\bysame, \emph{{$G$}-fonctions et transcendance}, J. Reine Angew. Math.
  \textbf{476} (1996), 95--125. \MR{1401697}

\bibitem[And00a]{AndreGevreyI}
\bysame, \emph{S\'{e}ries {G}evrey de type arithm\'{e}tique. {I}.
  {T}h\'{e}or\`emes de puret\'{e} et de dualit\'{e}}, Ann. of Math. (2)
  \textbf{151} (2000), no.~2, 705--740. \MR{1765707}

\bibitem[And00b]{AndreGevreyII}
\bysame, \emph{S\'{e}ries {G}evrey de type arithm\'{e}tique. {II}.
  {T}ranscendance sans transcendance}, Ann. of Math. (2) \textbf{151} (2000),
  no.~2, 741--756. \MR{1765708}

\bibitem[And04]{Andre}
\bysame, \emph{Sur la conjecture des {$p$}-courbures de {G}rothendieck--{K}atz
  et un probl\`eme de {D}work}, Geometric Aspects of Dwork Theory, vol. I, de
  Gruyter, Berlin, 2004, pp.~55--112. \MR{2023288}

\bibitem[Ang19]{Angelesco}
Aurel Angelesco, \emph{Sur deux extensions des fractions continues
  alg\'ebriques}, {Comptes Rendues Hebdomadaires des S\'eances de l'Acad\'emie
  des Sciences} \textbf{18} (1919), 262--263.

\bibitem[Ape79]{AperyHistoric}
Roger Apery, \emph{Irrationalit\'e de $\zeta 2$ et $\zeta 3$}, Ast\'{e}risque
  (1979), no.~61, 11--13, Luminy Conference on Arithmetic. \MR{3363457}

\bibitem[AR79]{AlladiRobinsonTwo}
Krishna Alladi and Michael~L. Robinson, \emph{On certain irrational values of
  the logarithm}, Number theory, {C}arbondale 1979 ({P}roc. {S}outhern
  {I}llinois {C}onf., {S}outhern {I}llinois {U}niv., {C}arbondale, {I}ll.,
  1979), Lecture Notes in Math., vol. 751, Springer, Berlin, 1979, pp.~1--9.
  \MR{564919}

\bibitem[AR80]{AlladiRobinson}
\bysame, \emph{Legendre polynomials and irrationality}, J. Reine Angew. Math.
  \textbf{318} (1980), 137--155. \MR{579389}

\bibitem[AZ90]{AZ}
Gert Almkvist and Doron Zeilberger, \emph{The method of differentiating under
  the integral sign}, J. Symbolic Comput. \textbf{10} (1990), no.~6, 571--591.
  \MR{1087980}

\bibitem[Bak64]{BakerCube}
Alan Baker, \emph{Rational approximations to {$\root 3\of 2$} and other
  algebraic numbers}, Quart. J. Math. Oxford Ser. (2) \textbf{15} (1964),
  375--383. \MR{171750}

\bibitem[Bak22]{Baker}
\bysame, \emph{Transcendental number theory}, Cambridge Mathematical Library,
  Cambridge University Press, Cambridge, 2022, With an introduction by David
  Masser, Reprint of the 1975 original. \MR{4404726}

\bibitem[BB85]{BertrandBeukers}
Daniel Bertrand and Frits Beukers, \emph{\'{E}quations diff\'{e}rentielles
  lin\'{e}aires et majorations de multiplicit\'{e}s}, Ann. Sci. \'{E}cole Norm.
  Sup. (4) \textbf{18} (1985), no.~1, 181--192. \MR{803199}

\bibitem[BB98]{BorweinPiAGM}
Jonathan~M. Borwein and Peter~B. Borwein, \emph{Pi and the {AGM}}, Canadian
  Mathematical Society Series of Monographs and Advanced Texts, vol.~4, John
  Wiley \& Sons, Inc., New York, 1998, A study in analytic number theory and
  computational complexity, Reprint of the 1987 original, A Wiley-Interscience
  Publication. \MR{1641658}

\bibitem[BBR90]{BezivinRobbaAMM}
Frits Beukers, Jean-Paul B\'{e}zivin, and Philippe Robba, \emph{An alternative
  proof of the {L}indemann-{W}eierstrass theorem}, Amer. Math. Monthly
  \textbf{97} (1990), no.~3, 193--197. \MR{1048429}

\bibitem[BC97a]{BombieriCohenGm}
Enrico Bombieri and Paula~B. Cohen, \emph{Effective {D}iophantine approximation
  on {$\bold{G}_m$}. {II}}, Ann. Scuola Norm. Sup. Pisa Cl. Sci. (4)
  \textbf{24} (1997), no.~2, 205--225. \MR{1487954}

\bibitem[BC97b]{BombieriCohenPade}
\bysame, \emph{Siegel's lemma, {P}ad\'{e} approximations and {J}acobians}, Ann.
  Scuola Norm. Sup. Pisa Cl. Sci. (4) \textbf{25} (1997), no.~1-2, 155--178,
  With an appendix by Umberto Zannier, Dedicated to Ennio De Giorgi.
  \MR{1655513}

\bibitem[BC22]{BostCharles}
Jean-Beno\^{i}t Bost and Fran\c{c}ois Charles, \emph{Projective and
  formal-analytic arithmetic surfaces}, 2022,
  \url{https://arxiv.org/abs/2206.14242v2}.

\bibitem[BCY04]{BertrandChirskiiYebbou}
Daniel Bertrand, Vladimir Chirskii, and Johan Yebbou, \emph{Effective estimates
  for global relations on {E}uler-type series}, Ann. Fac. Sci. Toulouse Math.
  (6) \textbf{13} (2004), no.~2, 241--260. \MR{2126743}

\bibitem[Ber99]{BertrandExpAndre}
Daniel Bertrand, \emph{On {A}ndr\'{e}'s proof of the {S}iegel--{S}hidlovsky
  theorem}, Colloque {F}ranco-{J}aponais: {T}h\'{e}orie des {N}ombres
  {T}ranscendants ({T}okyo, 1998), Sem. Math. Sci., vol.~27, Keio Univ.,
  Yokohama, 1999, pp.~51--63. \MR{1726524}

\bibitem[Ber12]{Bertrand}
\bysame, \emph{Le th\'{e}or\`eme de {S}iegel--{S}hidlovsky revisit\'{e}},
  Number theory, analysis and geometry, Springer, New York, 2012, pp.~51--67.
  \MR{2867911}

\bibitem[Beu79]{BeukersIntegral}
Frits Beukers, \emph{A note on the irrationality of {$\zeta (2)$} and {$\zeta
  (3)$}}, Bull. London Math. Soc. \textbf{11} (1979), no.~3, 268--272.
  \MR{554391}

\bibitem[Beu81]{BeukersPade}
\bysame, \emph{Pad\'e-approximations in number theory}, Pad\'e{} approximation
  and its applications, {A}msterdam 1980 ({A}msterdam, 1980), Lecture Notes in
  Math., vol. 888, Springer, Berlin-New York, 1981, pp.~90--99. \MR{649087}

\bibitem[Beu84]{BeukersPolylogs}
\bysame, \emph{The values of polylogarithms}, Topics in classical number
  theory, {V}ol. {I}, {II} ({B}udapest, 1981), Colloq. Math. Soc. J\'anos
  Bolyai, vol.~34, North-Holland, Amsterdam, 1984, pp.~219--228. \MR{781140}

\bibitem[Beu87]{Beukers}
\bysame, \emph{Irrationality proofs using modular forms}, Ast\'{e}risque
  (1987), no.~147-148, 271--283, 345, Journ\'{e}es arithm\'{e}tiques de
  Besan\c{c}on (Besan\c{c}on, 1985). \MR{891433}

\bibitem[Beu06]{BeukersFreeness}
\bysame, \emph{A refined version of the {S}iegel--{S}hidlovskii theorem}, Ann.
  of Math. (2) \textbf{163} (2006), no.~1, 369--379. \MR{2195138}

\bibitem[BG06]{BombieriGubler}
Enrico Bombieri and Walter Gubler, \emph{Heights in {D}iophantine {G}eometry},
  Cambridge New Mathematical Monographs, no.~4, Cambridge University Press,
  2006. \MR{2216774}

\bibitem[BGMN05]{StochasticLp}
Franck Barthe, Olivier Gu\'{e}don, Shahar Mendelson, and Assaf Naor, \emph{A
  probabilistic approach to the geometry of the {$l^n_p$}-ball}, Ann. Probab.
  \textbf{33} (2005), no.~2, 480--513. \MR{2123199}

\bibitem[BGS94]{BostGilletSoule}
Jean-Beno\^{\i}t Bost, Henri Gillet, and Christophe Soul\'{e}, \emph{Heights of
  projective varieties and positive {G}reen forms}, J. Amer. Math. Soc.
  \textbf{7} (1994), no.~4, 903--1027. \MR{1260106}

\bibitem[Bin16]{Binyamini}
Gal Binyamini, \emph{Multiplicity estimates: a {M}orse-theoretic approach},
  Duke Math. J. \textbf{165} (2016), no.~1, 95--128. \MR{3450743}

\bibitem[BLM13]{BoucheronLugosiMassart}
St\'{e}phane Boucheron, G\'{a}bor Lugosi, and Pascal Massart,
  \emph{Concentration inequalities}, Oxford University Press, Oxford, 2013, A
  nonasymptotic theory of independence, With a foreword by Michel Ledoux.
  \MR{3185193}

\bibitem[BM80]{BrownawellMasser}
W.~Dale Brownawell and David~W. Masser, \emph{Multiplicity estimates for
  analytic functions. {II}}, Duke Math. J. \textbf{47} (1980), no.~2, 273--295.
  \MR{575898}

\bibitem[BM83]{BombieriMueller}
Enrico Bombieri and Julia Mueller, \emph{On effective measures of irrationality
  for {$\root r\of{a/b}$} and related numbers}, J. Reine Angew. Math.
  \textbf{342} (1983), 173--196. \MR{703487}

\bibitem[Bom81]{BombieriG}
Enrico Bombieri, \emph{On {$G$}-functions}, Recent progress in analytic number
  theory, {V}ol. 2 ({D}urham, 1979), Academic Press, London-New York, 1981,
  pp.~1--67. \MR{637359}

\bibitem[Bom82]{BombieriDyson}
\bysame, \emph{On the {T}hue--{S}iegel--{D}yson theorem}, Acta Math.
  \textbf{148} (1982), 255--296. \MR{666113}

\bibitem[Bom83]{BombieriWeildec}
\bysame, \emph{On {W}eil's ``th\'eor\`eme de d\'ecomposition''}, Amer. J. Math.
  \textbf{105} (1983), no.~2, 295--308. \MR{701562}

\bibitem[Bom93]{BombieriGm}
\bysame, \emph{Effective {D}iophantine approximation on {${\bf G}_m$}}, Ann.
  Scuola Norm. Sup. Pisa Cl. Sci. (4) \textbf{20} (1993), no.~1, 61--89.
  \MR{1215999}

\bibitem[Bor14]{Borel}
\'Emile Borel, \emph{Introduction g\'eom\'etrique \`a quelques th\'eories
  physiques}, 1914.

\bibitem[Bos99]{BostL21}
Jean-Beno\^{\i}t Bost, \emph{Potential theory and {L}efschetz theorems for
  arithmetic surfaces}, Ann. Sci. \'{E}cole Norm. Sup. (4) \textbf{32} (1999),
  no.~2, 241--312. \MR{1681810}

\bibitem[Bos01]{BostFoliations}
\bysame, \emph{Algebraic leaves of algebraic foliations over number fields},
  Publ. Math. Inst. Hautes \'{E}tudes Sci. (2001), no.~93, 161--221.
  \MR{1863738}

\bibitem[Bos04]{BostGerms}
\bysame, \emph{Germs of analytic varieties in algebraic varieties: canonical
  metrics and arithmetic algebraization theorems}, Geometric aspects of {D}work
  theory. {V}ol. {I}, Walter de Gruyter, Berlin, 2004, pp.~371--418.
  \MR{2023294}

\bibitem[Bos20]{BostBook}
\bysame, \emph{Theta invariants of {E}uclidean lattices and
  infinite-dimensional {H}ermitian vector bundles over arithmetic curves},
  Progress in Mathematics, vol. 334, Birkh\"{a}user/Springer, 2020.
  \MR{4180991}

\bibitem[Boy98]{BoydMahlerSurvey}
David~W. Boyd, \emph{Mahler's measure and special values of {$L$}-functions},
  Experiment. Math. \textbf{7} (1998), no.~1, 37--82. \MR{1618282}

\bibitem[BR89]{BR}
Jean-Paul B\'ezivin and Philippe Robba, \emph{A new {$p$}-adic method for
  proving irrationality and transcendence results}, Ann. of Math. (2)
  \textbf{129} (1989), no.~1, 151--160. \MR{979603}

\bibitem[BR01]{BallRivoal}
Keith Ball and Tanguy Rivoal, \emph{Irrationalit\'{e} d'une infinit\'{e} de
  valeurs de la fonction z\^{e}ta aux entiers impairs}, Invent. Math.
  \textbf{146} (2001), no.~1, 193--207. \MR{1859021}

\bibitem[BV89]{BismutVasserot}
Jean-Michel Bismut and \'{E}ric Vasserot, \emph{The asymptotics of the
  {R}ay--{S}inger analytic torsion associated with high powers of a positive
  line bundle}, Comm. Math. Phys. \textbf{125} (1989), no.~2, 355--367.
  \MR{1016875}

\bibitem[BZ19]{magnetic}
David Broadhurst and Wadim Zudilin, \emph{A magnetic double integral}, J. Aust.
  Math. Soc. \textbf{107} (2019), no.~1, 9--25. \MR{3978030}

\bibitem[BZ20]{BrunaultZudilin}
Fran\c{c}ois Brunault and Wadim Zudilin, \emph{Many variations of {M}ahler
  measures---a lasting symphony}, Australian Mathematical Society Lecture
  Series, vol.~28, Cambridge University Press, Cambridge, 2020. \MR{4382435}

\bibitem[BZ22]{BrownZudilin}
Francis Brown and Wadim Zudilin, \emph{On cellular rational approximations to
  $\zeta(5)$}, 2022, \url{https://arxiv.org/abs/2210.03391v2}.

\bibitem[Car54]{Caratheodory}
Constantin Carath\'eodory, \emph{Theory of {F}unctions of a {C}omplex
  {V}ariable. {V}ol. 2}, Chelsea Publishing Co., New York, 1954, Translated by
  F. Steinhardt. \MR{0064861}

\bibitem[Cat82]{Catalan}
Eug\`ene~Charles Catalan, \emph{Recherches sur la constante~{$G$}, et sur les
  int\'egrales eul\'eriennes}, M\'emoires de l'Acad\'emie Imp\'eriale des
  Sciences de St. P\'etersbourg \textbf{Tome XXXI} (1882), 55 p.

\bibitem[CC83]{Chudnovsky4}
David~V. Chudnovsky and Gregory~V. Chudnovsky, \emph{Rational approximations to
  solutions of linear differential equations}, Proc. Nat. Acad. Sci. U.S.A.
  \textbf{80} (1983), no.~16, 5158--5162. \MR{714302}

\bibitem[CC85a]{ChudnovskyG}
\bysame, \emph{Applications of {P}ad\'{e} approximations to {D}iophantine
  inequalities in values of {$G$}-functions}, Number theory ({N}ew {Y}ork,
  1983--84), Lecture Notes in Math., vol. 1135, Springer, Berlin, 1985,
  pp.~9--51. \MR{803349}

\bibitem[CC85b]{ChudnovskyAlg}
\bysame, \emph{Applications of {P}ad{\'e} approximations to the {G}rothendieck
  conjecture on linear differential equations}, Number theory ({N}ew {Y}ork,
  1983--84), Lecture Notes in Math., vol. 1135, Springer, Berlin, 1985,
  pp.~52--100. \MR{803350}

\bibitem[CC85c]{ChudnovskyIso}
\bysame, \emph{Pad{\'e} approximations and {D}iophantine geometry}, Proc. Nat.
  Acad. Sci. U.S.A. \textbf{82} (1985), no.~8, 2212--2216. \MR{788857}

\bibitem[CDT21]{UDC}
Frank Calegari, Vesselin Dimitrov, and Yunqing Tang, \emph{The unbounded
  denominators conjecture}, 2021, \url{https://arxiv.org/abs/2109.09040}.

\bibitem[CDT24]{zeta5}
\bysame, \emph{Arithmetic holonomy bounds and the irrationality of the $2$-adic
  $\zeta(5)$}, 2024.

\bibitem[Che09]{Chen}
Huayi Chen, \emph{Maximal slope of tensor product of {H}ermitian vector
  bundles}, J. Algebraic Geom. \textbf{18} (2009), no.~3, 575--603.
  \MR{2496458}

\bibitem[Chu79]{ChudnovskyHermite}
Gregory~V. Chudnovsky, \emph{Formules d'{H}ermite pour les approximants de
  {P}ad\'{e} de logarithmes et de fonctions bin\^{o}mes, et mesures
  d'irrationalit\'{e}}, C. R. Acad. Sci. Paris S\'{e}r. A-B \textbf{288}
  (1979), no.~21, A965--A967. \MR{540368}

\bibitem[Chu80]{ChudnovskyShidlovsky}
\bysame, \emph{Rational and {P}ad\'{e} approximations to solutions of linear
  differential equations and the monodromy theory}, Complex analysis,
  microlocal calculus and relativistic quantum theory ({P}roc. {I}nternat.
  {C}olloq., {C}entre {P}hys., {L}es {H}ouches, 1979), Lecture Notes in Phys.,
  vol. 126, Springer, Berlin-New York, 1980, pp.~136--169. \MR{579747}

\bibitem[Chu83a]{ChudnovskyExt}
\bysame, \emph{Number theoretic applications of polynomials with rational
  coefficients defined by extremality conditions}, Arithmetic and geometry,
  {V}ol. {I}, Progr. Math., vol.~35, Birkh\"{a}user Boston, Boston, MA, 1983,
  pp.~61--105. \MR{717590}

\bibitem[Chu83b]{ChudnovskyThueSiegel}
\bysame, \emph{On the method of {T}hue--{S}iegel}, Ann. of Math. (2)
  \textbf{117} (1983), no.~2, 325--382. \MR{690849}

\bibitem[Chu05]{Chu}
Wenchang Chu, \emph{Harmonic number identities and {H}ermite--{P}ad\'{e}
  approximations to the logarithm function}, J. Approx. Theory \textbf{137}
  (2005), no.~1, 42--56. \MR{2179622}

\bibitem[Coh78]{CohenApery}
Henri Cohen, \emph{D\'emonstration de l'irrationalit\'e de~$\zeta(3)$
  (d'apr\`es {R}. {A}p\'ery)}, S\'eminaire de th\'eorie des nombres de Grenoble
  \textbf{6} (1977--78), VI.1--9.

\bibitem[Coo12]{Cooper}
Shaun Cooper, \emph{Sporadic sequences, modular forms and new series for
  {$1/\pi$}}, Ramanujan J. \textbf{29} (2012), no.~1-3, 163--183. \MR{2994096}

\bibitem[D{\`e}b85]{DebesBombieri}
Pierre D{\`e}bes, \emph{Quelques remarques sur un article de {B}ombieri
  concernant le th\'eor\`eme de d\'ecomposition de {W}eil}, Amer. J. Math.
  \textbf{107} (1985), no.~1, 39--44. \MR{778088}

\bibitem[D{\`e}b86]{DebesG}
\bysame, \emph{{$G$}-fonctions et th\'{e}or\`eme d'irreductibilit\'{e} de
  {H}ilbert}, Acta Arith. \textbf{47} (1986), no.~4, 371--402. \MR{884733}

\bibitem[Den50]{Deny}
Jacques Deny, \emph{Les potentiels d'\'{e}nergie finie}, Acta Math. \textbf{82}
  (1950), 107--183. \MR{36371}

\bibitem[DF87]{DiaconisFreedman}
Persi Diaconis and David Freedman, \emph{A dozen de {F}inetti-style results in
  search of a theory}, Ann. Inst. H. Poincar\'{e} Probab. Statist. \textbf{23}
  (1987), no.~2, 397--423. \MR{898502}

\bibitem[DGS94]{Dwork}
Bernard Dwork, Giovanni Gerotto, and Francis~J. Sullivan, \emph{An introduction
  to {$G$}-functions}, Annals of Mathematics Studies, vol. 133, Princeton
  University Press, Princeton, NJ, 1994. \MR{1274045}

\bibitem[DHKK22]{DavidHirataKohnoKawashima}
Sinnou David, Noriko Hirata-Kohno, and Makoto Kawashima, \emph{Linear forms in
  polylogarithms}, Ann. Sc. Norm. Super. Pisa Cl. Sci. \textbf{XXIII} (2022),
  no.~5, 1447--1490.

\bibitem[Dir37]{Dirichlet}
Peter Gustav~Lejeune Dirichlet, \emph{Beweis des {S}atzes, dass jede
  unbegrenzte arithmetische {P}rogression, deren erstes {G}lied und {D}ifferenz
  ganze {Z}ahlen ohne gemeinschaftlichen {F}actor sind, unendlich viele
  {P}rimzahlen enth\"{a}lt}, 1837.

\bibitem[dlVP49]{dVP}
Charles-Jean de~la Vall\'ee-Poussin, \emph{Le potentiel logarithmique ---
  balayage et repr\'esentation conforme}, 1949.

\bibitem[DT97]{DrmotaTichy}
Michael Drmota and Robert~F. Tichy, \emph{Sequences, discrepancies and
  applications}, Lecture Notes in Mathematics, vol. 1651, Springer-Verlag,
  Berlin, 1997. \MR{0470456}

\bibitem[DV01]{diVizio}
Lucia Di~Vizio, \emph{Sur la th\'eorie g\'eom\'etrique des {$G$}-fonctions.
  {L}e th\'eor\`eme de {C}hudnovsky \`a{} plusieurs variables}, Math. Ann.
  \textbf{319} (2001), no.~1, 181--213. \MR{1812823}

\bibitem[Dwo99]{DworkDiff}
Bernard~M. Dwork, \emph{On the size of differential modules}, Duke Math. J.
  \textbf{96} (1999), no.~2, 225--239. \MR{1666546}

\bibitem[Dys47]{Dyson}
Freeman~John Dyson, \emph{The approximation to algebraic numbers by rationals},
  Acta Math. \textbf{79} (1947), 225--240. \MR{23854}

\bibitem[DZ14]{DaugetZudilin}
Simon Dauguet and Wadim Zudilin, \emph{On simultaneous diophantine
  approximations to {$\zeta(2)$} and {$\zeta(3)$}}, J. Number Theory
  \textbf{145} (2014), 362--387. \MR{3253310}

\bibitem[Ell06]{Ellis}
Richard~S. Ellis, \emph{Entropy, large deviations, and statistical mechanics},
  Classics in Mathematics, Springer-Verlag, Berlin, 2006, Reprint of the 1985
  original. \MR{2189669}

\bibitem[Eul35]{Euler}
Leonhard Euler, \emph{De summis serierum reciprocarum}, 1735.

\bibitem[Fis04]{FischlerApery}
St{\'e}phane Fischler, \emph{Irrationalit\'e{} de valeurs de z\^eta (d'apr\`es
  {A}p\'ery, {R}ivoal, {$\dots$})}, Ast\'erisque (2004), no.~294, vii, 27--62.
  \MR{2111638}

\bibitem[FPLL11a]{PrietoLagomasino}
Ulises Fidalgo~Prieto and Guillermo~Tom\'as L\'{o}pez-Lagomasino,
  \emph{Nikishin systems are perfect}, Constructive Approximation \textbf{34}
  (2011), no.~3, 297--356. \MR{2852293}

\bibitem[FPLL11b]{PrietoLagomasino2}
\bysame, \emph{Nikishin systems are perfect. {T}he case of unbounded and
  touching supports}, J. Approx. Theory \textbf{163} (2011), no.~6, 779--811.
  \MR{2832123}

\bibitem[FR03]{FischlerRivoalPade}
St{\'e}phane Fischler and Tanguy Rivoal, \emph{Approximants de {P}ad\'{e} et
  s\'{e}ries hyperg\'{e}om\'{e}triques \'{e}quilibr\'{e}es}, J. Math. Pures
  Appl. (9) \textbf{82} (2003), no.~10, 1369--1394. \MR{2020926}

\bibitem[FR17]{FischlerRivoal}
\bysame, \emph{On the denominators of the {T}aylor coefficients of
  {$G$}-functions}, Kyushu J. Math. \textbf{71} (2017), no.~2, 287--298.
  \MR{3727222}

\bibitem[FR18]{FischlerRivoalIneq}
\bysame, \emph{Rational approximation to values of {$G$}-functions, and their
  expansions in integer bases}, Manuscripta Math. \textbf{155} (2018), no.~3-4,
  579--595. \MR{3763419}

\bibitem[FR21]{FischlerRivoalO}
\bysame, \emph{Linear independence of values of {$G$}-functions, {II}: outside
  the disk of convergence}, Ann. Math. Qu\'e. \textbf{45} (2021), no.~1,
  53--93. \MR{4229177}

\bibitem[Fug60]{Fuglede}
Bent Fuglede, \emph{The logarithmic potential in higher dimensions}, Mat.-Fys.
  Medd. Danske Vid. Selsk. \textbf{33} (1960), no.~1, 14. \MR{125248}

\bibitem[FW94]{FaltingsWustholz}
Gerd Faltings and Gisbert W\"{u}stholz, \emph{Diophantine approximations on
  projective spaces}, Invent. Math. \textbf{116} (1994), no.~1-3, 109--138.
  \MR{1253191}

\bibitem[FW08]{SelbergIntegral}
Peter~J. Forrester and S.~Ole Warnaar, \emph{The importance of the {S}elberg
  integral}, Bull. Amer. Math. Soc. (N.S.) \textbf{45} (2008), no.~4, 489--534.
  \MR{2434345}

\bibitem[Gal74]{Galochkin3}
Alexander~Ivanovich Galo\v{c}kin, \emph{Lower bounds of polynomials in the
  values of a certain class of analytic functions}, Mat. Sb. (N.S.)
  \textbf{95(137)} (1974), 396--417, 471. \MR{357338}

\bibitem[Gal75]{Galochkin1}
\bysame, \emph{Lower bounds of linear forms of the values of certain
  {$G$}-functions}, Mat. Zametki \textbf{18} (1975), no.~4, 541--552.
  \MR{414499}

\bibitem[Gal96]{Galochkin2}
\bysame, \emph{Lower bounds for linear forms of values of {$G$}-functions},
  Vestnik Moskov. Univ. Ser. I Mat. Mekh. (1996), no.~3, 23--29, 91, in
  Russian. \MR{1445270}

\bibitem[Gol69]{Goluzin}
Gennady~Mikhailovich Goluzin, \emph{Geometric {T}heory of {F}unctions of a
  {C}omplex {V}ariable}, Translations of Mathematical Monographs, Vol. 26,
  American Mathematical Society, Providence, R.I., 1969. \MR{0247039}

\bibitem[Gro07]{Gromov}
Misha Gromov, \emph{Metric structures for {R}iemannian and non-{R}iemannian
  spaces}, {E}nglish ed., Modern Birkh\"{a}user Classics, Birkh\"{a}user
  Boston, Inc., Boston, MA, 2007, Based on the 1981 French original, With
  appendices by M. Katz, P. Pansu and S. Semmes, Translated from the French by
  Sean Michael Bates. \MR{2307192}

\bibitem[GS92]{GilletSoule}
Henri Gillet and Christophe Soul\'{e}, \emph{An arithmetic {R}iemann-{R}och
  theorem}, Invent. Math. \textbf{110} (1992), no.~3, 473--543. \MR{1189489}

\bibitem[Had99]{HadamardProd}
Jacques Hadamard, \emph{Th\'eor\`eme sur les s\'eries enti\`eres}, Acta Math.
  \textbf{22} (1899), no.~1, 55--63. \MR{1554900}

\bibitem[Hat93]{HataDilog}
Masayoshi Hata, \emph{Rational approximations to the dilogarithm}, Trans. Amer.
  Math. Soc. \textbf{336} (1993), no.~1, 363--387. \MR{1147401}

\bibitem[Hat98]{HataLogsOpp}
\bysame, \emph{The irrationality of {$\log(1+1/q)\log(1-1/q)$}}, Trans. Amer.
  Math. Soc. \textbf{350} (1998), no.~6, 2311--2327. \MR{1390038}

\bibitem[Her74]{Hermite}
Charles Hermite, \emph{Sur la fonction exponentielle}, 1874.

\bibitem[Her93]{Hermite2}
\bysame, \emph{Sur la g\'en\'eralisation des fractions continues
  alg\'ebriques}, Annali di Matematica, $2^{\textrm{e}}$ s\'erie, \textbf{XXI}
  (1893), 289--308, Extrait d'une lettre \`a M. Pincherle.

\bibitem[Her17]{HermiteOeuvres}
\bysame, \emph{Oeuvres}, 1917.

\bibitem[Hil62]{Hille}
Einar Hille, \emph{Analytic function theory. {V}ol. {II}}, Introductions to
  Higher Mathematics, Ginn and Company, Boston, Mass.-New York-Toronto, 1962.
  \MR{201608}

\bibitem[Hoe63]{Hoeffding}
Wassily Hoeffding, \emph{Probability inequalities for sums of bounded random
  variables}, J. Amer. Statist. Assoc. \textbf{58} (1963), 13--30. \MR{144363}

\bibitem[HPHP11]{MR2900448}
Khodabakhsh\ Hessami~Pilehrood and Tatiana Hessami~Pilehrood, \emph{Bivariate
  identities for values of the {H}urwitz zeta function and supercongruences},
  Electron. J. Combin. \textbf{18} (2011), no.~2, Paper 35, 30. \MR{2900448}

\bibitem[Jac59]{Jacobi}
Carl Gustav~Jacob Jacobi, \emph{Untersuchungen \"{u}ber die
  {D}ifferentialgleichung der hypergeometrischen {R}eihe}, J. Reine Angew.
  Math. \textbf{56} (1859), 149--165. \MR{1579090}

\bibitem[Jag64]{Jager}
Hendrik Jager, \emph{A multidimensional generalization of the {P}ad\'{e} table.
  {I}--{VI}}, Indag. Math. \textbf{26} (1964), 193--249, Nederl. Akad.
  Wetensch. Proc. Ser. A {\bf 67}.

\bibitem[Kat72]{katzpcurvature}
Nicholas~M. Katz, \emph{Algebraic solutions of differential equations
  ({$p$}-curvature and the {H}odge filtration)}, Invent. Math. \textbf{18}
  (1972), 1--118. \MR{337959}

\bibitem[Kol59]{Kolchin}
Ellis~R. Kolchin, \emph{Rational approximation to solutions of algebraic
  differential equations}, Proc. Amer. Math. Soc. \textbf{10} (1959), 238--244.
  \MR{107641}

\bibitem[KZ01]{KontsevichZagier}
Maxim Kontsevich and Don Zagier, \emph{Periods}, Mathematics unlimited---2001
  and beyond, Springer, Berlin, 2001, pp.~771--808. \MR{1852188}

\bibitem[Lan66]{Lang}
Serge Lang, \emph{Introduction to transcendental numbers}, Addison-Wesley
  Publishing Co., Reading, Mass.-London-Don Mills, Ont., 1966. \MR{214547}

\bibitem[Lan72]{Landkof}
Naum~S. Landkof, \emph{Foundations of modern potential theory}, Die Grundlehren
  der mathematischen Wissenschaften, vol. Band 180, Springer-Verlag, New
  York-Heidelberg, 1972, Translated from the Russian by A. P. Doohovskoy.
  \MR{350027}

\bibitem[Led01]{Ledoux}
Michel Ledoux, \emph{The concentration of measure phenomenon}, Mathematical
  Surveys and Monographs, vol.~89, American Mathematical Society, Providence,
  RI, 2001. \MR{1849347}

\bibitem[Leg94]{Legendre}
Adrien-Marie Legendre, \emph{{E}l\'{e}ments de {G}\'{e}om\'{e}trie}, 1794.

\bibitem[Lin82]{Lindemann}
Ferdinand~von Lindemann, \emph{\"{U}ber die {Z}ahl {$\pi$}}, Math. Ann.
  \textbf{20} (1882), no.~2, 213--225. \MR{1510165}

\bibitem[Lys18]{Lysov}
Vladimir~Genrikhovich Lysov, \emph{On {D}iophantine approximants for the
  product of logarithms}, 2018,
  \url{https://keldysh.ru/papers/2018/prep2018_158.pdf} (in Russian),
  pp.~1--20.

\bibitem[Mah53]{MahlerPi}
Kurt Mahler, \emph{On the approximation of {$\pi$}}, Indag. Math. \textbf{15}
  (1953), 30--42, Nederl. Akad. Wetensch. Proc. Ser. A {\bf 56}. \MR{54660}

\bibitem[Mah68]{MahlerPerfect}
\bysame, \emph{Perfect systems}, Compositio Math. \textbf{19} (1968), 95--166.
  \MR{239099}

\bibitem[Mah76]{MahlerBook}
\bysame, \emph{Lectures on transcendental numbers}, Lecture Notes in
  Mathematics, vol. Vol. 546, Springer-Verlag, Berlin-New York, 1976.
  \MR{491533}

\bibitem[Mah19a]{MahlerLindemann}
\bysame, \emph{Ein {B}eweis der {T}ranszendenz der {$P$}-adischen
  {E}xponentialfunktion}, Doc. Math. (2019), 325--331, Reprint of the
  original~1933 paper. \MR{4605009}

\bibitem[Mah19b]{MahlerLog}
\bysame, \emph{On the approximation of logarithms of algebraic numbers}, Doc.
  Math. (2019), 527--555, Reprint of the original~1953 paper. \MR{4605020}

\bibitem[Mai06]{Maillet}
\'Edmond Maillet, \emph{Introduction \`a la th\'eorie des nombres transcendants
  et des propri\'et\'es arithm\'etiques des fonctions}, 1906.

\bibitem[Mai27]{WMaier}
Wilhelm Maier, \emph{Potenzreihen irrationalen {G}renzwertes}, J. Reine Angew.
  Math. \textbf{156} (1927), 93--148. \MR{1581091}

\bibitem[Mas16]{MasserBook}
David~W. Masser, \emph{Auxiliary polynomials in number theory}, Cambridge
  Tracts in Mathematics, vol. 207, Cambridge University Press, Cambridge, 2016.
  \MR{3497545}

\bibitem[Mat97]{Mattner}
Lutz Mattner, \emph{Strict definiteness of integrals via complete monotonicity
  of derivatives}, Trans. Amer. Math. Soc. \textbf{349} (1997), no.~8,
  3321--3342. \MR{1422615}

\bibitem[Mil82]{Milnor2}
John Milnor, \emph{Hyperbolic geometry: the first 150 years}, Bull. Amer. Math.
  Soc. (N.S.) \textbf{6} (1982), no.~1, 9--24. \MR{634431}

\bibitem[Mil83]{Milnor1}
\bysame, \emph{On polylogarithms, {H}urwitz zeta functions, and the {K}ubert
  identities}, Enseign. Math. (2) \textbf{29} (1983), no.~3-4, 281--322.
  \MR{719313}

\bibitem[Mil92]{MilmanDvoretzky}
Vitali Milman, \emph{Dvoretzky's theorem---thirty years later}, Geom. Funct.
  Anal. \textbf{2} (1992), no.~4, 455--479. \MR{1191569}

\bibitem[MS86]{MilmanSchechtman}
Vitali Milman and Gideon Schechtman, \emph{Asymptotic theory of
  finite-dimensional normed spaces}, Lecture Notes in Mathematics, vol. 1200,
  Springer-Verlag, Berlin, 1986, With an appendix by M. Gromov. \MR{856576}

\bibitem[Nag97]{Nagata2d}
Makoto Nagata, \emph{Regular singularities in {$G$}-function theory}, Analytic
  number theory ({K}yoto, 1996), London Math. Soc. Lecture Note Ser., vol. 247,
  Cambridge Univ. Press, Cambridge, 1997, pp.~321--336. \MR{1694999}

\bibitem[Nai82]{Nair}
Mohan Nair, \emph{A new method in elementary prime number theory}, J. London
  Math. Soc. (2) \textbf{25} (1982), no.~3, 385--391. \MR{657495}

\bibitem[Nes88]{Nesterenko2}
Yuri~V. Nesterenko, \emph{Estimates for the number of zeros of certain
  functions}, New advances in transcendence theory ({D}urham, 1986), Cambridge
  Univ. Press, Cambridge, 1988, pp.~263--269. \MR{972005}

\bibitem[Nes96]{NesterenkoThm}
\bysame, \emph{Modular functions and transcendence questions}, Mat. Sb.
  \textbf{187} (1996), no.~9, 65--96. \MR{1422383}

\bibitem[Nes08]{NesterenkoPadic}
\bysame, \emph{Algebraic independence of {$p$}-adic numbers}, Izv. Ross. Akad.
  Nauk Ser. Mat. \textbf{72} (2008), no.~3, 159--174. \MR{2432756}

\bibitem[Nes16]{Nesterenko}
\bysame, \emph{On {C}atalan's constant}, Trudy Mat. Inst. Steklova \textbf{292}
  (2016), no.~Algebra, Geometriya i Teoriya Chisel, 159--176. \MR{3628459}

\bibitem[Nes19]{NesterenkoMahler}
\bysame, \emph{Mahler and transcendence: {E}ffective constructions in
  transcendental number theory}, Doc. Math. (2019), 123--148. \MR{4604999}

\bibitem[Nev70]{Nevanlinna}
Rolf Nevanlinna, \emph{Analytic {F}unctions}, Die Grundlehren der
  mathematischen Wissenschaften, Springer-Verlag, New York-Berlin, 1970.
  \MR{0164038}

\bibitem[Nie09]{Nielsen}
Niels Nielsen, \emph{Der eulersche dilogarithmus und seine verallgemeinerungen:
  Eine monographie}, Nova Acta Leop. \textbf{90} (1909), 121--212.

\bibitem[Nik80]{Nikishin}
Evgeny~Mikhailovich Niki\v{s}in, \emph{Simultaneous {P}ad\'{e} approximants},
  Mat. Sb. (N.S.) \textbf{113(155)} (1980), no.~4(12), 499--519, 637.
  \MR{602272}

\bibitem[NP01]{NesterenkoPhilippon}
Yuri~V. Nesterenko and Patrice Philippon (eds.), \emph{Introduction to
  algebraic independence theory}, Lecture Notes in Mathematics, vol. 1752,
  Springer-Verlag, Berlin, 2001, With contributions from F. Amoroso, D.
  Bertrand, W. D. Brownawell, G. Diaz, M. Laurent, Yuri V. Nesterenko, K.
  Nishioka, Patrice Philippon, G. R\'{e}mond, D. Roy and M. Waldschmidt.
  \MR{1837822}

\bibitem[NS91]{NikishinSorokin}
Evgeny~Mikhailovich Niki\v{s}in and Vladimir~Nikolaevich Sorokin,
  \emph{Rational approximations and orthogonality}, Translations of
  Mathematical Monographs, vol.~92, American Mathematical Society, Providence,
  RI, 1991, Translated from the Russian by Ralph P. Boas. \MR{1130396}

\bibitem[Osg85]{Osgood4}
Charles~F. Osgood, \emph{Sometimes effective
  {T}hue--{S}iegel-{R}oth--{S}chmidt--{N}evanlinna bounds, or better}, J.
  Number Theory \textbf{21} (1985), no.~3, 347--389. \MR{814011}

\bibitem[Phi91]{Philippon}
Patrice Philippon, \emph{Sur des hauteurs alternatives. {I}}, Math. Ann.
  \textbf{289} (1991), no.~2, 255--283. \MR{1092175}

\bibitem[P{\'o}l23]{PolyaOriginal}
George P{\'o}lya, \emph{Sur les s\'eries enti\`eres a coefficients enti\`ers},
  Proc. London Math. Soc. (2) \textbf{21} (1923), 22--38. \MR{1575353}

\bibitem[P{\'o}l28]{PolyaDet}
\bysame, \emph{\"{U}ber gewisse notwendige {D}eterminantenkriterien f\"{u}r die
  {F}ortsetzbarkeit einer {P}otenzreihe}, Math. Ann. \textbf{99} (1928), no.~1,
  687--706. \MR{1512473}

\bibitem[Pom69]{PommerenkeDet}
Christian Pommerenke, \emph{Hankel determinants and meromorphic functions},
  Mathematika \textbf{16} (1969), 158--166. \MR{257327}

\bibitem[PR21]{PrevostRivoalD}
Marc Pr{\'e}vost and Tanguy Rivoal, \emph{Diagonal convergence of the remainder
  {P}ad{\'e} approximants for the {H}urwitz zeta function}, J. Number Theory
  \textbf{222} (2021), 346--361. \MR{4215819}

\bibitem[Pr{\'e}96]{Prevost}
Marc Pr{\'e}vost, \emph{A new proof of the irrationality of {$\zeta(2)$} and
  {$\zeta(3)$} using {P}ad\'{e} approximants}, J. Comput. Appl. Math.
  \textbf{67} (1996), no.~2, 219--235. \MR{1390181}

\bibitem[PS98]{FeldmanNesterenko}
A.~N. Parshin and I.~R. Shafarevich, \emph{Number theory. {IV}}, Encyclopaedia
  of Mathematical Sciences, vol.~44, Springer-Verlag, Berlin, 1998,
  Transcendental Numbers, {\it by Naum I. Fel'dman and Yuri V. Nesterenko}, A
  translation of {\it Number theory. 4 (Russian)}, Ross. Akad. Nauk, Vseross.
  Inst. Nauchn. i Tekhn. Inform., Moscow, Translation by N. Koblitz.
  \MR{1603604}

\bibitem[PZ84]{PerelliZannier}
Alberto Perelli and Umberto Zannier, \emph{On recurrent mod {$p$} sequences},
  J. Reine Angew. Math. \textbf{348} (1984), 135--146. \MR{733927}

\bibitem[Rie38]{Riesz}
Marcel Riesz, \emph{Int\'egrales de {R}iemann--{L}iouville et potentiels}, Acta
  sci. math. Szeged \textbf{9} (1938), 1--42.

\bibitem[Riv00]{Rivoal1}
Tanguy Rivoal, \emph{La fonction z\^{e}ta de {R}iemann prend une infinit\'{e}
  de valeurs irrationnelles aux entiers impairs}, C. R. Acad. Sci. Paris
  S\'{e}r. I Math. \textbf{331} (2000), no.~4, 267--270. \MR{1787183}

\bibitem[Riv06]{RivoalCatalan}
\bysame, \emph{Nombres d'{E}uler, approximants de {P}ad\'{e} et constante de
  {C}atalan}, Ramanujan J. \textbf{11} (2006), no.~2, 199--214. \MR{2267674}

\bibitem[Riv19]{RivoalSur}
\bysame, \emph{Les {$E$}-fonctions et {$G$}-fonctions de {S}iegel}, 2019,
  p.~77.

\bibitem[Rob68]{robinson}
Raphael~M. Robinson, \emph{An extension of {P}\'{o}lya's theorem on power
  series with integer coefficients}, Trans. Amer. Math. Soc. \textbf{130}
  (1968), 532--543. \MR{219706}

\bibitem[Rot55]{Roth}
Klaus~F. Roth, \emph{Rational approximations to algebraic numbers}, Mathematika
  \textbf{2} (1955), 1--20; corrigendum, 168. \MR{72182}

\bibitem[Roy90]{MR1081274}
Ranjan Roy, \emph{The discovery of the series formula for {$\pi$} by {L}eibniz,
  {G}regory and {N}ilakantha}, Math. Mag. \textbf{63} (1990), no.~5, 291--306.
  \MR{1081274}

\bibitem[RR91]{ApproximateIndependence}
Svetlozar Rachev and Ludgar R\"{u}schendorf, \emph{Approximate independence of
  distributions on spheres and their stability properties}, Ann. Probab.
  \textbf{19} (1991), no.~3, 1311--1337. \MR{1112418}

\bibitem[RT86]{RhinToffin}
Georges Rhin and Philippe Toffin, \emph{Approximants de {P}ad\'{e}
  simultan\'{e}s de logarithmes}, J. Number Theory \textbf{24} (1986), no.~3,
  284--297. \MR{866974}

\bibitem[RV96]{RhinViolaZeta2}
Georges Rhin and Carlo Viola, \emph{On a permutation group related to
  {$\zeta(2)$}}, Acta Arith. \textbf{77} (1996), no.~1, 23--56. \MR{1404975}

\bibitem[RV05]{RhinViolaDilog1}
\bysame, \emph{The permutation group method for the dilogarithm}, Ann. Sc.
  Norm. Super. Pisa Cl. Sci. (5) \textbf{4} (2005), no.~3, 389--437.
  \MR{2185863}

\bibitem[RV19]{RhinViolaDilog2}
\bysame, \emph{Linear independence of 1, {$\rm Li_1$} and {$\rm Li_2$}}, Mosc.
  J. Comb. Number Theory \textbf{8} (2019), no.~1, 81--96. \MR{3864310}

\bibitem[Sal07]{Salikhovlog3}
Vladislav~Kh. Salikhov, \emph{On the irrationality measure of {$\ln3$}}, Dokl.
  Akad. Nauk \textbf{417} (2007), no.~6, 753--755. \MR{2462856}

\bibitem[SB85]{Stienstra}
Jan Stienstra and Frits Beukers, \emph{On the {P}icard--{F}uchs equation and
  the formal {B}rauer group of certain elliptic {$K3$}-surfaces}, Math. Ann.
  \textbf{271} (1985), no.~2, 269--304. \MR{783555}

\bibitem[Sch36]{Schneider}
Theodor Schneider, \emph{\"{U}ber die {A}pproximation algebraischer {Z}ahlen},
  J. Reine Angew. Math. \textbf{175} (1936), 182--192. \MR{1581507}

\bibitem[Sch66]{Schwartz}
Laurent Schwartz, \emph{Th\'{e}orie des distributions}, Publications de
  l'Institut de Math\'{e}matique de l'Universit\'{e} de Strasbourg, vol. IX-X,
  Hermann, Paris, 1966, Nouvelle \'{e}dition, enti\'{e}rement corrig\'{e}e,
  refondue et augment\'{e}e. \MR{209834}

\bibitem[Shi59]{Shidlovsky1959}
Andrei~Borisovich Shidlovskii, \emph{A criterion for algebraic independence of
  the values of a class of entire functions}, Izv. Akad. Nauk SSSR Ser. Mat.
  \textbf{23} (1959), 35--66. \MR{102503}

\bibitem[Shi89]{Shidlovsky}
\bysame, \emph{Transcendental numbers}, De Gruyter Studies in Mathematics,
  vol.~12, Walter de Gruyter \& Co., Berlin, 1989, Translated from the Russian
  by Neal Koblitz, With a foreword by W. Dale Brownawell. \MR{1033015}

\bibitem[Sie21]{SiegelZeitschrift}
Carl~Ludwig Siegel, \emph{Approximation algebraischer {Z}ahlen}, Math. Z.
  \textbf{10} (1921), no.~3-4, 173--213. \MR{1544471}

\bibitem[Sie37]{SiegelThue}
\bysame, \emph{Die {G}leichung {$ax^n$}--{$by^n=c$}}, Math. Ann. \textbf{114}
  (1937), no.~1, 57--68. \MR{1513124}

\bibitem[Sie49]{SiegelBook}
\bysame, \emph{Transcendental {N}umbers}, Annals of Mathematics Studies, vol.
  No. 16, Princeton University Press, Princeton, NJ, 1949. \MR{32684}

\bibitem[Sor96]{SorokinN}
Vladimir~Nikolaevich Sorokin, \emph{On the measure of transcendency of the
  number {$\pi^2$}}, Mat. Sb. \textbf{187} (1996), no.~12, 87--120.
  \MR{1442212}

\bibitem[Sor16]{Sorokin}
\bysame, \emph{On {S}alikhov's integral}, Trans. Moscow Math. Soc. (2016),
  107--126. \MR{3643967}

\bibitem[Sta23]{Pascal}
\emph{In {P}ascal's triangle without the $1$s, what is the sum of squares of
  reciprocals?}, Mathematics Stack Exchange, 2023,
  \url{https://math.stackexchange.com/q/4828596} (version: 2023-12-17).

\bibitem[Sto74]{Stolarsky}
Kenneth Stolarsky, \emph{Algebraic numbers and {D}iophantine approximation},
  Pure and Applied Mathematics, vol. No. 26, Marcel Dekker, Inc., New York,
  1974. \MR{374041}

\bibitem[SZ90]{SchechtmanZinn}
Gideon Schechtman and Joel Zinn, \emph{On the volume of the intersection of two
  {$L^n_p$} balls}, Proc. Amer. Math. Soc. \textbf{110} (1990), no.~1,
  217--224. \MR{1015684}

\bibitem[Thu77]{ThueSelected}
Axel Thue, \emph{Selected mathematical papers of {A}xel {T}hue},
  Universitetsforlaget, Oslo--Bergen--Troms\o, 1977, with an introduction by
  Carl Ludwig Siegel.

\bibitem[Thu82]{Thurston}
William~P. Thurston, \emph{Three-dimensional manifolds, {K}leinian groups and
  hyperbolic geometry}, Bull. Amer. Math. Soc. (N.S.) \textbf{6} (1982), no.~3,
  357--381. \MR{648524}

\bibitem[Thu97]{ThurstonBook}
\bysame, \emph{Three-dimensional geometry and topology. {V}ol. 1}, Princeton
  Mathematical Series, vol.~35, Princeton University Press, Princeton, NJ,
  1997, Edited by Silvio Levy. \MR{1435975}

\bibitem[vdP80]{vdP2}
Alfred van~der Poorten, \emph{Some wonderful formulae {$\ldots \ $}footnotes to
  {A}p\'{e}ry's proof of the irrationality of {$\zeta (3)$}}, S\'{e}minaire
  {D}elange-{P}isot-{P}oitou, 20e ann\'{e}e: 1978/1979. {T}h\'{e}orie des
  nombres, {F}asc. 2 ({F}rench), Secr\'{e}tariat Math., Paris, 1980, pp.~Exp.
  No. 29, 7. \MR{582435}

\bibitem[vdP79]{Apery}
\bysame, \emph{A proof that {E}uler missed{$\ldots $}{A}p\'{e}ry's proof of the
  irrationality of {$\zeta (3)$}}, Math. Intelligencer \textbf{1} (1978/79),
  no.~4, 195--203, An informal report. \MR{547748}

\bibitem[Vio04]{ViolaEulerIntegrals}
Carlo Viola, \emph{The arithmetic of {E}uler's integrals}, Riv. Mat. Univ.
  Parma (7) \textbf{3*} (2004), 119--149. \MR{2128843}

\bibitem[Wan04]{JulieWang}
Julie Tzu-Yueh Wang, \emph{An effective {S}chmidt's subspace theorem over
  function fields}, Math. Z. \textbf{246} (2004), no.~4, 811--844. \MR{2045840}

\bibitem[Wei85]{Weierstrass}
Karl Weierstrass, \emph{Zu {L}indemann's {A}bhandlung: `\"{U}ber die
  {L}udolph'sche {Z}ahl".}, Sitzungsberichte der Königlich Preussischen
  Akademie der Wissen-schaften zu Berlin \textbf{5} (1885), no.~2, 1067--1085.

\bibitem[Wei77]{WeilHeckeLemma}
Andr\'{e} Weil, \emph{Remarks on {H}ecke's lemma and its use}, Algebraic number
  theory ({K}yoto {I}nternat. {S}ympos., {R}es. {I}nst. {M}ath. {S}ci., {U}niv.
  {K}yoto, {K}yoto, 1976), Japan Soc. Promotion Sci., Tokyo, 1977,
  pp.~267--274. \MR{480540}

\bibitem[Wir71]{Wirsing}
Eduard Wirsing, \emph{On approximations of algebraic numbers by algebraic
  numbers of bounded degree}, 1969 {N}umber {T}heory {I}nstitute ({P}roc.
  {S}ympos. {P}ure {M}ath., {V}ol. {XX}, {S}tate {U}niv. {N}ew {Y}ork, {S}tony
  {B}rook, {N}.{Y}., 1969), Proc. Sympos. Pure Math., vol. Vol. XX, Amer. Math.
  Soc., Providence, RI, 1971, pp.~213--247. \MR{319929}

\bibitem[Zag09]{Zagier}
Don Zagier, \emph{Integral solutions of {A}p\'{e}ry-like recurrence equations},
  Groups and symmetries, CRM Proc. Lecture Notes, vol.~47, Amer. Math. Soc.,
  Providence, RI, 2009, pp.~349--366. \MR{2500571}

\bibitem[Zan09]{ZannierLectureNotes}
Umberto Zannier, \emph{Lecture notes on {D}iophantine analysis}, Appunti.
  Scuola Normale Superiore di Pisa (Nuova Serie) [Lecture Notes. Scuola Normale
  Superiore di Pisa (New Series)], vol.~8, Edizioni della Normale, Pisa, 2009,
  With an appendix by Francesco Amoroso. \MR{2517762}

\bibitem[Zan14]{Siegel1929SNS}
Umberto Zannier (ed.), \emph{On some applications of {D}iophantine
  approximations}, Quaderni/Monographs, vol.~2, Edizioni della Normale, Pisa,
  2014, A translation of Carl Ludwig Siegel's ``\"{U}ber einige Anwendungen
  diophantischer Approximationen'' by Clemens Fuchs, With a commentary and the
  article ``Integral points on curves: Siegel's theorem after Siegel's proof''
  by Fuchs and Umberto Zannier. \MR{3309332}

\bibitem[Zha95]{Zhang}
Shouwu Zhang, \emph{Positive line bundles on arithmetic varieties}, J. Amer.
  Math. Soc. \textbf{8} (1995), no.~1, 187--221. \MR{1254133}

\bibitem[Zud01]{Zudilin57911}
Wadim Zudilin, \emph{One of the numbers {$\zeta(5)$}, {$\zeta(7)$},
  {$\zeta(9)$}, {$\zeta(11)$} is irrational}, Uspekhi Mat. Nauk \textbf{56}
  (2001), no.~4(340), 149--150. \MR{1861452}

\bibitem[Zud03]{ZudilinCatalan}
\bysame, \emph{An {A}p\'{e}ry-like difference equation for {C}atalan's
  constant}, Electron. J. Combin. \textbf{10} (2003), Research Paper 14, 10.
  \MR{1975764}

\bibitem[Zud14]{ZudilinHypergeometricTales}
\bysame, \emph{Two hypergeometric tales and a new irrationality measure of
  {$\zeta(2)$}}, Ann. Math. Qu\'e. \textbf{38} (2014), no.~1, 101--117.
  \MR{3249415}

\bibitem[Zud17]{ZudilinDet}
\bysame, \emph{A determinantal approach to irrationality}, Constr. Approx.
  \textbf{45} (2017), no.~2, 301--310. \MR{3619445}

\end{thebibliography}

\end{document}